\documentclass[hidelinks]{article}
\usepackage[left=3cm,right=3cm,top=2cm,bottom=2cm]{geometry} 
\usepackage{amsmath}
\usepackage{amsthm}
\usepackage{amssymb}
\usepackage{dsfont}
\usepackage{stackengine}
\usepackage{graphicx}
\usepackage{appendix}
\usepackage{bm}
\usepackage{calc}
\stackMath
\usepackage{hyperref}
\usepackage[svgnames]{xcolor}
\usepackage[style=numeric,backend=bibtex,natbib=true,maxnames=5]{biblatex}
\numberwithin{equation}{section}
\usepackage{accents}

\usepackage{caption}
\usepackage{float}
\usepackage{orcidlink}

\hypersetup{
	colorlinks = true,
	citecolor = teal, 
	linkcolor=red,
	urlcolor=blue}

\DeclareMathOperator{\Tr}{Tr}

\newcommand{\norm}[1]{\left\lVert#1\right\rVert}

\DeclareMathOperator{\im}{Im}
\DeclareMathOperator{\re}{Re}

\DeclareMathOperator{\Expv}{\mathbb{E}}

\newcommand{\other}[1]{\widetilde{#1}}
\newcommand{\otherhat}[1]{\widehat{#1}}

\newcommand{\I}{i}

\DeclareFontFamily{OT1}{pzc}{}
\DeclareFontShape{OT1}{pzc}{m}{it}{ <-> s*[1.1] pzcmi7t }{}
\DeclareMathAlphabet{\mathpzc}{OT1}{pzc}{m}{it}

\newcommand{\dom}{\mathbb{D}}
\newcommand{\fuldom}{\mathbb{D}^{\mathrm{bulk}}}

\newcommand{\wh}{\widehat}

\newcommand{\indset}[1]{[\![#1]\!]}

\newcommand{\maxK}{K}
\newcommand{\extk}{n}

\newcommand{\bandexp}{{\zeta_0}}
\newcommand{\etaexp}{{\delta_0}}
\newcommand{\sumJ}{{\sum J\cup J^*}}

\newcommand{\tinit}{t_{\mathrm{in}}}
\newcommand{\maxQ}{{Q_{\mathrm{max}}}}

\newcommand{\mom}{\mu}

\newcommand{\crit}{{t_*}}
\newcommand{\QV}{ Q }
\makeatletter
\newcommand{\sbullet}{%
	\hbox{\fontfamily{lmr}\fontsize{.6\dimexpr(\f@size pt)}{0}\selectfont\textbullet}}
\DeclareRobustCommand{\reg}[1]{\accentset{\sbullet}{#1}}
\makeatother

\newcommand{\trless}[1]{\mathring{#1}}
\newcommand\dring[1]{\reg{\phantom{#1}} \hspace{\widthof{~}-\widthof{$#1$}}\reg{\phantom{#1}}\hspace{0pt-\widthof{$#1$}-\widthof{~}} #1 }

\newcommand{\trSet}{\mathfrak{T}}
\newcommand{\trl}{\theta}

\newtheorem{theorem}{Theorem} %
\numberwithin{theorem}{section} 
\newtheorem{lemma}[theorem]{Lemma} %
\newtheorem{Def}[theorem]{Definition} %
\newtheorem{prop}[theorem]{Proposition} %
\newtheorem{claim}[theorem]{Claim} %
\newtheorem{remark}[theorem]{Remark} %
\newtheorem{corollary}[theorem]{Corollary} %
\newtheorem{example}[theorem]{Example} %

\DeclareFieldFormat*{title}{#1}
\AtEveryBibitem{\clearfield{month}}
\AtEveryBibitem{\clearfield{number}}
\AtEveryBibitem{\clearfield{publisher}}
\AtEveryBibitem{\clearfield{url}}
\DeclareFieldFormat*{eprint}{
	\hfil\penalty90\hfilneg\space \ifhyperref{
		\href{https://arxiv.org/abs/#1}{arXiv\addcolon#1} 
	}{
		arXiv\addcolon\nolinkurl{#1}
	} 
}
\DeclareFieldFormat*{doi}{
	\hfil\penalty90\hfilneg\space \ifhyperref{
		\href{https://doi.org/#1}{DOI\addcolon\addnbspace#1}
	}{
		DOI\addcolon\addnbspace\nolinkurl{#1}
	} 
}
\DeclareFontFamily{U}{mathx}{}
\DeclareFontShape{U}{mathx}{m}{n}{<-> mathx10}{}
\DeclareSymbolFont{mathx}{U}{mathx}{m}{n}
\DeclareMathAccent{\widehat}{0}{mathx}{"70}
\DeclareMathAccent{\widecheck}{0}{mathx}{"71}

\newcommand{\pis}[1]{\psi_{#1}^{\mathrm{iso}}}
\newcommand{\pav}[1]{\psi_{#1}^{\mathrm{av}}}

\usepackage{mathrsfs} 

\begin{document}
	\begin{minipage}{0.85\textwidth}	\vspace{0.5cm}
	\end{minipage}
	\begin{center}
		\large\bf The Zigzag Strategy for Random Band Matrices
	\end{center}
	\vspace{0.75cm}
	
	\renewcommand{\thefootnote}{\fnsymbol{footnote}}

	\noindent
	\mbox{}%
	\hfill%
	\begin{minipage}{0.25\textwidth}
		\centering
		{L\'aszl\'o Erd\H{o}s}\footnotemark[1]~\orcidlink{0000-0001-5366-9603}\\
		\footnotesize{\textit{lerdos@ist.ac.at}}
	\end{minipage}
	\hfill%
	\begin{minipage}{0.25\textwidth}
		\centering
		{Volodymyr Riabov}\footnotemark[1]~\orcidlink{0009-0007-4989-7524}\\
		\footnotesize{\textit{vriabov@ist.ac.at}}
	\end{minipage}
	\hfill%
	\mbox{}%
	\footnotetext[1]{Institute of Science and Technology Austria, Am Campus 1, 3400 Klosterneuburg, Austria.  Supported by the ERC Advanced Grant "RMTBeyond" No.~101020331.
	}
	
	\renewcommand*{\thefootnote}{\arabic{footnote}}
	\vspace{0.25cm}
	
	\begin{center}
		\begin{minipage}{0.91\textwidth}\footnotesize{ {\bf Abstract.}} 
			We prove that a very general class of $N\times N$ Hermitian random band matrices
			is in the delocalized phase
			when the band width $W$ exceeds the critical threshold, $W\gg \sqrt{N}$.
			 In this regime,  we show that, in the bulk spectrum,
			 the eigenfunctions are fully delocalized,  the eigenvalues follow the universal  Wigner-Dyson
			statistics, and quantum unique ergodicity holds
			for general diagonal  observables with an optimal convergence  rate. 
			Our results are valid for 
			general variance profiles, arbitrary single entry distributions,  in both real-symmetric and 
			complex-Hermitian symmetry classes. In particular, our work  substantially generalizes 
			the recent breakthrough result of Yau and Yin~\cite{yauyin}, obtained for a specific complex Hermitian Gaussian block band matrix.
			 The main technical input is the optimal  multi-resolvent local laws---both in the averaged and fully isotropic form.  We also generalize the $\sqrt{\eta}$-rule from~\cite{cipolloni2021eigenstate}  to
			exploit the additional effect of traceless observables. 
			Our analysis is based on the zigzag strategy, complemented with a new
			global-scale estimate derived using the static version
			of the master inequalities, while the zig-step and the a priori estimates on the deterministic approximations
			are proven dynamically.
		\end{minipage}
	\end{center}

	\vspace{3mm}
	
	{\small
		\footnotesize{\noindent\textit{Date}: \today}\\
		\footnotesize{\noindent\textit{Keywords and phrases}:  Anderson phase transition,
		Wigner-Dyson universality, Delocalization, Local law 
		}
		
		\vspace{3mm}
		
		\thispagestyle{headings} 
		\normalsize

\begingroup
\hypersetup{linkcolor=black}
\tableofcontents
\endgroup

\section{Introduction}

\subsection{General background}
The celebrated Anderson phase transition  asserts that a disordered quantum system with sufficiently 
complex disorder exhibits two distinct phases, depending on the dimension,  energy, and disorder strength. The \emph{localized} phase is characterized by a dense pure point spectrum, exponentially localized eigenfunctions,
dynamical localization and Poisson statistics of the eigenvalues in the large finite volume approximation---that is,   the neighboring eigenvalues
are independent and their local statistics form a Poisson point process.
In contrast, in the {\it delocalized} phase features an absolutely continuous spectrum, non-normalizable
generalized eigenfunctions supported in the entire space, and diffusive quantum dynamics. 
In particular, {\it quantum unique ergodicity (QUE)} (also called {\it eigenstate thermalization hypothesis (ETH)}
in the physics literature)  holds, asserting   that the quantum expectation of an observable on eigenfunctions is well approximated
by its microcanonical ensemble average. 
Neighboring eigenvalues  exhibit strong correlations in the form of  \emph{ level repulsion}, 
and their local statistics are universal, depending only on the basic
symmetry type of the model. Systems without time reversal symmetry follow the celebrated Wigner-Dyson-Mehta
eigenvalue statistics of the Gaussian Unitary Ensemble (GUE), while  those 
with time reversal symmetry  follow the statistics of 
its real counterpart, the Gaussian Orthogonal Ensemble (GOE).

Random Schrödinger operators are the primary mathematical models for the Anderson phase transition, beginning with 
Anderson's tight-binding model~\cite{Anderson58},
where the phase transition was first discovered.
The Hamilton operator  consists of the  standard discrete 
Laplacian on $\mathbb{Z}^d$ with an on-site i.i.d. random potential. 
This model can be extended to  continuous space $\mathbb{R}^d$, or to include magnetic fields and short range correlations
in the disorder -- none of these generalization change the fundamental features.
It is generally believed that in dimensions $d=1, 2$, the system is always in the localized phase, provided some
randomness is present, while for $d\ge 3$, the  system undergoes   a phase transition: at large disorder 
or 
near the spectral edges (where the density of states is very low), localization occurs, while at small disorder in the bulk spectrum the system is delocalized.

After fifty years of intensive mathematical research, the localized regime is now relatively well understood.
The first result in dimension 
 $d=1$ appeared in the work   by Goldsheid, Molchanov and Pastur \cite{gol1977pure}, followed by the groundbreaking paper of Fröhlich and Spencer
\cite{frohlich1983absence} in $d\ge 3$, which introduced the {\it multi-scale analysis} approach. The alternative {\it fractional moment method} was developed shortly after   by Aizenman and Molchanov 
~\cite{aizenman1993localization}. 
Most of the characteristic features of the localized phase mentioned earlier have since been rigorously established, with the notable exception of the two-dimensional case $d=2$, which remains an open problem. 

In contrast, the delocalized regime for the Anderson model remains unsolved; in fact, even its existence is not rigorously established,
apart from the special case of a tree graph \cite{klein1994absolutely, aizenman2006absolutely, aggarwal2022mobility, aggarwal2025mobility}, 
which corresponds to the infinite-dimensional limit $d=\infty$. 
One might naively expect the small disorder regime to behave as a mere perturbation of the pure Laplacian, which  is exactly solvable. 
However, 
 the corresponding perturbation expansion is highly divergent and requires a rigorous renormalization
analysis to all order, which is currently out of reach (for some partial results in the diffusive  {\it scaling limit}
see \cite{erdHos2008quantum, hernandez2024quantum}).

The main reason why the delocalized regime for the Anderson model is so difficult is that randomness occurs only on the 
diagonal of the Hamilton in the natural position space representation; more  randomness in the   off-diagonal matrix
elements would clearly enhance the spreading of the wave packet. 
Taking this simplification to the extreme leads to the   corresponding completely mean-field 
model: the \emph{Wigner random matrix}, an $N\times N$ matrix with independent entries of identical distribution (up to symmetry).  
Wigner matrices are in the delocalized phase, and over the last 15 years, 
 most hallmarks of the delocalized regime have been rigorously proven for them:  
complete delocalization of eigenvectors \cite{erdHos2009local},
the universality of the Wigner-Dyson-Mehta local statistics \cite{erdHos2011universality, tao2011random}, 
level repulsion \cite{erdHos2010wegner}, and
quantum unique ergodicity \cite{bourgade2017eigenvector, cipolloni2021eigenstate},  complemented by  many further   refinements of these features.
 More recently,   these results have been extended to include spatial inhomogeneity, short-range correlation, sparsity, heavy-tailed distributions, etc.
However, all such models remain mean-field in nature, in the sense that they involve of order $N^2$ roughly comparable
independent or weakly correlated random variables.

To bridge the gap between random Schrödinger operators with $N$ random degrees of freedom
and mean-field random matrices with $N^2$ degrees of freedom, 
 one considers   \emph{random band matrices (RBM)}.
Introduced in the physics literature around 1990 \cite{casati1990behavior, casati1990scaling, feingold1991spectral}, these are random matrices with independent entries supported in a band of width $W$ around the main diagonal (in the simplest $d=1$ case).
The  \emph{bandwidth}   parameter $W$
allows one to interpolate between the tridiagonal random Schrödinger operator ($W\sim 1$)
and the fully mean-field Wigner case ($W\sim N$). 
Based on numerical calculations,  Casati, Molinari and Izrailev~\cite{casati1990scaling} conjectured  that this model 
undergoes an Anderson transition around the critical bandwidth $W\sim \sqrt{N}$. 
Later, Fyodorov and Mirlin provided non-rigorous supersymmetric (SUSY) arguments \cite{fyodorov1991scaling} supporting this conjecture\footnote{
	It is sometimes called the Fyodorov--Mirlin conjecture, since they gave the first physical argument.
}.
More precisely, for $W\ll N^{1/2}$
the system is expected to be in the localized phase with a {\it localization length} $\ell\sim W^2$,
meaning the eigenfunctions decay exponentially on scale $\ell\ll N$. In the opposite regime, $W \gg N^{1/2}$,
the system is predicted to be fully delocalized, exhibiting Wigner-Dyson-Mehta statistics and quantum unique ergodicity (QUE).
  These statements are understood in the bulk spectrum.  A similar
transition occurs at the spectral edge, but with a different critical scale $W\sim N^{5/6}$, as was fully understood
by Sodin~\cite{sodin2010spectral} via a very sophisticated moment method.
 
 In higher dimensions $d\ge 2$, the configuration space is the discrete $d$-dimensional torus of side length $L$,   so that $N=L^d$,  and the bandwidth $W$ represents the range of the interaction with independent random amplitudes.
 According to~\cite{fyodorov1991scaling}, the conjectured threshold for delocalization in the bulk spectrum is $W\sim \sqrt{\log L}$ in $d=2$, and  $W\sim W_0(d)$ (finite) for $d\ge 3$.  At the spectral edge, the transition occurs at the critical scale  $W\sim L^{1-d/6}$, which has been fully proven in all dimensions in~\cite{liu2023edge}, extending the method of Sodin.  
 
 \subsection{Informal description of the main results}
 In this paper, we focus on the delocalized regime $W\gg \sqrt{N}$ in the bulk spectrum, and prove all its
essential features for  a very general class of   one-dimensional random band matrices (RBMs) on $N$ sites. That is, for  $H=H^*\in \mathbb{C}^{N\times N}$, 
 including both the real-symmetric and complex-Hermitian symmetry classes.  
 In the original spirit of Wigner's 
vision on universality, our goal is to develop robust 
methods that can handle general distributions well beyond Gaussian and general variance profiles
without special structure, with the only essential parameter being the bandwidth $W$. 
We choose the standard normalization $\Expv{H_{ab}}=0$, and  define the   \emph{variance profile matrix}
by $S_{ab}:=\Expv{|H_{ab}|^2}$, which we assume satisfies $\sum_a S_{ab}=1$ for all $b$, and decays
away from the diagonal (in the periodic distance) as 
\begin{equation}\label{eq:Sdecay}
     S_{ab} \le \frac{C}{W} \left(1+ \frac{|a-b|}{W}\right)^{-6}, \quad \forall a, b=1, 2,\ldots, N,
\end{equation}
with some further conditions. 
This normalization ensures that the limiting eigenvalue density is the standard Wigner semicircle law, $\varrho_{sc}(x)=
\frac{1}{2\pi}\sqrt{(4-x^2)_+}$. 

Our main results are optimal \emph{multi-resolvent local laws} in the bulk spectrum, both 
in averaged and isotropic sense, as well as improved local laws in the presence of traceless observables.
Most other features of the delocalized phase will follow as relatively standard corollaries of these local laws.

To describe the result informally, 
we consider general multi-resolvent chains of the form $G_{[1,k]}:= G_1A_1G_2A_2\ldots A_{k-1}G_k$
where $G_i= (H-z_i)^{-1}$ is the resolvent of  $H$ at a spectral parameter $z_i\in \mathbb{C}\setminus\mathbb{R}$ 
and $A_1, A_2, \ldots$
are  general diagonal observable matrices. We prove optimal local laws in the bulk spectrum for such chains
both in averaged and isotropic form:
\begin{equation}\label{eq:ave}
      \Tr \big[ \big(G_{[1,k]} - M_{[1,k]}\big)A_k\big] \lesssim
       \frac{1}{\ell \eta}\times \mbox{Size $\Big( \Tr\big[M_{[1,k]}A_k\big]\Big)$},
\end{equation}
\begin{equation}\label{eq:iso}
       \big\langle \bm u,  \big(G_{[1,k]} - M_{[1,k]} \big)\bm v\big\rangle
     \lesssim \frac{\| \bm u\|\|\bm v\|}{\sqrt{\ell \eta}}\times \mbox{Size $\Big( M_{[1,k]}\Big)$},
\end{equation}
with very high probability, where $M_{[1,k]}$ is the deterministic approximation of the chain that we will precisely identify. 
The notion of \emph{Size}, which we specify later in~\eqref{eq:isoM}--\eqref{eq:aveM}, encompasses the nontrivial spatial structure of the resolvent chain.
Here $\eta\sim |\im z_i|$ is the common imaginary part of the spectral parameters,
 and $\ell = \ell(\eta):= \min\{ N, W/\sqrt{\eta}\}$
is the \emph{localization length}  of the resolvent. Roughly speaking, this means that $G$ exhibits 
an off-diagonal decay on scale $\ell$; i.e.
 that $|G_{ab}|$ is negligible if
$|a-b|\gg \ell$. 
 The estimates~\eqref{eq:ave}--\eqref{eq:iso} hold down to the smallest possible spectral scale, $\eta\gg 1/N$,
slightly above the mean level spacing. Depending on $\eta$, there are two distinct regimes, that we describe next.

For $\eta\gg \eta_\mathrm{c}:=(W/N)^2$ the localization length is smaller than the system size, $\ell\ll N$.
In this regime, the  resolvents entries---and thus the {\it Size}
functional---exhibit off-diagonal decay on scale $\ell$, see~\eqref{eq:polyUps_notime}--\eqref{eq:expUps_notime}
more precisely later.
 The inverse of this threshold $t_\mathrm{c}:=1/\eta_\mathrm{c}= (N/W)^2$ corresponds to
the {\it Thouless time}; the time it takes for the classical random walk with step of order $W$ to equilibrate in the one-dimensional configuration space of size $N$. Our local law thus also describes quantum diffusion.  Observe that the basic condition 
for the delocalized phase,
$W\gg\sqrt{N}$, is equivalent to $t_\mathrm{c}\ll N$, guaranteeing that the quantum
diffusion is effectively equilibrated before the
Heisenberg time (inverse of the mean level spacing in the bulk 
spectrum\footnote{For convenience, we always scale $H$  such that $\| H\|$ remains order one for large $N$.}). Eventually, this physical mechanism results in 
the Anderson delocalization in our setting.

In contrast, when $\eta\ll \eta_\mathrm{c}$, the off-diagonal profile of the resolvent becomes essentially flat, as in the Wigner case, and   ${\it Size}(M)$ assumes its constant mean-field value without any relevant spatial variation. In this
regime, however,  a new feature emerges: 
the {\it Size} of $G_{[1,k]}$ (and its deterministic approximation $M_{[1,k]}$) becomes sensitive  to whether some of
the observables $A_i$ are traceless. We show that each traceless observable reduces {\it Size}
by a small factor $N\sqrt{\eta}/W$, thereby generalizing the {\it $\sqrt{\eta}$-rule} first discovered
in~\cite{cipolloni2021eigenstate} for Wigner matrices.   
As a simple corollary, this yields quantum unique ergodicity
for general diagonal observables $A$ with an optimal error bound of the form
\begin{equation}\label{eq:ETHinformal}
\Big| \langle \bm u, A \bm u \rangle - \frac{1}{N}\Tr A\Big| \lesssim \frac{\sqrt{N}}{W} \|A\|,
\end{equation}
with very high probability, for any normalized bulk eigenvector $\bm u$, satisfying $H\bm u =\lambda\bm u$.
For simplicity, we restrict our discussion to the bulk regime in this paper; however, all our methods extend naturally to the spectral edge as well---see~\cite{cipolloni2023eigenstate} for an application of the zigzag strategy at the edge.  

We now informally describe our approach, before turning to a more detailed discussion of previous results on one-dimensional RBM.

\subsection{Informal description of the main methods}\label{sec:informal}
A fundamental feature of random matrices $H$ is that their resolvents, $G(z)=(H-z)^{-1}$,
tend to become deterministic in the large $N$ limit---even if $\eta=|\im z|$ is very small.
In the ideal situation, this holds almost down to the local eigenvalue spacing, i.e. for  any $\eta\gg 1/N$
in the bulk spectrum.
Precise mathematical results of this type are called \emph{local laws}. When $\eta\sim 1$ they are referred to as \emph{global laws} instead.  Local laws are the cornerstone in the modern
mathematical analysis of the non-invariant random matrix ensembles\footnote{Invariant ensembles
and related models such as $\beta$-log gases are analyzed with very different tools and methods.}
not only because they directly imply several key 
results, such as eigenvalue rigidity and eigenvector delocalization, but more importantly,
they serve as the necessary {\it a priori} input for the general dynamical method
based on Dyson Brownian motion and the Green function comparison method~\cite{erdHos2017dynamical}.

Over the past two decades, both local laws and their proof techniques have developed extensively.    
However, they all relied on the central
idea that the resolvent $G$, and later more generally the resolvent chains $G_{[1,k]}$,  satisfy an approximate self-consistent
equation. Ignoring fluctuations and other error terms, the solution of the corresponding deterministic self-consistent
equation is exactly $M_{[1,k]}$, indicating the heuristic explanation for why $G_{[1,k]}$ is close to $M_{[1,k]}$.
Unfortunately, in the most interesting regimes, this self-consistent equation is unstable, and these
instabilities must be compensated by  additional cancellations in the fluctuating error terms.
We give a more detailed explanation of this core difficulty and the history of its gradual 
resolution in Section~\ref{sec:dev}. Here we 
only mention that in mean-field models (characterized, roughly,
by $\Expv{|H_{ab}|^2}\sim 1/N$) the instability arises from a  single mode, which can be analyzed separately,
often with the help of extra algebraic identities. In other words, 
the corresponding \emph{stability operator} has an order one \emph{spectral gap}, separating 
the single unstable eigenvalue from the rest.  In contrast, for RBMs, there are many unstable modes, the
stability operator is practically \emph{gapless}. 
This fundamental difference is at the heart of the technical challenges in proving local laws for RBMs.  

Another key feature of the equation for $G_{[1,k]}$ is that the fluctuation term in general depends on resolvent chains longer than length  $k$. This leads naturally to a hierarchy of equations involving chains of increasing length, potentially forming an infinite system. Structurally, this is reminiscent of the BBGKY hierarchy in many-body theory.
 As in any mathematically rigorous  treatment of the BBGKY hierarchy, 
truncating the   resolvent-chain hierarchy    is essential. In some simple cases, this can be achieved using
Ward identities,  but in most cases, truncation introduces additional error terms that must be carefully controlled.   
The most  robust (and sophisticated) approach developed so far is to construct   an entire \emph{hierarchy of master inequalities}
that control the fluctuations of resolvent chains of different length in a \emph{self-improving} way: a cruder initial bound
on the fluctuations of practically all chains can  be systematically refined  to yield sharp estimates via an iterative procedure.  
 The full hierarchy---including averaged and 
isotropic quantities in tandem---was first developed for the simplest
Wigner case in~\cite{Cipolloni2022Optimal}, and later refined in~\cite{cipolloni2022rank};
shorter hierarchies were extended to more general models~\cite{Cipolloni2022overlap, cipolloni2023gaussian}.

While these \emph{static} master inequalities robustly 
addressed the truncation problem, the instability and cancellation issue
was still handled on a case-by-case basis, relying heavily on the spectral gap and on model-specific information about
the single unstable  mode. A new \emph{dynamic}  framework,   the \emph{characteristic flow method},  resolved this problem in
a surprisingly elegant way. The idea is to embed
the problem into a flow that simultaneously moves the spectral parameter $z=z_t$ 
along the characteristic equation
  \begin{equation} \label{eq:char_flow1}
	\frac{\mathrm{d}z_t}{\mathrm{d}t} = -\frac{1}{2}z_t - m(z_t),
\end{equation}
and evolves
 the random matrix along the Ornstein-Uhlenbeck (OU) process 
\begin{equation}\label{OU}
\mathrm{d}H_t = -\frac{1}{2} H_t\mathrm{d}t  +\frac{1}{\sqrt{N}} \mathrm{d}\mathfrak{B}_t.
\end{equation}
The tandem of these two flows yields a major cancellation in 
the time derivative of $G_t(z_t)=(H_t-z_t)^{-1}$, 
the effect of which
is equivalent to the hidden cancellation in the static equations---see Section~\ref{sec:dev} for details.
Since $z_t$ approaches the real axis along the trajectories of \eqref{eq:char_flow1}, this mechanism effectively transfers local laws from larger to smaller   spectral scales, connecting the more (relatively) easily provable global laws at $\eta\sim 1$ to the challenging local law regime $\eta \gg 1/N$.  
The truncation problem still compels one to consider the full hierarchy of resolvent chains, leading to a \emph{dynamic} version of the master inequalities.

However, there is a price for   considering  $G_t(z_t)$ instead of $G(z_t)$:
the  OU process introduces a Gaussian component to $H_t$. To remove this component and return to the original ensemble $H$,   the characteristic flow method has to be complemented by a \emph{Green Function comparison (GFT)} 
argument.  Since GFT can typically remove only a small Gaussian component
 at a time, the procedure is applied iteratively in small steps.  
This gives rise to the \emph{zigzag} strategy,
where a short-time characteristic flow (\emph{zig step})  is used to decrease the $\eta$, and a GFT (\emph{zag step}) removes the added Gaussian component.

The first version of the characteristic flow appeared in a paper by Pastur~\cite{pastur},
 who related   the complex Burgers equation to the evolution of the resolvent along the
OU flow. The idea was later revived by Lee and Schnelli in~\cite{lee2015edge} to prove edge universality 
for deformed Wigner matrices and then by
 von Soosten and Warzel in~\cite{Soosten}
where it was used to prove a local law, albeit a non-optimal one.
The full power of the characteristic flow method has gradually been realized in the context of Dyson Brownian motion
 in  \cite{huang2019rigidity, adhikari2020dyson, adhikari2023local, aggarwal2024edge},
 and later for matrix models, e.g., in  \cite{bourgade2021extreme,
 landon2022almost, landon2024single}. Its  combination with a GFT argument
 to prove specific multi-resolvent local laws---that is, the full zigzag strategy---first appeared in \cite{cipolloni2024mesoscopic}
 in the setting of Hermitization of non-Hermitian matrices, where the characteristics were explicitly used to resolve
 the instability problem.
 The full dynamical version of the master inequalities, systematically controlling the fluctuations 
 of resolvent chains of arbitrary length, was introduced in~\cite{cipolloni2023eigenstate}. There, it was used to study
 Wigner matrices---a setting where the averaged master inequalities remain self-consistent, 
  forgoing the analysis of isotropic chains\footnote{Technically speaking, isotropic local laws were also derived in~\cite{cipolloni2023eigenstate} and used in the GFT step, but, due to the control of traceless observables in the optimal Hilbert–Schmidt norm, the isotropic laws followed directly from the averaged ones.
  	}.

The first use of isotropic master inequalities in tandem with the averaged ones, both in their dynamical form, 
appeared first in~\cite{cipolloni2024out}, which also introduced the term \emph{zigzag}.
Since then, the full power of the zigzag method was successfully exploited in several models and regimes that were previously inaccessible 
due to instabilities~\cite{Cipolloni2022overlap, cipolloni2023universality, cipolloni2024eigenvector, erdHos2024eigenstate, erdHos2024cusp}. By now, the zigzag strategy has proved to be  
a very powerful and robust tool for proving local laws across a wide class of mean-field matrix model. 
   
   \medskip
In this paper, we apply the zigzag strategy to RBM, 
 thus demonstrate its effectiveness well beyond the mean-field setting.  
The key difference between 
the mean-field case and RBM is 
that all relevant quantity exhibit a non-trivial spatial dependence, stemming from the basic fact that $|G_{ab}|$
off-diagonally decays on an intermediate length scale $\ell=\ell(\eta)$.  Heuristically, the correct formula for $\Expv{ |G_{ab}|^2}$
was already identified in~\cite{erdHos2013delocalization},
\begin{equation}\label{Th}
    \Expv{ |G_{ab}|^2} \approx \Theta_{ab}: = \Big(\frac{|m|^2 S}{1-|m|^2S} \Big)_{ab}.
\end{equation}
However, owing to the basic  instability of the static self-consistent equation, this formula could not be proven 
in the entire parameter regime $W\gg \sqrt{N}$, $\eta\gg 1/N$. 
Only partial results and suboptimal bounds were previously available.  The main obstacle
is that the operator $1-|m|^2S$ is gapless: besides the trivial almost zero mode ${\bf 1} :=N^{-1/2}(1,1,\ldots, 1)$, for which
$\| (1-|m|^2S){\bf 1} \| = |1-|m|^2| \sim \eta$,
the operator also has many other low-lying modes.  
A central achievement of this work is to show that the  zigzag strategy can effectively handle the cancellations 
required to control   \emph{all}  these unstable   modes simultaneously.  
A key feature of the zigzag strategy is that throughout the analysis, the fluctuations $G_{[1,k]}-M_{[1,k]}$ are controlled by the size of the deterministic approximation $M_{[1,k]}$ times an appropriate small factor. 
Consequently, establishing sharp bounds on $M_{[1,k]}$ is essential. 
In prior works, such deterministic estimates were typically obtained through independent, often model-specific arguments. 
One of the novelties of the present paper is that we incorporate these estimates into the dynamical framework of the zig-step itself,
effectively reducing the delicate $M$-bounds for long chains to understanding the associated linear propagator $\Theta$ from \eqref{Th}---the fundamental two-resolvent object.     
   In other words, $\Theta$ encodes all necessary information about the model, just as physics intuition predicts. 
   Once it is properly analyzed, the zigzag machinery enables us to rigorously establish local laws in the bulk spectrum for general
   RBMs across the entire delocalized regime.  
   In the next section we comment on previous results on random band matrices.

\subsection{History and related results  on random band matrices}\label{sec:historyRBM}
	Given the elegance of the SUSY methods, substantial efforts have been devoted
	to make the argument of Fyodorov and Mirlin \cite{fyodorov1991scaling} rigorous and establish the (de)localization transition. However, very serious technical complications force one to consider quite restricted models: not only the matrix entries must be Gaussian, but the variance profile of the band must be very specific, as first introduced in \cite{disertori2002density} where  the density of states was computed using SUSY. 
	Nevertheless, the precise transition at $W\sim N^{1/2}$ was rigorously established
	for these restricted models through a sequence of outstanding achievements
	by M. Shcherbina  and T. Shcherbina---first for the local two-point correlation function of
the  characteristic polynomial~\cite{shcherbina2018universality}, and later for the two-point eigenvalue correlation function~\cite{shcherbina2021universality}, after several other deep results \cite{shcherbina2014second, shcherbina2014universality, shcherbina2015universality, shcherbina2016transfer, shcherbina2017characteristic, shcherbina2022susy}.
  Rigorous SUSY methods appear less suitable for establishing other key features of delocalization, such as eigenfunction delocalization or resolvent decay. One notable exception is~\cite{bao2017delocalization}, where an optimal local law and eigenvector delocalization were established under the suboptimal condition
  $W\gg N^{6/7}$. Substantially greater progress has been achieved using other methods, both in the localized and delocalized regimes.

	In contrast to the case of random Schrödinger operators, progress on understanding the
	 localized regime in the bulk spectrum of one-dimensional RBMs has been much slower.
	The first result was due to Schenker ($W\ll N^{1/8}$ in \cite{schenker2009eigenvector}),
	 followed by $W\ll N^{1/7}$ in \cite{peled2019wegner},
	 by $W\ll N^{1/4}$ in \cite{cipolloni2024dynamical}
	and \cite{chen2022random}, and more recently the optimal threshold $W\ll \sqrt{N}$
	was announced~\cite{goldstein2022fluctuations}; 
	  all these results were for Gaussian block band matrices\footnote{Schenker's
	result extends to more general distribution but in a more restricted regime $W\ll N^\mu$ for some $\mu>0$.}.

On the delocalized side, progress has been more substantial, mainly due to intensive developments 
in mean-field random matrix theory over the past decade. New tools, such as resolvent methods
and Dyson Brownian motion techniques, have since been adapted to the RBM setting, gradually lowering the exponent $\alpha$
in threshold $W\gg N^\alpha$ for delocalization proofs.
Disregarding additional assumptions\footnote{Some of these results
	directly extend to higher dimensions, here we just list them in $d=1$.},  and just focusing only on the progress in terms of $\alpha$,
we mention: $\alpha = 6/7$ in \cite{erdHos2011general, erdHos2011quantum} 
and \cite{bao2017delocalization}, $\alpha = 5/4$ in  
\cite{erdHos2013delocalization}, and $\alpha = 3/4$ in \cite{BourgadeI, bourgadeII, YYIII}.  The threshold $\alpha=3/4$ appeared to  mark the limit of {\it static methods} based on   approximate self-consistent equations for
the resolvent.
 The emergence of dynamical methods broke this barrier:   
Dubova and Kang~\cite{dubova} first adapted the characteristic flow method to RBMs, achieving $W\gg N^{8/11}$, 
and finally the entire range $W\gg \sqrt{N}$ was covered by Yau and Yin in their recent fundamental paper~\cite{yauyin}, 
which was then extended to cover the edge regime as well~\cite{yang2025edge}.
After presenting our results in Section~\ref{sec:results}, we compare them with those of~\cite{yauyin}, and in Section~\ref{sec:ideas}, we further comment on the comparison of the methods in more details. 
Our work was initiated and developed independently of~\cite{dubova, yauyin}, but there are natural overlaps. 
Here we only remark that, in the broader context of earlier multi-resolvent zigzag proofs such as~\cite{cipolloni2023eigenstate, cipolloni2024out}, 
both~\cite{dubova} and~\cite{yauyin} essentially implement only the zig-step for the
hierarchy of averaged chains in the RBM setting. In contrast, we systematically combine averaged and isotropic chains in both zig and zag steps.
This full power of zigzag is essential to treat general variance profiles, observables
 and general non-Gaussian entry distributions. 

\medskip

So far we focused on  one-dimensional RBM, the main topic of this paper. For completeness, we briefly summarize the best-known results in higher dimensions.
Recall that the ultimate goal is to establish the analogue of the Anderson delocalization conjecture for RBMs:
namely, that for any $d\ge 3$, there exists a threshold $W_0=W_0(d)$
such that the RBM with $W\ge W_0$ is delocalized uniformly for any $N$.
In the translation-invariant Gaussian 
model, the best proven  result is in $d\ge 7$, where delocalization holds for any $W\ge N^\epsilon$,
as shown in a series of remarkable works~\cite{xu2024bulk, yang2021delocalization1, yang2022delocalization2}.
These results are based on a sophisticated version of the
 static expansion method via Feynman diagrams (called the \emph{T-expansion},
originally introduced in~\cite{erdHos2013delocalization}),
combined with a self-consistent  renormalization scheme (called the \emph{sum zero property}). 
The key instability of the static approach is mitigated by increasing the dimension, which explains the  condition $d\ge 7$.

Another class of higher-dimensional RBM models is inspired by the classical Wegner orbital model
\cite{wegner1979disordered}, in which
$W$ additional internal degrees of freedom (called \emph{spin})
are added to each site of the classical Anderson model.
The main role of the spin variable is to artificially enhance the decorrelation of  the wave function along the time evolution.
The resulting $N\times N$ matrix has a block structure, consisting of $n^d\times n^d$ blocks, each of size $W^d\times W^d$,
where $N=L^d$, $L=nW$. 
The diagonal blocks are independent GUE matrices, representing the completely mean-field internal spin dynamics.
The off-diagonal blocks introduce coupling between neighboring sites and are scaled by a small coupling constant  $\lambda$, which can be tuned to study the Anderson transition\footnote{
	This two-parameter setup, involving  $W$ and $\lambda$, helps decouple two physical effects: $W$ enhances decorrelation, while $\lambda$ is the usual disorder parameter
	in the Anderson model. In contrast, these two roles are intermingled in the single-parameter translation-invariant models.
}. 
The precise form of the off-diagonal blocks may vary: they can be Gaussian matrices (as in the original Wegner model), deterministic Laplacians on the discrete $L$-torus (Block Anderson model) or block Laplacians on the discrete $n=L/W$ torus (Anderson orbital model). 
 One feature of block RBMs, compared to the single-parameter translation invariant models, is that the localized regime is well understood via the fractional moment method, which establishes localization for small coupling: specifically, when $\lambda\ll W^{-d/2}$,  the model exhibits all expected features of Anderson localization~\cite{peled2019wegner}.
In the delocalized regime, similar results to the translation-invariant case were recently obtained using the static $T$-expansion method. In particular, delocalization was shown in~\cite{yang2025delocalization} for any dimension  $d\ge 7$, provided $W\ge L^\epsilon$ and $\lambda\gg W^{-d/4}$.
Once again, the restrictions to large dimension and non-optimal thresholds for $\lambda$ originate
from the inherent instability of the underlying static equation.

The recent emergence of dynamical methods has significantly changed the outlook for higher-dimensional models, as these techniques bypass the key instability inherent in all static approaches.
In particular, the method introduced in~\cite{dubova, yauyin} was extended to two dimensions in~\cite{dubova2025delocalization}, where delocalization in the bulk spectrum was proven for Gaussian block band matrices with $W\ge L^\epsilon$.
The edge regime in dimensions $d=1, 2$ was settled in~\cite{yang2025edge}, delocalization was established in the optimal range 
$W\gg L^{1-d/6}$.
Finally, for block RBMs in the weak disorder regime, the optimal threshold $\lambda\gg W^{-d/2}$
was achieved in dimensions $d=1,2$, for any $W\ge L^\epsilon$ in~\cite{truong2025localization}.

\subsection{Notations and conventions} 
Our model introduced precisely in the next Section~\ref{sec:model} will entail several parameters
 whose values are considered fixed throughout the paper and they are independent of the
 two basic parameters; the dimension $N$ of the matrix and the band-width $W$. We will call 
 them {\it model parameters} and they are listed in Definition~\ref{def:model_param} below.
 
 We investigate the regime, where the basic parameters $N$ and $W$ are very large
 (in fact they are always related by $\sqrt{N}\le W\le N$) and even if not explicitly stated, 
 every estimate involving $N$ is understood to hold for $N$ being sufficiently large, $N\ge N_0$, 
 with $N_0$ depending
 on the model parameters.
 
 For positive quantities  $f,g$ we use the notation $f\lesssim g$ to indicate that $f\le Cg$ for
 some constant $C$ whose exact value is irrelevant and it depends only on the model parameters.
 If $f\lesssim g$ and $g\lesssim f$, we denote this fact by $f\sim g$.
For any real number $r$ we use the standard notations
 $\lceil r\rceil := \min \{ n\in  \mathbb{N} \; : \; n\ge r\}$ and
 $\lfloor r\rfloor := \max \{ n\in  \mathbb{N} \; : \; n\le r\}$ for the upper  and lower integer parts of $r$.
 We also set $\langle r \rangle : =  (r^2+1)^{1/2}$. 
For integers $k\le n$, let $\indset{k,n}$ and $\indset{n}$ denote the index sets
\begin{equation}
	\indset{k,n} := [k,n]\cap \mathbb{Z} = \{k,\dots, n\}, \quad \indset{n} := \indset{1,n} = \{1,\dots, n\}.
\end{equation}
All  summations $\sum_a$ with unspecified limits run over $a\in\indset{N}$.
Our band matrices live on the  configuration space $\indset{N}$ that is  canonically identified with the discrete $N$-torus
$\mathbb{T}_N := \mathbb{Z}/(N\mathbb{Z})$. We use the notation  $|x-y|_N$ for the (periodic) distance 
 between $x, y \in \mathbb{T}_N$, for example $|1-N|_N=1$.

Matrices are denoted by uppercase letters (with the sole exception of the matrix $h$ in~\eqref{eq:kronecker}), their
 entries are indexed by lowercase Roman letters. The identity matrix is $I$. We denote vectors by bold-faced lowercase 
Roman letters ${\bm x}, {\bm y}\in\mathbb{C}^N$, their scalar product is
 $\langle {\bm x},{\bm y}\rangle:= \sum_{i} \overline{x}_i y_i$.  
 Vector and matrix norms, $\lVert {\bm x}\rVert$ and $\lVert A\rVert$, indicate the usual Euclidean norm and the corresponding induced matrix norm. 
 We  write $A_{\bm u \bm v}=\langle \bm u, A \bm v\rangle$. For any $N\times N$ matrix $A$ we use the notation $\langle A\rangle:= N^{-1}\mathrm{Tr}  A$ for its normalized trace. 
 
We introduce the concept of {\it with very high probability}  \emph{(w.v.h.p.)} meaning that for any fixed $C>0$, the probability of an $N$-dependent event is bigger than $1-N^{-C}$ for $N\ge N_0(C)$.
  We also introduce the  standard notion of \emph{stochastic domination}: given two families of non-negative random variables
\[
X=\left(X^{(N)}(u) : N\in\mathbb{N}, u\in U^{(N)} \right) \quad \mathrm{and}\quad Y=\left(Y^{(N)}(u) : N\in\mathbb{N}, u\in U^{(N)} \right)
\] 
indexed by $N$ (and possibly some parameter $u$  in a parameter space $U^{(N)}$), 
we say that $X$ is stochastically dominated by $Y$, if for all $\xi, C>0$ we have 
\begin{equation}
	\label{eq:stochdom}
	\sup_{u\in U^{(N)}} \mathbf{P}\left[X^{(N)}(u)>N^\xi  Y^{(N)}(u)\right]\le N^{-C}
\end{equation}
for large enough $N\ge N_0(\xi,C)$. In this case we use the notation $X\prec Y$ or $X= \mathcal{O}_\prec(Y)$.

\bigskip

{\it Acknowledgment.} 
We are exceptionally grateful to Giorgio Cipolloni, an early and integral member of this project since its inception in the summer of 2023.
Giorgio's insights were absolutely indispensable in shaping the conceptual foundations of this work, particularly through his idea of applying the zigzag strategy and multi-resolvent hierarchy to RBM.
Many fundamental realizations and early breakthroughs stem directly from his involvement, and the completion of the project would not have been possible without his initial participation.

\section{Definition of the model, spatial control functions}\label{sec:model}
We now define our random band matrix ensemble.

\begin{Def}[Random band matrix (RBM) with bandwidth $W$] \label{def:RBM} Let $W\in [1,N]$ be such that 
\begin{equation}\label{eq:WN}
   W^2\ge N^{1+\bandexp} 
 \end{equation}
 with some $\bandexp>0$.
Let $S\in \mathbb{R}^{N\times N}$ be a symmetric matrix with non-negative entries, $S_{ab}=S_{ba}\ge 0$
such that
\begin{equation} \label{eq:sumS=1}
	\sum_a S_{ab} =1, \quad b\in \indset{N}, 
\end{equation} 
\begin{equation} \label{eq:S_bound}
	S_{ab} \le \frac{C_W}{W},
\end{equation}
with some constant $C_W$.
An $N\times N$ real symmetric or complex Hermitian matrix
 $H=H^*\in \mathbb{C}^{N\times N}$ is called a \emph{random band matrix (RBM) with variance profile} $S$
 if 
 $$
     H_{ab} = \sqrt{S_{ab}} h_{ab}, \qquad a,b\in \indset{N},
 $$
 where $\{ h_{ab}\; : \; 1\le a\le b\le N\}$ is a collection of independent centered random variables with variance 1 
 with a uniform moment 
 bound, i.e. 
 \begin{equation}\label{eq:all_moments}
   \Expv{h_{ab}}=0, \quad  \Expv{|h_{ab}|^2}=1, \quad \max_{a,b} \Expv{ |h_{ab}|^p} \le \nu_p
 \end{equation}
 holds for any $p\in \mathbb{N}$ with some  constant $\nu_p$. In the complex-Hermitian case we
 additionally  assume\footnote{In Section~\ref{sec:real} we will explain how to remove this condition.}
 that $\Expv{ h_{ab}^2} =0$.
 The exponent $\bandexp$ and constants $C_W$, $\{\nu_p\}_{p\in \mathbb{N}}$ are independent of $N$ and $W$.
\end{Def}

We present two examples for RBMs that have earlier appeared in the literature.

\begin{example}[Translation-invariant RBM]\label{ex:tr_inv}
Our primary example  is the translation-invariant random band matrix with the variance profile given by
\begin{equation}\label{eq:trinvS}
 S_{ab} : = \frac{1}{W} f\Big( \frac{ |a-b|_N}{W}\Big),
\end{equation}
with some non-negative profile function $f: \mathbb{R}\to \mathbb{R}_+$. We assume
 that $\int_{\mathbb{R}} f(x) \mathrm{d}x=1$ and $f$ satisfies the following
decay and matching regularity conditions:
\begin{equation}\label{eq:fdecay}
  f(x) \le \frac{C_f}{\langle x\rangle^{D+2}}, \quad |\partial f(x)| \le  \frac{C_f}{\langle x\rangle^{D+3}},
  \quad |\partial^2 f(x)| \le  \frac{C_f}{\langle x\rangle^{D+4}}, \qquad \forall x\in \mathbb{R},
\end{equation}
with some fixed decay exponent $D\ge 6$.
Note that the condition $\int f=1$ does not exactly guarantee~\eqref{eq:sumS=1}; however, 
$\sum_a S_{ab}$ remains independent of $b$, hence the limiting density of states is still the Wigner semicircle---possibly
rescaled by an irrelevant constant very close to one, which we  can ignore.
\end{example}

\begin{example}[Block Band Matrix] \label{ex:block}
Suppose that $L:=N/W$ is an integer, assume that there is a symmetric,
non-negative $L\times L$ Toeplitz matrix $\sigma$, i.e. $\sigma_{ij}$ depends only on $|i-j|_L$,  such that 
$\sum_i \sigma_{ij}=1$ for all $j\in \indset{L}$,
$\sigma_{ij}=0$ for $|i-j|\le K_1$ with some constant $K_1$ and set
$$
    S_{ab}: = \frac{1}{W} \sigma_{\lceil a/W\rceil, \lceil b/W\rceil}.
$$
The simplest case is $\sigma_{ij}=1/3$ for $|i-j|_L \le 1$ and $\sigma_{ij}=0$ otherwise; this
corresponds to the  variance profile of the model in~\cite{yauyin}. 
\end{example}

The general Definition~\ref{def:RBM} of RBM does not yet enforce banded structure, and further conditions are necessary.
However, these will be formulated in terms of decay and regularity of certain
 resolvent kernels of $S$. Later we will verify 
that the variance matrix $S$ in our two examples satisfy them. Instead of formulating further conditions directly on $S$,
we  present these more general assumptions since they indicate what are really necessary  for our proof.
The conditions may be verified for more general $S$ 
using techniques from heat kernel analysis on metric measure spaces
but we will not investigate them in this paper.

We need a few preparatory notations.
For any $z\in\mathbb{C}\setminus [-2,2]$, let $m(z)\equiv m_{sc}(z)$ be the Stieltjes transform of the standard 
semicircle law $\varrho_{sc}(x)= \frac{1}{2\pi}\sqrt{[4-x^2]_+}$. Then $m(z)$ satisfies the quadratic equation
\begin{equation} \label{eq:Dyson}
- \frac{1}{m(z)} = z+ m(z), \qquad  \bigl(\im m(z)\bigr)(\im z)>0,
\end{equation}
 where the additional constraint that $\im m(z)$ has the same sign as $\im z$
  uniquely defines $m(z)$ out of the two solutions of~\eqref{eq:Dyson}.
 We list some important properties of $m(z)$ that can be checked with direct calculation. For any 
 (small) constant $\kappa$ and (large) constant $C_0$ we have
 \begin{equation}\label{eq:mbound}
    c \eta\le 1- |m(z)|^2 \le C\eta,  \quad  |m(z)^2-1|\ge c, \quad  \biggl\lvert \biggl( \frac{m(z)}{|m(z)|}\biggr)^2-1\biggr\rvert\ge c,
     \quad c\le  |m(z)|< 1, 
      \end{equation}
 for all $z\in \mathbb{C}\setminus\mathbb{R}$ with $|\re z|\le 2-\kappa$, where $\eta:= |\im z|\le C_0$. 
The positive constants $c, C$  depend only on $\kappa$ and $C_0$.

We fix a small $\etaexp>0$,  without loss of generality we assume
that  $\etaexp < \bandexp$, where $\bandexp$ is from~\eqref{eq:WN}.
Throughout the paper we work with spectral parameters $z$ that are bounded,
 well separated away from the edges $\pm2$ of
the semicircular distribution and away from the real line,   $|\im z|\ge N^{-1+\etaexp}$.
 To quantify these constraints, we fix a (large) constant $C_0$ and a (small) constant $\kappa$
and define our \emph{spectral domain} $\fuldom$ as
\begin{equation}\label{def:spectraldomain}
\fuldom \equiv \fuldom_{\kappa, \etaexp, C_0}:=\Big\{ z\in \mathbb{C} \; : \; N^{-1+\etaexp} \le |\im z|\; , \;  |z|\le C_0, \; 
\Big| \Big( \frac{m(z)}{|m(z)|}\Big)^2-1\Big|\ge \kappa \Big\}.
\end{equation}
The condition $| (m/|m|)^2-1|\ge \kappa$ guarantees that $z$ is separated away from $\pm 2$, 
i.e., we are in the \emph{bulk regime}\footnote{Typically the 
	concept of bulk regime is defined directly by $E=\re z\in [-2+\kappa, 2-\kappa]$, however the larger domain
	$z\in \fuldom$ better suits our proof since $\fuldom$ is invariant under the inverse characteristic flow 
	if the inessential upper bound $|z|\le C_0$ is ignored.}. 
In particular, if $z= E+i \eta$ with $\eta\ge N^{-1+\etaexp}$ and $E\in [-2+\kappa, 2-\kappa]$, then $z\in \fuldom$.
Throughout the paper, we consider $z\in\fuldom$, hence we always have $|\im m(z)|\sim 1$ and  $|m(z)|\sim 1$.

For two given spectral parameters $z_1, z_2\in \fuldom$ in the same half plane\footnote{Clearly $\Xi(z_1, z_2) = \Theta (z_1, \bar z_2)$, so the introduction of $\Xi$ may look superfluous if $\Theta$ were defined for any pair $z_1, z_2$, 
 but we prefer to introduce both functions and restrict their arguments to highlight that the behavior of these
 resolvents is very different depending whether the two spectral parameters are on the same half plane or not.},
  i.e. $(\im z_1)(\im z_2)>0$,  we introduce 
   \begin{equation}\label{def:ThetaXigen}
    \Theta(z_1, z_2 ):= \frac{m(z_1)\overline{ m(z_2)} S}{I- m(z_1)\overline{ m(z_2)}S}, 
    \qquad  \Xi(z_1, z_2 ):= \frac{m(z_1)m(z_2) S}{I- m(z_1)m(z_2) S}.
 \end{equation}
In the special $z_1=z_2=z$ case, we set  
 \begin{equation}\label{def:ThetaXi}
    \Theta(z):= \Theta(z,z)= \frac{|m(z)|^2 S}{I- |m(z)|^2 S}, \qquad  \Xi(z):=\Xi(z,z)= \frac{m(z)^2 S}{I- m(z)^2 S}.
  \end{equation}
       Note that all these families of matrices commute. 

      The matrix $\Theta(z)$, first introduced in~\cite{erdHos2013delocalization},
        plays a central role in the  analysis of random band matrices. It was already
      realized in~\cite{erdHos2013delocalization} that
       $\Theta(z)_{ab}$ approximates $|G(z)_{ab}|^2$ 
      in the sense that
      \begin{equation}\label{eq:heurTheta}
       \Theta(z)_{ab}\approx \Expv |G(z)_{ab}|^2 \approx \sum_c S_{ac} |G(z)_{cb}|^2,
      \end{equation}
      i.e. it encompasses\footnote{We point out that without some averaging, the
      relation $\Theta(z)_{ab}\approx |G(z)_{ab}|^2$ does not hold as $|G(z)_{ab}|^2$ itself is a random quantity
      with fluctuation comparable with its size.}  the key off-diagonal behavior of $|G(z)_{ab}|$. 
      Note that the phase of $G(z)_{ab}$ is not captured by $\Theta$; in fact, it is impossible to track this random phase
      precisely, nevertheless its effect can be tracked after some averaging.
      This is called the \emph{fluctuation averaging mechanism} \cite{erdHos2013averaging}
             that has been extensively exploited in
      many works.  For example, without absolute value 
      in the sum in~\eqref{eq:heurTheta}, we have
      $$
         \Xi(z)_{ab} \approx  \sum_c S_{ac} G(z)_{bc}G(z)_{cb},
     $$
and we will see later that $\Theta(z)$ and $\Xi(z)$ have very different properties;  in general $|\Xi_{ab}|\ll |\Theta_{ab}|$
and $\Theta_{ab}$ is supported on a much larger scale away from its diagonal than $\Xi_{ab}$.

 The Perron-Frobenius theorem and~\eqref{eq:sumS=1} show that the top eigenvalue of $S$ is 1
 with normalized eigenvector ${\bm 1}= N^{-1/2}(1,1,\ldots , 1)$, $S{\bm 1}={\bm 1}$,  therefore 
  we have $\mbox{Spec}(S) \subset (-1,1]$. Since $|m|<1$, the denominators in~\eqref{def:ThetaXi} are
  non-singular, hence $\Theta$ and $\Xi$ are well defined. Moreover, the Neumann series,
      $$
        \Theta = \sum_{k=1}^\infty  |m|^{2k} S^k, 
      $$
  converges and shows that $\Theta$ is entry-wise non-negative with 
  \begin{equation}\label{eq:maxT}
     \| \Theta(z) \| = \max \mbox{Spec}(\Theta(z)) = \frac{|m(z)|^2}{1-|m(z)|^2} =\frac{\im m(z)}{\im z},
  \end{equation}
  using again Perron-Frobenius and~\eqref{eq:Dyson} in the very last step.  We immediately obtain
  the important relation, called the \emph{sum rule},  
 \begin{equation} \label{eq:sumTheta}
 	\sum_{a}\Theta_{ab}(z) = \frac{\im m(z)}{\im z}.
 \end{equation}
 Note that for small $\im z>0$,  $\im m(z)$ is basically the semicircular density, in fact
 $\im m(E+\I 0) = \pi \rho_{sc}(E)$ for $|E|\le 2$.

       Informally, our main assumption is that   $\Theta_{ab}(z)$ lives on a certain scale $\ell$, i.e. it 
  has a strong decay if  $|a-b|_N\gg \ell$ and its maximum is of order $\max_{ab} \Theta_{ab}(z)\sim (\ell \eta)^{-1}$,
  consistently with the sum-rule~\eqref{eq:sumTheta}.
  Here $\ell$ is the 
 $\eta$-dependent {\it localization length-scale}  of the resolvent,   defined by
 \begin{equation}\label{eq:ell_def_notime}
    \ell : =\ell(\eta)= \min \Big\{ \frac{W}{\sqrt{\eta }}, N\Big\}, \quad \eta = |\im z|.
 \end{equation}
 We say that the resolvent is in the \emph{localized regime} if $\ell < N$, or equivalently, if $\eta > (W/N)^2$, while the complementary \emph{delocalized regime} is characterized by $\ell = N$.   
In contrast to $\Theta$, the function $\Xi$ is required to live on the scale of order $W$ (like $S$), and it plays only a secondary role in our analysis.

The more general functions $\Theta(z_1, z_2)$,  $\Xi(z_1, z_2)$ have  similar properties: although their entries
are typically not positive, but the Neumann series still converges and the sum rules in the form
\begin{equation} \label{eq:sumgenTheta}
 	\sum_{a}(\Theta(z_1,z_2 ))_{ab} = \frac{m(z_1)\overline{ m(z_2)} }{1- m(z_1)\overline{ m(z_2)}} ,
	\qquad \sum_{a}(\Xi(z_1,z_2 ))_{ab} =\frac{m(z_1)m(z_2)}{1- m(z_1)m(z_2)},
 \end{equation}
hold since they rely only on  $S{\bm 1}={\bm 1}$. 
If $(\im z_1)(\im z_2)>0$, i.e. $z_1, z_2$ are in the same half plane and away from the spectral edges, 
then $\Xi(z_1,z_2 )$ off-diagonally will live on scale $W$. The precise length scale of $\Theta(z_1,z_2 )$ is not essential
but it is not larger than $\max\{ \ell(\eta_1), \ell(\eta_2)\} = \ell(\eta)$, where $\eta:= \min\{ \eta_1, \eta_2\}$.

 To formulate our necessary working
   conditions on $\Theta$ and $\Xi$ precisely, we introduce the concept of {\it admissible control functions} 
 $(\Upsilon_\eta)_{xy}$,  an $\eta$-dependent deterministic matrix indexed by $x,y\in \indset{N}$,
 that will serve as a convenient upper bound\footnote{ We point out that the current $\Upsilon$ slightly differs from
the analogous $\Upsilon$  introduced in Eq. (2.33) of \cite{erdHos2013delocalization} for a similar purpose as now it does not
contain an additional $1/(N\eta)$ term. Therefore, the current $(\Upsilon_\eta)_{xy}$ bounds $(\Theta(z))_{xy}$
but typically does not exactly bound $| (G(z))_{xy}|^2$
(see Eq. (2.32) of \cite{erdHos2013delocalization}) in the tail regime, $|x-y|_N\gg \ell(\eta)$.}  for $(\Theta(z))_{xy}$
together with some consistency relations.
 The reader may think of  $(\Upsilon_\eta)_{xy}$ as essentially supported on 
the regime $|x-y|_N\lesssim \ell:=\ell(\eta)$ with maximal size of $1/(\ell\eta)$, analogously to $\Theta(z)$.
Therefore, in this context, we also refer to $\ell=\ell(\eta)$ and $\eta$ as {\it scale and size parameters}, respectively.

\begin{Def} [Admissible Control Functions] \label{def:adm_ups_notime}
Fix a large constant $C_0$ and small constants $\kappa$, $\etaexp$ and set $\fuldom=\fuldom_{\kappa,\etaexp, C_0}$ from~\eqref{def:spectraldomain}. 
Setting $\eta= |\im z|$ and recalling $\ell=\ell(\eta)$ 
from \eqref{eq:ell_def_notime},  note that $W\lesssim \ell \le N$ for $z\in \fuldom$.
We record the following straightforward monotonicity properties that follow from~\eqref{eq:ell_def_notime}:
\begin{equation}\label{eq:mon_notime}
 \eta_1\le \eta_2 \quad \Longrightarrow \quad  \ell_1\ge \ell_2, \quad 1\le \ell_1\eta_1\le \ell_2\eta_2, \qquad \ell_i : =\ell(\eta_i).
\end{equation}
We say that an $\eta$-dependent family of entry-wise positive symmetric
$N\times N$ matrices $(\Upsilon_\eta)$, $\eta\in [N^{-1+\etaexp}, C_0]$,
 is an \emph{admissible control function} if there exist  (large) constants $C_1, C_2$, a small constant $c_2>0$ 
 such that $\Upsilon_\eta$ satisfies the following properties:
\begin{itemize}
	\item[(i)]\emph{Majoration of $\Theta$ and $\Xi$}.
	\begin{equation} \label{eq:Ups_majorates_notime}
		\big| \bigl(\Theta(z_1, z_2)\bigr)_{xy} \big| \le C_1 (\Upsilon_{\eta})_{xy}, \quad
		 \bigl\lvert(\Xi(z_1, z_2))_{xy}\bigr\rvert \le C_1(\Upsilon_{\eta=1})_{xy} \quad z_1, z_2 \in \fuldom, \quad  x,y\in\indset{N}, 
	\end{equation} 
	with $\eta:= \min \{ |\im z_1|, |\im z_2| \}$. 
	
	\item [(ii)] \emph{Bounds}.
	\begin{equation}\label{eq:Ups_norm_bounds_notime}
		\max_{x,y} (\Upsilon_\eta)_{xy} \le C_1(\ell\eta)^{-1}, \quad \max_x \sum_a (\Upsilon_\eta)_{xa} \le C_1 \eta^{-1}.
	\end{equation}
	Moreover, for some (large) $D' > 0$, we also have the lower bound  for all $\eta\in [N^{-1+\etaexp}, C_0]$
	\begin{equation} \label{eq:ups_lower_bound_notime}
		\min_{x,y}(\Upsilon_\eta)_{xy} \ge N^{-2D'},
	\end{equation}
	and for $\eta\le (W/N)^2$  (delocalized regime) we have
 \begin{equation}\label{eq:Upsilon_deloc_notime}
       (\Upsilon_{\eta})_{xy}\sim  \frac{1}{N\eta}.
 \end{equation}

	\item[(iii)] \emph{Monotonicity}. For  any $N^{-1+\etaexp}\le \eta_1\le \eta_2\le C_0$ 
	\begin{equation} \label{eq:Ups_time_monot_notime} 
		(\Upsilon_{\eta_2})_{xy}  \le C_1   (\Upsilon_{\eta_1})_{xy}, \quad x,y\in\indset{N}.
	\end{equation}  
	
	\item[(iv)] \emph{Triangle and Convolution inequalities}. For  any $N^{-1+\etaexp}\le \eta_1\le \eta_2\le C_0$ 
	and $\ell_j:=\ell(\eta_j)$,   we have
	\begin{equation} \label{eq:triag_notime} 
		\max_a \bigl((\Upsilon_{\eta_2})_{xa} (\Upsilon_{\eta_1})_{ay}\bigr)  \le C_1  (\ell_2\eta_2)^{-1}(\Upsilon_{\eta_1})_{xy} ,  \quad x,y\in\indset{N},
	\end{equation} 
	\begin{equation} \label{eq:convol_notime} 
		\sum_a \sqrt{(\Upsilon_{\eta_2})_{xa} (\Upsilon_{\eta_1})_{ay}} \le C_1\frac{1}{\eta_2}\sqrt{\ell_2\eta_2(\Upsilon_{\eta_1})_{xy}}, \quad x,y\in\indset{N}.
	\end{equation}

	\item [(v)] \emph{Weighted convolution inequalities}.  For  any $N^{-1+\etaexp}\le \eta_1\le \eta_2\le C_0$
	 and $\ell_i:=\ell(\eta_i)$ 
	\begin{equation} \label{eq:suppressed_convol_notime} 
		\sum_{a} \frac{[|a-x|_N+W]\wedge \ell_1}{\ell_1} \sqrt{(\Upsilon_{\eta_2})_{xa} (\Upsilon_{\eta_i})_{ay}}  
		\le C_1 \frac{1}{\eta_2}\frac{\ell_2}{\ell_1} \sqrt{\frac{\ell_2\eta_2}{\ell_i\eta_i}} \sqrt{\ell_1\eta_1(\Upsilon_{\eta_1})_{xy}} , 
	\end{equation}
	for all $x,y\in\indset{N}$ and $i \in {1,2}$. 
		\item[(vi)] \emph{Regularity of $\Theta$: localized regime}.  
		Let $z_1, z_2\in \fuldom$, $\eta:= \min \{ |\im z_1|, |\im z_2| \}$
		and $\ell=\ell(\eta)$. For any 
		$\eta\ge  c_2(W/N)^2$,  we have
	\begin{equation} \label{eq:Theta_regularity_notime}
		\big| (\Theta(z_1, z_2))_{ab} - (\Theta(z_1, z_2))_{ac} \big|  \le C_1  \frac{(|b-c|_N+W)\wedge\ell}{\ell} \bigl((\Upsilon_\eta)_{ab}+(\Upsilon_\eta)_{ac}\bigr),  
 	\end{equation}
 	for all $a,b,c \in \indset{N}$.
 \item[(vii)]  \emph{Regularity of $\Theta$:  delocalized regime}. 
 Let $z_1, z_2\in \fuldom$, then\footnote{The sum is explicitly given in~\eqref{eq:sumgenTheta}.}
 \begin{equation} \label{eq:supercrit_Theta_notime} 
		\max_{ab}\biggl\lvert (\Theta(z_1, z_2))_{ab} - \frac{1}{N}\sum_c  (\Theta(z_1,z_2))_{ac} 
		\biggr\rvert \le C_1\frac{N}{W^2}, \quad a,b \in \indset{N}.
 \end{equation} 
\end{itemize}  
\end{Def}
Note that the choice~\eqref{eq:Upsilon_deloc_notime} in the delocalized regime, where $\ell\sim N$, 
is consistent
with the sum rule~\eqref{eq:sumTheta} and the bound $\Theta_{xy}\lesssim \Upsilon_{xy}$
from~\eqref{eq:Ups_majorates_notime}.
The condition~\eqref{eq:supercrit_Theta_notime}  for $z_1=z_2$ 
expresses the fact that if $\ell\sim N$, equivalently $\eta\lesssim (W/N)^2$, 
then $\Theta(z)_{ab}$ becomes essentially flat at the level $\sim 1/(N\eta)$, 
in agreement with the sum rule~\eqref{eq:sumTheta} using the identity from~\eqref{eq:maxT}.
For much smaller $\eta\ll (W/N)^2$, we even obtain that $\Theta(z)_{ab}$ becomes  constant with a negligible
error, $\Theta(z)_{ab}  = \frac{\im m(z)}{N\im z} (1+o(1))$. 
This is the main property of $\Theta$ in the delocalization regime. 
The estimate~\eqref{eq:supercrit_Theta_notime} is only effective if $\max \{ |\im z_1|, |\im z_2| \} \le C_2(W/N)^2$, otherwise it is weaker than the $\Theta \lesssim (\ell\eta)^{-1}$
bound following from~\eqref{eq:Ups_majorates_notime} and \eqref{eq:Ups_norm_bounds_notime}, since
$(\ell\eta)^{-1} \le N/W^2$ if $\eta\ge (W/N)^2$.

While the concept of admissible control functions is introduced for general scale and size parameters $\ell, \eta$
satisfying~\eqref{eq:mon_notime}, and a large part of our proof works in this generality, for simplicity 
in the rest of the paper we choose  $\ell=\ell(\eta)$  given
by \eqref{eq:ell_def_notime}. It is straightforward to check that~\eqref{eq:mon_notime} is satisfied with this choice.

Finally, we define the set of admissible variance profiles.
\begin{Def}[Admissible variance profile]\label{def:admS}  A symmetric, entry-wise
 non-negative matrix $S$ is an \emph{admissible variance profile} if
it satisfies~\eqref{eq:sumS=1}, \eqref{eq:S_bound} and  there exists a family of admissible control functions $\Upsilon_\eta$
with the properties given in Definition~\ref{def:adm_ups_notime}. 
\end{Def}

We mention two typical examples for $\Upsilon_\eta$. By direct elementary calculations, that we omit for brevity,
 they 
satisfy  conditions $(ii)-(v)$ of Definition \ref{def:adm_ups_notime}, i.e. those that involve only $\Upsilon$
but not $\Theta, \Xi$.

\begin{example}[Profile with Fast Polynomial Decay] \label{lemma:poly_Ups_notime} 
	For any fixed $D \ge 6$, the function $\Upsilon_\eta$, given by
	\begin{equation} \label{eq:polyUps_notime}
		(\Upsilon_\eta)_{xy} := \frac{1}{\ell\eta}\Big[1+ \frac{|x-y|_N}{\ell}\Big]^{-D}, \qquad \;\; \ell=\ell(\eta),
	\end{equation}
	satisfies the conditions  (ii)--(v) of Definition \ref{def:adm_ups_notime}.
\end{example}
 
\begin{example}[Profile with Exponential Decay]\label{lemma:exp_Ups_notime}
	Given  a constant $c_0>0$,  the function $\Upsilon_\eta$ defined as 
	\begin{equation} \label{eq:expUps_notime}
		(\Upsilon_\eta)_{xy} := \frac{1}{\ell\eta}\exp\biggl\{-c_0\frac{|x-y|_N}{\ell}\biggr\} + \frac{\ell}{\eta}N^{-D}, \quad D\ge 2,
		\quad \ell=\ell(\eta),
	\end{equation}
	satisfies the conditions  (ii)--(v) of Definition \ref{def:adm_ups_notime}.
	 provided constant $C$ in \eqref{eq:convol_notime} and \eqref{eq:suppressed_convol_notime} is replaced by $C\log N$. The resulting poly-log factors are irrelevant for the proof and can be absorbed into the $\prec$ notation used later.
\end{example}

Both examples \eqref{eq:polyUps_notime} and \eqref{eq:expUps_notime} indicate that $(\Upsilon_\eta)_{xy}$ is essentially supported near the diagonal $|x-y|_N \lesssim \ell$ with a maximal size of $(\ell\eta)^{-1}$, i.e. it properly mimics the
behavior of $\Theta(z)_{xy}$.
As $\eta$ decreases, the support of $\Upsilon_\eta$ grows, while its maximal value decreases, 
see \eqref{eq:mon_notime}.

Now we present sufficient conditions on the variance profile $S_{xy}$ that ensure 
that the above examples also satisfy  $(i)$, $(vi)$ and $(vii)$ of Definition \ref{def:adm_ups_notime}.
That is, the following proposition identifies a large class of admissible variance profiles. The proof is given in Section~\ref{sec:reg}.
\begin{prop}[Admissible $S$]\label{prop:admS} ~	
	
 \begin{itemize} \item[(i)] If $S$ is a translation-invariant variance profile  of the form~\eqref{eq:trinvS}
with a non-negative function $f$ satisfying~\eqref{eq:fdecay} for some $D\ge6$ as 
described in Example~\ref{ex:tr_inv},  then $S$ is an admissible variance profile
with the admissible control function  given in~\eqref{eq:polyUps_notime} with the same $D$.  If, additionally,  $f$ is compactly
supported, then $S$ is admissible with the control function given in~\eqref{eq:expUps_notime} with some small $c_0$
depending on the size of the support of $f$.
\item[(ii)] If $S$ is a block band matrix described in Example~\ref{ex:block}, then  $S$ is an admissible variance profile.
\end{itemize}
\end{prop}

Examples~\ref{ex:tr_inv} and \ref{ex:block}
indicate the typical situation for admissible variance profiles
as it was already described in~\cite{erdHos2013delocalization}. 
 In both cases the spectrum of
$S$ satisfies
\begin{equation}\label{eq:spectrumS}
      \mbox{Spec}(S) \subset \{ 1\} \cup \big[-1+ \delta, 1- c_1 \big(\tfrac{W}{N}\big)^2\big]
\end{equation}
for some fixed positive constants $\delta, c_1$, i.e. it has a gap of order $(W/N)^2$ below
the top Perron-Frobenius eigenvalue 1. Moreover, $S$ has many eigenvalues below the gap.
 The reader is encouraged to think of $S\approx 
  I+ W^2\Delta$,
where $\Delta\le 0$ is the standard one-dimensional discrete Laplacian on $\indset{N}$ with roughly regularly spaced
eigenvalues $\approx k^2/N^2$, $k=0, 1,2,3, \ldots$. Using the last identity in~\eqref{eq:maxT},   
 we have, for $\Theta=\Theta(z)$,
$$
     \Theta +I =  \frac{1}{I-|m|^2S}  \approx \frac{1}{|m|^2}\frac{1}{\frac{\eta}{|\im m|}  -  W^2 \Delta}.
$$
Ignoring the factors $|\im m|\sim 1$, $|m|\sim 1$ and using the plane waves as eigenfunctions of the Laplacian, we have
$$
  ( \Theta +I)_{xy}\approx \sum_{k\in\mathbb{Z}} \frac{1}{\eta + k^2(W/N)^2} e^{\I k(x-y)},
$$
which, by Fourier transform, gives rise to the exponential decay of the form
$$  
   \Theta_{xy} \sim \frac{1}{\ell\eta} \exp\biggl\{-c\frac{|x-y|_N}{\ell}\biggr\},
$$
with $\ell = W/\sqrt{\eta}$ in the regime where the gap $(W/N)^2$ is smaller than $\eta$
(localized regime with localization scale $\ell$).
In the opposite delocalized regime, $\eta\le (W/N)^2$, we have 
$$
   \Theta_{xy} \sim \frac{1}{N\eta}.
$$
This heuristics explains the choice of the length scale $\ell$  in~\eqref{eq:ell_def_notime}, the 
choice of $\Upsilon_\eta$ in~\eqref{eq:expUps_notime} and the estimate~\eqref{eq:supercrit_Theta_notime}.
The fast polynomial decay~\eqref{eq:polyUps_notime} arises when $S$ itself has a fast 
off-diagonal decay but  is not finite range on scale $W$. 

Although the heuristic calculation above used  the Fourier transform, in order to obtain a similar qualitative behavior
no explicit translation invariance  of $S$ is necessary (neither on scale one as in Example~\ref{ex:tr_inv}
nor on scale $W$ as in Example~\ref{ex:block}).  However, without some form of translation invariance,
estimating the correct off-diagonal behavior of $\Theta_{xy}$, hence guessing the correct form of $\Upsilon_\eta$ 
and checking its properties  becomes more involved and model dependent.

\medskip

We close this discussion by collecting the parameters of the model:

 \begin{Def}[Model parameters]\label{def:model_param}
 We consider the exponents $\bandexp,\etaexp$,
  the constants $C_W$, $\{ \nu_p\}_{p\in \mathbb{N}}$ in Definition~\ref{def:RBM}
 and 
the constants $\kappa,  C_0, C_1, C_2, c_2$ and $D'$ in Definition~\ref{def:adm_ups_notime}, as fixed \emph{model parameters}, independent of  the main parameters $N, W$.  In the rest of the paper there will be several constant whose precise 
value is irrelevant; they all depend on the model parameters only and hence independent of $N$ and $W$.
 \end{Def}
 
 \begin{remark} [Choice of Small Exponents and Large Parameters]
 	Throughout the proof, we work with several small positive exponents that we summarize here for 
	the reader's convenience.
	For clarity, we adopt consistent notation to reflect their relative sizes, though their exact values may vary locally.
 	
 	We fix three global exponents: $\bandexp$ from the bandwidth condition $W^2 \ge N^{1+\bandexp}$ (see~\eqref{eq:WN}), $\etaexp$ from the spectral scale lower bound $\eta \ge N^{-1+\etaexp}$ 
	(see~\eqref{def:spectraldomain}),
	 and $\xi_0$ from the high-probability tolerance in the definition of the stochastic domination~\eqref{eq:stochdom}
	 (we will use this when we  express our local laws in terms of $\prec$ in 
	 Theorems~\ref{th:local_laws}--\ref{th:local_laws_traceless}). 
	 Without loss of generality, we assume that they satisfy
 	$$
 	\xi_0 \ll \etaexp,\qquad  \etaexp\le \bandexp,
 	$$
 	where $\ll$ means smaller by a large but $N$-independent factor, for example, 
	$\xi_0 \le \etaexp/(10\maxK^2)$ with $\maxK$ the maximal chain length.
 	In addition, we use locally chosen exponents $\xi$, $\nu$, and $\nu'$ to fine tune thresholds in 
	stopping times (e.g., in~\eqref{eq:psi_choice} and~\eqref{eq:new_psi_choice}), satisfying
 	$$
 	\nu' \ll \nu\ll \xi \ll \xi_0.
 	$$
 	To prove local laws
	up to a maximal target chain length $K_{\rm target}$, along the proof 
	we need to consider much longer chains; their maximal chain length will be denoted by
	$\maxK$. Its relation to $K_{\rm target}$ depends on the setup.
	In the complex Hermitian case, we have $K\ge \max\{8, 2K_{\rm target}\}$
	while for the real-symmetric case we require $\maxK \ge \max\{ 8/\bandexp^2, 2K_{\rm target}\}$.   
	To prove local laws for chains  of length up to $K_{\rm target}$ with traceless observables, we will need $\maxK \ge 2K_{\rm target}/\etaexp$.  
 	
 	Finally, the number of zigzag steps, denoted by $\maxQ$, satisfies $\maxQ \le 1/\xi_0$.  
 \end{remark}

\section{Main results: Multi-resolvent  local laws and corollaries}\label{sec:results}
 
Fix a length $k\in \mathbb{N}$, a big constant $\other{C}_0$,  and a $k$-tuple of spectral parameters
$\bm z =(z_1,\dots, z_k) \in (\mathbb{C}\backslash\mathbb{R})^k$  such that
\begin{equation} \label{eq:admissible_z}
	\max_{j\in \indset{k}} \{ \eta_j \}\le \other{C}_0 \min_{j\in \indset{k}} \{ \eta_j\},  \quad \eta_j=|\im z_j|,  
\end{equation} 
 i.e. we assume that the imaginary parts of the spectral 
parameters are comparable. 
Let $(A_1,\dots, A_{k-1})$ be a vector of deterministic $N\times N$  diagonal matrices
that we will call \emph{observables}. 
Then, for $p < j \in \indset{k}$, we denote the resolvent chain with $j-p$ observables $A_i$ by
\begin{equation} \label{eq:resolvent_chains_notime} 
	G_{[p,j]}\equiv G_{[p,j]}(\bm z; A_{p}, \dots, A_{j-1}) 
	:= \biggl(\prod_{i = p}^{j-1} G_{i} A_i\biggr)G_{j}, \qquad G_j:= G(z_j) =(H-z_j)^{-1}.
\end{equation}
By convention, $G_{[j,j]} := G_{j}$ for $j \in \indset{k}$. 
Note that, in particular, the  subscript $[p, j]$ indicates the implicit dependence of $G_{[p,j]}$ on the  spectral parameters $(z_{p}, \dots, z_{j})$. 
For $p < j \in \indset{k}$, we denote the corresponding $N\times N$ deterministic approximation matrix $M_{[p, j]}$, 
defined recursively via \eqref{eq:otherM}--\eqref{eq:M_recursion}, as
\begin{equation} \label{eq:Mt_def_notime}
	M_{[p, j]} \equiv M_{[p, j]}(\bm z; A_p, \dots, A_{j-1}) := M(z_{p}, A_p, z_{p+1}, \dots, z_{j-1}, A_{j-1}, z_{j}).
\end{equation}
Similarly, by convention, $M_{[j,j]} := m (z_{j})$ for $j \in \indset{k}$, where $m(z) := m_{\mathrm{sc}}(z)$ is the solution to the Dyson equation \eqref{eq:Dyson}, and we identify $m(z)$ with the diagonal $N\times N$ matrix $m(z) I$. 
Note that $M_{[1,2]}$, the deterministic approximation of $G_1A_1G_2$, is given by
\begin{equation}\label{eq:M12}
    M_{[1,2]} = \bigl(1-m_1m_2\mathcal{S}\bigr)^{-1} \bigl[m_1m_2 A_1 \bigr], \qquad m_i=m(z_i),
\end{equation}
i.e. it is not simply the product of the deterministic approximations of $G_1$ and $G_2$.
Here $\mathcal{S}$ is the \emph{self-energy (super-)operator}, defined for any deterministic matrix $R$ by
\begin{equation}\label{def:mathcalS}
  \mathcal{S}  [R] : =  \Expv \bigl[ H \, R \, H\bigr].
\end{equation}
In Section~\ref{sec:calS} we will explain 
that for our matrices with independent entries $\mathcal{S}[R]$ is always diagonal 
under the additional assumption\footnote{This assumption holds in the complex-Hermitian symmetry
class (see Definition~\ref{def:RBM}). The
 necessary modifications without  $\Expv (h_{ab})^2 =0$, most importantly for the
real-symmetric case, will be explained separately in Section~\ref{sec:real}.} $\Expv (h_{ab})^2 =0$.
 Consequently, 
for diagonal observables $A_i$, the deterministic approximation $M_{[p, j]}$, defined in \eqref{eq:Mt_def_notime} is also a diagonal matrix.

For $G_\#(\,\cdot\,)$ and $M_\#(\,\cdot\,)$ with matching decorations $\#$ and arguments $(\,\cdot\,)$, we use the short-hand notation
\begin{equation}
	(G-M)_\#(\,\cdot\,) := G_\#(\,\cdot\,) - M_\#(\,\cdot\,).
\end{equation}

While our local laws hold for general diagonal  $A_i$'s, a particular class of observables will play
a central role in our analysis and we introduce separate notations for them. 
Once chains with such observables are understood, the extension to general observables
will be somewhat easier to prove. This special class is
when the diagonal of $A$ contains a row of the matrix of variances $S$. 
That is, when $A=S^x$ for some index $x\in\indset{N}$, where
 $S^x$ denotes the diagonal matrix with the entries of the $x$-th row of $S$ on the main diagonal,
\begin{equation} \label{eq:S_obs}
	\bigl(S^x\bigr)_{ab} := \delta_{ab} S_{xa}, \quad a,b\in\indset{N}.
\end{equation}
It turns out that the control functions $(\Upsilon_\eta)_{xy}$ are well suited to estimate resolvent chains
with these special observables, for example
$$
   \Tr \bigl[GS^xG^*S^y\bigr] = \sum_{ab} S_{xa} |G_{ab}|^2 S_{by}
 $$ 
 is well approximated by $(S\Theta)_{xy}$ with an error term $(\ell\eta)^{-1} (\Upsilon_\eta)_{xy}$.

 In order to effectively handle general observables,
 we now extend the definition of the control functions to include general diagonal observables $A_i$'s
as well as to genuinely isotropic objects to cover general test vectors $\bm u, \bm v$,
 and not only coordinate vectors. All definitions
 rely on the original control functions $\Upsilon_\eta$. We first introduce an appropriate norm to measure the
 size of the general observables, then we give the extension of the control function $\Upsilon$:
 \begin{Def}[Observable Norm, General Control Functions] For any $N\times N$ diagonal matrix $A$ we define its norm as
 \begin{equation}\label{eq:tri}
    |\! | \! | A  |\! | \! | := \min\Big\{ \sum_i a_i\; : \; |A_{qq}|\le \sum_i a_i  S_{iq},
     \;\; \forall q\in\indset{N} \Big\},
 \end{equation}
 where the minimum is taken over all $N$-tuples of non-negative numbers $\bm a =(a_1, \ldots, a_N)$.

 Furthermore, for any diagonal matrices $A, B$ and any vectors $\bm u, \bm v\in \mathbb{C}^N$ we define
 \begin{equation}\label{eq:genUp}
 (\Upsilon_\eta)_{\bm u\bm v} \equiv  (\Upsilon_\eta)_{\bm v\bm u} := \sum_{ij} |u_i|^2 (\Upsilon_\eta)_{ij} |v_j|^2,
 \end{equation}
  \begin{equation}\label{eq:genUpA}
 (\Upsilon_\eta)_{\bm u A} \equiv (\Upsilon_\eta)_{A\bm u} := \sum_{ij} |u_i|^2 (\Upsilon_\eta)_{ij} a_j,  
 \qquad (\Upsilon_\eta)_{ AB} \equiv (\Upsilon_\eta)_{BA} := \sum_{ij} a_i (\Upsilon_\eta)_{ij} b_j, 
 \end{equation}
 where the sequence $\{ a_i\} $ is the minimizer\footnote{The minimum clearly exists by compactness; if there
 are several minima, we may take any of them for this definition. If $A=S^x$ for some $x\in\indset{N}$,  then we choose
 the natural $\bm a= \bm e_x$ coordinate vector for definiteness, see Example~\ref{ex:S}.} on  which the minimum in~\eqref{eq:tri} is achieved, 
 and similarly 
 for $\{ b_i\}$. For coordinate vectors $\bm e_x$ we recover the old control functions, e.g. $(\Upsilon_\eta)_{\bm e_x,  \bm e_y}= (\Upsilon_\eta)_{xy}$.
 \end{Def}

Note that, by~\eqref{eq:Ups_norm_bounds_notime}, we have
$$  
   (\Upsilon_\eta)_{\bm u\bm v} \lesssim \frac{\| \bm u\|^2\|\bm v\|^2}{\ell\eta}, \qquad  
   (\Upsilon_\eta)_{\bm u A}\lesssim \frac{\| \bm u\|^2  |\! | \! | A  |\! | \! | }{\ell\eta}, \qquad
    (\Upsilon_\eta)_{ AB }\lesssim \frac{  |\! | \! | A  |\! | \! | \, |\! | \! | B  |\! | \! | }{\ell\eta}.
 $$
However in the localized regime, $\eta\ge (W/N)^2$,
 the general control functions $\Upsilon$ may be much smaller than these trivial upper bounds on them
since they are sensitive to the overlap
of the spatial structure of $\bm u, \bm v, A, B$, while the pure norm bounds are not.
In the delocalized regime, $\eta\le (W/N)^2$, we have
$\ell\sim N$ and thus $(\Upsilon_\eta)_{xy}\sim \frac{1}{N\eta}$
from~\eqref{eq:Upsilon_deloc_notime}, therefore, we have
\begin{equation}\label{eq:Ups_deloc}
 (\Upsilon_\eta)_{\bm u\bm v} \sim \frac{\| \bm u\|^2 \| \bm v\|^2}{N\eta}, \qquad
  (\Upsilon_\eta)_{\bm u  A} \sim \frac{\| \bm u\|^2  |\! | \! | A  |\! | \! | }{N\eta}, \qquad
   (\Upsilon_\eta)_{AB} \sim \frac{ |\! | \! | A  |\! | \! | \; |\! | \! | B|\! | \! |}{N\eta}.
\end{equation}

 \begin{example}\label{ex:S}
	In  the case of the special observable $A=S^{x}$ for some fixed $x\in \indset{N}$,  we have
	\begin{equation}
		\label{eq:UpsSS}
		\Upsilon_{S^x \bm v} = \Upsilon_{x \bm v}, \quad  \Upsilon_{S^x S^y} = \Upsilon_{x y}, \quad  \Upsilon_{S^x B} 
		= \Upsilon_{\bm e_x B} =: \Upsilon_{x B}, \quad \mbox{etc., and}
		\qquad  |\! | \! | S^x  |\! | \! | =1.
	\end{equation} 
\end{example}

We are now ready to formulate our first main result, the multi-resolvent local law in the bulk,
both in isotropic and averaged form, for general  diagonal
observables.

\begin{theorem}[Multi-Resolvent  Local Laws] \label{th:local_laws}
	Let $H$ be a real symmetric or complex Hermitian random band matrix with $\Expv H =0$ and with
	an admissible variance profile $S$. Fix an integer $k$, constants $C_0$,  $\other{C}_0$
	 and two small constants $\kappa,\etaexp > 0$.
	For each $j\in \indset{k}$ let  $z_j\in \fuldom_{\kappa, \etaexp, C_0}$, defined in~\eqref{def:spectraldomain},
	be a spectral parameter
	 satisfying  the compatibility condition~\eqref{eq:admissible_z}.
	 Set $\eta:= \min_j \eta_j$,  $\ell:=\ell(\eta)$ and $\Upsilon:= \Upsilon_\eta$.
	Let $A_j$, $j\in\indset{k}$,  be deterministic diagonal matrices and let $\bm u, \bm v$ be two deterministic vectors.  
	
	Then, the resolvent chain $G_{[1,k]}= G_{[1,k]}(\bm z; A_{1}, \dots, A_{k-1})$, defined in  \eqref{eq:resolvent_chains_notime},
	satisfies the following isotropic and averaged local laws,
	\begin{align}
			\biggl\lvert \bigl((G-M)_{[1,k]}\bigr)_{\bm u\bm v}\biggr\rvert &\prec \frac{1}{(\ell\eta)^{(k-1)/2 }} 
			\sqrt{ \Upsilon_{\bm u A_1}\Upsilon_{A_1A_2} \ldots \Upsilon_{A_{k-1}\bm v}} , \quad &&k\ge 2,  \label{eq:isolaw} \\
			\biggl\lvert \Tr\big[ \bigl((G-M)_{[1,k]}\bigr) A_k \big]\biggr\rvert &\prec \frac{1}{(\ell\eta)^{k/2 }}  
			\sqrt{\Upsilon_{A_1A_2} \ldots \Upsilon_{A_{k-1}A_k}  \Upsilon_{A_kA_1}}, \quad &&k\ge 2, \label{eq:avelaw}
	\end{align}
	where $M_{[1,k]}$ is its corresponding deterministic approximation from~\eqref{eq:Mt_def_notime}. The deterministic term  $M_{[1,k]}$ satisfies the bounds\footnote{Here we use the probabilistic $\prec$ concept 
	for deterministic $M$ objects to absorb
	additional poly-logarithmic factors that we do not wish to follow precisely.}
	\begin{align}
		\bigl\lvert \bigl(M_{[1,k]}\bigr)_{\bm u \bm v}\bigr\rvert 
		&\prec \frac{1}{(\ell\eta)^{k/2-1 }}  
		\sqrt{ \Upsilon_{\bm u A_1}\Upsilon_{A_1A_2} \ldots \Upsilon_{A_{k-1}\bm v} },  \quad &&k\ge 2, \label{eq:isoM}\\
		\bigl\lvert \Tr \big[ M_{[1,k]}  A_k\big]\bigr\rvert &\prec \frac{1}{(\ell\eta)^{k/2-1 }} 
		\sqrt{\Upsilon_{A_1A_2} \ldots \Upsilon_{A_{k-1}A_k}  \Upsilon_{A_kA_1} }, \quad &&k \ge 2 \label{eq:aveM}.
	\end{align}
	
	For the special case $k=1$, the resolvent $G \equiv G_1$ satisfies the local laws,
	\begin{equation}\label{eq:law1}
	\bigl\lvert \bigl(G-m)_{\bm u \bm v}\bigr\rvert \prec \sqrt{\Upsilon_{\bm u\bm v}}, \qquad 
	\bigl\lvert \Tr\big[ (G-m) A_1 \big]\bigr\rvert \prec \frac{1}{\ell\eta}  |\! | \! | A _1|\! | \! |,  
	\end{equation} 
	with the corresponding $M$-bounds given by
  	\begin{equation}\label{eq:M1}
	\bigl\lvert \bigl(m)_{\bm u \bm v}\bigr\rvert = \big| m\langle \bm u, \bm v\rangle \big| \le \| \bm u\|\|\bm v\|,
	\qquad \bigl\lvert \Tr \big[ m  A_1\big]\bigr\rvert  = |m| \bigl\lvert \Tr \big[ A_1\big]\bigr\rvert 
	\le  |\! | \! | A_1 |\! | \! |.
	\end{equation}
	\end{theorem}
	Notice that the fluctuation $G-M$ in 
	 isotropic sense is smaller by a factor of $1/\sqrt{\ell\eta}\ll 1$ than the corresponding $M$-term, while
	for the averaged law this factor is $1/(\ell\eta)$. 
	
 In the next main result 
we exploit an additional smallness factor $(N/W)\sqrt{\eta}$ 
from each traceless observable  $\langle A_j \rangle=0$ in the delocalized regime, $\eta\le (W/N)^2$. 

\begin{theorem}[Multi-Resolvent  Local Laws with Traceless Observables] \label{th:local_laws_traceless}
	Let $H$ be a real symmetric or complex Hermitian random band matrix with $\Expv H =0$ and with
	an admissible variance profile $S$. Fix an integer $k$,  (large) constants $C_0$, $\other{C}_0$
	 and two small constants $\kappa,\etaexp > 0$.
	For each $j\in \indset{k}$ let  $z_j\in \fuldom_{\kappa, \etaexp, C_0}$ be a spectral parameter
	 satisfying  the compatibility condition~\eqref{eq:admissible_z}.
	 Set $\eta:= \min_j \eta_j$, and assume that $\eta\le (W/N)^2$.
	
	Let $A_j$, $j\in\indset{k}$,  be deterministic diagonal matrices and let $\bm u, \bm v$ two deterministic vectors. 
	Let $n$ be the number of traceless matrices among the first $k-1$ observables $A_j$, and let $n^*$ be the number of traceless matrices among all $A_j$, that is
	$$
	     n:= \#\big\{ j\in \indset{k-1} \; : \; \langle A_j \rangle=0 \big\}, \qquad n^*: = \#\big\{ j\in \indset{k} \; : \; \langle A_j \rangle=0 \big\}.
	$$
	
	Then, the resolvent chain $G_{[1,k]}= G_{[1,k]}(\bm z; A_{1}, \dots, A_{k-1})$, defined in  \eqref{eq:resolvent_chains_notime},
	satisfies the following isotropic and averaged local laws, 
	\begin{align}
		\biggl\lvert \bigl((G-M)_{[1,k]}\bigr)_{\bm u\bm v}\biggr\rvert &\prec \frac{1}{(N\eta)^{k-1/2 }}
		\biggl(\frac{N\sqrt{\eta}}{W}\biggr)^{n }
		\|\bm u\| \, \|\bm v\| \prod_{j=1}^{k-1} |\! | \! | A_j  |\! | \! |,   \label{eq:tr_isolaw}\\
		\biggl\lvert \Tr\big[ (G-M)_{[1,k]} A_k \big]\biggr\rvert &\prec \frac{1}{(N\eta)^{k }} 
		\biggl(\frac{N\sqrt{\eta}}{W}\biggr)^{n^* }\prod_{j=1}^{k} |\! | \! | A_j  |\! | \! |, 		 \label{eq:tr_avelaw}
	\end{align}
	where $M_{[1,k]}$ is its corresponding deterministic approximation from~\eqref{eq:Mt_def_notime}.
	The deterministic term  $M_{[1,k]}$ satisfies the bounds
	\begin{align}
		\label{eq:tr_isoM}
		\bigl\lvert \bigl(M_{[1,k]}\bigr)_{\bm u \bm v}\bigr\rvert &\prec \frac{1}{(N\eta)^{k-1 }} 
		 \biggl(\frac{N\sqrt{\eta}}{W}\biggr)^{2 \lceil n/2 \rceil}  \|\bm u\| \, \|\bm v\| \prod_{j=1}^{k-1} |\! | \! | A_j  |\! | \! |,\\ 
		\label{eq:tr_aveM}
		\bigl\lvert \Tr \big[ M_{[1,k]}  A_k\big]\bigr\rvert &\prec \frac{1}{(N\eta)^{k-1 }}
		 \biggl(\frac{N\sqrt{\eta}}{W}\biggr)^{2 \lceil n^*/2 \rceil}\prod_{j=1}^{k} |\! | \! | A_j  |\! | \! |. 
	\end{align}
	\end{theorem}
Note that in regime $\eta\le(W/N)^2$, both the estimates on $M$'s and on the fluctuations $G-M$ 
have   no spatial dependence in accordance with~\eqref{eq:Ups_deloc}.
For the mean-field situation, $W=N$, all estimates in Theorems~\ref{th:local_laws}--\ref{th:local_laws_traceless}
correspond to the optimal bounds in the Wigner case from~\cite{Cipolloni2022Optimal}.

We now formulate several results that either follow directly from the multi-resolvent local laws---such as isotropic delocalization and quantum unique ergodicity---or follow by well-established procedures, such as Wigner–Dyson universality.

\begin{corollary}[Delocalization and Quantum Unique Ergodicity]\label{thm:deloc}
 Let $H$ be a real symmetric or complex Hermitian random band matrix with $\Expv H =0$ and with
an admissible variance profile $S$. Let \mbox{$\{\bm u_j\}_{j=1}^N \subset\mathbb{C}^N$} denote the orthonormal basis of eigenvectors of $H$ corresponding to eigenvalues $\{\lambda_j\}_{j=1}^N$. 
Fix two small constants $\xi, \kappa > 0$. Then, the bulk eigenvectors $\bm u_i$ with $|\lambda_i| \le 2-\kappa$ satisfy
the following:
\begin{itemize}
	\item[(i)] \textbf{Isotropic Delocalization}. For any fixed deterministic vector $\bm v \in \mathbb{C}^N$ with $\norm{\bm v} = 1$, the bound
	\begin{equation} \label{eq:deloc}
		\max_{\substack{i \in \indset{N} \,: \,  |\lambda_i| \le 2-\kappa }} \big| \langle \bm v, \bm u_i \rangle \big| \le \frac{N^\xi}{\sqrt{N}},
	\end{equation}
	holds with very high probability.
	
	\item[(ii)] \textbf{Quantum Unique Ergodicity}. For any fixed deterministic diagonal matrix $A\in\mathbb{C}^{N\times N}$, the estimate
	\begin{equation} \label{eq:QUE}
		\max_{i,j \in \indset{N} \,:\, |\lambda_i|, |\lambda_j| \le 2-\kappa  }\biggl\lvert \bigl\langle \bm u_i, A \,\bm u_j \bigr\rangle -\delta_{ij} \frac{\Tr[A]}{N} \biggr\rvert \le N^{\xi}\sqrt{\frac{N}{W^2}} \frac{|\! | \! | A  |\! | \! |}{N} \le N^{\xi}\sqrt{\frac{N}{W^2}}\norm{A},
	\end{equation}
	holds with very high probability.
\end{itemize}  
\end{corollary}

The proof of~\eqref{eq:deloc}  is standard, it follows  from the imaginary part of the single resolvent isotropic law~\eqref{eq:law1}
at $\eta:=N^{-1+\xi}$, the corresponding $m$-bound, \eqref{eq:M1},  and the spectral decomposition of $\im G$:
$$
   \big| \langle \bm v, \bm u_i \rangle \big|^2\lesssim   \max_{E\in (-2+\kappa, 2-\kappa)} \eta (\im G(E+i\eta) )_{\bm v \bm v}
   \prec N^\xi+ \frac{1}{N}.
$$
The proof of~\eqref{eq:QUE} follows similarly from the two-resolvent averaged local law~\eqref{eq:tr_avelaw}
with  observables $A_1=A_2^*= \trless{A} = A -N^{-1}\Tr[ A]$ in a standard way~\cite{cipolloni2021eigenstate}.

\begin{theorem}[Wigner-Dyson Universality]\label{thm:WD} 
Let $H$ be a real symmetric or complex Hermitian random band matrix with $\Expv H =0$ and with
	an admissible variance profile $S$. Fix any integer $k\in \mathbb{N}$ and a small constant $\kappa > 0$.
Then the  local eigenvalue correlation functions 
of the $H$ around a fixed energy $E \in [-2+\kappa, 2-\kappa]$ 
are universal in the large $N$ limit,  that is\footnote{The proof gives an effective error term of order $O(N^{-c})$ 
with some small $c>0$.}
\begin{equation}
	\int_{\mathbb{R}^k} f(\bm x) \bigl(\rho_{H}^{(k)}-\rho_{\mathrm{GUE/GOE}}^{(k)}\bigr)
	\Big( E+ \frac{\bm x}{N\varrho_{sc}(E)}\Big) \mathrm{d}\bm x = o(1), \quad N \to \infty,
\end{equation}
for any smooth compactly supported test function $f$ of $k$-variables,
 where $\rho_{H}^{(k)}$ denotes the $k$-th order correlation function of the eigenvalues 
 of $H$, and $\bm x := (x_1, \dots, x_k)$.
\end{theorem}
The proof of Theorem~\ref{thm:WD} is given in Section~\ref{sec:WD}. It will  rely on the single resolvent
isotropic and averaged local laws~\eqref{eq:law1} and on
the averaged two-resolvent local law, \eqref{eq:tr_avelaw} for $k= 2$, with both observables being
traceless ($n^*=2$); the gain factor $(N\sqrt{\eta}/W)^{2}$ will be essential for the proof.

\bigskip

We close this section by comparing our results with the recent paper~\cite{yauyin} by Yau and Yin,
in which the optimal threshold $W\gg N^{1/2}$ for delocalization was first rigorously established.

 First, the model we consider is more general in several 
aspects: (i) in addition to the specific block band matrix from~\cite{yauyin} (Example~\ref{ex:block}),
we also treat general translation-invariant variance profiles (Example~\ref{ex:tr_inv}); (ii) we allow for general single-entry distributions
with finite moments, not restricted to the Gaussian case as in~\cite{yauyin};
and (iii) we include both complex Hermitian and real symmetric symmetry classes. 
 
Second, we establish multi-resolvent local laws for chains of arbitrary length, both in the averaged and fully isotropic sense, with optimal error bounds that accurately capture the spatial dependence. We also allow for general diagonal observables.
 In contrast,~\cite{yauyin} focuses on averaged resolvent chains, and spatial dependence is explicitly tracked only for chains of length $k=2$. There are two key reasons why our approach naturally accommodates this greater generality: On the one hand, 
we perform a fully general global law analysis, initializing the characteristic flow at $\eta\sim 1$
level for a general RBM, while~\cite{yauyin} starts from the zero matrix thereby remains entirely within the Gaussian ensemble.  
On the other hand, we perform a general isotropic analysis, 
 which is crucial for executing the GFT argument (our zag step), and ultimately allows one to depart from Gaussian models.
This step requires   controlling isotropic chains with  optimal spatial dependence. 
In~\cite{yauyin}, the specific block-constant variance profile permits special factorizations that turn critical resolvent sums into products of full traces, thus circumventing the need for isotropic estimates. 

Third, we prove optimal local laws  involving any number of traceless observables; among other benefits this gives QUE for general diagonal observables with an optimal error bound that holds with very high probability, as stated in~\eqref{eq:QUE}.  

In the next Section~\ref{sec:ideas} we summarize the main ideas of our proof and at the end of each subsection
we will compare them with the analogous steps in~\cite{yauyin}.

\section{Main ideas of the proofs}\label{sec:ideas}
 
This section outlines the central ideas behind our proofs and situates them within the broader historical development of techniques for establishing local laws. 
In Section~\ref{sec:dev}, we begin by outlining the main technical challenges in proving local laws for mean-field models, along with the historical development of the methods used to overcome them. This serves as a foundation for Section~\ref{sec:ideasRBM}, where we describe how we address the analogous yet more severe obstacles that arise in the case of random band matrices.

\subsection{Development of local laws: Mean-field case}\label{sec:dev}

Local laws and their proof techniques have developed enormously in the past twenty years, 
with each major advancement leading to new insights and consequences. 
Early work focused on the single resolvent $G$ of a Wigner
matrix in the bulk regime, relying on the observation that both its average trace $\langle G\rangle = N^{-1}\Tr [G]$
and its diagonal elements $G_{ii}$ 
asymptotically satisfy a self-consistent relation. This relation is close to the defining equation $- m^{-1} =m+z$ for
the Stieltjes transform $m=m(z)$ of the semicircle law.  
Initially, this relation was obtained from the Schur complement formula. Later, an equation for the difference $G-m$ was derived directly,
which roughly takes the form
\begin{equation}\label{eq:G-m}
    \langle G-m\rangle  = m^2 \langle G-m\rangle + m\langle G-m\rangle^2 -m\langle\underline{HG}\rangle.
\end{equation}
Here, the \emph{underline} term, defined as  $\underline{HG}= HG +\langle G\rangle G$ in the Wigner case,
captures the main fluctuation of $G-m$. In the Gaussian case, it satisfies $\Expv \underline{HG}=0$, making it a kind of renormalization of $HG$. 

We now informally outline three main difficulties that arise when analyzing~\eqref{eq:G-m}. While these issues can be relatively easily resolved in the Wigner case, they become significantly more serious for more general ensembles—especially for random band matrices.

\begin{itemize}
\item[(i)] {\bf (Truncation of the hierarchy)}
Since the variance of the  fluctuating error term $\langle\underline{HG}\rangle$
involves  squares of the resolvent, more precisely $GG^*$, understanding a single resolvent through~\eqref{eq:G-m} 
requires control over products of two resolvents. 
In turn, analyzing $GG^*$, via an equation similar to~\eqref{eq:G-m}, would lead to terms involves four resolvents, and so on. 
Potentially, this phenomenon leads to an infinite hierarchy of equations, reminiscent of the BBGKY hierarchy
in many-body theories. 
In the Wigner case, this hierarchy is automatically truncated at the very first step via the Ward identity, $GG^* = \im G/\eta$, 
 which is applicable here  thanks to the constant variance profile.  However,   once the variance profile $S_{ab}$ becomes not-constant,
the corresponding term takes the form $G\Sigma G^*$, where $\Sigma$ is a nontrivial matrix directly related to $S$. This prevents direct application of the Ward identity. While a norm bound on $\Sigma$
can restore the $GG^*$ structure, it often results in overestimating $\Sigma$. 
This issue is particularly pronounced for band matrices,    where the maximum norm $\max_{ab} S_{ab}\lesssim 1/W$ of $S$ is much larger than 
its $\ell^1$-norm $N^{-2}\sum_{ab} S_{ab} = 1/N$.

\item[(ii)] {\bf (Higher-order moment expansion)} 
To establish high-probability bounds on $\langle G-m\rangle$, one requires strong concentration estimates on the fluctuating term, which is usually achieved by computing high moments. Historically, this involved intricate expansion techniques, often requiring
very involved bookkeeping of Feynman diagrams.  Recently, a more efficient   cumulant expansion technique
substantially reduced the combinatorial complexity (see Footnote~\ref{min} for one aspect
of this simplification). 

\item[(iii)] {\bf (Stability)} The key property of~\eqref{eq:G-m} is that the target quantity $\langle G-m\rangle$
also appears linearly on the right-hand side, multiplied by a prefactor (in this case, $m^2$).  Solving the equation then involves
bringing this term to the left-hand side and divide by $1-m^2$. This is not a coincidence for the Wigner case,
 rather it reflects a general pattern that also holds for analogous \emph{self-consistent equations} in more complex settings (including longer resolvent chains).
In the bulk spectrum, $1-m^2$ is separated away from zero, so  the equation is stable, 
but near the edges of the semicircle, $z\approx \pm 2$, $1-m^2$ becomes small,
hence~\eqref{eq:G-m} becomes potentially unstable. 
For more general models, $1-m^2$ is 
replaced by a \emph{stability operator}, which may have nontrivial unstable directions that require separate analysis. 
This instability then needs to be balanced by improved bounds on the fluctuation term $\langle\underline{HG}\rangle$
and possible other terms on the right-hand side of~\eqref{eq:G-m}---which is often a daunting task.

\end{itemize}

Beyond a single resolvent $G$, one is naturally led to study multi-resolvent chains of the form
 $G_{[1,k]}:= G_1A_1G_2A_2\ldots A_{k-1}G_k$ with
arbitrary deterministic observables. There are several motivations for considering such chains: multi-resolvent local laws
with nontrivial observables arise naturally in thermalization problems~\cite{cipolloni2021thermalisation}
and in functional central limit theorems~\cite{cipolloni2023functional}. Some version 
of the special case when $A_j=I$ and $k\le 3$ have appeared much earlier in CLT for mesoscopic  linear eigenvalue
 statistics~\cite{He2017WignerCLT} and
 related results \cite{erdHos2018fluctuations, he2020mesoscopic}.  
 
 Another key observation is that the typical size of $G_{[1,k]}$ in the interesting small $\eta$ regime heavily
 depends on whether some of the observables are traceless. Very roughly speaking, we
have $G_{[1,k]}\sim \eta^{-k+1}$ for general $A_j$'s (in particular if all $A_j=I$) but each traceless
observable reduces its size by $\sqrt{\eta}$. This so-called \emph{$\sqrt{\eta}$-rule} was first observed
along the proof of the general QUE for Wigner matrices in~\cite{cipolloni2021eigenstate},
and was later formalized in~\cite{Cipolloni2022Optimal}.
 
\subsubsection{Hierarchy of static master inequalities}\label{sec:staticmaster}
The systematic study of multi-resolvent local laws was initiated in~\cite{Cipolloni2022Optimal},
where the Wigner case was settled in full generality. More importantly, that
work introduced the concept of \emph{hierarchy of master inequalities}, which also plays a crucial role in the current paper.
We now explain it in a bit more detail in the context of Wigner matrices.
 
The first step  in analyzing the multi-resolvent chains $G_{[1,k]}$   is to obtain formulas and optimal bounds on their deterministic approximations $M_{[1,k]}$. 
This was already achieved  in~\cite{cipolloni2021thermalisation}, using explicit combinatorial identities based on non-crossing partitions, which are specific to the Wigner case. 
The point is that the natural recursion for $M_{[1,k]}$ (see~\eqref{eq:otherM}--\eqref{eq:M_recursion} later) is nonlinear, quite unstable, and carries many inherent cancellations. 
For example, consider
\begin{equation}\label{M12}
	M_{[1,2]} =  m_1m_2 \trless{A}_1 + \frac{m_1m_2}{1-m_1m_2} \langle A_1 \rangle, \qquad m_j = m(z_j),
\end{equation}
where 
$\trless{A}:= A-\langle A\rangle$ denotes the traceless part of $A$. 
If $z_1\approx \overline{z}_2$ then the denominator $|1-m_1m_2|\sim\eta$ becomes very small. 
However, this instability affects only the tracial part of $A_1$.
For longer chains, the recursion involves multiple denominators of the form $1-m_im_j$, and their stability depends on whether $z_i$ is close to $\overline z_j$. 
A naive estimation of the recursion leads to bounds on $M_{[1,k]}$ that are much worse than its true size; only by applying nontrivial algebraic identities one can  rewrite the formula for $M_{[1,k]}$ in a form which shows its correct size.
These identities also take the number of traceless observables into account.

The second step, once the correct bounds on $M_{[1,k]}$ are established, is to control  
the fluctuations \mbox{$G_{[1,k]}-M_{[1,k]}$} when all observables are traceless.
Indeed, it suffices to consider only such chains. Any observable can be decomposed as $A=\langle A\rangle I + \trless{A}$, which implies 
\begin{equation} \label{Adecomp}
	\ldots GAG\ldots = \big(\ldots GG \ldots\big)\langle A\rangle+\ldots G\trless{A}G\ldots.
\end{equation}
The first term on the right involves a resolvent square $G^2$, which can be reduced to a single $G$ using contour integral representations.
This results in a shorter chain, suitable for induction. 
Thus the main task reduces to understanding chains with traceless observables.
To this end, \cite{Cipolloni2022Optimal} introduces two sequences of random variables $\Psi^{\mathrm{iso}}_k$ and $\Psi^{\mathrm{av}}_k$, to control isotropic and averaged fluctuations, respectively. These quantities are renormalized by the expected size of the fluctuations, roughly speaking\footnote{
	The precise definition involved appropriate maxima over other parameters.
	}
$$
   \Psi^{\mathrm{av}}_k = (N\eta) \eta^{k/2-1}\langle (G_{[1,k]}-M_{[1,k]})A_k\rangle, \quad 
   \Psi^{\mathrm{iso}}_k = \sqrt{N\eta}\, \eta^{k/2} \big( G_{[1,k+1]}-M_{[1,k+1]} \big)_{\bm u \bm v}.
$$
Here, the factors $\eta^{k/2-1}, \eta^{k/2}$ indicate the typical (optimal) size of the deterministic $M$-term, and the (inverse of the) factors $N\eta$ and $\sqrt{N\eta}$ reflect the relative size of the fluctuation.
To estimate $G_{[1,k]}-M_{[1,k]}$, a self-consistent equation analogous to~\eqref{eq:G-m} have been derived, with the right hand side containing a fully underline term, together with error terms involving shorter chains. 
Assuming an initial very-high-probability bound  $\Psi^{\mathrm{iso}}_k\lesssim \psi^{\mathrm{iso}}_k$,
$\Psi^{\mathrm{av}}_k\lesssim \psi^{\mathrm{av}}_k$ for some deterministic control sequences $\psi^{\mathrm{iso}}_k$, $\psi^{\mathrm{av}}_k$,
one obtains a schematic estimate of the form:
\begin{equation}\label{eq:GkMk}
     \langle (G_{[1,k]}-M_{[1,k]})A_k\rangle \approx  -m\langle\underline{H G_{[1,k]}A_k} \rangle +\mathcal{O}\Bigg( 
     \frac{1}{N\eta^{k/2}} \Bigg\{ 1+ \sum_{j=1}^{k-1} \psi^{\mathrm{av}}_j\Big(1+ \frac{  \psi^{\mathrm{av}}_{k-j}}{N\eta}\Big)\Bigg\}     \Bigg).
\end{equation}
The underline term can be estimated in high moment sense using a cumulant expansion\footnote{ \label{min}
	Specifically, a \emph{minimalistic cumulant expansion} is employed: instead of computing $\Expv |\langle\underline{HG}\rangle|^{2p}$ directly, one writes
	$$
		\Expv | \langle G-M\rangle|^{2p} =   \Expv \langle \underline{HG}\rangle
		\langle G-M\rangle^{p-1} \overline{\langle G-M\rangle}^{p}, \qquad G = G_{[1,k]}A_k, \quad M = M_{[1,k]}A_k,
	$$
	and expands the underline. This technique preserves many instances of $\langle G-M\rangle$, reducing the complexity of the expansion.
	}. 
Higher-order cumulants (of order at least three) naturally bring in isotropic quantities, even in the estimates for the average chains.
Therefore, the average hierarchy cannot be closed independently of the isotropic one---except in the special Gaussian case, where all higher cumulants vanish.
 
The result is an inequality of the form
\begin{equation}\label{eq:Psik}
   \Psi_k^{\mathrm{av}}    \lesssim 1+ \frac{\sqrt{\psi_{2k}^{\mathrm{av}}}}{\sqrt{N\eta}} +  \sum_{j=1}^{k-1} \psi^{\mathrm{av}}_j\Big(1+ \frac{  \psi^{\mathrm{av}}_{k-j}}{N\eta}\Big) +  \frac{\psi_k^{\mathrm{iso}}}{\sqrt{N\eta}}+ \ldots,
\end{equation}
where the remaining terms (indicated by $\ldots$) are not relevant for this discussion. 
A similar inequality is derived for the isotropic quantity $\Psi_k^{\mathrm{iso}}$.
The goal is to exploit the \emph{self-improving} character of~\eqref{eq:Psik}: given $\Psi^{\mathrm{iso}}_k\lesssim \psi^{\mathrm{iso}}_k$, $\Psi^{\mathrm{av}}_k\lesssim \psi^{\mathrm{av}}_k$, the master inequalities yield a \emph{smaller} upper bounds on $\Psi^{\mathrm{iso/av}}_k$.
This improved bound can then be fed back into the same inequality, producing an iterative scheme.  
Starting from a rough a priori bound (with $\psi_k$'s), repeated application of the master inequalities produces the optimal estimate $\Psi_k^{\mathrm{av/iso}}\lesssim 1$.
Roughly speaking, the self-improving character of~\eqref{eq:Psik} means that $\psi_k^{\mathrm{av/iso}}\sim1$ is the globally stable fixed point of the iteration.

This self-improving mechanism relies crucially on two structural properties of the hierarchy:
\begin{itemize}
	\item[(i)] The terms on the right-hand side of~\eqref{eq:Psik} (and its isotropic version) contain only $\psi$'s with smaller or equal $k$,
	\item[(ii)] The terms containing $\psi_k$ themselves are suppressed by a small factor, such as $1/\sqrt{N\eta}\ll 1$.
\end{itemize}
The presence of $\sqrt{\psi_{2k}}$ apparently breaks this pattern. 
This term originates from the variance of the underline term, which behaves like a $2k$-chain (even after applying two Ward identities).
In fact, this is the manifestation of the truncation problem, described in Section~\ref{sec:dev}.  
To restore self-improvement,   \eqref{eq:Psik} must be complemented by a reduction inequality that estimates a $2k$-chain in terms of a product of two $k$-chains. 
Ideally, one would like to have $\Psi_{2k}^{\mathrm{av}}\le (\psi_k^{\mathrm{av}})^2$, so that the problematic term  $\sqrt{\psi_{2k}^{\mathrm{av}}}$ could be replaced with $\psi_k^{\mathrm{av}}$ thereby closing the hierarchy. 
In practice, however, such reduction estimates always incur a loss. 
For instance, in~\cite{Cipolloni2022Optimal} established that $\Psi_{2k}^{\mathrm{av}}\le (N\eta)^2 + (\psi_k^{\mathrm{av}})^2$ (at least for even $k$), which necessitated a two-step procedure: 
first, a suboptimal bound $\psi_k^{\mathrm{av/iso}}\lesssim \sqrt{N\eta}$ was proven for all $k$; 
 this bound was then used to control $\psi_{2k}$ and reinserted into~\eqref{eq:Psik}.  

In the present paper, we develop a more refined mechanism using a system of \emph{loss exponents}, which enables us to distribute the reduction loss at level $2k$ across intermediate chain lengths between $k$ and $2k$.
Most importantly, we need to control longer resolvent chains, 
even when the ultimate goal is to prove local laws for chains of small length (e.g., $k=2$).
This seems to be a common feature of BBGKY-type arguments: the cost of truncation or reduction may be affordable only if it happens at a high order of the hierarchy.

A fundamental guiding principle behind the master inequality strategy is that it is ultimately governed by the optimal bounds on the deterministic $M$-terms.  
This key point is somewhat implicit in the definition of the normalized quantities $\Psi_k$, which are normalized for optimal estimates on $G-M$, rather than $M$.
Nevertheless, one can heuristically observe that the hierarchy remains internally consistent even with the looser bounds $\psi_k^{\mathrm{av}} = N\eta$, $\psi_k^{\mathrm{iso}} = \sqrt{N\eta}$ if the $1$'s in \eqref{eq:Psik} were replaced by $N\eta$ 
 (corresponding to the typical size of $M_{[1,k]}$). 
This observation reflects the fact that the fluctuations remain controlled as long as $G_{[1,k]}$ does not exceed the size of $M_{[1,k]}$.
Upgrading such statement to genuine local laws requires the $1$'s in \eqref{eq:Psik} that actually come from the optimal bounds on the leading $M$-terms  coming from the underline.   
This principle makes the master inequalities particularly powerful: 
they basically reduce the problem of finding the optimal fluctuations of $G_{[1,k]}-M_{[1,k]}$ to  
 estimating longer $M$-terms.   
Roughly speaking, the $M$’s drive the dynamics of the entire hierarchy, and the $G-M$ fluctuations follow them.

\subsubsection{Dynamical master inequalities and zigzag}

Considering the three main difficulties listed in Section~\ref{sec:dev}, 
the master inequalities presented in Section~\ref{sec:staticmaster} provide an effective solution to the truncation problem, via reduction estimates that are pushed to sufficiently long chains.
The second difficulty—moment expansion—is still addressed using cumulant methods. As for the third issue, stability, it appears to vanish entirely when only traceless observables are considered. 
For this last point, the Wigner model in the bulk features two major simplifications that makes the critical stability issue look deceivingly harmless\footnote{
	This problem already resurfaced at the spectral edge for Wigner matrices and it required a new approach, see~\cite{cipolloni2023eigenstate}.
	}.
The first one is the spectral gap, common to all mean-field models, and the second one is a special structural property of the Wigner matrix.
We now explain both, as overcoming them is central to extending the theory to RBM.

Just as~\eqref{eq:G-m} expresses a self-consistent equation for the single resolvent, one can write similar equations for chains of arbitrary length  $k\ge 2$.
The associated stability operator is called the \emph{two-body stability operator}\footnote{
	The central role of the two-body stability operator for chains of any length can also be seen in the recursion formulas for $M$, 	see~\eqref{eq:otherM}--\eqref{eq:M_recursion}. 
	When $z_1=z_2$, it is called the \emph{one-body stability operator}; 
	this appears in the single-resolvent equation~\eqref{eq:G-m}.
	}.
For Wigner matrices, this linear operator takes the form
$$
   B_{12}: = I  - m_1m_2\langle\cdot \rangle,
$$
acting on matrices $\mathbb{C}^{N\times N}$. 
For example, the potentially singular term in~\eqref{M12} arises from inverting $B_{12}$.
This operator has a very convenient spectrum: 
it has one potentially small (if $z_1\approx\bar z_2$) eigenvalue $1-m_1m_2$ with eigenvector being the identity matrix. 
All other eigenvalues of $B_{12}$ are harmless $1$'s. In particular, $B_{12}$ acts as the identity on traceless matrices and is thus trivially invertible on that subspace.
In other words, $B_{12}$ has a large (order one) spectral gap above its single smallest eigenvalue, enabling one to treat the tracial and traceless components separately.
This gap structure is typical in mean-field models, including Wigner-type matrices with variance profiles satisfying $S_{ab}\sim 1/N$. 
We stress that random band matrices with $W\ll N$ do not have this large gap, which is the principal reason why band matrices are significantly harder to analyze than mean-field models.

However, there is another structural simplification that holds exclusively for Wigner matrices.
To explain this, we recall that in more general models, the two-body stability operator takes the form
$$
  B_{12} = I- M_1 \mathcal{S}[\,\cdot\,] M_2,
$$
where $\mathcal{S}$ is the self-energy operator, defined by $\mathcal{S}[R]= \Expv [H R H]$, and $M_j$ is the solution to the matrix Dyson equation
$$
   -M_j^{-1}= z_j- D+ \mathcal{S}[M_j], 
$$
subject to the positivity condition $(\im M_j)(\im z_j)>0$. 
In the Wigner case, $\mathcal{S}[\cdot ]=\langle \cdot \rangle$ and $M_j = m_j\cdot I$.  
For random band matrices with the condition $\sum_a S_{ab}=1$, one still has $M_j = m_j\cdot I$, but now $\mathcal{S}$ and hence $B_{12}=(1-m_1m_2 \mathcal{S})$ are nontrivial, although their bad eigendirection is still the identity.
For example, in case of RBM, the two-resolvent analogue of~\eqref{eq:G-m} for $G_{[1,2]}=G_1AG_2$ is\footnote{
	The precise definition of the renormalization $\underline{G_{1}A'HG_{2}}$ is given in~\eqref{eq:underline_def}.
	}
\begin{equation} \label{eq:2G}
	(G-M)_{[1,2]} = -\underline{G_{1}A'HG_{2}} 
	+ (G_1-m_1)A' 
	+G_{1} \mathcal{S}\bigl[(G_1-m_1)A'\bigr]G_{2} 
	+ G_{1}A'\mathcal{S}\bigl[G_{2}-m_{2}\bigr]G_{2}.
\end{equation}
where $A'  := m_2 B_{12}^{-1}[A]$.
In the Wigner case, where $S_{ab}=1/N$, for $A=I$, this simplifies to 
\begin{equation} \label{eq:2GW}
							(1-m_1m_2) (G-M)_{[1,2]} = -m_2\underline{G_{1}HG_{2}} 
				+ m_2 (G_1-m_1)
				+m_2G_{1} G_2 \big[ \langle G_1-m_1 \rangle  +  \langle G_{2}-m_{2}\rangle\big].
		\end{equation}
If $A$ is traceless, $\langle A \rangle=0$, the equation becomes
\begin{equation} \label{eq:2GWnotr}
					 (G-M)_{[1,2]} = -m_2\underline{G_{1}AHG_{2}} 
				+ m_2 (G_1-m_1)A
				+m_2G_{1} G_2 \langle(G_1-m_1)A\rangle 
				+ m_2G_{1}A G_2 \langle G_{2}-m_{2}\rangle.
		\end{equation}
These last two identities highlight a key dichotomy specific to the Wigner setting:
either the potentially bad stability factor $(1-m_1m_2)$ is absent, as in~\eqref{eq:2GWnotr}, or the largest term has a pure product of $G_1G_2$,
as in~\eqref{eq:2GW}, which is particularly amenable to algebraic reduction, either via contour integral representations or (in some cases) a Ward identity.
This algebraically resolves a major cancellation on the right hand side of~\eqref{eq:2GW} when $z_1\approx\bar z_2$.
Indeed, after taking the averaged trace, the left hand side is only $(1-|m|^2)\times 1/N\eta^2 \approx 1/N\eta$, while the last term is of order $1/N\eta^2$, and it cancels the underline term to leading order.

 However, in essentially any generalization of the Wigner matrix, the simplifications appearing in~\eqref{eq:2GW}–\eqref{eq:2GWnotr} no longer hold. Instead of a clean product $G_1 G_2$, one typically encounters  
$G_1\Sigma G_2$ 
 with a non-trivial observable $\Sigma \ne I$.  
This is already evident in models such as deformed Wigner matrices of the form $H+D$~ \cite{cipolloni2023gaussian}, and Wigner-type matrices~\cite{R2023bulk, erdHos2024eigenstate}.
For random band matrices, even for $A=I$ in~\eqref{eq:2G}, we have
\begin{equation}\label{def:Sa}
  G_{1} \mathcal{S}\bigl[(G_1-m_1)A'\bigr]G_{2} = m_2\sum_a (G_1-m_1)_{aa} G_1S^aG_2, \qquad (S^a)_{ij} := \delta_{ij}S_{ai},
\end{equation}
so this term generates a two-resolvent chain $G_1S^aG_2$ with a non-trivial diagonal observable $S^a$ in the middle.
Therefore, we need to handle chains with quite general observables directly and face the consequences of the bad eigendirection.

This instability in the self-consistent equations for resolvent chains represents a major obstacle: it conceals a crucial cancellation that is often not apparent at the algebraic level.
It is typically connected with the instability of the two-body operator $B_{12}$ when $z_1\approx\overline{z}_2$, but the same problem is already present in the single-resolvent equation~\eqref{eq:G-m} at the spectral edge where $|1-m^2|\ll 1$.

The {\it characteristic flow method} resolves this instability problem by studying the evolution of $G_t(z_t)=(H_t-z_t)^{-1}$, where the spectral parameter follows the characteristic equation~\eqref{eq:char_flow1}, and at the same time the matrix $H_t$ evolves along the Ornstein-Uhlenbeck process~\eqref{OU}.
A simple calculation using It\^o-calculus for the averaged trace of the resolvent $G:=G_t(z_t)=(H_t-z_t)^{-1}$ yields the stochastic differential equation
\begin{equation}\label{gmeq}
   \mathrm{d}\langle G-m\rangle = \Big[ \frac{1}{2} \langle G-m\rangle  + \Big( \langle G\rangle   + \frac{1}{2}z_t + 
   \frac{\mathrm{d}z_t}{\mathrm{d}t}\Big) \langle G^2\rangle\Big] \mathrm{d}t  + \mathrm{d} \mathcal{M}_t 
   = \Big[ \frac{1}{2}+   \langle G^2\rangle\Big]  \langle G-m\rangle \mathrm{d}t  + \mathrm{d} \mathcal{M}_t,
\end{equation}
where $\mathcal{M}_t$ is a martingale term. 
The key cancellation occurs in the inner brackets $(\dots)$: thanks to the carefully chosen trajectory of the spectral parameter $z_t$, governed by ~\eqref{eq:char_flow1}, reduces $ \langle G\rangle\sim m\sim1$ to $\langle G-m\rangle\lesssim 1/N\eta$ as the coefficient of the biggest term $ \langle G^2\rangle$.
Note that the equation~\eqref{gmeq} is a linear stochastic ODE with a first \emph{linear term} plus an inhomogeneous martingale term.
 
Similar
equations hold for longer chains, for example for the critical $G_{[1,2]}$ chain from~\eqref{eq:2GW} we have
\begin{equation}\label{2Gzig}
	\mathrm{d}\langle (G-M)_{[1,2]}\rangle = 2 \Big[ \frac{1}{2}+
	\langle M_{[1,2]}\rangle\Big] \langle (G-M)_{[1,2]} \rangle \mathrm{d}t
	+ \mathcal{F}  \mathrm{d}t + \mathrm{d} \mathcal{M}_t,
\end{equation}
where $\mathcal{F}$ is a  \emph{forcing term} with a typical size $|\mathcal{F}_t|\lesssim 1/N\eta_t^3$, consistent with the time derivative\footnote{
	Note from~\eqref{eq:char_flow1} that $\eta_t$ is roughly linear in $t$, so $\mathrm{d}\eta_t/\mathrm{d} t \sim -1$ and thus 	$|\mathrm{d}\eta_t^{-k}/\mathrm{d} t |\sim  \eta_t^{-k-1}$.
	}
of the target bound $\langle (G-M)_{[1,2]}\rangle_t\le 1/N\eta_t^2$. 
A key point is that this estimate on $\mathcal{F}$ can be obtained without knowing the sharp  target bound on $(G-M)_{[1,2]}$, only weaker bounds on $G_{[1,2]}$ and   sharp bounds on the shorter chain $G-m$ are used.
This structure makes the equation not only \emph{self-consistent}, but also \emph{self-improving}, eventually leading to the proof of the target bound after using an inductive iterative scheme. 
  
To analyze~\eqref{2Gzig}, we introduce the \emph{propagator}\footnote{
	Here the propagator is just a scalar factor, but we prepare the terminology for the more complicated RBM case.
	}
\begin{equation}\label{prop}
   \mathcal{P}_{s,t}  : =  \exp\Big(  \int_s^t  \Big[ \frac{1}{2}+ \langle  M_{[1,2]}\rangle_r \Big] \mathrm{d}r \Big)
   \lesssim  \frac{\eta_s}{\eta_t},
\end{equation}
 where in the last step we used that $\langle M_{[1,2]}\rangle \approx \rho/\eta$ with $\rho:=\im m$, and 
\begin{equation}\label{int}
	\int_s^t \frac{\rho_r\mathrm{d}r}{\eta_r}\lesssim \log \frac{\eta_s}{\eta_t},
\end{equation}
which follows from~\eqref{eq:char_flow1}.  
Using Duhamel's principle for~\eqref{2Gzig} and assuming the initial bound 
\begin{equation}\label{initbound}
 \langle (G-M)_{[1,2]} \rangle_s\lesssim \frac{1}{N\eta_s^2},
 \end{equation}
 at time $s$,
   we obtain the final bound
  \begin{equation}\label{GGzig}
    \langle (G-M)_{[1,2]}\rangle_t \approx \mathcal{P}_{s,t}^2 \langle (G-M)_{[1,2]}\rangle_s 
    + \int_s^t \mathcal{P}_{u,t}^2  \Big[ \mathcal{F}_u \mathrm{d} u +\mathrm{d} \mathcal{M}_u\Big]
   \lesssim \frac{1}{N\eta_t^2},
 \end{equation}
 at a later time $t$.
In the last step, for the forcing term we used
$$
   \Big| \int_s^t \mathcal{P}_{u,t}^2  \mathcal{F}_u \mathrm{d} u \Big|
    \lesssim  \int_s^t \Big(\frac{\eta_u}{\eta_t}\Big)^2 \frac{1}{N\eta_u^3}
  \mathrm{d} u \lesssim\frac{1}{N\eta_t^2},
$$
up to irrelevant logarithmic factors. A similar estimate holds for the martingale term.

Since $\eta_t$ decreases along the flow~\eqref{eq:char_flow1} (assuming $z$ is in the upper half plane), the equation~\eqref{2Gzig} propagates the initial bound~\eqref{initbound} for large $\eta_s$ to all smaller $\eta_t$ in a consistent way and without any intricate cancellation. 
As an additional bonus, no high moment cumulant expansion is necessary, unlike when estimating the fluctuating underline term in the static method in~\eqref{eq:GkMk}. 
In the dynamical approach, the fluctuation is encoded in the martingale term, where a simple Burkholder-Davis-Gundy (BDG) inequality can turn bounds on the second moment (quadratic variation) immediately into high probability bounds.

The only drawback of the characteristic flow method is that the  distribution of the
 original matrix changes along the flow~\eqref{OU}, even if its expectation and variance are preserved; compared to $H_s$,  $H_t$ pics up a Gaussian component of order $\sqrt{t-s}$:
$$
     H_t \stackrel{\rm d}=e^{(s-t)/2} H_s + \sqrt{1- e^{s-t}} \mathcal{H}_{\mathrm{G}},
$$
where $\mathcal{H}_{\mathrm{G}}$ is a standard GUE or GOE matrix, independent of $H_s$. 
This Gaussian component can, however, be removed using {\it Green function comparison (GFT)} arguments, at least provided the time increment $t-s$ is sufficiently small, and if isotropic estimates on resolvent chains are available.
Since only a small Gaussian component can be removed with each GFT step, we must iterate the procedure to reach the small $\eta$ regime. 
This gives rise   to the \emph{zigzag} strategy,  an alternating sequence of   
running the characteristic flow to reduce $\eta$ a bit, and using GFT to remove the resulting Gaussian component,
see the schematic Figure~\ref{fig:zigzag_fig} later.  
 
In this paper, we apply this general idea to RBM---a significantly more challenging model. 
The main novelty is that, unlike in mean-field models, all relevant quantities in RBM exhibit a non-trivial spatial dependence, stemming from the fundamental off-diagonal decay of $|G_{ab}|$ on an intermediate length scale $\ell=\ell(\eta)$. 
As a result,  the simple scalar ODE~\eqref{2Gzig} in the Wigner case is replaced by a matrix-valued equation,
 resembling a space-discretized  $1+1$ dimensional PDE. 
The evolution is governed by a  propagator $\mathcal{P}_{s,t}$, which now acts as a matrix operator tracking how the spatial profile of quantities such as $(GS^aG)_{bb} = \sum_c S_{ac} |G_{cb}|^2$ changes along the flow.
Heuristically, the correct leading-order behavior of these spatial profiles was already understood in~\cite{erdHos2013delocalization};
$$
    (GS^aG)_{bb} \approx \Theta_{ab} = \Big(\frac{|m|^2 S}{1-|m|^2S} \Big)_{ab}.
$$
However, due to the instability inherent in the RBM analogue of~\eqref{eq:2GW}, this behavior could not previously be rigorously proven across the optimal parameter range $W\gg \sqrt{N}$, $\eta\gg 1/N$.

 \subsection{Main ideas for band matrices}\label{sec:ideasRBM}

The core part of the analysis is to prove averaged and entry-wise local laws for the complex Hermitian symmetry class
in the special case when all observables in the resolvent chains
are of the form $A_i =S^{x_i}$ for some $x_i\in \indset{N}$. Such \emph{special} 
observables naturally arise 
along the evolution, see~\eqref{def:Sa}, and the hierarchy remains closed for such chains.
Extensions to general observables,  isotropic chains,  real symmetric matrices 
and improvements for traceless observables will be explained in Section~\ref{sec:extension}. 
 
 \subsubsection{Failure of the  mean-field approach}
 
To illustrate the main idea, we consider the evolution of the  averaged $k=2$ chain 
$$
  \Tr \big[ G_{[1,2]}S^{x_2}  \big]=  \Tr \big[ G_1 S^{x_1}G_2S^{x_2} \big], \qquad G_j=  G_t(z_{j, t}), \quad j=1, 2, 
$$
along the zig flow with the most critical choice of
spectral parameters $z:=z_1=\bar z_2$.  
The spectral parameters evolve along the same characteristic equation~\eqref{eq:char_flow1}, while
the random matrix $H_t$ evolves along the OU flow~\eqref{OU}, driven by a matrix Brownian motion $\mathcal{B}$ with variance $\mathrm{d}\mathcal{B}_{ab} \mathrm{d}\mathcal{B}_{ba} = S_{ab}\mathrm{d}t$. This choice
preserves the first two moments of $H_t$ along the flow. 
Denote the fluctuation by
$$
  \mathcal{X}  = \mathcal{X}_t(x_1, x_2):  = \Tr \big[(G_{[1,2]}-M_{[1,2]})S^{x_2}\big],
$$
where $M_{[1,2]} =M_{[1,2]}(x_1) $ is the deterministic approximation of $G_{[1,2]} =G_1 S^{x_1}G_2$, which is a
diagonal matrix with  $(M_{[1,2]})_{aa} = \Theta_{a x_1}(z_1)$.
We view $\mathcal{X}$ as a time-dependent (random) function of two space variables $x_1, x_2$.
Although we often omit the time variables from the notation, all quantities are time-dependent.
Similarly to \eqref{2Gzig},   we find that 
\begin{equation}\label{eq:12}
\mathrm{d} \mathcal{X}  = \bigoplus_{j=1}^2 \Big( \frac{1}{2}I  +  \Theta \Big)^{(j)}  \bigl[ \mathcal{X}\bigr]\mathrm{d}t 
		+ \mathrm{d}\mathcal{M}  + \mathcal{F} \mathrm{d}t,
\end{equation}
with some forcing and martingale terms, 
where the upper index $(j)$ indicates that the operator acts on the $j$-th spatial variable. The propagator corresponding to the (time-dependent) generator $I+\Theta$\footnote{The discrepancy between 
$I+\Theta$ and $\frac{1}{2}I +\Theta$ is irrelevant and will be ignored in this discussion.} is given by
 \begin{equation} \label{eq:propa}
	\mathcal{P}_{s,t} := \exp\biggl\{ \int_s^t (I + \Theta_r)\mathrm{d}r \biggr\} = \frac{|m_t|^2}{|m_s|^2}\frac{1-|m_s|^2 S}{1-|m_t|^2 S}, \quad s \le t, 
\end{equation}
where we used the relation $\mathrm{d}\Theta/\mathrm{d}t = (1+\Theta)\Theta$ along the 
characteristic flow~\eqref{eq:char_flow1}, which can be verified by an explicit calculation.
The propagator is a matrix with positive entries, and one can easily check that
\begin{equation}\label{eq:propapprox}
   \mathcal{P}_{s,t} \approx I + (t-s)\Theta_t \approx  I + \eta_s \Theta_t, \quad s \le t,
\end{equation}
up to irrelavant constants, 
and that it satisfies the same norm bound as its scalar version~\eqref{prop}:
\begin{equation}\label{Pnorm}
   \|   \mathcal{P}_{s,t} \| \lesssim \frac{\eta_s}{\eta_t}, \quad s \le t.
\end{equation}
We have the Duhamel formula, analogous to~\eqref{GGzig},
\begin{equation}\label{GSGzig}
    \mathcal{X}_t = \bigotimes_{j=1}^2 \mathcal{P}_{s,t}^{(j)} \bigl[\mathcal{X}_s \bigr]
    + \int_s^t \bigotimes_{j=1}^2 \mathcal{P}_{u,t}^{(j)} \Big[ \mathrm{d} \mathcal{M}_u +
    \mathcal{F}_u \mathrm{d} u \Big].
 \end{equation}
Using the norm bound~\eqref{Pnorm}, one could consistently propagate an estimate of the form $|\mathcal{X}_t|\lesssim 1/(N\eta_t)^2$ 
by~\eqref{GSGzig}, as was previously done in the mean-field case. However, the target estimate for RBM is different: in the regime $\eta\ge (W/N)^2$, it has a non-trivial spatial dependence
\begin{equation}\label{X12}
   | \mathcal{X}_t (x_1,x_2)| \lesssim \frac{1}{\ell_t\eta_t} (\Upsilon_t)_{x_1x_2},\qquad \mbox{with}\quad
    \Upsilon_{x_1x_2}  \lesssim \frac{1}{\ell\eta}.
\end{equation}
In particular, in this  regime  the target  $\ell^\infty$-bound on $\mathcal{X}$, 
\begin{equation}\label{target}
      \| \mathcal{X}_t\|_\infty \le \frac{1}{(\ell_t\eta_t)^2}, 
      \end{equation}
cannot be propagated using the norm bound~\eqref{Pnorm} in the Duhamel formula:
\begin{equation}\label{badest}
   \Big\| \bigotimes_{j=1}^2 \mathcal{P}_{s,t}^{(j)}  \mathcal{X}_s\Big\|_\infty \le \Big( \frac{\eta_s}{\eta_t}\Big)^2 
   \frac{1}{(\ell_s\eta_s)^2} 
      = \Big( \frac{\ell_t}{\ell_s}\Big)^2 \frac{1}{(\ell_t\eta_t)^2}  
\end{equation}
since the excess factor $(\ell_t/\ell_s)^2$ is typically very large. 

In our actual proof, we avoid using crude $\ell^\infty$-norm as long as $\eta \ge (W/N)^2$. Instead, we keep track of
the spatial dependence of every chain using the $\Upsilon$ factors, as introduced in~\eqref{X12}. In this more refined way,
the analogue of~\eqref{badest} becomes (see~\eqref{eq:bad_av_prop})
\begin{equation}\label{badest1}
  | \mathcal{X}_s (x_1,x_2)| \lesssim \frac{1}{\ell_s\eta_s} (\Upsilon_s)_{x_1x_2} \quad \Longrightarrow \quad
   \Big( \bigotimes_{j=1}^2 \mathcal{P}_{s,t}^{(j)}  \mathcal{X}_s\Big)_{x_1x_2} 
   \lesssim \frac{\ell_t}{\ell_s} 
   \frac{1}{\ell_t\eta_t} (\Upsilon_t)_{x_1x_2},
\end{equation}
so the excess factor is only $\ell_t/\ell_s$, but still too large.
Similar problem occurs with the inhomogeneous forcing and martingale terms in~\eqref{GSGzig}. In short, there are
too many propagators in~\eqref{GSGzig} compared with the $(\ell\eta)$-singularity of the target bound
$|\mathcal{X}_t|\lesssim (\ell_t\eta_t)^{-1} \Upsilon_t$.
We will refer to this new key obstacle as the \emph{excess propagator problem}. 

We point out that this problem arises only for the case when $z_1\approx\bar z_2$, i.e. when the
generator in~\eqref{eq:12} consists of two $\Theta$’s.  Whenever at least one $\Theta$ is replaced by
$\Xi$, the corresponding propagator is bounded, $\|  \mathcal{P}_{s,t} \|\lesssim 1$, so
the Duhamel formula at least yields a self-consistent bound. 
Moreover, this problem is present only in averaged chains, the equation corresponding isotropic chain $(GS^xG)_{ab}$ 
has only one propagator instead of two.
Therefore, in our actual proof,  we
need to handle the general averaged $k$-chains with alternating $\mbox{sgn} (\im z_j)$ separately. We call them \emph{saturated chains}. All other, \emph{non-saturated}, chains are  easier, so we focus on the saturated averaged chains now.

We stress that the mismatch in the bound~\eqref{badest} is not just a lousy estimate: 
the estimate~\eqref{Pnorm} saturates on the constant mode ${\bf 1}$. 
In the mean field models (e.g. Wigner type matrices, $W\sim N$), a standard way to resolve a single
problematic mode would be to separate it from the rest and use that $ \mathcal{P}_{s,t} \lesssim 1$
on all other modes, basically due to the spectral gap of $1-|m|^2S$. In the RBM problem, $\mathcal{P}_{s,t}$
is large on many other low momenta (large distance) modes. Thus,
 the Duhamel formula carries a delicate balance between the low-lying spectrum of $1-|m|^2S$
 and the scale (in terms of $x_1, x_2$) on which $\mathcal{X}_s(x_1, x_2)$ lives. In fact, this balance is
 the essence of the quantum diffusion phenomenon which heuristically governs the dynamics of RBM.
 Taking care of the scale $\mathcal{X}_s(x_1, x_2)$
 in terms of a consistent upper bound, e.g., in terms of $\Upsilon$, helps, but it is still not sufficient to restore consistency 
 as the excess factor in~\eqref{badest1} shows. This is because $\mathcal{X}_s(x_1, x_2)$ itself
 is an oscillatory quantity, and there is an additional cancellation when the propagator acts on it.
 The oscillatory character, in principle, may be detected by fully decomposing  $\mathcal{X}$
 into the low-lying eigenmodes of $1-|m|^2S$, which are essentially plane waves, but this
 seems a practically impossible task. Therefore, we rely on another approach,
 the idea of \emph{observable regularization}, which we explain the next Section~\ref{sec:rego}.
 
 We remark that the key problem of excess propagators in~\eqref{GSGzig} has also been identified in~\cite{yauyin},
 where a concept of \emph{sum zero operator} was introduced to tackle it. The basic mechanism of these two approaches
 is similar, as they both capitalize on the different length scale of the propagators and the objects they act upon in 
 the Duhamel formula~\eqref{GSGzig}, but their implementation is quite different. Most importantly, our regularization
 is a fully local  operation on observables   within the chain; 
 it can be applied independently at several sites, with the effects from multiple regularizations combining multiplicatively.   In particular, we heavily rely on \emph{two} simultaneous regularizations 
 in the most critical parts of our proof: namely, to handle the longest average chain 
 (see Section~\ref{sec:rego}), and for providing a particularly short proof for the optimal $M$-bounds in Section~\ref{sec:Mb}.

\subsubsection{Observable regularization}\label{sec:rego}

	Regularization, in general, is an operation on observables that aims to reduce their size under the action of the propagator $\mathcal{P}_{s,t}$. Historically, 
	the guiding idea behind regularization stems from mean-field models, such as Wigner matrices, where subtracting the tracial part of the observable already produced the desired effect. In more general mean-field models, regularization took more complicated form (see \cite{Cipolloni2022overlap, cipolloni2023gaussian} where the idea was first introduced for deformed Wigner matrices), but always amounted to projecting the observable onto a co-dimension one subspace of matrices orthogonal to the lowest mode of the corresponding stability operator. 
	In particular, this  \emph{mean-field regularization} was effective because of the spectral gap.
		
	As we alluded to previously, the stability operator for random band matrices is essentially gapless, and simply subtracting any finite number of modes does not suffice. Instead, we introduce a new form of observable regularization, called \emph{spatial regularization}. It
	  takes advantage of the spatial structure inherent to RBM, which manifests in the regularity of the propagator entries $(\mathcal{P}_{s,t})_{xy}$ in the space-variables $x,y$. Modulo the irrelevant identity term from \eqref{eq:propapprox}, the profile of $\mathcal{P}_{s,t}$ matches that of $\eta_s \Theta_t$, and is essentially given by a smooth\footnote{
		For translation-invariant random band matrices, this can be obtained directly using Fourier analysis.
		} 
	function on the scale $\ell_t=\ell(\eta_t)$. In particular, it acts as a kind of low-pass filter, mollifying all functions that live on shorter length scales (higher frequencies). Therefore, we construct the spatial
	regularization of an observable to be a short-range object that averages out to zero.
	 We emphasize that the mean-field regularization (called simply regularization in earlier papers on mean-field models)
	and the new spatial regularizations are very different operations, even though they eventually serve similar purposes.
	For RBMs, both types of regularization appear in our proofs, but spatial regularization plays a much more central and ubiquitous role. Therefore, for brevity, we refer to it as simply as \emph{regularization} throughout this paper.
	The concept of mean-field regularization will only appear in Section~\ref{sec:traceless}, which focuses on traceless observables.
	We now explain the mechanism and implementation of this new form of regularization.

Consider a long resolvent chain $G_1S^{x_1}G_2 S^{x_2}G_3 S^{x_3} \ldots G_k$. 
To regularize the observable $S^{x_j}$, depending only on $x_j$, we use the preceding index $x_{j-1}$ as a reference,
and define     
$$
	\reg{S}^{x_j}: = S^{x_j} - \delta_{x_j, x_{j-1}}I = \sum_b (\delta_{x_j b} - \delta_{x_j x_{j-1}})S^b.
$$
The additional index $x_{j-1}$ is used as an anchor, which, together with 
the decay properties of the resolvent $G$, localizes $S^{x_{j-1}}G_{j}\reg{S}^{x_j}$ in space.
Note that observable regularization is a local procedure, depending only on the preceding observable\footnote{
	In mean-field models, the (mean-field) 
	regularization was a completely local, single-site operation on each observable. For more general ensembles, the only additional dependence came from the spectral parameters of adjacent resolvents, but the key point remained: \emph{each} observable in the chain could be regularized independently.
	}. 
In particular,  it lends itself to \emph{tensorization}, allowing   \emph{every second} observable to be regularized independently. 

To deal with the excess linear term (which eventually produces an excess propagator via Duhamel),   we can split the chain with $S^{x_j}$ into a sum of a chain with a regularized $ \reg{S}^{x_j}$
and another term with the identity matrix.  This second term is subject to 
the usual Ward identity treatment, which reduces the length of chain by one, meaning it can be handled inductively.
 For the rest of this informal discussion we focus on the regularized part.

In the concrete $k=2$ case, we define $\reg{\mathcal{X}}^{k}_t$ to be
\begin{equation} \label{eq:Greg}
	\reg{\mathcal{X}}_t \equiv \reg{\mathcal{X}}^{k}_t(x_1, x_2) :=\Theta_t^{(2)}\biggl[\Tr\bigl[ 
	G_1 {S}^{x_1} G_2 - M_{[1,2]}\big]_t
	\reg{S}^{x_2}\bigr]\biggr],
\end{equation}
thus the  linear term $\Theta_t^{(j)}[\mathcal{X}_t\big]$ for $j=2$ in~\eqref{eq:12} becomes essentially $\reg{\mathcal{X}}_t $
up to less relevant terms  involving shorter chains.  
We now claim that $\reg{\mathcal{X}}^{k}_t$ is smaller than the original $\Theta_t^{(2)}[\mathcal{X}_t\big]$
when acted upon by the propagator $\mathcal{P}_{s,t}$.
\begin{equation}\label{dX}
    \reg{\mathcal{X}}_s  = \sum_a\Tr \Big[ \bigl(G_1 {S}^{x_1} G_2 - M_{[1,2]}\bigr) S^a\Big]_s
     \big( (\Theta_s)_{ax_2} -  (\Theta_s)_{x_1x_2}\big), \quad s \le  t,
\end{equation}
where both terms in the summation are localized on the scale $\ell_s$, that is $|a-x_1|\sim |a-x_2| \sim \ell_s$. 
Crucially, the difference of the two $\Theta$ profiles becomes small when averaged in $x_2$ by the propagator $\mathcal{P}_{s,t}$, which lives on a larger scale $\ell_t$.  The gain from this cancellation is proportional 
to the ratio of length scales\footnote{In fact, the actual gain is $\ell_s\ell_t\eta_t/W^2$, see~\eqref{2prop}
below, which in some regimes is even smaller than $\ell_s/\ell_t$.} $\ell_s/\ell_t$,
  compensating for the large factor in \eqref{badest1}. 
This cancellation is not present in the unregularized term $\Theta_t^{(2)}[\mathcal{X}_t\big]$, where the difference 
in the second factor in~\eqref{dX} is replaced by $\Theta_{ax_2}$.

To rigorously capture the effect of regularization, we derive a separate evolution equation for $\reg{\mathcal{X}}_t$,
and extract the desired improvement from the action of its associated propagator.  
Ignoring irrelevant identity terms, the Duhamel formula analogous to~\eqref{GSGzig} takes the form
\begin{equation}\label{GSGregzig}
    \reg{\mathcal{X}}_t = \bigotimes_{j=1}^2 \mathcal{P}_{s,t}^{(j)} \circ  \reg{\Theta}^{(2)}_t \big[ \mathcal{X}_s \big]
    + \int_s^t \bigotimes_{j=1}^2 \mathcal{P}_{u,t}^{(j)}  \circ  \reg{\Theta}^{(2)}_t\Big[\mathrm{d} \mathcal{M}_u 
    +  \mathcal{F}_u \mathrm{d} u + \other{\mathcal{F}}_u \mathrm{d} u\Big].
 \end{equation}
Here, $\reg{\Theta}$ is the regularization of the generator, in this case $\reg{\Theta}_{ab} = \Theta_{ab} -\Theta_{ax_1}$.
Note that we picked up an additional forcing term $\other{\mathcal{F}}$ from commuting the generator of $\mathcal{X}$ from~\eqref{eq:12} through the regularization. 
 
The key difference between~\eqref{GSGzig} and the unregularized evolution~\eqref{GSGregzig} lies in the appearance of the  regularized operator~$\reg{\Theta}^{(2)}_t = \mathcal{P}_{u,t}^{(2)} \circ \reg{\Theta}^{(2)}_u$, which encodes a cancellation effect arising from the spatial regularity of the $\Theta$-profile, as described above.  

For instance, considering the contribution from the initial condition (i.e., the first term in~\eqref{GSGregzig}), we have the bound (see~\eqref{eq:one_ring})
\begin{equation}\label{badest2}
  | \mathcal{X}_s (x_1,x_2)| \lesssim \frac{1}{\ell_s\eta_s} (\Upsilon_s)_{x_1x_2} \quad \Longrightarrow \quad
   \Big( \bigotimes_{j=1}^2 \mathcal{P}_{s,t}^{(j)} \circ  \reg{\Theta}_t^{(2)} \big[ \mathcal{X}_s\big]\Big)_{x_1x_2} 
   \lesssim
   \frac{\ell_t}{\ell_s}\,  \frac{\ell_s\ell_t\eta_t}{W^2} \,\frac{1}{\eta_t} \frac{1}{\ell_t\eta_t} (\Upsilon_t)_{x_1x_2}.
\end{equation}
Compared to~\eqref{badest1}, the additional factor $1/\eta_t$ reflects the natural size of the bare $\Theta_t$;
however, due to regularization, the effective size of $\reg{\Theta}_t$ is reduced by a crucial improvement factor $\ell_s\ell_t\eta_t/W^2$.
Since this satisfies $\ell_s\ell_t\eta_t/W^2\le \ell_s/\ell_t$ (in fact in the regime $\eta_t\ge (W/N)^2$ they are equal),
 this factor compensates for the excess factor $\ell_t/\ell_s$ in~\eqref{badest1}. As a result, we recover
  the consistent estimate 
\begin{equation}\label{direct}
     \big| \reg{\mathcal{X}}_t(x_1,x_2)\big|\le \frac{1}{\eta_t}  \frac{1}{\ell_t\eta_t}(\Upsilon_t)_{x_1x_2}.
\end{equation}
While we have explained the regularization gain only for the initial condition term, a similar improvement holds for both the martingale and forcing terms in~\eqref{GSGregzig}. 

Having established the bound~\eqref{direct}, we can now treat the regularized component $\reg{\mathcal{X}}_t$ of the linear term $\Theta^{(2)}[\mathcal{X}]$ in~\eqref{eq:12} as an inhomogeneous forcing term, using its direct estimate~\eqref{direct}. 
In this way, we essentially reduced the two generators in~\eqref{eq:12} to a single $(\frac{1}{2}I+\Theta^{(1)})$,
and the Duhamel’s formula is now applied only for $\mathcal{P}^{(1)}$ rather than $\mathcal{P}^{(1)}\otimes \mathcal{P}^{(2)}$. 
In the analog of the problematic estimate~\eqref{badest1}, this effectively removes the excess factor $\ell_t/\ell_s$.

The essence of this improvement is already visible in the following toy example.  Suppose a time-dependent function $f_s(x)$
lives on scale $\ell_s$ and
satisfies  $\| f_s\|_\infty\le C/(\ell_s\eta_s)$ 
 with some constant $C$. Then, using the scale and regularity properties of $\Theta_t$, 
an elementary calculation yields
\begin{equation}\label{2prop}
 \|\Theta_t[f_s]\|_\infty \le \frac{1}{\eta_s} \frac{C}{\ell_t\eta_t}, \qquad  \|\reg{\Theta}_t[f_s]\|_\infty \le
  \frac{\ell_s\ell_t\eta_t}{W^2} \times 
  \frac{1}{\eta_s} \frac{C}{\ell_t\eta_t}.
\end{equation}
Recalling from~\eqref{eq:propapprox} that the true propagator $\mathcal{P}_{s,t}$ is essentially $\eta_s\Theta_t$, it is easy to reconstruct the gain from regularization on the actual propagator
from~\eqref{2prop}.  In the regime $\eta_s\ge \eta_t\ge (W/N)^2$, the improvement factor satisfies
$$
     \frac{\ell_s\ell_t\eta_t}{W^2} = \frac{\ell_s}{\ell_t} = \Big(\frac{\eta_t}{\eta_s}\Big)^{1/2},
$$
confirming  that the gain comes from the ratio of length scales on which $f_s$ and $\Theta_t$ live.
In the regime $\eta_t\le \eta_s\le (W/N)^2$, the gain factor becomes even better---at least a full power $\eta_t/\eta_s$.
 
So far we showed that regularizing a single observable restores the consistency of the basic power counting in Duhamel’s formula. 
However, the full proof must still address the truncation of the hierarchy, which introduces a loss—at least at the longest chain length $K$.
For technical reasons, this loss is redistributed across chains of length between $K/2$ and $K$, meaning that mere consistency is no longer sufficient in this range. 
To compensate, we require an additional improvement for chains of length $ \ge K/2$. 
Since regularization is a local procedure, it can be applied to every second observable, yielding independent improvement factors of $\ell_s\ell_t\eta_t/W^2$ from each of them if necessary. 
It turns out that  \emph{two regularizations} already suffice to balance the redistributed loss---hence a minimal chain length of $4$ is required, implying  $K/2\ge4$ as a necessary lower bound on the maximal chain length. 
The actual proof has many other subtleties,
 but multiple observable regularization remains the most crucial ingredient.

 In this informal summary, we focused on the averaged chains where the effect of regularization plays the most important role. 
 However, the full argument necessarily involves isotropic chains as well. The averaged and isotropic master inequalities are fundamentally coupled, albeit for different reasons.
 In the isotropic analysis, the main challenge is the quadratic variation for long chains, which must be reduced to shorter, averaged chains where better control is available. 
 This way the isotropic master inequalities can be closed without regularization.
 Conversely, in the averaged analysis, isotropic quantities emerge due to the general form of the variance matrix  $S$. To illustrate this, consider one term from $\mathcal{F}$ in~\eqref{eq:12} of the form (see~\eqref{eq:1_k+1_bound1_1})
 \begin{equation}\label{fact}
     \Tr \Big[ \mathcal{S} \big[ G_1-m_1\big] \big(  G_{[1,2]}S^{x_2} G_1-M \big)  \Big]
     =\sum_q \Tr \big[ (G_1-m_1)S^q \big] \big(  G_{[1,2]}S^{x_2} G_1-M \big)_{qq},
 \end{equation}
 which cannot be expressed purely in terms of tracial quantities. 
 If $S$ were factorized as $S_{ab} = \sum_c \other{S}_{ac}\other{S}_{cb}$ with $\other{S}$ having
 similar scale and decay properties to $S$, then ~\eqref{fact} could be rewritten as
 $$
    \sum_c \Tr \big[ (G_1-m_1)\other{S}^c \big] \Tr \big[ \other{S}^c \big(  G_{[1,2]}S^{x_2} G_1-M \big)\big],
 $$
 essentially allowing to run the averaged analysis alone.
  This factorization holds for the special block-constant variance profile used in~\cite{yauyin}, explaining why that proof foregoes substantial isotropic analysis. Although a weak isotropic estimate without spatial dependence is still derived in~\cite{yauyin} 
  separately via the Schur complement, such bounds would likely be insufficient to control terms like~\eqref{fact} without factorization. No such simplification is available for general variance profiles $S$.
 
 Beyond these structural considerations, there is a conceptual reason to include a full isotropic analysis: it is essential for the Green Function Comparison (GFT) argument (zag-step), which allows to depart from purely Gaussian ensembles. Even if one aims to estimate only averaged quantities, isotropic chains still emerge from terms involving third and higher cumulants in any known GFT method, and require quite precise control. The model in~\cite{yauyin} was fully Gaussian, so this complication did not arise.

 \subsubsection{Deterministic $M$-bound}\label{sec:Mb}
A core element of the zigzag strategy is to establish sharp estimates on the deterministic $M$-terms, since these ultimately steer all fluctuation bounds.
As discussed in Section~\ref{sec:staticmaster}, their defining recursion exhibits the same stability issues as the static master inequalities: it features many cancellations that are not easy to discern. 
An added challenge in the RBM setting is that not only the size of the $M$-terms, but also their spatial structure must be accurately followed.
 Instead of identifying the cancellations in the $M$-terms directly, we propose a significantly shorter, purely \emph{dynamical} approach to deriving the $M$-bounds, leveraging the observable regularization idea and the associated zig machinery developed in this work.   
 
Let $X_t=X_t^k: = \Tr\bigl[M_{[1,k],t}(\bm x') S^{x_k}\bigr]$. 
Recall that $\mathcal{X}_t^k$, analyzed in the previous section, is just the fluctuation of the averaged $k$-chain around $X_t^k$.
The deterministic quantity $X_t^k$ satisfies a similar evolution equation of the form
\begin{equation}\label{eq:av_M_evol1}
	\frac{\mathrm{d}}{\mathrm{d}t}X_t^k = 
	\bigoplus_{j=1}^k \Big( \frac{1}{2}I  +  \Theta \Big)^{(j)} 
	\bigl[X_t^k \bigr] +  F_{t}^k,
\end{equation}
which was already implicitly used in deriving equation~\eqref{eq:12} for $\mathcal{X}_t^k$.

Unlike its stochastic counterpart, equation~\eqref{eq:av_M_evol1} has no martingale term, and, moreover, the forcing term $F_t^k$ is much simpler.
However, in one crucial respect, it is \emph{worse}: the issue of excess propagators is more pronounced, because the natural size for $X$ is larger by a factor of $\ell_t \eta_t$ than that of $\mathcal{X}$.
Concretely, the target $\ell^\infty$-bound becomes
\begin{equation} 
	\| X_t\| \lesssim \frac{1}{(\ell_t\eta_t)^{k-1}}, 
\end{equation}
which replaces the estimate~\eqref{target} (written for $k=2$). 
Taking into account the fact that $X_s$ lives on scale $\ell_s$, the naive propagation yields
\begin{equation} 
	\Big\| \bigotimes_{j=1}^k \mathcal{P}_{s,t}^{(j)} [X_s] \Big\| \lesssim  \Big(\frac{\ell_t}{\ell_s}\Big)^2
	\frac{1}{(\ell_t\eta_t)^{k-1}}, 
\end{equation}
i.e., it is off by a factor $(\ell_t/\ell_s)^2$, compared with the single factor $\ell_t/\ell_s$ in~\eqref{badest1}.
To compensate for it, we need \emph{two} regularizations.
Fortunately, any $k\ge 4$ chain is sufficiently long to accommodate two independent regularizations. 
For $k=3$, the chain is not saturated (since it is odd), so not all propagators are critical. 
The case $k=2$ must be handled separately; here $M_{[1,2]}$ is essentially $\Theta$, for which we give an explicit estimate in Section~\ref{sec:reg}.
 
Once the two regularizations are in place, the analysis of~\eqref{eq:av_M_evol1} becomes much simpler than that of the corresponding fluctuation term $\mathcal{X}^k_t$. No reduction is needed, since there is no martingale term---whose quadratic variation would otherwise generate longer chains---and the forcing terms are simpler and do not involve excess length.
We emphasize, however, that two regularizations are strictly necessary for obtaining the correct $M$-bounds. By contrast, in the analysis of $\mathcal{X}^k_t$, a single regularization would suffice if optimal estimates for all chain lengths were available; the need for suboptimal reduction inequalities to truncate the hierarchy is what necessitates an additional regularization.  
 
We also note that optimal $\ell^\infty$-bounds for averaged $M$-chains were derived independently in~\cite[Section 3]{yauyin} using a more intricate expansion technique that tracks delicate cancellations arising from multiple length scales in their analogue of~\eqref{eq:av_M_evol1}. The precise spatial structure of long $M$-chains was not captured there;  this may be partly the reason why accurate spatial dependence in the local laws was only established for $2$-chains in~\cite[Theorem 2.21]{yauyin}.

 \subsubsection{Initial global law and completion of zigzag}\label{sec:informalinit}
For purely Gaussian models, the zig step---described in its most essential form in Section~\ref{sec:rego}---is sufficient. One can start from the initial condition $H_0 = 0$ in the OU flow and generate the desired Gaussian ensemble after time of order one. This simplification of the characteristic flow method in the RBM context was first introduced in~\cite{dubova} and adopted in~\cite{yauyin}. However, for non-Gaussian models, the zig step must be complemented by two further ingredients of the zigzag strategy: an initial global law at $\eta\sim 1$, and a GFT (zag) step to remove the Gaussian component.

In all previous applications of the zigzag strategy to mean field models, the initial global estimate was relatively routine and often relegated to appendices (see, e.g.,~\cite[Appendix A.2]{cipolloni2023eigenstate}, \cite[Appendix B]{Cipolloni2022Optimal}, \cite[Appendix A]{cipolloni2022rank}).
The main reason is that norm bounds on the resolvent $\| G(z)\|\le 1/|\im z|\sim 1$ are affordable, and all stability operators remain bounded. 
This automatically removes the two main complications (truncation of the
hierarchy and stability) of a typical local law proof, while reducing the higher moment expansions to elementary power counting.

While the same simplifications hold for RBM at the level of bounds, a substantial additional difficulty arises: a useful global law must capture the spatial structure of $G_{[1,k]} - M_{[1,k]}$.
The  standard bootstrapping procedure together with a cumulant expansion for the
underline term,  inherited from mean field models,  yields only
$\ell^\infty$-bounds, e.g., for a single resolvent, they are
\begin{equation}\label{global}
	\big| (G-m)_{ab}\big|\lesssim \frac{1}{(\ell\eta)^{1/2}} \lesssim W^{-1/2}, 
	\qquad   \Tr \big[  (G-m)S^x\big]  \lesssim W^{-3/4}, 
\end{equation}
where the first bound is optimal, the second is off by a factor $W^{-1/4}$.
Similar bounds can be derived for longer chains, but they all
lack the correct dependence on $a, b, x_1, x_2, \ldots , x_k$. 

To address this, in Section~\ref{sec:global}, we recover  the optimal spatial dependence by running a version of the static master inequalities described in Section~\ref{sec:staticmaster}.  
We first  perform an iterative refinement for the isotropic $k=1$ and $k=2$ chains to recover the correct spatial factor. 
Roughly speaking, instead of operating with the proper spatial factor $\Upsilon$ right away,
we instead use its uniformly-relaxed  version $\Upsilon+\omega$, where the parameter $\omega$ 
controls the size of the flat tail,   which we gradually decrease. 
The crude $\ell^\infty$-bound corresponds to the starting value 
$\omega=\| \Upsilon\|_\infty\sim
W^{-1}$ and after several iterations of the static master inequalities 
we reduce $\omega$ to $0$, where the optimal spatial control is achieved. 
Finally, the optimal bounds on $k=1$ and $k=2$ chains are used to initialize the system of static master inequalities for chains of arbitrary length (isotropic and averaged).
Their optimal spatial decay is obtained immediately, without any further  profile-recovering iterations, 
in line with the general philosophy that master inequalities transfer the quality of estimates from two-resolvent chains to longer chains.

The final component of the zigzag strategy is the zag step, the GFT argument. 
In most mean-field applications, this straightforward---with  notable exceptions~\cite{campbell2024spectral},~\cite{cipolloni2023eigenstate}
and \cite{erdHos2024cusp}, which anticipate two complications that arise simultaneously in RBM.  
First, typical GFT arguments use the single resolvent local law as an external input; instead, here we need to establish it internally through a self-consistent procedure, as in~\cite{campbell2024spectral,erdHos2024cusp}. 
Second, we must use the more refined one-by-one replacement method to prove GFT, rather than the more commonly used continuity method along the OU flow. 
This is because capturing the precise spatial structure within GFT is delicate. In fact, even the standard estimate on
each one-by-one replacement is not sufficient—replacing each matrix entry does not yield an error of $o(N^{-2})$---but by carefully tracking the spatial decay, we can ensure that the \emph{cumulative error} across all replacements remains small.
A similar idea was implemented in~\cite{cipolloni2023eigenstate} (also adapted in~\cite{erdHos2024cusp})
in a simpler setting, where multi-resolvent local laws for Wigner matrices were established with optimal control in the Hilbert–Schmidt norm of the observables that required nontrivial estimates on the off-diagonal structure.

 \subsubsection{Extensions}\label{sec:extension}
 
 We briefly discuss three extensions of our main proof; 
   general observables and isotropic laws (Section~\ref{sec:general}), traceless observables (Section~\ref{sec:traceless}), and the
 real-symmetric case (Sections~\ref{sec:real}).  Surprisingly, they are not just simple technical 
 extensions, but pose genuine difficulties that require separate, and in places quite extensive, proofs.
 
\medskip

{\bf General observables and test vectors.} 
In Section~\ref{sec:general}, we start with extending our chains to include general diagonal observables (beyond the special ones 
of the form $S^x$) and upgrade entry-wise estimates to full isotropic laws.
The key difficulty is that the generalized control functions, $\Upsilon_{\bm u \bm v}$ and 
$\Upsilon_{AB}$ from~\eqref{eq:genUp}--\eqref{eq:genUpA}, are naturally linear in $A, B$ and  
quadratic in $\bm u, \bm v$,  yet they appear under square root in the fluctuation bounds.
Roughly speaking, this is because the spatial dependence of a chain 
is designed by schematically writing, e.g.,
$$
   \langle\bm u, GAGBG\bm v\rangle \sim (\bm u G \sqrt{A})(\sqrt{A}G\sqrt{B})(\sqrt{B}G\bm v)
   \lesssim \sqrt{\Upsilon_{\bm u A}\Upsilon_{AB} \Upsilon_{B\bm v}},
$$
where  $\Upsilon$ naturally controls the \emph{variance} of each constituent  factor. For example, the 
size of $(\bm u G \sqrt{A})$ is controlled by the square root of the size of $\langle\bm u, GAG\bm u\rangle$ which is
$\Upsilon_{\bm u A}$. This square root, however, hinders optimal  concatanations of such estimates;
we will call it the \emph{$\sqrt{\Upsilon}$-problem.}
Most importantly,  the analogue of~\eqref{eq:convol_notime} does not hold in the generalized form (take $\eta_1=\eta_2$
for simplicity),
\begin{equation}\label{con}
 \sum_a \sqrt{\Upsilon_{\bm u a} \Upsilon_{a\bm v} } = \sum_a \sqrt{\sum_{i,j}  |u_i|^2 \Upsilon_{ia} \Upsilon_{aj} |v_j|^2}
 \not \le \sqrt{\frac{\ell}{\eta}} \sqrt{ \sum_{i,j}  |u_i|^2 \Upsilon_{ij}|v_j|^2}
   = \sqrt{\frac{\ell}{\eta}} \sqrt{ \Upsilon_{\bm u \bm v}},
\end{equation}
even though the inequality holds for $\bm u,\bm v$ being coordinate vectors. This is not just a minor technical 
problem, but connected with the  phase cancellation in  sums like
$  G_{\bm u \bm v} = \sum_{ij} u_i G_{ij} v_j$ or $(GAG)_{aa} =\sum_i G_{ai}A_{ii} G_{ia}$ 
which can only be correctly controlled by computing the variance,
leading to the square root problem.

Following our main proof, we realized that in most instances only a weaker form of the convolution inequality~\eqref{con}
without any square roots is needed, which still holds for the general $\Upsilon$-factors. We identified the terms, where~\eqref{con} 
were really needed, and paid 
the price of an additional large factor $\sqrt{N/\ell}$ in the right hand side which makes~\eqref{con}
valid by a simple Schwarz estimate. In short, resolving the $\sqrt{\Upsilon}$-problem costs an extra $\sqrt{N/\ell}$ factor.
We then gain back this factor again by an additional observable regularization.

\medskip

{\bf Traceless observables.}
Similarly to Wigner matrices, traceless observables make resolvent chains and their fluctuations smaller, but
only in the regime $\eta\ll (W/N)^2$.  As a complete
analogue of the $\sqrt{\eta}$-rule in the Wigner case~\cite{cipolloni2021eigenstate}, in Section~\ref{sec:traceless} 
we show that an improvement factor of order $\theta:=\sqrt{\eta}(N/W)\ll 1$ is gained
for each traceless observable. The source of this improvement is 
the \emph{mean-field regularization}, i.e. the mechanism that
replacing one special observable $S^x$ with its traceless version,
$\trless{S}^x: = S^x - \frac{1}{N}I$, effectively removes the constant mode of the propagator acting on $x$,
and thus the gain is proportional to the ratio of the two smallest eigenvalues.   The modified propagator is roughly
$$
\widetilde{\mathcal{P}}_{s, t} 	\approx I + \eta_s \trless{\Theta}_t, \qquad \trless{\Theta}_t: =\Theta_t - \frac{1}{N}\Tr \big[\Theta_t | {\bf 1}\rangle\langle {\bf 1}|\big],
$$
and in the relevant regime, $\eta_t\le\eta_s\le (W/N)^2$, it is bounded, $\| \widetilde{\mathcal{P}}_{s, t}\|\lesssim 1$,
in contrast to the large $\eta_s/\eta_t$ bound~\eqref{Pnorm}. Inserting
this improved propagator estimate for each traceless
observable into the analysis of the Duhamel’s formula for the fluctuation terms $G-M$,
we obtain a consistent improvement factor 
$\theta_t =\sqrt{\eta_t}(N/W)$.

The actual proof follows our main proof  with  several additional steps. We run
a double induction scheme, where the main induction step is on the number of traceless observables, while
the second induction is a bootstrap where the full $\theta$ factor is gained in small steps of order $(N\eta)^{-1/2}$.
In both steps the maximal chain length decreases; the stronger bound is proven for shorter chains than
the induction hypothesis. This bootstrap approach is convenient: since we already have bounds for chains
of arbitrary length with  fewer traceless observables (in particular, we have bounds
for chains with general observables as a starting point), we can use them  to  estimate various
terms in the corresponding master inequalities without worrying about closability of  the hierarchy in 
terms of chain length\footnote{
We note that this direct approach differs from the analogous proofs
in the Wigner case, where the traceless hierarchy was closable hence no direct inflation of chain lengths were present
(see, e.g.~\cite{Cipolloni2022Optimal} for the static hierarchy or \cite{cipolloni2023eigenstate} for the dynamical hierarchy).
A tighter control of the traceless master inequalities in the current RBM problem may also lead to a similar, 
somewhat  more conceptual proof.}. 
The most important term, the quadratic variation of the martingale, roughly
in the form of $Q:=\Tr \big[ G_{[1,k+1]}S^x G_{[1,k+1]}^*\big]$, however, still needs a reduction 
not because of its doubled length, but because it also doubled the number $n$ of traceless observables;
note that closability in $n$ is still needed.
The difficulty is that reduction cannot be performed at traceless observables without losing the traceless effect,
and one reduction at $S^x$ in the middle of $Q$ is not  sufficient for self-improvement. This forces us
to study a class of symmetrized chains of the form $Q$ separately, by analyzing their own evolution equations.
Since $Q$ already has one non-traceless observable, $S^x$, its quadratic variation will have three,
giving us enough freedom for the necessary two reductions without destroying any traceless observable.
\medskip

{\bf Real symmetric case.} All our results and proofs are extended from the complex Hermitian case to the real symmetry class
in Section~\ref{sec:real}.  This means that the self-energy operator $\mathcal{S}[R]:= \Expv[HRH]$,
 besides the usual diagonal tracial term
$(\mathscr{S}[R])_{ab} = \delta_{ab}\Tr [S^a R]$, has an additional offdiagonal term 
$(\mathscr{T}[R])_{ab} = (S^\mathrm{od})_{ab} R_{ba}$, where $ S^\mathrm{od}$ is the offdiagonal part of $S$.
In mean field models, the effect of this term is typically negligible just by power counting the size,
 since $(S^\mathrm{od})_{ab}$ is order $1/N$
and the offdiagonal elements of resolvent chains, the typical argument $R$ of $\mathscr{T}$, is small.
In contrast, our RBM analysis is very sensitive to the spatial decay and structurally the offdiagonal
 $(\mathscr{T}[R])_{ab}$ contains a very different spatial relation between $a$ and $b$ than the diagonal 
 $(\mathscr{S}[R])_{ab}$. This creates new terms, called \emph{skewed forcing terms}, that require
 separate treatment for general observables and for genuinely isotropic chains. In short, 
 the $\sqrt{\Upsilon}$-problem occurs in more terms than before and they are more critical.
 We resolve this problem by increasing and maximal chain length $K$ and
 adjusting the loss exponents to accommodate the $\sqrt{N/\ell}$ loss
 from the $\sqrt{\Upsilon}$-problem.

 \medskip

We remark that~\cite{yauyin} focused only on the complex Hermitian case,
 on special observables without general isotropic law and without analyzing traceless observables;
  therefore none of these complications
 were present.  Nevertheless,~\cite{yauyin} contains a weak version of the  upper bound
 on the averaged two-resolvent chain with two traceless observables
  since it is 
 necessary both for Wigner-Dyson universality and for their weak version of the QUE.
 First, it is shown in
 ~\cite[Theorem 2.4]{yauyin}\footnote{For consistency,
 here we use our notation, note that our $\ell$ and $S^x$ correspond to $W\ell$ and $E_x$ in~\cite{yauyin},
 respectively.} that
 $$
    \Big| \Expv \Tr \Big[ \big( GS^xG^* - M_{[1,2]}\big) S^y\Big] \Big| \lesssim \frac{1}{\ell \eta} \frac{1}{(\ell \eta)^2},
$$
which carries an additional $(\ell\eta)^{-1}$ improvement over the  local law~\eqref{eq:avelaw} (without spatial decay)
due to taking expectation. The same bound holds for the traceless observables $\trless{S}$ as well. 
Direct calculation shows that the deterministic term exhibits the $[\sqrt{\eta}(N/W)]^2$ 
factor, yielding the estimate
$$
    \Big| \Expv \Tr \big( G\trless{S}^xG^* \trless{S}^y\big)\Big| \le \Tr \big( M_{[1,2]} S^y\big) + \frac{C}{(\ell\eta)^3}
    \lesssim \frac{1}{(\ell\eta)^2}\Big[ \Big( \frac{\sqrt{\eta}N}{W}\Big)^2+   \frac{1}{\ell\eta} \Big].
$$
 Note that the second error term is suboptimal in the relevant $\eta\le (W/N)^2$ regime and the bound holds only in
 expectation sense.

\section{Zigzag strategy}\label{sec:zigzagsec}

 The key idea of the zigzag strategy is to embed the multi-resolvent local law problem 
at a fixed constellation of spectral parameters ${\bm z}$ into a flow of spectral parameters ${\bm z}_t$.
Simultaneously with moving the spectral parameters, the matrix $H$ itself is
also embedded into a stochastic matrix flow $H_t$ in such a way that the 
resulting flow of resolvents $G_t=(H_t-z_t)^{-1}$ exhibits a remarkable cancellation.
This feature allows us to propagate a bound from a larger $\eta$ to a smaller $\eta$
at the expense of adding a small Gaussian component; this is the zig-step.
It is complemented by the zag-step, where the Gaussian component is removed
by a moment matching argument and  by an initial global bound in order to start the procedure.
Each zig-step is followed by a zag-step and the procedure is repeated finitely many times
until the smallest $\eta\sim N^{-1+\etaexp}$ is achieved.

First we introduce the necessary ingredients in Section~\ref{sec:zigzag_setup}. 
In  Section~\ref{sec:zz} we then formalize  the three key steps
of the zigzag strategy, namely the initial global bound, the zig-step and the zag-step
as separate propositions and conclude from them the proof of the local laws, Theorem~\ref{th:local_laws}.
In the subsequent sections then we prove these three key propositions.
The most important zig-step is proven in Section~\ref{sec:masters_sec}--\ref{sec:reg_sec}
with auxiliary results deferred to Section~\ref{sec:aux}. In separate Section~\ref{sec:traceless} we extend 
the zig-step to traceless observables.
The initial global bound is proven in Section~\ref{sec:global}, while the zag-step in Section~\ref{sec:GFT}. 
Along the proof we need various bounds on the deterministic approximation $M$, these are
proven in Section~\ref{sec:M}. 

Furthermore, in the main proof we only consider special observables of the form $A_i := S^{x_i}$
and its traceless version $\trless{S}^{x_i}: = S^{x_i}- \frac{1}{N}I$,
and instead of full isotropic laws with arbitrary test vectors $\bm u, \bm v$ we consider only
coordinate vectors (so-called entry-wise local laws). Once the proof in this special case is
completed, in a separate Section~\ref{sec:gen_tr} we will explain how to remove these
restrictions.

\subsection{Setup and notations for the zigzag strategy}\label{sec:zigzag_setup}

\subsubsection{Self-energy operator $\mathcal{S}$} \label{sec:calS}
We  recall the definition of  the  self-energy  super-operator $\mathcal{S}$ of the random matrix $H$
from~\eqref{def:mathcalS}. For matrices with independent entries, 
  $\mathcal{S}$ admits the decomposition
\begin{equation}\label{eq:Sdecomp}
	\mathcal{S}[R] := \Expv\bigl[H \, R \, H\bigr] = \mathscr{S}[R] + \mathscr{T}[R],  \quad R \in \mathbb{C}^{N\times N},
\end{equation}
where $\mathscr{S}[R]$ and $\mathscr{S}[R]$ represent the diagonal and off-diagonal components of $\mathcal{S}$, respectively, acting on $R$ as follows:
\begin{equation}\label{eq:entries}
	\bigl(\mathscr{S}[R]\bigr)_{ab} := \delta_{ab} \sum_j S_{aj} R_{jj} =\delta_{ab} \Tr \big[S^a R\big]  ,
	 \qquad \bigl(\mathscr{T}[R]\bigr)_{ab} := \bigl(S^\mathrm{od}\bigr)_{ab} \, R_{ba} , \quad a,b \in \indset{N}.
\end{equation}
Here the entries of the $N\times N$ off-diagonal matrix $S^\mathrm{od}$ are given by
\begin{equation}
	 \bigl(S^\mathrm{od}\bigr)_{ab} := \delta_{a\neq b}\Expv [(H_{ab})^2] = \delta_{a\neq b}S_{ab}\Expv \bigl[(h_{ab})^2\bigr], \quad a,b\in\indset{N}.
\end{equation}
As their names suggest, $\mathscr{T}$ acts exclusively on the off-diagonal entries of $R$, while $\mathscr{S}$ acts only on the diagonal entries, moreover $\mathscr{S}[R]$ is  diagonal, while $\mathscr{T}[R]$ is an off-diagonal matrix.
 In coordinate-free notation, the off-diagonal super-operator $\mathscr{T}$ can be expressed as
\begin{equation}\label{eq:calT}
	\mathscr{T}[R] = S^\mathrm{od}\odot R^\mathfrak{t},
\end{equation}
where $\odot$ denotes the Hadamard (entry-wise) product, and $(\cdot)^\mathfrak{t}$ denotes transposition. 
In the rest of the paper, except Section~\ref{sec:real}, we will always assume that $\mathscr{T} \equiv 0$.
In Section~\ref{sec:real} we will explain the necessary modifications if $\mathscr{T} \not\equiv 0$,
in particular we handle the real symmetry class.

\subsubsection{Characteristic flow}\label{sec:charflow} 
The characteristic flow $z_t$ is the solution to the simple ODE on the complex plane:
\begin{equation} \label{eq:char_flow}
	\frac{\mathrm{d}z_t}{\mathrm{d}t} = -\frac{1}{2}z_t - m(z_t).
\end{equation}
 Instead of an initial condition $z_0=z$ we will usually specify a final time $T$
together with a target spectral parameter $z_T\in \mathbb{C}$, we solve the equation \eqref{eq:char_flow} backwards
in time until $t=0$, and then we consider the flow~\eqref{eq:char_flow} running
forward  from  $z_{0}$ at time $t=0$
up to the final point $z_T$ at $t=T$. We always assume that $z_T\in \fuldom=\fuldom_{\kappa, \etaexp, C_0}$ for some fixed
parameters $\kappa, \etaexp, C_0$, see~\eqref{def:spectraldomain}.

We set
\begin{equation}\label{def:etat}
   \eta_t := |\im z_t|
\end{equation}
as our basic parameter controlling the distance of $z_t$ from the real axis. Since 
$z_T\in \fuldom=\fuldom_{\kappa, \etaexp, C_0}$, we have 
$ N^{-1+\etaexp} \le \eta_T=|\im z_T|\lesssim1$. In particular,
\begin{equation}\label{eq:etabound}
N^{-1+\etaexp}\le \eta_t\lesssim 1, \quad  t \in [0, T], 
\end{equation}   
 since $\eta_t$ is monotone decreasing. 
Note that $\im z_t$ has the same sign as $\im z_T$
for all $t\le T$.

A direct calculation shows  that if the  final condition $z_T$ is 
chosen such that $z_T\in \fuldom_{\kappa, \etaexp, C_0}$  
for some small fixed $\kappa,\etaexp>0$, then
\begin{equation}\label{eq:etatasymp}
\eta_t - \eta_T \sim (T -t), \qquad 0\le t\le T,
\end{equation}
where the implicit constant in $\sim$ depends on $\kappa$ and $T$.  In particular, by choosing $T\sim 1$
sufficiently large,
we achieve that for the initial condition we have $\eta_0=|\im z_0| \ge 1$.

Furthermore,  if a $k$-tuple of spectral parameters  $(z_1, z_2,\ldots, z_k)$ satisfies~\eqref{eq:admissible_z}
at a final time $T$, i.e. 
$$
 \max_j |\im z_j |\le \other{C}_0 \min_j |\im z_j |,  
 $$
then the same relation holds for all times
\begin{equation}\label{eq:comparability}
 \max_j |\im z_{j,t} |\le \other{C}_0' \min_j |\im z_{j,t}|,  \quad t\in [0, T],
 \end{equation}
 with a possibly larger constant $\other{C}_0'$ that is irrelevant for the proof.
 Here $z_{j,t}$ is the solution to \eqref{eq:char_flow} with final condition $z_{j,T} = z_j$.
 Finally, setting $m_t:= m(z_t)$, notice that
$$
      \frac{\mathrm{d}}{\mathrm{d}t} m_t= \frac{1}{2} m_t.
$$
From this relation it follows that $m_t/|m_t|$ is constant along the time evolution, in particular, the 
backward flow maps the domain $\fuldom$ into itself, in the sense that
 \begin{equation}\label{eq:edge}
   z_{j, T} \in \fuldom_{\kappa,\etaexp, C_0} \quad \Longrightarrow \quad z_{j,t }\in \fuldom_{\kappa, \etaexp, C_0'}, \quad t\in [0,T].
 \end{equation}

 Set $\Xi_t:= \Xi(z_{1,t}, z_{2,t})$  and $\Theta_t:= \Theta(z_{1,t}, z_{2,t})$, then a by direct calculation we have
\begin{equation} \label{eq:dTheta}
    \frac{\mathrm{d}}{\mathrm{d}t}\Theta_t= (I+\Theta_t)\Theta_t, \qquad   
    \frac{\mathrm{d}}{\mathrm{d}t}\Xi_t= (I+\Xi_t)\Xi_t.
 \end{equation}

  If $\eta_t$ is given by \eqref{def:etat},  we set 
 \begin{equation}\label{def:ell_def}
 \ell_t := \ell(\eta_t) = \min\Big\{ \frac{W}{\sqrt{\eta_t}}, N \Big\}. 
 \end{equation}
   We will often encounter the product $\ell_t\eta_t  \gtrsim \min\{W\sqrt{\eta_t}, N\eta_t\}  \gtrsim N^{\etaexp} $, where the first inequality is a consequence of \eqref{def:ell_def}, 
 and the second inequality follows from  our basic assumptions \eqref{eq:WN} and \eqref{eq:etabound}. 
 An easy  calculation
 using \eqref{eq:etatasymp} shows the following
{\it integration rules}
 \begin{equation} \label{eq:int_rules}
 	\int_{\tinit}^t \frac{\mathrm{d}s}{\eta_s} \lesssim \log N,\quad \int_{\tinit}^t \frac{\mathrm{d}s}{\eta_s^{1+\gamma}} \lesssim \frac{1}{\eta_t^{\gamma}} , \quad \int_{\tinit}^t \frac{1}{(\ell_s\eta_s)^\gamma}\frac{\mathrm{d}s}{\eta_s} \lesssim \frac{1}{(\ell_t\eta_t)^\gamma}, \quad \gamma > 0
 \end{equation}
 for any times $0\le \tinit\le t\le T$ recalling that $\eta_T \gtrsim 1/N$.
  
\subsubsection{Stochastic matrix flow}
 
Next, we define a stochastic matrix flow. The simplest way to formulate it is to 
write $H$ in the form of an entry-wise (Kronecker) matrix product:
\begin{equation}\label{eq:kronecker}
   H= \sqrt{S} \odot h,
\end{equation}
where $\sqrt{S}$ is the entry-wise square root, i.e. $(\sqrt{S})_{ab}: = \sqrt{S_{ab}}$ and $h \in \mathbb{C}^{N\times N}$
denotes the {\it normalized matrix}\footnote{By somewhat deviating from standard conventions, in this paper $h_{ab}$ denotes the matrix elements of the normalised matrix $h$; while the matrix elements of $H$
 will be denoted by $H_{ab}$.}, defined 
as $h_{ab} : = H_{ab}/\sqrt{S_{ab}}$ if $S_{ab}\neq 0$. The definition of $h_{ab}$ for index pairs $(a,b)$ with $S_{ab}=0$
is irrelevant, for completeness we set them to be  independent 
standard real or complex Gaussians according to the symmetry class of $H$. In this way $h$ is a Wigner matrix,
i.e. $\Expv h_{ab}=0$, $\Expv |h_{ab}|^2=1$,
$\Expv |h_{ab}|^p \le \nu_p$ and
$\Expv h_{ab}^2=0$ in the complex case.

Fix a time  $\tinit\in[0,T]$, and let $h_t$ be the solution to the  stochastic differential equation 
	\begin{equation} \label{eq:zigOU}
		\mathrm{d}h_t = -\frac{1}{2}h_t \mathrm{d}t + \mathrm{d}\mathfrak{B}_t, \quad t \in [\tinit,T], \quad h_{\tinit} = h,
	\end{equation}
where $\mathfrak{B}_t$ is the standard Brownian motion in the space of $N\times N$ matrices
	of the same symmetry class as $H$. Finally, we set\footnote{We remark that one may define $H_t$ directly as
	the solution of 
$$
\mathrm{d}H_t = -\frac{1}{2}H_t \mathrm{d}t + \Sigma^{1/2}\bigl[\mathrm{d}\mathfrak{B}_t\bigr], \quad t \in [\tinit,T], \quad H_{\tinit} = H,
$$
	where $\Sigma$ is the covariance tensor 
	of the random matrix $H$, defined by its action on $N\times N$ matrices $R$ as
	$\bigl(\Sigma[R]\bigr)_{ab} := S_{ab} R_{ab}$.  
The covariance tensor $\Sigma$ is positive definite as a super-operator on the space of $N\times N$ matrices equipped with the Hilbert-Schmidt scalar product, and $\Sigma^{1/2}$ denotes its operator square root.}
$$
H_t: = \sqrt{S}\odot h_t.
$$

Fix a chain length $k\in \mathbb{N}$.
For every $i \in \indset{k}$, we define the time-dependent resolvents $G_{i,t}$ as 
\begin{equation}
	G_{i,t} \equiv G_t(z_{i,t}) :=  \bigl(H_t - z_{i,t}\bigr)^{-1},
\end{equation}
where $z_{i,t}$ is the solution to the characteristic flow \eqref{eq:char_flow} with initial condition $z_i$.
We denote ${\bm z}_t:= (z_{1,t}, z_{2,t}, \ldots, z_{k,t})$ the vector of these solutions.
Now we introduce the time dependent versions of the resolvent chains and their deterministic approximations
analogously to~\eqref{eq:resolvent_chains_notime}--\eqref{eq:Mt_def_notime}.

For a vector of observables  $(A_1,\dots, A_{k-1})$ and for $p < j \in \indset{k}$, we set
\begin{equation} \label{eq:resolvent_chains}
	G_{[p,j],t}  
	\equiv G_{[p,j],t}(\bm z_t; A_{p}, \dots, A_{j-1}) 
	:= \biggl(\prod_{i = p}^{j-1} G_{i,t} A_i\biggr)G_{j,t}, \quad \mbox{and} \quad G_{[j,j],t} := G_{j,t}.
\end{equation} 
The  subscript $[p, j]$ indicates the implicit dependence of $G_{[p,j],t}$ on the time-dependent spectral parameters 
$(z_{p,t}, \dots, z_{j,t})$.  
The corresponding deterministic approximation matrix is  
\begin{equation} \label{eq:Mt_def}
	M_{[p, j],t} \equiv M_{[p, j],t}(\bm z_t; A_p, \dots, A_{j-1}) := M(z_{p,t}, A_p, z_{p+1,t}, \dots, z_{j-1,t}, A_{j-1}, z_{j,t}).
\end{equation}
Similarly, by convention, $M_{[j,j],t} := m (z_{j,t})$ for $j \in \indset{k}$. 
 
\subsubsection{Specialized observables $A=S^x$}

Now we specialize these notations for the observables $S^x$, defined in~\eqref{eq:S_obs},
 since our main proof will be carried out first in this case.
In Section~\ref{sec:traceless}  
we show how to account for the additional smallness factor $(N/W)\sqrt{\eta}$ for each  traceless observable $\trless{S}^x
=S^x-\frac{1}{N}I$ instead of $S^x$
in the delocalized regime.

Fix a length $k \in \mathbb{N}$, and a vector of spectral parameters $\bm z := (z_1,\dots, z_k) \in(\fuldom)^k$. 
Let $\bm x\in \indset{N}^k$ denote a vector of $k$ external indices, then we consider observables $A_i := S^{x_i}$ for all $i \in \indset{k}$, and denote, for all $p \le j \in \indset{k}$,
\begin{equation} \label{eq:Gk_def} 
	G_{[p,j],t} \equiv G_{[p,j],t}(x_{p}, \dots, x_{j-1}) := G_{[p,j],t}\bigl(\bm z_t, S^{x_p}, \dots, S^{x_{j-1}} \bigr),
\end{equation}  
\begin{equation} \label{eq:Mk_def}
	M_{[p,j],t}  
	\equiv M_{[p,j],t}(x_{p}, \dots, x_{j-1}) := M_{[p,j],t}\bigl(\bm z_t, S^{x_p}, \dots, S^{x_{j-1}} \bigr).
\end{equation}
In the remaining of Section~\ref{sec:zigzagsec} and in Sections~\ref{sec:masters_sec}--\ref{sec:reg_sec}
we will work with these special observables.

Before starting the proof, we introduce further notations.
For any vector $\bm x \in \indset{N}^{k}$, let  $\bm x' \in \indset{N}^{k-1}$ denote its projector onto the first $k-1$ coordinates, that is,
\begin{equation} \label{eq:xprime}
	\bm x' := (x_1, \dots, x_{k-1}).
\end{equation} 
Similarly, we define $\bm x''$, $\bm x'''$, and in general $\bm x^{(p)}$ to be the projection of $\bm x$ onto the first $k-2$, $k-3$, and $k-\ell$, respectively, coordinates, that is 
\begin{equation} \label{eq:xprimes}
	\bm x'' := (x_1, \dots, x_{k-2}), \quad \bm x''' := (x_1, \dots, x_{k-3}),\quad \bm x^{(p)} := (x_1, \dots, x_{k-p}), \quad \ell \in \indset{k-1}.
\end{equation}
Moreover, for a pair of index vectors $\bm x \in \indset{N}^{k}$ and $\bm a \in \indset{N}^{j}$, we denote
\begin{equation}
	(\bm x,\bm a) := (x_1,\dots, x_{k}, a_1,\dots, a_j) \in \indset{N}^{k+j}.
\end{equation}

Let $f(\bm x)$ be a $\mathbb{C}$-valued function of a variable $\bm x \in \indset{N}^k$, and let 
$(\mathcal{A}_j)_{j=1}^{k}$ be a family of linear operators on $\mathbb{C}^N$ (i.e. $N\times N$ matrices), 
then we denote the action of $\bigotimes_{j=1}^k\mathcal{A}_j$ and $\bigoplus_{j=1}^k\mathcal{A}_j$ on $f$ by
\begin{equation} \label{eq:tensor_action}
	\begin{split}
		\biggl(\bigotimes_{j=1}^k \mathcal{A}_j\biggr)\bigl[f\bigr](\bm x) &:= \sum_{\bm a \in \indset{N}^k} \prod_{i=1}^k (\mathcal{A}_i)_{x_i a_i}f(\bm a),\\
		\biggl(\bigoplus_{j=1}^k \mathcal{A}_j\biggr)\bigl[f\bigr](\bm x) &:= \sum_{j=1}^k \sum_{\bm a \in \indset{N}^k}(\mathcal{A}_j)_{x_j a_j} \prod_{i\neq j} \delta_{x_i a_i} f(\bm a).
	\end{split}
\end{equation}
Moreover, for a non-empty subset $J \subset \indset{k}$ we define the partial trace operators $\mathcal{T}_J$, acting on $\mathbb{C}$-valued function $f(\bm x)$ of a variable $\bm x \in \indset{N}^k$, as
\begin{equation} \label{eq:partial_trace}
	\mathcal{T}_J\bigl[f\bigr]\bigl((x_i)_{i\notin J}\bigr) := \sum_{\bm a \in \indset{N}^k} \prod_{i\notin J} \delta_{x_i a_i} f(\bm a).
\end{equation}

\subsubsection{Time-dependent control functions, size-functions and $M$-bounds}
Now we define the time dependent control functions that are straightforward reformulations
of Definition~\ref{def:adm_ups_notime}. 

\begin{Def} [Time-Dependent Admissible Control Functions] \label{def:adm_ups}
Fix a terminal time $T > 0$. For all $t\in[0, T]$, and with $z_t $ being the solution of~\eqref{eq:char_flow}
we set $\eta_t=|\im z_t|$ and $\ell_t=\ell(\eta_t)$ to be time-dependent scale and size parameters.
They satisfy $W\lesssim \ell_t \le N$, $N^{-1} \lesssim \eta_t \lesssim 1$, and  the monotonicity properties
\begin{equation} \label{eq:elleta_assumes}
	\ell_s \le \ell_t, \quad \eta_s \ge \eta_t, \quad \ell_s\eta_s \ge \ell_t\eta_t \ge 1, \quad 0 \le s \le t \le T.
\end{equation}
Let $\Upsilon_\eta$ be an admissible control function by Definition~\ref{def:adm_ups_notime}. 
$\Upsilon_t:=\Upsilon_{\eta_t}$, where $\eta_t := |\im z_t|$, is called a (time dependent) admissible control function.
\end{Def}

We recall the  two typical examples for $\Upsilon_t=\Upsilon_{\eta_t}$ from~\eqref{eq:polyUps_notime}
and~\eqref{eq:expUps_notime}. 
They indicate that $(\Upsilon_t)_{xy}$ is essentially supported near the diagonal $|x-y|_N \lesssim \ell_t$ with a maximal size of $(\ell_t\eta_t)^{-1}$. As time increases, the support of $\Upsilon_t$ grows, while its maximal value decreases, see \eqref{eq:elleta_assumes}.

We use the admissible control function 
as a building block to construct the \emph{size functions} $\mathfrak{s}$ to control the resolvent chains~\eqref{eq:Gk_def}
and their deterministic approximation~\eqref{eq:Mk_def},
especially to encode their
spatial behavior as a function of $\bm x$ and $a, b$ (in case of the isotropic chain).

Fix an even integer $\maxK\ge 8$ to be the maximal chain length we consider. Let $\widehat{z}_1, \dots, \widehat{z}_\maxK \in \fuldom = \fuldom_{\kappa, \etaexp, C_0}$ be a fixed set of spectral parameters, 
satisfying \eqref{eq:admissible_z}, in particular $|\im \widehat{z}_j| \ge N^{-1+\etaexp}$.
Consider the trajectories $\widehat{z}_{j,t}$ of the characteristic \eqref{eq:char_flow} for $0 \le t \le T$, satisfying the final condition $\widehat{z}_{j,T} = \widehat{z}_j$ for all $j \in \indset{\maxK}$. 
We define a discrete time-dependent set of spectral parameters $\dom_t$ as
\begin{equation} \label{eq:domt}
	\dom_t :=  \bigl\{ \widehat{z}_{j,t}, (\widehat{z}_{j, t})^*\,: \, j \in \indset{\maxK}\bigr\} \subset 
	(\fuldom_{\kappa, \etaexp, C_0'})^{2K},  
\end{equation}
where $(\widehat{z}_{j, t})^*$ denotes the complex conjugate or $\widehat{z}_{j, t}$.   The inclusion
follows from~\eqref{eq:edge}. 
Recall from \eqref{eq:comparability} that for each fixed $t$, 
the  absolute values of the imaginary parts of $\widehat{z}_{1,t}, \widehat{z}_{2,t}, \ldots, \widehat{z}_{K,t}$ are comparable.

For $k \in \mathbb{N}$ and $\bm z_t := (z_{1,t},\dots, z_{k,t})\in \dom_t^k$, we set
$$
  \eta_t: = \min_j |\im z_{j, t}|, \quad \ell_t =\ell(\eta_t),
$$
i.e. we choose the smallest $\eta$ to represent all spectral parameters, therefore $\ell_t$ represents the
largest length scale. 
We define the \emph{isotropic}\footnote{Here and in the sequel we use the adjective \emph{isotropic} for quantities  that depend on two matrix entries $a, b$, hence in the standard RMT 
literature they are called \emph{entry-wise}. Our slight change  in terminology is justified
since in a separate argument  we will eventually prove general isotropic bounds using these entry-wise quantities.}
 and \emph{averaged size functions} $\mathfrak{s}_k^\mathrm{iso/av}$, depending on
 time $t$ via $\eta_t$,
\begin{equation} \label{eq:sfunc_def}
	\begin{split}
		\mathfrak{s}_{1,t}^\mathrm{iso}(a,b) &:= \sqrt{(\Upsilon_t)_{ab}}, \quad\mathfrak{s}_{1,t}^\mathrm{av}(x) := \frac{1}{\ell_t\eta_t}, \\
		\mathfrak{s}_{k,t}^\mathrm{iso}(a,\bm x',b) &:= \frac{1}{(\ell_t\eta_t)^{(k-1)/2}}\sqrt{(\Upsilon_t)_{ax_1}(\Upsilon_t)_{x_{k-1}b}\prod_{j=2}^{k-1}(\Upsilon_t)_{x_{j-1}x_{j}}}, \quad k \ge 2,\\
		\mathfrak{s}_{k,t}^\mathrm{av}(\bm x) &:= \frac{1}{(\ell_t\eta_t)^{k/2}}\sqrt{(\Upsilon_t)_{x_kx_1}\prod_{j=2}^{k}(\Upsilon_t)_{x_{j-1}x_{j}}} = \frac{\mathfrak{s}_{k,t}^\mathrm{iso}(x_k,\bm x',x_k)}{\sqrt{\ell_t\eta_t}}, \quad k \ge 2,
	\end{split}
\end{equation} 
for all $k \in \mathbb{N}$, $\bm x \in \indset{N}^k$ and $a,b \in \indset{N}$. 

We  now state the bounds on the deterministic approximations $M$ in terms of the size functions
 that will be proven in Section \ref{sec:M_bounds}: 
\begin{lemma}[$M$-Bounds] \label{lemma:M_bounds}
	Fix an integer $k \in \mathbb{N}$. Then, for any $\bm x\in\indset{N}^k$ and $\bm z_t \in \dom_t^k$, 
	  the deterministic approximation $M_{[1,k],t}(\bm x')$, defined in \eqref{eq:Mk_def} satisfy
	\begin{equation} \label{eq:M_bound_av}
		\bigl\lvert  \Tr \big[ M_{[1,k],t}(\bm z_t, \bm{x}') S^{x_k} \big] \bigr\rvert \lesssim  (\log N)^{C_k} \times(\ell_t\eta_t)\,\mathfrak{s}_{k,t}^\mathrm{av}(\bm x), 
	\end{equation}
	\begin{equation} \label{eq:M_bound}
		\bigl\lvert \bigl(M_{[1,k],t}(\bm z_t, \bm{x}') \bigr)_{ab}\bigr\rvert \lesssim \delta_{ab} \,(\log N)^{C_k+1} \times \sqrt{\ell_t\eta_t}\,\mathfrak{s}_{k,t}^\mathrm{iso}(a, \bm x', b),
	\end{equation}
	where $\eta_t= \min_j \eta_{j, t}$, $\ell_t=\ell(\eta_t)$.
\end{lemma}
The additional $\ell_t\eta_t$ factors are present because we designed the size functions to control directly the fluctuations
in the local laws and not the $M$-terms themselves, see the main statements of Theorem~\ref{th:local_laws}.

Since we use the $M$-bound to establish estimates on the $(G-M)$ quantities up to an $N^{\xi}$ tolerance, the $(\log N)$-factors  in \eqref{eq:M_bound_av}--\eqref{eq:M_bound} can safely be ignored. In fact, there is an alternative static approach,
 which relies purely on the recursion \eqref{eq:M_recursion}, that allows to establish the $M$ bounds without the poly-logarithmic factors. However, the static prove requires cumbersome graphical expansions, so we only present the simpler dynamical proof in Section \ref{sec:M_bounds}.

We remark that we immediately defined the size functions and the $M$-bounds
in their time-dependent versions since we use them in this form,
but they depend on time only via $\eta_t$, so these are really 
time-independent bounds for any fixed $\bm z\in (\fuldom)^k$ under the condition \eqref{eq:admissible_z}.

\subsubsection{Formulation of the local laws via $\Psi$-functions}

Resolvent chains and their deterministic $M$ approximations have natural sizes and spatial structure, expressed
via the size functions as Lemma~\ref{lemma:M_bounds} demonstrates.
 We define the random $\Psi^\mathrm{iso}$ and $\Psi^\mathrm{av}$-functions that are essentially the ratios 
of the fluctuations of the resolvent chains and their size functions, i.e. these are dimensionless versions of 
the isotropic and averaged local laws. 
Moreover, $\Psi$'s have only two parameters; the length of the chain $k$ and the time $t$, thus 
they encode uniformity in all parameters  within a single scalar quantity.
Proving local laws amounts to proving that $\Psi$'s are essentially of order 1 with high probability, 
 more precisely  that $\Psi\prec 1$.   Introducing the $\Psi$'s substantially shortens the formulas
 and clarifies the essence of the proofs.
 
Their precise definition is as follows. 
Recall that $\maxK \ge 8$ is a fixed even integer, corresponding to the maximal length of the resolvent chain we consider.

For $k \in \indset{K}$, we define the {\it isotropic} and {\it averaged loss exponents} $\alpha_k$ and $\beta_k$, respectively, as
\begin{equation} \label{eq:loss_exponents}
	\alpha_k := \begin{cases}
		0, &\quad k \le \maxK/2,\\
		\frac{1}{2}\sqrt{2k/K-1}, &\quad k \ge \maxK/2+1,  
	\end{cases}
	\quad 
	\beta_k := \begin{cases} 
		\alpha_{k+1}, &\quad 
		k \le \maxK-1,\\
		1/2 + \alpha_{\maxK/2+1}, &\quad k=\maxK.
	\end{cases}
\end{equation}
For $k\in\indset{\maxK}$, we define\footnote{The somewhat complicated definition via 
the set $\dom_t$ is necessary to keep track of all possible combinations of spectral parameters in various chains 
 that will arise in  our analysis. The reader is encouraged to think of the simplest case where there is one fixed
 trajectory $z_t$, $\dom_t= \{ z_t, \bar z_t\}$,
  and every coordinate of $\bm z_{t}$ is either $z_t$ or $\bar z_t$. We stress that even in this simple case
 the additional freedom to consider the complex conjugate trajectory $\bar z$ together with $z$ is necessary.}   the
 random control quantities $\Psi_{k,t}^\mathrm{iso}$ and $\Psi_{k,t}^\mathrm{av}$ as
 \begin{equation} \label{eq:Psi_def}
	\begin{split}
		\Psi_{k,t}^\mathrm{iso} :=
		 \max_{ \bm z_{t} \in \dom_t^k} 
		\max_{a,b\in\indset{N}} \max_{\bm x'\in \indset{N}^{k-1}}  \frac{\bigl\lvert \bigl(G_{[1,k],t} - M_{[1,k],t}\bigr)_{ab} (\bm x') \bigr\rvert}{(\ell_t\eta_t)^{\alpha_k}\mathfrak{s}_{k,t}^\mathrm{iso}(a,\bm x',b)},\\
		\Psi_{k,t}^\mathrm{av} := 
		 \max_{  \bm z_{t} \in \dom_t^k} 
		\max_{\bm x\in \indset{N}^k} \frac{\bigl\lvert \Tr \bigl[ \bigl(G_{[1,k],t} - M_{[1,k],t}\bigr) (\bm x')  S^{x_k}\bigr]\bigr\rvert}{(\ell_t\eta_t)^{\beta_k} \mathfrak{s}_{k,t}^\mathrm{av}(\bm x)},
	\end{split}
\end{equation}
where $\dom_t$ is defined in \eqref{eq:domt}. 
 For $k=1$, we consider the argument $\bm x'$ to be dropped by convention. Therefore, our goal is to prove that $\Psi_{k,t}^{\mathrm{iso/av}} \prec 1$.

Note that the loss exponents satisfy $0\le \alpha_k\le 1/2$, $0\le \beta_k<1$, and both $\alpha_k, \beta_k$ are  non-decreasing in $k$. Moreover, $k\mapsto \alpha_k$ is concave where it is positive. 
These exponents capture the gradual deterioration in the error term of the local law  as the chain length $k$ increases. That is, optimal error estimates would correspond to $\alpha_k = \beta_k = 0$ for all $k$; however, in our setup, 
we only have $\alpha_k=0$ for $k\le K/2$ and $\beta_k=0$ for $k\le K/2-1$. 
This deterioration for larger values of $k$ reflects the cost of using reduction inequalities, which are necessary for truncating the hierarchy at length $\maxK$. While we will only use reduction inequality for chains beyond the length $K$, the resulting loss is distributed evenly over the range $K/2 \le k \le K$; the loss exponents govern this distribution.
The square root function used in~\eqref{eq:loss_exponents} is chosen for convenience, though other concave choices would also be possible.
We emphasize that this deterioration is merely a technical artifact of the proof. Since the maximal length $K$ can ultimately be chosen arbitrary
large, and $\alpha_k=0$ for $k\le K/2$ and $\beta_k=0$ for $k\le K/2-1$,
the final result still yield optimal multi-resolvent local laws for chains of any length.

\medskip

Armed with this notation, we can reformulate the local laws in terms of $\Psi$'s. 
Fix a maximal chain length $\maxK \ge 8$. 
Recall the definition of $\dom_t$ from \eqref{eq:domt}, and that for all $\bm z_t \in \dom_t^k$,
the comparability of the imaginary parts
remains valid along the whole trajectory. More precisely, by~\eqref{eq:comparability}, 
\begin{equation}
	|\im z_{i,t}|\sim  |\im z_{j,t}|, \quad i,j\in\indset{k}, \quad 0 \le t \le T.
\end{equation}
Moreover, the entire trajectory stays far away from the spectral edges~\eqref{eq:edge}. 
For the initial condition we have $|\im z_{j,t=0}|\sim 1$ for all $j$.
The maximal chain length $\maxK$ and the entire set of trajectories $\dom_t$ are considered fixed for the
rest of Section~\ref{sec:zigzagsec}.

\begin{Def}[Local Laws in terms of $\Psi$] \label{def:local_laws} 
	Let $\xi > 0$ be a tolerance exponent, and let $\mathcal{I} \subset [0,T]$ be a closed non-empty (but possibly degenerate) time interval . 
	Let $H_t$ be a real symmetric or complex Hermitian random matrix, depending on the parameter $t \in \mathcal{I}$, with $\Expv H_t = 0$ and admissible variance profile $S$. 
	We say that the (multi-resolvent) local laws hold for the resolvent $G_{i,t}:=(H_t-z_{i,t})^{-1}$ with tolerance $\xi$ 
	uniformly in $t\in \mathcal{I}$ if and only if the bounds
	\begin{equation}\label{eq:LLwithPsi}
		\max_{t\in\mathcal{I}}\Psi_{k,t}^{\mathrm{av/iso}} \le N^{\xi}, \quad k \in \indset{\maxK},
	\end{equation}
	hold with very high probability, where $\Psi_{k,t}^{\mathrm{av/iso}}$ are the quantities given in \eqref{eq:Psi_def} with $G_{[p,j],t}$ defined according to \eqref{eq:resolvent_chains} and \eqref{eq:Gk_def}.
	 If $|\im z_{i,t}|\gtrsim 1$, then we talk about global laws instead of local laws.
\end{Def}

\subsection{Multi-resolvent local laws via zigzag}\label{sec:zz}

Before we formulate the three key ingredients of the zigzag strategy as separate propositions, we
point out a technical point that will play an important role later when we extend local laws to include 
more general observables and to their full isotropic version with arbitrary test vectors instead of 
coordinate vectors only. The full set of conditions on the admissible control functions
listed  in Definition~\ref{def:adm_ups_notime}
are necessary only for the zig-step. The global laws, the zag-step and even some technical estimates
within the zig-step do not require
\eqref{eq:convol_notime}--\eqref{eq:suppressed_convol_notime}, instead they can be proven under
the following weaker condition:
for all $z_1, z_2\in \fuldom$ with $\eta_1\le \eta_2$ we have
\begin{equation} \label{eq:true_convol_notime} 
	\sum_a (\Upsilon_{\eta_2})_{xa} (\Upsilon_{\eta_1})_{ay} \le C_1\frac{1}{\eta_2} (\Upsilon_{\eta_1})_{xy}, \quad x,y\in\indset{N}.
\end{equation}
In the sequel we will always indicate if a statement holds under the weaker condition \eqref{eq:true_convol_notime} instead of \eqref{eq:convol_notime}--\eqref{eq:suppressed_convol_notime}.  This weaker condition
directly refers to the admissible control functions $\Upsilon_\eta$ but, consequently, 
it also weakens the admissibility condition of the variance profile. Whenever we refer to
the weaker condition \eqref{eq:true_convol_notime}, we also mean that $S$ is admissible in
the sense of Definition~\ref{def:admS} with the weaker admissibility conditions on $\Upsilon_\eta$.

In the following three propositions the maximal chain length $\maxK$ and the entire trajectory $\bm z_t$, $t\in [0,T]$
with components of comparable imaginary parts  are fixed.

\begin{prop}[Global Laws] \label{prop:global_laws} Let $H$ be a real symmetric or complex Hermitian random matrix with $\Expv H = 0$ and admissible 
	variances profile $S$ but we assume only the weaker
condition \eqref{eq:true_convol_notime} instead of \eqref{eq:convol_notime}--\eqref{eq:suppressed_convol_notime}.  
	 Let $\bm z_0 =(z_{1,0}, z_{2,0}, \ldots, z_{K,0})\subset(\fuldom)^K$  
	  be a 
	 $\maxK$-tuple of spectral parameters satisfying the comparability condition \eqref{eq:admissible_z}, 
	  $\eta_{j,0} \ge c$ and $|z_{j,0}|\le C$, $j\in\indset{K}$,  then the 
	  multi-resolvent global laws hold for the $K$-tuple of resolvents $G_{j,0} := (H_0 - z_{j,0})^{-1}$, $j\in\indset{\maxK}$, 
	  i.e. for the one element time interval  $\mathcal{I}:=\{0\}$. 
\end{prop}
We prove Proposition~\ref{prop:global_laws} in Section~\ref{sec:global3}, after gradually developing the necessary ingredients in Section~\ref{sec:global}.

\begin{prop}[Zig-Step] \label{prop:zig}
	Let $H$ be a real symmetric or complex Hermitian random matrix with $\Expv H = 0$ and admissible 
	variances profile $S$.  Fix the set of trajectories $\dom_t$, defined in \eqref{eq:domt}.  
	 Fix a time  $\tinit\in[0,T]$, and let $H_t$ be the solution to the  stochastic differential equation \eqref{eq:zigOU}. 
	Assume that for some tolerance exponent $\xi \in (0, \etaexp/100)$, the multi-resolvent local laws hold for the resolvent $G_{j,\tinit} := (H_{\tinit} - z_{j,\tinit})^{-1}$ with tolerance $\xi$ at time $t = t_{\tinit}$ in the sense of Definition \ref{def:local_laws}. Then, for any fixed positive $\nu > 0$, the local laws also hold for $G_{j,t} := (H_t - z_{j,t})^{-1}$ with tolerance $\xi+\nu$ uniformly in  $t \in [\tinit, T]$.
\end{prop}

We prove Proposition \ref{prop:zig} in Section \ref{sec:masters_sec}.
To remove the Gaussian component introduced during the zig-step in Proposition \ref{prop:zig} above, we use the following Green Function Comparison theorem, called as zag-step in the zigzag context.
\begin{prop}[Zag-Step] \label{prop:zag}
	Assume only the weaker  condition \eqref{eq:true_convol_notime} instead of \eqref{eq:convol_notime}--\eqref{eq:suppressed_convol_notime}.  
	Let $H := \bigl(\sqrt{S_{ij}}h_{ij}\bigr)_{i,j=1}^N$ and $V := \bigl(\sqrt{S_{ij}}v_{ij}\bigr)_{i,j=1}^N$ be two $N\times N$ real symmetric or complex Hermitian random matrices with $\Expv H = \Expv V = 0$ and admissible variance profile $S$.	
	Fix the set of trajectories $\dom_t$, defined in \eqref{eq:domt}.  
	Fix a time $t \in [0, T]$, and assume that, for some tolerance exponent $\xi \in (0, \etaexp/100)$, the local laws hold for $\other{G}_{j,t} := (V-z_{j,t})^{-1}$ with tolerance $\xi$ at time $t$ in the sense of Definition \ref{def:local_laws}. Then, the local laws also hold for $G_{j,t} := (H-z_{j,t})^{-1}$ with tolerance $\xi+\nu$  at time $t$ for any fixed $\nu >0$, provided the following moment-matching conditions are satisfied:
	\begin{itemize}
		\item [(a)]  The first three are moments matched exactly,
		\begin{equation} \label{eq:moment_match}
			\Expv\bigl[ h_{ij}^s h_{ji}^{r-s} \bigr] = \Expv\bigl[ v_{ij}^s v_{ji}^{r-s} \bigr], \quad i,j\in\indset{N}, \quad s \in \{0,\dots, r\}, \quad r \in \{0,1,2,3\}.
		\end{equation}
		
		\item [(b)] The fourth and fifth moments are matched approximately, 
		\begin{equation} \label{eq:m4_cond}
			\Expv\bigl[ h_{ij}^s h_{ji}^{r-s} \bigr] = \Expv\bigl[ v_{ij}^s v_{ji}^{r-s} \bigr] + \mathcal{O}(\lambda), \quad i,j\in\indset{N}, \quad s \in \{0,\dots, r\}, \quad r \in \{4,5\},
		\end{equation}
		where the implicit constant in $\mathcal{O}(\lambda)$ is uniform in $i, j$, and 
		the small parameter $\lambda \equiv \lambda(\xi,t)$ is given by
		\begin{equation} \label{eq:lambda_assume}
			\lambda(\xi,t) := N^{\xi/2} \ell_t^{-1} W \sim N^{\xi/2} \max\bigl\{\sqrt{\eta_t}, W/N \bigr\}.
		\end{equation} 
	\end{itemize}
\end{prop}
We prove Proposition \ref{prop:zag} in Section \ref{sec:zag_proof_for_real}, following the preparatory results established in Section~\ref{sec:GFT}.
Equipped with Propositions~\ref{prop:zig} and~\ref{prop:zag}, which encapsulate the cardinal steps of the zigzag strategy, we now combine them to prove the multi-resolvent local laws of Theorem \ref{th:local_laws}.

Note that we formulated all three propositions in a time dependent form primarily to facilitate the explanation of the zigzag proof. However, intrinsically time plays a role only in the zig-step, while
the global law and the zag-step are, in essence, time-independent statements. 
Specifically, the global law holds for any $\bm z\in (\fuldom)^k$ satisfying the comparability condition and with $\eta\ge c$, while the zag-step applies to any $\bm z\in (\fuldom)^k$ under the comparability condition \eqref{eq:admissible_z}.

\subsection{Proof of Theorem~\ref{th:local_laws} with special observables $S^x$}

\begin{proof} [Proof of Theorem \ref{th:local_laws}]  
	Fix a target tolerance exponent $\xi_0 \in (0,\etaexp/100)$. Let $\bm z$ be the fixed target spectral parameters.
	 Recall that the final time $T\sim 1$ is fixed, and $\bm z_t$ is the trajectory of \eqref{eq:char_flow} satisfying
	  $\bm z_T = \bm z$
	and $\eta:=\min |\im z_j| = \min |\im (z_T)_j|$. We partition the time interval $t\in [0,T]$ into finitely many sub-intervals $\{[t_{q-1}, t_{q}]\}_{q=1}^{\maxQ}$, with $\maxQ$ satisfying the bound $\maxQ \le 1/\xi_0$.    We will apply
	Proposition~\ref{prop:zig} on each sub-interval; at step~$q$, we set the initial $\tinit := t_{q-1}$ and deduce the
	local law at time $t_q$.  We choose the lengths $t_q' := t_{q}-t_{q-1}$ of these intervals in such a way that the Gaussian component of size $\mathcal{O}(\sqrt{t_q'})$ introduced by the process~\eqref{eq:zigOU}
	in the zig-step  can be removed using Proposition~\ref{prop:zag} in the zag-step. 
	
	Fix $\xi := \xi_0/2$, and let $\maxQ \in \mathbb{N}$ be the smallest integer such that 
	\begin{equation} \label{eq:Q_cond}
		N^{-\maxQ\xi} \lesssim N^{-\xi/2}\max\bigl\{\sqrt{\eta}, W/N\}.
	\end{equation}
	It follows from \eqref{eq:etabound} and \eqref{eq:WN} that $\max\bigl\{\sqrt{\eta}, W/N\} \gtrsim N^{-1/2+\bandexp/2}$, hence $\maxQ \le 1/\xi_0$. 
	Consider the sequence of times $\{t_q\}_{q=0}^{\maxQ}$, defined as
	\begin{equation} \label{eq:t_seq}
		t_0 := 0, \quad t_q := T - N^{-q\xi}T, \quad q \in \indset{\maxQ-1}, \quad t_{\maxQ} := T.
	\end{equation} 
	It is straightforward to check, using \eqref{eq:ell_def_notime}, \eqref{eq:etatasymp} and $\eta=\eta_T$,
	that $\eta_{t_q} \sim  \eta + N^{-q\xi}$ and 
	\begin{equation} \label{eq:tq'_bound}
		 t_{q}' = t_{q}-t_{q-1} \sim N^{-(q-1)\xi} \lesssim \lambda(\xi, t_q), \quad q\in\indset{\maxQ},
	\end{equation} 
	where $\lambda(\xi, t_q)$ is defined in \eqref{eq:lambda_assume}.  
	We now construct the sequence of random matrices~$\{H^{q}\}_{q=0}^{\maxQ}$ using the following backward induction:
	\begin{itemize}
		\item Set $H^{\maxQ} := H$ to be the target random matrix;
		
		\item For each $q \in \indset{\maxQ}$, the bound \eqref{eq:tq'_bound} ensures (via a standard moment-matching argument, see Lemma \ref{lemma:moment_match}  below for more details) the existence of a random matrix $H^{q-1}$ in the same symmetry class as $H^q$, such that the Gaussian-divisible matrix $V^q$, defined by 
		\begin{equation} \label{eq:Vq}
			V^q := \sqrt{\mathrm{e}^{-t_q'}} H^{q-1} + \sqrt{1-\mathrm{e}^{-t_q'}} \mathcal{H}_{\mathrm{G}}, 
		\end{equation}
		satisfies the  conditions \eqref{eq:moment_match} and \eqref{eq:m4_cond}
		on the normalized matrix elements of $h^q_{ij}=H^q_{ij}/\sqrt{S_{ij}}$ and $v^q_{ij}=V^q_{ij}/\sqrt{S_{ij}}$.  Here $\mathcal{H}_{G}$ denotes a Gaussian random band matrix with zero mean and variance profile $S$, in the same symmetry class as $H^q$, and independent of~$H^{q-1}$.
	\end{itemize}
	We proceed by forward induction in $q \in \indset{0,\maxQ}$. For the base case $q=0$,
	the global law, Proposition \ref{prop:global_laws}, implies that the local laws in the sense of Definition \ref{def:local_laws} hold for $G_{j,0} := \bigl(H^0 - z_{j,0}\bigr)^{-1}$ with tolerance $\xi$ at time $t_0  = 0$.  
	
	For the induction step $q\in\indset{\maxQ}$, fix a positive tolerance exponent $0 < \nu  < \maxQ^{-1}\xi/10$. Assume that the local laws  hold for $G_{j,q-1} := \bigl(H^{q-1} - z_{j,t_{q-1}}\bigr)^{-1}$ with tolerance~$\xi + (2q-2)\nu$ at time~$t_{q-1}$. Consider the matrix process \eqref{eq:zigOU} $H_t$ with $\tinit := t_{q-1}$ and initial data $H_{\tinit} := H^{q-1}$, then
	 $H_{t_q}$ is equal in distribution   to $V^q$ , defined in \eqref{eq:Vq}. Therefore, Proposition~\ref{prop:zig} implies that the local laws hold for $\other{G}^q_j:= \bigl(V^q - z_{j,t_q}\bigr)^{-1}$ with tolerance~$\xi + (2q-1)\nu$ at time $t_q$. 
	Since the normalized matrix elements of $V^q$ and $H^q$ satisfy the moment-matching conditions~\eqref{eq:moment_match} and \eqref{eq:m4_cond}, by Proposition~\ref{prop:zag}, the local laws also hold for $G^{q}_j := \bigl(H^{q} - z_{j,t_{q}}\bigr)^{-1}$ with tolerance~$\xi + 2q\nu$ at time~$t_{q}$.  
	
	\begin{figure}[H]
		\centering
		\begin{minipage}{0.6\textwidth}
			\centering
			\includegraphics[width=\textwidth]{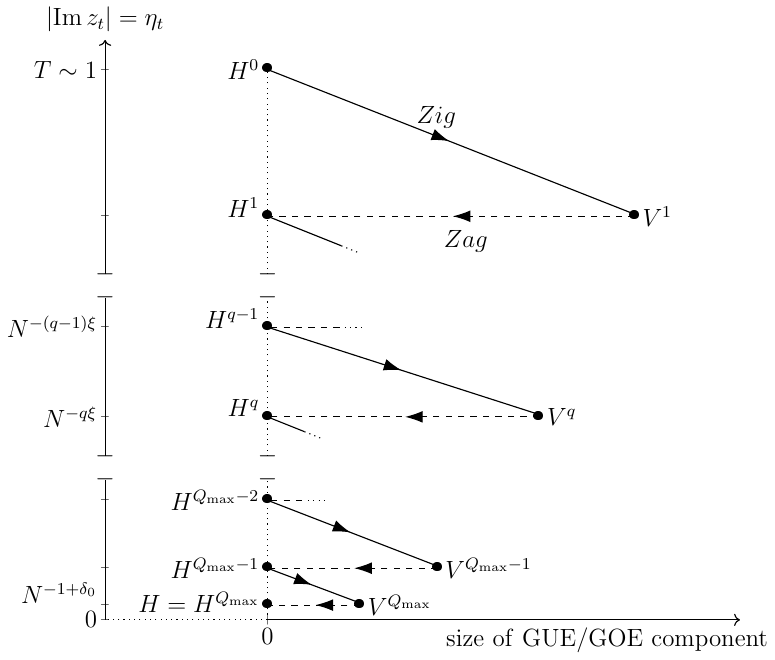}
		\end{minipage}%
		\hfill
		\begin{minipage}{0.4\textwidth}
			\centering
			\captionsetup{width=.85\textwidth, justification=justified}
			\caption[Zigzag caption]{
				Schematic representation of the Zigzag induction.
				The random matrices $H^q, V^q$, as defined in \eqref{eq:Vq}, positioned within an abstract coordinate system. The horizontal axis corresponds to the size of the Gaussian component, while the vertical axis represents the lower bound on $|\im z_t| = \eta_t$ for which the local laws are established.
				Solid arrows indicate applications of Proposition~\ref{prop:zig} (referred to as \textit{Zig} steps), while dashed arrows denote applications of Proposition~\ref{prop:zag} (\textit{Zag} steps).
			}
			\label{fig:zigzag_fig}
		\end{minipage}
	\end{figure}
	
	After $\maxQ$ steps in this induction, 
	 the local laws holds for $G_{j,T} := (H_T-z_{j,T})^{-1} = (H-z_j)^{-1}$ with tolerance $\xi + 2\maxQ\nu \le \xi_0$ at time $T$. This concludes the proof of Theorem~\ref{th:local_laws}.
\end{proof}

In the proof above we used the following elementary moment matching lemma. We omit the proof as it is very 
similar to that of the  Lemma 16.2 of \cite{erdHos2017dynamical}.
\begin{lemma}[Moment Matching] \label{lemma:moment_match}
	Let $h$ be a real- or complex-valued random variable with zero mean and unit variance, and let $\lambda \in (0,1)$ be a fix constant. 
	Then, there exists a real- or complex-valued, respectively, random variable $h'$, 
	 such that the Gaussian-divisible random variable  $v$, defined as 
	\begin{equation}
		v = \sqrt{\mathrm{e}^{-\lambda}} h' + \sqrt{1-\mathrm{e}^{-\lambda}} h_{\mathrm{G}},
	\end{equation}
	with $h_G$ being a standard Gaussian variable independent of $h'$,
	satisfies the moment-matching conditions
	\begin{equation}
		\begin{alignedat}{2}
			\Expv\bigl[ h^s (\overline{h})^{r-s} \bigr] &= \Expv\bigl[ v^s (\overline{v})^{r-s} \bigr], \quad &&s \in \{0,\dots, r\}, \quad r \in \{0,1,2,3\},\\
			\Expv\bigl[ h^s (\overline{h})^{r-s} \bigr] &= \Expv\bigl[ v^s (\overline{v})^{r-s} \bigr] + \mathcal{O}(\lambda), \quad &&s \in \{0,\dots, r\}, \quad r \in \{4,5\}.
		\end{alignedat}
	\end{equation}
\end{lemma}

\section{Master inequalities: Proof of Proposition \ref{prop:zig}} \label{sec:masters_sec}
The main tool is the following system of dynamical \emph{master inequalities}, which represent a self-improving mechanism for controlling the size of our key quantities $\Psi_{k,t}^\mathrm{iso/av}$.
\begin{prop}[Dynamical Master Inequalities] \label{prop:masters}
	Let $\{ \psi_{k,0}^\mathrm{iso}, \psi_k^\mathrm{iso}, \psi_{k,0}^\mathrm{av}, \psi_k^\mathrm{av}\}_{k=1}^\maxK$ be a set of time-independent deterministic control parameters satisfying 
	$\psi_k^\mathrm{iso/av} > \psi_{k,0}^\mathrm{iso/av} \ge 1$.
	Define the stopping time $\tau \in [\tinit,T]$ as
	\begin{equation}\label{eq:tau_def}
		\tau :=  \inf\biggl\{ t\in[\tinit,T] \, :\, 
		\max_{k\in \indset{\maxK}} \Psi_{k,t}^\mathrm{av} / \pav{k} \ge 1 \biggr\} \wedge \inf\biggl\{ t\in[\tinit,T] \, :\, \max_{k \in \indset{\maxK}}  \Psi_{k,t}^\mathrm{iso} / \pis{k} \ge 1\biggr\} ,
	\end{equation}
%
%
	where the random quantities $\Psi_{k,s}^{\mathrm{av/iso}}$ are given by \eqref{eq:Psi_def}.  Suppose that at time $t=\tinit$, 
	\begin{equation} \label{eq:zig_init}
		\Psi_{k,\tinit}^\mathrm{av/iso} \prec \psi_{k,0}^\mathrm{av/iso}, \quad k\in\indset{\maxK}.
	\end{equation}   Then $\Psi_{k,t}^{\mathrm{av/iso}}$ satisfy the bounds
	\begin{subequations} \label{eq:masters}
		\begin{equation}
			\max_{\tinit\le t\le \tau}\Psi_{k,t}^\mathrm{iso} \prec \varphi_{k,\tau}^\mathrm{iso},   \quad k\in\indset{\maxK}, \label{eq:iso_masters}
		\end{equation}
		\begin{equation}
			\max_{\tinit\le t\le \tau}\Psi_{k,t}^\mathrm{av} \prec \varphi_{k,\tau}^\mathrm{av},   \quad k\in\indset{\maxK}, \label{eq:av_masters}
		\end{equation}
	\end{subequations} 
	where, for $t \in [0,T]$, the time-dependent  control parameters $\varphi_{k,t}^\mathrm{av/iso}$ are defined as 	
	\begin{subequations} \label{eq:phi_masters_def} 
		\begin{equation}
			\varphi_{k,t}^\mathrm{iso} := 
			\psi_{k,0}^\mathrm{iso}  + \sum_{j=2}^{k-1} \frac{\pis{j}}{(\ell_{t}\eta_{t})^{\alpha_k-\alpha_j}} +\sum_{j=2}^k \pis{j} \biggl( \frac{\sqrt{N\eta_t}}{\ell_t\eta_t} \pis{k-i+1} + \frac{\pis{k-j+2}}{\sqrt{\ell_t\eta_t}}\biggr)  + \frac{\pav{k}}{(\ell_{t}\eta_{t})^{1/2+\alpha_k-\beta_k 
			}} + \varphi_{k,t}^\mathrm{iso, qv}, \label{eq:phi_iso_masters}
		\end{equation}
		\begin{equation}
			\varphi_{k,t}^\mathrm{av} := \psi_{k,0}^\mathrm{av} +  \sum\limits_{j=2}^{k-1} \frac{\pav{j}}{(\ell_{t}\eta_{t})^{\beta_k-\beta_j}} \biggl(1+\frac{\pis{k+2-j}}{(\ell_{t}\eta_{t})^{1/2-\alpha_{k+2-j}}}\biggr) + \frac{\pis{2}\pav{k}}{\sqrt{\ell_{t}\eta_{t}}}  +\varphi_{k,t}^\mathrm{av,qv} + \varphi_{k,t}^\mathrm{av,f}. \label{eq:phi_av_masters}
		\end{equation}
	\end{subequations}
	Here, the auxiliary control parameters $\varphi_{k,t}^\mathrm{iso, qv}$ and $\varphi_{k,t}^\mathrm{av,qv}$ (arising from quadratic variation estimates), and $\varphi_{k,t}^\mathrm{av,f}$ (arising from using a reduction inequality in a forcing term), are given by
	\begin{equation} \label{eq:phi_isoqv}
		\varphi_{k,t}^\mathrm{iso,qv}:= \frac{\psi_k^\mathrm{iso}\sqrt{\psi_2^\mathrm{iso}}}{(\ell_t\eta_t)^{1/4}}  + \sum_{j=2}^{\lfloor (k+1)/2 \rfloor} \frac{\pis{k-j+1}}{(\ell_t\eta_t)^{\alpha_{k}-\alpha_{k-j+1}}}\sqrt{1+ \frac{\pis{2j}}{(\ell_t\eta_t)^{1/2-\alpha_{2j}}}} + \other{\varphi}_{k,t}^\mathrm{\,iso,qv}
	\end{equation}
	\begin{equation} \label{eq:other_phi_isoqv}
		\other{\varphi}_{k,t}^\mathrm{\,iso,qv}:=\begin{cases}
			\sqrt{1+ \frac{\pis{2k}}{(\ell_t\eta_t)^{1/2-\alpha_{2k}}}},\quad& k \in\indset{\maxK/2},\\
			\sqrt{\frac{\pis{k-1}\pis{k+1}}{(\ell_t\eta_t)^{2\alpha_{k}- \alpha_{k-1}-\alpha_{k+1}} }}, \quad& k \in\indset{\maxK/2+1,\maxK-1},\\
			\pis{\maxK/2}\sqrt{1+\frac{\pav{\maxK}}{(\ell_t\eta_t)^{1-\beta_\maxK}}}, \quad& k = K,
		\end{cases}
	\end{equation}
	\begin{equation} \label{eq:phi_avqv}
		\varphi_{k,t}^\mathrm{av,qv} := \begin{cases}
			\sqrt{1 + \frac{\pis{2k+2}}{(\ell_t\eta_t)^{1/2-\alpha_{2k+2}}}} , \quad & k \in \indset{\maxK/2-1},\\
			\pis{k+1}, \quad& k \in \indset{\maxK/2, \maxK-1},\\
			\pis{\maxK/2+1}\sqrt{1 + \frac{\pav{\maxK}}{(\ell_t\eta_t)^{1-\beta_{\maxK} }}}, \quad &k=\maxK,  
	\end{cases}
	\end{equation}
	\begin{equation} \label{eq:phi_av_force}
		\varphi_{k,t}^\mathrm{av,f} :=   \begin{cases}
			 \frac{\sqrt{N_t\eta_t}}{\ell_t\eta_t} \pav{1}\pis{k+1}, \quad& k =1,\\
			 \frac{\sqrt{N_t\eta_t}}{\ell_t\eta_t} \pav{1}\pis{k+1}  + \frac{\pav{1}}{(\ell\eta)^{\beta_{k}}}, \quad& k \in\indset{2,\maxK-1},\\
			 \frac{\sqrt{N_t\eta_t}}{\ell_t\eta_t} \pav{1}\pis{\maxK/2}\pis{\maxK/2+1}, \quad &k=\maxK.
		\end{cases}
	\end{equation}
	 Here the non-negative loss exponents $\{\alpha_k, \beta_k\}_{k=1}^\maxK$ are defined in \eqref{eq:loss_exponents}. 
\end{prop}
After some preparation, we prove Proposition~\ref{prop:masters} in Section~\ref{sec:masters_proof}.

It is easy to check from \eqref{eq:loss_exponents} that all powers of $\ell_\tau \eta_\tau \ge N^{\etaexp}$ in the 
denominators on the right-hand side of \eqref{eq:masters}--\eqref{eq:phi_av_force} are non-negative. 
 Moreover, $\sqrt{N\eta_\tau}/(\ell_\tau \eta_\tau) \lesssim N^{-\etaexp/2}$.   
The structure of the
master inequalities is such that the right hand side of \eqref{eq:masters}  contains control parameters
with index at most $K$.  Moreover, they are {\it self-improving} in the sense that, 
for any fixed $k$, every term on the right hand side of the isotropic equation \eqref{eq:iso_masters}-\eqref{eq:phi_iso_masters} 
 contains either  $\psi$'s with indices strictly less than $k$,
or they are suppressed by a small positive power of $1/(\ell_t\eta_t)$, with the exception of the very first term $\psi_{k,0}^\mathrm{iso}$, which controls the initial condition. 
Similarly, apart from the initial $\psi_{k,0}^\mathrm{av}$,
 every term on the right hand side of the averaged equation \eqref{eq:av_masters}-\eqref{eq:phi_av_masters} 
contains either $\psi^\mathrm{iso}_j$, $j\le K$,  or $\psi_j^\mathrm{ave}$ with strictly smaller index $j <k$, or
comes suppressed by a small negative power of $N$.
This self-improving structure allows us to easily iterate the system of master inequalities, eventually achieving  
$\Psi^\mathrm{av/iso}_{k,t}\prec 1$  for all $k\le K$ and $t\le T$, 
which is formalized in the following proof.  

\begin{proof}[Proof of Proposition \ref{prop:zig}]
	Recall that $\ell_t\eta_t \gtrsim N^{\etaexp}$ for all $t \in [\tinit, T]$ by \eqref{eq:WN} and \eqref{eq:etabound}. 
	By assumption, the local laws hold for $G_t$ with tolerance $\xi \in (0, \etaexp/100)$ at time $t=\tinit$ in the sense of  Definition~\ref{def:local_laws}, hence we choose
	\begin{equation}
		\psi_{k,0}^\mathrm{av/iso} := N^\xi, \quad k \in\indset{\maxK}.
	\end{equation}
	Let $\nu> 0$ be a fixed tolerance exponent. Without loss of generality, we can assume that $\nu < \xi/(10K^2)$, and choose the family of control parameters
	\begin{equation} \label{eq:psi_choice}
		\pis{k}:= N^{\xi+k\nu},\quad \pav{k}:=N^{\xi+(k+1)\nu}, \quad k\in \indset{\maxK}.
	\end{equation}
	In particular, since $\psi^{\mathrm{av/iso}}_k > \psi^{\mathrm{av/iso}}_{k,0}$, the stopping time $\tau$, defined in \eqref{eq:tau_def}, satisfies $\tau > \tinit$ with very high probability. 
	Then it is straightforward to check from \eqref{eq:loss_exponents} that the master inequalities \eqref{eq:masters} imply
		\begin{equation}\label{eq:it}
		\max_{\tinit\le t \le \tau}\Psi_{k,t}^\mathrm{av/iso} \prec \psi_{k,0}^\mathrm{av/iso} + N^{-\etaexp/(2\maxK^2)}\psi_k^\mathrm{av/iso}, \quad k \in \indset{\maxK}.
	\end{equation}
	  Here we repeatedly used the explicit formulas for $\alpha_k, \beta_k$ from~\eqref{eq:loss_exponents},
	in particular explicit lower bounds on the exponents like $\alpha_{k+1}-\alpha_k \le \alpha_{K/2+1} \le 1/\sqrt{2K}$
	to effectively control the positive powers of $1/(\ell_t\eta_t)\lesssim N^{-\etaexp}$ in the right hand side of \eqref{eq:phi_iso_masters}. This gives rise to the prefactor of $N$ to a small negative power $-\etaexp/(2K^2)$ 
	in~\eqref{eq:it}. 	 
	In particular, by definition of the stopping times \eqref{eq:tau_def}, we deduce that $\tau = T$ with very high probability. Therefore, by definition of $\Psi_{k,t}^\mathrm{av/iso}$ in \eqref{eq:Psi_def}, we conclude that the 
	local laws hold for $G_t$ with tolerance $\xi + \nu$, uniformly in $t \in [\tinit, T]$ in the sense of Definition~\ref{def:local_laws}.
	 This concludes the proof of Proposition \ref{prop:zig}. 	
\end{proof}

The remainder of this section is dedicated to proving the master inequalities \eqref{eq:masters}. 
We will present the proof in a simplified situation, explaining why this entails no loss of generality at the end of Section~\ref{sec:masters_sec}.
While the definition of the control quantities  $\Psi$,
as well as Proposition~\ref{prop:zig}, were formulated for an arbitrary set of spectral parameters
with comparable imaginary parts, for notational simplicity throughout the proof, we assume a 
fixed $z_t$ trajectory. Specifically, for every $j\in \indset{\maxK}$, we assume $z_{j, t}=z_t$
or $z_{j, t}=\bar z_t$. This simplifies the first maximum in the definition  $\Psi$'s in \eqref{eq:Psi_def} 
to  $\max_{\bm z_t\in \{ z_t, \bar z_t\}^k}$.

 The starting point of the proof of~\eqref{eq:masters}
 under the simplification $z_{j, t} \in \{ z_t, \bar z_t\}$ 
 are the evolution equations for averaged and isotropic resolvent chains along the combination of the characteristic flow \eqref{eq:char_flow}
and the Ornstein-Uhlenbeck process \eqref{eq:zigOU}. These evolution equations are the content of the following lemma.
\begin{lemma}[Evolution Equations] \label{lemma:av_iso_evol}
	\textbf{(i)~Averaged Chains}.
	Define the quantity
	\begin{equation} \label{eq:G-M_kav}
		\mathcal{X}^{k}_{t} \equiv \mathcal{X}^{k}_{\bm z_t,t}(\bm x) := \Tr\bigl[(G - M)_{[1,k],t}(\bm x')S^{x_k}\bigr],
	\end{equation}
	and view it as a function of $\bm x \in \indset{N}^k$. 
	Then $\mathcal{X}^{k}_{t}(\bm x)$ satisfies the following evolution equation,
	\begin{equation} \label{eq:k_av_evol}
		\mathrm{d}\mathcal{X}^{k}_t = \biggl(\frac{k}{2}I+\bigoplus_{j=1}^k \mathcal{A}_{j,t}\biggr)\bigl[ \mathcal{X}^{k}_t\bigr]\mathrm{d}t 
		+ \mathrm{d}\mathcal{M}^\mathrm{av}_{[1,k],t}  + \mathcal{F}^\mathrm{av}_{[1,k],t} \mathrm{d}t,
	\end{equation}
	where, the $\mathcal{A}_{j,t}$, $\mathrm{d}\mathcal{M}^\mathrm{av}_{[1,k],t}$, and $\mathcal{F}^\mathrm{av}_{[1,k],t}$ are defined as follows.
	With the cyclic convention $z_{k+1,t} := z_{1,t}$, the linear operators $\mathcal{A}_{j,t}$ on $\mathbb{C}^N$ 
	are given by\footnote{Note that $\Theta(z) =\Theta(\bar z)$, while  $\Xi(z) = \overline{\Xi(\bar z})$,
	hence the definition of $\Theta_t=\Theta(z_{j, t})$ is
	insensitive whether $z_{j,t} = z_t$ or $z_{j,t} = \bar z_t$, unlike for $\Xi(z_{j, t})$.}
	\begin{equation} \label{eq:lin_prop_ops}
		\mathcal{A}_{j,t} :=  \frac{m(z_{j,t})m(z_{j+1,t})S}{1- m(z_{j,t})m(z_{j+1,t})S}
		 = \begin{cases}
			\Theta_t, \quad &\text{if}~ z_{j+1} = \overline{z_j},\\
			\Xi(z_{j,t}), \quad &\text{if} ~z_{j+1} = z_j 
		\end{cases}.
	\end{equation}
	The martingale differential term $\mathrm{d}\mathcal{M}^\mathrm{av}_{[1,k],t}$ denotes 
	\begin{equation} \label{eq:av_k_mart}
		\mathrm{d}\mathcal{M}^\mathrm{av}_{[1,k],t} \equiv 	\mathrm{d}\mathcal{M}^\mathrm{av}_{[1,k],t}(\bm x)  := \sum_{ab}\sqrt{S_{ab}}\Tr\bigl[ \partial_{ab}G_{[1,k],t}(\bm x')  S^{x_k}\bigr] \mathrm{d}\mathfrak{B}_{ab,t},
	\end{equation}
	where $\partial_{ab}:= \partial_{H_{ab}}$ denotes the complex partial derivative with respect to $H_{ab}$.
	
	Finally, $\mathcal{F}^\mathrm{av}_{[1,k],t} \equiv \mathcal{F}^\mathrm{av}_{[1,k],t}(\bm x) $ comprises the forcing terms: For $k=1$, $\mathcal{F}^\mathrm{av}_{[1,1],t}$ is given by  
	\begin{subequations}\label{eq:F_def}
		\begin{equation}  \label{eq:F1_def}
			\mathcal{F}^\mathrm{av}_{[1,1],t} = 
			\mathcal{F}^\mathrm{av}_{[1,1],t}(x_1) :=  \Tr\bigl[\mathcal{S}[G_{1,t}-m_{1,t}]\bigl( (G-M)_{[1,1],t}^{(x_1)}\bigr)\bigr],  
		\end{equation}
		while for $k \ge 2$, we define  
		\begin{equation} \label{eq:Fk_def}
			\begin{split}
				\mathcal{F}^\mathrm{av}_{[1,k],t} :=&~ \sum_{1\le j\le k} \Tr\bigl[G_{[1,j],t}\mathcal{S}[G_{j,t}-m_{j,t}] G_{[j,k],t}S^{x_k}\bigr]\\
				&+ \sum_{1 \le i < j \le k} \Tr\bigl[ \mathcal{S}\bigl[(G-M)_{[i,j],t}\bigr] (G-M)^{(x_k)}_{[j,i],t} \bigr] \\
				&+ \sum_{\substack{1\le i < j \le k \\ j-i \le k-2}} \Tr\bigl[\mathcal{S}\bigl[(G-M)_{[i,j],t}\bigr]M^{(x_k)}_{[j,i],t}\bigr] + \sum_{\substack{1\le i < j \le k \\ 2 \le j-i}} \Tr\bigl[\mathcal{S}\bigl[(G-M)^{(x_k)}_{[j,i],t}\bigr]M_{[i,j],t}\bigr],
			\end{split}
		\end{equation}	
	\end{subequations}
	where for $1\le i \le j \le k$, we denote
	\begin{equation} \label{eq:sub_chains} 
			G_{[j,i],t}^{(x_k)} \equiv G_{[j,i],t}^{(x_k)}(x_1, \dots, x_{i-1}, x_j,\dots, x_{k-1})  := G_{[j, k],t}S^{x_k} G_{[1,i],t}, 
	\end{equation}
	and use $M^{(x_k)}_{[j,i],t}$ to denote the corresponding deterministic approximations\footnote{ 
	We warn the reader that for the special $i=j$ case, $G_{[j,j],t}^{(x_k)}$ is a long chain and
	it is very different from $G_{[j,j],t}=G_{j,t}$ from~\eqref{eq:Gk_def} consisting of  a single resolvent.}.

	\textbf{(ii)~Isotropic Chains}.
	Fix entry indices $a,b \in \indset{N}$. Define the quantity
	\begin{equation} \label{eq:G-M_kiso}
		\mathcal{Y}^{k}_{t} \equiv \mathcal{Y}^{k}_{\bm z_t,t}(a,\bm x',b) := \bigl((G - M)_{[1,k],t}(\bm x')\bigr)_{ab},
	\end{equation}
	and view it as a function of $\bm x' \in \indset{N}^{k-1}$.
	Then $\mathcal{Y}^{k}_{t}$ satisfies the following evolution equation,
	\begin{equation} \label{eq:k_iso_evol}
		\mathrm{d}\mathcal{Y}^{k}_{t} = \biggl(\frac{k}{2}I+\bigoplus_{j=1}^{k-1} \mathcal{A}_{j,t}\biggr)\bigl[ \mathcal{Y}^{k}_{t}\bigr]\mathrm{d}t 
		+ \mathrm{d}\mathcal{M}^\mathrm{iso}_{[1,k],t}  + \mathcal{F}^\mathrm{iso}_{[1,k],t} \mathrm{d}t,
	\end{equation}
	where, $\mathcal{A}_{j,t}$ are defined in \eqref{eq:lin_prop_ops}, the martingale term $\mathrm{d}\mathcal{M}^\mathrm{iso}_{[1,k],t}$ denotes 
	\begin{equation} \label{eq:iso_k_mart}
		\mathrm{d}\mathcal{M}^\mathrm{iso}_{[1,k],t} \equiv 	\mathrm{d}\mathcal{M}^\mathrm{iso}_{[1,k],t}(a,\bm x',b)  := \sum_{cd}\sqrt{S_{cd}}\,\partial_{cd}\bigl( G_{[1,k],t}(\bm x')\bigr)_{ab} \mathrm{d}\mathfrak{B}_{cd,t},
	\end{equation}
	while $\mathcal{F}^\mathrm{iso}_{[1,k],t} \equiv \mathcal{F}^\mathrm{iso}_{[1,k],t}(a,\bm x',b) $ comprises the forcing terms:  For $k=1$, $\mathcal{F}^\mathrm{iso}_{[1,1],t}$ is given by  
	\begin{subequations}\label{eq:iso_F_def}
		\begin{equation} \label{eq:iso_F1_def}
			\mathcal{F}^\mathrm{iso}_{[1,1],t} = 
			\mathcal{F}^\mathrm{iso}_{[1,1],t}(a,b) := \bigl(G_{1,t}\mathcal{S}\bigl[G_{1,t}-m_{i,t}\bigr]G_{1,t}\bigr)_{ab},  
		\end{equation}
		 while for $k \ge 2$, we define  
		\begin{equation} \label{eq:iso_Fk_def}
			\begin{split}
				\mathcal{F}^\mathrm{iso}_{[1,k],t} :=&~ \sum_{1\le j\le k} \bigl(G_{[1,j],t}\mathcal{S}\bigl[G_{j,t}-m_{i,t}\bigr]G_{[j, k],t}\bigr)_{ab} \\
				&+ \sum_{1\le i < j \le k} \sum_q \bigl((G-M)_{[i,j],t}\bigr)_{qq} \bigl(G_{[1,i],t}S^qG_{[j, k],t}-M_{[1,i],[j,k],t}^{(q)}\bigr)_{ab}\\
				&+ \sum_{\substack{1 \le i < j \le k\\j-i\le k-2}} \sum_q \bigl((G-M)_{[i,j],t}\bigr)_{qq} \bigl(M_{[1,i],[j,k],t}^{(q)}\bigr)_{ab}\\
				&+ \sum_{\substack{1\le i < j \le k \\ 2 \le j-i}} \sum_q \bigl(M_{[i,j],t}\bigr)_{qq} \bigl(G_{[1,i],t}S^qG_{[j, k],t}-M_{[1,i],[j,k],t}^{(q)}\bigr)_{ab}\\
				&+\delta_{ab}\sum_q m_{1,t}m_{k,t} \bigl(I+\mathcal{A}_{k,t}\bigr)_{aq} 
				\mathcal{X}_{t}^k(\bm x',q).
			\end{split}
		\end{equation}
	\end{subequations}
	Here, for all $q\in\indset{N}$, the deterministic term $M_{[1,i],[j,k],t}^{(q)}$ 
	is the deterministic approximation to the chain $G_{[1,i],t}S^qG_{[j, k],t}$, given by
	\begin{equation} \label{eq:M_with_q}
		M_{[1,i],[j,k],t}^{(q)} := M\bigl(z_{1,t}, S^{x_1}, z_{2,t}, \dots, S^{x_{i-1}}, z_{i,t}, S^q, z_{j,t}, S^{x_{j+1}},\dots, z_{k-1,t}, S^{x_{k-1}}, z_{k,t} \bigr).
	\end{equation}
\end{lemma}
We defer the straightforward proof of Lemma \ref{lemma:av_iso_evol} to Section \ref{sec:evols}. 
The evolution equations  \eqref{eq:k_av_evol} and \eqref{eq:k_iso_evol} consist of three main components: 
a {\it linear part}, given by the  first term on their right-hand sides involving the linear 
operators $\mathcal{A}$ (along with trivial constants); 
a {\it martingale part} $ \mathrm{d}\mathcal{M}$; and a {\it forcing term} $\mathcal{F} \mathrm{d}t$. 
By treating the linear operators as generators of a semigroup, we can apply Duhamel's principle to rewrite 
these differential equations in integral form.
The resulting semigroup---also called the {\it propagator}--then acts on the initial condition,
the martingale term, and the forcing term. Crucially, 
the spatial decay  of each of these components (as functions of $\bm x$) is controlled by our size functions $\mathfrak{s}^{\mathrm{av/iso}}$, introduced in~\eqref{eq:sfunc_def}. 
Therefore, we now turn to analyzing the action of the propagators on such functions.

 Integrating the linear part of the (non-linear) evolution equations \eqref{eq:k_av_evol} and \eqref{eq:k_iso_evol} using Duhamel's principle, we encounter three kinds of propagators: The \emph{saturated} propagator $\mathcal{P}_{s,t}$, defined as
\begin{equation} \label{eq:Psat_def}
	\mathcal{P}_{s,t} := \exp\biggl\{ \int_s^t (I + \Theta_r)\mathrm{d}r \biggr\} = \frac{|m_t|^2}{|m_s|^2}\frac{1-|m_s|^2 S}{1-|m_t|^2 S}, \quad 0 \le s \le t\le T,
\end{equation}
the unsaturated propagator $\mathcal{Q}_{s,t}$, given by
\begin{equation} \label{eq:Q_desat_def}
	\mathcal{Q}_{s,t} := \exp\biggl\{ \int_s^t (I + \Xi_r)\mathrm{d}r \biggr\} = \frac{m_t^2}{m_s^2}\frac{1-m_s^2 S}{1-m_t^2 S}, \quad 0 \le s \le t\le T,
\end{equation}
and its complex conjugate $\overline{\mathcal{Q}}_{s,t}$. 
The second relations in \eqref{eq:Psat_def}--\eqref{eq:Q_desat_def} follow from \eqref{def:ThetaXi} and \eqref{eq:dTheta}.

Therefore, for every $k\in\indset{\maxK}$, we consider the class of linear operators $\other{\mathcal{P}}^{k}_{s,t}$ on the space $(\mathbb{C}^{N})^{\otimes k}$ of the form
\begin{equation} \label{eq:general_prop}
	\other{\mathcal{P}}^{k}_{s,t} := \bigotimes_{j=1}^{k}\other{\mathcal{P}}^{(j)}_{s,t}, \quad \text{with}\quad  \other{\mathcal{P}}_{s,t}^{(j)} \in \bigl\{ \mathcal{P}_{s,t}, \mathcal{Q}_{s,t}, \overline{\mathcal{Q}}_{s,t}, 
	I \bigr\}, \quad 0 \le s \le t\le T,
\end{equation}
where the tilde  in the notation indicates that for each factor we take one of the possible four choices.

The following lemma asserts that the target estimates on the isotropic and averaged resolvent chains are stable under the action of linear propagators belonging to the class \eqref{eq:general_prop}, with one exception: for averaged chains
in the case where all propagators are saturated.
\begin{lemma} [Propagator Estimate] \label{lemma:good_props}
	Let $s,t \in [0,T]$ be a pair of times satisfying $s\le t$, and let $\other{\mathcal{P}}^{k}_{s,t}$ be a linear operator of the form \eqref{eq:general_prop} for some $k \in\indset{\maxK}$. Let $\varphi >0 $ be a positive control parameter. 
	Then for any function $f_{ab}(\bm x)$ of $\bm x \in \indset{N}^{k}$ and all indices $a,b \in \indset{N}$,
	\begin{equation} \label{eq:iso_prop_bound}
		\forall \bm x~~ \bigl\lvert f_{ab}(\bm x) \bigr\rvert \lesssim \varphi \, \mathfrak{s}_{k+1,s}^\mathrm{iso}(a,\bm x, b) \quad \Longrightarrow \quad 
		\forall \bm x~~ \bigl\lvert \other{\mathcal{P}}^{k}_{s,t}\bigl[f_{ab} \bigr]  (\bm x)\bigr\rvert \lesssim \varphi \,\sqrt{\frac{\ell_t\eta_t }{\ell_s\eta_s}}\, \mathfrak{s}_{k+1,t}^\mathrm{iso}(a,\bm x , b).
	\end{equation} 
	
	If we assume additionally that there exists at least one
	 index $j \in \indset{k}$, such that $\other{\mathcal{P}}_{s,t}^{(j)} \neq \mathcal{P}_{s,t}$, then for any function $g(\bm x)$ of $\bm x \in \indset{N}^{k}$,
	\begin{equation} \label{eq:av_prop_bound}
		\forall \bm x~~ \bigl\lvert g(\bm x) \bigr\rvert \lesssim \varphi \, \mathfrak{s}_{k,s}^\mathrm{av}(\bm x) \quad \Longrightarrow \quad 
		\forall \bm x~~ \bigl\lvert \other{\mathcal{P}}^{k}_{s,t}\bigl[g \bigr]  (\bm x)\bigr\rvert \lesssim \varphi \, \frac{\ell_t\eta_t}{\ell_s\eta_s}  \, \mathfrak{s}_{k,t}^\mathrm{av}(\bm x).
	\end{equation}
	
	Furthermore, estimates \eqref{eq:iso_prop_bound} and \eqref{eq:av_prop_bound} also hold with the linear operator $\other{\mathcal{P}}^{k}_{s,t}$ replaced by its entry-wise absolute value\footnote{That is, $|\other{\mathcal{P}}^{k}_{s,t}|_{\bm x\bm y} := |(\other{\mathcal{P}}^{k}_{s,t})_{\bm x\bm y}|$ for all $\bm x, \bm y \in \indset{N}^k$. } $|\other{\mathcal{P}}^{k}_{s,t}|$.
\end{lemma}  
We prove Lemma \ref{lemma:good_props} in Section \ref{sec:props} below. 
We now briefly comment on how Lemma \ref{lemma:good_props} will be applied and highlight its main limitation.

\bigskip
{\bf Smoothing and consistency.} 
 Since $\ell_t\eta_t\le \ell_s\eta_s$,
both estimates~\eqref{eq:iso_prop_bound}--\eqref{eq:av_prop_bound} represent an  improvement
over the initial bound on $f$ in terms of the scalar prefactor in front of $\mathfrak{s}$; we will call it the {\it smoothing effect}.
In fact, both bounds allow propagating the  target estimates 
 $\mathcal{X}^{k}_{t} \prec (\ell_t\eta_t)^{\beta_k}\mathfrak{s}_{k,t}^\mathrm{av}$ 
and $\mathcal{Y}^{k}_{t}\prec (\ell_t\eta_t)^{\alpha_k}\mathfrak{s}_{k,t}^\mathrm{iso}$ using Duhamel's principle---a property we refer to as {\it consistency}.  
For example, consistency is straightforward for the propagation of the initial condition at time $\tinit$ to a later time $t\ge \tinit$:
\begin{equation}\label{eq:ini}
 \bigl\lvert \other{\mathcal{P}}_{\tinit,t}^k \bigl[\mathcal{X}^{k}_{\tinit}\bigr]\bigr\rvert \prec (\ell_{\tinit}\eta_{\tinit})^{\beta_k}\other{\mathcal{P}}_{\tinit,t}^k\bigl[\mathfrak{s}_{k,\tinit}^\mathrm{av}\bigr]
 \lesssim \frac{\ell_t\eta_t}{(\ell_{\tinit}\eta_{\tinit})^{1-\beta_k}} \mathfrak{s}_{k,t}^\mathrm{av} \lesssim (\ell_t\eta_t)^{\beta_k}\mathfrak{s}_{k,t}^\mathrm{av},
\end{equation}
where we used  \eqref{eq:av_prop_bound}, $\ell_t\eta_t\le \ell_{\tinit}\eta_{\tinit}$, and $\beta_k\le1$.  
For the other terms on the right-hand side of \eqref{eq:k_av_evol} behave similarly. Heuristically, these terms can be compared to $|\mathcal{X}^{k}_{s}|/\eta_s$, where the additional factor $1/\eta_s$ arises morally from the differentiation of a resolvent. Applying 
the integration rules~\eqref{eq:int_rules}, we then obtain the consistent estimate
\begin{equation}\label{eq:fo}
  \int_{\tinit}^t \other{\mathcal{P}}_{s,t}^k  \Big[\frac{ |\mathcal{X}^{k}_{s}|}{\eta_s}\Big] \mathrm{d} s \prec
   \int_{\tinit}^t \frac{(\ell_s\eta_s)^{\beta_k}}{\eta_s}   \other{\mathcal{P}}_{s,t}^k   \bigl[ \mathfrak{s}_{k,s}^\mathrm{av}\bigr] \mathrm{d} s
   \lesssim  \int_{\tinit}^t \frac{\ell_t\eta_t}{(\ell_s\eta_s)^{1-\beta_k}}  \, \mathfrak{s}_{k,t}^\mathrm{av}  \frac{\mathrm{d} s}{\eta_s} 
   \lesssim 
   (\ell_t\eta_t)^{\beta_k}\mathfrak{s}_{k,t}^\mathrm{av}.
 \end{equation}  

The consistency of \eqref{eq:iso_prop_bound} for the isotropic case is analogous. 
Note that consistency is a necessary but not sufficient condition
 to guarantee closability  of the resulting master inequalities, we will also need to establish its self-improving structure
 as explained after Proposition~\ref{prop:masters}.
  
 We remark that the actual proofs of the consistency are much more involved. 
First, the forcing term is more intricate than the heuristic expression $|\mathcal{X}^{k}_{s}|/\eta_s$;
 it contains products of shorter chains. However, using the corresponding estimates---encoded in the definition of 
the stopping time~\eqref{eq:tau_def}---along with the structural relations among the size functions $\mathfrak{s}_{k,t}$
(which follow directly from  the properties of the admissible control 
functions~\eqref{eq:Ups_majorates_notime}--\eqref{eq:suppressed_convol_notime}), we ultimately show
that the forcing term admits an upper bound comparable to that on $|\mathcal{X}^{k}_{s}|/\eta_s$. 
Second, the martingale term must be estimated stochastically via its quadratic variation 
and the Burkholder-Davis-Gundy (BDG) inequality. Besides this technical complication, the 
main feature of the martingale term is that its quadratic variation contains chains of length up to $2k+2$,
which are not directly controlled by the stopping time when $k\ge K/2$. 
To ensure the closability of the master inequalities, these chains first need to be broken up into shorter chains of length at most $K$. 
 To this end, we employ three distinct {\it reduction strategies}  (see the explanation around~\eqref{eq:QG^2_bound} later). 
These are used pragmatically depending on the context,  but all of them inherently result in a deterioration compared to the optimal bound.
Typically the martingale term arising from the longest chain is the most critical and its reduction incurs the highest cost. As already mentioned after~\eqref{eq:Psi_def}, we
account for this deterioration by carefully choosing  the loss exponents.

\medskip

{\bf Limitation of Lemma~\ref{lemma:good_props}.}
The smoothing effect is present for all the propagators $\other{\mathcal{P}}_{s,t}^k$
with one exception: in the averaged case  when \textbf{all} $k$ propagators in $\other{P}_{s,t}^k$ are saturated---that is, $\other{\mathcal{P}}_{s,t}^{(j)} = \mathcal{P}_{s,t}$ for all $j\in\indset{k}$. In this special case  the implication in \eqref{eq:av_prop_bound} is violated. Specifically, 
rather than gaining the small factor $(\ell_t\eta_t)/(\ell_s\eta_s)$, the propagated estimate 
deteriorates by a large factor  $\ell_t/\ell_s$, namely,
\begin{equation} \label{eq:bad_av_prop}
	\forall \bm x~~ \bigl\lvert g(\bm x) \bigr\rvert \lesssim \varphi \, \mathfrak{s}_{k,s}^\mathrm{av}(\bm x) \quad \Longrightarrow \quad 
	\forall \bm x~~ \bigl\lvert \mathcal{P}^{\otimes k}_{s,t}\bigl[g \bigr]  (\bm x)\bigr\rvert \lesssim \frac{\ell_t}{\ell_s}  \, \varphi \,  \mathfrak{s}_{k,t}^\mathrm{av}(\bm x),
\end{equation}
and it is straightforward to verify that \eqref{eq:bad_av_prop} is optimal in general.
Since typically $\ell_t/\ell_s\gg 1$,  this estimate is not even sufficient to consistently control the initial term
in the Duhamel's formula.
We emphasize that this problem arises only for the averaged chains; for isotropic chains the smoothing effect is always effective.  

To remedy this limitation of Lemma~\ref{lemma:good_props} 
and restore consistency to include this special case,  we need a new idea. 
It follows from \eqref{eq:lin_prop_ops} that all $k$ propagators in $\other{P}_{s,t}^k$ can be saturated if and only if the resolvent chain $G_{[1,k],t}$ has even length and $G(z)$ and $G^*(z)$ alternate in the chain, that is,
\begin{equation} \label{eq:sat_cahins}
	G_{i+1,t} = G_{i,t}^*, \quad i\in\indset{k-1}.
\end{equation}
By analogy with the propagators, we call such chains \emph{saturated}\footnote{ 
	Saturated chains are actually the most important ones since they are the biggest and they are inevitable. For example, the basic $2$-chain $\langle GS^xG^*S^x\rangle$ arising in the quadratic variation of  $\langle G-m\rangle$ is saturated.
},
they give rise to the special case for which Lemma~\ref{lemma:good_props} does not apply.
 Equivalently, saturated chains correspond to $\bm z_t \in \{\bm z_{t,\mathrm{sat}},  \overline{\bm z_{t,\mathrm{sat}}}\}$, where 
\begin{equation} \label{eq:alt_z}
	\bm z_{t,\mathrm{sat}} := (z_t, \overline{z_t}, z_t, \overline{z_t}, \dots, z_t, \overline{z_t}) \in (\mathbb{C}\backslash\mathbb{R})^k, \quad k \in \indset{\maxK} \text{ -- even}.
\end{equation}
Since the large factor $\ell_t/\ell_s$ in \eqref{eq:bad_av_prop} is generally not compatible with the target estimate on the averaged resolvent chains that operates with powers of $\ell\eta$
 (see \eqref{eq:Psi_def}), at least one\footnote{In the proof, we actually remove all of them.} of the $k$ linear terms in the evolution equation \eqref{eq:k_av_evol} for saturated averaged chains has to be estimated
independently before applying Duhamel's principle. 
To this end, we study the linear term as a separate 
quantity\footnote{We remark that similar separate estimates of the linear terms could also be obtained
for the non-saturated chains, hence their propagators could have also been removed, however the notation becomes slightly more cumbersome, as the complex conjugates of the $z$'s must be explicitly tracked.  Here we follow the
simpler route of estimating the linear terms separately only when it is really necessary, i.e. for saturated chains,
where the difference between $z$ and $\bar z$ is absent. }
 and derive a self-consistent estimate for it which does not involve $\pav{k}$ without additional suppressing factors. This is the content of the crucial Proposition \ref{prop:lin_term} below.

Let $\Theta_t^{(j)} := I^{\otimes (j-1)}\otimes\Theta_t\otimes I^{\otimes (k-j)}$ be the linear operator that acts by $\Theta_t$ on the $j$-th component in the tensor-product space $(\mathbb{C}^N)^{\otimes k}$, for $j \in \indset{k}$. 
Recall from~\eqref{eq:lin_prop_ops} that for saturated chains $\mathcal{A}_{j,t}= \Theta_t$ for all $j$. 
\begin{prop}[Linear Term Estimate] \label{prop:lin_term}
	Let $k\in\indset{\maxK}$ be an even integer. Assume that \eqref{eq:zig_init} holds, and let $\tau$ be the stopping time defined in \eqref{eq:tau_def}, then, for any $j \in \indset{k}$ and any $\bm x \in \indset{N}^k$, the linear terms satisfy the estimate
	\begin{equation} \label{eq:lin_term_est}
		\max_{\tinit\le s\le t \wedge \tau}\frac{\eta_s\bigl\lvert \Theta_s^{(j)} \bigl[ \mathcal{X}^{k}_{\bm z_s, s}\bigr](\bm x)\bigr\rvert}{(\ell_s\eta_s)^{\beta_k}\mathfrak{s}^\mathrm{av}_{k,s}(\bm x)} \prec  
		\varphi_{k,t\wedge\tau}^\mathrm{av}, \quad \bm z_s \in \{\bm z_{s,\mathrm{sat}},  \overline{\bm z_{s,\mathrm{sat}}}\}, 
	\end{equation}
	uniformly in $t\in [\tinit,T]$, where $\varphi_{k,t\wedge\tau}^\mathrm{av}$ is as defined in \eqref{eq:phi_av_masters}, and the alternating vector of spectral parameters $\bm z_{t,\mathrm{sat}}$ is defined in \eqref{eq:alt_z}.
\end{prop}
We defer the proof of Proposition \ref{prop:lin_term} to Section \ref{sec:reg_sec}. Using \eqref{eq:lin_prop_ops}, the linear terms in \eqref{eq:k_av_evol} for saturated chains will be treated as additional forcing terms
(estimated by~\eqref{eq:lin_term_est}), 
hence avoiding the excess of saturated propagators. In fact, the same reasoning is behind treating the term $\sum_q m_{1,t}m_{k,t} \bigl( I+\mathcal{A}_{k,t} \bigr)_{aq}
\mathcal{X}_{t}^k(\bm x',q)$ in the last line of \eqref{eq:iso_Fk_def} as a forcing term instead of a linear term. This is done despite the fact that, using the identity
\begin{equation}
	\sum_q m_{1,t}m_{k,t} \bigl( I+\mathcal{A}_{k,t} \bigr)_{aq} 
	\mathcal{X}_{t}^k(\bm x',q) = \sum_{q}\bigl(\mathcal{A}_{k,t}\bigr)_{aq}\mathcal{Y}^k_t(q,\bm x', q),
\end{equation} 
it could have been included in the linear part of \eqref{eq:k_iso_evol}. This approach, where such linear terms are reclassified as forcing terms, with their estimation delegated to a separate equation (in this case, the averaged law), was also adopted in prior work \cite{Cipolloni2022Optimal}.

The final ingredient that goes into the proof of the master inequalities \eqref{eq:masters} are the estimates on the martingale and forcing terms that appear on the right-hand sides of \eqref{eq:k_av_evol} and \eqref{eq:k_iso_evol}.

\begin{lemma}[Martingale Estimates] \label{lemma:mart_est} 
Assume only the weaker condition \eqref{eq:true_convol_notime} instead of \eqref{eq:convol_notime}--\eqref{eq:suppressed_convol_notime}.  
	Let $k\in\indset{\maxK}$.
	For all $s, t \in [0,T]$, satisfying $s\le t$, let $\otherhat{\mathcal{P}}^k_{s,t}$ be a linear operator on $(\mathbb{C}^N)^k$ with matrix elements $(\otherhat{\mathcal{P}}^k_{s,t})_{\bm x \bm a} := (\otherhat{\mathcal{P}}^k_{s,t})_{x_1, \dots, x_k, a_1, \dots, a_k}$ for $\bm x,\bm a \in \indset{N}^k$ satisfying $|(\otherhat{\mathcal{P}}^k_{s,t})_{\bm x \bm a}|\le N^C$.
	
	Let $\tau$ be the stopping time defined in \eqref{eq:tau_def},  and let $\mathrm{d}\mathcal{M}^\mathrm{av}_{[1,k],s}(\bm x)$ be the martingale differential defined in \eqref{eq:av_k_mart}. Then, for all $t\in [\tinit,T]$, 
	\begin{equation} \label{eq:mart_bound}
		\begin{split}
			\biggl\lvert \int_{\tinit}^{t\wedge\tau} \otherhat{\mathcal{P}}^k_{s,t\wedge\tau}\bigl[\mathrm{d}\mathcal{M}^\mathrm{av}_{[1,k],s}\bigr](\bm x) \biggr\rvert
			\prec&~ \biggl(\int_{\tinit}^{t\wedge\tau}  \bigl(\varphi_{k,s}^\mathrm{av,qv}\bigr)^2 (\ell_s\eta_s)^{2\beta_k} \biggl\lvert\sum_{\bm a\in\indset{N}^k} \bigl\lvert(\otherhat{\mathcal{P}}^k_{s,t\wedge\tau})_{\bm x\bm a}\bigr\rvert \mathfrak{s}_{k,s}^\mathrm{av}(\bm a) \biggr\rvert^2 \frac{\mathrm{d}s}{\eta_s}
			\biggr)^{1/2}\\
			&+ \mathfrak{s}_{k,t\wedge\tau}^\mathrm{av}(\bm x),
		\end{split}
	\end{equation}
	where the time-dependent control parameter $\varphi_{k,s}^\mathrm{av,qv}$ is defined in \eqref{eq:phi_avqv}.
	
	Assume additionally\footnote{
	In \eqref{eq:iso_mart_bound}, we formulated the estimate for 
	 resolvent chain of length  $k+1$ for notational similarity with its counterpart~\eqref{eq:mart_bound}
	 so that both contain $k$ propagators. In its actual usage  in \eqref{eq:k_iso_evol} we will use 
	 \eqref{eq:iso_mart_bound} for $k$ instead of $k+1$.}
that $k\in \indset{\maxK-1}$. Fix $a,b \in \indset{N}$, let $\mathrm{d}\mathcal{M}^\mathrm{iso}_{[1,k+1],s}(\bm x) \equiv \mathrm{d}\mathcal{M}^\mathrm{iso}_{[1,k+1],s}(a,\bm x,b)$ be the quantity defined in \eqref{eq:iso_k_mart}. Then, for all $t\in [t_0,T]$,
	\begin{equation} \label{eq:iso_mart_bound}
		\begin{split}
			\biggl\lvert \int_{\tinit}^{t\wedge\tau} \otherhat{\mathcal{P}}^{k}_{s,t\wedge\tau}\bigl[\mathrm{d}\mathcal{M}^\mathrm{iso}_{[1,k+1],s}\bigr](\bm x) \biggr\rvert
			\prec&~ \biggl(\int_{\tinit}^{t\wedge\tau}  \frac{\bigl(\varphi_{k+1,s}^\mathrm{iso,qv}\bigr)^2}{(\ell_s\eta_s)^{-2\alpha_{k+1}}}  \biggl\lvert\sum_{\bm c\in\indset{N}^{k}} \bigl\lvert(\otherhat{\mathcal{P}}^{k}_{s,t\wedge\tau})_{\bm x\bm c}\bigr\rvert \mathfrak{s}_{k+1,s}^\mathrm{iso}(a,\bm c,b) \biggr\rvert^2 \frac{\mathrm{d}s}{\eta_s}
			\biggr)^{1/2}\\
			&+ \mathfrak{s}_{k+1,t\wedge\tau}^\mathrm{iso}(a,\bm x,b),
		\end{split}
	\end{equation}
	where the time-dependent control parameter $\varphi_{k+1,s}^\mathrm{iso,qv}$ is defined in \eqref{eq:phi_isoqv}.
\end{lemma} 
 We will  use this lemma for $\otherhat{\mathcal{P}}^{k} =\other{\mathcal{P}}^{k}$ or 
$\mathcal{P}_{s,t\wedge\tau}^{\otimes k}$ sometimes
with a small twist as in~\eqref{eq:Gcirc_solve} or \eqref{eq:2ring_mart_bound}.
We also remark that the $\mathfrak{s}^\mathrm{ave}$ and $\mathfrak{s}^\mathrm{iso}$
terms in the last lines of~\eqref{eq:mart_bound}--\eqref{eq:iso_mart_bound} just conveniently overestimate
a tiny term of order $N^{-kD'}$ coming from exceptional events in the stochastic domination bounds.

\begin{lemma}[Forcing Term Estimates] \label{lemma:forcing}
	Assume only the weaker condition \eqref{eq:true_convol_notime} instead of \eqref{eq:convol_notime}--\eqref{eq:suppressed_convol_notime}.  
	Let $\tau$ be the stopping time defined in \eqref{eq:tau_def}, and let $s \in [\tinit,\tau]$. 
	Then, for all $k\in \indset{\maxK}$, for any $\bm z_s \in \{z_s, \overline{z}_s\}^k$, the quantity $\mathcal{F}^\mathrm{av}_{[1,k],t}(\bm x)$, defined in \eqref{eq:F1_def}--\eqref{eq:Fk_def}, satisfies
	\begin{equation} \label{eq:av_forcing_bound}
		\frac{\eta_s\bigl\lvert \mathcal{F}^\mathrm{av}_{[1,k],s}(\bm x) \bigr\rvert}{(\ell_s\eta_s)^{\beta_k}\mathfrak{s}_{k,s}^\mathrm{av}(\bm x)} 
		\lesssim \sum\limits_{j=2}^{k-1} \frac{\pav{j}}{(\ell_s\eta_s)^{\beta_k-\beta_j}} \biggl(1+\frac{\pis{k+2-j}}{(\ell_s\eta_s)^{1/2-\alpha_{k+2-j}}}\biggr) + \frac{\pis{2}\pav{k}}{\sqrt{\ell_s\eta_s}} 
		+ \varphi_{k,s}^\mathrm{av,f},
	\end{equation} 
	for all $\bm x\in\indset{N}^k$, where the time-dependent control parameter $\varphi_{k,s}^\mathrm{av,f}$ is defined in \eqref{eq:phi_av_force}.
 
	Moreover, for $s\in[\tinit,\tau]$ the quantity $\mathcal{F}^\mathrm{iso}_{[1,k],s}(a,\bm x',b)$, defined in \eqref{eq:iso_F1_def}--\eqref{eq:iso_Fk_def}, satisfies
	\begin{equation} \label{eq:iso_forcing_bound}
		\frac{\eta_s\bigl\lvert \mathcal{F}^\mathrm{iso}_{[1,k],s}(a,\bm x',b) \bigr\rvert}{(\ell_s\eta_s)^{\alpha_k}\mathfrak{s}_{k,s}^\mathrm{iso}(a,\bm x',b)} 
		\lesssim \sum_{j=2}^k \pis{j}  \biggl( \frac{\sqrt{N\eta_s}}{\ell_s\eta_s} \pis{k-i+1} + \frac{\pis{k-j+2}}{\sqrt{\ell_s\eta_s}}\biggr)  + \sum_{j=2}^{k-1} \frac{\pis{j}}{(\ell_s\eta_s)^{\alpha_k-\alpha_j}} + \frac{\pav{k}}{(\ell_s\eta_s)^{1/2 +\alpha_k-\beta_k 
			}} ,
	\end{equation} 
	for all $k\in \indset{\maxK}$, all $\bm z_s \in \{z_s, \overline{z}_s\}^k$, $\bm x' \in \indset{N}^{k-1}$, and all $a,b \in \indset{N}$.
\end{lemma}

We prove Lemmas \ref{lemma:mart_est} and \ref{lemma:forcing} in Sections \ref{sec:mart} and \ref{sec:forcing}, respectively.

\subsection{Completing the proof of master inequalities} \label{sec:masters_proof}
 Equipped with Lemmas \ref{lemma:av_iso_evol}--\ref{lemma:forcing} and Proposition \ref{prop:lin_term}, we are ready to prove the Proposition \ref{prop:masters}.
\begin{proof}[Proof of Master Inequalities (Proposition \ref{prop:masters})]
	We begin by proving \eqref{eq:av_masters}. Fix a chain length $k\in\indset{\maxK}$ and an a vector of spectral parameters $\bm z_t \in \{z_t, \overline{z}_t\}^k$. 
	
	Since saturated averaged chains require special treatment, we first consider non-saturated chains, that is, we assume that $\bm z_t \notin \{ \bm z_{t,\mathrm{sat}}, \overline{\bm z_{t,\mathrm{sat}}}\}$. In this case, we apply Duhamel's principle to the evolution equation \eqref{eq:k_av_evol}, and rewrite $\bigoplus_{j=1}^k \mathcal{A}_{j,t} + k/2$ as $\bigoplus_{j=1}^k (I+\mathcal{A}_{j,t}) - k/2$, to obtain
	\begin{equation} \label{eq:k_av_solve}
		\mathcal{X}^k_{\tau} = \mathrm{e}^{k(\tinit-\tau)/2} \other{\mathcal{P}}^k_{\tinit, \tau} \bigl[\mathcal{X}^k_{\tinit}\bigr] + \int_{\tinit}^\tau \mathrm{e}^{k(s-\tau)/2} \other{\mathcal{P}}^k_{s, \tau} \bigl[\mathrm{d}\mathcal{M}_{[1,k],s}^\mathrm{av}\bigr] + \int_{\tinit}^\tau \mathrm{e}^{k(s-\tau)/2} \other{\mathcal{P}}^k_{s, \tau} \bigl[\mathcal{F}_{[1,k],s}^\mathrm{av}\bigr] \mathrm{d}s,
	\end{equation}
	where $\tau$ is the stopping time defined in \eqref{eq:tau_def}, and with $\mathcal{A}_{j,t}$ being the linear operators from \eqref{eq:lin_prop_ops}, the propagator $\other{\mathcal{P}}^k_{s, \tau}$ is given by
	\begin{equation} \label{eq:k_av_propagator}
		\other{\mathcal{P}}^k_{s, t} := \exp\biggl\{ \int_{s}^t \bigoplus_{j=1}^k \bigl(I + \mathcal{A}_{j,r}\bigr)\mathrm{d}r \biggr\}, \quad \tinit \le s\le t\le T.
	\end{equation}
	Clearly, $\other{\mathcal{P}}^k_{s, t}$ belongs to the class \eqref{eq:general_prop}. Moreover, since the chain corresponding to $\bm z_t$ is not saturated by assumption, $\other{\mathcal{P}}^k_{s, t}$ contains at least one non-saturated propagator. Hence, \eqref{eq:av_prop_bound} is applicable.
	We now estimate the terms on the right-hand side of \eqref{eq:k_av_solve} one-by-one. Since $s,\tinit \le T \lesssim 1$, the exponential factors $\mathrm{e}^{k(\tinit-\tau)}$ $\mathrm{e}^{k(s-\tau)}$ are bounded by a constant, and hence can safely be ignored. 
	
	It follows from \eqref{eq:Psi_def}, and \eqref{eq:av_prop_bound} in Lemma \ref{lemma:good_props}, that the first term on the right-hand side of \eqref{eq:k_av_solve} admits the bound
	\begin{equation} \label{eq:k_av_init}
		\bigl\lvert  \other{\mathcal{P}}^k_{\tinit, \tau} \bigl[\mathcal{X}^k_{\tinit}\bigr](\bm x)\bigr\rvert \lesssim \Psi_{k,\tinit}^\mathrm{av} (\ell_{\tinit}\eta_{\tinit})^{\beta_k} \frac{\ell_\tau \eta_\tau}{\ell_{\tinit}\eta_{\tinit}} \mathfrak{s}_{k,\tinit}^\mathrm{av}(\bm x) \lesssim \Psi_{k,\tinit}^\mathrm{av} (\ell_\tau\eta_\tau)^{\beta_k} \mathfrak{s}_{k,\tinit}^\mathrm{av}(\bm x), \quad \bm x\in\indset{N}^k,
	\end{equation}
	where we used $\beta_k \le 1$ and the monotonicity property
	\begin{equation} \label{eq:elleta_monot}
		\ell_s\eta_s \ge \ell_t\eta_t, \quad 0\le s \le t\le T.
	\end{equation}
	
	Next, we estimate the second term on the right-hand side of \eqref{eq:k_av_solve}. It follows from \eqref{eq:av_prop_bound} that 
	\begin{equation}
		\sum_{\bm a\in\indset{N}^k} \bigl\lvert (\other{\mathcal{P}}^k_{s,t})_{\bm x\bm a}\bigr\rvert \mathfrak{s}_{k,s}^\mathrm{av}(\bm a)  \lesssim \frac{\ell_t\eta_t}{\ell_s\eta_s}\mathfrak{s}_{k,t}^\mathrm{av}(\bm x), \quad 0 \le s\le t\le T.
	\end{equation}
	Therefore, using \eqref{eq:mart_bound} in Lemma \ref{lemma:mart_est}, we obtain, for all $\bm x \in \indset{N}^k$,
	\begin{equation} \label{eq:k_av_mart_int}
		\begin{split}
			\biggl\lvert\int_{\tinit}^{t\wedge\tau}  \other{\mathcal{P}}^k_{s,t\wedge\tau}\bigl[\mathrm{d}\mathcal{M}^\mathrm{av}_{[1,k],s}\bigr](\bm x) \biggr\rvert
			&\prec \biggl(1+\int_{\tinit}^{\tau}  \bigl(\varphi_{k,s}^\mathrm{av,qv}\bigr)^2 (\ell_s\eta_s)^{2\beta_k} \biggl\lvert\frac{\ell_\tau\eta_\tau}{\ell_s\eta_s}\biggr\rvert^2 \frac{\mathrm{d}s}{\eta_s}
			\biggr)^{1/2} \mathfrak{s}_{k,t}^\mathrm{av}(\bm x)\\
			&\prec \varphi_{k,\tau}^\mathrm{av,qv}(\ell_\tau\eta_\tau)^{\beta_k}\, \mathfrak{s}_{k,t}^\mathrm{av}(\bm x),
		\end{split}
	\end{equation}
	where in the second step we used \eqref{eq:elleta_monot}, $\varphi_{k,t}^\mathrm{av,qv} \ge 1$ by \eqref{eq:phi_avqv}, $\beta_k < 1$, and the integration rules \eqref{eq:int_rules}.

	Finally, we estimate the contribution of the third term on the right-hand side of \eqref{eq:k_av_solve}. It follows from \eqref{eq:av_forcing_bound} in Lemma \ref{lemma:forcing}, and \eqref{eq:av_prop_bound} in Lemma \ref{lemma:good_props}, that
	\begin{equation} \label{eq:k_av_forcing_integral}
		\begin{split}
			\biggl\lvert\int_{\tinit}^\tau \other{\mathcal{P}}^k_{s, \tau} \bigl[\mathcal{F}_{[1,k],s}^\mathrm{av}\bigr] \mathrm{d}s\biggr\rvert \lesssim&~ 
			\biggl(\sum\limits_{j=2}^{k-1} \frac{\pav{j}}{(\ell_\tau\eta_\tau)^{\beta_k-\beta_j}} \biggl(1+\frac{\pis{k+2-j}}{(\ell_\tau\eta_\tau)^{1/2-\alpha_{k+2-j}}}\biggr) + \frac{\pis{2}\pav{k}}{\sqrt{\ell_\tau\eta_\tau}} 
			+ \varphi_{k,\tau}^\mathrm{av,f}\biggr)\\
			&~\times\biggl(\int_{\tinit}^\tau (\ell_s\eta_s)^{\beta_k}\frac{\ell_\tau\eta_\tau}{\ell_s\eta_s} \frac{\mathrm{d}s}{\eta_s}\biggr)\,\mathfrak{s}_{k,\tau}^\mathrm{av}.
		\end{split}
	\end{equation}
	It follows from \eqref{eq:int_rules} that the integral on the right-hand side of \eqref{eq:k_av_forcing_integral} is bounded by $(\ell_\tau\eta_\tau)^{\beta_k}$. Therefore, combining \eqref{eq:k_av_solve}, \eqref{eq:k_av_init}, \eqref{eq:k_av_mart_int}, \eqref{eq:k_av_forcing_integral}, we conclude by \eqref{eq:Psi_def} that the master inequality \eqref{eq:av_masters} holds for non-saturated chains. 
	
	Next, we prove \eqref{eq:av_masters} for saturated chains. The proof follows the same outline as the non-saturated case, except we treat the linear term $\Theta_t^{\oplus k}[\mathcal{X}^k_t]$ in \eqref{eq:k_av_evol} as a forcing term that we estimate by using \eqref{eq:lin_term_est} in Proposition \ref{prop:lin_term} since  $\Theta_t^{\oplus k} =\sum_{j=1}^k \Theta_t^{(j)}$. Since all linear terms are removed from the equation before Duhamel, the propagator $\mathrm{e}^{k(s-t)/2}\other{\mathcal{P}}^k_{s,t}$ is replaced by the operator $\mathrm{e}^{k(t-s)/2} I^{\otimes k}$. It follows from \eqref{eq:lin_term_est} that 
	\begin{equation} \label{eq:lin_forcing}
		\int_{\tinit}^\tau \mathrm{e}^{k(\tau-s)/2}\Theta_s^{\oplus k}[\mathcal{X}^k_s](\bm x) \mathrm{d}s \prec \int_{\tinit}^\tau \varphi_{k,s}^\mathrm{av} (\ell_s\eta_s)^{\beta_k}  \mathfrak{s}_{k,s}^\mathrm{av}(\bm x) \frac{\mathrm{d}s}{\eta_s} \prec   \varphi_{k,\tau}^\mathrm{av} (\ell_\tau\eta_\tau)^{\beta_k}  \mathfrak{s}_{k,\tau}^\mathrm{av}(\bm x),
	\end{equation}
	where in the last step we used \eqref{eq:sfunc_def}, \eqref{eq:elleta_monot}, and \eqref{eq:int_rules}. The other terms are estimated analogously to their counterparts in non-saturated case above. Therefore, \eqref{eq:av_masters} is established.
	
	Finally, the proof of the isotropic master inequality \eqref{eq:iso_masters} follows the exact same outline as the proof of \eqref{eq:av_masters} for the non-saturated chains, except we use the bound \eqref{eq:iso_prop_bound} to estimate the action of the propagators, \eqref{eq:iso_mart_bound} to bound the martingale term, and \eqref{eq:iso_forcing_bound} for the forcing term. We leave the straightforward details to the reader.   
	This concludes the proof of Proposition \ref{prop:masters} for the simplified case $z_{j, t} \in \{ z_t, \bar z_t\}$.

	Now we explain the minor modifications needed for the case of $k$ general spectral parameters $\bm z_t \in \mathbb{D}t^k$.
	Recall that $\eta_t$ was defined as the smallest among all $\eta{j,t}$, and the control quantities are monotone in $\eta$ (see~\eqref{eq:Ups_time_monot_notime}).
	Therefore, the $\Upsilon$ terms corresponding to the smallest $\eta_t$ indeed control all other $\Upsilon$ terms.
	
	Since $\eta_{j,t} \sim \eta_{i,t}$ for all $i, j \in \indset{\maxK}$, it follows that $\ell(\eta_{j,t}) \sim \ell(\eta_{i,t})$ as well. Thus, any estimate involving $\eta$ and $\ell$ factors remains valid when $\eta_{i,t}$ is replaced by its minimum, up to irrelevant constant factors.

	There is only one point in the previous proof where the special choice $z_{j,t} \in { z_t, \bar z_t }$ is explicitly used,
	namely, in the definition of the propagator 
	$$
	\mathcal{A}_{j,t} :=  \frac{m(z_{j,t})m(z_{j+1,t})S}{1- m(z_{j,t})m(z_{j+1,t})S}
	$$
	in~\eqref{eq:lin_prop_ops}. However, in the general case, $\mathcal{A}_{j,t}$ is either $\Theta(z_{j,t}, z_{j+1,t})$
	or $\Xi(z_{j,t}, z_{j+1,t})$. Inspecting the proof above, we see that the precise form of $\mathcal{A}_{j,t}$ as given in~\eqref{eq:lin_prop_ops} was not essential, and aside from the identities in~\eqref{eq:dTheta}, we relied only on general bounds. These bounds have already been formulated for the general setting in Definitions~\ref{def:adm_ups_notime} and~\ref{def:admS}, and the proof holds verbatim in the case of general spectral parameters. 
	This concludes  the proof of Proposition \ref{prop:masters} in the general case.
\end{proof}

\section{Observable regularization: Proof of Proposition \ref{prop:lin_term}} \label{sec:reg_sec}
In this section, we prove the estimate \eqref{eq:lin_term_est} on the linear term in \eqref{eq:k_av_evol} for saturated chains. This is a crucial ingredient for propagating the local laws along the zig flow
using Proposition \ref{prop:zig}, as it enables the closure of the hierarchy of master inequalities \eqref{eq:masters}. However, as we discussed previously in Section \ref{sec:masters_sec}, only the saturated chains (those satisfying \eqref{eq:sat_cahins}) of even length require special treatment. 
Hence, in the entire Section~\ref{sec:reg_sec} we assume that $k\in\indset{\maxK}$ is even, we fix, without loss of generality, $\bm z_t := \bm z_{t,\mathrm{sat}}$, as defined in \eqref{eq:alt_z}, and abbreviate
\begin{equation} \label{eq:Xk_shorthand}
	\mathcal{X}^{k}_t(\bm x) \equiv\mathcal{X}^{k}_{\bm z_{t,\mathrm{sat}}, t}(\bm x) := \Tr\bigl[(G - M)_{[1,k],t}(\bm x')S^{x_k}\bigr].
\end{equation}

The key idea for proving Proposition \ref{prop:lin_term} is observable regularization.
Recall the definition of the observable $S^y$ from \eqref{eq:S_obs}. Let $x,y$ be two indices in $\indset{N}$. We define  the regularization of $S^y$ with respect to the index $x$, denoted $\reg{S}^{x,y}$, as
\begin{equation} \label{eq:circ_above}
	\reg{S}^{x,y} := S^y - \delta_{xy}I = \sum_b \bigl(\delta_{by}-\delta_{xy}\bigr)S^b,
\end{equation}
where $I$ is an $N\times N$ identity matrix. In the sequel, we always regularize $S^{x_j}$ for $j \ge 2$ with respect to the preceding external index $x_{j-1}$. Hence, to simplify the notation where no ambiguity arises, we abbreviate
\begin{equation} \label{eq:Sring_def}
	\reg{S}^{x_j} \equiv \reg{S}^{x_j} (x_{j-1},x_j) := \reg{S}^{x_{j-1}, x_j} .
\end{equation} 

Owing to the cyclicity of trace,  the linear terms appearing on the left-hand side of \eqref{eq:lin_term_est} are structurally identical for all $j\in\indset{k}$, so we only treat the linear term corresponding to $j=k$,
\begin{equation} \label{eq:the_lin_term}
	\Theta_t^{(k)}\bigl[\mathcal{X}_t^k\bigr],
\end{equation}
to which we henceforth refer to as \emph{the linear term}, in full detail. All other linear terms are treated completely analogously, up to a cyclic shift of the vectors $\bm x$ and $\bm z_t$.

The general mechanism of observable regularization is used in the following way. Decomposing the observable $S^{x_k}$ into $\reg{S}^{x_k}$ and $\delta_{x_{k-1}x_k}I$ yields the identity
\begin{equation} \label{eq:G_circ_decomp}
	\Tr\bigl[(G-M)_{[1,k],t} S^{x_k} \bigr] = \Tr\bigl[(G-M)_{[1,k],t} \reg{S}^{x_k} \bigr] + \delta_{x_{k-1}x_k}\Tr\bigl[(G-M)_{[1,k],t}  \bigr].
\end{equation}
Now, the first term on the right-hand side of \eqref{eq:G_circ_decomp}, under the mollifying action of $\Theta_t$, will 
be smaller due to the presence of the regularized observable (see Propositions \ref{prop:circ_bound} and \ref{prop:2circ_bound} below); while the second term on the right-hand side of \eqref{eq:G_circ_decomp} can be reduced
 to a shorter chain of length $k-1$ by Ward identity, and subsequently estimated using the stopping time, without 
 jeopardising  the closability and the self-improving structure of the master inequalities. 

In the sequel, we treat short ($k\le K/2-1$) and long ($k\ge K/2$)
chains somewhat differently.  
To extract the mollifying effect of $\Theta_t$ on  the 
chain with a regularized $\reg{S}^{x_k}$, we derive separate evolution
equations for the mollified quantities, see~\eqref{eq:Z_elov_short}--\eqref{eq:Z_elov_long} later.
As it turns out (see discussion below \eqref{eq:noring}), mollification improves the propagation estimate by a factor $\sqrt{\eta_t/\eta_s}$,
which is sufficient to handle shorter chains. 
In fact the mollification of a single regularization compensates for the excess
propagator in~\eqref{eq:bad_av_prop}, which was the original reason for treating the linear term \eqref{eq:the_lin_term}
separately.

However, for long chains, the evolution equation for the mollified quantity---particularly bounding the quadratic variation of its martingale term---requires a reduction
inequality which incurs a loss. We account for this loss by increasing the loss
exponent:  from $\beta_k=0$ for $k\le K/2-1$  up to $\beta_K> 1/2$. This means, however, that 
that a suboptimal bound is fed into the Duhamel's formula. To compensate for this, we need an additional 
gain from further regularization.
It turns out that a second regularization, providing an additional  $\sqrt{\eta_t/\eta_s}$  factor, is 
already sufficient. Hence for longer chains 
we regularize two observables. 
Since the regularization $\reg{S}^{x_j}$ also involves
the previous index $x_{j-1}$, it is convenient to keep the two regularizations at least two sites apart in the chain to extract their effects independently. This means that only chains of length at least four 
can accommodate two regularizations, leading to the constraint $K/2\ge 4$ and thereby explaining the lower bound $K\ge 8$
on the maximal chain length in our proof.

\subsection{Short chains}
For chains of even length $k \in \indset{\maxK/2-1}$, it suffices to regularize a single observable, therefore, 
we define $\reg{\mathcal{X}}^{k}_t$ to be
\begin{equation} \label{eq:Gcirc}
	\reg{\mathcal{X}}^{k}_t \equiv \reg{\mathcal{X}}^{k}_t(\bm x) :=\Theta_t^{(k)}\biggl[\Tr\bigl[ (G - M)_{[1,k],t}(\bm x') \reg{S}^{x_k}\bigr]\biggr],
\end{equation}
where we recall the convention \eqref{eq:Sring_def}. Note that besides regularization of the observable, the definition of the quantity $\reg{\mathcal{X}}^{k}_t$ involves mollification by $\Theta_t$.
Using \eqref{eq:sumS=1}, it follows immediately that the linear term admits the decomposition
\begin{equation} \label{eq:short_lin_term_decomp}
	\Theta_t^{(k)} \bigl[ \mathcal{X}^{k}_t\bigr] = \reg{\mathcal{X}}^{k}_t + \mathcal{R}_{[1,k],t}, \quad k\in \indset{\maxK/2-1}, \quad k\text{--even},
\end{equation}
where, for even $k\in \indset{\maxK/2-1}$, the remainder term $\mathcal{R}_{[1,k],t}\equiv \mathcal{R}_{[1,k],t}(\bm x)$ is defined as
\begin{equation} \label{eq:remainder_short}
	\mathcal{R}_{[1,k],t}(\bm x) := (\Theta_t)_{x_{k-1}x_k}\mathcal{T}_{k}\bigl[\mathcal{X}^{k}_t\bigr](\bm x').
\end{equation} 
For the two terms in~\eqref{eq:short_lin_term_decomp} we have the following estimate:
\begin{prop}[Regularized Linear Term for Short Chains] \label{prop:circ_bound}
	 Assume that \eqref{eq:zig_init} holds,   and let $\tau$ be the  stopping time defined in \eqref{eq:tau_def}, then for all even $k \in \indset{\maxK/2-1}$ and $\bm x\in \indset{N}^k$, the mollified quantity $\reg{\mathcal{X}}^{k}_t(\bm x)$, defined in \eqref{eq:Gcirc}, admits the bound
	\begin{equation} \label{eq:circ_bound}
		\max_{\tinit \le s \le  t\wedge\tau} \frac{\eta_s\bigl\lvert \reg{\mathcal{X}}^{k}_s(\bm x) \bigr\rvert}{\mathfrak{s}^\mathrm{av}_{k,s}(\bm x)} \prec  \varphi_{k,t\wedge\tau}^\mathrm{av},
	\end{equation}
	uniformly for all $t \in [\tinit,T]$, where $\varphi_{k,t\wedge\tau}^\mathrm{av}$ is defined in \eqref{eq:phi_av_masters}.

	Moreover,  assuming only the weaker condition \eqref{eq:true_convol_notime} instead of \eqref{eq:convol_notime}--\eqref{eq:suppressed_convol_notime},   the remainder term  $\mathcal{R}_{[1,k],t}(\bm x)$,  defined in \eqref{eq:remainder_short} for all even integers $k \in \indset{\maxK/2-1}$, admits the bound
	\begin{equation} \label{eq:short_remainder_est}
		\max_{\tinit \le s \le  t\wedge\tau} \frac{\eta_s\bigl\lvert \mathcal{R}_{[1,k],s}(\bm x) \bigr\rvert}{(\ell_s\eta_s)^{\beta_k}\mathfrak{s}_{k,s}^\mathrm{av}(\bm x)} \lesssim \frac{\pav{k-1}}{(\ell_{t\wedge\tau}\eta_{t\wedge\tau})^{\beta_k - \beta_{k-1}}},
	\end{equation}
	uniformly in $t \in [\tinit,T]$.
\end{prop}
We prove Proposition \ref{prop:circ_bound} in Section \ref{sec:reg_evols}. The estimate of the more critical 
$\reg{\mathcal{X}}^{k}_t$ requires to write up a separate evolution equation and use Duhamel's principle.
The error term $\mathcal{R}$ will be estimated directly.

 Equipped with Proposition \ref{prop:circ_bound}, we prove Proposition \ref{prop:lin_term} for short chains ($k \le \maxK/2-1$).
\begin{proof}[Proof of Proposition \ref{prop:lin_term} for short chains] Assume that $k\in\indset{\maxK/2-1}$ is even. Then \eqref{eq:lin_term_est} follows immediately from \eqref{eq:short_lin_term_decomp}, \eqref{eq:phi_av_masters}, and Proposition \ref{prop:circ_bound}.
\end{proof}

\subsection{Long chains}
For chains of length $k \in [\maxK/2, \maxK]$, one regularization is not sufficient. Therefore, we regularize two observables
and we need to extract the improvements from both regularizations multiplicatively. 
However, since both regularizations need to be mollified in space by a $\Theta$, and only one $\Theta$ is 
present in \eqref{eq:lin_term_est}, we first consider the evolution equation of  
$\Theta_t^{(k)}\bigl[ \mathcal{X}^{k}_t \bigr]$. The linear term of this evolution equation will then contain two $\Theta$'s,
for which an improved estimate will be established.

It follows immediately from \eqref{eq:k_av_evol} and \eqref{eq:dTheta} that the linear term \eqref{eq:the_lin_term} satisfies the evolution equation
\begin{equation} \label{eq:long_lin_tem_evol}
	\begin{split}
		\mathrm{d}\biggl(\Theta_t^{(k)}\bigl[ \mathcal{X}^{k}_t \bigr]\biggr) =&~ \biggl(I+\Theta_t^{(k)}+\frac{k}{2}I + \sum_{i\neq k-2}\Theta_t^{(i)}\biggr)\circ \Theta_t^{(k)}\bigl[ \mathcal{X}^{k}_t \bigr]\mathrm{d}t\\ 
		&+ \Theta_t^{(k-2)}\circ \Theta_t^{(k)} \bigl[ \mathcal{X}^{k}_t \bigr]\mathrm{d}t + \Theta_t^{(k)}\biggl[\mathrm{d}\mathcal{M}^\mathrm{av}_{[1,k],t}  + \mathcal{F}^\mathrm{av}_{[1,k],t} \mathrm{d}t\biggr].
	\end{split}
\end{equation}
Note that out of the $k+1$ linear terms we separated the one with $\Theta_t^{(k-2)}$, as this term will be
estimated directly as a forcing term. In this way, the remaining $k$ linear terms can be treated
by Duhamel's formula.  We remark that more linear terms could have been estimated
in this way, but not all of them since along estimating $\Theta_t^{(k-2)}\circ \Theta_t^{(k)}$ it will be
important that the two indices $k-2$ and $k$ are separated from each other in the cyclic order.

To estimate this new forcing term, 
for $ k \in [\maxK/2, \maxK]$, we define the mollified quantity $\dring{\mathcal{X}}^{k}_t\equiv \dring{\mathcal{X}}^{k}_t(\bm x)$ to have two regularized observables at positions $k$ and $k-2$, that is
\begin{equation} \label{eq:G2circ}
	\dring{\mathcal{X}}^{k}_t(\bm x) :=\Theta_t^{(k-2)}\circ \Theta_t^{(k)}\biggl[\Tr\bigl[ (G - M)_{[1,k],t}\bigl(S^{x_1},\dots, S^{x_{k-3}}, \reg{S}^{x_{k-2}}, S^{x_{k-1}} \bigr) \reg{S}^{x_k}\bigr]\biggr],
	\end{equation}
where we recall the convention \eqref{eq:Sring_def}. Note that the observables we regularize are not neighboring, as this choice simplifies the subsequent analysis (see the proof of Lemma \ref{lemma:reg_props} in Section \ref{sec:props} below). Since $k \ge \maxK/2 \ge 4$ by assumption, such choice is possible.

The definition of the regularized observables \eqref{eq:circ_above} implies the identity
\begin{equation} \label{eq:long_lin_term_decomp}
	\Theta_t^{(k-2)}\circ \Theta_t^{(k)}\bigl[ \mathcal{X}^{k}_t\bigr] = \dring{\mathcal{X}}^{k}_t + \mathcal{R}_{[1,k],t}, \quad \maxK/2 \le k \le \maxK,
\end{equation}
where, for $k \in [\maxK/2,\maxK]$, the remainder term $\mathcal{R}_{[1,k],t}\equiv \mathcal{R}_{[1,k],t}(\bm x)$ is given by
\begin{equation} \label{eq:remainder_long}
	\begin{split}
		\mathcal{R}_{[1,k],t}(\bm x) :=&~  (\Theta_t)_{x_{k-1}x_k} \mathcal{T}_{k}\circ \Theta_t^{(k-2)} \bigl[\mathcal{X}^{k}_t\bigr](\bm x') \\
		&+ (\Theta_t)_{x_{k-2}x_{k-3}} \mathcal{T}_{k-2}\circ \Theta_t^{(k)} \bigl[\mathcal{X}^{k}_t\bigr](\bm x''',x_{k-1},x_k)\\ 
		&- (\Theta_t)_{x_{k-1}x_k}(\Theta_t)_{x_{k-2}x_{k-3}} \mathcal{T}_{k-2, k}\bigl[\mathcal{X}^{k}_t\bigr](\bm x''', x_{k-1}),
	\end{split}
\end{equation}
where we recall that $\bm x''' := (x_1,\dots, x_{k-3}) \in \indset{N}^{k-3}$, and the operators $(\Theta_t)^{(j)}$ act on  the index $x_j$.
For these two terms in~\eqref{eq:long_lin_term_decomp} we have the following estimate.
\begin{prop}[Regularized Linear Term for Long Chains] \label{prop:2circ_bound}
	Assume that \eqref{eq:zig_init} holds,   and let $\tau$ be the  stopping time defined in \eqref{eq:tau_def}, then the mollified quantity $\dring{\mathcal{X}}^{k}_t(\bm x)$ defined in \eqref{eq:G2circ} admits the bound
	\begin{equation} \label{eq:two_ring_ll}
		\max_{\tinit \le s \le  t\wedge\tau} \frac{\eta_s^2\bigl\lvert \dring{\mathcal{X}}^{k}_s(\bm x) \bigr\rvert}{(\ell_s\eta_s)^{\beta_k} \mathfrak{s}^\mathrm{av}_{k,s}(\bm x)} 
		\prec  \varphi_{k,t\wedge\tau}^\mathrm{av},
	\end{equation}
	uniformly in $t\in[\tinit,T]$, where $\varphi_{k,t\wedge\tau}^\mathrm{av}$ is defined in \eqref{eq:phi_av_masters}.
	
	Moreover,  assuming only the weaker condition \eqref{eq:true_convol_notime} instead of \eqref{eq:convol_notime}--\eqref{eq:suppressed_convol_notime},   the remainder term  $\mathcal{R}_{[1,k],t}(\bm x)$,  defined in \eqref{eq:remainder_long} for even $k \in \indset{\maxK/2, \maxK}$, admits the bound
	\begin{equation} \label{eq:remainder_est}
		\max_{\tinit \le s \le  t\wedge\tau} \frac{\eta_s^2\bigl\lvert \mathcal{R}_{[1,k],s}(\bm x) \bigr\rvert}{(\ell_s\eta_s)^{\beta_k}\mathfrak{s}_{k,s}^\mathrm{av}(\bm x)} \lesssim \frac{\pav{k-1}}{(\ell_{t\wedge\tau}\eta_{t\wedge\tau})^{\beta_k - \beta_{k-1}}} + \frac{\pav{k-2}}{(\ell_{t\wedge\tau}\eta_{t\wedge\tau})^{\beta_k - \beta_{k-2}}},
	\end{equation}
	uniformly in $t\in[\tinit,T]$.
\end{prop}
We prove Proposition \ref{prop:2circ_bound} in Section \ref{sec:reg_evols}.
Similarly to the proof of Proposition~\ref{prop:circ_bound} for short chains, the more critical 
$\dring{\mathcal{X}}^{k}_s$ term will be estimated by writing up its own evolution equation 
and using Duhamel's formula. The simpler $\mathcal{R}$ will be estimated directly.

\begin{proof} [Proof of Proposition \ref{prop:lin_term} for long chains]	
	Assume that $k \in \indset{\maxK/2,\maxK}$. 
	It follows from \eqref{eq:long_lin_term_decomp}, Proposition \ref{prop:2circ_bound} and \eqref{eq:remainder_est} that, uniformly for all $t \in [\tinit,T]$,
	\begin{equation} \label{eq:lin_term_ll}
		\max_{\tinit\le s\le t\wedge\tau} \frac{\eta_s^2\bigl\lvert \Theta_s^{(k-2)}\circ \Theta_s^{(k)} \bigl[ \mathcal{X}^{k}_s\bigr]\bigr\rvert}{(\ell_s\eta_s)^{\beta_k}\mathfrak{s}_{k,s}^{\mathrm{av}}(\bm x)} \prec
	\varphi_{k,t\wedge\tau}^\mathrm{av}.
	\end{equation}

	Similarly to the proof of \eqref{eq:av_masters}, we apply Duhamel's principle to \eqref{eq:long_lin_tem_evol},
	considering the first line of   \eqref{eq:long_lin_tem_evol} as the linear term,
	 while treating $\Theta_s^{(k-2)}\circ \Theta_s^{(k)} \bigl[ \mathcal{X}^{k}_s\bigr]$ as an additional forcing term; 
	  we obtain
	 \begin{equation} \label{eq:long_lin_tem_solve}
	 	\begin{split}
	 		\Theta_{t\wedge\tau}^{(k)}\bigl[ \mathcal{X}^{k}_{t\wedge\tau} \bigr] =&~ \Theta_{t\wedge\tau}^{(k)}\Biggl[ \mathrm{e}^{\frac{2-k}{2}(\tinit-t\wedge\tau)} \other{\mathcal{P}}_{\tinit,t\wedge\tau}^{k} \bigl[ \mathcal{X}^{k}_{\tinit} \bigr] 
	 		+\int\limits_{\tinit}^{t\wedge\tau} \mathrm{e}^{\frac{2-k}{2}(s-t\wedge\tau)}\other{\mathcal{P}}^k_{s,t\wedge\tau} \biggl[\mathrm{d}\mathcal{M}^\mathrm{av}_{[1,k],s}  + \mathcal{F}^\mathrm{av}_{[1,k],s} \mathrm{d}s\biggr]\Biggl]\\
	 		&+\int\limits_{\tinit}^{t\wedge\tau} \biggl(\mathrm{e}^{\frac{2-k}{2}(s-t\wedge\tau)}\mathcal{P}_{s,t\wedge\tau}^{(k)}\circ \other{\mathcal{P}}_{s,t\wedge\tau}^{k}\biggr)\biggl[\Theta_s^{(k-2)}\circ \Theta_s^{(k)} \bigl[ \mathcal{X}^{k}_s\bigr]\biggr]\mathrm{d}s,
	 	\end{split}
	 \end{equation}
	 where $\other{\mathcal{P}}^k_{s,t} := \bigotimes_{i\neq k-2} \mathcal{P}_{s,t}^{(i)}$, and we used the identity $\mathcal{P}_{s,t}^{(k)}\circ\Theta_{s}^{(k)} = \Theta_{t}^{(k)}$. Note that the propagator $\other{\mathcal{P}}^k_{s,t}$ belongs to the class \eqref{eq:general_prop} with one of its components, namely $\other{\mathcal{P}}_{s,t}^{(k-2)} = I$ being non-saturated. In particular, $\other{\mathcal{P}}^k_{s,t}$ satisfies the conditions of \eqref{eq:av_prop_bound}. Hence using Lemmas~\ref{lemma:good_props}, \ref{lemma:mart_est}, \ref{lemma:forcing}, and \eqref{eq:zig_init}, we conclude that the contribution of terms in the first line of \eqref{eq:long_lin_tem_solve} is stochastically dominated by
	\begin{equation} \label{eq:line1_law}
		\varphi_{k,t\wedge\tau}^\mathrm{av}(\ell_{t\wedge\tau}\eta_{t\wedge\tau})^{\beta_k}\Theta_{t\wedge\tau}^{(k)}\bigl[\mathfrak{s}_{k,t\wedge\tau}^\mathrm{av}(\bm x)\bigr] \lesssim  \varphi_{k,t\wedge\tau}^\mathrm{av}\eta_{t\wedge\tau}^{-1}(\ell_{t\wedge\tau}\eta_{t\wedge\tau})^{\beta_k} \mathfrak{s}_{k,t\wedge\tau}^\mathrm{av}(\bm x).
	\end{equation}
	On the other hand, it follows from \eqref{eq:av_prop_bound} and \eqref{eq:lin_term_ll} that the contribution of terms in the second line of \eqref{eq:long_lin_tem_solve} is stochastically dominated by
	\begin{equation} \label{eq:line2_law}
		\int\limits_{\tinit}^{t\wedge\tau}\varphi_{k,s}^\mathrm{av} \frac{(\ell_s\eta_s)^{\beta_k}}{\eta_s^2}\frac{\ell_{t\wedge\tau}\eta_{t\wedge\tau}}{\ell_s\eta_s} \mathcal{P}_{s,t\wedge\tau}^{(k)} \bigl[\mathfrak{s}_{k,t\wedge\tau}^\mathrm{av}(\bm x)\bigr]\mathrm{d}s 
		\lesssim \varphi_{k,t\wedge\tau}^\mathrm{av}\eta_{t\wedge\tau}^{-1}(\ell_{t\wedge\tau}\eta_{t\wedge\tau})^{\beta_k} \mathfrak{s}_{k,t\wedge\tau}^\mathrm{av}(\bm x).
	\end{equation} 
	 Here we used the integration rules \eqref{eq:int_rules} and the monotonicity of $\ell_s\eta_s$ in $s$. Combining the estimates \eqref{eq:line1_law} and \eqref{eq:line2_law} on the first and second lines of
	 \eqref{eq:long_lin_tem_solve}, we obtain the desired \eqref{eq:lin_term_est}.   
	This concludes the proof of Proposition \ref{prop:lin_term}.  
\end{proof}

\subsection{Equations for chains with regularized observables: Proof of Proposition \ref{prop:circ_bound} and \ref{prop:2circ_bound}} \label{sec:reg_evols}
In this section we analyze the evolution equations for the mollified quantities $\reg{\mathcal{X}}^k_{t}$ and $\dring{\mathcal{X}}^k_{t}$, defined in \eqref{eq:Gcirc} and \eqref{eq:G2circ}, respectively. The intrinsic regularizations in $\reg{\mathcal{X}}^k_{t}$ and $\dring{\mathcal{X}}^k_{t}$ modify their evolution equations compared to that of their unregularized counterpart $\mathcal{X}^k_{t}$ in \eqref{eq:k_av_evol}.  To express these equations efficiently, we introduce the linear operators $\reg{\Theta}_t^{(j)}$. 
For all $j \in \indset{2,k}$, let $\reg{\Theta}_t^{(j)}$ denote a linear operator acting on functions $f(\bm x)$ of $\bm x \in \indset{N}^k$ by 
\begin{equation} \label{eq:Thetaring}
	\reg{\Theta}_t^{(j)}\bigl[f\bigr](\bm x) := \sum_{\bm a \in \indset{N}^k} \prod_{i \neq j} \delta_{x_ia_i} \bigl(\reg{\Theta}_t^{x_{j-1}}\bigr)_{x_{j}a_j} f(\bm a), \quad \bigl(\reg{\Theta}_t^{x}\bigr)_{ab} :=
	 (\Theta_t)_{ab}-(\Theta_t)_{ax}, \quad a,b,x \in \indset{N}.
\end{equation}
We will call $\reg{\Theta}_t^x$ the {\it regularized operator} or {\it regularized version} of $\Theta_t$ with respect to $x$.  While the circle notation 
is similar to the regularization $\reg{S}^{x,y}$ of the observable $S^y$ defined in \eqref{eq:circ_above},
these two definitions are different. The observable regularization just subtracted an identity matrix 
while the regularization of $\Theta_t$ subtracts another $\Theta_t$ anchored at the fix index $x$.
Note that the resulting matrix $\reg{\Theta}_t$ is non-symmetric, but this fact  will play no  role. However, the two
regularizations are closely related since $\reg{\Theta}_t$ naturally arises in the equation 
for chains with regularized observable as the following lemma shows.

Before proceeding with the proof of Proposition \ref{prop:circ_bound} and \ref{prop:2circ_bound}, we collect the necessary technical Lemmas \ref{lemma:circ_evol}--\ref{lemma:M_Ward}, that will be proven in later sections.

\begin{lemma} [Equation for Chains with Regularized Observables] \label{lemma:circ_evol}
	Let  $\tau \in [\tinit, T]$ be a stopping time, and let $k \in \indset{\maxK}$ be an even integer. Let the quantity $\mathcal{Z}^k_t  := \reg{\mathcal{X}}_{t}^{k}$, defined in \eqref{eq:Gcirc}, if $k\in \indset{\maxK/2-1}$ and $\mathcal{Z}^k_t  := \dring{\mathcal{X}}_{t}^{k}$, defined in \eqref{eq:G2circ}, if $k\in \indset{\maxK/2,\maxK}$, then
	\begin{equation} \label{eq:Gcirc_solve}
		\begin{split}
			\mathcal{Z}_{t\wedge\tau}^{k} =&~ \mathrm{e}^{k(\tinit-t\wedge\tau)/2}\mathcal{P}^{\otimes k}_{\tinit,t\wedge\tau}\circ \mathcal{U}^{(k)}_{t\wedge\tau}  \bigl[ 
			\mathcal{X}_{\tinit}^{k,\mathrm{av}}
			\bigr] 
			+\int_{\tinit}^{t\wedge\tau} \mathrm{e}^{k(s-t\wedge\tau)/2}\mathcal{P}^{\otimes k}_{s,t\wedge\tau}\circ \mathcal{U}^{(k)}_{t\wedge\tau}  \bigl[ \mathrm{d}\mathcal{M}^\mathrm{av}_{[1,k],s}  \bigr]\\
			&+ \int_{\tinit}^{t\wedge\tau} \mathrm{e}^{k(s-t\wedge\tau)/2}\mathcal{P}^{\otimes k}_{s,t\wedge\tau} \circ \mathcal{U}^{(k)}_{t\wedge\tau} \bigl[ \mathcal{F}^\mathrm{av}_{[1,k],s}+ \other{\mathcal{F}}^\mathrm{av}_{[1,k],s} \bigr]  \mathrm{d}s,
		\end{split}
	\end{equation}
	where the linear operator $\mathcal{U}_{t}^{(k)}$ is given by
	\begin{equation} \label{eq:U_operator}
		\mathcal{U}^{(k)}_{t} := \begin{cases}
			\reg{\Theta}_{t}^{(k)}, \quad &k \in \indset{\maxK/2-1},\\
			\reg{\Theta}_{t}^{(k-2)} \circ  \reg{\Theta}_{t}^{(k)}, \quad &k \in \indset{\maxK/2,\maxK},\\
		\end{cases}
	\end{equation}	
	and $\other{\mathcal{F}}^\mathrm{av}_{[1,k],t}\equiv \other{\mathcal{F}}^\mathrm{av}_{[1,k],t}(\bm x)$ is an additional forcing term, arising from $\delta_{xy}I$ in \eqref{eq:circ_above}, the definition of which depends on $k$. For $k\in\indset{\maxK/2-1}$, the additional forcing term $\other{\mathcal{F}}^\mathrm{av}_{[1,k],t}(\bm x)$ is given by
	\begin{equation} \label{eq:extra_forcing1}
		\other{\mathcal{F}}^\mathrm{av}_{[1,k],t}(\bm x) := (\Theta_t)_{x_{k-1}x_k}  \bigl(\mathcal{T}_{k}\bigl[\mathcal{X}^{k}_t\bigr](\bm x'',x_k)  +  \mathcal{T}_{k}\bigl[\mathcal{X}^{k}_t\bigr](\bm x')\bigr).
	\end{equation}
	while for $k \in \indset{\maxK/2,\maxK}$, it is defined as
	\begin{equation} \label{eq:extra_forcing2}
		\begin{split}
			\other{\mathcal{F}}^\mathrm{av}_{[1,k],t}(\bm x) :=&~  
			(\Theta_t)_{x_{k-3} x_{k-2}} \bigl(\mathcal{T}_{k-2} \bigl[\mathcal{X}^{k}_t\bigr](\bm x^{(4)}, x_{k-2}, x_{k-1}, x_k)  + \mathcal{T}_{k-2}\bigl[\mathcal{X}^{k}_t\bigr](\bm x''', x_{k-1},x_k)\bigr)\\
			& +(\Theta_t)_{x_{k-1}x_k} \bigl(\mathcal{T}_{k}\bigl[\mathcal{X}^{k}_t\bigr](\bm x'', x_k)  +  \mathcal{T}_{k}\bigl[\mathcal{X}^{k}_t\bigr](\bm x')\bigr).
		\end{split}
	\end{equation}
	Here we recall the conventions \eqref{eq:xprime} and \eqref{eq:xprimes}.
\end{lemma}
We prove Lemma \ref{lemma:circ_evol} in Section \ref{sec:evols}.

The following lemma provides the crucial estimate for the composition of the saturated propagator $\mathcal{P}_{s,t}$ with the regularized operator $\reg{\Theta}_t^{(j)}$.
\begin{lemma} [Regularized Propagator Estimates] \label{lemma:reg_props}
	Let $s,t \in [0,T]$ be a pair of times satisfying $s\le t$, and let $k\in\indset{\maxK}$ be even.
	Then, for any function $f(\bm x)$ of $\bm x \in \indset{N}^{k}$, and $\varphi > 0$,
	\begin{equation} \label{eq:one_ring}
		\forall \bm x~~\bigl\lvert f(\bm x) \bigr\rvert \lesssim \varphi \, \mathfrak{s}_{k,s}^\mathrm{av}(\bm x) \quad \Longrightarrow \quad 
		\forall \bm x~~ \biggl\lvert \mathcal{P}_{s,t}^{\otimes k}\circ \reg{\Theta}_t^{(k)} \bigl[f \bigr]  (\bm x)\biggr\rvert \lesssim \varphi \, \frac{1}{\eta_t}  \, \mathfrak{s}_{k,t}^\mathrm{av}(\bm x).
	\end{equation}
	
	Moreover, if $k \ge 4$, then for any function $f(\bm x)$ of $\bm x \in \indset{N}^{k}$, and $\varphi > 0$,
	\begin{equation} \label{eq:two_rings}
		\forall \bm x~~\bigl\lvert f(\bm x) \bigr\rvert \lesssim \varphi \, \mathfrak{s}_{k,s}^\mathrm{av}(\bm x) \quad \Longrightarrow \quad 
		\forall \bm x~~ \biggl\lvert \mathcal{P}_{s,t}^{\otimes k}\circ \reg{\Theta}_t^{(k-2)}\circ \reg{\Theta}_t^{(k)} \bigl[f \bigr]  (\bm x)\biggr\rvert \lesssim \varphi \, \frac{1}{\eta_t^2}\frac{\ell_t\eta_t}{\ell_s\eta_s}  \, \mathfrak{s}_{k,t}^\mathrm{av}(\bm x),
	\end{equation}
	where $\reg{\Theta}_t^{(j)}$ are defined in \eqref{eq:Thetaring}. Furthermore, estimates \eqref{eq:one_ring} and \eqref{eq:two_rings} hold with $\reg{\Theta}_t$ replaced by its entry-wise absolute value  $|\reg{\Theta}_t|$.  
\end{lemma}
We prove Lemma \ref{lemma:reg_props} in Section \ref{sec:props}. We note that if all $\reg{\Theta}_t^{(j)}$ in \eqref{eq:one_ring} and \eqref{eq:two_rings}
were replaced with the corresponding $\Theta_t^{(j)}$ (without regularization), the estimates would degrade by a factor of $\ell_t/\ell_s$ and $\eta_s/\eta_t$, respectively, to
\begin{equation}\label{eq:noring}
\biggl\lvert \mathcal{P}_{s,t}^{\otimes k}\circ {\Theta}_t^{(k)} \bigl[f \bigr]  (\bm x)\biggr\rvert \lesssim \varphi \, \frac{1}{\eta_t}\, \frac{\ell_t}{\ell_s}  \, \mathfrak{s}_{k,t}^\mathrm{av}(\bm x),
\qquad\biggl\lvert \mathcal{P}_{s,t}^{\otimes k}\circ {\Theta}_t^{(k-2)}\circ{\Theta}_t^{(k)} \bigl[f \bigr]  (\bm x)\biggr\rvert \lesssim \varphi \, \frac{1}{\eta_t^2}\frac{\ell_t}{\ell_s}  \, \mathfrak{s}_{k,t}^\mathrm{av}(\bm x),
\end{equation}
 see \eqref{eq:bad_av_prop}. More precisely,  along the proof of Lemma \ref{lemma:reg_props} 
 we will show that each change of a ${\Theta}_t^{(j)}$ to
 $\reg{\Theta}_t^{(j)}$ improves the bound by a factor $\ell_s\ell_t\eta_t/W^2$, see 
 \eqref{eq:Thetaring_action_basic}--\eqref{eq:Thetaring_action_basic1}
 later. 
 The individual improvement factor satisfies $\ell_s\ell_t\eta_t/W^2 \sim \ell_s/\ell_t\sim
 \sqrt{\eta_t/\eta_s}$ when $\eta_t \ge (W/N)^2$, thus it exactly compensates the $\ell_t/\ell_s$ factor
 owing to the excess propagator
 in the saturated propagator estimate~\eqref{eq:bad_av_prop}.  In the complementary  $\eta_t \ll (W/N)^2 $ regime 
 the improvement factor is, in fact, much smaller than $ \ell_s/\ell_t$ (since $\ell_t^2\eta_t/W^2 \ll 1$). 

The following lemma, proven in Section \ref{sec:forcing}, provides an estimate for the additional forcing term in \eqref{eq:Gcirc_solve}.
\begin{lemma} [Estimates on Additional Forcing Terms] \label{lemma:extra_forcing}
	 Assume only the weaker condition \eqref{eq:true_convol_notime} instead of \eqref{eq:convol_notime}--\eqref{eq:suppressed_convol_notime}.  
	Let $\tau$ be the stopping time defined in \eqref{eq:tau_def}, and let $s \in [0,\tau]$. 
	The quantity  $\other{\mathcal{F}}^\mathrm{av}_{[1,k],s}(\bm x)$, defined in \eqref{eq:extra_forcing1} for all even integers $k \in \indset{\maxK/2-1}$, and in \eqref{eq:extra_forcing2} for even $k \in \indset{\maxK/2, \maxK}$, admits the bound
	\begin{equation} \label{eq:extra_forcing}
		\frac{\eta_s\bigl\lvert \other{\mathcal{F}}^\mathrm{av}_{[1,k],s}(\bm x) \bigr\rvert}{(\ell_s\eta_s)^{\beta_k}\mathfrak{s}_{k,s}^\mathrm{av}(\bm x)} \lesssim \frac{\pav{k-1}}{(\ell_s\eta_s)^{\beta_k - \beta_{k-1}}}, \quad k \in \indset{\maxK}~\text{-- even}.
	\end{equation}
\end{lemma}

\begin{lemma}[Divided Difference Identity for $M$] \label{lemma:M_Ward}
	For any integer length $k \ge 2$ and any index $j \in \indset{k-1}$, let $M_{[1,k]}\rvert_{A_j=I} := M(z_1, A_1, \dots, z_{j-1}, I, z_{j}, A_j, z_{j+1}, \dots, A_{k-1}, z_k)$ denote the matrix $M_{[1,k]}$ with $A_j :=I$.
	Then $M_{[1,k]}\rvert_{A_j=I}$   satisfies the divided-difference identity, 
	\begin{equation} \label{eq:M_Ward}
		\begin{split}
			M_{[1,k]}\bigr\rvert_{A_j=I} = \frac{1}{z_{j}-z_{j+1}}\biggl(~&M(z_1, A_1, \dots, z_{j-1},A_{j-1}, z_{j},A_{j+1},z_{j+2},\dots, A_{k-1},z_k) \\
			&- M(z_1,A_1,\dots z_{j-1},A_{j-1},z_{j+1}, A_{j+1}, z_{j+2}, \dots, A_{k-1},z_k) \biggr),
		\end{split}
	\end{equation}
	assuming $z_j \neq z_{j+1}$.
\end{lemma}
We prove Lemma \ref{lemma:M_Ward} in Section \ref{sec:M_rec_analysis}.  Before proceeding with the proof of Proposition \ref{prop:circ_bound}, we record two important properties of the time-dependent control function $\Upsilon_t$, that follow immediately from Definition \ref{def:adm_ups_notime}, and will be used repeatedly in the sequel. First, it follows from \eqref{eq:Ups_majorates_notime} that
\begin{equation} \label{eq:Ups_majorates}
	(\Theta_t)_{xy} \le C (\Upsilon_t)_{xy}, \quad \bigl\lvert(\Xi_t)_{xy}\bigr\rvert \le C (\Upsilon_{0})_{xy} \quad x,y\in\indset{N}, \quad 0 \le t \le T.
\end{equation}
Next, \eqref{eq:triag_notime} implies that 
\begin{equation} \label{eq:triag}
	\max_a (\Upsilon_s)_{xa} (\Upsilon_t)_{ay} \le C (\ell_s\eta_s)^{-1}(\Upsilon_t)_{xy},  \quad x,y\in\indset{N}, \quad 0 \le s \le t \le T.
\end{equation}
Moreover, \eqref{eq:true_convol_notime} implies that
\begin{equation} \label{eq:true_convol}
	\sum_a (\Upsilon_s)_{xa} (\Upsilon_t)_{ay} \le C\frac{1}{\eta_s} (\Upsilon_t)_{xy}, \quad x,y\in\indset{N}, \quad 0 \le s \le t \le T.
\end{equation}
Finally, combining  Schwarz inequality and \eqref{eq:true_convol}, we obtain, for all $0\le s \le t \le T$,
\begin{equation} \label{eq:Schwarz_convol}
	\sum_q \sqrt{(\Upsilon_s)_{x_1q}(\Upsilon_s)_{x_2q} (\Upsilon_t)_{qy_1} (\Upsilon_t)_{qy_2}} 
	\le \frac{1}{\eta_s}C  \sqrt{ (\Upsilon_t)_{x_1y_1} (\Upsilon_t)_{x_2y_2}}, \quad x_1,x_2,y_1,y_2\in\indset{N}.
\end{equation}

\begin{proof} [Proof of Proposition \ref{prop:circ_bound}]
	First, we prove \eqref{eq:circ_bound}.
	It follows directly from equation \eqref{eq:Gcirc_solve} and Lemmas \ref{lemma:mart_est}, \ref{lemma:forcing}, and \ref{lemma:extra_forcing}, together with the bound \eqref{eq:one_ring} that, for all $t\in[\tinit,T]$,
	\begin{equation}
		\max_{0 \le s \le t\wedge\tau} \frac{\eta_s\bigl\lvert \reg{\mathcal{X}}^{k}_s(\bm x) \bigr\rvert}{\mathfrak{s}^\mathrm{av}_{k,t}(\bm x)} \prec  \varphi_{k,t\wedge\tau}^\mathrm{av},
	\end{equation}
	where we used the fact that $\beta_k = 0$ for $k \in \indset{\maxK/2-1}$ and $\alpha_j =0 $ for $j \in \indset{\maxK/2}$.
	
	Next, we prove \eqref{eq:short_remainder_est}. 
	Recall from \eqref{eq:alt_z} that $z_{1,s} = z_s$ and $z_{k,s} = \overline{z}_s$. 
	Ward identity implies that 
	\begin{equation} \label{eq:Gk_Ward}
		\sum_{a} \Tr\bigl[G_{[1,k],s} S^a\bigr] = \Tr\bigl[G_{[1,k],s}\bigr] = \frac{1}{2\mathrm{i}\eta_s}\biggl(\Tr\bigl[G_{s} S^{x_1}G_{[2,k-1],s}\bigr] - \Tr\bigl[G_s^* S^{x_1}G_{[2,k-1],s}\bigr]\biggr),		
	\end{equation} 
	where $G_s := G_s(z_s)$.
	Therefore, combining the definition of $\reg{\mathcal{X}}^{k}_s$ from \eqref{eq:Gcirc}, \eqref{eq:Gk_Ward} and the divided difference identity \eqref{eq:M_Ward} for the deterministic approximation terms $M_{[1,k],s}$, we deduce from \eqref{eq:Psi_def}, and \eqref{eq:tau_def} that  
	\begin{equation} \label{eq:ward_local_law}
		\bigl\lvert \mathcal{T}_{k}\bigl[\mathcal{X}^{k}_s\bigr](\bm x') \bigr\rvert \lesssim \eta_s^{-1}(\ell_s\eta_s)^{\beta_{k-1}} \mathfrak{s}_{k-1,s}^\mathrm{av}(\bm x') \Psi_{k-1,s}^\mathrm{av} \lesssim \eta_s^{-1}(\ell_s\eta_s)^{\beta_{k-1}} \pav{k-1} \mathfrak{s}_{k-1,s}^\mathrm{av}(\bm x'),
	\end{equation}	
	uniformly in $s \in [\tinit,\tau]$.
	In particular, for $k=2$, \eqref{eq:extra_forcing} follows immediately. On the other hand, for $k\ge 4$, we conclude \eqref{eq:short_remainder_est} from \eqref{eq:ward_local_law} by triangle inequality \eqref{eq:triag}. 	
	This concludes the proof of Proposition \ref{prop:circ_bound}.
\end{proof}

\begin{proof}[Proof of Proposition \ref{prop:2circ_bound}]
	First, we prove \eqref{eq:two_ring_ll}. 
	It follows from Lemma \ref{lemma:mart_est} and the bound \eqref{eq:two_rings} that
	\begin{equation} \label{eq:2ring_mart_bound}
		\begin{split}
			\int_{\tinit}^{t\wedge\tau} \mathcal{P}_{s,t\wedge\tau}^{\otimes k}\circ\reg{\Theta}_{t\wedge\tau}^{(k-2)} \circ  \reg{\Theta}_{t\wedge\tau}^{(k)} \bigl[\mathrm{d}\mathcal{M}^\mathrm{av}_{[1,k],s}\bigr](\bm x) 
			&\prec \biggl(\int_{\tinit}^{t\wedge\tau} \bigl(\varphi_{k,s}^\mathrm{av,qv}\bigr)^2  \frac{1}{\eta_{t\wedge\tau}^4}\frac{(\ell_{t\wedge\tau}\eta_{t\wedge\tau})^2}{(\ell_s\eta_s)^{2-2\beta_k}} 
			 \frac{\mathrm{d}s}{\eta_s}  \biggr)^{1/2}  \mathfrak{s}_{k,t}^\mathrm{av}(\bm x)\\
			&\prec \frac{(\ell_t\eta_t)^{\beta_k}}{\eta_{t\wedge\tau}^2} \varphi_{k,t\wedge\tau}^\mathrm{av,qv} \,\mathfrak{s}_{k,t}^\mathrm{av}(\bm x),
		\end{split}
	\end{equation}
	where in the last step we used the fact that $\beta_k < 1$, and $\varphi_{k,s}^\mathrm{av,qv}$ is monotone non-decreasing in $s\in[0,T]$.
	Hence, \eqref{eq:two_ring_ll} follows directly from the evolution equation \eqref{eq:Gcirc_solve} and Lemmas \ref{lemma:forcing}, and \ref{lemma:extra_forcing}, together with the bounds \eqref{eq:two_rings} and \eqref{eq:2ring_mart_bound}. 
	
	Finally, we prove \eqref{eq:remainder_est}. We estimate the terms on the right-hand side of \eqref{eq:remainder_long} one by one. For the first term, it follows from \eqref{eq:Ups_majorates} and  \eqref{eq:ward_local_law}, that
	\begin{equation} \label{eq:long_rem_term1}
		\begin{split}
			\biggl\lvert\frac{\eta_t^2(\Theta_t)_{x_{k-1}x_k}}{(\ell_t\eta_t)^{\beta_{k-1}}}  \mathcal{T}_{k}\circ \Theta_t^{(k-2)} \bigl[\mathcal{X}^{k}_t\bigr](\bm x')\biggr\rvert 
			&\lesssim \frac{\eta_t^2(\Upsilon_t)_{x_{k-1}x_k}}{(\ell_t\eta_t)^{\beta_{k-1}}}\sum_{a}(\Upsilon_t)_{x_{k-2}a}\bigl\lvert \mathcal{T}_{k} \bigl[\mathcal{X}^{k}_t\bigr](\bm x''',a,x_{k-1})\bigr\rvert\\
			&\lesssim \frac{\pav{k-1}\, \mathfrak{s}_{k-1,t}^\mathrm{av}(\bm x') (\Upsilon_t)_{x_{k-1}x_k} }{(\ell_{t\wedge\tau}\eta_{t\wedge\tau})^{\beta_k - \beta_{k-1}}}
			\lesssim \frac{\pav{k-1}}{(\ell_{t\wedge\tau}\eta_{t\wedge\tau})^{\beta_k - \beta_{k-1}}}\,\mathfrak{s}_{k,t}^\mathrm{av}(\bm x),
		\end{split}
	\end{equation}
	where in the second inequality we used  \eqref{eq:triag}, Schwarz inequality together with \eqref{eq:Schwarz_convol}, \eqref{eq:sfunc_def}, and in the last step we used \eqref{eq:triag}. 
	The term in the second line of \eqref{eq:remainder_long} admits the same bound. 
	
	For the third and final term, similarly to \eqref{eq:Gk_Ward}, we combine two Ward identities for the resolvent $G_t$ and use Lemma \ref{lemma:M_Ward} twice to deduce that
	\begin{equation}
		\bigl\lvert \mathcal{T}_{k-2, k}\bigl[\mathcal{X}^{k}_t\bigr](\bm x''', x_{k-1})\bigr\rvert \lesssim \frac{1}{\eta_t^2}(\ell_t\eta_t)^{\beta_{k-2}}\pav{k-2}\mathfrak{s}_{k-2,t}^\mathrm{av}(\bm x''', x_{k-1}).
	\end{equation}
	Hence, by triangle inequality \eqref{eq:triag}, we obtain the bound
	\begin{equation} \label{eq:long_rem_term3}
		\bigl\lvert \eta_t^2(\Theta_t)_{x_{k-1}x_k}(\Theta_t)_{x_{k-2}x_{k-3}} \mathcal{T}_{k-2, k}\bigl[\mathcal{X}^{k}_t\bigr](\bm x''', x_{k-1})\bigr\rvert \lesssim \pav{k-2}(\ell_t\eta_t)^{\beta_{k-2}}\mathfrak{s}_{k,t}^\mathrm{av}(\bm x).
	\end{equation}
	Therefore, estimates \eqref{eq:long_rem_term1} and \eqref{eq:long_rem_term3} imply \eqref{eq:remainder_est} for $k\in\indset{\maxK/2,\maxK}$. 
	This concludes the proof of Proposition \ref{prop:2circ_bound}.
\end{proof}

\section{Proof of technical lemmas from Sections \ref{sec:masters_sec} and \ref{sec:reg_sec}}\label{sec:aux}

We record the following properties of $\Upsilon_t$ and $\Upsilon_s$ for all $0 \le s \le t \le T$, that follow immediately from Definition~\ref{def:adm_ups_notime}: By \eqref{eq:Ups_norm_bounds_notime}, we have 
\begin{equation}\label{eq:Ups_norm_bounds}
	\max_{xy} (\Upsilon_t)_{xy} \le C_1(\ell_t\eta_t)^{-1}, \quad \max_x \sum_a (\Upsilon_t)_{xa} \le C_1 \eta_t^{-1}, \quad 0\le t\le T.
\end{equation}
It follows from \eqref{eq:convol_notime} that, for all $s \in  [0,t]$,
\begin{equation} \label{eq:convol}
	\sum_a \sqrt{(\Upsilon_s)_{xa} (\Upsilon_t)_{ay}} \le C_1\frac{1}{\eta_s}\sqrt{\ell_s\eta_s(\Upsilon_t)_{xy}}, \quad x,y\in\indset{N}.
\end{equation}
Finally, \eqref{eq:Ups_time_monot_notime} implies that
\begin{equation} \label{eq:Ups_time_monot}
	(\Upsilon_s)_{xy}  \le C_1   (\Upsilon_t)_{xy}, \quad x,y\in\indset{N}.
\end{equation}
\subsection{Action of the propagators: Proof of Lemmas \ref{lemma:good_props} and \ref{lemma:reg_props}} \label{sec:props}
We always assume that the times $s,t$ satisfy $0 \le s \le t \le T$.
Before proceeding with the proof of Lemmas \ref{lemma:good_props} and \ref{lemma:reg_props}, we collect some necessary properties of the propagators $\mathcal{P}_{s,t}$ and $\mathcal{Q}_{s,t}$, defined in \eqref{eq:Psat_def} and \eqref{eq:Q_desat_def}, respectively. Observe that
\begin{equation} \label{eq:P_decomp}
	\mathcal{P}_{s,t} = \mathrm{e}^{t-s}I + \bigl(\mathrm{e}^{t-s}-1\bigr)\Theta_t, \quad \mathcal{Q}_{s,t} = \mathrm{e}^{t-s}I + \bigl(\mathrm{e}^{t-s}-1\bigr)\Xi_t.
\end{equation}
Using the fact that  $\mathrm{e}^{t-s}-1  \sim \eta_s-\eta_t \lesssim \eta_s$, $\eta_{0}\ge 1$,
\eqref{eq:Ups_time_monot_notime}  and \eqref{eq:Ups_majorates}, we obtain the entry-wise bounds
\begin{equation} \label{eq:prop_xy_estimates}
	(\mathcal{P}_{s,t})_{ab} \lesssim \delta_{ab} + \eta_s(\Upsilon_t)_{ab}, \quad \bigl\lvert(\mathcal{Q}_{s,t})_{ab}\bigr\rvert \lesssim \delta_{ab} + \eta_s(\Upsilon_{0})_{ab}.
\end{equation}
Note that $\mathcal{P}_{s,t}$ is entry-wise positive.

Using \eqref{eq:prop_xy_estimates},  
we conclude that
\begin{equation} \label{eq:P_ss_bound}
	\sum_{a} (\mathcal{P}_{s,t})_{ua} \sqrt{(\Upsilon_s)_{xa}(\Upsilon_s)_{ay}} \lesssim \sqrt{(\Upsilon_s)_{xu}(\Upsilon_s)_{uy}} + \sqrt{(\Upsilon_t)_{xu}(\Upsilon_t)_{uy}} \lesssim \sqrt{(\Upsilon_t)_{xu}(\Upsilon_t)_{uy}},
\end{equation}
for all $x,y,u \in \indset{N}$, where in the first step we wrote $(\Upsilon_t)_{ua}$ from the propagator as 
$\sqrt{(\Upsilon_t)_{ua}}\sqrt{(\Upsilon_t)_{ua}}$ and used \eqref{eq:Schwarz_convol};   while in the last step
we used time-monononicity of $\Upsilon_t$ from \eqref{eq:Ups_time_monot}. 
In other words, \eqref{eq:P_ss_bound} shows that  the propagator $\mathcal{P}_{s,t}$ acts on $\Upsilon_s$'s
practically as a delta function on the spatial components, and the time index $s$ on $\Upsilon_s$ 
is changed to $t$ corresponding to the mechanism that the action of the propagator $\mathcal{P}_{s,t}$ creates 
a function living on the scale $\ell_t$.

Completely analogously, we deduce that the unsaturated propagator $\mathcal{Q}_{s,t}$ satisfies
\begin{equation} \label{eq:Q_ss_bound}
	\sum_{a} \bigl\lvert (\mathcal{Q}_{s,t})_{ua} \bigr\rvert \sqrt{(\Upsilon_{r_1})_{xa}(\Upsilon_{r_2})_{ay}} \lesssim \sqrt{(\Upsilon_{r_1})_{xu}(\Upsilon_{r_2})_{uy}},  \quad u,x,y \in \indset{N}, \quad r_1,r_2 \in [0,T].
\end{equation} 
Notice that $|\mathcal{Q}_{s,t}|$ practically acts as a delta function on $\Upsilon$'s,  this is because $\Upsilon$'s live
on a scale $\ell$, while the effective range of $|\mathcal{Q}_{s,t}|$ is $W\lesssim \ell$. 
On the other hand, using \eqref{eq:triag}--\eqref{eq:convol}, and \eqref{eq:Ups_time_monot}, similarly to the proof of~\eqref{eq:P_ss_bound}, 
we have
\begin{equation} \label{eq:P_st_bound}
	\sum_{a} (\mathcal{P}_{s,t})_{ua} \sqrt{(\Upsilon_t)_{xa}(\Upsilon_s)_{ay}} \lesssim \sqrt{\frac{\ell_s\eta_s}{\ell_t\eta_t}} \sqrt{(\Upsilon_t)_{xu}(\Upsilon_t)_{uy}}, \quad u,x,y \in \indset{N}.
\end{equation}
We stress that at this point we needed the stronger convolution estimate in the form~\eqref{eq:convol_notime}
and its weaker version~\eqref{eq:true_convol_notime} would not suffice.

Note that the bound \eqref{eq:P_st_bound} is worse by a factor of $\sqrt{(\ell_s\eta_s)/(\ell_t\eta_t)}$ compared to \eqref{eq:Q_ss_bound} with $r_1 := t, r_2:= s$, and the time in the superscript of $\Upsilon$ on the right-hand side of \eqref{eq:P_st_bound} are replaced by $t$.  Observe also the difference between \eqref{eq:P_ss_bound} and \eqref{eq:P_st_bound}; in the former both $\Upsilon$'s
live on the shorter scale $\ell_s$, while in the latter there is a mismatch of scales of the two $\Upsilon$ factors. 
We record that if both $\Upsilon$'s are at time $t$, then  we have (see also~\eqref{eq:last_prop_bound} later) 
\begin{equation} \label{eq:P_tt_bound}
	\sum_{a} (\mathcal{P}_{s,t})_{ua} \sqrt{(\Upsilon_t)_{xa} (\Upsilon_t)_{ay}} \lesssim \frac{\eta_s}{\eta_t}\sqrt{(\Upsilon_t)_{xu}(\Upsilon_t)_{uy}}, \quad u,x,y \in \indset{N}.
\end{equation}

So far we analysed the action of $\mathcal{P}$ and $\mathcal{Q}$.
Note that the remaining trivial propagator, given by the identity matrix $I$, that can appear in \eqref{eq:general_prop}, is bounded from above by $|\mathcal{Q}|$, and hence can be treated the same way as the non-saturated propagator $\mathcal{Q}$ using \eqref{eq:Q_ss_bound}.
Armed with \eqref{eq:P_ss_bound}--\eqref{eq:P_st_bound}, we are ready to prove Lemmas \ref{lemma:good_props} and Lemma \ref{lemma:reg_props}.

\begin{proof} [Proof of Lemma \ref{lemma:good_props}]
	First, we prove \eqref{eq:iso_prop_bound}. Recall the definition of $\mathfrak{s}_{k+1,s}^\mathrm{iso}$ from \eqref{eq:sfunc_def}, 
	\begin{equation}
		\mathfrak{s}_{k+1,s}^\mathrm{iso}(a,\bm x,b)  = \frac{1}{(\ell_s\eta_s)^{k/2}}\sqrt{\prod_{j=1}^{k+1}(\Upsilon_s)_{x_{j-1}x_{j}}},
	\end{equation}
	where we use the convention $x_0 := a$ and $x_{k+1} := b$ for brevity.
	Note that every index $x_j$ for $j\in\indset{k}$ occurs in exactly two $\sqrt{(\Upsilon_s)_{\dots}}$-factors. 
	We apply propagators $\other{\mathcal{P}}^{(j)}_{s,t}$ sequentially. It follows from \eqref{eq:P_ss_bound} and \eqref{eq:Q_ss_bound} with $r_1=r_2=s$ that, if $f_{ab}(\bm x)$ satisfies the assumption of \eqref{eq:iso_prop_bound}, then for all $a,b\in \indset{N}$ and $\bm x \in \indset{N}^k$,
	\begin{equation}
		\biggl\lvert \other{\mathcal{P}}^{(1)}_{s,t}\bigl[f_{ab} \bigr]  (\bm x)\biggr\rvert \lesssim  \frac{\varphi}{(\ell_s\eta_s)^{k/2}}\sqrt{(\Upsilon_t)_{x_0x_1}(\Upsilon_t)_{x_1x_2}\prod_{j=3}^{k+1}(\Upsilon_s)_{x_{j-1}x_{j}}}.
	\end{equation}
	Note that the power of $(\ell_s\eta_s)$ is the same as in the definition of $\mathfrak{s}_{k+1,s}^\mathrm{iso}$
	that bounds $f_{ab}(\bm x)$. 
	For all subsequent propagators indexed by $p \in \indset{2,k}$, using \eqref{eq:P_st_bound} and \eqref{eq:Q_ss_bound} with $r_1 = t, r_2=s$, we iteratively obtain   
	\begin{equation} \label{eq:prop_appl_bound}
		\biggl\lvert \other{\mathcal{P}}^{(p)}_{s,t} \circ \dots \circ \other{\mathcal{P}}^{(1)}_{s,t}\bigl[f_{ab} \bigr]  (\bm x)\biggr\rvert \lesssim \frac{\varphi }{(\ell_s\eta_s)^{\frac{k-p+1}{2}}(\ell_t\eta_t)^{\frac{p-1}{2}}}\sqrt{\prod_{i=1}^{p+1}(\Upsilon_t)_{x_{i-1}x_{i}}\prod_{j=p+2}^{k+1}(\Upsilon_s)_{x_{j-1}x_{j}}}, 
	\end{equation}
	where we treat $\other{\mathcal{P}}^{(j)}_{s,t}$ as acting on the $j$-th component of the tensor-product space 
	$(\mathbb{C}^N)^{\otimes k}$. Note that in the worst case
	 each step may pick up a large factor $\sqrt{(\ell_s\eta_s)/(\ell_t\eta_t)}$
	from  \eqref{eq:P_st_bound}.
	In particular, \eqref{eq:iso_prop_bound} follows from \eqref{eq:prop_appl_bound} and \eqref{eq:sfunc_def}.
	
	Next, we prove \eqref{eq:av_prop_bound}. If $k=1$, $\mathfrak{s}_{1,s}^\mathrm{av}(x_1) = (\ell_s\eta_s)^{-1}$ is a constant, hence for any $g(\bm x)$ satisfying the assumption of \eqref{eq:av_prop_bound}, it holds that
	\begin{equation}
		\bigl\lvert \other{\mathcal{P}}^{(1)}_{s,t}\bigl[g \bigr] (x_1)\bigr\rvert \lesssim \sum_a \bigl\lvert 
		(\mathcal{Q}_{s,t})_{x_1a} \bigr\rvert \frac{1}{\ell_s\eta_s} \lesssim \frac{1}{\ell_s\eta_s},
	\end{equation}
	where in the last step we used \eqref{eq:prop_xy_estimates} and \eqref{eq:Ups_norm_bounds}. In particular, \eqref{eq:av_prop_bound} holds with $k=1$. Therefore, it remains to consider $k\ge 2$.
	
	Recall the definition of $\mathfrak{s}_{k,s}^\mathrm{av}$ from \eqref{eq:sfunc_def}, and note that 
	\begin{equation}
		\mathfrak{s}_{k,s}^\mathrm{av}(\bm x)  = (\ell_s\eta_s)^{-1/2} \mathfrak{s}_{k,s}^\mathrm{iso}(x_k,\bm x',x_k).
	\end{equation}
	Owing to the cyclic structure of $\mathfrak{s}_{k,s}^\mathrm{av}$, we can assume without loss of generality, that the index $j$ of the unsaturated propagator in $\other{\mathcal{P}}^{k}_{s,t}$ is $j=k$.
	Therefore, using \eqref{eq:iso_prop_bound} for $f_{x_kx_k}(\bm x'):= g(\bm x)$, we obtain 
	\begin{equation} \label{eq:k-1_props}
		\biggl\lvert \other{\mathcal{P}}^{(k-1)}_{s,t} \circ \dots \circ \other{\mathcal{P}}^{(1)}_{s,t}\bigl[f_{x_kx_k} \bigr]  (\bm x)\biggr\rvert \lesssim \frac{\ell_t\eta_t}{\ell_s\eta_s}\frac{\varphi }{(\ell_t\eta_t)^{\frac{k}{2}}}\sqrt{(\Upsilon_t)_{x_{k}x_{1}}\prod_{i=2}^{k}(\Upsilon_t)_{x_{i-1}x_{i}}} = \varphi\,\frac{\ell_t\eta_t}{\ell_s\eta_s}\,\mathfrak{s}_{k,t}^\mathrm{av}(\bm x).
	\end{equation}
	Since the final propagator $\other{\mathcal{P}}^{(k)}_{s,t} \neq \mathcal{P}_{s,t}$, using \eqref{eq:Q_ss_bound} with $r_1=r_2=t$, we conclude \eqref{eq:av_prop_bound} from \eqref{eq:k-1_props}.
	This concludes the proof of Lemma \ref{lemma:good_props}.
\end{proof}

For the proof of the regularized propagator estimates, Lemma~\ref{lemma:reg_props}, we need regularity
estimates for $\Theta_t$, expressed in terms of upper bounds on $\reg{\Theta}_t$. These immediately follow
from the corresponding regularity conditions~\eqref{eq:Theta_regularity_notime} and~\eqref{eq:supercrit_Theta_notime};
we list them here for further reference:
\begin{lemma}[Regularity of Theta] \label{lemma:Theta_regular}
	Let $\crit$ be the \emph{critical time} along the characteristic flow~\eqref{eq:char_flow}
	 satisfying $\eta_{\crit} = (W/N)^2$.  
	If $t \in [0, \crit \wedge T]$, then the entries of $\reg{\Theta}_t^x$, defined in \eqref{eq:Thetaring}, satisfy
	\begin{equation} \label{eq:Theta_regularity}
		\bigl\lvert(\reg{\Theta}_t^x)_{ab}\bigr\rvert  \lesssim  \frac{\bigl(|b-x|_N+W\bigr)\wedge\ell_t}{\ell_t} \bigl((\Upsilon_t)_{ab}+(\Upsilon_t)_{ax}\bigr), \quad a,b,x \in \indset{N}.
	\end{equation}
	On the other hand, if $t \in [\crit, T]$, then the entries of $\Theta_t$ satisfy
	\begin{equation} \label{eq:supercrit_Theta}
		(\Theta_t)_{ab} \sim \frac{1}{N\eta_t}, \quad  \biggl\lvert (\Theta_t)_{ab} - \frac{\im m_t}{N\eta_t} \biggr\rvert \lesssim \frac{N}{W^2}, \quad a,b \in \indset{N}.
	\end{equation}
\end{lemma}

We recall the weighted convolution property of the admissible control function $\Upsilon$.
It follows immediately from \eqref{eq:suppressed_convol_notime} that, for all $s \in  [0,t]$ and $r \in \{s,t\}$, 
\begin{equation} \label{eq:suppressed_convol}
	\sum_{a} \frac{\bigl(|a-x|_N + W\bigr)\wedge \ell_t}{\ell_t} \sqrt{(\Upsilon_s)_{xa} (\Upsilon_r)_{ay}} \le C_1 \frac{1}{\eta_s} \frac{\ell_s}{\ell_t} \sqrt{\frac{\ell_s\eta_s}{\ell_r\eta_r}} \sqrt{\ell_t\eta_t (\Upsilon_t)_{xy}}, \quad x,y\in\indset{N}.
\end{equation}
 
\begin{proof}[Proof of Lemma \ref{lemma:reg_props}]	
	Recall from \eqref{eq:sfunc_def} that, in the definition of $\mathfrak{s}_k^{\mathrm{av}}$, exactly two decaying factors $\sqrt{(\Upsilon_s)_{\dots}}$  depend on the external index $x_j$ for any given $j \in \indset{k}$.
	We claim the following estimate on the action of a single $\reg{\Theta}_t^{(j)}$, defined in \eqref{eq:Thetaring} for $j \in \indset{2,k}$, on the product of two decaying factors,
	\begin{equation} \label{eq:Thetaring_action_basic}
		\biggl\lvert\reg{\Theta}_t^{(j)}\biggl[\sqrt{(\Upsilon_s)_{x_{j-1}x_j}(\Upsilon_r)_{x_jx_{j+1}}}\biggr]\biggr\rvert \lesssim \frac{1}{\eta_s} \frac{\ell_s\ell_t\eta_t}{W^2} \sqrt{\frac{\ell_s\eta_s}{\ell_r\eta_r}} \sqrt{(\Upsilon_t)_{x_{j-1}x_j}(\Upsilon_t)_{x_jx_{j+1}}},
	\end{equation}
	for all $s \in [0,t]$ and all $r \in \{s,t\}$, where we adhere to the cyclic convention $x_{k+1} := x_1$ and $x_0 := x_k$.
 It is useful to compare this bound with the analogous estimate (obtained easily from~\eqref{eq:Ups_majorates}, \eqref{eq:triag}, and 
\eqref{eq:convol})
 	\begin{equation} \label{eq:Thetaring_action_basic1}
		\biggl\lvert\Theta_t^{(j)}\biggl[\sqrt{(\Upsilon_s)_{x_{j-1}x_j}(\Upsilon_r)_{x_jx_{j+1}}}\biggr]\biggr\rvert \lesssim \frac{1}{\eta_s} 
		 \sqrt{\frac{\ell_s\eta_s}{\ell_r\eta_r}} \sqrt{(\Upsilon_t)_{x_{j-1}x_j}(\Upsilon_t)_{x_jx_{j+1}}},
	\end{equation}
where $\Theta_j$ is not regularized; notice the improvement factor  $\ell_s\ell_t\eta_t/W^2$ owing to the regularization
which is typically much smaller than 1 (see discussion below~\eqref{eq:noring}).
	
	We defer the derivation of the prototypical estimate \eqref{eq:Thetaring_action_basic} until the end of the proof. In the subsequent analysis, we use $r=s$ and $r=t$ to bound the actions of $\reg{\Theta}_t^{(k)}$ and $\reg{\Theta}_t^{(k-2)}$, respectively. 
	
.
	
	First, we prove \eqref{eq:one_ring}.  
	It follows from \eqref{eq:Thetaring_action_basic} with $r=s$ that for any $f(\bm x)$ satisfying the assumption of \eqref{eq:one_ring}, we have 
	\begin{equation} \label{eq:Thetaring_action}
		\biggl\lvert \reg{\Theta}_t^{(k)} \bigl[f \bigr]  (\bm x)\biggr\rvert \lesssim \varphi \, \frac{1}{\eta_s} \frac{\ell_s\ell_t\eta_t}{W^2}  \frac{1}{(\ell_s\eta_s)^{k/2}}\sqrt{(\Upsilon_t)_{x_kx_1}(\Upsilon_t)_{x_{k-1}x_{k}}\prod_{j=2}^{k-1}(\Upsilon_s)_{x_{j-1}x_{j}}}.
	\end{equation}
	Using \eqref{eq:P_st_bound}, we deduce, similarly to \eqref{eq:k-1_props}, that
	\begin{equation} \label{eq:1ring_k-2}
		\biggl\lvert \mathcal{P}^{(k-2)}_{s,t} \circ \dots \circ \mathcal{P}^{(1)}_{s,t} \circ \reg{\Theta}_t^{(k)} \bigl[f \bigr]  (\bm x)\biggr\rvert \lesssim \varphi \, \frac{\eta_t}{\eta_s^2}\frac{\ell_t^2\eta_t}{W^2}  \mathfrak{s}_{k,t}^\mathrm{av}(\bm x).
	\end{equation}
	To estimate the action of the remaining two propagators, we observe that, by Schwarz inequality together with \eqref{eq:P_decomp} and \eqref{eq:true_convol},
	\begin{equation} \label{eq:last_prop_bound}
		\mathcal{P}_{s,t}^{(j)}\bigl[\mathfrak{s}_{k,t}^\mathrm{av}(\bm x)\bigr] \lesssim  \frac{\eta_s}{\eta_t} \mathfrak{s}_{k,t}^\mathrm{av}(\bm x), \quad j \in \indset{k},
	\end{equation}
	where we used that $\eta_s\ge \eta_t$.
	Hence, acting by $\mathcal{P}_{s,t}^{(k)}\circ{P}_{s,t}^{(k-1)}$ on both sides of  \eqref{eq:1ring_k-2} and using \eqref{eq:last_prop_bound} for $j=k-1, k$, we deduce that 
	\begin{equation}
		\biggl\lvert \mathcal{P}^{\otimes k}_{s,t} \circ \reg{\Theta}_t^{(k)} \bigl[f \bigr]  (\bm x)\biggr\rvert \lesssim \varphi \, \frac{\eta_s^2}{\eta_t^2}\frac{\eta_t}{\eta_s^2}\frac{\ell_t^2\eta_t}{W^2} \mathfrak{s}_{k,t}^\mathrm{av}(\bm x) = \frac{1}{\eta_t}\frac{\ell_t^2\eta_t}{W^2}\mathfrak{s}_{k,t}^\mathrm{av}(\bm x) \lesssim \frac{1}{\eta_t}\mathfrak{s}_{k,t}^\mathrm{av}(\bm x),
	\end{equation}
	where in the last step we used $\ell_t^2 \eta_t \lesssim W^2$ .This concludes the proof of \eqref{eq:one_ring}.
	
	Next, we prove \eqref{eq:two_rings} under the assumption $k \ge 4$. 
	It follows form \eqref{eq:P_st_bound} and  \eqref{eq:Thetaring_action} that
	\begin{equation} \label{eq:thringP_est}
		\biggl\lvert \mathcal{P}^{(k-1)}_{s,t}\circ \reg{\Theta}_t^{(k)} \bigl[f \bigr]  (\bm x)\biggr\rvert \lesssim \varphi \, \frac{1}{\eta_s} \frac{\ell_s\ell_t\eta_t}{W^2}\frac{\sqrt{\ell_s\eta_s}}{\sqrt{\ell_t\eta_t}}  \frac{1}{(\ell_s\eta_s)^{k/2}}\sqrt{(\Upsilon_t)_{x_kx_1}\prod_{j=2}^{k-2}(\Upsilon_s)_{x_{j-1}x_{j}}\prod_{j=k-1}^{k}(\Upsilon_t)_{x_{j-1}x_{j}} }.
	\end{equation} 
	Therefore, combining \eqref{eq:thringP_est} and \eqref{eq:Thetaring_action_basic} with $r=t$, we deduce that
	\begin{equation} \label{eq:2thringP_est}
		\biggl\lvert \reg{\Theta}_t^{(k-2)} \circ \mathcal{P}^{(k-1)}_{s,t}\circ \reg{\Theta}_t^{(k)} \bigl[f \bigr]  (\bm x)\biggr\rvert \lesssim  
		\frac{\varphi}{\eta_s^2} \frac{\ell_s^2}{W^4} \frac{\ell_t\eta_t}{(\ell_s\eta_s)^{k/2-1}}\sqrt{(\Upsilon_t)_{x_kx_1}\prod_{j=2}^{k-3}(\Upsilon_s)_{x_{j-1}x_{j}}\prod_{j=k-2}^{k}(\Upsilon_t)_{x_{j-1}x_{j}} }.
	\end{equation}
	By \eqref{eq:P_st_bound}, applying the propagators acting on the  indices $1, \dots, (k-4)$ to both sides of \eqref{eq:2thringP_est}, and using $\ell_t^2\eta_t \lesssim W^2$, $\ell_s \le \ell_t$, we obtain
	\begin{equation} \label{eq:2ring_k-4}
		\biggl\lvert \mathcal{P}^{(k-4)}_{s,t} \circ \dots \circ \mathcal{P}^{(1)}_{s,t} \circ \reg{\Theta}_t^{(k-2)}\circ \mathcal{P}^{(k-1)}_{s,t}\circ \reg{\Theta}_t^{(k)} \bigl[f \bigr]  (\bm x)\biggr\rvert 
		\lesssim \varphi\, \frac{\eta_t^3}{\eta_s^3} \frac{1}{\eta_t^2} \frac{\ell_t\eta_t}{\ell_s\eta_s} \mathfrak{s}_{k,t}^\mathrm{av}(\bm x).
	\end{equation}
	By \eqref{eq:last_prop_bound}, the remaining three propagators $\mathcal{P}^{(k)}_{s,t}\circ\mathcal{P}^{(k-2)}_{s,t}\circ\mathcal{P}^{(k-3)}_{s,t}$ contribute a factor $\eta_s^3/\eta_t^3$, hence yielding \eqref{eq:two_rings}, since $\reg{\Theta}_t^{(k-2)}$ commutes with $\mathcal{P}_{s,t}^{(k-1)}$.  
	
	Therefore, it remains to prove \eqref{eq:Thetaring_action_basic}.  The proof is slightly different for $t\le t_*$ 
	and $t\ge t_*$,  since it uses different inputs \eqref{eq:Theta_regularity} and \eqref{eq:supercrit_Theta}, respectively.
	Without loss of generality we assume that $j=k$, other cases are analogous.
	\medskip
	
	{\bf Case 1.  $t \le \crit$.}  
	In this case it follows from \eqref{eq:Theta_regularity} that 
	\begin{equation} \label{eq:Thetaring_est}
		\begin{split}
			\biggl\lvert\reg{\Theta}_t^{(k)}\biggl[\sqrt{(\Upsilon_s)_{x_{k-1}x_k}(\Upsilon_r)_{x_kx_{k+1}}}\biggr]\biggr\rvert 
			\lesssim&~ \sum_{b}\frac{|b-x_{k-1}|_N\wedge\ell_t}{\ell_t} (\Upsilon_t)_{x_{k}b} \sqrt{(\Upsilon_s)_{x_{k-1}b}(\Upsilon_r)_{bx_{k+1}}}\\
			&+(\Upsilon_t)_{x_kx_{k-1}}  \sum_{b}\frac{|b-x_{k-1}|_N\wedge\ell_t}{\ell_t} \sqrt{(\Upsilon_s)_{x_{k-1}b}(\Upsilon_r)_{bx_{k+1}}}.
		\end{split}
	\end{equation}
	To estimate the first sum  on the right-hand side of \eqref{eq:Thetaring_est}, we use \eqref{eq:triag} for $\sqrt{(\Upsilon_t)_{x_{k}b}(\Upsilon_r)_{bx_{k+1}}}$ and \eqref{eq:suppressed_convol},
	\begin{equation} \label{eq:Thetaring_est1}
		\begin{split}
			\sum_{b}\frac{|b-x_{k-1}|_N\wedge\ell_t}{\ell_t} (\Upsilon_t)_{x_{k}b} \sqrt{(\Upsilon_s)_{x_{k-1}b}(\Upsilon_r)_{bx_{k+1}}} \lesssim \frac{1}{\eta_s} \frac{\ell_s}{\ell_t} \sqrt{\frac{\ell_s\eta_s}{\ell_r\eta_r}} \sqrt{(\Upsilon_t)_{x_{k-1}x_k}(\Upsilon_t)_{x_kx_{k+1}}}.
		\end{split}
	\end{equation} 
	Using \eqref{eq:suppressed_convol} followed by \eqref{eq:triag} to estimate the second sum on the right-hand side of \eqref{eq:Thetaring_est}, we obtain
	\begin{equation}\label{eq:Thetaring_est2}
		\begin{split} 
			(\Upsilon_t)_{x_kx_{k-1}}  \sum_{b}\frac{|b-x_{k-1}|_N\wedge\ell_t}{\ell_t} \sqrt{(\Upsilon_s)_{x_{k-1}b}(\Upsilon_r)_{bx_{k+1}}} 
			&\lesssim \frac{1}{\eta_s} \frac{\ell_s}{\ell_t} \sqrt{\frac{\ell_s\eta_s}{\ell_r\eta_r}}\sqrt{(\Upsilon_t)_{x_{k-1}x_k}(\Upsilon_t)_{x_kx_{k+1}}}.
		\end{split}
	\end{equation} 
	It follows from \eqref{eq:Thetaring_est1}, \eqref{eq:Thetaring_est2}, and $\ell_t^2\eta_t \sim W^2$, that \eqref{eq:Thetaring_action_basic} holds for $t \le \crit$. 
	
	\medskip
	{\bf Case 2.  $t \ge \crit$.} In this case, by definition of $\mathcal{P}_{\crit,t}$ in \eqref{eq:Psat_def}, the matrix $\reg{\Theta}_t^x$, defined in \eqref{eq:Thetaring}, satisfies
	\begin{equation} \label{eq:Theta_prop}
		\reg{\Theta}_t^{x} = \mathcal{P}_{\crit,t}\reg{\Theta}_\crit^{x}.
	\end{equation}
	Hence, plugging the decomposition \eqref{eq:P_decomp} into \eqref{eq:Theta_prop}, we obtain
	\begin{equation} \label{eq:Theta_pre_decom}
		\reg{\Theta}_t^{x} = \mathrm{e}^{t-\crit} \reg{\Theta}_\crit^{x} + \bigl(\mathrm{e}^{t-\crit}-1\bigr)\Theta_t\reg{\Theta}_\crit^{x}.
	\end{equation}
	Let $\Pi := \bm1\bm1^*$ be the projector onto the constant vector $\bm 1 := N^{-1/2}(1,\dots, 1)\in\mathbb{C}^N$. It follows from \eqref{eq:sumTheta} that $\Pi\Theta_t = \Theta_t \Pi = \frac{\im m_t}{N\eta_t} \Pi$. Moreover, \eqref{eq:Thetaring} implies that $\Pi \reg{\Theta}^x_t = 0$ for all $x \in \indset{N}$.
	 Therefore, we deduce from \eqref{eq:Theta_pre_decom} that
	\begin{equation} \label{eq:Theta_decom}
		\reg{\Theta}_t^{x} = \mathrm{e}^{t-\crit} \reg{\Theta}_\crit^{x} + \bigl(\mathrm{e}^{t-\crit}-1\bigr)\bigl(\Theta_t - \Pi\Theta_t\bigr)\reg{\Theta}_\crit^{x}.
	\end{equation}
	Therefore, it suffices to estimate the actions of $\reg{\Theta}_\crit^{x}$ and 
	$(\Theta_t - \Pi\Theta_t)\reg{\Theta}_\crit^{x}$ on the $x=x_k$ argument of $\reg{\Theta}_t^{(k)}$ in \eqref{eq:Thetaring_action_basic} (recall we set $j=k$).

	Analogously to the proof of  \eqref{eq:Thetaring_action_basic} for $t \le \crit$ above, using \eqref{eq:suppressed_convol}
	and the first estimate in \eqref{eq:supercrit_Theta}, we deduce that
	\begin{equation} \label{eq:crit_Thetatirg}
		\biggl\lvert\reg{\Theta}_\crit^{(k)}\biggl[\sqrt{(\Upsilon_s)_{x_{k-1}x_k}(\Upsilon_r)_{x_kx_{k+1}}}\biggr]\biggr\rvert 
		\lesssim \frac{1}{\eta_s}\frac{\ell_s}{\ell_\crit} \frac{1}{\ell_\crit\eta_\crit} \sqrt{\frac{\ell_s\eta_s}{\ell_r\eta_r}}
		 \sim  \frac{\ell_s}{W^2\eta_s} \sqrt{\frac{\ell_s\eta_s}{\ell_r\eta_r}}, \quad r\in\{s,t\}.
	\end{equation}
	Furthermore, the second estimate in \eqref{eq:supercrit_Theta} and \eqref{eq:crit_Thetatirg} imply, for $r\in \{s,t\}$,
	\begin{equation} \label{eq:crit_Thetatirg2}
		\biggl\lvert\sum_{ab}\bigl(\Theta_t - \Pi\Theta_t\bigr)_{x_kb}\bigl(\reg{\Theta}_\crit^{(k)}\bigr)_{ba} \sqrt{(\Upsilon_s)_{x_{k-1}a}(\Upsilon_r)_{ax_{k+1}}} \biggr\rvert \lesssim \sum_b \frac{N}{W^2} 
		 \frac{\ell_s}{W^2\eta_s} \sqrt{\frac{\ell_s\eta_s}{\ell_r\eta_r}} \lesssim   \frac{N^2}{W^2}\frac{\ell_s}{W^2\eta_s} \sqrt{\frac{\ell_s\eta_s}{\ell_r\eta_r}}.
	\end{equation}
	Recall that the scalar factors in \eqref{eq:Theta_decom} satisfy $\mathrm{e}^{t-\crit} \sim 1$ and $(\mathrm{e}^{t-\crit}-1) \sim t-\crit \lesssim \eta_\crit$ by \eqref{eq:etatasymp}. 
	Therefore, combining \eqref{eq:Theta_decom}, \eqref{eq:crit_Thetatirg} and \eqref{eq:crit_Thetatirg2}, we conclude that, for $t \ge \crit$, $r\in \{s,t\}$,
	\begin{equation} \label{eq:crit_Thetatirg3}
		\begin{split}
			\biggl\lvert\reg{\Theta}_t^{(k)}\biggl[\sqrt{(\Upsilon_s)_{x_{k-1}x_k}(\Upsilon_r)_{x_kx_{k+1}}}\biggr]\biggr\rvert &\lesssim \biggl(1 + \eta_\crit\frac{N^2}{W^2}\biggr)\frac{\ell_s}{W^2\eta_s} \sqrt{\frac{\ell_s\eta_s}{\ell_r\eta_r}}\\
			&\lesssim \frac{1}{\eta_s} \frac{\ell_s\ell_t\eta_t}{W^2} \sqrt{\frac{\ell_s\eta_s}{\ell_r\eta_r}} \sqrt{(\Upsilon_t)_{x_{k-1}x_k}(\Upsilon_t)_{x_kx_{k+1}}},
		\end{split}
	\end{equation}
	where we used the first bound in \eqref{eq:supercrit_Theta},
	together with~\eqref{eq:Ups_majorates} to have $(\Upsilon_t)_{ab} \gtrsim (\Theta_t)_{ab} \sim (N\eta_t)^{-1} 
	\sim 1/(\ell_t\eta_t)$, and that $\eta_\crit = (W/N)^2$ by definition of $\crit$. Hence, \eqref{eq:Thetaring_action_basic} also holds for $t \in [\crit, T]$.	
	This concludes the proof of \eqref{eq:Thetaring_action_basic}.
	
	Since the proof of \eqref{eq:Thetaring_action_basic} and subsequently \eqref{eq:one_ring} and \eqref{eq:two_rings} relies only on entry-wise bounds on $\reg{\Theta}_t$ from Lemma \ref{lemma:Theta_regular},  \eqref{eq:Theta_regularity} on $\reg{\Theta}_t$, estimates \eqref{eq:one_ring} and \eqref{eq:two_rings} also hold with $\reg{\Theta}_t$ replaced by its entry-wise absolute value $|\reg{\Theta}_t|$.
	This concludes the proof of Lemma \ref{lemma:reg_props}.
\end{proof}

\subsection{Estimates on martingale and forcing terms: Proof of Lemmas \ref{lemma:mart_est}, \ref{lemma:forcing}, \ref{lemma:extra_forcing}} \label{sec:mart_forcing}
We emphasize that in Section \ref{sec:mart_forcing}, we use only the weaker convolution estimate \eqref{eq:true_convol_notime} (in its time-dependent forms \eqref{eq:true_convol} and \eqref{eq:Schwarz_convol}) instead of \eqref{eq:convol_notime}--\eqref{eq:suppressed_convol_notime}.

In this section we use the following bounds on resolvent chain of length up to $\maxK$. It follows from \eqref{eq:Psi_def} and \eqref{eq:tau_def} that for all $k\in\indset{\maxK}$, $\bm z_s \in \{z_s, \overline{z}_s\}^k$ and $\bm x \in \indset{N}^k$,
\begin{equation} \label{eq:G-M_psi_bounds}
	\begin{alignedat}{3}
		\biggl\lvert \bigl(&(G-M)_{[1,k]}(\bm z_s, \bm x')\bigr)_{ab} \biggr\rvert &&\le \pis{k}(\ell_s\eta_s)^{\alpha_k}~\mathfrak{s}_{k,s}^\mathrm{iso}(a,\bm x',b), \quad a,b\in\indset{N}, \quad &&s\le \tau, \\
		\biggl\lvert\Tr \bigl[&(G-M)_{[1,k]}(\bm z_s, \bm x') S^{x_k} \bigr]  \biggr\rvert &&\le \pav{k}(\ell_s\eta_s)^{\beta_k}~\mathfrak{s}_{k,s}^\mathrm{av}(\bm x), \quad &&s\le \tau.
	\end{alignedat}
\end{equation}
Moreover, combining \eqref{eq:G-M_psi_bounds} with \eqref{eq:M_bound}, we obtain, for all $k\in\indset{\maxK}$, $\bm z_s \in \{z_s, \overline{z}_s\}^k$ and $\bm x \in \indset{N}^k$,
\begin{equation} \label{eq:G_psi_bounds}
	\begin{alignedat}{3}
		\biggl\lvert \bigl(&G_{[1,k]}(\bm z_s, \bm x')\bigr)_{ab} \biggr\rvert &&\le \bigl(\delta_{ab}\sqrt{\ell_s\eta_s}+\pis{k}(\ell_s\eta_s)^{\alpha_k}\bigr)~\mathfrak{s}_{k,s}^\mathrm{iso}(a,\bm x',b), \quad a,b\in\indset{N}, \quad &&s\le \tau,\\
		\biggl\lvert\Tr \bigl[&G_{[1,k]}(\bm z_s, \bm x') S^{x_k} \bigr]  \biggr\rvert &&\le \bigl(\ell_s\eta_s+\pav{k}(\ell_s\eta_s)^{\beta_k}\bigr)~\mathfrak{s}_{k,s}^\mathrm{av}(\bm x), \quad &&s\le \tau.
	\end{alignedat}
\end{equation}

Furthermore, we record the following deterministic estimates that we use repeatedly in the sequel.  
Since $S = |m_0|^{-2}\Theta_0 - S\Theta_0$ and the matrix $S\Theta_0$ entry-wise non-negative, we have
\begin{equation} \label{eq:SUps_comvol}
	S_{xy} \lesssim (\Theta_{0})_{xy} \lesssim (\Upsilon_{0})_{xy}, \quad \bigl(S\Upsilon_t\bigr)_{xy} \lesssim (\Upsilon_t)_{xy}.
\end{equation}

Throughout this section, we repeatedly encounter products of size functions $\mathfrak{s}_{k,t}^\mathrm{iso/av}$, defined in \eqref{eq:sfunc_def} with a shared index (denoted here by $q$) among their arguments. In some cases, the shared index $q$ is also summed. Recall from \eqref{eq:sfunc_def} that, up to a scalar factor, the size functions $\mathfrak{s}_{k,t}^\mathrm{iso/av}$ are given by products of decaying factors $\sqrt{(\Upsilon_t)_{\dots}}$. 
Therefore, such products can be estimated locally (involving only the decaying factors with the shared index $q$ among their indices) using inequalities \eqref{eq:triag} and   \eqref{eq:true_convol}
 (if $q$ is summed). 
 
More precisely, typically, at least four $\sqrt{\Upsilon}$-factors containing the $q$ index are involved in the product,
in which case we use  \eqref{eq:true_convol} via~\eqref{eq:Schwarz_convol}.  There is one exceptional case:
 terms of the form $(G\mathcal{S}[G-m]G)$ give rise to convolutions with only two factors containing $q$. To deal with such terms, we use Schwarz inequality together with \eqref{eq:true_convol} to obtain, for all $0\le s \le t \le T$,
\begin{equation} \label{eq:sqrt_convol}
	\sum_q \sqrt{(\Upsilon_s)_{xq} (\Upsilon_t)_{qy}} \le \sqrt{N} \biggl(\sum_q (\Upsilon_s)_{xq} (\Upsilon_t)_{qy} \biggr)^{1/2} \le \frac{1}{\eta_s}\sqrt{C} \sqrt{ N\eta_s} \sqrt{ (\Upsilon_t)_{xy}}, \quad x,y\in\indset{N}.
\end{equation}
Using \eqref{eq:sqrt_convol}, we obtain the following bound for the   convolution of two isotropic size functions, 
\begin{equation} \label{eq:iso_convol}
	\sum_q \mathfrak{s}_{i+1,t}^\mathrm{iso}(a,\bm x, q) \mathfrak{s}_{j+1,t}^\mathrm{iso}(q,\bm y, b) \lesssim \eta_t^{-1} 
	(N  \eta_t)^{1/2}\,\mathfrak{s}_{i+j+1}^\mathrm{iso}(a, \bm x, \bm y, b), \quad \bm x\in \indset{N}^i, \quad \bm y \in \indset{N}^j.
\end{equation}
The general heuristic that $(\Upsilon_t)_{xy}$ is essentially supported on $|x-y|_N \lesssim \ell_t$ suggests that the factor $(N\eta_t)^{1/2}$ should be replaced by $(\ell_t\eta_t)^{1/2}$, hence using Schwarz inequality in \eqref{eq:sqrt_convol} incurs a price of $\sqrt{N/\ell_t}$ in~\eqref{eq:iso_convol}.
 In the sequel we find that, whenever \eqref{eq:sqrt_convol}  has to be used, this price is affordable.  

Note that if the summation in $q$  is replaced by maximum over $q\in \indset{N}$, instead we obtain
\begin{equation} \label{eq:iso_concat}
	\max_{q} \mathfrak{s}_{i+1,t}^\mathrm{iso}(a,\bm x, q) \mathfrak{s}_{j+1,t}^\mathrm{iso}(q,\bm y, b) \lesssim (\ell_t\eta_t)^{-1/2}\,\mathfrak{s}_{i+j+1}^\mathrm{iso}(a, \bm x, \bm y, b), \quad \bm x\in \indset{N}^i, \quad \bm y \in \indset{N}^j,
\end{equation}
where we used \eqref{eq:triag}.

If we represent each $\sqrt{(\Upsilon_t)_{xy}}$-factor by an edge connecting vertices $x$ and $y$, we obtain the following diagrammatic depiction of \eqref{eq:iso_convol} (ignoring some overall common scalar factors 
in $\mathfrak{s}^\mathrm{iso}$ and dropping out all indices not adjacent to $q$).

\vspace{2pt}
\begingroup \hfil  \includegraphics[scale=.7]{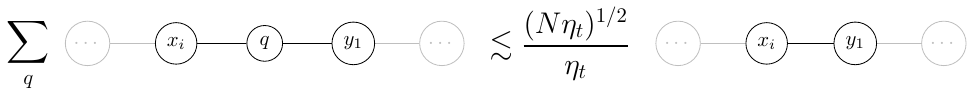} \hfill \endgroup

The following example is more complicated, but represents the typical situation when at least four $\sqrt{\Upsilon}$-factors contain the summation index $q$. 
Consider the   convolution of an averaged and an isotropic size functions, where the index $q$ is repeated four times in total;  using \eqref{eq:Schwarz_convol}, we obtain the bound  
\begin{equation} \label{eq:av_insertion}
	\sum_q \mathfrak{s}_{i+1,t}^\mathrm{av}(\bm x, q) \mathfrak{s}_{j+1,t}^\mathrm{iso}(q,\bm y, q) \lesssim \eta_t^{-1} (\ell_t\eta_t)^{-1/2}\,\mathfrak{s}_{i+j}^\mathrm{av}(\bm x, \bm y), \quad \bm x\in \indset{N}^i, \quad \bm y \in \indset{N}^j.
\end{equation}
The corresponding  diagram is depicted below. 

\vspace{3pt}
\begingroup \hfil  \includegraphics[scale=.7]{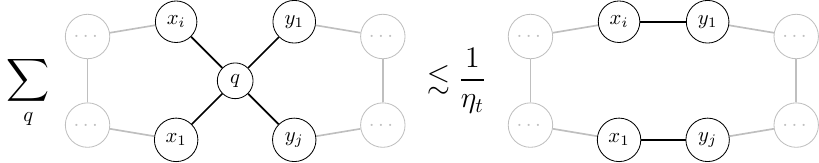} \hfill \endgroup

Note that due to the symmetric nature of $\mathfrak{s}_{j+1,t}^\mathrm{iso}(q,\bm y, q)$, the vector $\bm y$ on the right-hand side of \eqref{eq:av_insertion} can be replaced by $\bm y^* := (y_{j}, \dots, y_1)$ obtained by reversing the order of the elements in $\bm y$. Similarly to \eqref{eq:iso_concat} above, if the summation in $q$ is replaced by maximum over $q\in \indset{N}$, the estimate \eqref{eq:av_insertion} improves by a factor of $\ell_t^{-1}$: 
\begin{equation} \label{eq:av_concat}
	\max_q \mathfrak{s}_{i+1,t}^\mathrm{av}(\bm x, q) \mathfrak{s}_{j+1,t}^\mathrm{iso}(q,\bm y, q) \lesssim (\ell_t\eta_t)^{-3/2}\,\mathfrak{s}_{i+j}^\mathrm{av}(\bm x, \bm y), \quad \bm x\in \indset{N}^i, \quad \bm y \in \indset{N}^j.
\end{equation}
The improvement by a factor of $\ell_t^{-1}$ in \eqref{eq:av_concat} compared to \eqref{eq:av_insertion} is in line with the general heuristic that $(\Upsilon_t)_{xy}$ is essentially supported on $|x-y|_N \lesssim \ell_t$.  

As the final example, we consider the following double convolution involving two isotropic size functions and the variance profile $S$,
\begin{equation} \label{eq:iso_S}
\sum_{q'q} S_{q'q}  \bigl(\mathfrak{s}_{k+1}^\mathrm{iso}(q',\bm x, q)\bigr)^2  \lesssim \sum_{q} \bigl(\mathfrak{s}_{k+1}^\mathrm{iso}(q,\bm x, q)\bigr)^2  
\lesssim \eta_t^{-1}  (\mathfrak{s}_{k}^\mathrm{av}(\bm x)\bigr)^2, \quad \bm x \in \indset{N}^k.
\end{equation}
Here we used $S_{av} \lesssim (\Upsilon_0)_{ab}$ by the first bound in \eqref{eq:SUps_comvol},  then \eqref{eq:Schwarz_convol}   as in \eqref{eq:av_insertion} to first perform the $q'$ summation; the remaining $q$ summation is analogous to \eqref{eq:av_insertion}.

In Section~\ref{sec:mart} and \ref{sec:forcing}, we perform 
estimates similar to \eqref{eq:iso_convol}--\eqref{eq:iso_S} using  the following strategy: first, we  remove all but four decay factors involving $q$ using \eqref{eq:triag}, then we sum the remaining two factors using \eqref{eq:true_convol} and \eqref{eq:Schwarz_convol}, and finally we   assemble the resulting decaying factors into a size function using \eqref{eq:sfunc_def}. We will call such arguments {\it local contractions}
of the size functions; as above, they can easily be proven using the properties of 
the control function $\Upsilon$  described in Definition~\ref{def:adm_ups}.

\subsubsection{Martingale estimates} \label{sec:mart}
It follows immediately from \eqref{eq:ups_lower_bound_notime}, that, with $D'$ being the constant from \eqref{eq:ups_lower_bound_notime},
\begin{equation} \label{eq:ups_lower_bound}
	\min_{xy} (\Upsilon_t)_{xy} \ge N^{-2D'}.
\end{equation}

To bound the martingale terms in \eqref{eq:k_av_evol} and \eqref{eq:k_iso_evol}, we use the following path-wise Burkholder-Davis-Gundy inequality (see~\cite{martingale}, Appendix B.6, Eq. (18)). For any continuous martingale $\mathcal{M}_t$, finite stopping time $\tau$, and any fixed $x,y > 0$,
\begin{equation} \label{eq:mart_ineq}
	\mathbb{P}\biggl(\sup_{0\le s \le \tau} \biggl\lvert \int_0^{s} \mathrm{d}\mathcal{M}_r \biggr\rvert \ge x, \quad  \biggl[ \int_0^{\cdot} \mathrm{d}\mathcal{M}_r \biggr]_{\tau} \le y \biggr) \le 2\mathrm{e}^{-x^2/(2y)},
\end{equation}
where $[\cdot]_t$ denotes the quadratic variation process.
 
\begin{proof}[Proof of Lemma \ref{lemma:mart_est}] 
	\textbf{Averaged Chains.} Proof of \eqref{eq:mart_bound}.
	
	It follows from \eqref{eq:mart_ineq} together with a simple dyadic argument with $x := 2^j\log N$, $y := 2^{2j}$ for $j \in [-100LD'\log N, 100KD' \log N]$, 
	that 
	\begin{equation}
		\begin{split}
			\biggl\lvert \int_{\tinit}^{t\wedge\tau} \otherhat{\mathcal{P}}_{s,t\wedge\tau}\bigl[\mathrm{d}\mathcal{M}^\mathrm{av}_{[1,k],s}\bigr](\bm x) \biggr\rvert
			\prec&~ \biggl(\int_{\tinit}^{t\wedge\tau} \sum_{ab} S_{ab} \biggl\lvert \otherhat{\mathcal{P}}^k_{s,t\wedge\tau}\biggl[ \Tr\bigl[\partial_{ab}G_{[1,k],s}S^{x_k}\bigr]\biggr](\bm x) \biggr\rvert^2\mathrm{d}s\biggr)^{1/2}\\ &+ N^{-kD'-k/2},
		\end{split}
	\end{equation}
	where $D'$ is the constant from \eqref{eq:ups_lower_bound} and \eqref{eq:ups_lower_bound_notime}.
	
	For any chain of length $k\in\indset{\maxK}$,  we define the positive quantity $\QV^{(j)}_{k,s}(\bm z_s, \bm x)$ as 
	\begin{equation} \label{eq:Qkj}
		\QV^{(j)}_{k,s}(\bm z_s,\bm x) :=  \sum_{ab} S_{ab} \bigl\lvert \bigl(G_{[j,k],s}S^{x_k} G_{[1,j],s}\bigr)_{ab} \bigr\rvert^2, \quad j\in\indset{k},
	\end{equation}
	where $G_{[j,k],s}$ and $G_{[1,j],s}$ are resolvent chains defined according to \eqref{eq:Gk_def} for a set of external indices $\bm x\in\indset{N}^k$ and an arrangement of spectral parameters $\bm z_s \in \{z_s, \overline{z}_s\}^k$.
	It follows from Schwarz inequality that
	\begin{equation} \label{eq:QV_Q_est}
		\sum_{ab} S_{ab} \biggl\lvert \otherhat{\mathcal{P}}^k_{s,t\wedge\tau}\biggl[ \Tr\bigl[\partial_{ab}G_{[1,k],s}S^{x_k}\bigr]\biggr](\bm x) \biggr\rvert^2 \lesssim \sum_{j=1}^k \biggl(\sum_{\bm a\in\indset{N}^k} \bigl\lvert(\otherhat{\mathcal{P}}^k_{s,t\wedge\tau})_{\bm x\bm a}\bigr\rvert \sqrt{\QV^{(j)}_{k,s}(\bm z_s,\bm a)} \biggr)^2.
	\end{equation}
	Note that by \eqref{eq:ups_lower_bound}, \eqref{eq:sfunc_def} and $\ell_t\eta_t \le N$, $\mathfrak{s}^\mathrm{av}_{k,t\wedge\tau}(\bm x) \ge N^{-kD'-k/2}$.
	Therefore, to prove \eqref{eq:mart_bound}, it suffices to show that for any $\bm z_s \in \{z_s, \overline{z}_s\}^k$ and any $j \in \indset{k}$, the quantity $\QV^{(j)}_s(\bm z_s,\bm x)$, defined in \eqref{eq:Qkj}, satisfies
	\begin{equation} \label{eq:Q_est}
		\QV^{(j)}_{k,s}(\bm z_s,\bm x) \lesssim \bigl(\varphi_{k,s}^\mathrm{av,qv}\bigr)^2  (\ell_s\eta_s)^{2\beta_k}\frac{\bigl(\mathfrak{s}_{k,s}^\mathrm{av}(\bm x)\bigr)^2}{\eta_s}, \quad \bm x\in \indset{N}^k, \quad s \le \tau,
	\end{equation} 
	where $\varphi_{k,s}^\mathrm{av,qv}$ is the control parameter defined in \eqref{eq:phi_avqv}.
	
	We now proceed with the proof of \eqref{eq:Q_est}. Throughout this proof, we use the shorthand notation $\QV_k^{(j)}(\bm x) := \QV^{(j)}_{k,s}(\bm z_s, \bm x)$, $\eta := \eta_s$, $\ell := \ell_s$, and omit the explicit dependence of the control quantities $\mathfrak{s}_{\cdot,s}$ as well as the resolvent chains $G_{[\dots], s}$ on the time variable $s\le \tau$.
	
	Our goal is to express $\QV_{k}^{(j)}$ in terms of the resolvent chains of length at most $\maxK$, which can then be bounded using \eqref{eq:G_psi_bounds}, in a way that minimizes the number of $\psi$ factors. There will be 
	three strategies for such reduction\footnote{Reduction strategies are the cornerstones of proving local laws via
	master inequalities, they first appeared in \cite{Cipolloni2022Optimal} and 
	 have been systematically used later \cite{cipolloni2022rank, Cipolloni2022overlap, cipolloni2023eigenstate}. Their precise forms slightly varied but
	 they always followed one of these three strategies.}	
	that we briefly describe heuristically before proceeding with the technical proof.
	
	{\bf First reduction strategy.} While $\QV_{k}^{(j)}$ itself is not a resolvent chain, the definition \eqref{eq:Qkj} represents $\QV_{k}^{(j)}$ as a double sum of squares of resolvent chains of length $k+1$. Since typically $a\neq b$, \eqref{eq:G_psi_bounds} (ignoring the loss exponents $\alpha$) implies the symbolic bound
	\begin{equation} \label{eq:QG^2_bound}
		\QV_k = \sum_{ab} S_{ab} \bigl\lvert (G_{[k+1]})_{ab} \bigr\rvert^2 \lesssim (\pis{k+1})^2 \times \mathfrak{s}^2, \quad k+1 \le \maxK.
	\end{equation}
	Here $\mathfrak{s}$ represents a size function that we do not follow precisely for the moment as we focus on
	the structure of the $\psi$'s.
	
	{\bf Second reduction strategy.} 
	Note that $\QV_{k}^{(j)}$ can also be written as a sum of resolvent chains of length $2k+2$ (in two different ways), 
	\begin{equation} \label{eq:av_Q_rep}
		\begin{split}
			\QV_{k}^{(j)}(\bm x) &= \sum_b \bigl((G_{[1,j]})^*S^{x_k}(G_{[j,k]})^* S^b G_{[j,k]}S^{x_k} G_{[1,j]}\bigr)_{bb}\\
			&= \sum_a \bigl( G_{[j,k]}S^{x_k} G_{[1,j]} S^a (G_{[1,j]})^*S^{x_k}(G_{[j,k]})^*\bigr)_{aa}.
		\end{split}
	\end{equation}
	Symbolically, \eqref{eq:G_psi_bounds} applied to \eqref{eq:av_Q_rep} yields (ignoring the loss exponents $\alpha$) 
	\begin{equation} \label{eq:QG_bound}
		\QV_k = \sum_{b}   (G_{[k+1]}^*S^bG_{[k+1]})_{bb} =\sum_{b}   ( G_{[2k+2]})_{bb}   \lesssim \biggl(1 + \frac{\pis{2k+2}}{\sqrt{\ell\eta}}\biggr) \times \mathfrak{s}^2, \quad 2k+2\le\maxK.
	\end{equation}
	The advantage of 	\eqref{eq:QG_bound} over \eqref{eq:QG^2_bound} is that the right-hand side of \eqref{eq:QG_bound} contains only a single $\psi$, which is further suppressed by $(\ell\eta)^{-1/2}$, as compared to the two on the right-hand side of \eqref{eq:QG^2_bound}. However, the drawback is that \eqref{eq:QG_bound} is only applicable for $k\in \indset{\maxK/2-1}$.
	Neither \eqref{eq:QG^2_bound} nor \eqref{eq:QG_bound} is applicable when $k = \maxK$, since the chains involved are too long to apply \eqref{eq:G_psi_bounds}. Instead, for $k=\maxK$,  we need a third method:
	
	{\bf Third reduction strategy.} 
	Using a simple trace inequality $\Tr[AB] \le \Tr[A]\Tr[B]$ for positive-definite $N\times N$ matrices $A,B \ge 0$, to obtain the following  \emph{reduction inequality}. For all $N\times N$ matrices $X,Y$, we have
	\begin{equation} \label{eq:reduction}
		\begin{split}
			\bigl(X^* S^{x_q} Y^* S^a Y S^{x_q} X\bigr)_{bb} &= \Tr \bigl[(S^{x_q})^{1/2} Y^* S^a Y (S^{x_q})^{1/2}(S^{x_q})^{1/2} X \lvert b \rangle \langle b \rvert X^*(S^{x_q})^{1/2} \bigr]\\
			&\le \Tr \bigl[(S^{x_q})^{1/2} Y^* S^a Y (S^{x_q})^{1/2}\bigr]~\Tr \bigl[(S^{x_q})^{1/2} X \lvert b \rangle \langle b \rvert X^*(S^{x_q})^{1/2} \bigr]\\
			&= \Tr \bigl[Y^* S^a Y S^{x_q} \bigr]~ \bigl(X^* S^{x_q} X\bigr)_{bb}.
		\end{split}
	\end{equation}
	Hence, applying \eqref{eq:reduction} to the $2k+2$ chain on the right-hand side of \eqref{eq:av_Q_rep}, symbolically we obtain (ignoring the loss exponents $\alpha$ and $\beta$, see \eqref{eq:QK_final} below for more details)
	\begin{equation} \label{eq:Qred_bound}
		\QV_\maxK \le \sum_b \Tr \bigl[G_{[\maxK]} S^{b} \bigr]~ \bigl(G_{[\maxK+2]}\bigr)_{bb} \lesssim \sqrt{\ell\eta}\biggl(1 + \frac{\pav{\maxK}}{\ell\eta}\biggr) (\pis{\maxK/2+1})^2 \times \mathfrak{s}^2.
	\end{equation}
	 Note that the highest order control parameter 
	$\psi^\mathrm{av}_K$ comes with a small prefactor $1/\sqrt{\ell\eta}$ and to the first power, which,
	after taking the square root (owing to the quadratic variation) results in an estimate 
	$\psi^\mathrm{av}_K \le \sqrt{\psi^\mathrm{av}_K}/(\ell\eta)^{1/4} +\ldots$ which is closable assuming $\psi^\mathrm{iso}$
	for shorter chains is already controlled. 
	
	In conclusion, we will use the first strategy \eqref{eq:av_Q_rep} and \eqref{eq:QG_bound} for $k \in \indset{\maxK/2-1}$,
	the second strategy \eqref{eq:Qkj} and \eqref{eq:QG^2_bound} for $k\in\indset{\maxK/2,\maxK-1}$, and 
	the third strategy \eqref{eq:av_Q_rep}, \eqref{eq:reduction} with \eqref{eq:Qred_bound} to treat $k=\maxK$.
	After this heuristic sketch, we now start the actual proofs.

	First, we consider $k \in \indset{\maxK/2 -1}$. Since in this case $2k+2\le \maxK$, the representation \eqref{eq:av_Q_rep} is immediately useful. Using \eqref{eq:G_psi_bounds} for the $(2k+2)$-chain $(G_{\dots})_{bb}$ on the right-hand side in the first line of \eqref{eq:av_Q_rep}, we deduce from \eqref{eq:M_bound} and \eqref{eq:tau_def}  that, for all $j \in \indset{k}$,
	\begin{equation}\label{eq:bb_short}
		\bigl\lvert  \bigl((G_{[1,j]})^*S^{x_k}(G_{[j,k]})^* S^b G_{[j,k]}S^{x_k} G_{[1,j]}\bigr)_{bb} \bigr\rvert \lesssim \sqrt{\ell\eta}\biggl( 1+ \frac{\pis{2k+2}}{(\ell\eta)^{1/2-\alpha_{2k+2}}}\biggr) \mathfrak{s}_{2k+2}^\mathrm{iso}(b, \bm y_j^*, b, \bm y_j, b),
	\end{equation}
	where $\bm y_j := (x_j, \dots, x_k, x_1, \dots, x_{j-1})$ denotes a cyclic shift of  $\bm x$, and $\bm y_j^* := (x_{j-1}, \dots, x_1, x_k, \dots, x_j)$ is obtained by reversing the order of the elements in  $\bm y_j$.
	It follows from \eqref{eq:sfunc_def} that for $\bm x^* := (x_k, x_{k-1}, \dots, x_1)$,
	\begin{equation} \label{eq:sfunc_doubling}
		\mathfrak{s}_{2k+2}^\mathrm{iso}(b,\bm x^*, a, \bm x, b) =  (\ell\eta)^{-1/2}\bigl(\mathfrak{s}_{k+1}^\mathrm{iso}(a,\bm x ,b)\bigr)^2, \quad k \in \indset{\maxK}.
	\end{equation}	
	Therefore, by \eqref{eq:av_Q_rep}, \eqref{eq:bb_short}  
	and \eqref{eq:sfunc_doubling}, we conclude that
	\begin{equation} \label{eq:Qk_short} 
			\QV_k^{(j)}(\bm x)  
			\lesssim \sum_b \bigl(\mathfrak{s}^\mathrm{iso}_{k+1}(b,\bm y_j, b)\bigr)^2\biggl(1 + \frac{\pis{2k+2}}{(\ell\eta)^{1/2-\alpha_{2k+2}}}\biggr)  \lesssim \frac{\bigl(\mathfrak{s}_{k}^\mathrm{av}(\bm x)\bigr)^2}{\eta}\biggl(1 + \frac{\pis{2k+2}}{(\ell\eta)^{1/2 -\alpha_{2k+2}}} \biggr),
	\end{equation}
	where in the second step we used \eqref{eq:sfunc_def} and the convolution inequality \eqref{eq:true_convol}. 
	Estimate \eqref{eq:Qk_short} implies \eqref{eq:Q_est} for $k \in\indset{\maxK/2-1}$ by \eqref{eq:phi_avqv}.

	Next, we consider $k \in \indset{\maxK/2, \maxK -1}$. Since in this case $2k+2 > \maxK$ but $k+1 \le \maxK$, we go back to the definition \eqref{eq:Qkj}, that expresses $\QV_k^{(j)}(\bm x)$ in terms of $(k+1)$-chain $(G_{\dots})_{ab}$. Similarly to \eqref{eq:Qk_short}, by using \eqref{eq:G_psi_bounds}, we deduce that
	\begin{equation} \label{eq:Qk_long}
		\begin{split}
			\QV_k^{(j)}(\bm x) 
			&\lesssim \bigl(\pis{k+1}(\ell\eta)^{\alpha_{k+1}} \bigr)^2 \sum_{ab} S_{ab}  \bigl(\mathfrak{s}_{k+1}^\mathrm{iso}(a,\bm x, b)\bigr)^2 
			+ \ell\eta\sum_{b}S_{bb}\bigl(\mathfrak{s}_{k+1}^\mathrm{iso}(b,\bm x, b)\bigr)^2\\
			&\lesssim\frac{\bigl(\mathfrak{s}_{k}^\mathrm{av}(\bm x)\bigr)^2}{\eta}(\ell\eta)^{2\alpha_{k+1}}\bigl(\pis{k+1}\bigr)^2,
		\end{split}
	\end{equation}
	where in the last step  we used the local contraction bound \eqref{eq:iso_S}. 	
	 Estimate \eqref{eq:Qk_long} together with \eqref{eq:phi_avqv} implies \eqref{eq:Q_est} for $k \in \indset{\maxK/2, \maxK -1}$.
	
	Finally, we consider the case $k=\maxK$.  
	We apply the reduction inequality \eqref{eq:reduction} to the $(2k+2)$-chain on the right-hand side of \eqref{eq:av_Q_rep} by choosing the splitting index $q$ in such a way that $Y$ is a resolvent chain of length $\maxK/2$, and $X$ is a resolvent chain of length $\maxK/2 + 1$. The precise choice of $q$ depends on the index $j$ of $\QV_K^{(j)}$. 
	For $j \in \indset{\maxK/2+1}$, denoting $q \equiv q(j) := \maxK/2+j-1$, we conclude that, for all $j \in \indset{\maxK/2+1}$,
	\begin{equation} \label{eq:QK_intermediate}
		\begin{split}
			\QV_K^{(j)}(\bm x) &\le \sum_b \Tr \bigl[G_{[j,q]}S^{x_{q}}(G_{[j,q]})^* S^b\bigr]\bigl((G_{[1,j]})^*S^{x_\maxK}(G_{[q+1,\maxK]})^* S^{x_{q}} G_{[q+1,\maxK]}S^{x_\maxK}G_{[1,j]}\bigr)_{bb} \\
			&=  \sum_b \Tr \bigl[G_{[j,q]}S^{x_{q}}(G_{[j,q]})^* S^b\bigr] \sum_a S_{ax_{q}} \bigl\lvert\bigl(G_{[q+1,\maxK]}S^{x_\maxK}G_{[1,j]}\bigr)_{ab}\bigr\rvert^2,
		\end{split}
	\end{equation}
	We now bound the terms appearing on the right-hand side of \eqref{eq:QK_intermediate}.
	Using \eqref{eq:sfunc_def}, \eqref{eq:G_psi_bounds} and \eqref{eq:sfunc_doubling}, for every $b \in \indset{N}$, we obtain the estimate 
	\begin{equation} \label{eq:QK_intermediate1}
		\begin{split}
			\Tr \bigl[G_{[j,q]}S^{x_{q}}G_{[j,q]}^* S^b\bigr] &\lesssim  \bigl(\ell\eta + (\ell\eta)^{\beta_\maxK}\pav{\maxK}\bigr)\mathfrak{s}^\mathrm{av}_{\maxK}(x_j, \dots,x_{q-1}, x_q, x_{q-1}, \dots, x_j, b)\\
			&\lesssim \bigl(\mathfrak{s}^\mathrm{iso}_{\maxK/2}(b, x_j, \dots x_{q-1}, x_q)\bigr)^2 \biggl(1 + \frac{\pav{\maxK}}{(\ell\eta)^{1-\beta_\maxK}}\biggr).
		\end{split}
	\end{equation}
	On the other hand, similarly to \eqref{eq:Qk_long}, using \eqref{eq:SUps_comvol}, we deduce from \eqref{eq:G_psi_bounds} that
	\begin{equation} \label{eq:QK_intermediate1_1}
		\begin{split}
			\sum_a S_{ax_{q}} \biggl\lvert\bigl(G_{[q+1,\maxK]}S^{x_\maxK}G_{[1,j]}\bigr)_{ab}\biggr\rvert^2 \lesssim&~ \ell\eta S_{bx_q} \bigl( \mathfrak{s}^\mathrm{iso}_{\maxK/2+1}(b,\bm x_{[q+1,j-1]},b)\bigr)^2\\
			&+ 
			\bigl((\ell\eta)^{\alpha_{\maxK/2+1}} \pis{\maxK/2+1}\, \mathfrak{s}^\mathrm{iso}_{\maxK/2+1}(x_{q}, \bm x_{[q+1,j-1]},b)\bigr)^2,
		\end{split}
	\end{equation}
	where we use the shortcut notation $\bm x_{[q+1,j-1]} := (x_{q+1},\dots,x_K,x_1,\dots, x_{j-1}) \in \indset{N}^{\maxK/2}$.
	
	To combine the estimates \eqref{eq:QK_intermediate1} and \eqref{eq:QK_intermediate1_1}, we observe that, by definition of $\mathfrak{s}^\mathrm{iso}$ in \eqref{eq:sfunc_def}, and  \eqref{eq:triag}, \eqref{eq:true_convol}, \eqref{eq:SUps_comvol}, we have the following local contraction bounds
	\begin{equation} \label{eq:QK_diag_term}
		\sum_b \bigl(\mathfrak{s}^\mathrm{iso}_{\maxK/2}(b, x_j, \dots x_{q-1}, x_q)\bigr)^2 S_{bx_q} \bigl( \mathfrak{s}^\mathrm{iso}_{\maxK/2+1}(b,\bm x_{[q+1,j-1]},b)\bigr)^2 \lesssim  \bigl(\mathfrak{s}_{\maxK}^\mathrm{av}(\bm x)\bigr)^2,
	\end{equation}
	\begin{equation} \label{eq:QK_offdiag_term}
		\sum_b \bigl(\mathfrak{s}^\mathrm{iso}_{\maxK/2}(b, x_j, \dots x_{q-1}, x_q)\bigr)^2 \bigl(\mathfrak{s}^\mathrm{iso}_{\maxK/2+1}(x_{q}, \bm x_{[q+1,j-1]},b)\bigr)^2 \lesssim  \ell\bigl(\mathfrak{s}_{\maxK}^\mathrm{av}(\bm x)\bigr)^2.
	\end{equation}
	Therefore, we conclude from \eqref{eq:QK_intermediate}, \eqref{eq:QK_diag_term} and \eqref{eq:QK_offdiag_term} that
	\begin{equation} \label{eq:QK_final}
		\QV_\maxK^{(j)}(\bm x) \le \frac{\bigl(\mathfrak{s}_{\maxK}^\mathrm{av}(\bm x)\bigr)^2}{\eta} (\ell\eta)^{1 + 2\alpha_{\maxK/2+1}}\biggl(1 + \frac{\pav{\maxK}}{(\ell\eta)^{1-\beta_\maxK}}\biggr)\bigl(\pis{\maxK/2+1} \bigr)^2, \quad j \in \indset{\maxK/2+1},
	\end{equation}
	which implies \eqref{eq:Q_est} for $k=\maxK$ and $j \in \indset{\maxK/2+1}$ by \eqref{eq:phi_avqv}.
	Completely analogously, for $j \in \indset{\maxK/2+2, \maxK}$, we perform the reduction by cutting the resolvent chain at $q \equiv q(j) :=  j - \maxK/2 $ to obtain
	\begin{equation} \label{eq:QK_intermediate2}
		\QV_\maxK^{(j)}(\bm x) \le \sum_b \Tr \biggl[G_{[j,q]}^{(x_\maxK)} S^{x_{q}}\bigl(G_{[j,q]}^{(x_\maxK)} \bigr)^* S^b\biggr] \sum_a S_{ax_{q}} \bigl\lvert\bigl(G_{[q,j]}\bigr)_{ab}\bigr\rvert^2, \quad j \in  \indset{\maxK/2+2, \maxK}.
	\end{equation}
	Therefore, estimate \eqref{eq:QK_final} also holds for $j \in  \indset{\maxK/2+2, \maxK}$. This concludes the proof of \eqref{eq:Q_est}, and hence that of \eqref{eq:mart_bound}.
	
	\vspace{5pt}
	\noindent\textbf{Isotropic Chains}. Proof of \eqref{eq:iso_mart_bound}. Recall that $k \in \indset{\maxK-1}$.
	As in the averaged case, it follows from \eqref{eq:mart_ineq}
	that 
	\begin{equation}
		\begin{split}
			\biggl\lvert \int_{\tinit}^{t\wedge\tau} \otherhat{\mathcal{P}}^{k}_{s,t\wedge\tau}\bigl[\mathrm{d}\mathcal{M}^\mathrm{iso}_{[1,k+1],s}\bigr](\bm x) \biggr\rvert
			\prec&~ \biggl(\int_{\tinit}^{t\wedge\tau} \sum_{cd} S_{cd} \biggl\lvert \otherhat{\mathcal{P}}^k_{s,t\wedge\tau}\biggl[ \partial_{cd}\bigl(G_{[1,k+1],s}\bigr)_{ab} \biggr](\bm x) \biggr\rvert^2\mathrm{d}s\biggr)^{1/2}\\ &+ N^{-(k+1)D'- k/2},
		\end{split}
	\end{equation}
	where $D'$ is the constant from \eqref{eq:ups_lower_bound}, and we abbreviate $\mathrm{d}\mathcal{M}^\mathrm{iso}_{[1,k+1],s}(\bm x) \equiv \mathrm{d}\mathcal{M}^\mathrm{iso}_{[1,k+1],s}(a,\bm x,b)$
	as $a, b$ are fixed.
	
	Similarly to \eqref{eq:Qkj}, for all $j \in \indset{k+1}$, we define the quantity $\QV^{(j)}_{k+1,s}(\bm z_s;a,\bm x,b) $ as
	\begin{equation} \label{eq:iso_Qkj}
		\QV^{(j)}_{k+1,s}(\bm z_s;a,\bm x,b) :=  \sum_{cd} S_{cd} \bigl\lvert \bigr(G_{[1,j],s}\bigr)_{ad} \bigl(G_{[j,k+1],s}\bigr)_{cb}  \bigr\rvert^2, \quad j\in\indset{k+1}.
	\end{equation}
	Analogously to \eqref{eq:QV_Q_est}, Schwarz inequality implies that 
	\begin{equation} \label{eq:iso_QV_Q_est}
		\sum_{cd} S_{cd} \biggl\lvert \otherhat{\mathcal{P}}^{k}_{s,t}\biggl[ \partial_{cd}\bigl(G_{[1,k+1],s}\bigr)_{ab}\biggr](\bm x) \biggr\rvert^2 \lesssim \sum_{j=1}^{k} \biggl(\sum_{\bm c\in\indset{N}^{k}} \bigl\lvert(\otherhat{\mathcal{P}}^{k}_{s,t})_{\bm x\bm c}\bigr\rvert \sqrt{\QV^{(j)}_{k+1,s}(\bm z_s;a,\bm c,b)} \biggr)^2.
	\end{equation}
	Since by \eqref{eq:ups_lower_bound} and \eqref{eq:sfunc_def}, $\mathfrak{s}^\mathrm{iso}_{k+1,t\wedge\tau}(a,\bm x,b) \ge N^{-(k+1)D'-k/2}$, to prove \eqref{eq:iso_mart_bound}, it suffices to show that for any $\bm z_s \in \{z_s, \overline{z}_s\}^{k+1}$, all $a,b\in\indset{N}$, the quantities $\QV^{(j)}_{k+1,s}(\bm z_s;a,\bm x',b)$, defined in \eqref{eq:iso_Qkj}, satisfy
	\begin{equation} \label{eq:isoQ_est}
		\QV^{(j)}_{k+1,s}(\bm z_s;a,\bm x,b) \lesssim \bigl(\varphi_{k+1,s}^\mathrm{iso,qv}\bigr)^2(\ell_s\eta_s)^{2\alpha_{k+1}}\frac{\bigl(\mathfrak{s}_{k+1,s}^\mathrm{iso}(a,\bm x,b)\bigr)^2}{\eta_s}, \quad \bm x\in\indset{N}^{k}, \quad s\le\tau,
	\end{equation}
	where $\varphi_{k+1,s}^\mathrm{iso,qv}$ is defined in \eqref{eq:phi_isoqv}.
	As in the proof of \eqref{eq:Q_est} above, we use the shortcut notation $\QV^{(j)}_{k+1}(a,\bm x,b) := \QV^{(j)}_{k+1,s}(\bm z_s;a,\bm x,b)$ and omit the dependence of $\ell, \mathfrak{s}, \eta, m, \Theta, G_{\dots}$ on $s \le \tau$. 
	
	Due to symmetry, it suffices to consider $j \le  \lfloor (k+2)/2 \rfloor$. Since $\maxK$ is even, for all $j \in \indset{k+1}$ and $k \in \indset{\maxK-1}$ satisfying $j \le  (k+1)/2$, we also have $2j \le \maxK$, hence it follows from  \eqref{eq:G_psi_bounds} and \eqref{eq:sfunc_doubling} that
	\begin{equation}
		\sum_{d} S_{cd} \bigl\lvert \bigr(G_{[1,j]}\bigr)_{ad} \bigr\rvert^2 = \bigr((G_{[1,j]})^* S^cG_{[1,j]}\bigr)_{aa} \lesssim \biggl(1+ \frac{\pis{2j}}{(\ell\eta)^{1/2-\alpha_{2j}}} \biggr)  \bigl(\mathfrak{s}_{j}^\mathrm{iso}(a,x_1,\dots, x_{j-1},c)\bigr)^2.
	\end{equation}
	Therefore, for $j \in \indset{2, \lfloor (k+1)/2 \rfloor}$, it follows from \eqref{eq:G_psi_bounds} applied to $(G_{[j,k+1]})_{cb}$ in \eqref{eq:iso_Qkj} that
	\begin{equation} \label{eq:isoQj_est}
		\begin{split}
			\QV^{(j)}_{k+1}(a,\bm x,b) \lesssim&~ \sum_{c}\biggl(1+ \frac{\pis{2j}}{(\ell\eta)^{1/2-\alpha_{2j}}} \biggr)  \bigl(\mathfrak{s}_{j}^\mathrm{iso}(a,x_1,\dots, x_{j-1},c)\bigr)^2\\
			&\qquad \times \biggl(\delta_{bc}(\ell\eta)^{1/2} + \pis{k-j+2}(\ell\eta)^{\alpha_{k-j+2}} \biggr)^2\bigl(\mathfrak{s}_{k-j+2}^\mathrm{iso}(c,x_j,\dots, x_{k},b)\bigr)^2\\
			\lesssim&~ \biggl(1+ \frac{\pis{2j}}{(\ell\eta)^{1/2-\alpha_{2j}}} \biggr) \biggl(\frac{\pis{k-j+2}}{(\ell\eta)^{\alpha_{k+1}-\alpha_{k-j+2}}} \biggr)^2 \frac{\bigl((\ell\eta)^{\alpha_{k+1}}\mathfrak{s}_{k+1}^\mathrm{iso}(a,\bm x,b)\bigr)^2}{\eta},
		\end{split}
	\end{equation}
	where in the second step  we used \eqref{eq:Schwarz_convol} (or \eqref{eq:true_convol}),   similarly to \eqref{eq:av_insertion}. Note that by \eqref{eq:phi_isoqv},  \eqref{eq:isoQj_est} implies \eqref{eq:isoQ_est} for all indices $j \in \indset{2, \lfloor (k+1)/2 \rfloor}$.
	
	It remains to estimate $\QV^{(1)}_{k+1}(a,\bm x,b)$. We begin by performing the following decomposition,
	\begin{equation} \label{eq:isoQ_decomp}
		\begin{split}
			\QV^{(1)}_{k+1}(a,\bm x,b) = \sum_{c} \bigr(G_{1} S^c G_{1}^* -\Theta^c\bigr)_{aa}  \bigl\lvert \bigl(G_{[1,k+1]}\bigr)_{cb}  \bigr\rvert^2 + \sum_{d} (I+\Theta)_{ad} ~\other{\QV}_{k+1}(d,\bm x,b),
		\end{split}
	\end{equation}
	where we subtracted $\Theta^c$ from $G_{1} S^c G_{1}^* $, added it back by using
	 the identity $\Theta = |m|^2(I+\Theta)S$ and introduced the auxiliary quantity
	\begin{equation} \label{eq:otherQ}
		\other{\QV}_{k+1}(d,\bm x,b) := |m|^2\sum_c S_{dc}  \bigl\lvert \bigl(G_{[1,k+1]}\bigr)_{cb}  \bigr\rvert^2.
	\end{equation}
	It follows from \eqref{eq:G-M_psi_bounds} and \eqref{eq:G_psi_bounds} that the first term on the right-hand side of \eqref{eq:isoQ_decomp} admits the bound
	\begin{equation} \label{eq:Q1_subleading}
		\biggl\lvert\sum_{c} \bigr(G_{1} S^c G_{1}^* -\Theta^c\bigr)_{aa}  \bigl\lvert \bigl(G_{[1,k+1]}\bigr)_{cb}  \bigr\rvert^2\biggr\rvert \lesssim \frac{\pis{2}\bigl(\pis{k+1}\bigr)^2}{\sqrt{\ell\eta}}\frac{\bigl((\ell\eta)^{\alpha_{k+1}}\mathfrak{s}_{k+1}^\mathrm{iso}(a,\bm x,b)\bigr)^2}{\eta}.
	\end{equation}
	As to the second term on the right-hand side of \eqref{eq:isoQ_decomp} the quantity $\other{\QV}_{k+1}$, defined in \eqref{eq:otherQ}, requires different treatment depending on the relative size of  $k$  and $\maxK$, similarly to the proof of \eqref{eq:Q_est}.   
	
	First, consider the case $k \in \indset{\maxK/2-1}$. In this case $2k+2 \le \maxK$, hence, we view $\other{\QV}_{k+1}$ as $(2k+2)$-chain and use \eqref{eq:sfunc_def}, \eqref{eq:G_psi_bounds}, \eqref{eq:sfunc_doubling} to obtain, for all $k \in \indset{\maxK/2-1}$, 
	\begin{equation} \label{eq:otherQ_short}
		\begin{split}
			\other{\QV}_{k+1}(d,\bm x, b) &=  |m|^2 \bigl((G_{[1,k+1]})^* S^d G_{[1,k+1]}\bigr)_{bb}  \lesssim \biggl(1+ \frac{\pis{2k+2}}{(\ell\eta)^{1/2-\alpha_{2k+2}}} \biggr)  \bigl(\mathfrak{s}_{k+1}^\mathrm{iso}(d,\bm x,b)\bigr)^2. 
		\end{split}
	\end{equation} 
	
	Next, we consider the case $k\in\indset{\maxK/2,\maxK-2}$. In this case the $(2k+2)$-chain is too long to be estimated directly, but can be viewed as a combination of $k$- and a $(k+2)$-chains,
	\begin{equation} \label{eq:otherQ_long}
		\begin{split}
			\other{\QV}_{k}(d,\bm x,b)  &\lesssim   \sum_q S_{x_1q} \bigl\lvert\bigl(G_{[2,k]}\bigr)_{qb} \bigl(G_{1}^*S^dG_{[1,k]}\bigr)_{qb}\bigr\rvert  
			\lesssim \frac{\pis{k}\pis{k+2}}{(\ell\eta)^{2\alpha_{k+1}- \alpha_{k}-\alpha_{k+2}} } \bigl((\ell\eta)^{\alpha_{k+1}}\mathfrak{s}_{k+1}^\mathrm{iso}(d,\bm x,b)\bigr)^2, 
		\end{split}
	\end{equation}
	for all $k\in\indset{\maxK/2,\maxK-2}$, where we used local contraction argument
	based upon  \eqref{eq:triag}, \eqref{eq:sfunc_def}, \eqref{eq:G_psi_bounds}, and \eqref{eq:SUps_comvol}.

	Finally, we consider the case $k = \maxK-1$. In this case, similarly to \eqref{eq:QK_intermediate}, we need to perform a reduction. Letting $q:=\maxK/2$, it follows from \eqref{eq:reduction} that
	\begin{equation}
		\begin{split}
			\other{\QV}_{\maxK}(d,\bm x,b)  &= |m|^2 \bigl((G_{[1,\maxK]})^* S^d G_{[1,\maxK]}\bigr)_{bb}\\
			&\le   \Tr \bigl[(G_{[1,q]})^* S^d G_{[1,q]} S^{x_q} \bigr]  \bigl( (G_{[q,\maxK]})^* S^{x_q} G_{[q,\maxK]}\bigr)_{bb}\\
			&\le   \Tr \bigl[(G_{[1,q]})^* S^d G_{[1,q]} S^{x_q} \bigr] \sum_c S_{x_qc} \bigl\lvert\bigl(G_{[q,\maxK]}\bigr)_{cb}\bigr\rvert^2.
		\end{split}
	\end{equation}
	Similarly to \eqref{eq:QK_final}, using \eqref{eq:G_psi_bounds} and recalling that $\alpha_\maxK = 1/2$ from \eqref{eq:loss_exponents}, we obtain
	\begin{equation} \label{eq:otherQ_K}
			\other{\QV}_{\maxK}(d,\bm x,b)  \lesssim  \bigl(\pis{\maxK/2}\bigr)^2\biggl(1+\frac{\pav{\maxK}}{(\ell\eta)^{1-\beta_\maxK}}\biggr) \bigl((\ell\eta)^{\alpha_{\maxK}}\mathfrak{s}_{k}^\mathrm{iso}(d,\bm x',b)\bigr)^2. 
	\end{equation}

	It follows from \eqref{eq:true_convol}, \eqref{eq:isoQ_decomp}, \eqref{eq:Q1_subleading} the estimates \eqref{eq:otherQ_short}, \eqref{eq:otherQ_long} and \eqref{eq:otherQ_K}, that
	\begin{equation} \label{eq:Q1_est}
		\QV^{(1)}_{k+1}(a,\bm x,b) \lesssim \biggl(\frac{\pis{2}\bigl(\pis{k+1}\bigr)^2}{\sqrt{\ell\eta}} + \bigl(\other{\varphi}_{k+1}^\mathrm{\,iso,qv}\bigr)^2\biggr)\frac{\bigl((\ell\eta)^{\alpha_{k+1}}\mathfrak{s}_{k+1}^\mathrm{iso}(a,\bm x,b)\bigr)^2}{\eta},
	\end{equation} 
	where $\other{\varphi}_{k+1}^\mathrm{\,iso,qv} := \other{\varphi}_{k+1,s}^\mathrm{\,iso,qv}$ is defined in \eqref{eq:other_phi_isoqv}.
	Combining \eqref{eq:phi_isoqv}, \eqref{eq:isoQj_est} and \eqref{eq:Q1_est}, we obtain the desired \eqref{eq:isoQ_est}. 
	This concludes the proof of Lemma \ref{lemma:mart_est}.
\end{proof}

\subsubsection{Forcing term estimates} \label{sec:forcing}
 Before proceeding with the proof of Lemma \ref{lemma:forcing}, we record an additional bound on a special class of deterministic $M$-terms that is used along the proof. 
 \begin{lemma}\label{lemma:extra_M}
	Fix $k \in \mathbb{N}$ and let $\bm z_t \in \{z_t, \overline{z}_t\}^k$, $\bm x\in \indset{N}^k$, then, for any $j \in \indset{k}$,
	\begin{equation} \label{eq:resum_M_bound}
		\sum_q \bigl\lvert \big(M_{[j,j],t}^{(x_k)}\big)_{qq} \bigr\rvert \lesssim (\log N)^{C_k'}\frac{1}{\eta_t}\times\ell_t\eta_t\, \mathfrak{s}_{k,t}^{\mathrm{av}}(\bm x),
	\end{equation}
	for some positive $k$-dependent constant $C_k'$, where we recall that for a fixed $k \in \mathbb{N}$,  and any $j \in \indset{k}$,
	\begin{equation} \label{eq:special_M}
		M_{[j,j],t}^{(x_k)} \equiv M_{[j,j],t}^{(x_k)}(\bm x') := M\bigl(z_{j,t}, S^{x_j}, \dots, z_{k,t}, S^{x_{k}}, z_{1,t}, S^{x_1}, \dots, z_{j-1,t}, S^{x_{j-1}}, z_{j,t} \bigr).
	\end{equation}
\end{lemma}
In the sequel, we safely ignore the $\log N$ in \eqref{eq:resum_M_bound}. We prove Lemma \ref{lemma:extra_M} is Section \ref{sec:M_bounds}. 

Note that Lemma~\ref{lemma:extra_M} is essentially a deterministic analog of a Ward estimate. For a heuristic, consider that
$$
\Tr \big(M_{[j,j],t}^{(x_k)}\big) =\sum_q \big(M_{[j,j],t}^{(x_k)}\big)_{qq},
$$
approximates $\Tr [G_j S^{x_j} G_{j+1} \ldots  S^{x_{j-1}} G_j]\sim \frac{1}{\eta_{j,t}}
 \Tr [G_j S^{x_j} G_{j+1} \ldots  G_{j-1}S^{x_{j-1}}]$, where we replaced $G_j^2$ by $G_j/\eta_j$, 
  ignoring the difference between $G_j$ and $G_j^*$.

We proceed to prove Lemma \ref{lemma:forcing}.  
\begin{proof}[Proof of Lemma \ref{lemma:forcing}]
	Throughout the proof, we suppress the explicit dependence of the control quantities $\mathfrak{s}_{\cdot,s}$ as well as resolvent chains $G_{[\dots], s}$ and the corresponding deterministic approximations $M_{[\dots], s}$ on the time variable $s$ (or $t$ when we refer to~\eqref{eq:F1_def}--\eqref{eq:Fk_def}).
	
	\vspace{5pt}
	\noindent\textbf{Averaged Chains}. Proof of \eqref{eq:av_forcing_bound}.
	
	First, for $k \ge 2$, we consider the terms of the form (first line in \eqref{eq:Fk_def}) 
	\begin{equation} \label{eq:1_k+1_terms}
		\Tr\bigl[\mathcal{S}[G_{j}-m_{j}] G_{[j,k]}S^{x_k}G_{[1,j]}\bigr] = \sum_{q}\Tr\bigl[(G_{j}-m_{j})S^q\bigl] \bigl(G_{[j,k]}S^{x_k}G_{[1,j]}\bigr)_{qq}, \quad j \in \indset{k}.
	\end{equation}
	Note that the resolvent chain $G_{[j,k]}S^{x_k}G_{[1,j]}$ has length $k+1$. 
	
	In the case $k \in \indset{\maxK-1}$, the $(k+1)$-chain can be estimated directly.  Subtracting the corresponding deterministic approximation, denoted by $M_{[j,j]}^{(x_k)}$, we obtain,
	\begin{equation}\label{eq:1_k+1_bound1_1}
		\begin{split}
			\bigl\vert \Tr\bigl[\mathcal{S}[G_{j}-m_{j}] \bigl(G_{[j,k]}S^{x_k}G_{[1,j]} - M_{[j,j]}^{(x_k)}\bigr)\bigr] \bigr\rvert &= \sum_q \Tr\bigl[(G_{j}-m_{j})S^q\bigr]\bigl(G_{[j,k]}S^{x_k}G_{[1,j]} - M_{[j,j]}^{(x_k)}\bigr)_{qq}\\
			&\lesssim \frac{\pav{1}}{\ell\eta}    \sum_q (\ell\eta)^{\alpha_{k+1}}\pis{k+1}\mathfrak{s}_{k+1}^\mathrm{iso}(q,\bm x,q) \\
			&\lesssim \frac{(\ell\eta)^{\beta_{k}}\mathfrak{s}_{k}^\mathrm{av}(\bm x)}{\eta}\frac{\sqrt{N\eta}}{\ell\eta}\pav{1}\pis{k+1}  \\ 
		\end{split}
	\end{equation}
	where in the first inequality we used \eqref{eq:G-M_psi_bounds}, $\mathfrak{s}_{1}^\mathrm{av}=1/(\ell\eta)$, 
	and in the second step we used \eqref{eq:sqrt_convol}, similarly to \eqref{eq:iso_convol}. 
	Recall also that $\alpha_{k+1}=\beta_k$ for $k\le K-1$, see \eqref{eq:loss_exponents}. 
	 On the other hand, by \eqref{eq:G-M_psi_bounds} and \eqref{eq:resum_M_bound}, the contribution from $M_{[j,j]}^{(x_k)}$ admits the bound
	\begin{equation}\label{eq:1_k+1_bound1_2}
		\bigl\vert \Tr\bigl[\mathcal{S}[G_{j}-m_{j}] M_{[j,j]}^{(x_k)} \bigr] \bigr\rvert = \sum_q \Tr\bigl[(G_{j}-m_{j})S^q\bigr]\bigl( M_{[j,j]}^{(x_k)}\bigr)_{qq}  
		\lesssim \frac{\mathfrak{s}_{k}^\mathrm{av}(\bm x)}{\eta} \pav{1}.
	\end{equation}
	Therefore, combining \eqref{eq:1_k+1_bound1_1} and \eqref{eq:1_k+1_bound1_2}, we obtain
	\begin{equation}\label{eq:1_k+1_bound1}
		\bigl\vert (\ref{eq:1_k+1_terms}) \bigr\rvert \lesssim \frac{(\ell\eta)^{\beta_{k}}\mathfrak{s}_{k}^\mathrm{av}(\bm x)}{\eta}
		\biggl( \frac{\sqrt{N\eta}}{\ell\eta} \pav{1}\pis{k+1}  + \frac{\pav{1}}{(\ell\eta)^{\beta_{k}}}\biggr).
	\end{equation}

	In the case $k=\maxK$, we cut the $(\maxK+1)$-chain $G_{[j,\maxK]}S^{x_\maxK}G_{[1,j]}$ into two chains of length $\maxK/2$ and $\maxK/2+1$. Similarly to the estimate on the quadratic variation terms in \eqref{eq:QK_final}, we need to consider two cases: First, if $j \in\indset{\maxK/2 +1}$, then we define the cutting index $q \equiv q(j) := \maxK/2+j-1$, to obtain, by \eqref{eq:G-M_psi_bounds} and \eqref{eq:G_psi_bounds}, that
	\begin{equation} \label{eq:1_k+1_reduction1}
		\begin{split}
			\bigl\lvert\Tr\bigl[\mathcal{S}[G_{j}-m_{j}] G_{[j,k]}S^{x_\maxK}G_{[1,j]}\bigr] \bigr\rvert \lesssim \frac{\pav{1}}{\ell\eta}\sum_{ab} S_{ax_q}\bigl\lvert \bigl(G_{[q+1,j]}^{(x_\maxK)}\bigr)_{ab} \bigl(G_{[j,q]}\bigr)_{ba} \bigr\rvert, \quad j \in\indset{\maxK/2 +1}.
		\end{split}
	\end{equation}
	Therefore, for $k=K$ and  $j \in\indset{\maxK/2 +1}$,  it follows from \eqref{eq:G-M_psi_bounds} and \eqref{eq:G_psi_bounds} and
	local contraction arguments, similarly to \eqref{eq:iso_convol} and \eqref{eq:av_insertion}, that 
	\begin{equation} \label{eq:1_k+1_bound2}
		\begin{split}
			\bigl\vert (\ref{eq:1_k+1_terms}) \bigr\rvert \lesssim&~ \pav{1} \biggl(1+ \frac{\pis{\maxK/2}+\pis{\maxK/2+1}}{(\ell\eta)^{1/2-\alpha_{\maxK/2+1}}}\biggr)
			\sum_{b} S_{bx_q} \mathfrak{s}_{\maxK/2}^\mathrm{iso}(b,x_{j},\dots, x_{q-1},b)  \mathfrak{s}_{\maxK/2+1}^\mathrm{iso}(b,\bm x_{[q+1,j-1]},b)\\
			&+\frac{\pav{1}\pis{\maxK/2}\pis{\maxK/2+1}}{(\ell\eta)^{1-\alpha_{\maxK/2}-\alpha_{\maxK/2+1}}} \sum_{ab} S_{ax_q} \mathfrak{s}_{\maxK/2}^\mathrm{iso}(b,x_{j},\dots, x_{q-1},a)\mathfrak{s}_{\maxK/2+1}^\mathrm{iso}(a,\bm x_{[q+1,j-1]},b)\\
			\lesssim&~ \frac{(\ell\eta)^{\beta_\maxK}\mathfrak{s}_{\maxK}^\mathrm{av}(\bm x)}{\eta} 
			\frac{\sqrt{N\eta}}{\ell\eta} \pav{1}\pis{\maxK/2}\pis{\maxK/2+1} ,
		\end{split} 
	\end{equation}
	 
	where we use the shortcut notation $\bm x_{[q+1,j-1]} := (x_{q+1},\dots,x_K,x_1,\dots, x_{j-1}) \in \indset{N}^{\maxK/2}$.
	
	On the other hand, for  $k=K$ and $j \in\indset{\maxK/2 +2, \maxK}$, we define $q \equiv q(j) :=  j- \maxK/2$. 
	Hence, completely analogously to \eqref{eq:1_k+1_reduction1} and \eqref{eq:1_k+1_bound2}, for $j \in\indset{\maxK/2 +2, \maxK}$,
	\begin{equation} \label{eq:1_k+1_bound3}
		\begin{split}
			\bigl\lvert\Tr\bigl[\mathcal{S}[G_{j}-m_{j}] G_{[j,\maxK]}S^{x_\maxK}G_{[1,j]}\bigr] \bigr\rvert &\lesssim \frac{\pav{1}}{\ell\eta}\sum_{ab} S_{ax_q}\bigl\lvert \bigl(G_{[q+1,j]}\bigr)_{ab} \bigl(G_{[j,q]}^{(x_\maxK)}\bigr)_{ba} \bigr\rvert\\
			&\lesssim  
			\frac{(\ell\eta)^{\beta_\maxK}\mathfrak{s}_{\maxK}^\mathrm{av}(\bm x)}{\eta} 
			\frac{\sqrt{N\eta}}{\ell\eta} \pav{1}\pis{\maxK/2}\pis{\maxK/2+1}. 
		\end{split}
	\end{equation} 
	holds for all $j \in \indset{\maxK}$. 
	
	Next, we estimate the quadratic terms (second line of \eqref{eq:Fk_def}), that is (replacing $j \mapsto i+j-1$)
	\begin{equation} \label{eq:quad_terms}
		\Tr\bigl[ \mathcal{S}\bigl[(G-M)_{[i,i+j-1]}\bigr] (G-M)^{(x_k)}_{[i+j-1,i]} \bigr], \quad i\in \indset{k}, \quad j \in \indset{2, k+1-i}.
	\end{equation}
	Using \eqref{eq:G-M_psi_bounds}, and the local contraction~\eqref{eq:av_insertion}, we conclude that
	\begin{equation} \label{eq:quad_estimate}
		\begin{split}
			\bigl\lvert (\ref{eq:quad_terms}) \bigr\rvert &\lesssim \pav{j}\pis{k+2-j}  (\ell\eta)^{\beta_{j}+\alpha_{k+2-j}} \sum_{a}\mathfrak{s}_{j}^\mathrm{av}(x_i,\dots x_{i+j-2},a)\mathfrak{s}_{k+2-j}^\mathrm{iso}(a,\bm x_{[i+j-1,i-1]},a)\\
			&\lesssim \frac{\mathfrak{s}_{k}^\mathrm{av}(\bm x)}{\eta} (\ell\eta)^{\beta_{j}}\pav{j} \frac{\pis{k+2-j}}{(\ell\eta)^{1/2-\alpha_{k+2-j}}}.
		\end{split}
	\end{equation}
	
	Finally, we bound the terms  in the third line of \eqref{eq:Fk_def}, that is (setting $j \mapsto i+j-1$)
	\begin{equation} \label{eq:sub_long}
		\Tr\bigl[\mathcal{S}\bigl[(G-M)_{[i,i+j-1]}\bigr]M^{(x_k)}_{[i+j-1,i]}\bigr],  \quad i \in \indset{k}, \quad  j \in \indset{2, (k-1) \wedge (k+1-i)},
	\end{equation}
	and (setting $j \mapsto i+k-j+1$, such that the length of the $(G-M)_{\dots}$ term is $j$)
	\begin{equation} \label{eq:sub_short}
		\Tr\bigl[\mathcal{S}\bigl[(G-M)^{(x_k)}_{[i+k-j+1,i]}\bigr]M_{[i,i+k-j+1]}\bigr], \quad i \in \indset{k}, \quad  j \in \indset {i+1,  k-1}.
	\end{equation}
	Similarly to \eqref{eq:quad_estimate}, using \eqref{eq:G-M_psi_bounds} and \eqref{eq:M_bound}, we deduce that
	\begin{equation} \label{eq:sub_final}
		\bigl\lvert  (\ref{eq:sub_long})\bigr\rvert \lesssim \frac{\mathfrak{s}_{k}^\mathrm{av}(\bm x)}{\eta} (\ell\eta)^{\beta_{j}}\pav{j} \quad \text{and} \quad \bigl\lvert  (\ref{eq:sub_short})\bigr\rvert \lesssim \frac{\mathfrak{s}_{k}^\mathrm{av}(\bm x)}{\eta} (\ell\eta)^{\beta_{j}}\pav{j}.
	\end{equation}
	Note that by \eqref{eq:sub_long}--\eqref{eq:sub_short}, the index $j$ in \eqref{eq:sub_final} only takes values in $\indset{2,k-1}$, hence $\pav{k}$ never occurs in~\eqref{eq:sub_final}.
	
	 For $k=1$, the only term on right-hand side of \eqref{eq:F1_def} is estimated analogously to \eqref{eq:1_k+1_bound1_1}, hence
	\begin{equation}
		\frac{\eta\bigl\lvert \mathcal{F}^\mathrm{av}_{[1,1]}(\bm x) \bigr\rvert}{\mathfrak{s}_1^\mathrm{av}(\bm x)} 
		\lesssim \frac{\sqrt{N\eta}}{\ell\eta} \pav{1}\pis{k+1}.
	\end{equation}	
	For $k\ge 2$,   summing \eqref{eq:1_k+1_bound1}, \eqref{eq:quad_estimate} and \eqref{eq:sub_final} according to the definition of $\mathcal{F}^\mathrm{av}_{[1,k],t}$ in \eqref{eq:Fk_def}, we conclude that, for all $k \in \indset{\maxK-1}$ and all $\bm x \in \indset{N}^k$,
	\begin{equation}
		\frac{\eta\bigl\lvert \mathcal{F}^\mathrm{av}_{[1,k]}(\bm x) \bigr\rvert}{(\ell\eta)^{\beta_k}\mathfrak{s}_k^\mathrm{av}(\bm x)} 
		\lesssim   \frac{\sqrt{N\eta}}{\ell\eta} \pav{1}\pis{k+1}  + \frac{\pav{1}}{(\ell\eta)^{\beta_{k}}} 
		+ \sum_{j=2}^{k-1} \frac{\pav{j}}{(\ell\eta)^{\beta_k-\beta_j}} \biggl(1+\frac{\pis{k+2-j}}{(\ell\eta)^{1/2-\alpha_{k+2-j}}}\biggr) + \frac{\pis{2}\pav{k}}{\sqrt{\ell\eta}}.
	\end{equation}
	Using \eqref{eq:1_k+1_bound2}--\eqref{eq:1_k+1_bound3}
	 instead of \eqref{eq:1_k+1_bound1},  
	 we obtain, for $k=\maxK$ and all $\bm x \in \indset{N}^k$,
	\begin{equation}
		\frac{\eta\bigl\lvert \mathcal{F}^\mathrm{av}_{[1,\maxK]}(\bm x) \bigr\rvert}{(\ell\eta)^{\beta_\maxK}\mathfrak{s}_\maxK^\mathrm{av}(\bm x)} \lesssim  
		 \frac{\sqrt{N\eta}}{\ell\eta} \pav{1}\pis{\maxK/2}\pis{\maxK/2+1}   
		+\sum_{j=2}^{\maxK-1} \frac{\pav{j}}{(\ell\eta)^{\beta_\maxK-\beta_j}} \biggl(1+\frac{\pis{\maxK+2-j}}{(\ell\eta)^{1/2-\alpha_{\maxK+2-j}}}\biggr) 
		+ \frac{\pis{2}\pav{\maxK}}{\sqrt{\ell\eta}}.
	\end{equation}
	In particular, \eqref{eq:av_forcing_bound} holds for all $k \in \indset{\maxK}$ by \eqref{eq:phi_av_force}. 
	Bringing back the time variable for $\mathcal{F}, \ell, \eta$,  this concludes the proof of \eqref{eq:av_forcing_bound}.
	
	\vspace{5pt}
	\noindent\textbf{Isotropic Chains}. Proof of \eqref{eq:iso_forcing_bound}. 
	Throughout the proof, we suppress the explicit dependence of the control quantities $\mathfrak{s}_{\cdot,s}$
	 as well as resolvent chains $G_{[\dots], s}$ and the corresponding deterministic approximations $M_{[\dots], s}$ on the time variable $s\le \tau$.
	
	Fist, we bound the terms appearing in \eqref{eq:iso_F1_def} and in the first line on the right-hand side of \eqref{eq:iso_Fk_def},  that is (time
	variable omitted)
	\begin{equation} \label{eq:k+1_1_iso_term}
		\bigl(G_{[1,i]}\mathcal{S}\bigl[(G_{i}-m_{i})S^q\bigr] G_{[i, k]}\bigr)_{ab} = \sum_q \Tr\bigl[(G_{i}-m_{i})S^q\bigr] \bigl(G_{[1,i]}\bigr)_{aq}\bigl(G_{[i, k]}\bigr)_{qb}, \quad i \in \indset{k}.
	\end{equation} 
	It follows from \eqref{eq:G-M_psi_bounds} that $\Tr\bigl[(G_{i}-m_{i})S^q\bigr] \lesssim (\ell\eta)^{-1}\pav{1}$. Hence, using \eqref{eq:G_psi_bounds}, we obtain, for all $i\in \indset{k}$,  
	\begin{equation} \label{eq:k+1_1_iso_bound}
		\begin{split}
			\bigl\lvert (\ref{eq:k+1_1_iso_term}) \bigr\rvert 
			\lesssim&~  \frac{\pav{1}\pis{i}\pis{k-i+1}}{(\ell\eta)^{1-\alpha_i-\alpha_{k-i+1}}}\sum_q  \mathfrak{s}_{i}^\mathrm{iso}(a,x_1,\dots,x_{i-1},q)\, \mathfrak{s}_{k-i+1}^\mathrm{iso}(q,x_i,\dots,x_{k-1},b)\\
			&+\pav{1} \biggl(\delta_{ab} + \frac{\pis{k-i+1}}{(\ell\eta)^{1/2-\alpha_{k-i+1}}}\biggr) \mathfrak{s}_{i}^\mathrm{iso}(a,x_1,\dots,x_{i-1},a)\, \mathfrak{s}_{k-i+1}^\mathrm{iso}(a,x_i,\dots,x_{k-1},b)\\
			&+\frac{\pav{1}\pis{i}}{(\ell\eta)^{1/2-\alpha_i}}\, \mathfrak{s}_{i}^\mathrm{iso}(a,x_1,\dots,x_{i-1},b)\, \mathfrak{s}_{k-i+1}^\mathrm{iso}(b,x_i,\dots,x_{k-1},b)\\
			\lesssim&~  \biggl(\frac{\pav{1}}{\sqrt{\ell\eta}}+\frac{\sqrt{N\eta}}{\ell\eta}\frac{\pav{1}\pis{i}\pis{k-i+1}}{(\ell\eta)^{\alpha_k-\alpha_i-\alpha_{k-i+1}}}\biggr) \frac{(\ell\eta)^{\alpha_k}\mathfrak{s}_{k}^\mathrm{iso}(a,\bm x',b)}{\eta}  \\
			&\lesssim   \frac{\sqrt{N\eta}}{\ell\eta} \pav{1}\pis{i}\pis{k-i+1} \frac{(\ell\eta)^{\alpha_k}\mathfrak{s}_{k}^\mathrm{iso}(a,\bm x',b)}{\eta},
		\end{split}
	\end{equation}
	  where in the penultimate step we used inequalities \eqref{eq:iso_convol}, \eqref{eq:iso_concat}, and in the last step we used that $\alpha_{i}+\alpha_{k-i+1}\le \alpha_k $ by \eqref{eq:loss_exponents} for all $i\in\indset{k}$.
	
	Next, we bound the quadratic terms in the second line of \eqref{eq:iso_Fk_def}, namely (setting $j \mapsto i+j-1$)
	\begin{equation} \label{eq:iso_quad_term}
		\sum_q \bigl((G-M)_{[i,i+j-1]}\bigr)_{qq} \bigl(G_{[1,i]}S^q G_{[i+j-1, k]}-M_{[1,i],[i+j-1,k]}^{(q)}\bigr)_{ab}, \quad i\in\indset{k},  j \in \indset{2,k-i+1}.
	\end{equation}
	It follows from \eqref{eq:G-M_psi_bounds} that 
	\begin{equation} \label{eq:iso_quad_bound}
		\begin{split}
			\bigl\lvert (\ref{eq:iso_quad_term})\bigr\rvert &\lesssim \sum_q\pis{j}\pis{k-j+2} (\ell\eta)^{\alpha_{j}+\alpha_{k-j+2}}  \mathfrak{s}_{j}^\mathrm{iso}(q,x_i,\dots, x_{i+j-2}, q)\\
			&\qquad\quad\times\mathfrak{s}_{k-j+2}^\mathrm{iso}(a,x_1,\dots,x_{i-1},q,x_{i+j-1},\dots,x_{k-1},b)\\ 
			&\lesssim \frac{\pis{j}\pis{k-j+2}}{(\ell\eta)^{1/2 + \alpha_k - \alpha_{j}-\alpha_{k-j+2}}}\,\frac{(\ell\eta)^{\alpha_k}\mathfrak{s}_{k}^\mathrm{iso}(a,\bm x',b)}{\eta} \lesssim \frac{\pis{j}\pis{k-j+2}}{\sqrt{\ell\eta}}\,\frac{(\ell\eta)^{\alpha_k}\mathfrak{s}_{k}^\mathrm{iso}(a,\bm x',b)}{\eta},
		\end{split}
	\end{equation}
	where in the second step we used \eqref{eq:Schwarz_convol}, similarly to \eqref{eq:av_insertion}, and in the last step we used $\alpha_{j}+\alpha_{k-j+2}\le \alpha_k $ that follows from \eqref{eq:loss_exponents} for all $j \in \indset{2,k-1}$.
	
	For the terms in the third line of \eqref{eq:iso_Fk_def}, (setting $j \mapsto i+j-1$)
	\begin{equation} \label{eq:iso_sublin_1}
		\sum_q \bigl((G-M)_{[i,i+j-1]}\bigr)_{qq} \bigl(M_{[1,i],[i+j-1,k]}^{(q)}\bigr)_{ab}, \quad i \in \indset{k}, j \in \indset{2, (k-1)\wedge(k-i+1)},
	\end{equation}
	we deduce, using \eqref{eq:G-M_psi_bounds},  \eqref{eq:Schwarz_convol}, and \eqref{eq:M_bound}, that (similarly to \eqref{eq:av_insertion})
	\begin{equation} \label{eq:iso_sublin_1_bound}
		\bigl\lvert (\ref{eq:iso_sublin_1})  \bigr\rvert  
		\lesssim  \frac{\pis{j}}{(\ell\eta)^{\alpha_{k}-\alpha_j}}\frac{(\ell\eta)^{\alpha_{k}}\mathfrak{s}_{k}^\mathrm{iso}(a,\bm x',b)}{\eta}.
	\end{equation}
	Similarly, for the terms in the fourth line of \eqref{eq:iso_Fk_def},
	\begin{equation} \label{eq:iso_sublin_2}
		\sum_q \bigl(M_{[i,i+k-j+1]}\bigr)_{qq} \bigl(G_{[1,i]}S^qG_{[i+k-j+1, k]}-M_{[1,i],[i+k-j+1,k]}^{(q)}\bigr)_{ab}, \quad i \in \indset{k},  j \in \indset {i+1, k-1}
	\end{equation}
	we deduce that 
	\begin{equation} \label{eq:iso_sublin_2_bound}
		\bigl\lvert (\ref{eq:iso_sublin_2})  \bigr\rvert 
		\lesssim \frac{\pis{j}}{(\ell\eta)^{\alpha_k-\alpha_j}} \frac{(\ell\eta)^{\alpha_k}\mathfrak{s}_{k}^\mathrm{iso}(a,\bm x',b)}{\eta}.
	\end{equation}
	
	Finally, we estimate the term in the fifth line of \eqref{eq:iso_Fk_def},
	\begin{equation} \label{eq:iso_av_term}
		\begin{split}
			\delta_{ab}\sum_q m_{1}m_{k} \bigl(I+\mathcal{A}_{k}\bigr)_{aq}\mathcal{X}_{t}^k(\bm x',q),
		\end{split}
	\end{equation}
	where $\mathcal{A}_{k}:= \mathcal{A}_{k,s}$ 	is defined in \eqref{eq:lin_prop_ops}. It follows from \eqref{eq:lin_prop_ops} and \eqref{eq:Ups_majorates} that 
	\begin{equation} \label{eq:Aab_bound}
		\bigl\lvert m_{1,s}m_{k,s}\bigl(I+\mathcal{A}_{k,s}\bigr)_{aq}\bigr\rvert \lesssim \delta_{aq} + (\Upsilon_s)_{aq}, \quad a,q \in \indset{N}, \quad 0\le s \le T.
	\end{equation}
	Hence, using \eqref{eq:G-M_psi_bounds}, \eqref{eq:Aab_bound}, \eqref{eq:triag} and \eqref{eq:Schwarz_convol}, we obtain
	by local contraction,  similarly to \eqref{eq:av_insertion},
	\begin{equation} \label{eq:iso_av_bound}
			\bigl\lvert (\ref{eq:iso_av_term}) \bigr\rvert \lesssim \delta_{ab} (\ell\eta)^{\beta_k}\pav{k} \sum_q \Upsilon_{aq}\mathfrak{s}_k^\mathrm{av}(\bm x',q) \lesssim \frac{\pav{k}}{(\ell\eta)^{1/2+\alpha_k-\beta_k}}  \frac{(\ell\eta)^{\alpha_k}\mathfrak{s}_{k}^\mathrm{iso}(a,\bm x',b)}{\eta}.
	\end{equation} 
	Hence, combining \eqref{eq:k+1_1_iso_bound}, \eqref{eq:iso_quad_bound}, \eqref{eq:iso_sublin_1_bound}, \eqref{eq:iso_sublin_2_bound} and \eqref{eq:iso_av_bound}, re-introducing the
	time variable, we conclude \eqref{eq:iso_forcing_bound}. This concludes the proof of Lemma \ref{lemma:forcing}. 
\end{proof}

\begin{proof} [Proof of Lemma \ref{lemma:extra_forcing}]
	First, assume that $k\in\indset{\maxK/2-1}$. The two terms on the right-hand side of \eqref{eq:extra_forcing1} are estimated analogously to \eqref{eq:short_remainder_est} using \eqref{eq:ward_local_law}.
	Hence, hence we conclude \eqref{eq:extra_forcing} for $k \in \indset{\maxK/2-1}$.
	
	Next, assume that  $k \in \indset{\maxK/2, \maxK}$. Note that the second line of \eqref{eq:extra_forcing2} is identical to the right-hand side of \eqref{eq:extra_forcing1}. The remaining terms in the first line on the right-hand side of \eqref{eq:extra_forcing2} are structurally identical to those in \eqref{eq:extra_forcing1} and hence we estimated analogously. This concludes the proof of Lemma \ref{lemma:extra_forcing}.	
\end{proof}

\subsection{Evolution equations: Proof of Lemmas \ref{lemma:av_iso_evol} and \ref{lemma:circ_evol}} \label{sec:evols}
In this section we compute the time-differentials of the $(G-M)_{[1,k],t}$ quantities along the combination of the characteristic flow \eqref{eq:char_flow} and the Ornstein-Uhlenbeck process \eqref{eq:zigOU}.
But before we analyze the random resolvent chains, we record the evolution equation for their deterministic $M$ term counterparts in the following lemma. 
\begin{lemma}[Time-Derivative of $M$] \label{lemma:dM}
	For any $k \in \mathbb{N}$, any $\bm z_{t} \in \{z_t, \overline{z}_t\}^{k}$, and any $\bm x\in\indset{x}$, the derivative of $M_{[1,k],t} := M_{[1,k],t}(\bm z_t, \bm x)$, defined in \eqref{eq:Mt_def}, satisfies the identity
	\begin{equation} \label{eq:dM}
		\begin{split}
			\frac{\mathrm{d}}{\mathrm{d}t}M_{[1,k],t} =&~ \frac{k}{2}M_{[1,k],t} +\sum_{1 \le i < j \le k} \sum_q \bigl(M_{[i,j],t}\bigr)_{qq} M_{[1,i],[j,k],t}^{(q)},
		\end{split}
	\end{equation}
	where we recall the definition of $M_{[1,i],[j,k],t}^{(q)}$ from \eqref{eq:M_with_q}.
\end{lemma}
We defer the proof of Lemma \ref{lemma:dM} to Section \ref{sec:M_rec_analysis}, and proceed to prove Lemma \ref{lemma:av_iso_evol}.
\begin{proof}[Proof of Lemma \ref{lemma:av_iso_evol}]
	As the starting point of our computation, we apply It\^{o}'s formula to the definition of $G_{[1,k],t}$ in
	 \eqref{eq:Gk_def},  together with the characteristic flow~\eqref{eq:char_flow}, obtaining 
	\begin{equation} \label{eq:dG_eq}
		\mathrm{d}G_{[1,k],t} = \frac{k}{2}G_{[1,k],t} + \sum_{1\le i < j \le k} G_{[1,i]} \mathcal{S}\bigl[ G_{[i,j],t}\bigr]G_{[j,k],t} + \sum_{i=1}^k  G_{[1,i]} \mathcal{S}\bigl[ G_{i,t}-m_{i,t}\bigr]G_{[i,k],t} + \mathrm{d}\mathcal{M}_{[1,k],t}.
	\end{equation}
	Combining \eqref{eq:dG_eq} with \eqref{eq:dM} and collecting the terms, we conclude that the evolution equation \eqref{eq:k_iso_evol} for the isotropic objects $\mathcal{Y}^k_{t}$, defined in \eqref{eq:G-M_kiso} holds.
	
	Multiplying \eqref{eq:k_iso_evol} by $S^{x_k}$ and taking the trace, we obtain the evolution equation \eqref{eq:k_av_evol} for the averaged objects $\mathcal{X}_t^k$, defined in \eqref{eq:G-M_kav}. Note that the term in the last line of \eqref{eq:iso_Fk_def} becomes $\mathcal{A}_{k,t}[\mathcal{X}_t^k]$ on the right-hand side of \eqref{eq:k_av_evol},
	using the trivial identities $|m|^2 S (I + \Theta) =\Theta$ and $m^2S (I+\Xi) = \Xi$.
	This concludes the proof of Lemma \ref{lemma:av_iso_evol}.
\end{proof}

Next, we prove Lemma \ref{lemma:circ_evol}.
\begin{proof} [Proof of Lemma \ref{lemma:circ_evol}]
	Consider linear operators $\mathcal{C}_t$
	and $\mathcal{D}$ on $(\mathbb{C}^N)^{\otimes 2}$ acting on functions $f(x,y)$ of $x,y \in \indset{N}$ by
	\begin{equation}
		\mathcal{D}[f](x,y) := \sum_b \bigl(\delta_{yb}-\delta_{xy}\bigr) f(x,b), \quad \mathcal{C}_t[f](x,y) := \sum_{a} (\Theta_t)_{xa} f(a,y) + \sum_{b} (\Theta_t)_{yb} f(x,b).
	\end{equation}
	A direct computation reveals that their commutator $[\mathcal{D},\mathcal{C}_t]$ satisfies
	\begin{equation} \label{eq:DC_commut}
		[\mathcal{D},\mathcal{C}_t][f](x,y) = \mathcal{D}\biggl[ (\Theta_t)_{xy} \sum_b f(y,b) + (\Theta_t)_{xy} \sum_b f(x,b) \biggr].
	\end{equation}
	
	Now fix an even $k \in \indset{\maxK/2-1}$, then $\mathcal{Z}_t^k = \reg{\mathcal{X}}_t^k$, defined in \eqref{eq:Gcirc}. Denoting by $\mathcal{D}^{(k)}:=I^{\otimes (k-2)}\otimes\mathcal{D}$,
	we can rewrite 	$\mathcal{Z}_t^k = \Theta_t^{(k)}\circ \mathcal{D}^{(k)}[\mathcal{X}_t^k]$. Furthermore, with  $\mathcal{C}_t^{(k)}:=\Theta_t^{(k-1)}+\Theta_t^{(k)}$, we deduce from  \eqref{eq:k_av_evol}, \eqref{eq:Sring_def}, \eqref{eq:Gcirc}, and \eqref{eq:dTheta}, that 
	\begin{equation}
		\begin{split}
			\mathrm{d}\mathcal{Z}_t^k &= \mathrm{d}\bigl\{\Theta_t^{(k)}\circ \mathcal{D}^{(k)}[\mathcal{X}_t^k] \bigr\} = \mathrm{d}\bigl\{\Theta_t^{(k)}\bigr\}\circ \mathcal{D}^{(k)}[\mathcal{X}_t^k] + \Theta_t^{(k)}\circ \mathcal{D}^{(k)}[\mathrm{d}\mathcal{X}_t^k]\\
			&= \biggl(1+\frac{k}{2}+\Theta_t^{(k)} + \Theta_t^{\oplus k}\biggr)[\mathcal{Z}_t^k]\mathrm{d}t +\Theta_t^{(k)}\circ [\mathcal{D}^{(k)},\mathcal{C}_t^{(k)}] [\mathcal{X}_t^k]\mathrm{d}t + \mathcal{U}_t^{(k)}\biggl[
			\mathrm{d}\mathcal{M}^\mathrm{av}_{[1,k],t}  + \mathcal{F}^\mathrm{av}_{[1,k],t} \mathrm{d}t\biggr],
		\end{split}
	\end{equation}
	where we used  that  
	$\Theta^{(k)}_t\circ\mathcal{D}^{(k)} = \reg{\Theta}_t^{(k)} = \mathcal{U}^{(k)}_t$ for $k \in \indset{\maxK/2-1}$.
	Here $\mathcal{U}^{(k)}_t$ is defined in \eqref{eq:U_operator}. We also
	recall that $\reg{\mathcal{X}}_t^k$ is defined only for saturated chains of even length, hence $\mathcal{A}_{t,j} = \Theta_t$ for all $j\in\indset{k}$ on the right-hand side of the evolution equation \eqref{eq:k_av_evol} for the corresponding $\mathcal{X}^k_t$.
	It follows from \eqref{eq:DC_commut} that 
	\begin{equation}
		[\mathcal{D}^{(k)},\mathcal{C}_t^{(k)}] [\mathcal{X}_t^k]  = \mathcal{D}^{(k)} \bigl[ \other{\mathcal{F}}^\mathrm{av}_{[1,k],t}\bigr],
	\end{equation}
	where $\other{\mathcal{F}}^\mathrm{av}_{[1,k],t}$ is defined in \eqref{eq:extra_forcing1}. Hence, we obtain the evolution equation
	\begin{equation} \label{eq:Z_elov_short}
		\mathrm{d}\mathcal{Z}_t^k = \biggl(1+\frac{k}{2}+\Theta_t^{(k)} + \Theta_t^{\oplus k}\biggr)[\mathcal{Z}_t^k]\mathrm{d}t +\mathcal{U}_t^{(k)}\biggl[
		\mathrm{d}\mathcal{M}^\mathrm{av}_{[1,k],t}  + \mathcal{F}^\mathrm{av}_{[1,k],t} \mathrm{d}t + \other{\mathcal{F}}^\mathrm{av}_{[1,k],t}\mathrm{d}t\biggr], \quad 
		k\in\indset{\maxK/2-1}.
	\end{equation}
	
	Completely analogously, using the commutator identity $[\mathcal{D}\otimes\mathcal{D}, 
	\mathcal{C}_t\oplus\mathcal{C}_t] = \mathcal{D}\otimes[\mathcal{D}, \mathcal{C}_t] + [\mathcal{D}, \mathcal{C}_t] \otimes \mathcal{D} $, we deduce that, for all even $k\in\indset{\maxK/2,\maxK}$,
	\begin{equation} \label{eq:Z_elov_long}
		\mathrm{d}\mathcal{Z}_t^k = \biggl(2+\frac{k}{2}+\Theta_t^{(k-2)}+\Theta_t^{(k)}  + \Theta_t^{\oplus k}\biggr)[\mathcal{Z}_t^k]\mathrm{d}t +\mathcal{U}_t^{(k)}\biggl[
		\mathrm{d}\mathcal{M}^\mathrm{av}_{[1,k],t}  + \mathcal{F}^\mathrm{av}_{[1,k],t} \mathrm{d}t + \other{\mathcal{F}}^\mathrm{av}_{[1,k],t}\mathrm{d}t\biggr], \quad 
	\end{equation}
	where, for $k \in \indset{\maxK/2,\maxK}$, the additional forcing term $\other{\mathcal{F}}^\mathrm{av}_{[1,k],t}$ is defined in \eqref{eq:extra_forcing2}.
	
	By definition of $\mathcal{P}_{s,t}$ in \eqref{eq:Psat_def} and $\reg{\Theta}^{(j)}_s$ in \eqref{eq:Thetaring}, we have $\mathcal{P}_{s,t}\circ\reg{\Theta}^{(j)}_s = \reg{\Theta}^{(j)}_t$, and hence
	\begin{equation} \label{eq:U_prop}
		\mathcal{P}_{s,t}^{(k)}\circ\mathcal{U}^{(k)}_s = \mathcal{U}^{(k)}_t,\quad k \in\indset{\maxK/2-1}, \quad\text{and} \quad \mathcal{P}_{s,t}^{(k-2)}\circ\mathcal{P}_{s,t}^{(k)}\circ\mathcal{U}^{(k)}_s = \mathcal{U}^{(k)}_t,\quad k \in\indset{\maxK/2,\maxK}.
	\end{equation}
	Therefore, applying Duhamel's principle to \eqref{eq:Z_elov_short} and \eqref{eq:Z_elov_long}, and using \eqref{eq:U_prop}, together
	with the identity $\mathcal{Z}_t^k = \mathcal{U}^{(k)}_t\big[ \mathcal{X}_t^k \big]$,
	we obtain \eqref{eq:Gcirc_solve} for all $k\in\indset{K}$. This concludes the proof of Lemma \ref{lemma:circ_evol}.
\end{proof}

\section{Global laws: Proof of Proposition \ref{prop:global_laws}}\label{sec:global} 
In this section we prove the multi-resolvent global laws, that will be applied to $z_0$ at the initial time $t=0$. Recall that spectral parameter $z_0$ under consideration satisfies $\eta_0 := |\im z_0| \ge c$ and $|z_0|\le C$. Thus we
have  $\ell_0\sim W$ by~\eqref{def:ell_def}.  This choice  of $\eta_0, \ell_0$
will remain valid throughout Section~\ref{sec:global}, and we drop the subscript $0$ from $\eta_0, \ell_0, G_{0}$.
Recall that the statement of Proposition \ref{prop:global_laws} is equivalent to
\begin{equation} \label{eq:av_global}
		\max_{\bm z_0\in \{z_0, \bar z_0\}^k} \max_{\bm x\in \indset{N}^k} \frac{\bigl\lvert \Tr \bigl[ \bigl(G_{[1,k],0} - M_{[1,k],0}\bigr) (\bm x')  S^{x_k}\bigr]\bigr\rvert}{  \mathfrak{s}_{k,0}^\mathrm{av}(\bm x)} \prec 1, \quad k \in \indset{\maxK}, 
	\end{equation}
	\begin{equation} \label{eq:iso_global}
		\max_{\bm z_0 \in \{z_0, \bar z_0\}^k} \max_{a,b\in\indset{N}} \max_{\bm x'\in \indset{N}^{k-1}}  \frac{\bigl\lvert \bigl(G_{[1,k],0} - M_{[1,k],0}\bigr)_{ab} (\bm x') \bigr\rvert}{\mathfrak{s}_{k,0}^\mathrm{iso}(a,\bm x',b)} \prec 1, \quad k \in \indset{\maxK},
	\end{equation}
	with $G_{[p,j],0}$ defined according to \eqref{eq:resolvent_chains} and \eqref{eq:Gk_def}, and the size functions $\mathfrak{s}_{k,0}^{\mathrm{av/iso}}$ are defined in \eqref{eq:sfunc_def}.

The proof of Proposition \ref{prop:global_laws} is carried out in three steps. First, we establish a priori estimates on a single resolvent without any spatial component, that is, every entry is controlled by the same parameter. 
We call them {\it flat bounds}. In the key 
second step, we use an inductive procedure to gradually improve the flat bound from the first step to reflect the optimal spatial decay of the resolvent entries $G_{ab}$ in $|a-b|$. This step will also require controlling chains of length $k=2$, but
not longer. Finally, in the third step, we derive a system of static master inequalities for both isotropic and averaged chains of arbitrary length, that we gradually solve by induction in chain length $k$ using the estimates from steps one and two as input, thus completing the proof of Proposition \ref{prop:global_laws}. 

Before carrying out these three steps, we mention that global laws are generally  easier than local laws
since the precise $\eta$-dependence of the estimates is not critical. Moreover, reduction is not a serious issue since
every excess resolvent can affordably be bounded by $\| G\|\le \eta^{-1}\sim 1$. 
Indeed, for mean-field models the proof of the global law is typically a routine argument spanning only $2$--$3$ pages.
In the setup of band matrices, however, 
establishing global laws is still not an easy task, due to the delicate spatial dependence. In particular, proving the decay of $G_{ab}$
in $|a-b|_N$ is highly nontrivial.

We emphasize that in Section \ref{sec:global}, we use only the weaker convolution estimate \eqref{eq:true_convol_notime} instead of \eqref{eq:convol_notime}--\eqref{eq:suppressed_convol_notime}.

\subsection{Step 1. Flat bound}
In this section, we establish the initial bound on a single resolvent $G$, which is the content of the following lemma.
\begin{lemma}[Flat A Priori Bounds] \label{lemma:apriori} Let $z = E + \I \eta \in \mathbb{H}$ be a spectral parameter satisfying  $|z|\le C$ and $\eta \ge c$. Then the resolvent $G(z) := (H-z)^{-1}$ satisfies 
	\begin{equation} \label{eq:apriori_flat}
		\bigl\lvert \bigl(G(z)-m(z)\bigr)_{ab} \bigr\rvert \prec W^{-1/2},
	\end{equation}
	\begin{equation} \label{eq:apriori_av}
		\bigl\lvert \Tr \bigl[\bigl(G(z)-m(z)\bigr) S^x\bigr] \bigr\rvert \prec W^{-3/4},
	\end{equation}
	uniformly in $a,b,x \in \indset{N}$, recalling that $m(z)_{ab}= \delta_{ab} m(z)$.
\end{lemma}
We prove Lemma \ref{lemma:apriori} in two steps: First, we establish the entry-wise bound \eqref{eq:apriori_flat} using a routine bootstrapping procedure to reduce $\eta$, presented in Section \ref{sec:flat_entrywise}; then, we use the entry-wise bound as an input to obtain the averaged bound \eqref{eq:apriori_av} in Section \ref{sec:apriori_av}. 
Note that the entry-wise flat bound~\eqref{eq:apriori_flat} is optimal, but
the averaged bound \eqref{eq:apriori_av} only captures part of the fluctuation averaging effect
of the trace, and is off by a factor $W^{-1/4}$ from its optimal form \mbox{$\lvert \Tr [(G(z)-m(z)) S^x] \rvert \prec W^{-1}$},
 that we will eventually prove in the subsequent Section~\ref{sec:global3}.

\subsubsection{Flat entry-wise bound: Proof of \eqref{eq:apriori_flat}} \label{sec:flat_entrywise}
For an $N\times N$ matrix $X$, we defined the norm $\norm{X}_{\max}$ as 
\begin{equation}\label{def:maxnorm}
	\norm{X}_{\max} := \max_{a,b} |X_{ab}|.
\end{equation}
The proof of the entry-wise bound \eqref{eq:apriori_flat} is based on the following lemma. 
\begin{lemma} [Gap in the Values of $\norm{G-m}_{\max}$] \label{lemma:iso_gap}
	Let $\xi, \nu > 0$ be small fixed tolerance exponents, let $c \sim 1$ be a small constant, and let $z := E+\I \eta \in \mathbb{H}$ be a spectral parameter satisfying $c \le \eta \le N^{\xi}$ and $|E| \le C$. Then,
	\begin{equation} \label{eq:iso_gap}
		\bigl\lVert G(z)-m(z)\bigr\rVert_{\max} \lesssim N^{-\xi}~ \text{w.v.h.p.} \quad \Longrightarrow\quad \bigl\lVert G(z)-m(z)\bigr\rVert_{\max} \le \frac{N^{\nu}}{\sqrt{W}}~ \text{w.v.h.p.}
	\end{equation}
\end{lemma}
Once Lemma \ref{lemma:iso_gap} is established, \eqref{eq:apriori_flat} follows by the usual continuity argument presented below.
\begin{proof}[Proof of the entry-wise a priori bound \eqref{eq:apriori_flat}]
	First, observe that the assumption of \eqref{eq:iso_gap} holds trivially at $\eta = N^{\xi}$ owing to the norm bound $\norm{G(E+\I \eta)} \le \eta^{-1}$ and $|m(E+\I\eta)| \lesssim \eta^{-1}$. 
	
	Moreover, it is straightforward to check that the map $\eta \mapsto \norm{G(E+\I \eta)- m(E+\I\eta)}_{\max}$ is $\eta^{-2}$-Lipschitz regular. Therefore, we choose a fine $N^{-1}$-grid $\{\other{\eta}_j := (N^{\xi} - N^{-1}j)\vee c\}$ of $\eta$'s in the interval $[c, N^{\xi}]$, with cardinality at most $N^{1+\xi}$, and prove that 
	\begin{equation}
		\bigl\lVert G(E+\I \other{\eta}_j)-m(E+\I\other{\eta}_j)\bigr\rVert_{\max} \le \frac{N^{\nu}}{\sqrt{W}}~ \text{w.v.h.p.}
	\end{equation}
	for every $\other{\eta}_j$ in the grid by induction in $j$ using Lipschitz regularity of $\norm{G-m}_{\max}$ and \eqref{eq:iso_gap}. In particular, we obtain \eqref{eq:apriori_flat} for $\eta=c$ with any fixed $c\sim 1$. This concludes the proof of \eqref{eq:apriori_flat}.
\end{proof}
Therefore, it remains to prove Lemma \ref{lemma:iso_gap}. To this end, we employ our main tool for estimating the resolvent at large $\eta$---the second order Gaussian renormalization, introduced in \cite{cipolloni2021eigenstate}. For functions $f(H), g(H)$ of a random matrix $H$, the re normalization of the matrix $f(H)Hg(H)$, denoted by the underline, is defined as
\begin{equation} \label{eq:underline_def}
	\underline{f(H)H g(H)}:= f(H)Hg(H) - \Expv_{\other{H}}\bigl[\bigl(\partial_{\other{H}}f\bigr)(H)\other{H}g(H) + f(H)\other{H}\bigl(\partial_{\other{H}}g\bigr)(H)\bigr],
\end{equation}
where $\partial_{\other{H}}$ denotes the derivative in the direction of a random matrix $\other{H}=\other{H}^*$, that has independent centered Gaussian entries with matrix of variances $S$, independent of $H$, and $\Expv_{\other{H}}$ denotes the expected value with respect to $\other{H}$. This normalization  is designed in such a way that $\Expv \underline{f(H)H g(H)}=0$
if the entries of $H$ were Gaussian.
In particular, the definition \eqref{eq:underline_def} with $f(H) := I$ and $g(H):=G(z)$ yields 
\begin{equation}
	\underline{HG} = HG +\mathcal{S}[G]G.
\end{equation}
Therefore, \eqref{eq:Dyson} and the resolvent identity $HG = I + zG$ imply the well-knows quadratic equation for $G-m$ (that we refer to as the \emph{static equation}),
\begin{equation} \label{eq:static1G}
	G-m = -m \, \mathcal{B}^{-1}\bigl[\underline{HG}\bigr] + m \,\mathcal{B}^{-1}\bigl[ \mathcal{S}[G-m](G-m)\bigr],
\end{equation}
where the (stability) super-operator $\mathcal{B} \equiv \mathcal{B}(z)$ is defined by its action on $R \in \mathbb{C}^{N\times N}$,
\begin{equation} \label{eq:B_stab_def}
	\mathcal{B}(z)\bigl[R\bigr] := R - m(z)^2\mathcal{S}\bigl[R\bigr].
\end{equation}
Therefore, the action of the operator $\mathcal{B}^{-1}$ on $R\in \mathbb{C}^{N\times N}$ is given by
\begin{equation} \label{eq:Binv_action}
	\bigl(\mathcal{B}^{-1}[R]\bigr)_{xy} = R_{xy}, \quad \bigl(\mathcal{B}^{-1}[R]\bigr)_{xx} = \sum_{c} (I + \Xi)_{xc} R_{cc}, \quad x,y\in \indset{N}, \quad x\neq y,
\end{equation}
where $\Xi := \Xi(z)$ was defined in~\eqref{def:ThetaXi}. 
In particular,  for any $z\in\mathbb{H}$ satisfying $\im z \ge c \sim 1$, we have
\begin{equation} \label{eq:stab_bound}
	\norm{\mathcal{B}(z)^{-1}[R]}_{\max} \lesssim \norm{R}_{\max}.
\end{equation}

We make use of the following standard technical lemma.
\begin{lemma} \label{lemma:high-moment}
	Let $\mathcal{X}$ be a random variable satisfying $|\mathcal{X}| \le N^D$ for some constant $D > 0$, and let $\varphi$ be a deterministic control parameter satisfying $\varphi \ge N^{-D}$. Assume that $\mathcal{X}$ satisfies the equation\footnote{
	Here, we choose the notation $\underline{\mathcal{X}}$ by analogy with the static equation \eqref{eq:static1G} and its generalization \eqref{eq:Gk_static}. However, in the context of Lemma \ref{lemma:high-moment}, $\underline{\mathcal{X}}$ is unrelated to the underline notation used for Gaussian renormalization, defined in \eqref{eq:underline_def}.}
	\begin{equation}
		\mathcal{X} = \underline{\mathcal{X}} + \mathcal{E},
	\end{equation}
	where $\underline{\mathcal{X}}$ and $\mathcal{E}$ are random variables satisfying 
	\begin{equation}
		\lvert \mathcal{E} \rvert \lesssim \varphi, \quad \text{w.v.h.p.}, \quad \bigl\lvert \Expv \bigl[ \underline{\mathcal{X}}\cdot  \overline{\mathcal{X}} \lvert \mathcal{X}\rvert^{2p-2} \bigr]  \bigr\rvert \lesssim C_p \varphi^{2p} + N^{-\nu} \Expv \bigl[\lvert \mathcal{X}\rvert^{2p}\bigr],
	\end{equation}
	for any integer $p\in\mathbb{N}$ and some small exponent $\nu > 0$. Here, $\overline{\mathcal{X}}$ denotes the complex conjugate of $\mathcal{X}$. Then, $\mathcal{X}$ satisfies the bound
	\begin{equation}
		\lvert \mathcal{X} \rvert \prec \varphi.
	\end{equation}
\end{lemma}
The proof of Lemma \ref{lemma:high-moment} is a straightforward application of the definition of stochastic domination, and we leave its details to the reader.

\begin{proof} [Proof of Lemma \ref{lemma:iso_gap}]
	Throughout the proof, since the spectral parameter $z := E + \I \eta$ is fixed, we omit it from the argument of $G$ and $m$. Assume that for a deterministic control parameter $\psi \ge 1$, the norm $\norm{G-m}_{\max}$ satisfies 
	\begin{equation} \label{eq:max_ansatz}
		\norm{G-m}_{\max} \lesssim W^{-1/2}\psi, \quad \text{w.v.h.p.}
	\end{equation}
	In particular, it follows from the assumption of \eqref{eq:iso_gap} that the second term on the right-hand side of \eqref{eq:static1G} satisfies
	\begin{equation} \label{eq:1G_flat_quad_bound}
		\bigl\lVert \mathcal{B}^{-1}\bigl[\mathcal{S}[G-m](G-m)\bigr]\bigr\rVert_{\max} \lesssim N^{-\xi} W^{-1/2}\psi, \quad \text{w.v.h.p.,}
	\end{equation}
	where we used \eqref{eq:stab_bound}. 
	Therefore, to use Lemma \ref{lemma:high-moment} for $\mathcal{X} = \mathcal{Y}_0:= (G-m)_{ab}$, it suffices to establish a bound on $\lvert\Expv[(\mathcal{B}^{-1}[\underline{HG}])_{ab}\overline{(G-m)_{ab}}\lvert (G-m)_{ab}\rvert^{2p-2}]\rvert$, which is the content of the following claim.
	\begin{claim} \label{claim:flat1G_underline}
		Fix $a,b\in\indset{N}$, and let $\mathcal{Y}_0 := (G-m)_{ab}$. Then, under the assumptions of Lemma \ref{lemma:apriori}, and additionally assuming~\eqref{eq:max_ansatz},
		 for any integer $p\in\mathbb{N}$ and $0 < \nu < 1/100$, we have
		\begin{equation} \label{eq:flat1G_underline}
			\biggl\lvert\Expv\biggl[\bigl(\mathcal{B}^{-1}[\underline{HG}]\bigr)_{ab}\overline{\mathcal{Y}_0 } \lvert \mathcal{Y}_0  \rvert^{2p-2}\biggr]\biggr\rvert 
			\lesssim \biggl(\frac{N^{\nu}}{\sqrt{W}}\bigl(\sqrt{\psi}+  W^{-1/4}\psi\bigr) \biggr)^{2p} + N^{-\nu} \Expv \bigl[|\mathcal{Y}_0 |^{2p}\bigr]. 
		\end{equation}
	\end{claim}
	We prove Claim \ref{claim:flat1G_underline} in Section \ref{sec:cum_expand}, after presenting the general formalism of cumulant expansion. It follows from \eqref{eq:static1G}, \eqref{eq:1G_flat_quad_bound}, Lemma \ref{lemma:high-moment}, and Claim \ref{claim:flat1G_underline}, that
	\begin{equation}
		\norm{G-m}_{\max} \prec W^{-1/2}\bigl(\sqrt{\psi}+ W^{-1/4}\psi\bigr).
	\end{equation}
	Hence, the conclusion of \eqref{eq:iso_gap} follows by a standard iteration argument
	that gradually reduces $\psi$ from its initial value $\psi= W^{1/2}N^{-\xi}$, guaranteed by
	the assumption in \eqref{eq:iso_gap}, to its final value $\psi=N^{\nu}$. 
	  This concludes the proof of Lemma \ref{lemma:iso_gap}.
\end{proof}

\subsubsection{A priori averaged bound: Proof of \eqref{eq:apriori_av}} \label{sec:apriori_av}
\begin{proof}[Proof of the averaged a priori bound \eqref{eq:apriori_av}]
	We consider the index $x\in\indset{N}$ to be fixed. It follows from \eqref{eq:static1G} that 
	\begin{equation} \label{eq:1G_av_static}
		\Tr\bigl[(G-m)A_1\bigr] = -\, \Tr\bigl[ \underline{HG} A_1'\bigr] + \,\Tr\bigl[  \mathcal{S}[G-m](G-m)A_1'\bigr] ,
	\end{equation}
	 holds for any deterministic matrix $A_1$, 
	where  $A_1' := m\,\mathcal{B}^{-1}[A_1]$. Let $A_1 := S^x$, then the entry-wise bound \eqref{eq:apriori_flat} proven above and \eqref{eq:Binv_action} imply that
	\begin{equation} \label{eq:apriori_av_quad}
		\bigl\lvert  \Tr\bigl[  \mathcal{S}[G-m](G-m)A_1'\bigr] \bigr\rvert  
		\lesssim \norm{G-m}_{\max}^2  \sum_{du}\bigl\lvert (I+\Xi)_{du}\bigr\rvert S_{ux} \prec W^{-1}.
	\end{equation}
	On the other hand, the underline term in \eqref{eq:1G_av_static} satisfies the following bound.
	\begin{claim} \label{claim:1Gav_underline}
		Let $A_1=S^x$ and $\mathcal{X}_1 := \Tr[(G-m)A_1]$, then, for any integer $p\in\mathbb{N}$, and any $\nu \in (0,1/100)$, we have
		\begin{equation} \label{eq:apriori_av_underline}
			\biggl\lvert \Expv\biggl[\Tr\bigl[ \underline{HG} A_1'\bigr]\cdot \overline{\mathcal{X}_1} \lvert \mathcal{X}_1 \rvert^{2p-2}\biggr] \biggr\rvert \lesssim \biggl(\frac{N^{\nu}}{W^{3/4}}\biggr)^{2p} + N^{-\nu} \Expv\bigl[\lvert \mathcal{X}_1 \rvert^{2p}\bigr].
		\end{equation}
	\end{claim}
	We prove Claim \ref{claim:1Gav_underline} in Section \ref{sec:cum_expand}. 
	Together with \eqref{eq:1G_av_static}, \eqref{eq:apriori_av_quad} and \eqref{eq:apriori_av_underline}, Lemma \ref{lemma:high-moment} implies the desired \eqref{eq:apriori_av}. In particular, together with the proof of \eqref{eq:apriori_flat} in Section \ref{sec:flat_entrywise} above, this concludes the proof of Lemma \ref{lemma:apriori}.
\end{proof}

\subsection{Step 2. Spatial decay for one- and two-resolvent chains} \label{sec:glob_decay}
	In this section, we gradually improve the entry-wise bound \eqref{eq:apriori_flat} to capture the optimal spatial decay of the resolvent entries $G_{ab}$. 	
	The leading contribution to the fluctuations of $G-m$, coming from the underline term in \eqref{eq:static1G}, is given by
	\begin{equation}\label{eq:GSG}
		\biggl(\sum_{c}S_{ac} |G_{cb}|^2\biggr)^{1/2} = \bigl( (G^*S^aG)_{bb}\bigr)^{1/2}.
	\end{equation}
	Therefore, to recover the optimal spatial decay of $(G-m)_{ab}$, we have to append the two-resolvent chains to our analysis. Estimating the left hand side of~\eqref{eq:GSG} via $G_{ab}$ would not 
	give a self-improving relation, but the right hand side of~\eqref{eq:GSG} can be
	estimated more effectively by treating it as a two-resolvent chain and analyzing its own evolution equation.  Hence, our goal is to prove the following lemma.
	\begin{lemma}[Isotropic Global Laws with Spatial Components for $1G$ and $2G$ Chains]\label{lemma:1G_global_decay}
		Let $z = E + \I \eta \in \mathbb{H}$ be a spectral parameter satisfying  $|z| \le C$ and $\eta \ge c$. Then, for any choice of $z_1,z_2 \in \{z,\overline{z}\}$, the resolvents $G_j(z) := (H-z_j)^{-1}$ with $j \in \{1,2\}$ satisfy 
		\begin{equation} \label{eq:1G_global_decay}
			\bigl\lvert (G_{j}-m_{j})_{ab} \bigr\rvert \prec \mathfrak{s}_1^\mathrm{iso}(a,b), \quad \bigl\lvert \bigl((G-M)_{[1,2]}(x)\bigr)_{ab} \bigr\rvert \prec \mathfrak{s}_2^\mathrm{iso}(a,x,b), \quad a,b,x \in \indset{N},
		\end{equation}
		where $G_{[1,2]} \equiv G_{[1,2]}(x) := G_1 S^{x}G_2$, and $M_{[1,2]} \equiv M_{[1,2]}(x) := \bigl(1- m_1m_2\mathcal{S}\bigr)^{-1}[m_1m_2S^{x}]$ denotes the corresponding deterministic approximation. Here $\mathfrak{s}_k (\cdot) := \mathfrak{s}_{k,0} (\cdot)$ are defined in \eqref{eq:sfunc_def} with $t=0$.
	\end{lemma}
	That is, the isotropic global law \eqref{eq:iso_global} in Proposition \ref{prop:global_laws} holds for $k \in \{1,2\}$.
	
	To gradually recover the spatial decay from the flat bound \eqref{eq:apriori_flat}, we introduce a family of isotropic \emph{modified size functions}  $\mathfrak{s}_k(\omega;\,\cdot \,)$ parameterized by $\omega \in [0, W^{-1}]$ that interpolate between the target size functions $\mathfrak{s}_k^\mathrm{iso}(\,\cdot \,)$ and a constant function, and 
	we will iteratively decrease $\omega$ from $W^{-1}$ to $N^{-2D'}$.
	For $\omega \in [0, W^{-1}]$, define the $\omega$-dependent modified size functions $\mathfrak{s}_1(\omega; a,b)$ and $\mathfrak{s}_2(\omega; a,x,b)$ as
	\begin{equation} \label{eq:omegas_def}
		\mathfrak{s}_1(\omega; a,b) := \sqrt{\Upsilon_{ab} + \omega}, \quad 
		\mathfrak{s}_2(\omega; a,x,b) := \frac{1}{\sqrt{W}}\sqrt{\bigl(\Upsilon_{ax} + \omega\bigr)\bigl(\Upsilon_{xb} + \omega\bigr)}, \quad a,b,x\in \indset{N},
	\end{equation}
	where $\Upsilon := \Upsilon_0$ is an admissible control function as in Definition \ref{def:adm_ups}, i.e., with $\ell \sim W$ and $\eta\sim 1$.
	
	We prove Lemma \ref{lemma:1G_global_decay} by iteration, using the following lemma as a basic step.
	\begin{lemma} [$\omega$-Bootstrap] \label{lemma:omega_lemma}
		Let $\zeta > 0$ be a fixed tolerance exponent satisfying $\zeta < 1/100$, and let  $\omega \in [0, W^{-1}]$ and $\psi \in [N^{10\zeta}, W^{1/2}]$ be a pair of deterministic control parameters. Let $z:=E+\I\eta\in\mathbb{H}$ be a spectral parameter satisfying $|E|\le C$ and $\eta \sim 1$. Assume that, for $j=1,2$,  the bounds
		\begin{equation} \label{eq:omega_assume}
			\bigl\lvert (G_j-m_j)_{ab} \bigr\rvert \lesssim N^{2\zeta}\mathfrak{s}_1(\omega; a,b), \quad \bigl\lvert \bigl((G-M)_{[1,2]}(x)\bigr)_{ab} \bigr\rvert \lesssim \psi\, \mathfrak{s}_2(\omega; a,x,b),
		\end{equation}
		hold with very high probability, with some deterministic $\psi\in [1, \sqrt{W}]$,
		 for all ${z_1, z_2} \in \{z, \overline{z}\}$, and $a,b,x \in \indset{N}$. Then the improved bounds
		\begin{equation} \label{eq:omega_conclude}
			\bigl\lvert (G_j-m_j)_{ab} \bigr\rvert \lesssim N^{\zeta}\mathfrak{s}_1\bigl(N^{-\zeta}\omega; a,b\bigr), \quad \bigl\lvert \bigl((G-M)_{[1,2]}(x)\bigr)_{ab} \bigr\rvert \lesssim \bigl(N^{10\zeta} + N^{-\zeta}\psi\bigr)\, \mathfrak{s}_2\bigl(N^{-\zeta}\omega; a,x,b\bigr),
		\end{equation}
		also hold with very high probability, for all ${z_1, z_2} \in \{z, \overline{z}\}$, and $a,b,x \in \indset{N}$. 
	\end{lemma}
	
	We defer the proof of Lemma \ref{lemma:omega_lemma} to the end of this section. Equipped with Lemma \ref{lemma:omega_lemma}, we prove Lemma \ref{lemma:1G_global_decay}.
	\begin{proof}[Proof of Lemma \ref{lemma:1G_global_decay}]
		Fix a $\zeta \in (0, 1/100)$ to be the target tolerance exponent in the stochastic domination bounds \eqref{eq:1G_global_decay}. 
		First, observe that, by \eqref{eq:Ups_norm_bounds}, for $\omega_0 = W^{-1}$ we have
		\begin{equation} \label{eq:omega0_s}
			\mathfrak{s}_1(\omega_0; a,b) \sim W^{-1/2}, \quad \mathfrak{s}_2(\omega_0; a,x,b) \sim W^{-3/2}.
		\end{equation}
		On the other hand, \eqref{eq:ups_lower_bound} implies that for any $\omega \le \omega_{\mathrm{f}} := N^{-2D'}$, the modified size functions $\mathfrak{s}_{k}(\omega; \, \cdot\,)$ satisfy
		\begin{equation} \label{eq:omega_final_s}
			\mathfrak{s}_{k}(\omega; \,\cdot\,) \sim \mathfrak{s}_{k}^\mathrm{iso}(\, \cdot \,), \quad k \in \{1,2\}, \quad \omega \le  N^{-2D'},
		\end{equation}
		where $\mathfrak{s}_{k}^\mathrm{iso} := \mathfrak{s}_{k,0}^\mathrm{iso}$ are the size functions defined in \eqref{eq:sfunc_def}. Moreover, or any $z_1,z_2 \in \{z, \overline{z}\}$ and all $a,b,x\in\indset{N}$, using Schwarz inequality and Ward identity, we obtain 
		\begin{equation} \label{eq:2G_flat}
			\begin{split}
				\bigl\lvert \bigl(G_{[1,2]}(x) \bigr)_{ab}\bigr\rvert 
				&\le \sum_{c} S_{xc} \lvert G_{ac}(z_1) G_{cb}(z_2) \rvert \lesssim\frac{1}{W}\sqrt{\sum_{c} \lvert G_{ac}(z_1)\rvert^2\sum_{d} \lvert G_{db}(z_2) \rvert^2}\\
				&= \frac{1}{W\eta} \sqrt{\bigl(\im G(z_1)\bigr)_{aa}\bigl(\im G(z_2)\bigr)_{bb}} \le \frac{1}{W\eta^2},
			\end{split}
		\end{equation}
		where we used \eqref{eq:S_bound}.
		Therefore, it follows from the entry-wise bound in \eqref{eq:apriori_flat} and \eqref{eq:omega0_s} that the estimates
		\begin{equation} \label{eq:omega_induct_base}
			\bigl\lvert (G_j-m_j)_{ab} \bigr\rvert \le N^{\zeta}\,\mathfrak{s}_1(\omega_0; a,b), \quad \bigl\lvert \bigl((G-M)_{[1,2]}(x)\bigr)_{ab} \bigr\rvert \lesssim \sqrt{W} \,\mathfrak{s}_2(\omega_0; a,x,b),
		\end{equation}
		hold with very high probability, for all $ z_1, z_2  \in \{z, \overline{z}\}$, and $a,b,x \in \indset{N}$. That is, the resolvent $G$ satisfies the assumption of Lemma \ref{lemma:omega_lemma} with $\omega := \omega_0 = W^{-1}$ and  $\psi:= \sqrt{W}$. 
		
		Using~\eqref{eq:omega_induct_base} to initialize the iteration, we use Lemma \ref{lemma:omega_lemma} inductively to conclude that, for any  $l \in \indset{2D'/\zeta}$, the bounds \eqref{eq:omega_assume} holds with $\omega := N^{-l\zeta}$ and $\psi := lN^{10\zeta} + N^{-l\zeta}\sqrt{W}$.
		Hence, it follows from \eqref{eq:omega_final_s} that, with very high probability,
		\begin{equation} 
			\bigl\lvert (G_{j}-m_{j})_{ab} \bigr\rvert \lesssim N^{\zeta}\mathfrak{s}_1^\mathrm{iso}(a,b), \quad \bigl\lvert \bigl((G-M)_{[1,2]}(x)\bigr)_{ab} \bigr\rvert \lesssim N^{10\zeta} \mathfrak{s}_2^\mathrm{iso}(a,x,b), \quad a,b,x \in \indset{N},
		\end{equation}
		Therefore, since $\zeta>0$ can be chosen arbitrarily, \eqref{eq:1G_global_decay} holds. This concludes the proof of Lemma~\ref{lemma:1G_global_decay}.
	\end{proof}
	
	Therefore, it remains to prove Lemma \ref{lemma:omega_lemma}.
	It follows from \eqref{eq:Ups_majorates} and the second estimate in \eqref{eq:omega_assume}, and the assumption $\psi \le \sqrt{W}$,  that
	\begin{equation} \label{eq:2G_omega_ansatz}
		\bigl\lvert (G^*S^xG)_{bb} \bigr\rvert \le \Theta_{xb} + \bigl\lvert (G^*S^xG - \Theta^x)_{bb} \bigr\rvert \lesssim \Upsilon_{xb} + \psi \mathfrak{s}_2(\omega; b,x,b) \lesssim \Upsilon_{xb} + \omega, \quad \text{w.v.h.p.},\\
	\end{equation}
	Similarly, it follows from \eqref{eq:omega_assume}, that
	\begin{equation} \label{eq:omega_chain_ansatz}
		 \lvert G_{ab}  \rvert \le \delta_{ab}+ N^{2\zeta} \mathfrak{s}_1(\omega; a,b), \quad \bigl\lvert (G^*S^xG)_{ab} \bigr\rvert \le \delta_{ab}\Upsilon_{xa} + \psi \mathfrak{s}_2(\omega; b,a,b), \quad \text{w.v.h.p}.
	\end{equation}
	
	\begin{proof}[Proof of Lemma \ref{lemma:omega_lemma}]
		It follows from \eqref{eq:Ups_norm_bounds}, \eqref{eq:triag}, and \eqref{eq:true_convol}, that the structural blocks of $\mathfrak{s}_k(\omega;\,\cdot\,)$ satisfy
		\begin{equation}\label{eq:omega_s_triag}
			\max_a \sqrt{(\Upsilon_{xa} + \omega)(\Upsilon_{ay} + \omega)} \lesssim W^{-1/2}\sqrt{\Upsilon_{xy} + \omega}, \quad x,y \in \indset{N}, \quad \omega \in [0, W^{-1}],
		\end{equation}
		\begin{equation} \label{eq:omega_s_convol}
			\sum_a \Upsilon_{ba} \sqrt{(\Upsilon_{xa} + \omega)(\Upsilon_{ay} + \omega)} \lesssim \sqrt{(\Upsilon_{xb} + \omega)(\Upsilon_{by} + \omega)}, \quad b,x,y \in \indset{N}.
		\end{equation}
		That is, the analogs of \eqref{eq:triag} and \eqref{eq:true_convol} hold for the modified control function $\Upsilon + \omega$. 
		Moreover, for any constant $0 \le \lambda \le 1$, we have
		\begin{equation} \label{eq:omega_ineqs}
			\lambda^{k/2} \mathfrak{s}_k(\omega; \,\cdot\,) \le  \mathfrak{s}_k(\lambda \omega; \,\cdot\,) \le \mathfrak{s}_k(\omega; \,\cdot\,), \quad k \in \{1,2\}.
		\end{equation}
		
		\textbf{Decay for a single resolvent}.
		First, we prove the estimate on the single resolvent in \eqref{eq:omega_conclude}. To this end, we use \eqref{eq:static1G}. Note that for $a=b$, the estimate \eqref{eq:apriori_flat} is already sharp and it suffices to  consider the case $a \neq b$. Hence, $\mathcal{B}^{-1}[X]_{ab} = X_{ab}$ and the second term on the right-hand side of \eqref{eq:static1G} is bounded by
		\begin{equation} \label{eq:1G_decay_quad_term}
			\bigl\lvert m \,\mathcal{B}^{-1}\bigl[ \mathcal{S}[G-m](G-m)\bigr]_{ab} \bigr\rvert \lesssim W^{-3/4}N^{2\zeta}\mathfrak{s}_1(\omega; a,b) \lesssim W^{-3/4}N^{5\zeta/2}\mathfrak{s}_1(N^{-\zeta}\omega; a,b), \quad \text{w.v.h.p.},
		\end{equation}
		where we used \eqref{eq:apriori_av}, the assumed bounds \eqref{eq:omega_assume}, and \eqref{eq:omega_ineqs}.
		
		To apply Lemma \ref{lemma:high-moment}, we show that the underline term in \eqref{eq:static1G} satisfies the following bound.
		\begin{claim} \label{claim:omega_underline1} Under the assumptions of Lemma \ref{lemma:omega_lemma}, for any integer $p \in \mathbb{N}$, we have, for all $a\neq b$,
			\begin{equation} \label{eq:1G_decay_high_moment}
				\biggl\lvert\Expv\biggl[\mathcal{B}^{-1}[\underline{HG}]_{ab}\overline{G_{ab}} \lvert G_{ab} \rvert^{2p-2}\biggr]\biggr\rvert \lesssim \biggl( N^{\zeta/2} \, \mathfrak{s}_1(N^{-\zeta}\omega; a,b)\biggr)^{2p} + N^{-\zeta}\Expv\bigl[ \lvert G_{ab} \rvert^{2p}\bigr].
			\end{equation}
		\end{claim}
		We prove Claim \ref{claim:omega_underline1} in Section \ref{sec:cum_expand}. 
		Together with \eqref{eq:static1G}, \eqref{eq:1G_decay_quad_term}  and \eqref{eq:1G_decay_high_moment}, Lemma \ref{lemma:high-moment} implies 
		\begin{equation} \label{eq:1G_prec_from_high_moment}
			\bigl\lvert (G-m)_{ab} \bigr\rvert \lesssim N^{\zeta}\mathfrak{s}_1\bigl(N^{-\zeta}\omega; a,b\bigr), \quad a\neq b \in \indset{N},\quad  \text{w.v.h.p.}
		\end{equation}
		Therefore, since the $a=b$ term can be estimated by \eqref{eq:apriori_flat}, the first estimate in \eqref{eq:omega_conclude} is established.
		
		\textbf{Decay for a two-resolvent chain}.
		We now prove the second bound in \eqref{eq:omega_conclude}. To estimate the two-resolvent chain, we use the following generalization of the static equation \eqref{eq:static1G}, that follows from \eqref{eq:underline_def} and \eqref{eq:M_recursion}. Let $A' \equiv A'(z_1,z_2,x) := m_2 (1-m_1m_2\mathcal{S})^{-1}[S^x]$, then
		\begin{equation} \label{eq:2G_static}
			\begin{split}
				(G-M)_{[1,2]} =&~ -\underline{G_{1}A'HG_{2}} 
				+ (G_1-m_1)A' 
				+G_{1} \mathcal{S}\bigl[(G_1-m_1)A'\bigr]G_{2} 
				+ G_{1}A'\mathcal{S}\bigl[G_{2}-m_{2}\bigr]G_{2}.
			\end{split}
		\end{equation}
		It follows from \eqref{eq:Ups_majorates} that
		\begin{equation} \label{eq:A'_bound}
			\bigl\lvert A'_{ab} \bigr\rvert \lesssim \delta_{ab}\Upsilon_{xa}, \quad a,b\in \indset{N}.
		\end{equation}
		
		We begin by estimating the underline-free terms on the right-hand side of \eqref{eq:2G_static}. It follows from \eqref{eq:omega_assume} and  \eqref{eq:A'_bound} that the second term on the right-hand side of \eqref{eq:2G_static} satisfies
		\begin{equation} \label{eq:2G_forcing1}
			\bigl\lvert \bigl((G_1-m_1)A'\bigr)_{ab} \bigr\rvert \lesssim  N^{2\zeta} \Upsilon_{xa} \mathfrak{s}_{1}(\omega; a,b) \lesssim  W^{-1/2} N^{5\zeta/2}   \mathfrak{s}_{2}(N^{-\zeta}\omega; a,x,b),
		\end{equation} 
		where we used  \eqref{eq:omegas_def}, \eqref{eq:omega_s_triag}, and \eqref{eq:omega_ineqs}.

		By \eqref{eq:apriori_flat}, the third term on the right-hand side of \eqref{eq:2G_static} admits the bound
		\begin{equation} \label{eq:2G_forcing2}
			\begin{split}
				\bigl\lvert \bigl(G_{1} \mathcal{S}\bigl[(G_1-m_1)A'\bigr]G_{2} \bigr)_{ab} \bigr\rvert &= \biggl\lvert\sum_{q} \bigl(G_{1}S^qG_{2}\bigr)_{ab} (G_1-m_1)_{qq} A'_{qq} \biggr\rvert\\
				&\lesssim  \frac{N^{\zeta}}{\sqrt{W}}\sum_q \biggl(\bigl\lvert \bigl(M_{[1,2]}(q)\bigr)_{ab} \bigr\rvert + \psi \mathfrak{s}_2(\omega; a, q, b)\biggr) \bigl\lvert A'_{qq}\bigr\rvert\\
				&\lesssim  \biggl(N^{\zeta}+\frac{N^{2\zeta}\psi}{\sqrt{W}}\biggr)  \mathfrak{s}_2(N^{-\zeta}\omega; a, x, b),
			\end{split}
		\end{equation}
		with very high probability, where we used \eqref{eq:M_bound}, \eqref{eq:true_convol}, \eqref{eq:omegas_def}, and \eqref{eq:omega_s_convol}. Finally, by \eqref{eq:apriori_av}, the fourth term on the right-hand side of \eqref{eq:2G_static} admits the bound
		\begin{equation} \label{eq:2G_forcing3}
			\begin{split}
				\bigl\lvert \bigl(G_{1}A'\mathcal{S}\bigl[G_{2}-m_{2}\bigr]G_{2}\bigr)_{ab}\bigr\rvert &= \biggl\lvert \sum_q (G_1)_{aq} A'_{qq} \Tr\bigl[(G_2-m_2)S^q\bigr] (G_2)_{qb}\biggr\rvert\\
				&\lesssim \frac{N^\zeta}{W^{3/4}} \sum_q \bigl\lvert A'_{qq}\bigr\rvert |(G_{1})_{aq}||(G_{2})_{qb}|\\
				&\lesssim \frac{N^\zeta}{W^{3/4}} \sum_q \Upsilon_{xq} \bigl(\delta_{aq} + N^{2\zeta} \mathfrak{s}_1(\omega; a,q)\bigr)\bigl(\delta_{qb} + N^{2\zeta} \mathfrak{s}_1(\omega; q,b)\bigr)\\
				&\lesssim \frac{N^{6\zeta}}{W^{1/4}} \mathfrak{s}_2(N^{-\zeta}\omega; a, x, b),
			\end{split}
		\end{equation}
		where, to obtain the final bound, we used \eqref{eq:omega_chain_ansatz}, \eqref{eq:omega_s_triag}, \eqref{eq:omega_s_convol}, together with \eqref{eq:omegas_def} and \eqref{eq:omega_ineqs}.
		
		Finally, the underline term in \eqref{eq:2G_static} satisfies the following claim.
		\begin{claim} \label{claim:omega_underline2} Let $\mathcal{Y}_1 := ((G-M)_{[1,2]}(x))_{ab}$. Then, under the assumptions of Lemma \ref{lemma:omega_lemma}, for any integer $p \in \mathbb{N}$, we have 
			\begin{equation} \label{eq:omega_underline2}
				\biggl\lvert\Expv\biggl[\bigl(\underline{G_{1}A'HG_{2}}\bigr)_{ab} \overline{\mathcal{Y}_1} \lvert \mathcal{Y}_1 \rvert^{2p-2}\biggr]\biggr\rvert 
				\lesssim \biggl(N^{4\zeta}\sqrt{\psi} \, \mathfrak{s}_2(N^{-\zeta}\omega; a,x,b)\biggr)^{2p} +N^{-\zeta}\Expv\bigl[\lvert \mathcal{Y}_1\rvert^{2p}\bigr].
			\end{equation}
		\end{claim}
		We prove Claim \ref{claim:omega_underline2} in Section \ref{sec:cum_expand}. The second estimate in \eqref{eq:omega_conclude} follows now
		 from Lemma \ref{lemma:high-moment}, \eqref{eq:2G_static}, \eqref{eq:2G_forcing1}--\eqref{eq:2G_forcing3}, and \eqref{eq:omega_underline2}, similarly to \eqref{eq:1G_prec_from_high_moment}. 
		 Note that the control parameter $\psi$ comes either with a small prefactor, as in~\eqref{eq:2G_forcing2}, 
		 or with a square root, as in~\eqref{eq:omega_underline2}; both give rise to the iterable  $N^{-\zeta}\psi$ 
		 factor in~\eqref{eq:omega_conclude}. 		 This concludes the proof of Lemma \ref{lemma:omega_lemma}.
	\end{proof}
	
	\subsection{Step 3. Chains of arbitrary length} \label{sec:global3}
	In this section, we prove Proposition \ref{prop:global_laws} for chains of arbitrary length $k \in \indset{\maxK}$. Recall that the isotropic estimate \eqref{eq:iso_global} was already established for $k\in\{1,2\}$ in Section \ref{sec:glob_decay} above ( see Lemma \ref{lemma:1G_global_decay}). In the present section, we show that \eqref{eq:1G_global_decay} is a sufficient input to prove optimal decaying estimate for both isotropic and averaged chains of length $k\in\indset{\maxK}$ directly, without performing the $\omega$-bootstrap, which was only necessary for proving Lemma \ref{lemma:1G_global_decay} in Step 2.
	Furthermore, there is no need to handle $k$ and $k+1$ chains together
	unlike the simultaneous treatment of the $k=1$ and $k=2$ chains in Section~\ref{sec:glob_decay}. 
	However, an induction on $k$ will still be used. 
	
	In the sequel we consider size functions $\mathfrak{s}_k^{\mathrm{iso/av}}$ given by
	\begin{equation} \label{eq:global_size_funcs}
		\mathfrak{s}_k^{\mathrm{iso/av}}(\,\cdot\,) := \mathfrak{s}_{k,t=0}^{\mathrm{iso/av}}(\,\cdot\,), \quad k\in \mathbb{N},
	\end{equation}
	with $\ell_0 \sim W$ and $\eta_0 \sim 1$. Moreover, the corresponding control function $\Upsilon := \Upsilon_{t=0}$, as in Definition \ref{def:adm_ups}, satisfies the following analogs of \eqref{eq:triag} and \eqref{eq:Schwarz_convol},
	\begin{equation} \label{eq:global_triag}
		\max_a \sqrt{\Upsilon_{xa} \Upsilon_{ay}} \lesssim W^{-1/2}\sqrt{ \Upsilon_{xy}}, 
		\quad x,y\in\indset{N}, 
	\end{equation}
	\begin{equation} \label{eq:global_convol}
		\sum_a \sqrt{\Upsilon_{x_1a} \Upsilon_{ay_1}\Upsilon_{x_2a} \Upsilon_{ay_2}}  \lesssim  \sqrt{\Upsilon_{x_1y_1}\Upsilon_{x_2y_2}}, \quad x_1,x_2,y_1,y_2\in\indset{N}.  
	\end{equation}
	
	We denote the resolvent chains $G_{[p,j]} \equiv G_{[p,j],0}(\bm z; A_{p},\dots, A_{j-1})$ and their corresponding deterministic approximations $M_{[p,j]} \equiv M(\bm z; A_{p},\dots, A_{j-1})$ as in  \eqref{eq:resolvent_chains} and \eqref{eq:Mk_def}, respectively, for $k\in\mathbb{N}$, $p\le j \in \indset{k}$, and  $\bm z \in \{z,\overline{z}\}^k$.  
	 
	As in Section \ref{sec:glob_decay} above, our main tool for estimating $(G-M)_{[1,k]}$ is the multi-resolvent analog of the static equation \eqref{eq:static1G}. We will prove both averaged and isotropic bounds. Now we setup the precise
	objects that the master inequalities will control. 
	
	 In the averaged case\footnote{
		We present the more general matrix-valued version of \eqref{eq:Gk_static_av}, used in the isotropic case, in \eqref{eq:Gk_static} below. 
		} for $k \ge 2$, we deduce from \eqref{eq:underline_def} and \eqref{eq:M_recursion} that
	\begin{equation} \label{eq:Gk_static_av}
		\begin{split}
			\Tr\bigl[(G-M)_{[1,k]}A_k\bigr] =&~ -\Tr\biggl[\underline{HG_{[1,k]}}A_k'\biggr] + \Tr\bigl[(G-M)_{[2,k]}A_k'A_1\bigr]   \\
			&+ \sum_{j=2}^{k-1}\Tr\biggl[\mathcal{S}\bigl[(G-M)_{[1,j]}\bigr]G_{[j,k]}A_k'\biggr] +\sum_{j=2}^{k-1}\Tr\biggl[\mathcal{S}\bigl[(G-M)_{[j,k]}A_k'\bigr]M_{[1,j]}\biggr] \\
			&+ \Tr\bigl[\mathcal{S}\bigl[G_1-m_1\bigr]G_{[1,k]}A_k'\bigr] + \Tr\bigl[\mathcal{S}\bigl[(G_k-m_k)A_k'\bigr]G_{[1,k]}\bigr],
		\end{split}
	\end{equation}
	where $A_k' := m_1\bigl(1 -m_1 m_k \mathcal{S}\bigr)^{-1}[A_{k}]$. As can be seen from the second term on the right-hand side of \eqref{eq:Gk_static_av}, up to applying the stability operator $\bigl(1 -m_1 m_k \mathcal{S}\bigr)^{-1}$ to $A_j$, this equation also produces observables that are products of the original $A_i$'s. To account for this, we define the following hierarchy of observable classes $\{\mathcal{V}_k\}_{k=1}^{\maxK}$. Let $\mathcal{V}_1$ be a set of diagonal matrices depending on a single external index $x_1$, given by
	\begin{equation} \label{eq:V1_class}
		\mathcal{V}_1 := \bigl\{ A : x_1 \mapsto  S^{x_1}, \quad x_1 \in \indset{N} \bigr\}.
	\end{equation} 
	To construct $\mathcal{V}_{k}$ with $k \ge 2$, we first defined the set of all possible operators that map $A_k$ to $A_k'$ in \eqref{eq:Gk_static_av}, that is
	\begin{equation}
		\mathfrak{A} := \bigl\{\, \other{m}_1\bigl(1 -\other{m}_1 \other{m}_2 \mathcal{S}\bigr)^{-1} \, :\, \other{m}_1,\other{m}_2 \in \bigl\{m(z), m(\overline{z}) \bigr\} \,\bigr\}.
	\end{equation}
	Then, for $k \in \indset{2, \maxK }$, we define
	\begin{equation} \label{eq:Vk_class}
		\mathcal{V}_{k} := \bigl\{ A: \bm x \mapsto \mathcal{A}\bigl[A_{k-1}(x_2, \dots, x_{k})\bigr]S^{x_{1}}, \quad \bm x = (x_1,\dots, x_{k}) \in \indset{N}^{k} \, : \, A_{k-1} \in \mathcal{V}_{k-1}, \, \mathcal{A} \in \mathfrak{A}  \bigr\}.
	\end{equation}
	For example, if $A_1 = S^{x_1}$ and $A_2 = S^{x_2}$, then the product $A_2' A_1$ belongs to $\mathcal{V}_2$. 
	The emergence of the new classes $\mathcal{V}_{k}$ of somewhat more involved observables is 
	a feature of the averaged proof (for $k\ge 2$) stemming from the product $A_k'A_1$ in~\eqref{eq:Gk_static_av}. 
	While they complicate the notation, we stress that conceptually they do not involve any extra difficulty. The proof 
	only uses their spatial decay properties that are natural and can be easily established.
	
	Every  resolvent chain we consider contains at most one distinguished observable belonging to $\mathcal{V}_k$ with $k\ge 2$, while all other observables are from $\mathcal{V}_1$, i.e., the usual $S^{x_i}$. Recall that the spectral parameter $z := E+ \I\eta$ with $\eta \sim 1$ and $|E| \le C$ and the maximal chains length $\maxK \in \mathbb{N}$ are fixed.  
	Also recall the definitions~\eqref{eq:resolvent_chains}--\eqref{eq:Mt_def} that will be used for $t=0$ (and $t$ is omitted
	from the notation).
	 We define the corresponding control quantities $\Phi_{k}^\mathrm{av}$ and $\Phi_{k+1}^\mathrm{iso}$, for any $k$, as
	\begin{subequations}\label{eq:static_Phi_def}
		\begin{equation}  \label{eq:static_Phi_av_def}
			\Phi_{k}^\mathrm{av} := \max_{\bm z \in \{z, \overline{z}\}^k} \max_{\extk\in \indset{\maxK - k+1}} \max_{A_k \in \mathcal{V}_\extk} \max_{\bm x \in \indset{N}^{k+\extk-1}} \frac{\bigl\lvert \Tr\bigl[(G-M)_{[1,k]} (A_1^k, \dots, A_{k}^k)\bigr] \bigr\rvert}{\mathfrak{s}_{k+\extk-1}^\mathrm{av}(\bm x)},
		\end{equation}
		\begin{equation} \label{eq:static_Phi_iso_def}
			\Phi_{k+1}^\mathrm{iso} := \max_{a,b}\max_{\bm z \in \{z, \overline{z}\}^{k+1}} \max_{\extk\in \indset{(\maxK - k)\vee 1}} \max_{j\in\indset{k}}  \max_{A_j \in \mathcal{V}_\extk}  \max_{\bm x \in \indset{N}^{k+\extk-1}} \frac{\bigl\lvert \bigl((G-M)_{[1,k+1]}(A_1^j, \dots, A_{k}^j)\bigr)_{ab} \bigr\rvert}{\mathfrak{s}_{k+\extk}^\mathrm{iso}(a,\bm x,b)},
		\end{equation}
	\end{subequations}
	where, for all $i,j \in \indset{k}$, $\bm x \in \indset{N}^{k+\extk-1}$, and an observable function $A_j$ dependent on $\extk$ indices, we define the observables\footnote{
	Here, with a mild abuse of notation (moving the index $j$ from $A_j$ to a superscript),
	 we use the superscript $j$ of $A_{i}^{j}$ to signify the position of the distinguished observable $A_j(x_j,\dots, x_{j+\extk-1})$ in the resolvent chain. } $A_i^j \equiv A_i^j(\bm x)$ as
	\begin{equation} \label{eq:static_observables}
		A_i^j(\bm x) := \begin{cases}
			S^{x_i}, \quad & i < j,\\
			A_j(x_j, \dots, x_{j+\extk-1}), \quad & i = j,\\
			S^{x_{i+\extk-1}}, \quad & i > j.\\
		\end{cases}
	\end{equation}
	
	In the averaged case, by cyclicity, it suffices to consider the case with the last observable $A_k$ being the distinguished observable. However, in the isotropic case such simplification is not available, and we have to allow the distinguished observable $A_j$ to occupy any position $j \in \indset{k}$ in the chain $G_{[1,k+1]}$ hence the complicated definition of $\Phi_{k+1}^\mathrm{iso}$ in \eqref{eq:static_Phi_iso_def}. Note that the arguments of the size functions $\mathfrak{s}_{k+\extk-1}^\mathrm{av}$ and $\mathfrak{s}_{k+\extk}^\mathrm{iso}$ in the denominators of \eqref{eq:static_Phi_av_def} and \eqref{eq:static_Phi_iso_def}, respectively, contain all $k+\extk -1$ indices of $\bm x$: $k-1$ indices coming from the observables $A_i^j = S^{x_{\dots}}$ with $i\neq j$, and $\extk$ indices coming from the distinguished observable $A_j$. Finally, we note that the maxima in the definition of $\Phi_{k}^\mathrm{iso/av}$ in \eqref{eq:static_Phi_def} are taken over a set of polynomial cardinality in $N$, hence they are compatible with taking a union bound in the definition of stochastic domination.
	
	For ease of notation, we define
	\begin{equation}
		\Phi^\mathrm{iso}_1 := \Psi_{1,0}^\mathrm{iso} = \max_{a,b\in\indset{N}} \frac{\bigl\lvert (G-m)_{ab} \bigr\rvert}{\mathfrak{s}_{1}^\mathrm{iso}(a,b)},
	\end{equation}

	The goal of this section is to prove the following propositions.
	\begin{prop}[Static Master Inequalities for Isotropic Global Laws] \label{prop:iso_static_masters}
		Fix an integer $k \in \indset{\maxK}$.
		Let $\pis{k+1}$ be a deterministic control parameter satisfying $\pis{k+1}\ge 1$.
		Assume that 
		\begin{equation} \label{eq:iso_static_assume}
			\Phi_{k'}^\mathrm{iso} \prec 1, \quad k' \in \indset{k}, \quad    \Phi_{k+1}^\mathrm{iso} \prec \psi_{k+1}^\mathrm{iso}.
		\end{equation}	
		Then $\Phi_{k+1}^{\mathrm{iso}}$  satisfies the improved bound
		\begin{equation}
			\Phi_{k+1}^\mathrm{iso} \prec 1 +  \frac{\pis{k+1}}{W^{1/4}} + \sqrt{\pis{k+1}}. \label{eq:static_iso_masters}
		\end{equation}
	\end{prop}
	\begin{prop}[Averaged Global Laws] \label{prop:static_masters}
		Fix an integer $k \in \indset{\maxK}$.
		Assume that 
		\begin{equation} \label{eq:av_static_assume}
			\Phi_{k'}^\mathrm{av} \prec 1, \quad k' \in \indset{k-1}, \quad \text{and} \quad \Phi_{k'}^\mathrm{iso} \prec 1, \quad k' \in \indset{k+1}. 
		\end{equation}	
		Then $\Phi_{k}^{\mathrm{av}} \prec 1$ satisfies the stochastic domination bound
		\begin{equation}
			\Phi_{k}^{\mathrm{av}} \prec 1. \label{eq:static_av_masters}
		\end{equation}
	\end{prop}
	We defer the proof of Propositions \ref{prop:iso_static_masters} and \ref{prop:static_masters} to the end of this section. Propositions \ref{prop:iso_static_masters} and \ref{prop:static_masters} comprise an "economical" version of master inequalities, that allow us to prove global laws by strong induction in the chain length $k$.
	Furthermore, notice that the isotropic master inequalities involve only isotropic quantities and thus can be
	solved without any averaged bounds. Averaged bounds are then proven separately. In typical previous applications
	of the static master inequalities \cite{Cipolloni2022Optimal, cipolloni2022rank} averaged and isotropic bounds were proven 
	in tandem. Now, being in the global regime $\eta\sim 1$ simplifies the procedure and allows decoupling the
	isotropic master inequalities from the averaged ones. Note, however, that the  bound for the averaged $k$-chain 
	requires a bound on the isotropic $(k+1)$-chain as an input.
	
	\begin{proof}[Proof of Proposition \ref{prop:global_laws}]
		By definition of $\Phi_{k}^\mathrm{iso/av}$ from \eqref{eq:static_Phi_def}, to prove \eqref{eq:av_global} and \eqref{eq:iso_global}, it suffices to show that $\Phi_{k}^\mathrm{iso/av} \prec 1$ for $k\in\indset{\maxK}$. 
		
		Clearly, assuming that \eqref{eq:iso_static_assume} holds, \eqref{eq:static_iso_masters} and a standard iteration argument implies that $\Phi_{k+1}^\mathrm{iso} \prec 1$.		
		Hence, using strong induction with $\Phi_1^\mathrm{iso} \prec 1$ (that follows from Lemma \ref{lemma:1G_global_decay}) as the base, we prove that $\Phi_{k+1}^\mathrm{iso} \prec 1$ for all $k \in \indset{\maxK}$. 
		
		Analogously, using strong induction and the isotropic estimates $\Phi_{k'}^\mathrm{iso} \prec 1$ with $k'\in \indset{\maxK+1}$ as an input, we prove that $\Phi_k^\mathrm{av} \prec 1$ for all $k \in \indset{\maxK}$. This concludes the proof of Proposition \ref{prop:global_laws}.
	\end{proof}
	Therefore, it remains to prove Propositions \ref{prop:iso_static_masters} and \ref{prop:static_masters}.	
	Before proceeding with the proof, we collect the preliminary bounds on the corresponding $M$-terms.
	\begin{claim}[Global $M$ Bounds with Distinguished Observables] \label{claim:staticM} Let $\extk \in \indset{\maxK}$, then, for any $A \in \mathcal{V}_\extk$ defined in \eqref{eq:V1_class} and \eqref{eq:Vk_class}, we have
		\begin{equation} \label{eq:V_class_bound}
			\bigl\lvert \bigl( A(\bm y) \bigr)_{ab} \bigr\rvert \lesssim \delta_{ab}\sqrt{W} \, \mathfrak{s}_{\extk+1}^\mathrm{iso}(a,\bm y, b), \quad a,b \in \indset{N}, \quad \bm y \in \indset{N}^\extk,
		\end{equation}
		\begin{equation} \label{eq:V_class_av}
		 \sum_{q}\bigl\lvert \bigl( A(\bm y) \bigr)_{qq} \bigr\rvert \lesssim W \, \mathfrak{s}_{\extk}^\mathrm{av}(\bm y), \quad \bm y \in \indset{N}^\extk.  
		\end{equation}
		Furthermore, for any $k \in \indset{\maxK}$ and $\extk \in \indset{(\maxK-k)\vee 1}$, 
		\begin{equation} \label{eq:staticM}
			\bigl\lvert \bigl(M_{[1,k+1]}(A_1^j, \dots, A_{k}^j)\bigr)_{ab} \bigr\rvert \lesssim \delta_{ab}\sqrt{W}\,\mathfrak{s}_{k+\extk}^\mathrm{iso}(a,\bm x,b), \quad a,b \in \indset{N}, \quad j \in \indset{k},
		\end{equation}
		where $A_{i}^j$ with $i,j \in \indset{k}$ are defined in \eqref{eq:static_observables}.
	\end{claim}
	We prove  Claim \ref{claim:staticM} in Section \ref{sec:M_bounds}. 
	Equipped with Claim \ref{claim:staticM}, we  prove Proposition \ref{prop:static_masters},
	\begin{proof}[Proof of Proposition \ref{prop:static_masters}]
		It follows from the definition of $\Phi_{k'}^{\mathrm{av/iso}}$ in  \eqref{eq:static_Phi_def}, \eqref{eq:av_static_assume},  that
		\begin{equation} \label{eq:static_G-M_ansatz}
			\begin{split}
				\bigl\lvert \Tr\bigl[(G-M)_{[1,k']} (A_1^{k'}, \dots, A_{k'}^{k'})\bigr] \bigr\rvert &\prec \mathfrak{s}_{k'+\extk-1}^\mathrm{av}(\bm x), \quad\extk \in \indset{\maxK-k' + 1}, \quad k' \in \indset{k-1}, \\
				\bigl\lvert \bigl((G-M)_{[1,k'+1]}(A_1^j, \dots, A_{k'}^j)\bigr)_{ab} \bigr\rvert &\prec  \mathfrak{s}_{k'+\extk}^\mathrm{iso}(a,\bm x,b), \quad \extk \in \indset{(\maxK-k')\vee 1}, \quad k'\in\indset{k},
			\end{split}
		\end{equation}
		for all $a,b \in \indset{N}$,  $\bm x\in\indset{N}^{k'+\extk-1}$, $j \in \indset{k'}$, and $A_{i}^j$ defined in \eqref{eq:static_observables}. Hence, under the same assumptions, \eqref{eq:1G_global_decay}, \eqref{eq:staticM} and \eqref{eq:static_G-M_ansatz} imply that 
		\begin{equation} \label{eq:static_G_bound}
			\lvert G_{ab}  \rvert \prec \delta_{ab} + \mathfrak{s}_{1}^\mathrm{iso}(a,b),\quad	\bigl\lvert \bigl(G_{[1,k'+1]}(A_1^j, \dots, A_{k'}^j)\bigr)_{ab} \bigr\rvert \prec \bigl(\delta_{ab}\sqrt{W} + 1\bigr)\,\mathfrak{s}_{k'+\extk}^\mathrm{iso}(a,\bm x,b), \quad k'\in\indset{k}.
		\end{equation}

		Fix   $\extk\in \indset{\maxK - k + 1}$. Let $\bm x := (x_1,\dots, x_{k+\extk-1}) \in \indset{N}^{k+\extk-1}$, and $\bm y := (x_k, \dots, x_{k+\extk-1}) \in \indset{N}^{\extk}$, such that $(\bm x^{(\extk)}, \bm y) = \bm x$.
		Let  $A_k$ be a fixed observable function, belonging to $\mathcal{V}_\extk$, defined in \eqref{eq:V1_class} and \eqref{eq:Vk_class}.  We abbreviate $A_k \equiv A_k(\bm y)$. Let $\bm z  \in \{z, \overline{z} \}^k$ be a fixed vector of spectral parameters. Our goal is to estimate $\lvert \Tr[ (G-M)_{[1,k]}(\bm x^{(\extk)})\, A_k ] \rvert$. Since $A_k \in \mathcal{V}_{\extk}$, it follows from \eqref{eq:Binv_action} and \eqref{eq:V_class_bound} that 
		\begin{equation} \label{eq:Ak'_bound}
			\bigl\lvert \bigl(A_{k}'\bigr)_{ab} \bigr\rvert \lesssim \delta_{ab} \sqrt{W} \mathfrak{s}_{\extk+1}^\mathrm{iso} 
			(b,\bm y, b).
		\end{equation} 
		
		We consider two cases: $k=1$ and $k \ge 2$, for which we use equations \eqref{eq:1G_av_static} and \eqref{eq:Gk_static_av}, respectively. Since the underline terms in  \eqref{eq:1G_av_static} and \eqref{eq:Gk_static_av}
		are the same, we only draw the distinction between the two cases when estimating the underline-free terms.
		
		\textbf{Case 1.} Assume that $k=1$. Note that in this case $\bm x = \bm y$ and $k+\extk-1 = \extk$. 
		Hence, from \eqref{eq:apriori_av} and \eqref{eq:1G_global_decay}, we deduce  that the only underline-free term on the right-hand side of \eqref{eq:1G_av_static} satisfies 
		\begin{equation} \label{eq:1av_quad_term_static}
			\begin{split}
				\bigl\lvert \Tr\bigl[  \mathcal{S}[G-m](G-m)A_1'\bigr] \bigr\rvert &= \biggl\lvert\sum_{q} \Tr\bigl[(G-m)S^q\bigr] (G-m)_{qq} \bigl(A_1'\bigr)_{qq} \biggr\rvert\\
				&\prec W^{-3/4} W^{-1/2}\sum_{q} \bigl\lvert \bigl(A_1'\bigr)_{qq} \bigr\rvert \prec W^{-1/4}  \mathfrak{s}_{\extk}^\mathrm{av}(\bm x).
			\end{split}
		\end{equation} 
		where in the last step we used  \eqref{eq:V_class_av}.
		
		\textbf{Case 2.} Assume now that $k \ge 2$. Then our goal is to estimate the underline-free terms on the right-hand side of \eqref{eq:Gk_static_av}. Since $A_1 \in \mathcal{V}_1$ and $A_k \in \mathcal{V}_\extk$, the matrix $A_k'A_1 \in \mathcal{V}_{\extk+1}$ by definition \eqref{eq:Vk_class}. Hence, it follows from \eqref{eq:av_static_assume} and \eqref{eq:static_G-M_ansatz} that the second term on the right-hand side of \eqref{eq:Gk_static_av} admits the bound
		\begin{equation} \label{eq:static_k-1_term}
			\bigl\lvert \Tr\bigl[(G-M)_{[2,k]}A_k'A_1\bigr] \bigr\rvert \prec \mathfrak{s}_{(k-1)+(\extk+1)-1}^\mathrm{av}(x_2,\dots, x_{k+\extk-1},x_1) =  \mathfrak{s}_{k+\extk-1}^\mathrm{av}(\bm x).
		\end{equation}
		Similarly to \eqref{eq:quad_estimate} and \eqref{eq:sub_final} with $\eta \sim 1,$ $\ell \sim W$ and $\beta_j = 0$, using \eqref{eq:av_static_assume},  \eqref{eq:V_class_bound}, \eqref{eq:static_G-M_ansatz}, and \eqref{eq:static_G_bound}, we deduce that the third term on the right-hand side of \eqref{eq:Gk_static_av}
		\begin{equation} \label{eq:static_quad_terms}
			\begin{split}
				\biggl\lvert\sum_{i=2}^{k-1}\Tr\biggl[\mathcal{S}\bigl[(G-M)_{[1,i]}\bigr]G_{[i,k]}A_k'\biggr]\biggr\rvert \prec&~ W\sum_{i=2}^{k-1}\sum_q  \mathfrak{s}^\mathrm{av}_{i}(x_1, \dots, x_{i-1},q)  \\
				&\qquad\times \mathfrak{s}_{k-i+1}(q, x_i, \dots,x_{k-1}, q) \,\mathfrak{s}^\mathrm{iso}_{\extk+1}(q,\bm y,q)\\
				&\prec  \mathfrak{s}_{k+\extk-1}^\mathrm{av}(\bm x),
			\end{split}
		\end{equation}
		where in the last step we used \eqref{eq:global_size_funcs}, \eqref{eq:global_triag}--\eqref{eq:global_convol}. Completely analogously, using \eqref{eq:staticM}, and the fact that $A_k'S^q \in \mathcal{V}_{\extk+1}$ as a function of $(\bm y,q) \in \indset{N}^{\extk+1}$,  the fourth term on the right-hand side of \eqref{eq:Gk_static_av} satisfies
		\begin{equation} \label{eq:static_sublin_terms}
			\biggl\lvert \sum_{i=2}^{k-1}\Tr\biggl[\mathcal{S}\bigl[(G-M)_{[i,k]}A_k'\bigr]M_{[1,i]}\biggr]  \biggr\rvert \prec  \mathfrak{s}_{k+\extk-1}^\mathrm{av}(\bm x).
		\end{equation}
		Similarly to \eqref{eq:1_k+1_bound1}, using  the averaged bound from \eqref{eq:static_G-M_ansatz} with $k'=1$, \eqref{eq:static_G_bound} with $k'=k$, and \eqref{eq:Ak'_bound}, we deduce that the fifth term on term on the right-hand side of \eqref{eq:Gk_static_av} admits the bound
		\begin{equation} \label{eq:static_1_k_term}
			\bigl\lvert \Tr\bigl[\mathcal{S}\bigl[G_1-m_1\bigr]G_{[1,k]}A_k'\bigr] \bigr\rvert \prec  \mathfrak{s}_{k+\extk-1}^\mathrm{av}(\bm x).
		\end{equation}
		Finally, the last term on the right-hand side of \eqref{eq:Gk_static_av} can be bounded in the same way, using \eqref{eq:static_G-M_ansatz} for $k'=1$, \eqref{eq:static_G_bound} for $k'=k$, and $A_k'S^q \in \mathcal{V}_{\extk+1}$, to obtain
		\begin{equation} \label{eq:static_k_1_term}
			\bigl\lvert \Tr\bigl[\mathcal{S}\bigl[(G_k-m_k)A_k'\bigr]G_{[1,k]}\bigr] \bigr\rvert \prec  \mathfrak{s}_{k+\extk-1}^\mathrm{av}(\bm x).
		\end{equation}
		
		Next, we estimate the contribution of the underline term for all $k$, using the following claim.
		\begin{claim} \label{claim:underline_k_av} Let $\mathcal{X}_{k} := \Tr\bigl[(G-M)_{[1,k]}A_k\bigr]$, then, under the assumptions of Proposition \ref{prop:static_masters}, for any integer $p\in\mathbb{N}$ and $0 < \nu < 1/100$, we have
			\begin{equation} \label{eq:underline_k_av} 
				\biggl\lvert \Expv \biggl[ \Tr\bigl[\underline{HG_{[1,k]}}A_k'\bigr] \cdot \overline{\mathcal{X}_{k}} \lvert \mathcal{X}_{k}\rvert^{2p-2} \biggr]  \biggr\rvert \lesssim 
				C_p \biggl( N^{2\nu} \mathfrak{s}_{k+\extk-1}^\mathrm{av}(\bm x)\biggr)^{2p} + N^{-\nu}\Expv \bigl[ \lvert \mathcal{X}_{k}\rvert^{2p} \bigr], 
			\end{equation}
			for some positive constant $C_p$ depending only on $p$.
		\end{claim}
		We prove Claim \ref{claim:underline_k_av} in Section \ref{sec:cum_expand}. It follows from 
		\eqref{eq:1G_av_static}, 
		 \eqref{eq:1av_quad_term_static} for $k=1$ and \eqref{eq:Gk_static_av}, \eqref{eq:static_k-1_term}--\eqref{eq:static_k_1_term} for $k\ge 2$, Claim \ref{claim:underline_k_av}, and Lemma \ref{lemma:high-moment}, that \eqref{eq:static_av_masters} holds. This concludes the proof of Proposition \ref{prop:static_masters}.
	\end{proof}

	Next, we prove Proposition~\ref{prop:iso_static_masters}. 	Note that we  present the proof of Proposition~\ref{prop:iso_static_masters} after the proof of Proposition~\ref{prop:static_masters} only because 
	the latter is conceptually cleaner, 
	but we do not use Proposition~\ref{prop:static_masters} or any of its assumptions along the proof
	of  Proposition~\ref{prop:iso_static_masters}.  
	\begin{proof}[Proof of Proposition \ref{prop:iso_static_masters}]
		We first collect the bounds on the resolvent chains under the assumption \eqref{eq:iso_static_assume}. It follows from the definition of $\Phi_{k'}^{\mathrm{iso}}$ in \eqref{eq:static_Phi_def} and \eqref{eq:iso_static_assume} that
		\begin{equation} \label{eq:iso_static_G-M_ansatz}
			\begin{split}
				\bigl\lvert \bigl((G-M)_{[1,k'+1]}(A_1^j, \dots, A_{k'}^j)\bigr)_{ab} \bigr\rvert &\prec  \pis{k'+1}\,\mathfrak{s}_{k'+\extk}^\mathrm{iso}(a,\bm x,b),   \quad k'\in\indset{k},
			\end{split}
		\end{equation}
		for all $a,b \in \indset{N}$, $\extk \in \indset{(\maxK-k')\vee 1}$, $\bm x\in\indset{N}^{k'+\extk-1}$, $j \in \indset{k'}$, and $A_{i}^j$ defined in \eqref{eq:static_observables}, where we set $\pis{k'} = 1$ for $k' \in \indset{k}$. Under the same assumptions, \eqref{eq:1G_global_decay}, \eqref{eq:staticM}, \eqref{eq:iso_static_G-M_ansatz} imply that 
		\begin{equation} \label{eq:iso_static_G_bound}
				\lvert G_{ab}  \rvert \prec \delta_{ab} + \pis{1}\,\mathfrak{s}_{1}^\mathrm{iso}(a,b), \quad \bigl\lvert \bigl(G_{[1,k'+1]}(A_1^j, \dots, A_{k'}^j)\bigr)_{ab} \bigr\rvert \prec \bigl(\delta_{ab}\sqrt{W} + \pis{k'+1}\bigr)\,\mathfrak{s}_{k'+\extk}^\mathrm{iso}(a,\bm x,b),
		\end{equation}
		for all $k' \in \indset{k}$.
		
		Fix an integer $\extk\in \indset{(\maxK - k)\vee 1}$. Let $\bm x := (x_1,\dots, x_{k+\extk-1}) \in \indset{N}^{k+\extk-1}$, and let $\bm x_1 := (x_1,\dots, x_{j-1})$, $\bm y := (x_j, \dots, x_{j+\extk-1}) \in \indset{N}^{\extk}$, $\bm x_2 := (x_{j+\extk},\dots, x_{k+\extk-1})$ such that $(\bm x_1, \bm y, \bm x_2) = \bm x$.
		Let $\bm z  \in \{z, \overline{z} \}^{k+1}$ be a fixed vector of spectral parameters. Let  $A_j$ be a fixed observable function, belonging to $\mathcal{V}_\extk$, defined in \eqref{eq:V1_class} and \eqref{eq:Vk_class}.  We abbreviate $A_j \equiv A_j(\bm y)$.  
		
		Our goal is to estimate $\lvert \bigl((G-M)_{[1,k+1]}(A_1^j,\dotsm A_k^j)\bigr)_{ab} \rvert$ with $A_i^j$ defined in \eqref{eq:static_observables}. To this end, we use the following static equation.  
		Denoting $A_j' := m_{j+1}(1-m_{j+1}m_{j}\mathcal{S})^{-1}[A_j]$, it follows from \eqref{eq:underline_def} 
		that
		\begin{equation} \label{eq:Gk_static}
			\begin{split}
				(G-M)_{[1,k+1]} =&~ -\underline{G_{[1,j]}A_{j}'HG_{[j+1,k+1]}}\\ 
				&+ G_{[1,j]}A_{j}'A_{j+1}G_{[j+2,k+1]} - M(z_1, \dots, z_j, A_j'A_{j+1}, z_{j+2}, \dots, z_{k+1})\\
				&+\sum_{i=1}^{j-1}G_{[1,i]} \mathcal{S}\bigl[(G-M)_{[i,j]}A_{j}'\bigr]G_{[j+1,k+1]}\\
				&+ \sum_{i=j+2}^{k+1} G_{[1,j]}A_{j}'\mathcal{S}\bigl[(G-M)_{[j+1,i]}\bigr]G_{[i,k+1]}\\
				&+\sum_{i=1}^{j-1} \sum_q \bigl(G_{[1,i]} S^qG_{[j+1,k+1]} - M_{[1,i],[j+1,k+1]}^{(q)}\bigr) \bigl(M_{[i,j]}A_{j}'\bigr)_{qq}\\
				&+ \sum_{i=j+2}^{k+1} \sum_q \bigl(G_{[1,j]}A_{j}'S^qG_{[i,k+1]} - M(z_1,\dots, z_j, A_{j}'S^q, z_{i}, \dots)\bigr) \bigl(M_{[j+1,i]}\bigr)_{qq}\\
				&+G_{[1,j]} \mathcal{S}\bigl[(G_j-m_j)A_{j}'\bigr]G_{[j+1,k+1]} + G_{[1,j]}A_{j}'\mathcal{S}\bigl[G_{j+1}-m_{j+1}\bigr]G_{[j+1,k+1]},
			\end{split}
		\end{equation}
		where in the special case $j=k$, the term in the second line on the right-hand side of \eqref{eq:Gk_static} is understood as $(G-M)_{[1,k]}A_{k}'$, and the empty summations are zero by convention.
		Here, in order to add and subtract the corresponding $M$-terms, we used the following identity
		\begin{equation} \label{eq:M_re-expansion_identity}
			\begin{split}
				M_{[1,k+1]} =&~ M(z_1, A_1, \dots, z_j, A_j'A_{j+1}, z_{j+2}, \dots, A_{k-1},z_{k+1})\\
				&+ \sum_{i=1}^{j-1}  M\bigl(z_1, \dots, z_i, \mathcal{S}\bigl[M_{[i,j]}A_{j}'\bigr], z_{j+1},\dots, z_{k+1}\bigr),\\
				&+ \sum_{i=j+2}^{k+1}  M\bigl(z_1,\dots, z_j, A_{j}'\mathcal{S}\bigl[M_{[j+1,i]}\bigr], z_{i}, \dots, z_{k+1}\bigr),
			\end{split}
		\end{equation}
		where in the case $j=k$, $M(z_1, A_1, \dots, z_j, A_j'A_{j+1}, z_{j+2}, \dots, A_{k-1},z_{k+1})$ is replaced by $M_{[1,k]}A_{k}'$. We give a detailed proof of the identity \eqref{eq:M_re-expansion_identity} in Section \ref{sec:M_rec_analysis}.
		
		Equation \eqref{eq:Gk_static} is a generalization of \eqref{eq:static1G}, and, unlike \eqref{eq:Gk_static_av}, it provides a matrix-valued expansion for $G_{[1,k+1]}$ (without taking the trace against another observable). In \eqref{eq:Gk_static}, we choose to expand  the resolvent $G_{j+1}$ following the distinguished observable $A_j$ to preserve the number of distinguished observables. Indeed, if a different resolvent were expanded, it would produce a second special observable, not belonging to $\mathcal{V}_1$, on the right-hand side of \eqref{eq:Gk_static}. 
		The complication of~\eqref{eq:Gk_static} compared with~\eqref{eq:Gk_static_av}
		stems from expanding the resolvent in the middle of the chain. Note that in~\eqref{eq:Gk_static_av} the 
		last resolvent was expanded and by cyclicity of the trace any observable could be brought to the last position.
		Such convenience is not available for isotropic chains.
		
		To apply Lemma \ref{lemma:high-moment}, we first estimate the underline-free terms in \eqref{eq:Gk_static}. Using \eqref{eq:iso_static_G-M_ansatz} and $\psi_k^\mathrm{iso}=1$, we obtain, for all $j \in \indset{k-1}$,
		\begin{equation} \label{eq:k-1_term_iso}
			\bigl\lvert \text{ $a,b$ matrix entry of second line of (\ref{eq:Gk_static})} \bigr\rvert \prec \mathfrak{s}_{k+\extk}^\mathrm{iso}(a,\bm x, b).
		\end{equation}
		Similarly, using \eqref{eq:global_size_funcs}, \eqref{eq:global_triag},  \eqref{eq:iso_static_G-M_ansatz}, and \eqref{eq:V_class_bound}, we deduce that the bound \eqref{eq:k-1_term_iso} also holds for $j=k$. 
		
		Next, we observe that all terms in the third, fourth, fifth and sixth line on the right-hand side of \eqref{eq:Gk_static} contain only chains of length up to $k$. Therefore, similarly to \eqref{eq:static_quad_terms}--\eqref{eq:static_sublin_terms}, we obtain
		\begin{equation}
			\bigl\lvert \text{$a,b$ matrix entry of lines 3--6 of (\ref{eq:Gk_static})} \bigr\rvert \prec \mathfrak{s}_{k+\extk}^\mathrm{iso}(a,\bm x, b).
		\end{equation}
		For the last line of~\eqref{eq:Gk_static},
		 using \eqref{eq:iso_static_G-M_ansatz} and \eqref{eq:iso_static_G_bound}, we deduce that 
		\begin{equation}
			\begin{split}
				\bigl\lvert \bigl(G_{[1,j]} \mathcal{S}\bigl[(G_j-m_j)A_{j}'\bigr]G_{[j+1,k+1]}\bigr)_{ab} \bigr\rvert &\prec \sum_{q}  \sqrt{W}\biggl(1 + \frac{\pis{k+1}}{\sqrt{W}}\biggr) \mathfrak{s}_{\extk+1}^\mathrm{iso}(q,\bm y,q)   \mathfrak{s}_{k+1}^\mathrm{iso}(a,\bm x_1,q, \bm x_2,b)\\
				&\prec \biggl(1 + \frac{\pis{k+1}}{\sqrt{W}}\biggr) \mathfrak{s}_{k+\extk}^\mathrm{iso}(a,\bm x, b),
			\end{split}
		\end{equation}
		where in the second step we used  \eqref{eq:global_size_funcs} and \eqref{eq:global_convol}. Finally, using \eqref{eq:apriori_av}, \eqref{eq:iso_static_G-M_ansatz}, \eqref{eq:V_class_bound}, and \eqref{eq:iso_static_G_bound}, we obtain
		\begin{equation}
			\begin{split}
				\bigl\lvert \bigl(G_{[1,j]}A_{j}'\mathcal{S}\bigl[G_{j+1}-m_{j+1}\bigr]G_{[j+1,k+1]}\bigr)_{ab} \bigr\rvert \prec&~ W^{3/4}\sum_{q}  \mathfrak{s}_{j}^\mathrm{iso}(a,\bm x_1,q)  \mathfrak{s}_{\extk+1}^\mathrm{iso}(q,\bm y,q) 
				 \mathfrak{s}_{k-j+1}^\mathrm{iso}(q,\bm x_2,b)  \\
				&\qquad\qquad\times \bigl(\delta_{aq} + W^{-1/2}\bigr)\bigl(\delta_{qb} + W^{-1/2}\bigr)\\
				\prec&~ W^{-1/4}  \mathfrak{s}_{k+\extk}^\mathrm{iso}(a,\bm x, b),
			\end{split}
		\end{equation}
		where we used \eqref{eq:global_size_funcs}, \eqref{eq:global_triag} and \eqref{eq:global_convol}. Therefore, for all $a,b\in\indset{N}$ and $j \in \indset{k}$,
		\begin{equation} \label{eq:Gk_static_est}
			\bigl((G-M)_{[1,k+1]}\bigr)_{ab} =  -\bigl(\underline{G_{[1,j]}A_{j}'HG_{[j+1,k+1]}}\bigr)_{ab} + \mathcal{E}_{ab}, \quad \bigl\lvert \mathcal{E}_{ab} \bigr\rvert \prec \bigl(1+ W^{-1/2}\pis{k+1}\bigr)\mathfrak{s}_{k+\extk}^\mathrm{iso}(a,\bm x, b).
		\end{equation}
		
		The underlined term in \eqref{eq:Gk_static_est} satisfies the following claim.
		\begin{claim} \label{claim:underline_k_iso}
			Let $j \in \indset{k}$, $k \in \indset{\maxK}$, and $\mathcal{Y}_{k}  := \bigl((G-M)_{[1,k+1]}\bigr)_{ab}$. Then, under the assumptions of Proposition \ref{prop:iso_static_masters}, for any integer $p\in\mathbb{N}$ and $0 < \nu < 1/100$, we have 
			\begin{equation} \label{eq:underline_k_iso}
				\begin{split}
					\biggl\lvert \Expv \biggl[ \bigl(\underline{G_{[1,j]}A_{j}'HG_{[j+1,k+1]}}\bigr)_{ab}\cdot \overline{\mathcal{Y}_{k}} \lvert \mathcal{Y}_{k}\rvert^{2p-p} \biggr]  \biggr\rvert \lesssim&~ 
					C_p \biggl( N^{2\nu}\biggl(\sqrt{\pis{k+1}}+\frac{\pis{k+1}}{W^{1/4}} \biggr) \mathfrak{s}_{k+\extk}^\mathrm{iso}(a,\bm x,b)\biggr)^{2p}\\
					& + N^{-\nu}\Expv \bigl[ \lvert \mathcal{Y}_{k}\rvert^{2p} \bigr].
				\end{split}		
			\end{equation}
			for some positive constant $C_p$ depending only on $p$.
		\end{claim}
		We prove Claim \ref{claim:underline_k_iso} in Section \ref{sec:cum_expand}. It follows from \eqref{eq:Gk_static_est},  Claim \ref{claim:underline_k_iso}, and Lemma \ref{lemma:high-moment}, that \eqref{eq:static_iso_masters} holds. This concludes the proof of Proposition \ref{prop:iso_static_masters}.
	\end{proof}

\subsection{Cumulant expansions} \label{sec:cum_expand}
Our key tool in this section is the standard cumulant expansion formula.
We recall that the joint cumulant $\kappa(X_1, X_2, \ldots, X_k)$ 
of several (possibly identical) random variables $X_1, X_2, \ldots  X_p$ is defined via the cumulant generating function, i.e.
by 
$$
\kappa(X_1, X_2, \ldots, X_k) : = (-i)^p \partial_{t_1}\partial_{t_2}\ldots \partial_{t_p} 
\log \Expv e^{i\sum_{j=1}^p t_j X_j}\Big|_{t_1=\ldots =t_p=0}.
$$

\begin{lemma}[Cumulant Expansion] \label{lemma:cumexp} (Section II in \cite{Khorunzhy1999}, Lemma 3.1 in \cite{He2017WignerCLT}) Let $H=H^*$ be a random matrix with entries $H_{ab} := \sqrt{S_{ab}}h_{ab}$. Assume that $\{h_{ab}\}_{a,b=1}^N$ have finite moments of all orders, and let $\kappa_{ab}^{q,q'}$ denote the joint cumulant
 of $q$ copies of $h_{ab}$ and $q'$ copies of $\overline{h_{ab}} = h_{ba}$. Then, for any smooth 
  function $f:\mathbb{C}^{N\times N}\to \mathbb{C}$, and any fixed $L \in \mathbb{N}$, we have 
	\begin{equation}\label{eq:cumulant}
		\Expv\bigl[H_{ab}\cdot\,f(H)\bigr] = \sum_{r=1}^{L} \frac{1}{r!}\bigl(S_{ab}\bigr)^{(r+1)/2} \sum_{q+q'=r}\kappa_{ab}^{q+1,q'} \Expv\bigl[\partial^q_{ab}\partial^{q'}_{ba}f(H) \bigr] + (S_{ab})^{L/2+1}\Omega_{L+1},
	\end{equation}
	where we recall that $\partial_{ab}$ denotes the partial derivative in the direction of the  $(a,b)$ matrix entry, and the error term $\Omega_{L}$ admits the bound
	\begin{equation} \label{eq:cumexp_error}
		\lvert \Omega_L\rvert \lesssim \Expv \bigl[ |h_{ab}|^{L+2} \mathds{1}_{|h_{ab}|\ge N^{\xi}} \bigr] \bigl\lVert f^{(L+1)}\bigr\rVert_\infty + \Expv \bigl[ |h_{ab}|^{L+2}  \bigr] \bigl\lVert \lvert f^{(L+1)}(H) \rvert \mathds{1}_{|h_{ab}|\le N^{\xi}}\bigr\rVert_{\infty}.
	\end{equation}
\end{lemma}

To present the cumulant expansion more compactly, we use the following notation, introduced in 
\cite{Cipolloni2022Optimal}\footnote{We point out a small typo in~\cite{Cipolloni2022Optimal}; from Eq. (4.15) onwards
the cardinality of the union $|J\cup J^*|$ was incorrectly used instead of $|J|+|J^*|$ and similarly 
the definition of $\sum J\cup J^*$ included the elements $\bm j\in J\cap J^*$ only once instead of the correct~\eqref{eq:J_sizes}.  The same typo appeared
in~\cite{Cipolloni2022overlap}  from Eq. (5.11) onwards.}. For a vector $\bm j := (j_1,j_2) \in (\mathbb{Z}_{\ge 0})^2$ and two (possibly empty) subsets $J,J^* \subset (\mathbb{Z}_{\ge 0})^2 \backslash \{\bm 0\}$, where $\bm 0 := (0,0) \in (\mathbb{Z}_{\ge 0})^2$, we define
\begin{equation} \label{eq:J_sizes}
	|\bm j| := j_1 + j_2, \quad \sumJ := \sum_{\bm j\in J}|\bm j| + \sum_{\bm j\in J^*}|\bm j|.
\end{equation}
We use $\bm l, \bm j$ to denote vectors in $(\mathbb{Z}_{\ge 0})^2$. Integer $p$ is always used to denote half of the order of the moment under consideration, and we assume that $J, J^*$ satisfy $|J| \le p-1$ and $|J^*| \le p$, respectively.

We begin by considering the case of averaged chains and proving Claims \ref{claim:1Gav_underline}, \ref{claim:underline_k_av}. 
\subsubsection{Averaged chains} \label{sec:cumexp_av}
To illustrate the overall strategy, we prove the more general Claim \ref{claim:underline_k_av} first.
\begin{proof}[Proof of Claim \ref{claim:underline_k_av}]
	To estimate the left-hand side of \eqref{eq:underline_k_av}, we use the cumulant expansion formula from Lemma \ref{lemma:cumexp} and the bounds $\kappa^{q,q'}_{ab} \le C_{q,q'}$ from~\eqref{eq:all_moments} to have
	\begin{equation} \label{eq:k_av_cumexp}
		\begin{split}
			\bigl\lvert \text{l.h.s. of } (\ref{eq:underline_k_av}) \bigr\rvert
			\lesssim&~ \sum_{\sumJ=1} \Expv \biggl[ Z_{k}^{\mathrm{av}}(\bm 0, J, J^*) \lvert\mathcal{X}_k\rvert^{2p-2}\biggr]\\
			&+ \sum_{r=2}^{L(p)}\sum_{|\bm l| + \sumJ = r } \Expv\biggl[Z_k^{\mathrm{av}}(\bm l, J, J^*) \lvert\mathcal{X}_k\rvert^{2p-|J|-|J^*|-1}\biggr] + C_p\bigl(W^{-k/2}N^{-kD'}\bigr)^{2p},
		\end{split}
	\end{equation}
	where $L(p) \in \mathbb{N}$, $\bm l\in (\mathbb{Z}_{\ge 0})^2$, $J,J^* \subset (\mathbb{Z}_{\ge 0})^2 \backslash \{\bm 0\}$ with $J, J^*$ satisfy $|J| \le p-1$, $|J^*| \le p$, and the quantities $Z_k^\mathrm{av}(\bm l, J, J^*)$ are defined as
	\begin{equation} \label{eq:Zkav_def}
		Z_k^\mathrm{av}(\bm l, J, J^*) := \sum_{cd} \bigl(S_{cd}\bigr)^{(|\bm l| + \sumJ +1)/2} \bigl\lvert \partial_{cd}^{\bm l} \bigl(G_{[1,k]}\bigr)_{dc}\bigl(A_k'\bigr)_{cc} \bigr\rvert \biggl\lvert \prod_{\bm j \in J} \partial_{cd}^{\bm j} \mathcal{X}_k \biggr\rvert \biggl\lvert \prod_{\bm j \in J^*} \partial_{cd}^{\bm j} \overline{\mathcal{X}_k} \biggr\rvert, 
	\end{equation}
	where $\partial_{cd}^{\bm j} := \partial_{cd}^{j_1} \partial_{dc}^{j_2}$.
	 Note that the term with $|\bm l| = 1$ and $J=J^* = \emptyset$ is canceled by definition of the underline in \eqref{eq:underline_def}. Hence, we put the terms involving second-order cumulants, that is, with $\bm l = \bm 0$ and $\sumJ = 1$ into a separate summation. It is straightforward to check, using $\norm{G}\lesssim \eta^{-1}\lesssim 1$ and \eqref{eq:cumexp_error}, that by choosing the expansion order $L(p)$ large enough, the error term $\Omega_{L(p)+1}$ can indeed be bounded by $C_p(W^{-k/2}N^{-kD'})^{2p}$, where $D'$ is the constant from \eqref{eq:ups_lower_bound}.
	
	Clearly, if $Z_k^\mathrm{av}(\bm l, J, J^*) \lesssim \varphi^{1+|J| +|J^*|}$ with very high probability for some $\bm l,   J , J^*$, and some deterministic $\varphi > N^{-D}$, then, by Young's inequality, 
	\begin{equation} \label{eq:enough_to_bound_Z}
		\Expv\bigl[Z_k^{\mathrm{av}}(\bm l, J, J^*) \lvert\mathcal{X}_k\rvert^{2p-|J| -|J^*|-1}\bigr] \lesssim \bigl(N^{\nu}\varphi\bigr)^{2p} + N^{-\nu}\Expv\bigl[\lvert\mathcal{X}_k\rvert^{2p}\bigr].
	\end{equation}
	Note that, by Young's inequality, we achieved a small factor $N^{-\nu}$ in front of 
	$\Expv\bigl[\lvert\mathcal{X}_k\rvert^{2p}\bigr]$, hence,  when~\eqref{eq:enough_to_bound_Z} will be
	used to bound $\Expv\bigl[\lvert\mathcal{X}_k\rvert^{2p}\bigr]$ on the left hand side,  this error
	can be incorporated. We will use this procedure several times and we call it {\it self-consistent high moment bound}.  
	Therefore, it suffices to estimate all $Z_k^\mathrm{av}(\bm l, J, J^*)$. The quantity $Z_{k}^\mathrm{av}(\bm l, J, J^*)$, defined in \eqref{eq:Zkav_def}, contains $|J | + |J^*|$ chains destroyed by derivatives, therefore, we need to recover exactly $(|J| + |J^*|+1)$ decaying factors $\mathfrak{s}_{k+\extk-1}^\mathrm{av}(\bm x)$, with the additional one accounting for the original underlined factor.
	
	First, we consider the terms involving second order cumulants, that is $\bm l =0$ and $\sumJ = 1$ in \eqref{eq:k_av_cumexp}. In particular, we analyze the case  $J=\{(1,0)\}$, $J^* = \emptyset$ in full detail, since the other cases are treated analogously. It follows from \eqref{eq:V_class_bound}, \eqref{eq:static_G_bound} that 
	\begin{equation}\label{eq:Zkav_10}
		\begin{split}
			Z_k^\mathrm{av}(\bm 0, \{(1,0)\}, \emptyset) &\lesssim  \sum_{j=1}^k\sum_{cd} S_{cd}  \bigl\lvert \bigl(G_{[1,k]}\bigr)_{dc}\bigl(A_k'\bigr)_{cc} \bigr\rvert \bigl\lvert \bigl(G_{[j,k]}A_kG_{[1,j]}\bigr)_{dc} \bigr\rvert\\
			&\prec  \sum_{j=1}^k\sum_{cd} S_{cd} \bigl(\delta_{cd}\sqrt{W} + 1\bigr)\,\mathfrak{s}_{k}^{\mathrm{iso}}(c, \bm x^{(\extk)}, d) \sqrt{W}\,\mathfrak{s}^{\mathrm{iso}}_{\extk+1}(c,\bm y, c)\\
			&\qquad\qquad \times\bigl(\delta_{cd}\sqrt{W} + 1\bigr)\,  \mathfrak{s}_{k+1}^\mathrm{iso}(d,x_{j},\dots, x_k,\bm y, x_1, \dots, x_{j-1},c)\\
			&\prec   \mathfrak{s}^\mathrm{av}_{k+\extk-1}(\bm x),
		\end{split}
	\end{equation}
	with very high probability, where $\bm x := (x_1, \dots, x_{k+\extk -1}) \in \indset{N}^{k + \extk-1}$, and $\bm y := (x_k, \dots, x_{k+\extk-1}) \in \indset{N}^{\extk}$. Here, in the last step, we used \eqref{eq:global_size_funcs}, \eqref{eq:global_triag}, and \eqref{eq:global_convol}, together with \eqref{eq:SUps_comvol}. Completely analogously, we show that
	\begin{equation} \label{eq:Zkav_cum2}
		Z_k^\mathrm{av}(\bm 0, J, J^*) \prec  \mathfrak{s}^\mathrm{av}_{k+\extk-1}(\bm x),  \quad \sumJ = 1.
	\end{equation}
	
	Next, we estimate the terms containing higher-order cumulants in \eqref{eq:k_av_cumexp}, that is $|\bm l| + \sumJ \ge 2$. 	
	To this end, we first estimate $(S_{cd})^{|\bm j|/2}\partial_{cd}^{\bm j} \mathcal{X}_k$. For all $\bm j  \in (\mathbb{Z}_{\ge 0})^2 \backslash\bm 0$ with $|\bm j| =: u \ge 1$,
	\begin{equation} \label{eq:partial_X_1}
		\bigl\lvert \partial_{cd}^{\bm j} \mathcal{X}_k \bigr\rvert \lesssim \sum_{1\le a_1\le\dots\le a_u\le k} \sum_{\substack{(c_i,d_i)\in \{(c,d),(d,c)\}\\ i \in \indset{u}}} \biggl(\prod_{i=1}^{u-1} \bigl\lvert (G_{[a_i, a_{i+1}]})_{ c_i d_{i+1}} \bigr\rvert\biggr) \bigl\lvert (G_{[a_u, k]}A_kG_{[1, a_{1}]})_{ c_u d_{1}} \bigr\rvert,
	\end{equation}
	where the first summation runs over integers $\{a_i\}_{i=1}^u$ satisfying $1\le a_1\le\dots\le a_u\le k$. For brevity, we denote 
	\begin{equation}
		\mathcal{G}_i := G_{[a_i, a_{i+1}]}, \quad i \in \indset{u-1}, \quad \mathcal{G}_u := G_{[a_u, k]}A_kG_{[1, a_{1}]},
	\end{equation}
	and their respective lengths by $ k_i:=a_{i+1}-a_i+1$ for $i\le u-1$ and $k_u=k+a_1-a_u+1$.
	 Clearly, $\sum_{i=1}^u k_i = k + u$. 
		
	Chains $\mathcal{G}_{i}$ with $k_i \ge 2$ in \eqref{eq:partial_X_1} carry the dependence on the original indices $\bm x$. We estimate every chain $\mathcal{G}_i$ in \eqref{eq:partial_X_1} by \eqref{eq:static_G-M_ansatz}, which results in $\sqrt{\Upsilon}$-decay factors coupling the original indices $\bm x$ to (potentially different) summation indices $c$ and $d$. 
	More precisely, for $u:=|\bm j|\ge 1$,
	\begin{equation} \label{eq:partial_X2}
		\bigl\lvert \partial_{cd}^{\bm j} \mathcal{X}_k \bigr\rvert \prec \sum_{1\le a_1\le\dots\le a_u\le k} \sum_{\{(c_i,d_i)\}} \prod_{i=1}^{u} \bigl(\sqrt{W}\delta_{c_id_{i+1}} + 1\bigr)\, \mathfrak{s}_i(c_i, d_{i+1}),
	\end{equation}
	where $\mathfrak{s}_i(a,b):=\mathfrak{s}_{k_i}^\mathrm{iso}(a,\bm x_i, b)$ for $i \in \indset{u-1}$,  $\mathfrak{s}_u(a,b) := \mathfrak{s}_{k_u+\extk-1}^\mathrm{iso}(a,\bm x_u, b)$, and  we denote the (potentially empty) vector of external indices on which $\mathcal{G}_i$ depends by $\bm x_i$. Here we adhere to the cyclic convention $d_{u+1} := d_1$, and we use the sum over $\{(c_i,d_i)\}$ to abbreviate the summations over $(c_i,d_i) \in \{(c,d), (d,c)\}$ for $i\in\indset{u}$ as in \eqref{eq:partial_X_1}. The bound \eqref{eq:partial_X2} also holds with $\mathcal{X}_k$ replaced by its complex conjugate. 
	
	Moreover, the $\bm l$-derivative in~\eqref{eq:Zkav_def} satisfies an analogous bound
	\begin{equation} \label{eq:partial_l_X}
		\bigl\lvert \partial_{cd}^{\bm l} \bigl(G_{[1,k]}\bigr)_{dc} \bigr\rvert \prec \sum_{1\le a_1\le\dots\le a_{u-1}\le k} \sum_{\{(c_i,d_i)\}} \prod_{i=1}^{u} \bigl(\sqrt{W}\delta_{c_{i}d_{i+1}} + 1\bigr)\, \other{\mathfrak{s}}_i(c_{i},d_{i+1}) ,
	\end{equation}
	where $u := |\bm l|+1 \ge 1$, $d_{u+1} := d_1$, and the decay factors $\other{\mathfrak{s}}_i$ correspond to chains $G_{[a_{i-1}, a_{i}]}$ with $a_0 := 1$ and $a_u := k$.  By a simple parity argument, for any arrangement of $(c_i,d_i)$,  the product of the $\delta$ terms in \eqref{eq:partial_X2}--\eqref{eq:partial_l_X} can be overestimated by 
	\begin{equation} \label{eq:delta_est}
		\prod_{i=1}^{u} \bigl(\sqrt{W}\delta_{c_{i}d_{i+1}} + 1\bigr) \le W^{u/2} \bigl(\delta_{cd} + 
		\mathds{1}_{u\text{--odd}} + W^{-1/2} \bigr). 
	\end{equation}

	 Our goal is to combine the  decay factors to recover the target spatial decay $\mathfrak{s}^\mathrm{av}_{k+\extk-1}(\bm x)$ in~\eqref{eq:underline_k_av}. 
	To this end, we bridge the distance between the indices $c$ and $d$ with the aid of 
	$\sqrt{S_{cd}}$-powers available in~\eqref{eq:Zkav_def}. In practice, we will replace all
	$d$ indices by $c$'s\footnote{
		Since the summation over $c$ in \eqref{eq:Zkav_def} involves an additional $\lvert (A_k')_{cc}\rvert$, in the sequel we opt for replacing $d$'s by $c$'s, but otherwise the roles of $c$ and $d$ are interchangeable.
	},  since the $\sqrt{S_{cd}}$ factors force them to be close.
	  We use two mechanisms to accomplish this: the triangle inequality \eqref{eq:global_triag}, and the convolution inequality \eqref{eq:global_convol}. 
	
	Indeed, since $S_{cd} \lesssim \Upsilon_{cd}$ by \eqref{eq:SUps_comvol}, it follows from \eqref{eq:triag} 
	(with $\ell\sim W$, $\eta\sim 1$) that
	\begin{equation} \label{eq:triag_replacement}
		\sqrt{S_{cd}}\,\mathfrak{s}_{i}(d, \,\cdot\,) \lesssim W^{-1/2}\mathfrak{s}_{i}(c, \,\cdot\,), \quad \sqrt{S_{cd}}\,\mathfrak{s}_{i}(\,\cdot\,, d) \lesssim W^{-1/2}\mathfrak{s}_{i}(\,\cdot\,, c),   \quad i \in \indset{u}.
	\end{equation}
	The estimate \eqref{eq:triag_replacement} can be used  as many times as we have $\sqrt{S_{cd}}$-factors available.
	On the other hand, \eqref{eq:global_convol} implies that 
	\begin{equation} \label{eq:convol_replacement}
		\sum_{d} S_{cd} \,\mathfrak{s}_i(d,c)\mathfrak{s}_j(c,d) \lesssim \mathfrak{s}_i(c,c)\mathfrak{s}_j(c,c), \quad \sum_d S_{cd}\, \mathfrak{s}_i(d,d) \lesssim \mathfrak{s}_i(c,c), \quad i,j\in\indset{u}.
	\end{equation}
	The estimate \eqref{eq:convol_replacement} also holds for any permutation of the arguments $(c,d)$ of $\mathfrak{s}$, however, it can only be used once (when $d$ is summed up in~\eqref{eq:Zkav_def}) and requires a full power of $S_{cd}$. Both \eqref{eq:triag_replacement} and \eqref{eq:convol_replacement} also hold for  $\other{\mathfrak{s}}_i$.
	
	Observe that by  \eqref{eq:global_size_funcs}  and \eqref{eq:global_triag}, we have, for any fixed $c\in\indset{N}$, 
	\begin{equation} \label{eq:s_triag_concat}
		\begin{split}
			\mathfrak{s}_{n_1}^\mathrm{iso}(c,\bm y_1, c)\,\mathfrak{s}_{n_2}^\mathrm{iso}(c,\bm y_{2}, \,\cdot\,) &\lesssim W^{-1/2}  \mathfrak{s}_{n_1+n_2-1}^\mathrm{iso}(c,\bm y_1, \bm y_2,  \,\cdot\,), \quad n_1,n_2\in \mathbb{N},\quad  \bm y_{i} \in \indset{N}^{n_i-1}\\
			\mathfrak{s}_{n_1}^{\mathrm{iso}}(c,\bm y_1, c) &\lesssim W^{-1/2}\mathfrak{s}_{n_1-1}^{\mathrm{av}}(\bm y_1), \quad n_1\in\mathbb{N}, \quad \bm y_1 \in \indset{N}^{n_1-1}.
		\end{split}
	\end{equation} 
	Hence, once all arguments $d$ are replaced by $c$, using \eqref{eq:s_triag_concat} repeatedly, we conclude that 
	\begin{equation} \label{eq:glob_concatenations}
		\prod_{i=1}^{u} \mathfrak{s}_i(c,c) \lesssim W^{-u/2} \mathfrak{s}_{k+n-1}^\mathrm{av}(\bm x), \quad \prod_{i=1}^{u} \other{\mathfrak{s}}_i(c,c) \lesssim W^{-(u-1)/2} \mathfrak{s}_k^\mathrm{iso}(c, \bm x^{(n)},c)
	\end{equation}

	Observe that both replacement mechanisms  \eqref{eq:triag_replacement}--\eqref{eq:convol_replacement} require a single $\sqrt{S_{cd}}$ for each index $d$ replaced by $c$. If $c\neq d$, it is straightforward to check that, across all $\bm j$- and $\bm l$-derivatives on the right-hand side of \eqref{eq:Zkav_def}, there are exactly $r := |\bm l| + \sumJ + 1$ instances of the index $d$ (each $\partial_{cd}$ or $\partial_{dc}$ creates one $d$, and the remnant of the underlined chain contains an additional $d$). Therefore, for $d\neq c$, using \eqref{eq:triag_replacement} $(r-1)$ times and \eqref{eq:convol_replacement} once, we replace all $d$'s by $c$'s, with the additional gain of $W^{-(r-1)/2}$ from \eqref{eq:triag_replacement}, obtaining
	\begin{equation} \label{eq:partial_X_plug_in}
		\begin{split}
			Z_{k}^{\mathrm{av}}(\bm l, J, J^*) \prec&~  C_r 
			\sum_{c} W^{-r/2+1} \mathfrak{s}_k^\mathrm{iso}(c, \bm x^{(n)},c) \bigl\lvert (A_k')_{cc} \bigr\rvert \bigl(\mathfrak{s}_{k+n-1}^\mathrm{av}(\bm x)\bigr)^{|J|+|J^*|}\\
			&\quad \times \bigl(\mathds{1}_{|\bm l| \text{--odd} } + W^{-1/2}\bigr) \prod_{\bm j \in J,J^*}  \bigl(\mathds{1}_{|\bm j| \text{--even} } + W^{-1/2}\bigr)\\
			&+C_r \sum_{c} W^{-r/2+1/2} \mathfrak{s}_k^\mathrm{iso}(c, \bm x^{(n)},c) \bigl\lvert (A_k')_{cc} \bigr\rvert \, \bigl(\mathfrak{s}_{k+n-1}^\mathrm{av}(\bm x)\bigr)^{|J|+|J^*|},
		\end{split}
	\end{equation}
	where the second sum over $c$ comes from the regime $d=c$, and is estimated using \eqref{eq:partial_X2}--\eqref{eq:delta_est}, \eqref{eq:glob_concatenations}. Here the constant $C_r > 0$ depends only on $r$, and we omit it in the sequel.
	
	Estimating $A_k'$ using \eqref{eq:Ak'_bound} and using \eqref{eq:global_convol}, we deduce from \eqref{eq:partial_X_plug_in}, that
	\begin{equation} \label{eq:partial_X_convol1}
		\frac{Z_{k}^{\mathrm{av}}(\bm l, J, J^*)}{\bigl(\mathfrak{s}_{k+n-1}^\mathrm{av}(\bm x)\bigr)^{|J|+|J^*|+1}} \prec W^{-r/2+3/2} \bigl(\mathds{1}_{|\bm l| \text{--odd} } + W^{-1/2}\bigr) \prod_{\bm j \in J,J^*}  \bigl(\mathds{1}_{|\bm j| \text{--even} } + W^{-1/2}\bigr) + W^{-r/2+1},
	\end{equation}
	where we recall that $r := |\bm l| + \sum J +1$.
	Note that for $|J|+|J^*| \ge 1$ and  $|\bm l|+\sumJ \ge 3$, estimate \eqref{eq:partial_X_convol1} immediately implies the desired
	\begin{equation} \label{eq:Zkav_high}
		Z_{k}^{\mathrm{av}}(\bm l, J, J^*) \prec  \bigl( \mathfrak{s}_{k+\extk-1}^\mathrm{av}(\bm x)\bigr)^{|J| + |J^*|+1}.
	\end{equation}
	Therefore, it suffices to consider $|\bm l|+\sumJ = 2$. We show that for any combination of $|\bm l|, J, J^*$, at least one of $|\bm j|$ is odd or $|\bm l|$ is even. As a result, the right-hand side of \eqref{eq:partial_X_convol1} will be bounded by a constant:
	\begin{itemize}
		\item If $|\bm l| = 2$, $\sumJ = 0$ or $|\bm l| = 0$, $\sumJ = 2$, then $|\bm l| + 1$ is odd, hence
		\begin{equation} \label{eq:Zkav_3_1} 
				Z_{k}^{\mathrm{av}}(\bm l, J, J^*) \prec   \bigl(\mathfrak{s}_{k+\extk-1}^\mathrm{av}(\bm x)\bigr)^{|J|+|J^*|+1}.
		\end{equation}
		\item If $|\bm l| = 1$ and $\sumJ = 1$, then $|\bm j| = 1$ is odd for the only $\bm j \in J, J^*$, hence
		\begin{equation} \label{eq:Zkav_3_2} 
				Z_{k}^{\mathrm{av}}(\bm l, J, J^*) \prec    \bigl(\mathfrak{s}_{k+\extk-1}^\mathrm{av}(\bm x)\bigr)^{2} = \bigl(\mathfrak{s}_{k+\extk-1}^\mathrm{av}(\bm x)\bigr)^{|J| +|J^*|+1}. 
		\end{equation}
	\end{itemize}
	 Therefore, we conclude that \eqref{eq:Zkav_high} holds for all $|\bm l|+\sumJ \ge 2$. 
	Combining \eqref{eq:enough_to_bound_Z},  \eqref{eq:Zkav_cum2} and \eqref{eq:Zkav_high}  for all $|\bm l|+\sumJ \ge 2$, we obtain desired \eqref{eq:underline_k_av}. This concludes the proof of Claim \ref{claim:underline_k_av}.
\end{proof}

Next, we prove Claim \ref{claim:1Gav_underline}. Our proof follows the same general approach as that of Claim \ref{claim:underline_k_av} above, with one key difference: since no spatial decay of $|(G-m)_{ab}|$ is known at this stage, we rely on the Ward identity
\begin{equation} \label{eq:Ward}
	GG^* = G^*G = \eta^{-1} \,\im G,
\end{equation}
and the trivial bound for the resolvent $\norm{G(z)}\lesssim \eta^{-1}$, that, provided $\eta \gtrsim 1$, implies
\begin{equation} \label{eq:G_trivial}
	\norm{G}_{\max} \lesssim 1.
\end{equation}

\begin{proof} [Proof of Claim \ref{claim:1Gav_underline}]
	Recall that $A_1 := S^x$, $A_1' := m\,\mathcal{B}^{-1}[A_1]$ 
	and the cumulant expansion \eqref{eq:k_av_cumexp} for $k=1$, namely
	\begin{equation} \label{eq:1G_av_cuexp}
		\begin{split}
			\biggl\lvert\Expv\biggl[\Tr\bigl[ \underline{HG} A_1'\bigr] \cdot\overline{\mathcal{X}_1} \lvert \mathcal{X}_1 \rvert^{2p-2}\biggr]\biggr\rvert 
			\lesssim&~ \sum_{\sumJ=1} \Expv \biggl[ Z_{1}^{\mathrm{av}}(\bm 0, J, J^*) \lvert\mathcal{X}_1\rvert^{2p-2}\biggr] + W^{-2p}\\
			&+ \sum_{r=2}^{L(p)}\sum_{|\bm l| + \sumJ = r } \Expv\biggl[Z_1^{\mathrm{av}}(\bm l, J, J^*) \lvert\mathcal{X}_1\rvert^{2p-|J| -|J^*|-1}\biggr],
		\end{split}
	\end{equation}
	where the quantities $Z_{1}^\mathrm{av}$ are defined in \eqref{eq:Zkav_def}.
	It follows from \eqref{eq:G_trivial} and \eqref{eq:S_bound}, that for all $a,b\in\indset{N}$,
	\begin{equation}
		\bigl\lvert \bigl(GS^xG\bigr)_{ab} \bigr\rvert \lesssim W^{-1}, \quad \sum_d \bigl\lvert\bigl(GS^xG\bigr)_{cd}\bigr\rvert^2 =  \bigl(GS^xGG^*S^xG^*\bigr)_{cc} \lesssim W^{-2}.
	\end{equation}
	Therefore, similarly to \eqref{eq:partial_X_1}, we obtain
	\begin{equation} \label{eq:partial_X1}
		\bigl\lvert \partial^{\bm j} \mathcal{X}_1\bigr\rvert \lesssim \bigl\lvert\bigl(GS^xG\bigr)_{cd}\bigr\rvert + \mathds{1}_{|\bm j| \ge 2} W^{-1}, \quad \bigl\lvert \partial^{\bm l} G_{dc}\bigr\rvert \lesssim \bigl\lvert G_{dc} \bigr\rvert + \mathds{1}_{|\bm l| \ge 1}.
	\end{equation}
	Moreover, it follows from \eqref{eq:S_bound} and \eqref{eq:stab_bound} that, for $A_1 := S^x$,
	\begin{equation} \label{eq:other_S_bounds}
		\lVert A_1'\rVert_{\max} \lesssim W^{-1}, \quad \sum_{c}\bigl\lvert (A_1')_{cc}\bigr\rvert \lesssim 1.
	\end{equation}
	Hence, similarly to \eqref{eq:Zkav_high}--\eqref{eq:Zkav_3_2}, by simple power counting, it follows from \eqref{eq:Ward}, the bounds \eqref{eq:partial_X1} and \eqref{eq:other_S_bounds} that, for $|\bm l| + \sumJ \ge 2$,
	\begin{equation}
		Z_{1}^\mathrm{av}(\bm l, J, J^*)  \lesssim W^{-(1+ |J| + |J^*|) }~.
	\end{equation}
	Therefore, it remains to estimate the terms with $\bm l = \bm 0$ and $\sumJ = 1$.
	
	Combining Schwarz inequality, Ward identity \eqref{eq:Ward}, \eqref{eq:sumS=1}, and \eqref{eq:G_trivial},  we obtain
	\begin{equation} \label{eq:Schwarz+Ward}
		\sum_{c} \bigl(S_{ac}\bigr)^{1/2}  |G_{cb}|  \le \biggl(\sum_{c} S_{ac}\biggr)^{1/2}\biggl(\sum_{c}  |G_{cb}|^2 \biggr)^{1/2} \lesssim \frac{1}{\eta}\lesssim 1.
	\end{equation}
	Hence, for $\bm l = \bm 0$, $J:=\{(1,0)\}$ and $J^* = \emptyset$, by \eqref{eq:S_bound}, \eqref{eq:other_S_bounds}, \eqref{eq:Schwarz+Ward},  we have
	\begin{equation}
		Z_{1}^\mathrm{av}(\bm 0, \{(1,0)\}, \emptyset) \lesssim \biggl\lvert\sum_{cd} \bigl(A_1'\bigr)_{cc} S_{cd} \bigl(GS^xG \bigr)_{dc}  \biggr\rvert \lesssim W^{-3/2},
	\end{equation}
	and the other terms  with $\sumJ = 1$ admit the same bound. 
	
	Therefore, using \eqref{eq:enough_to_bound_Z}, we obtain the desired \eqref{eq:apriori_av_underline}. 
	Note that the suboptimal $W^{-3/4}$ factor (instead of $W^{-1}$) in~\eqref{eq:apriori_av} comes solely
	from the second cumulant term $\sum J\cup J^*=1$ owing to the missing decay bound on a single resolvent $G$.
	This concludes the proof of Claim \ref{claim:1Gav_underline}.
\end{proof}

Next, we turn to the case of isotropic chains.
\subsubsection{Isotropic chains}
Following the structure of Section \ref{sec:cumexp_av} above, we first prove the most general Claim \ref{claim:underline_k_iso}. We then use the same method, with appropriate modifications, to establish the remaining Claims \ref{claim:flat1G_underline} and \ref{claim:omega_underline1}--\ref{claim:omega_underline2}.
\begin{proof} [ Proof of Claim \ref{claim:underline_k_iso}] 
	Throughout the proof, we consider the indices $a,b\in\indset{N}$ and $j \in \indset{k}$ to be fixed and omit the dependence of various quantities on them.
	As in \eqref{eq:k_av_cumexp}, to estimate the left-hand side of \eqref{eq:underline_k_iso}, we use the cumulant expansion formula from Lemma \ref{lemma:cumexp},
	\begin{equation} \label{eq:k_iso_cumexp}
		\begin{split}
			\bigl\lvert \text{l.h.s. of } (\ref{eq:underline_k_iso}) \bigr\rvert
			\lesssim&~ \sum_{\sumJ=1} \Expv \biggl[ Z_{k}^{\mathrm{iso}}(\bm 0, J, J^*) \lvert\mathcal{Y}_k\rvert^{2p-2}\biggr] + \bigl(W^{-k/2}N^{-(k+1)D'}\bigr)^{2p}\\
			&+ \sum_{r=2}^{L(p)}\sum_{|\bm l| + \sumJ = r } \Expv\biggl[Z_k^{\mathrm{iso}}(\bm l, J, J^*) \lvert\mathcal{Y}_k\rvert^{2p-|J| -|J^*|-1}\biggr],
		\end{split}
	\end{equation}
	where $L(p) \in \mathbb{N}$, $\bm l\in (\mathbb{Z}_{\ge 0})^2$, $J,J^* \subset (\mathbb{Z}_{\ge 0})^2 \backslash \{\bm 0\}$ with $J, J^*$ satisfy $|J| \le p-1$, $|J^*| \le p$, and the quantities $Z_k^\mathrm{iso}(\bm l, J, J^*)$ are defined as, similarly to \eqref{eq:Zkav_def},
	\begin{equation} \label{eq:Zkiso_def}
		Z_k^\mathrm{iso}(\bm l, J, J^*) := \sum_{cd} \bigl(S_{cd}\bigr)^{(r+1)/2} \bigl\lvert \partial_{cd}^{\bm l} \bigl\{\bigl(G_{[1,j]}A_j'\bigr)_{ac}\bigl(G_{[j+1,k+1]}\bigr)_{db}\bigr\} \bigr\rvert \biggl\lvert \prod_{\bm j \in J} \partial_{cd}^{\bm j} \mathcal{Y}_k \biggr\rvert \biggl\lvert \prod_{\bm j \in J^*} \partial_{cd}^{\bm j} \overline{\mathcal{Y}_k} \biggr\rvert,
	\end{equation}
	where $r := |\bm l| +\sumJ$, and we recall $\partial_{cd}^{\bm j} := \partial_{cd}^{j_1} \partial_{dc}^{j_2}$. 
	
	First, we bound the terms involving second order cumulants, that is $\bm l =0$ and $\sumJ = 1$. It follows from \eqref{eq:Ak'_bound}, \eqref{eq:iso_static_G_bound} that 
	\begin{equation} \label{eq:Zkiso_2cum1}
		\begin{split}
			Z_k^\mathrm{iso}(\bm 0, \{(1,0)\}, \emptyset) \lesssim&~ \sum_{q=1}^k \sum_{cd}  S_{cd}  \bigl\lvert   \bigl(G_{[1,j]}\bigr)_{ac} \bigl(A_j'\bigr)_{cc}\bigl(G_{[j+1,k+1]}\bigr)_{db} \bigr\rvert \bigl\lvert \bigl(G_{[1,q]}\bigr)_{ac} \bigl(G_{[q,k+1]}\bigr)_{db}  \bigr\rvert \\
			\lesssim&~ W^{5/2}\sum_{q=1}^{k+1} \sum_{cd} S_{cd} \bigl(\delta_{ac} + W^{-1/2} \pis{j} \bigr) \mathfrak{s}_{j}^\mathrm{iso}(a, \bm x_1, c) \mathfrak{s}^\mathrm{iso}_{\extk+1}(c,\bm y, c)\\
			&\qquad\qquad\times  \bigl(\delta_{db} +  W^{-1/2}\pis{k-j+1}   \bigr)\mathfrak{s}^\mathrm{iso}_{k-j+1}(d, \bm x_2, b)\\
			&\qquad\qquad\times\biggl(\delta_{ac} + \frac{\pis{q} }{\sqrt{W}}\biggr) \mathfrak{s}^\mathrm{iso}_{k_1(q)}(a, \other{\bm x}_1, c)\biggl(\delta_{db} + \frac{\pis{k+2-q} }{\sqrt{W}}\biggr) \mathfrak{s}^\mathrm{iso}_{k_2(q)}(d, \other{\bm x}_2, b),
		\end{split}
	\end{equation}
	where $\bm x_1$ and $\bm x_2$ are the vectors of external indices on which the resolvent chains $G_{[1,j]}$ and $G_{[j+1,k+1]}$, respectively, depend. Here $q\in\indset{k+1}$ denotes the index   of the resolvent $G_q$ in $G_{[1,k+1]}$ that was differentiated by $\partial_{cd}$, while $\other{\bm x}_1 \in \indset{N}^{k_1(q)-1}$ and $\other{\bm x}_2\in \indset{N}^{k_2(q)-1}$ denote the vectors of external indices on which $G_{[1,q]}$ and $G_{[q,k+1]}$, respectively depend. The index vectors  $\bm x_1$, $\bm x_2$, $\other{\bm x}_1$, $\other{\bm x}_2$ satisfy $(\bm x_1, \bm y, \bm x_2) = \bm x = (\other{\bm x}_1, \other{\bm x}_2)$. In particular, $k_1(q) + k_2(q) = k + \extk +1$.
	
	Using \eqref{eq:global_size_funcs}, \eqref{eq:global_triag}--\eqref{eq:global_convol}, and recalling that $\pis{k'} = 1$ for all $k'\in\indset{k}$, it is straightforward to check that \eqref{eq:Zkiso_2cum1} implies
	\begin{equation} \label{eq:Zkiso_2cum1final} 
		Z_k^\mathrm{iso}(\bm 0, \{(1,0)\}, \emptyset) \prec   \pis{k+1} \bigl(\mathfrak{s}^\mathrm{iso}_{k+n}(a, \bm x, b)\bigr)^2.
	\end{equation}
	All other terms with $\bm l =0$ and $\sumJ = 1$ are estimated analogously. 
	
	Next, we estimate the contribution from the higher-order cumulants. In preparation, we bound $\partial_{cd}^{\bm j} \mathcal{Y}_k$  and $\partial_{cd}^{\bm l} \{(G_{[1,j]}  )_{ac}(G_{[j+1,k+1]})_{db} \}$. 
	 Reasoning as in \eqref{eq:partial_X_1}--\eqref{eq:partial_l_X}, and using the last inequality in~\eqref{eq:iso_static_assume},
	we deduce that
	\begin{equation} \label{eq:partial_Y}
		\begin{split}
			\bigl\lvert \partial_{cd}^{\bm j} \mathcal{Y}_k \bigr\rvert \prec&~ \sum_{1\le a_1\le\dots\le a_u\le k+1} \sum_{\{(c_i,d_i)\}} \mathfrak{s}_0(a,d_1) \, \mathfrak{s}_u(c_u, b)\prod_{i=1}^{u-1} \mathfrak{s}_i(c_i, d_{i+1}) \\
			&\times W^{u/2} \bigl(W^{-1/2} \pis{k+1}  + \other{\delta}_1 \pis{k+1} + \bigl(\other{\delta}_2 + \other{\delta}_1 \mathds{1}_{|\bm j| \ge 2}\bigr)\sqrt{W} \bigr),
		\end{split}
	\end{equation}
	where $u := |\bm j|$, the decay factors $\mathfrak{s}_i$ correspond to chains $G_{[a_{i}, a_{i+1}]}$ with $a_0 := 1$ and $a_{u+1} := k$, and $\other{\delta}_1 \equiv\other{\delta}_1(c,d,a,b)$, $\other{\delta}_2\equiv\other{\delta}_2(c,d,a,b)$ are defined as
	\begin{equation} \label{eq:other_deltas}
		\other{\delta}_1 := \delta_{ca}+\delta_{cb} + \delta_{da} +\delta_{db}, \quad \other{\delta}_2 := (\delta_{ca}+\delta_{cb}) (\delta_{da} +\delta_{db}).
	\end{equation}
	The terms involving $\other{\delta}_i$ arise when diagonal  chains emerge after differentiation. 
	This occurs when one or both of the indices $c,d$ coincide with $a,b$, resulting in a restricted summation, and the subscript $i$ of $\other{\delta}_i$ denotes the number of restricted summations.
	In \eqref{eq:partial_Y}, we overestimated all other $\delta_{c_id_{i+1}}$~by~$1$. The bound \eqref{eq:partial_Y} 
	also holds with $\mathcal{Y}_k$ replaced by its complex conjugate.
	
	Similarly to \eqref{eq:partial_l_X}, we also have
	\begin{equation} \label{eq:partial_l_Y}
		\begin{split}
			\bigl\lvert \partial_{cd}^{\bm l} \bigl\{(G_{[1,j]})_{ac}(G_{[j+1,k+1]})_{db} \bigr\} \bigr\rvert \prec&~ \sum_{\{a_i\}} \sum_{\{(c_i,d_i)\}} \other{\mathfrak{s}}_0(a,d_1) \, \other{\mathfrak{s}}_u(c_u, b)\prod_{i=1}^{u-1} \other{\mathfrak{s}}_i(c_i, d_{i+1}) \\
			&\times W^{u/2} \bigl(W^{-1/2}  + \other{\delta}_1 + \bigl(\other{\delta}_2 + \other{\delta}_1 \mathds{1}_{|\bm l| \ge 1}\bigr)\sqrt{W} \bigr),
		\end{split}
	\end{equation}
	where $u := |\bm l| + 1$. Here $\{a_i\}_{i=1}^u \ni j+1$ is a monotone non-decreasing sequence of integers in $\indset{k+1}$ containing $j+1$, and $\other{\mathfrak{s}}_i$ are the decay factors corresponding to the chains  $\mathcal{G}_i$, given by
	\begin{equation}
		\mathcal{G}_i := \begin{cases}
			G_{[a_i, a_{i+1}]}, \quad &i \neq q,\\
			G_{[a_i, j]}, \quad &i=q,
		\end{cases}
	\end{equation}
	where $q$ is the smallest index such that $a_{q+1} = j+1$. Note that the estimate \eqref{eq:partial_l_Y} contains no $\pis{k+1}$ control parameter, since all resolvent chains involved in the $\partial^{\bm l}_{cd}$-derivative have length at most $k$ and~\eqref{eq:iso_static_assume} 
	was assumed in Claim \ref{claim:underline_k_iso}.
	
	Comparing \eqref{eq:triag_replacement} to \eqref{eq:convol_replacement}, we observe that each restricted summation in \eqref{eq:Zkiso_def} improves the bound by a factor of $W^{-1}$. Hence, analogously to \eqref{eq:Zkav_high}, using \eqref{eq:triag_replacement} and \eqref{eq:convol_replacement} to replace all indices $d$ by $c$, using \eqref{eq:s_triag_concat}, and analyzing the cases $\other{\delta}_1 = 0$, $\other{\delta}_1 \neq 0$ with $\other{\delta}_2 = 0$, and $\other{\delta}_2 \neq 0$ separately, we conclude that, 	 
	\begin{equation} \label{eq:Zkiso_high_final}
		Z_k^\mathrm{iso}(\bm l, J, J^*) \prec  \biggl( \sqrt{\pis{k+1}} +  \frac{\pis{k+1}}{W^{1/4}}\biggr)^{|J| + |J^*| + 1} \times \bigl(\mathfrak{s}_{k+\extk}^\mathrm{iso}(a,\bm x,b)\bigr)^{|J| + |J^*| + 1}.
	\end{equation}
	Here we used \eqref{eq:global_size_funcs}  
	and \eqref{eq:V_class_bound}. Hence, the desired \eqref{eq:underline_k_iso} follows from \eqref{eq:k_iso_cumexp}, \eqref{eq:Zkiso_2cum1final}, \eqref{eq:Zkiso_high_final}, by the same reasoning as in \eqref{eq:enough_to_bound_Z}. This concludes the proof of Claim \ref{claim:underline_k_iso}.
\end{proof}

Next, we prove Claim \ref{claim:flat1G_underline}, which is used in establishing the a priori flat bound \eqref{eq:apriori_flat}. The proof follows the same outline as the proof of Claim \ref{claim:underline_k_iso} above, with one key difference:  the a priori bound~\eqref{eq:iso_static_assume}  is not available but we do not aim to preserve the $|a-b|$ decay of $(G-m)_{ab}$. Instead, we resort to simpler estimates using \eqref{eq:S_bound}, the Ward identity \eqref{eq:Ward}, and \eqref{eq:G_trivial}.
\begin{proof} [Proof of Claim \ref{claim:flat1G_underline}] 
	Recall from \eqref{eq:B_stab_def} that action of the operator $\mathcal{B}^{-1}$ on $R\in \mathbb{C}^{N\times N}$ is given by \eqref{eq:Binv_action}.
	Therefore, we consider two cases: $a\neq b$ and $a = b$.
	
	First, assume that $a\neq b$, then $\big(\mathcal{B}^{-1}[R]\big)_{ab} = R_{ab}$, then, similarly to \eqref{eq:k_iso_cumexp}, we obtain
	\begin{equation} \label{eq:1Gflat_offdiag_cuexp}
		\begin{split}
			\biggl\lvert\Expv\biggl[\big(\mathcal{B}^{-1}[\underline{HG}]\big)_{ab}\overline{\mathcal{Y}_0} \lvert \mathcal{Y}_0 \rvert^{2p-2}\biggr]\biggr\rvert 
			\lesssim&~ \sum_{\sumJ = 1} \Expv\biggl[ Z_{ab}^\mathrm{iso}(\bm 0, J, J^*) \lvert \mathcal{Y}_0\rvert^{2p-2}\biggr]  
			+ N^{-pD'}\\
			&+ \sum_{2\le |\bm l| + \sumJ\le L } \Expv\biggl[ Z_{ab}^\mathrm{iso}(\bm l, J, J^*) \lvert \mathcal{Y}_0\rvert^{2p-1 - |J| -|J^*|}\biggr], \\
		\end{split}
	\end{equation}
	where $D'$ is the constant from \eqref{eq:ups_lower_bound}, and, similarly to \eqref{eq:Zkiso_def}, $Z_{ab}^\mathrm{iso}(\bm l, J, J^*)$ are defined as 
	\begin{equation} \label{eq:Ziso}
		Z_{ab}^\mathrm{iso}(\bm l, J, J^*) := \biggl\lvert\sum_{c} \bigl(S_{ac}\bigr)^{(1 + |\bm l| + \sumJ)/2} \bigl(\partial_{ac}^{\bm l} G_{cb} \bigr) \biggl(\prod_{\bm j \in J} \partial_{ac}^{\bm j} \mathcal{Y}_0\biggr) \biggl(\prod_{\bm j \in J^*} \partial_{ac}^{\bm j} \overline{\mathcal{Y}_0}\biggr) \biggr\rvert.
	\end{equation}
	
	We now estimate the $Z^\mathrm{iso}_{ab}$ quantities. 
	For the terms with $|\bm l| +\sumJ = 1$ and $\bm l = \bm 0$, it suffices to consider $J^* = \emptyset$, as the other case $J = \emptyset$ is analogous. For $J = \{(1,0)\}$, using \eqref{eq:S_bound} and \eqref{eq:G_trivial}, we obtain
	\begin{equation}
		Z_{ab}^\mathrm{iso}\bigl({\bm 0}, \{(1,0)\}, \emptyset\bigr) = \biggl\lvert \sum_c S_{ac} G_{cb} G_{aa}G_{cb} \biggr\rvert \lesssim \frac{1}{W}\sum_c  \lvert G_{cb}\rvert^2 = \frac{(\im G)_{bb}}{W\eta}  \lesssim \frac{1}{W\eta^2} \lesssim \frac{1}{W},
	\end{equation}
	where we used $\partial_{ac}  \mathcal{Y}_0 = \partial_{ac} G_{ab}$ and the Ward identity. Similarly, for $J = \{(0,1)\}$, 
	we bound
	\begin{equation} \label{eq:1Gflat_offdiag_Z1}
		Z_{ab}^\mathrm{iso}\bigl({\bm 0}, \{(0,1)\}, \emptyset\bigr) = \biggl\lvert \sum_c S_{ac} G_{cb} G_{ac}G_{ab} \biggr\rvert \lesssim \frac{|G_{ab}|}{W}\sum_c  \lvert G_{ac}G_{cb}\rvert \le \frac{\sqrt{(\im G)_{aa}(\im G)_{bb}}}{W\eta}  \lesssim \frac{1}{W}.
	\end{equation}

	Hence, using \eqref{eq:G_trivial}, for any $\bm j \in (\mathbb{Z}_{\ge 0})^2\backslash\bm 0$, we obtain, for some positive constant $C_{|\bm j|}$ dependent only on the number of derivatives $|\bm j|$, 
	\begin{equation} \label{eq:partialG_bounds1}
		\bigl\lvert \partial_{ac}^{\bm j}  G_{yb}\bigr\rvert \le \sum_{\substack{(a_i, c_i) \in \{(a,c),(c,a)\}\\ i \in \indset{\,|\bm j|\, }}} \bigl\lvert G_{xa_1} G_{c_1a_2}\dots G_{c_{|\bm j|-1}a_{|\bm j|}} G_{c_{|\bm j|}b} \bigr\rvert \le C_{|\bm j|} \bigl(|G_{ab}| + |G_{cb}|\bigr), \quad y\in\indset{N},
	\end{equation}
	and an analogous bound holds with $G$ replaced by $\overline{G}$. The contribution coming from the $|G_{ab}|$ factors is estimated using a self-consistent high-moment bound with a Young's inequality 
	 (see~\eqref{eq:enough_to_bound_Z} for the same procedure), while the $|G_{cb}|$ factors are estimated using \eqref{eq:max_ansatz} for $c\neq b$ and \eqref{eq:G_trivial} for $c=b$. 
	Since $\partial_{cd}\mathcal{Y}_0 = \partial_{cd}G_{ab}$, using \eqref{eq:partialG_bounds1}, for $|\bm l| + \sumJ \ge 2$ and $a \neq b$, we obtain
	\begin{equation} \label{eq:1Gflat_offdiag_Z2}
		\begin{split}
			Z_{ab}^\mathrm{iso}(\bm l, J, J^*) &\lesssim \sum_{c} \bigl(S_{ac}\bigr)^{(1 + |\bm l| + \sumJ)/2}  \biggl(|G_{ab}|^{1+|J| +|J^*|} + |G_{cb}|^{1+|J| +|J^*|}\biggr) \\
			&\lesssim W^{-(|\bm l| + \sumJ-1)/2}\biggl(|\mathcal{Y}_0|^{1+|J| +|J^*|} + W^{-1}  + W^{-1/2}\bigl(W^{-1/2}\psi\bigr)^{|J| + |J^*|}\biggr),
		\end{split}
	\end{equation} 
	with very high probability, where we used \eqref{eq:S_bound}, \eqref{eq:max_ansatz}, and \eqref{eq:Schwarz+Ward}. 
	Since $|\bm l| + \sumJ \ge 2$ and $\sumJ \ge |J| + |J^*|$,  we have $|\bm l| + \sumJ-1 \ge (|J| +|J^*|)/2$, hence 
	\begin{equation} \label{eq:1Gflat_offdiag_final} 
		Z_{ab}^\mathrm{iso}(\bm l, J, J^*)  \lesssim W^{-1/2}|\mathcal{Y}_0|^{1+|J| +|J^*|} + W^{-(1+|J| +|J^*|)/2}  + W^{-1/2}\bigl(W^{-3/4}\psi\bigr)^{|J| + |J^*|}, 
	\end{equation}
	with very high probability for all $|\bm l| + \sumJ \ge 2$. Hence, from \eqref{eq:1Gflat_offdiag_cuexp}, using \eqref{eq:1Gflat_offdiag_Z1} and \eqref{eq:1Gflat_offdiag_final}, we conclude that \eqref{eq:flat1G_underline} holds for $a \neq b$.
	
	Next, we consider the case $a=b$. It follows from the bound \eqref{eq:stab_bound} and \eqref{eq:Binv_action} that
	\begin{equation} \label{eq:diag_stab_flat_bound}
		\biggl\lvert\Expv\biggl[\bigl(\mathcal{B}^{-1}[\underline{HG}]\bigr)_{aa}\overline{\mathcal{Y}_0} \lvert \mathcal{Y}_0 \rvert^{2p-2}\biggr]\biggr\rvert  \lesssim \max_{d} \biggl\lvert\Expv\biggl[\bigl(\underline{HG}\bigr)_{dd}\overline{\mathcal{Y}_0} \lvert \mathcal{Y}_0 \rvert^{2p-2}\biggr]\biggr\rvert,
	\end{equation}
	where we recall that $\mathcal{Y}_0 = (G-m)_{ab} = (G-m)_{aa}$ since $a=b$ in the case under consideration. In particular, it suffices to estimate the quantities $Z^\mathrm{iso}_{dd}$ for all $d \in \indset{N}$
	in the formula analogous to~\eqref{eq:1Gflat_offdiag_cuexp}.

	First, we consider $\sumJ = 1$ and $|\bm l| = 0$. For $ J =\{(1,0)\}$, similarly to \eqref{eq:1Gflat_offdiag_Z1}, we have
	\begin{equation} \label{eq:1Gflat_diag_Z1}
		Z_{dd}^\mathrm{iso}\bigl(\bm 0, \{(1,0)\}, \emptyset\bigr) = \biggl\lvert \sum_c S_{dc} G_{cd} G_{ad}G_{ca} \biggr\rvert \lesssim \frac{1}{W}\sum_c  \lvert G_{ca} G_{cd} \rvert \le \frac{\sqrt{(\im G)_{aa}(\im G)_{dd}}}{W}  \lesssim \frac{1}{W}.
	\end{equation}
	The other terms involving second-order cumulants are estimated analogously.	
	
	Next, we consider terms with higher-order cumulants, $|\bm l| + \sumJ \ge 2$. Similarly to \eqref{eq:partialG_bounds1}, the partial derivatives of $G$ satisfy, $\bm j =(j_1,j_2) \in (\mathbb{Z}_{\ge 0})^2\backslash \bm 0$,
	\begin{equation}
		\bigl\lvert \partial_{dc}^{\bm j} G_{aa} \bigr\rvert \lesssim \mathds{1}_{|\bm j| \ge 2} + |G_{ac}| + |G_{ca}|, \quad \text{w.v.h.p.},
	\end{equation}
	and the same bound holds with $G$ replaced by $\overline{G}$.
	Therefore, for all $|\bm l| + \sumJ \ge 2$, recalling $\mathcal{Y}_0 = (G-m)_{aa}$, we obtain, using \eqref{eq:sumS=1}, \eqref{eq:max_ansatz}, \eqref{eq:Schwarz+Ward}, and Young's inequality,
	\begin{equation} \label{eq:1Gflat_diag_Z2}
		\begin{split}
			Z_{dd}^\mathrm{iso}(\bm l, J, J^*) &\lesssim \sum_{c} \bigl(S_{dc}\bigr)^{(1+|\bm l| + \sumJ)/2}  \bigl( \mathds{1}_{|\bm l| \ge 1} + |G_{cd}| \bigr) \bigl(  |G_{ac}| + |G_{ca}|\bigr)^{ n_1}\\
			&\lesssim W^{-(1+|J| +|J^*|)/2} + \psi^{n_1} W^{-(|\bm l| + \sumJ+n_1-1)/2} \bigl(\mathds{1}_{|\bm l|\ge1} + W^{-1/2}\bigr),
		\end{split}
	\end{equation}
	with very high probability, where $n_1 := |\{\bm j \in J, J^* \,:\, |\bm j| = 1\}|$. 
	To finish the estimate, we use $|\bm l| + \sumJ \ge 2$,  $\sumJ  + n_1 \ge 2 |J| + |J^*|$, and consider the following cases:
	\begin{itemize}
		\item  If $|\bm l| \neq 1$, or $|\bm l| = 1$ and $\sumJ \ge 2$, then
		\begin{equation}
			Z_{dd}^\mathrm{iso}(\bm l, J, J^*) \lesssim W^{-(1+|J| +|J^*|)/2}\bigl(1 + W^{-1/4}\psi\bigr)^{n_1}, \quad \text{w.v.h.p.}
		\end{equation}
		\item  If $|\bm l| = 1$ and $\sumJ = 1$, then
		\begin{equation}
			Z_{dd}^\mathrm{iso}(\bm l, J, J^*) \lesssim \bigl(W^{-1/2} \sqrt{\psi}\bigr)^{1+|J| +|J^*|}, \quad \text{w.v.h.p.}
		\end{equation}
	\end{itemize}
 	Therefore, for  $|\bm l| + \sumJ \ge 2$, we conclude that
 	\begin{equation} \label{eq:Z1iso_diag_high}
 		Z_{dd}^\mathrm{iso}(\bm l, J, J^*) \lesssim \biggl(W^{-1/2}\bigl( \sqrt{\psi} + W^{-1/4}\psi\bigr)\biggr)^{1+|J| +|J^*|}.
 	\end{equation}

	Therefore, using Young's inequality, we deduce from \eqref{eq:diag_stab_flat_bound}, \eqref{eq:Z1iso_diag_high}, and \eqref{eq:1Gflat_diag_Z1}, that the desired \eqref{eq:flat1G_underline} also holds for $a=b$. This concludes the proof of Claim \ref{claim:flat1G_underline}.
\end{proof}

Next, we prove Claim \ref{claim:omega_underline1}. Note that in this proof and the proof of Claim \ref{claim:omega_underline2} immediately below, we deal with modified size functions $\mathfrak{s}_{k}(\omega;\, \cdot\,)$ for $k \in \{1,2\}$. 
\begin{proof}[Proof of Claim \ref{claim:omega_underline1}]
	To prove \eqref{eq:1G_decay_high_moment}, we again use the cumulant expansion \eqref{eq:1Gflat_offdiag_cuexp};
	note that $\mathcal{Y}_0=G_{ab}$ for $a\ne b$.  
	We now estimate the contributions of $Z_{ab}^\mathrm{iso}$, defined in \eqref{eq:Ziso}, coming from each cumulant order, starting with the second cumulants.
	
	Cumulants of the second order correspond to $|\bm l| = 0$ and $\sumJ = 1$ in the expansion \eqref{eq:1Gflat_offdiag_cuexp}. We only treat the case $J = \{(1,0)\}$ and $J = \{(0,1)\}$ with $J^* = \emptyset$ in detail, since the other two terms are analogous. 
	The quantity $Z_{ab}^\mathrm{iso} (0,\{(1,0)\}, \emptyset ) $, defined in \eqref{eq:Ziso}, satisfies
	\begin{equation} \label{eq:omega_2cum}
		\begin{split}
			Z_{ab}^\mathrm{iso}\bigl(\bm 0,\{(1,0)\}, \emptyset\bigr) &= \biggl\lvert\sum_{c} S_{ac} \bigl(G_{cb}\bigr)^2 G_{aa}\biggr\rvert \le |G_{aa}|\, (G^*S^aG\bigr)_{bb} \lesssim \Upsilon_{ab} + \omega \\
			&\lesssim N^{\zeta}\bigl(\mathfrak{s}_1(N^{-\zeta}\omega; a,b)\bigr)^2,
		\end{split}
	\end{equation}
	with very high probability, where in the penultimate step we used \eqref{eq:G_trivial}, \eqref{eq:omegas_def} and \eqref{eq:2G_omega_ansatz}, while in the last step we used \eqref{eq:omega_ineqs}.

	Similarly, it follows by Schwarz inequality and \eqref{eq:2G_omega_ansatz}, that $Z_{ab}^\mathrm{iso} ({\bm 0},\{(0,1)\}, \emptyset )$ satisfies
	\begin{equation} \label{eq:omega_2cum2}
		\begin{split}
			Z_{ab}^\mathrm{iso} (\bm 0,\{(0,1)\}, \emptyset ) &= \biggl\lvert\sum_{c} S_{ac} G_{cb} G_{ac} G_{ab}\biggr\rvert \le |G_{ab}| \sqrt{ (G^*S^aG\bigr)_{bb} (GS^aG^*\bigr)_{aa} }\\
			&\lesssim W^{-1/2}|G_{ab}|N^{\zeta/2}\,\mathfrak{s}_1(N^{-\zeta}\omega; a,b) \lesssim W^{-1/2}N^{3\zeta}\bigl(\mathfrak{s}_1(N^{-\zeta}\omega; a,b)\bigr)^2,
		\end{split}
	\end{equation}
	with very high probability, where we used \eqref{eq:omega_assume} with $a\neq b$,  \eqref{eq:2G_flat}, \eqref{eq:2G_omega_ansatz}, and \eqref{eq:omega_ineqs}. Therefore, the contribution from the terms involving second order cumulants in \eqref{eq:1Gflat_offdiag_cuexp} is bounded by 
	\begin{equation}
		C N^{\zeta}\bigl(\mathfrak{s}_1(N^{-\zeta}\omega; a,b)\bigr)^2 \Expv\bigl[|G_{ab}|^{2p-2}\bigr].
	\end{equation}
	
	Next, we estimate the contributions coming from terms involving higher-order cumulants.  
	As in \eqref{eq:partial_X_1}, our goal is to recover a full $\mathfrak{s}_1(N^{-\zeta}\omega; a,b)$ decaying factor from each $\partial^{\bm j}_{cd}$ in \eqref{eq:Ziso}. 
	To this end, we use the following consequence of \eqref{eq:SUps_comvol} and the triangle inequality \eqref{eq:omega_s_triag}, for all $c\neq b$,
	\begin{equation} \label{eq:omega_triag}
		\bigl(S_{ac}\bigr)^{1/2}|G_{cb}| \lesssim N^{2\zeta} \sqrt{\Upsilon_{ac} (\Upsilon_{cb} + \omega)} \lesssim W^{-1/2}N^{2\zeta} \sqrt{\Upsilon_{a b} + \omega } \lesssim W^{-1/2}N^{5\zeta/2} \mathfrak{s}_1(N^{-\zeta}; a,b),
	\end{equation}
	with very high probability, where we used \eqref{eq:omegas_def}, \eqref{eq:omega_assume},  and \eqref{eq:omega_ineqs}.
	
	Therefore, for $|\bm l| + \sumJ \ge 2$ with $|J| + |J^*| \ge 1$, plugging \eqref{eq:partialG_bounds1} into \eqref{eq:Ziso}, using \eqref{eq:omega_triag} $(|J| + |J^*|-1)$ times and \eqref{eq:2G_omega_ansatz}, we obtain
	\begin{equation} \label{eq:Ziso_dest_chains}
		\begin{split} 
			Z_{ab}^\mathrm{iso}(\bm l, J, J^*) \lesssim&~ \sum_{c} \bigl(S_{ac}\bigr)^{(1 + |\bm l| + \sumJ)/2} \bigl(|G_{ab}|^{1+|J| +|J^*|} + |G_{cb}|^{1+|J| +|J^*|}\bigr)\\
			\lesssim&~ W^{-(|\bm l| + \sumJ-1)/2}|G_{ab}|^{1+|J| +|J^*|}  + \bigl(\Upsilon_{ab}\bigr)^{(1 + |\bm l| + \sumJ)/2} \\
			&+\bigl(N^{2\zeta} \bigr)^{|J| + |J^*|-1} W^{-(|\bm l| + \sumJ-1)/2}\bigl(N^{\zeta/2}\mathfrak{s}_1(N^{-\zeta}\omega; a,b)\bigr)^{1+|J| +|J^*|},
		\end{split}
	\end{equation}
	with very high probability, where we used \eqref{eq:SUps_comvol}, \eqref{eq:S_bound},  and \eqref{eq:sumS=1}. 
	
	On the other hand, if $|J| + |J^*| = 0$ and $|\bm l| \ge 2$, by Schwarz inequality, \eqref{eq:2G_omega_ansatz}, and \eqref{eq:omega_ineqs}, we have
	\begin{equation} \label{eq:Ziso_no_dest_chains}
		\begin{split}
			Z_{ab}^\mathrm{iso}(\bm l, \emptyset, \emptyset) \lesssim&~ \sum_{c} \bigl(S_{ac}\bigr)^{(|\bm l|+1)/2} \bigl(|G_{ab}|  + |G_{cb}| \bigr)\\
			\lesssim&~  |G_{ab}|  W^{-(|\bm l|-1)/2}  +  W^{-(|\bm l|-1)/2}  N^{\zeta/2}\mathfrak{s}_1(N^{-\zeta}\omega; a,b) ,
		\end{split}
	\end{equation}
	with very high probability. 
	Therefore, combining \eqref{eq:omega_ineqs}, \eqref{eq:Ziso_dest_chains} and \eqref{eq:Ziso_no_dest_chains}, we conclude that
	\begin{equation} \label{eq:omega_Lcum}
		\begin{split}
			Z_{ab}^\mathrm{iso}(\bm l, J, J^*) \lesssim&~ W^{-1/2}|G_{ab}|^{1+|J| +|J^*|} +\bigl(\mathfrak{s}_1(N^{-\zeta}\omega; a,b)\bigr)^{1+|J| +|J^*|}, \quad \text{w.v.h.p.}
		\end{split}
	\end{equation}
	
	Hence, it follows from \eqref{eq:1Gflat_offdiag_cuexp}, \eqref{eq:omega_2cum}, \eqref{eq:omega_2cum2}, and \eqref{eq:omega_Lcum}, that \eqref{eq:1G_decay_high_moment} holds for all $a\neq b$. This concludes the proof of Claim \ref{claim:omega_underline1}.
\end{proof}

We conclude the section by proving Claim \ref{claim:omega_underline2}. 
\begin{proof}[Proof of Claim \ref{claim:omega_underline2}]  Recall that $\mathcal{Y}_1 := ((G-M)_{[1,2]}(x))_{ab}$.
	To prove \eqref{eq:omega_underline2}, we use cumulant expansion \eqref{eq:k_iso_cumexp} with $k=1$  and $A_1' := A'$. 
	
	When estimating $Z_{1}^\mathrm{iso}(\bm l, J, J^*)$, defined in \eqref{eq:Zkiso_def}, our goal is to recover $(|J| + |J^*| + 1)$ decaying factors $\mathfrak{s}_2(N^{-\zeta}\omega;a,x,b)$. 
	First, we consider the terms with $\sumJ = 1$ and $|\bm l| = 0$. We analyze the term corresponding to $J = \{ (1,0)\}$ in detail, that is
	\begin{equation} \label{eq:Z_2_cum2_sum}
		Z_{1}^\mathrm{iso}(\bm 0, \{(1,0)\}, \emptyset) = \biggl\lvert\sum_{c} \sum_d A'_{cc}  S_{cd}  (G_1)_{ac}(G_2)_{db}   \biggl((G_1)_{ac} \bigl(G_{[1,2]}\bigr)_{db}  +   \bigl(G_{[1,2]}\bigr)_{ac} (G_2)_{db} \biggr) \biggr\rvert.
	\end{equation} 
	Using \eqref{eq:M_bound} and \eqref{eq:omega_assume} for every $G_{[1,2]}$ and $G_{i}$ factor, we obtain that the contribution to \eqref{eq:Z_2_cum2_sum} coming from $c\neq a$ and $d\neq b$ admits the bound
	\begin{equation} \label{eq:Z_2_cum2_offdiag}
		\begin{split}
			&\biggl\lvert\sum_{c\neq a} \sum_{d\neq b}  A'_{cc}  S_{cd}  (G_1)_{ac}(G_2)_{db}   \biggl((G_1)_{ac} \bigl(G_{[1,2]}\bigr)_{db}  +   \bigl(G_{[1,2]}\bigr)_{ac} (G_2)_{db} \biggr) \biggr\rvert\\
			&\quad\lesssim\frac{N^{6\zeta}\psi}{\sqrt{W}}\sqrt{\Upsilon_{xb}+\omega}\sum_{c}\Upsilon_{xc} \bigl(\Upsilon_{ac}+\omega\bigr) \sum_d S_{cd} \sqrt{\Upsilon_{db}+\omega}  \sqrt{\Upsilon_{dx}+\omega}\\ 
			&\qquad+   \frac{N^{6\zeta}\psi}{\sqrt{W}}\sqrt{\Upsilon_{ax}+\omega} \sum_{c}\Upsilon_{xc}  \sqrt{\Upsilon_{ac}+\omega}\sqrt{\Upsilon_{xc}+\omega}\sum_d S_{cd}\bigl(\Upsilon_{db}+\omega\bigr)\\
			&\quad\lesssim \frac{N^{6\zeta}\psi}{\sqrt{W}}\mathfrak{s}_2(\omega; a,x,b)\sum_{c}\Upsilon_{xc} \sqrt{\Upsilon_{ac}+\omega} \sqrt{\Upsilon_{cb}+\omega} \lesssim   N^{6\zeta}\psi\bigl(\mathfrak{s}_2(\omega; a,x,b) \bigr)^2,
		\end{split}
	\end{equation}
	with very high probability, where in the penultimate step we used \eqref{eq:SUps_comvol}, \eqref{eq:omega_s_triag},  and the assumption $\omega \lesssim W^{-1}$, and in the last step we used \eqref{eq:omega_s_convol}.
	Similarly, for $c=a$ and $d \neq b$, we obtain
	\begin{equation} \label{eq:Z_2_cum2_diag}
		\begin{split}
			&\biggl\lvert\sum_{d\neq b}  A'_{aa}  S_{ad}  (G_1)_{aa}(G_2)_{db}   \biggl((G_1)_{aa} \bigl(G_{[1,2]}\bigr)_{db}  +   \bigl(G_{[1,2]}\bigr)_{aa} (G_2)_{db} \biggr) \biggr\rvert\\
			&\quad\lesssim \frac{N^{2\zeta}\psi}{\sqrt{W}}\Upsilon_{xa} \sqrt{\Upsilon_{xb}+\omega} \sum_{d }    S_{ad}  \sqrt{(\Upsilon_{db}+\omega)(\Upsilon_{dx}+\omega)} +  \Upsilon_{ax} \bigl(\Upsilon_{xa} + \omega\bigr) \sum_{d}  S_{ad}  \bigl(\Upsilon_{db}+\omega\bigr)  \\
			&\quad\lesssim  N^{2\zeta}\psi\,\bigl(\mathfrak{s}_2(\omega; a,x,b)\bigr)^2,
		\end{split}
	\end{equation}
	where we used \eqref{eq:M_bound}, \eqref{eq:omega_s_triag} and \eqref{eq:omega_s_convol}. Contribution to the sum in \eqref{eq:Z_2_cum2_sum} coming from $d=b$ is estimated similarly. Hence, we obtain
	\begin{equation} 
		Z_{1}^\mathrm{iso}(\bm 0, \{(1,0)\}, \emptyset) \lesssim N^{6\zeta}\psi\bigl(\mathfrak{s}_2(\omega; a,x,b) \bigr)^2, \quad \text{w.v.h.p.} 
	\end{equation}
	Other terms involving second order cumulants are estimated analogously. Therefore, using \eqref{eq:omega_ineqs}, we conclude that
	\begin{equation} \label{eq:Z2_l=0j=1}
		\sum_{\sumJ = 1} Z_{1}^\mathrm{iso}(0, J, J^*)  \lesssim N^{6\zeta}\psi\bigl(\mathfrak{s}_2(N^{-\zeta}\omega; a,x,b) \bigr)^2, \quad \text{w.v.h.p.}
	\end{equation}
	
	Next, we estimate the contribution of the terms involving higher-order cumulants. Similarly to \eqref{eq:partial_X_1} and \eqref{eq:partial_Y}, we obtain  
	\begin{equation} \label{eq:partial_G_bounds2}
		\bigl\lvert \partial_{cd}^{\bm j}  \bigl(G_{[1,2]}\bigr)_{ab}\bigr\rvert 
		\le C_{|\bm j|}\sum_{q=1}^{|\bm j|+1}\sum_{\substack{(c_i, d_i) \in \{(c,d),(d,c)\}\\ i \in \indset{\,|\bm j|\, }}} \bigl\lvert (\mathcal{G}_1^q)_{ac_1} (\mathcal{G}_2^q)_{d_1c_2}\dots (\mathcal{G}_{|\bm j|}^q)_{d_{|\bm j|-1}c_{|\bm j|}} (\mathcal{G}_{|\bm j|+1}^q)_{c_{|\bm j|}b} \bigr\rvert,
	\end{equation}
	where, for $j,q\in\mathbb{N}$, the matrices  $\mathcal{G}_{j}^q$ denote
	\begin{equation}
		\mathcal{G}_j^q := \begin{cases}
			G_1, \quad& j < q,\\
			G_{[1,2]}, \quad &j=q,\\
			G_2, \quad &j > q.
		\end{cases}
	\end{equation}
	and the same bound holds for $G_{[1,2]}$ replaced by $\overline{G_{[1,2]}}$. 
	
	Note that only $\mathcal{G}_{j}^q$ with $j \in \{1, q, |\bm j| + 1\}$ in \eqref{eq:partial_G_bounds2} carry the dependence on the original indices $a,x,b$. We estimate them by \eqref{eq:omega_chain_ansatz}, which results in $\sqrt{\Upsilon + \omega}$-decay factors coupling the original indices $a,x,b$ to (potentially different) summation indices $c$ and $d$.
	We bring these factors together using $\sqrt{S_{cd}} \lesssim \sqrt{\Upsilon_{cd}}$ from \eqref{eq:Zkiso_def} and \eqref{eq:omega_s_triag} to bridge the gap between $c$ and $d$, thus recovering the target spatial decay $\mathfrak{s}_2(N^{-\zeta}\omega; a,x,b)$.

	Note that, similarly to \eqref{eq:partial_Y}, using \eqref{eq:triag}, the first bound in \eqref{eq:SUps_comvol}, \eqref{eq:omega_chain_ansatz} and \eqref{eq:omega_s_triag}, \eqref{eq:omega_ineqs},
	 we obtain, for all $\bm j \in (\mathbb{Z}_{\ge 0})^2 \backslash \{\bm 0\}$,
	\begin{equation} \label{eq:partial_G_bounds2_triag}
		\bigl(S_{cd}\bigr)^{|\bm j|/2}\bigl\lvert \partial_{cd}^{\bm j}  \bigl(G_{[1,2]}\bigr)_{ab}\bigr\rvert 
		\lesssim \frac{C_{|\bm j|}}{W^{|\bm j|/2}} \biggl(\frac{N^{3\zeta}\psi}{\sqrt{W}} + \other{\delta}_1 N^\zeta\psi + \bigl(\other{\delta}_2 + \other{\delta}_1 \mathds{1}_{|\bm j| \ge 2}\bigr)\sqrt{W}\biggr)  \mathfrak{s}_2(N^{-\zeta}\omega; a,x,b),
	\end{equation}
	with very high probability, where $\other{\delta}_i$ are defined in \eqref{eq:other_deltas}. Here, we also 
	used \eqref{eq:triag} to estimate the products of $\Upsilon$'s without the additive tail bound $\omega$, that arise from the $\delta_{ab}$ in \eqref{eq:omega_chain_ansatz}. 

	Completely analogously, using the bound \eqref{eq:A'_bound} for $A'$, we deduce that
	\begin{equation}
		\begin{split}
			\bigl(S_{cd}\bigr)^{|\bm l|/2}\bigl\lvert \partial_{cd}^{\bm l}  \bigl((G_1)_{ac}(G_{2})_{db}\bigr) \bigr\rvert 
			\lesssim&~ C_{|\bm l|+1} W^{-|\bm l|/2}  \mathfrak{s}_1(N^{-\zeta}\omega; a,c)\mathfrak{s}_1(N^{-\zeta}\omega; d,b)\\
			&\times \bigl(N^{5\zeta} + \other{\delta}_1 N^{5\zeta/2}\sqrt{W} + \bigl(\other{\delta}_2 + \other{\delta}_1 \mathds{1}_{|\bm l| \ge 1}\bigr)W\bigr) ,
		\end{split}
	\end{equation}
	\begin{equation} \label{eq:partial_l_2G_bounds}
		\begin{split}
			\bigl(S_{cd}\bigr)^{(|\bm l|+1)/2}\bigl\lvert \partial_{cd}^{\bm l}  \bigl((G_1)_{ac}(G_{2})_{db}\bigr) \bigr\rvert \bigl\lvert (A')_{cc} \bigr\rvert
			\lesssim&~ C_{|\bm l|+1} W^{-(|\bm l|+1)/2} \mathfrak{s}_2(N^{-\zeta}\omega; a,x,b)\\
			&\times \biggl(\frac{N^{5\zeta}}{\sqrt{W}} + \other{\delta}_1 N^{5\zeta/2}+ \bigl(\other{\delta}_2 + \other{\delta}_1 \mathds{1}_{|\bm l| \ge 1}\bigr)\sqrt{W} \biggr) ,
		\end{split}
	\end{equation}
	with very high probability, for all $\bm l \in (\mathbb{Z}_{\ge 0})^2$.
		
  We can now put together the estimate of \eqref{eq:Zkiso_def} with $k=1$ and $A_1 := A'$. Using \eqref{eq:partial_G_bounds2_triag}--\eqref{eq:partial_l_2G_bounds}, and recalling the general heuristic that each unrestricted $c$ and $d$ summation brings a factor of $W$, similarly to \eqref{eq:Zkiso_high_final}, we assert that $Z_1^\mathrm{iso}(\bm l, J, J^*) $ satisfies  
	\begin{equation} \label{eq:Z2iso_high_final}
		Z_1^\mathrm{iso}(\bm l, J, J^*) \lesssim \bigl(\sqrt{N^{\zeta}\psi} + W^{-1/4}N^\zeta\psi \bigr)^{|J| + |J^*| + 1} \bigl( \mathfrak{s}_2(N^{-\zeta}\omega; a,x,b)\bigr)^{|J| + |J^*| + 1}, \quad \text{w.v.h.p.}
	\end{equation}
	 It is straightforward to check that \eqref{eq:Z2iso_high_final} indeed holds using precisely the same procedure as in \eqref{eq:partial_X2}--\eqref{eq:Zkav_high} and \eqref{eq:partial_Y}--\eqref{eq:Zkiso_high_final}. 
	 In particular, we estimate each chain in \eqref{eq:partial_G_bounds2} by \eqref{eq:omega_chain_ansatz}, and use \eqref{eq:omega_s_triag} enough times so that at most two indices $d$ are left among the arguments of $\sqrt{\Upsilon+\omega}$ across all $\bm j$- and $\bm l$-derivatives in \eqref{eq:Zkiso_def} while preserving at least a full power of $S_{cd}$. Hence, for $\other{\delta}_1 = 0$, we conclude \eqref{eq:Z2iso_high_final} by using $S_{cd}$ to perform the $d$ summation, and \eqref{eq:A'_bound} to perform the $c$ summation by \eqref{eq:omega_s_convol}. The cases $\other{\delta}_1 \ge 1$ with $\other{\delta}_2=0$, and $\other{\delta}_2 \ge 1$ are estimated similarly.

	Hence, the desired \eqref{eq:omega_underline2} follows by combining \eqref{eq:k_iso_cumexp} with $k:=1$, \eqref{eq:Z2_l=0j=1}, and \eqref{eq:Z2iso_high_final} by the same reasoning as in \eqref{eq:enough_to_bound_Z}. 
	This concludes the proof of Claim \ref{claim:omega_underline2}.
\end{proof}

\section{Green function comparison: Proof of Proposition \ref{prop:zag}} \label{sec:GFT}
In this section we remove the Gaussian component that was added to the random matrix during the zig-step. Since the time-parameter $t$ remains fixed throughout the proof, to condense the notation, we omit the subscript $t$ from the quantities $z_t, m_t, \ell_t, \eta_t$, the deterministic approximations $M_t$ and the size functions $\Upsilon_t$, $\mathfrak{s}_t^{\mathrm{iso/av}}$.

The basic idea of the proof is a one-by-one replacement strategy analogous to 
Lindeberg's original proof for the central limit theorem. This method was introduced to random matrix theory
first by Chatterjee \cite{Chatterjee}, 
followed by the fundamental work of Tao and Vu  \cite{tao2011random} on their Four Moment Theorem. 
It was then systematically applied for resolvents \cite{erdHos2011universality}, 
 coined the Green Function Comparison Theorem (GFT),
as summarized in Chapter 16 of \cite{erdHos2017dynamical}. 
In our application, we transfer the local law from $V$ to $H$ by replacing each entry $v_{ij}$ of $V$ with $h_{ij}$ one-by-one
and monitor the change of $\Psi_{k,t}^{\mathrm{av/iso}}$, defined in \eqref{eq:Psi_def}, in each step.
Finally we sum up the effects of each replacement. Technically,
high probability bounds will be regularly swapped with high moment bounds since the latter can be effectively
controlled using the matched moments. 
 
After setting up the notations in Section~\ref{sec:prelim}, we will first prove  Proposition~\ref{prop:zag}
for  the isotropic single resolvent local laws ($k=1$) in Section~\ref{sec:1Grep}. 
Armed with the control of a single resolvent, in Section~\ref{sec:k_iso_rep} we prove isotropic comparisons
for chains of any length $k\le K$ by induction on $k$. Finally, comparison for averaged chains of arbitrary 
length  is proven in Section~\ref{sec:k_av_rep} using the isotropic comparison results.

We stress that in most versions of the GFT in the literature the single resolvent isotropic 
local law is a given input, proven beforehand
by other methods.  Lacking alternatives for band matrices, here we even need to prove the
single resolvent result self-consistently via the zigzag strategy, similarly to~\cite{campbell2024spectral}
(a partially self-consistent way to prove GFT for a single resolvent has already appeared in
\cite{knowles2017anisotropic}).

We also remark that besides the one-by-one replacement strategy, there is another comparison
approach, the continuity of the resolvents along the Ornstein-Uhlenbeck flow
(see e.g. Chapter 15 of \cite{erdHos2017dynamical}). It is technically simpler than the one-by-one replacement strategy
and it is applicable for a broader class of models (especially for random matrices with correlated entries),
but it is less precise as it matches only the first two moments. To achieve the required precision for band matrices,
 we need to match the third moment as well, hence we resort to the one-by-one replacement.
In most applications of the zigzag strategy, the continuity method was sufficient,  with the notable exception 
of \cite{cipolloni2023eigenstate}, where optimally controlling the error in terms of the Hilbert-Schmidt norm in the multi-resolvent local law
required more precise analysis, that only the one-by-one strategy could provide.

We emphasize that in  Section \ref{sec:GFT}, we use only the weaker convolution estimate \eqref{eq:true_convol_notime} (in its time-dependent forms \eqref{eq:true_convol},  \eqref{eq:Schwarz_convol}, and \eqref{eq:sqrt_convol}) instead of \eqref{eq:convol_notime}--\eqref{eq:suppressed_convol_notime}.

\subsection{Preliminaries}\label{sec:prelim}
For the entire proof, we fix a bijective map $\phi \equiv \phi_N$,
\begin{equation}
	\phi : \{(i,j) \in \indset{N}^2\,:\, i \le j\} \to \indset{\gamma(N)}, \quad \gamma(N) := N(N+1)/2.
\end{equation}
The bijection $\phi$ defines the order in which we replace the independent entries of the $N\times N$ random matrix $H$. For all $\gamma := \phi\bigl((i,j)\bigr) \in \indset{\gamma(N)}$, we define the random matrices with replaced entries as
\begin{equation}
	H^{(0)} := V, \quad H^{(\gamma)} := H^{(\gamma - 1)} + \sqrt{S_{ij}}\bigl(\Delta^{(\gamma)}_{H}-\Delta_{V}^{(\gamma)}\bigr),
\end{equation}
where the matrices $\Delta^{(\gamma)}_{H}$ and $\Delta^{(\gamma)}_{V}$ are given by
\begin{equation}
	\bigl(\Delta^{(\gamma)}_{H}\bigr)_{ab} := \frac{h_{ij}\, E^{ij} + h_{ji}\, E^{ji} }{1 + \delta_{ij}}, \quad 
	\bigl(\Delta^{(\gamma)}_{V}\bigr)_{ab} := \frac{v_{ij}\, E^{ij}  + v_{ji}\, E^{ji}}{1 + \delta_{ij}}.
\end{equation}
Here, for all $a,b,i,j \in \indset{N}$, we denote $E_{ab}^{ij} := \delta_{ai}\delta_{bj}$, and $E^{ij} \equiv E^{ij}(N) := (E_{ab}^{ij})_{a,b=1}^N$. Note that $H^{(\gamma(N))}=H$. 

Furthermore, for all $\gamma := \phi\bigl((i,j)\bigr) \in \indset{\gamma(N)}$, we define the intermediate matrices 
$\widecheck{H}^{(\gamma)}$ with the $(i,j)$-entry of $H^{(\gamma )}$ replaced by zero, i.e. we have
\begin{equation}
	\widecheck{H}^{(\gamma)} := H^{(\gamma - 1)} - \sqrt{S_{ij}}\Delta_{V}^{(\gamma)} =
	H^{(\gamma)} - \sqrt{S_{ij}}\Delta_{H}^{(\gamma)}.
\end{equation} 
Similarly, we denote the corresponding resolvents by 
\begin{equation}  \label{eq:place_resolvents}
	G^{(\gamma-1)}(z) := \bigl(H^{(\gamma-1)} - z\bigr)^{-1}, \quad \widecheck{G}^{(\gamma)}(z) := \bigl(\widecheck{H}^{(\gamma)} - z\bigr)^{-1}, \quad G^{(\gamma)}(z) := \bigl(H^{(\gamma)} - z\bigr)^{-1},
\end{equation}
for all $z \in \mathbb{C}\backslash\mathbb{R}$ and $\gamma \in \indset{\gamma(N)}$.

The resolvents $G^{(\gamma)}$ and $G^{(\gamma -1)}$ admit the following  finite-order expansions around the intermediate resolvent $\widecheck{G}^{(\gamma)}$. 
\begin{lemma}
	For any $L \in \mathbb{Z}_{\ge 0}$ and $i,j\in \indset{N}$,
	the resolvents $G^{(\gamma)} \equiv G^{(\gamma)}(z)$ and $G^{(\gamma -1)} \equiv G^{(\gamma-1)}(z)$ 
	with $\gamma= \phi((i,j))$ satisfy
	\begin{equation} \label{eq:resolvent_expand}
		\begin{split}
			G^{(\gamma)} &= \sum_{q=0}^{L} \bigl(S_{ij}\bigr)^{q/2}\frac{1}{(1+\delta_{ij})^q} \frac{1}{q!} \bigl( h_{ij}\partial_{ij} + h_{ji}\partial_{ji} \bigr)^q \widecheck{G}^{(\gamma)}  +\bigl(S_{ij}\bigr)^{(L+1)/2}\bigl(-\widecheck{G}^{(\gamma)} \Delta^{(\gamma)}_H\bigr)^{L+1} G^{(\gamma)},\\
			G^{(\gamma-1)} &= \sum_{q=0}^{L} \bigl(S_{ij}\bigr)^{q/2}\frac{1}{(1+\delta_{ij})^q\,q!}
			 \bigl( v_{ij}\partial_{ij} +  v_{ji} \partial_{ji}  \bigr)^q  \widecheck{G}^{(\gamma)} 
			 +\bigl(S_{ij}\bigr)^{(L+1)/2} \bigl(-\widecheck{G}^{(\gamma)} \Delta^{(\gamma)}_V\bigr)^{L+1} G^{(\gamma-1)},
		\end{split}
	\end{equation}
	where $\widecheck{G}^{(\gamma)} \equiv \widecheck{G}^{(\gamma)}(z)$ is defined in \eqref{eq:place_resolvents}. 
	Here $\partial_{ij}$ denotes the partial derivative in the direction of the $(i,j)$ matrix entry of $\widecheck{H}^{(\gamma)}$, that is  for any smooth function $F$ we define\footnote{Note that this definition of $\partial_{ij}$ slightly differs from 
	the one used in all previous sections (e.g. in the definition of the martingale term
	 \eqref{eq:av_k_mart} or in the cumulant expansion~\eqref{eq:cumulant}). Since the current definition used only
	 locally in Section~\ref{sec:GFT}, where the previous derivative plays no role, this should not lead to any confusion.}
	\begin{equation}
		\partial_{ij}F\bigl(h_{ij}, h_{ji}, v_{ij}, v_{ji}, \widecheck{H}^{(\gamma)}\bigr) := \partial_{\varepsilon}  F\bigl(h_{ij}, h_{ji}, v_{ij}, v_{ji}, \widecheck{H}^{(\gamma)} + \varepsilon E^{ij}\bigr)\bigr\rvert_{\varepsilon=0}, \quad \gamma := \phi\bigl((i,j)\bigr).
	\end{equation} For example, $\partial_{ij}\widecheck{G}^{(\gamma)} = -\widecheck{G}^{(\gamma)}E^{ij}\widecheck{G}^{(\gamma)}$.
\end{lemma}
Note that (e.g., the first line of) \eqref{eq:resolvent_expand} is equivalent to the simple resolvent expansion
$$
G^{(\gamma)} = \widecheck{G}^{(\gamma)} - \widecheck{G}^{(\gamma)} (\sqrt{S_{ij}} \Delta_{H}^{(\gamma)})
\widecheck{G}^{(\gamma)} + \widecheck{G}^{(\gamma)} (\sqrt{S_{ij}} \Delta_{H}^{(\gamma)})
\widecheck{G}^{(\gamma)} (\sqrt{S_{ij}} \Delta_{H}^{(\gamma)})
\widecheck{G}^{(\gamma)} - \ldots
$$
but written with a derivative formalism that will allow us to draw parallels with the cumulant expansion~\eqref{eq:cumulant}
and thus with the analysis in Section~\ref{sec:cum_expand}. 
The main point is that the expansions \eqref{eq:resolvent_expand}   of $G^{(\gamma-1)}$
and $G^{(\gamma)}$ around the common resolvent $\widecheck{G}^{(\gamma-1)}$ 
explicitly highlight the dependence of these resolvents on the $(i,j)$ and $(j,i)$ matrix elements---the only  entries in which $H^{(\gamma-1)}$
and $H^{(\gamma)}$ differ. 

In the sequel, we consider the expansions \eqref{eq:resolvent_expand} with a finite $N$-independent order $L$ to be chosen later.  Note that the last terms in \eqref{eq:resolvent_expand} contain the non-checked  resolvents.

\subsection{Single resolvent isotropic comparison} \label{sec:1Grep}
Throughout the proof, we will estimate the resolvents $G^{(\gamma)}$ in the high-moment sense. To this end, for a random variable $\mathcal{X}$ and a positive integer $p\in \mathbb{N}$, we denote
\begin{equation}
	\norm{\mathcal{X}}_{p} := \Expv\bigl[|\mathcal{X}|^p\bigr]^{1/p}.
\end{equation}
The goal of this section is to prove the entry-wise bound a single resolvent of the random matrix $\other{H}$,
\begin{equation} \label{eq:1G_GFT_goal}
	\bigl\lVert \bigl(G^{(\gamma(N))}-m\bigr)_{ab}\bigr\rVert_{p} \lesssim C_p N^{\xi}\sqrt{\Upsilon_{ab}}, \quad a,b\in\indset{N}, \quad p \in \mathbb{N},
\end{equation}
for some positive constants $C_p$ depending on $p \in \mathbb{N}$ and the model parameters, provided
we have basically the same bound initially, i.e. 
\begin{equation} \label{eq:1G_GFT_init}
	\bigl\lVert \bigl(G^{(0)}-m\bigr)_{ab}\bigr\rVert_{p} \lesssim C_p N^{\xi}\sqrt{\Upsilon_{ab}}, \quad a,b\in\indset{N}, \quad p \in \mathbb{N},
\end{equation}
for some (possibly different) constants $C_p$. Here and in the sequel, the precise value of $C_p$ may vary from line to line
but it is independent of $N$, the indices $a, b$ 
 and of the running superscript $\gamma$. Note that the tolerance exponent $\xi$
does not deteriorate from \eqref{eq:1G_GFT_init} to \eqref{eq:1G_GFT_goal}.

The basic step in the replacement process is the content of the following lemma.
\begin{lemma}[Replacement Lemma] \label{lemma:1G_replacement}
	Assume that \eqref{eq:1G_GFT_init} holds for some $0 < \xi < \etaexp/40$.
	There exist constants $\{\other{C}_p\}_{p\in\mathbb{N}}$, depending only on $p$ and the model parameters, such that the following holds.
	
	Let $\gamma_0 \in \indset{\gamma(N)}$, and  let $\{\mu_p\}_{p\in\mathbb{N}}$, independent of $\gamma_0$,
	 be a set of positive constants, such that the resolvents $G^{(\gamma)}$ with $\gamma \le \gamma_0 - 1$ satisfy
	\begin{equation} \label{eq:1G_rep_assume}
		\max_{\gamma \in \indset{\gamma_0-1}}	\norm{\bigl(G^{(\gamma)}-m\bigr)_{ab}}_{2p}^{2p} \le \mu_p \bigl(N^{\xi} \sqrt{\Upsilon_{ab}}\bigr)^{2p}, \quad a,b \in \indset{N},
	\end{equation}
	for any integer $p \in\mathbb{N}$.
	 Then, the resolvent $G^{(\gamma_0)}$ satisfies
	\begin{equation} \label{eq:1G_rep_goal}
		\norm{\bigl(G^{(\gamma_0)}-m\bigr)_{ab}}_{2p}^{2p} \le \bigl(\other{C}_p+ 2N^{-\xi}\mu_p\bigr) \bigl(N^{\xi}\sqrt{\Upsilon_{ab}}\bigr)^{2p}, \quad a,b \in \indset{N},
	\end{equation}
	for any integer $p \in \mathbb{N}$.
\end{lemma}
\begin{proof} [Proof of \eqref{eq:1G_GFT_goal}]
	Once Lemma \ref{lemma:1G_replacement} is established, the desired entry-wise bound \eqref{eq:1G_GFT_goal} follows by a straightforward induction in $\gamma \in \indset{\gamma(N)}$ by choosing $C_p= (2\widetilde C_p)^{1/2p}$
	and using Lemma \ref{lemma:1G_replacement} with $\mu_p= C_p^{2p}$.
	 This concludes the proof of \eqref{eq:1G_GFT_goal}.
\end{proof}

Therefore, it remains to prove Lemma \ref{lemma:1G_replacement}. In preparation, we record the following claim that we prove at the end of Section \ref{sec:1Grep}.
\begin{claim} \label{claim:rep_check}
	Assume that \eqref{eq:1G_rep_assume} holds for all $p$.	
	Then, the resolvents $\widecheck{G}^{(\gamma)}$ with $\gamma \le \gamma_0$ satisfy
	\begin{equation} \label{eq:1G_rep_check}
		\max_{\gamma \in \indset{\gamma_0}}	\norm{\bigl(\widecheck{G}^{(\gamma)}-m\bigr)_{ab}}_{2p}^{2p} 
		\le \bigl(1+N^{-\xi}\bigr)\max_{\gamma \in \indset{\gamma_0}}	\norm{\bigl(G^{(\gamma-1)}-m\bigr)_{ab}}_{2p}^{2p}  + C_p\bigl(N^{\xi}\sqrt{\Upsilon_{ab}}\bigr)^{2p},
	\end{equation}
	for all $a,b\in\indset{N}$, with some positive constants $C_{p}$ depending only on $p$ and the model parameters.
	In particular, we have the stochastic domination bound
	\begin{equation} \label{eq:rep_G-m_prec}
		\max_{\gamma \in \indset{\gamma_0}}\bigl\lvert \bigl(\widecheck{G}^{(\gamma)}-m\bigr)_{ab} \bigr\rvert \prec N^{\xi} \sqrt{\Upsilon_{ab}}, \quad a,b\in\indset{N}.
	\end{equation}
\end{claim}

Equipped with Claim \ref{claim:rep_check}, we carry out the proof of Lemma \ref{lemma:1G_replacement}.
\begin{proof} [Proof of Lemma \ref{lemma:1G_replacement}]
Fix an integer $p \in \mathbb{N}$. For all $a,b\in\indset{N}$ and $\gamma \in \indset{\gamma_0}$, consider the random quantities $\Psi_{ab}^{\gamma, p}$, defined as 
\begin{equation} \label{eq:replace_high_moment}
	\Psi_{ab}^{\gamma, p} := \Expv_{h_{ij}} \biggl[\bigl\lvert \bigl(G^{(\gamma)}-m\bigr)_{ab} \bigr\rvert^{2p}\biggr] - \Expv_{v_{ij}} \biggl[\bigl\lvert \bigl(G^{(\gamma-1)}-m\bigr)_{ab} \bigr\rvert^{2p}\biggr],
\end{equation}
where $\Expv_h$ denotes the expectation with respect to the random variable $h$. Then, by telescoping summation, the $2p$-th moment of $(G^{(\gamma)}-m)_{ab}$ satisfies 
\begin{equation} \label{eq:telescope}
	\biggl\lvert \bigl\lVert \bigl(G^{(\gamma_0)}-m\bigr)_{ab} \bigr\rVert_{2p}^{2p} - \bigl\lVert \bigl(G^{(0)}-m\bigr)_{ab} \bigr\rVert_{2p}^{2p} \biggr\rvert \le 
	\sum_{\gamma \le \gamma_0} 	 \Expv\biggl[\bigl\lvert\Psi_{ab}^{\gamma, p}\bigr\rvert  \biggr].
\end{equation}
Therefore, our goal is to estimate the quantities $\Psi_{ab}^{\gamma, p}$ in a way that is summable over $\gamma := \phi\bigl((i,j)\bigr)$. To this end, we employ the expansions \eqref{eq:resolvent_expand}.

For a random variable $h$, let $\mom_{p,q}(h) := \Expv [ h^p (\overline{h})^q]$ denote the joint moment of $h$ and $\overline{h}$ of order $(p,q)$, and $\mom_{0,0}(h) := 1$ by convention. Then, expanding every resolvent $G^{(\gamma)}$ and $G^{(\gamma-1)}$ in \eqref{eq:replace_high_moment} using \eqref{eq:resolvent_expand}, we obtain
\begin{equation} \label{eq:moment_expansion}
	\Psi_{ab}^{\gamma, p} = \sum_{r = 1}^{2p L}  \frac{1}{(1 + \delta_{ij})^{r}} \sum_{q + q' = r} \biggl(\mom_{q,q'}(h_{ij}) - \mom_{q,q'}(v_{ij})\biggr) \other{Z}_{ab}^{\gamma,p,L}(q,q') + \Omega_{ab}^{\gamma,p}(L),
\end{equation}
where the $(q,q')$ summation runs over all $q,q' \in \mathbb{Z}_{\ge 0}$ satisfying $q+q' = r$, the residual term $\Omega_{ab}^{\gamma,p}(L)$ in \eqref{eq:moment_expansion} comprises all terms in the expansion that contain the non-checked resolvents $G^{(\gamma)}$ or $G^{(\gamma-1)}$, and $\other{Z}_{ab}^{\gamma,p}(q,q')$ comprises all terms in the expansion that contain exactly $q$ copies of $\partial_{ij}$ and $q'$ copies of $\partial_{ji}$. More precisely, the random quantities $\other{Z}_{ab}^{\gamma,p,L}(q,q')$ are given by
\begin{equation} \label{eq:Zreplace_def}
	\other{Z}_{ab}^{\gamma,p,L}(q,q') := C^{p,L}_{q,q'}\sum_{(J,J^*) \in \mathfrak{J}^{p,L}_{q,q'}} Z_{ab}^{\gamma}(J,J^*) \bigl(\widecheck{\mathcal{Y}}^{(\gamma)}\bigr)^{p - |J|} \bigl(\overline{\widecheck{\mathcal{Y}}^{(\gamma)}}\bigr)^{p - |J^*|},
	\quad \widecheck{\mathcal{Y}}^{(\gamma)} : = (\widecheck{G}^{(\gamma)}-m)_{ab}, 
\end{equation}
where $r:=q+q'$, the set $\mathfrak{J}^{p,L}_{q,q'}$ is defined as
\begin{equation}\label{def:J}
	\begin{split}
		\mathfrak{J}^{p,L}_{q,q'} &:= \biggl\{ (J, J^*) \in (\mathfrak{J}^{p,L})^2 \,:\, \sum_{\bm j \in J, J^*} j_1 = q, \,  \sum_{\bm j \in J, J^*} j_2 = q'
		\biggr\},\\
		\mathfrak{J}^{p,L} &:= \biggl\{ J \subset (\mathbb{Z}_{\ge 0})^2\backslash\{\bm 0\}\, :\, |J|\le p, \, \max_{\bm j \in J}|\bm j| \le L
		\biggr\},
	\end{split}
\end{equation}
the random quantities $Z_{ab}^{\gamma}$ are given by
\begin{equation} \label{eq:replaceZ}
	Z_{ab}^{\gamma}(J,J^*): = \bigl(S_{ij}\bigr)^{r/2} \biggl(\prod_{\bm j \in J} \partial_{ij}^{\bm j} \widecheck{\mathcal{Y}}^{(\gamma)} \biggr) \biggl(\prod_{\bm j \in J^*} \partial_{ij}^{\bm j} \overline{\widecheck{\mathcal{Y}}^{(\gamma)}} \biggr),
	\quad \gamma = \phi\bigl((i,j)\bigr),
\end{equation}
and the combinatorial constant $C^{p,L}_{q,q'}$ in \eqref{eq:Zreplace_def} depends only on $p,L,q,q'$ (its precise value is irrelevant). We recall the notation  $ \partial_{ij}^{\bm j}  =  \partial_{ij}^{j_1}   \partial_{ji}^{j_2}$ for $\bm j= (j_1, j_2)$
(the index $j$ and the exponent $\bm j$ are independent letters in \eqref{eq:replaceZ}).

We stress the similarity between \eqref{eq:Zreplace_def} and the summands on the right-hand side of 
\eqref{eq:1Gflat_offdiag_cuexp}, as well as between \eqref{eq:replaceZ}
and \eqref{eq:Ziso} and the more general \eqref{eq:Zkiso_def}, as this analogy will be our guiding principle in estimating $Z_{ab}^{\gamma,p}$. 

As in Section \ref{sec:cum_expand}, our goal is to estimate the $\other{Z}_{ab}^{\gamma,p,L}$ quantities in a way that is summable over $\gamma := \phi\bigl((i,j)\bigr)$
within the telescopic sum~\eqref{eq:telescope}.  Technically, here the summation over $\gamma$, or, equivalently, over
 $(i,j)$,  plays a very similar role of the $(c,d)$-summation in \eqref{eq:Zkiso_def} even if their origin is different.
  Each $\bm j$-derivative in \eqref{eq:replaceZ} produces a collection of decay factors that connect the original indices $(a,b)$ to the fresh summation indices $(i,j)$. These factors can be recombined into the target decay $\mathfrak{s}_1(a,b) = \sqrt{\Upsilon_{ab}}$ using the triangle inequality \eqref{eq:triag} and the convolution inequality \eqref{eq:true_convol} 
 	(in the form \eqref{eq:Schwarz_convol}). Recall, however, that to apply \eqref{eq:Schwarz_convol} effectively, we need at least four $\sqrt{\Upsilon}$-factors to connect the summation index $q$ to other indices independent of $q$
   (in our case  this summation index  will be
 one of the index pair $(i,j)$).  
  We call such summations \emph{effective}. For example, the summation
\begin{equation} \label{eq:non-effective}
	\sum_j  \sqrt{\Upsilon_{jj}\Upsilon_{jj} \Upsilon_{xj}\Upsilon_{jy}} \sim \frac{1}{\ell\eta}\sum_j  \sqrt{ \Upsilon_{xj}\Upsilon_{jy}}
\end{equation}
is \emph{not} effective; It can only be estimated using \eqref{eq:sqrt_convol}, incurring an additional $\sqrt{N/\ell}$ cost compared to \eqref{eq:Schwarz_convol}. 
Note that in Section \ref{sec:cum_expand}, for example in~\eqref{eq:Zkav_10}, the $c$-summation was automatically effective due to the presence of $A_k'$, which connected $c$ to external indices via two $\sqrt{\Upsilon}$'s, and only the effectiveness of the $d$-summation had to be ensured. In the present section, however, no such simplification is available, and we have to establish the effectiveness of both the $i$- and $j$-summations. 
We will see that in most cases it is possible; from the structure of the expansion we will find
four resolvents with independent indices. In a few exceptional case 
when it is not possible, we will carefully show that the additional $\sqrt{N/\ell}$ is affordable
by estimating one resolvent more precisely by decomposing it into its deterministic leading term and 
fluctiation.  Now we present the formal proof.

It follows from \eqref{eq:rep_G-m_prec} that for any $\xi' \in (\xi, \etaexp/20)$, we have,  
\begin{equation} \label{eq:rep_G_bound}
	\max_{\gamma \in \indset{\gamma_0}}\bigl\lvert  \widecheck{G}^{(\gamma)}_{cd} \bigr\rvert \lesssim \bigl(\delta_{cd}\sqrt{\ell\eta} + N^{\xi'}\bigr) \sqrt{\Upsilon_{cd}}, \quad \text{w.v.h.p.},
\end{equation}
uniformly in $c,d\in\indset{N}$, with the implicit constant independent of $\gamma$.
Therefore,  analogously to \eqref{eq:partial_Y}, using \eqref{eq:triag} and \eqref{eq:rep_G_bound}, we obtain,
\begin{equation} \label{eq:rep_partial_bound}
	\begin{split}
		\bigl\lvert \partial_{ij}^{\bm j} \widecheck{\mathcal{Y}}^{(\gamma)}\bigr\rvert &\lesssim C_{1} \sum_{(c,d)\in \{(i,j),(j,i)\}} \sqrt{\Upsilon_{ac}\Upsilon_{db} } \,\bigl( N^{2\xi'}  + \other{\delta}_1 N^{\xi'}\sqrt{\ell\eta} + \other{\delta}_2 \ell\eta \bigr), \quad |\bm j| =1\\
		\bigl\lvert \partial_{ij}^{\bm j} \widecheck{\mathcal{Y}}^{(\gamma)}\bigr\rvert &\lesssim  C_{|\bm j|} \sum_{c,d\in \{i,j\}} \sqrt{\Upsilon_{ac}\Upsilon_{db} } \,\bigl( N^{2\xi'}  + \other{\delta}_1 \ell\eta \bigr), \quad |\bm j| \ge 2.
	\end{split}
\end{equation}
with very high probability,   uniformly in $\gamma:= \phi\bigl((i,j)\bigr) \in \indset{\gamma_0}$ and $a,b\in\indset{N}$. 
Here, we define, similarly to \eqref{eq:other_deltas},
\begin{equation} \label{eq:other_delta_rep}
	\other{\delta}_1 \equiv \other{\delta}_1(i,j,a,b,) := \delta_{ia}+\delta_{ib} + \delta_{ja}+\delta_{jb}, \quad \other{\delta}_2 \equiv \other{\delta}_2(i,j,a,b,) := (\delta_{ia}+\delta_{ib})(\delta_{ja}+\delta_{jb}).
\end{equation} 
 The constant $C_{|\bm j|}$ depends only on the order of derivative $|\bm j|$.
Note that all but at most two resolvents are estimated trivially as $| \widecheck{G}^{(\gamma)}_{cd} |\lesssim 1$
by \eqref{eq:rep_G_bound}.  Only  two decaying factors $\Upsilon$ are kept to guarantee the necessary connection between the index pairs $(a,b)$ and $(i,j)$.

The  terms with the delta-functions in the right hand side of~\eqref{eq:rep_partial_bound} may look complicated
and we will have similar estimates in the sequel. However, at the first reading, 
 the reader is advised to disregard all the terms containing the 
 summation-restricting $\other{\delta}_{1}$ and $\other{\delta}_{2}$ factors, and focus on the terms without  $\other{\delta}$'s,
 representing the generic situation. It turns out that 
  in the critical small-$\eta$ regime, the entropic gain from a restricted summation 
  vastly outweighs the loss factors of $(\ell\eta)$-powers
 associated with $\other{\delta}$'s. To maintain the integrity of the proof, however, we  still meticulously
 track all the terms involving $\other{\delta}$'s and show that their contribution is always negligible.

Furthermore, recalling that $S_{ij} \lesssim (\Upsilon_0)_{ij}$ by \eqref{eq:SUps_comvol}, the bounds \eqref{eq:rep_partial_bound} together with \eqref{eq:triag} imply a simpler estimate 
\begin{equation} \label{eq:rep_partial_bound2}
	\bigl(S_{ij}\bigr)^{|\bm j|/2}\bigl\lvert \partial_{ij}^{\bm j} \widecheck{\mathcal{Y}}^{(\gamma)}\bigr\rvert \lesssim  C_{|\bm j|} \frac{1}{W^{|\bm j|/2}} \sqrt{\Upsilon_{ab}}    \biggl( \frac{N^{2\xi'}}{\sqrt{\ell\eta}} +   \other{\delta}_1 N^{\xi'} +  \bigl(\other{\delta}_2  + \other{\delta}_1 \mathds{1}_{|\bm j|\ge 2} \bigr)\sqrt{\ell\eta}
	\biggr),
\end{equation}
with very high probability, uniformly in $\gamma := \phi\bigl((i,j)\bigr) \in \indset{\gamma_0}$ and $a,b \in \indset{N}$.
To preserve a full power of $S_{ij}$ and make the summation in $j$ effective, we use \eqref{eq:rep_partial_bound} to estimate two (or one if $|J|+|J^*| = 1$) $\bm j$-derivatives in \eqref{eq:replaceZ},  and use \eqref{eq:rep_partial_bound2} to estimate the remaining $\bm j$-derivatives, obtaining
\begin{equation} \label{eq:repZ_bound}
	\bigl\lvert Z_{ab}^{\gamma}(J,J^*) \bigr\rvert 
	\lesssim  C_{p}\,\psi_{ab}^{(\gamma)}(J,J^*),
\end{equation}
for some constant $C_{p} > 0$, where 
$$
n := |J|+|J^*|, \quad r := \sumJ,
$$
 and the deterministic control parameter $\psi_{ab}^{(\gamma)}(J,J^*)$ is defined as 
\begin{equation} \label{eq:psi1_rep_gamma}
	\psi_{ab}^{(\gamma)}(J,J^*):=  \frac{1}{W^{(r-2)/2}} S_{ij}  \sqrt{\Upsilon_{aj}\Upsilon_{jb}} \bigl( N^{2\xi'}  + \other{\delta}_1 \ell\eta \bigr), \quad \text{if}~n = |J|+|J^*| =1,
\end{equation}
while for $ n \ge 2$,  we define it as
\begin{equation} \label{eq:psi_rep_gamma}
		\psi_{ab}^{(\gamma)}(J,J^*):= \frac{\ell\eta}{W^{\frac{r-2}{2}}} S_{ij} 
		\bigl(\Upsilon_{ai}\Upsilon_{jb} + \Upsilon_{aj}\Upsilon_{ib}\bigr)  \bigl(\Upsilon_{ab}\bigr)^{\frac{n-2}{2}} \times 
		  \biggl( \frac{N^{2n\xi'}}{(\ell\eta)^{\frac{n}{2}}} + \other{\delta}_1 (\ell\eta)^{\frac{n-n_1}{2}} N^{n_1\xi'}  + \other{\delta}_2(\ell\eta)^{\frac{n}{2}}\biggr).
\end{equation}
Here we denote 
$$
n_1 := |\{ \bm j \in J  \,: \, |\bm j| = 1 \}| + |\{ \bm j \in J^* \,: \, |\bm j| = 1 \}|.
$$ 
We have separated the case $n=1$ because, as we will see in the sequel, it is the only special case
that results in one non-effective summation and needs more attention since~\eqref{eq:sqrt_convol} will be used.

Plugging \eqref{eq:repZ_bound} into \eqref{eq:Zreplace_def}, taking the expectation of both sides, and summing over $\gamma \in \indset{\gamma_0}$ yields
\begin{equation} \label{eq:otherZ_bound}
	\begin{split}
		\sum_{\gamma \le \gamma_0} \Expv \biggl[\bigl\lvert\other{Z}_{ab}^{\gamma,p,L}(q,q')\bigr\rvert\biggr] &\lesssim C_p \sum_{\gamma\le \gamma_0}\sum_{(J,J^*) \in \mathfrak{J}^{p,L}_{q,q'}} \psi_{ab}^{(\gamma)}(J,J^*)
		\bigl\lVert\widecheck{\mathcal{Y}}^{(\gamma)}  \bigr\rVert_{2p-n}^{2p-n}\\
		&\lesssim C_p  \sum_{(J,J^*)  \in \mathfrak{J}^{p,L}_{q,q'}} 
		\max_{\gamma'\le\gamma_0}\bigl\lVert\widecheck{\mathcal{Y}}^{(\gamma')}  \bigr\rVert_{2p-n}^{2p-n} \sum_{\gamma\le \gamma(N)} \psi_{ab}^{(\gamma)}(J,J^*),
	\end{split}
\end{equation}
with very high probability.  Therefore, it remains to bound the sum of $\psi_{ab}^{(\gamma)}(J,J^*)$. 

First, consider the case $n = 1$, then, summing the definition \eqref{eq:psi1_rep_gamma} over all $\gamma := \phi\bigl((i,j)\bigr)$ we obtain, for all $a,b \in \indset{N}$ and $J,J^* \subset (\mathbb{Z}_{\ge 0})^2\backslash\{\bm 0\}$ satisfying $|J|+|J^*| = 1$,
\begin{equation} \label{eq:sum_psi1_rep_bound}
	\begin{split}
		\sum_{\gamma\le \gamma(N)} \psi_{ab}^{(\gamma)}(J,J^*) \le&~ \sum_{i,j}\frac{1}{W^{\frac{r-2}{2}}} S_{ij}  \sqrt{\Upsilon_{aj}\Upsilon_{jb}} \bigl( N^{2\xi'}  + \other{\delta}_1 \ell\eta \bigr) 
		\lesssim \sqrt{\Upsilon_{ab}}  \frac{N^{2\xi'}\sqrt{N\eta}}{W^{\frac{r-2}{2}}\eta},
	\end{split}
\end{equation}
where we where we used \eqref{eq:other_delta_rep}, invariance of $\psi_{ab}^{(\gamma)}(J,J^*)$ under swapping $i$ and $j$, and the estimate \eqref{eq:Schwarz_convol} and \eqref{eq:sqrt_convol} together with the first bound in \eqref{eq:SUps_comvol} to sum the decaying factors in $i$ and $j$. 

Similarly, for  $J,J^* \subset (\mathbb{Z}_{\ge 0})^2\backslash\{\bm 0\}$ satisfying $|J|+|J^*| \ge 2$, summing the definition \eqref{eq:psi_rep_gamma} over all $\gamma := \phi\bigl((i,j)\bigr)$, we obtain
\begin{equation} \label{eq:sum_psi_rep_bound}
	\begin{split}
		\sum_{\gamma\le \gamma(N)} \psi_{ab}^{(\gamma)}(J,J^*) \lesssim&~ \frac{\ell\eta}{W^{\frac{r-2}{2}}} \bigl(\Upsilon_{ab}\bigr)^{\frac{n-2}{2}} \sum_{i,j}S_{ij} 
		\Upsilon_{ai}\Upsilon_{jb}
		\biggl( \frac{N^{2n\xi'}}{(\ell\eta)^{\frac{n}{2}}} + \other{\delta}_1 (\ell\eta)^{\frac{n-n_1}{2}} N^{n_1\xi'}  + \other{\delta}_2(\ell\eta)^{\frac{n}{2}}\biggr)\\
		\lesssim&~ \sqrt{\Upsilon_{ab}}\biggl(\frac{N^{2n\xi'}}{W^{\frac{r-2}{2}}(\ell\eta)^{\frac{n-2}{2}}\eta} + \frac{N^{n_1\xi'}(\ell\eta)^{\frac{n-n_1}{2}}}{W^{\frac{r-2}{2}}} + \frac{1}{W^{\frac{r-n}{2}}}\biggr),
	\end{split}
\end{equation}
where we used \eqref{eq:other_delta_rep}, and \eqref{eq:true_convol} to perform the $i$ and $j$ summations.

It follows from  the moment matching condition \eqref{eq:moment_match}, that all terms with $r \in \{0,1,2,3\}$ in the expansion \eqref{eq:moment_expansion} vanish identically, while the terms with $r=4,5$ come with a prefactor of $\mathcal{O}(\lambda)$ by \eqref{eq:m4_cond}. 
By a straightforward case-by-case analysis for the possible combinations of $n$ and $n_1$ given a fixed $r\ge 4$, we deduce that the scalar factors in \eqref{eq:sum_psi1_rep_bound}--\eqref{eq:sum_psi_rep_bound} satisfy
\begin{equation} \label{eq:other_Z_r_bound}
	 \frac{N^{2\xi'}\sqrt{N\eta}}{W^{\frac{r-2}{2}}\eta} +  \frac{N^{2n\xi'}}{W^{\frac{r-2}{2}}(\ell\eta)^{\frac{n-2}{2}}\eta} + \frac{N^{n_1\xi'}(\ell\eta)^{\frac{n-n_1}{2}}}{W^{\frac{r-2}{2}}} + \frac{1}{W^{\frac{r-n}{2}}} 
	\lesssim  1 + \frac{N^{2\xi'}\sqrt{N\eta}}{W^{\frac{r-2}{2}}\eta},
\end{equation}
for all $(J,J^*) \in \mathfrak{J}^{p,L}_{q,q'}$ with $q+q' = r \ge 4$, uniformly in $p \in \mathbb{N}$. 
The dominant term involving $\sqrt{N\eta}$, that arose from bounding the ineffective 
summations via \eqref{eq:sqrt_convol},   comes from the special $n=1$ case.
Also note 
that in the estimate~\eqref{eq:sum_psi_rep_bound}  for the $n\ge 2$ cases the largest term is the first one,
coming from the generic index distributions; the other two terms involving $\delta$-functions give smaller contributions.

Hence, from  \eqref{eq:moment_expansion}, \eqref{eq:otherZ_bound}--\eqref{eq:other_Z_r_bound}, we obtain
\begin{equation} \label{eq:psi_rep_bound1}
	 \sum_{\gamma \le \gamma_0} 	 \Expv\biggl[\bigl\lvert\Psi_{ab}^{\gamma, p}\bigr\rvert  \biggr]  \lesssim  C_p \biggl(1 +   \lambda  \frac{N^{2\xi'}\sqrt{N\eta}}{W\eta}  \biggr)  \sum_{n=1}^{2p} \bigl(\Upsilon_{ab}\bigr)^{n/2} \max_{\gamma'\le\gamma_0}\bigl\lVert\widecheck{\mathcal{Y}}^{(\gamma')}  \bigr\rVert_{2p-n}^{2p-n} + \sum_{\gamma \le \gamma_0}\Expv\biggl[\bigl\lvert\Omega_{ab}^{\gamma,p}(L)\bigr\rvert\biggr],
\end{equation}
where we used that $N^{2\xi'}\sqrt{N\eta}/(W^2\eta) \lesssim 1$, by \eqref{eq:WN},   $\eta \ge N^{-1+\etaexp}$ and $\xi' < \etaexp/20$. Moreover,  $\lambda  N^{2\xi'}\sqrt{N\eta}/(W\eta) \lesssim 1$ by the choice of $\lambda$ in \eqref{eq:lambda_assume},
hence the scalar factor in the big brackets on the right-hand side of \eqref{eq:psi_rep_bound1} is bounded by a constant. 
The overestimate  by a factor $\sqrt{N/\ell}$, incurred due to~\eqref{eq:sqrt_convol}, remains affordable thanks to an additional factor $\lambda$, which is small in the critical small $\eta$ regime. The presence of $\lambda$ reflects that while the moments of $h_{ij}$ and $v_{ij}$ beyond the third order do not  match exactly, they are still sufficiently close for our applications.  This is because the zag-step only removes the Gaussian component introduced in the preceding zig-step, and crucially, for smaller  $\eta$ this added Gaussian component is itself small.

We now estimate the error term $\Omega_{ab}^{\gamma,p}(L)$.
The non-checked resolvent  are estimated trivially by operator norm
$\norm{G^{(\gamma)}} \le \eta^{-1}$ while for all other resolvents we can use 
$\lVert\widecheck{G}^{(\gamma)}\rVert_{\max}\lesssim 1$, $\gamma\le \gamma_0$,  with very high probability from  \eqref{eq:rep_G_bound}.
  By choosing the expansion order $L \ge 4D' + 12$, the term containing the original resolvent
  $G^{(\gamma)}$ on the right-hand side of \eqref{eq:resolvent_expand} admits the bound
\begin{equation}
	\bigl\lVert\bigl(-\widecheck{G}^{(\gamma)} \Delta^{(\gamma)}_H\bigr)^{L+1} G^{(\gamma)}\bigr\rVert_{\max} \lesssim \bigl(|h_{ij}| + |v_{ij}|\bigr)^{L+1} W^{-(L+1)/2}\eta^{-1} \lesssim \bigl(|h_{ij}| + |v_{ij}|\bigr)^{L+1} N^{-3}\sqrt{\Upsilon_{ab}},
\end{equation}
with very high probability, 
where we used \eqref{eq:S_bound} to estimate $S_{ij}$. Therefore, it is straightforward to check that
\begin{equation} \label{eq:rep_Omega_est}
	\begin{split}
		\Expv\biggl[ \bigl\lvert \Omega_{ab}^{\gamma,p}(L) \bigr\rvert\biggr] 
		&\lesssim C_p \sum_{q=1}^{2p} \bigl(N^{-3}\sqrt{\Upsilon_{ab}}\bigr)^q \sum_{J \in \mathfrak{J}^{2p-q}} \Expv\biggl[\bigl\lvert Z^{\gamma}_{ab}(J,\emptyset) \bigr\rvert \bigl\lvert\widecheck{\mathcal{Y}}^{(\gamma)} \bigr\rvert^{2p-n-q}\biggr]\\
		&\lesssim C_p N^{-3}\sum_{q=1}^{2p} \sum_{n=1}^{2p-q}(\Upsilon_{ab})^{(q+n)/2}   \bigl\lVert\widecheck{\mathcal{Y}}^{(\gamma)} \bigr\rVert_{2p-n-q}^{2p-n-q},
	\end{split}
\end{equation} 
uniformly in $\gamma \in \indset{\gamma_0}$, where in the second step we used \eqref{eq:rep_partial_bound2} to bound $\bigl\lvert Z^{\gamma}_{ab}(J,\emptyset) \bigr\rvert$ by $C(\Upsilon_{ab})^{|J|/2}$ with very high probability. In particular, the  sum  involving $\Omega_{ab}^{\gamma,p}(L)$ in \eqref{eq:psi_rep_bound1} is negligible.

By a standard argument relying on Young's inequality, we conclude from \eqref{eq:telescope} and  \eqref{eq:psi_rep_bound1} that
\begin{equation} \label{eq:1G_rep}
	\begin{split}
		\bigl\lVert \bigl(G^{(\gamma_0)}-m\bigr)_{ab} \bigr\rVert_{2p}^{2p} &\le  \bigl\lVert \bigl(G^{(0)}-m\bigr)_{ab} \bigr\rVert_{2p}^{2p}  + \mathcal{O}_p(1) \bigl(N^{\xi}   \sqrt{\Upsilon_{ab}} \bigr)^{2p} + N^{-\xi} \max_{\gamma'\le\gamma_0}\bigl\lVert\widecheck{\mathcal{Y}}^{(\gamma')}  \bigr\rVert_{2p}^{2p}\\
		&\le \bigl\lVert \bigl(G^{(0)}-m\bigr)_{ab} \bigr\rVert_{2p}^{2p} +  \other{C}_p\bigl(N^{\xi}\sqrt{\Upsilon_{ab}}\bigr)^{2p} +   2N^{-\xi}\max_{\gamma' \in \indset{\gamma_0}}	\bigl\lVert\mathcal{Y}^{(\gamma'-1)}\bigr\rVert_{2p}^{2p}.
	\end{split}
\end{equation}
where we used Claim \ref{claim:rep_check} to bound $\lVert\widecheck{\mathcal{Y}}^{(\gamma')} \rVert_{2p}^{2p}$
in the last step. Hence, \eqref{eq:1G_rep_goal} holds by \eqref{eq:1G_rep}, \eqref{eq:1G_GFT_init}, and \eqref{eq:1G_rep_assume}. This concludes the proof of Lemma \ref{lemma:1G_replacement}.
\end{proof}

We close this section by proving Claim \ref{claim:rep_check}.
The proof of Claim \ref{claim:rep_check} uses the same ideas as that of Lemma \ref{lemma:1G_replacement}, but in a much cruder form since we do not require the error term to be 
summable over $\gamma$. Therefore, we present the proof only briefly.
\begin{proof}[Proof of Claim \ref{claim:rep_check}] 
	Completely analogously to \eqref{eq:resolvent_expand}, the resolvent $\widecheck{G}^{(\gamma)}$ can be expanded around $G^{(\gamma-1)}$, that is
	\begin{equation} \label{eq:Gcheck_expand}
		\widecheck{G}^{(\gamma)} = \sum_{q=0}^{L-1} \bigl(G^{(\gamma-1)}\Delta^{(\gamma)}_V\bigr)^q G^{(\gamma-1)} + \bigl(S_{ij}\bigr)^{(L+1)/2}\bigl(G^{(\gamma-1)}\Delta^{(\gamma)}_V\bigr)^L \widecheck{G}^{(\gamma)}, \quad \gamma = \phi\bigl((i,j)\bigr)
	\end{equation}
	For all $a,b\in\indset{N}$ and  $\gamma \in \indset{\gamma_0}$, define, similarly to \eqref{eq:replace_high_moment},
	\begin{equation} \label{eq:Psi_replace_check}
		\widecheck{\Psi}_{ab}^{\gamma, p} :=  \bigl\lvert \bigl(\widecheck{G}^{(\gamma)}-m\bigr)_{ab} \bigr\rvert^{2p}  - \Expv_{v_{ij}} \biggl[\bigl\lvert \bigl(G^{(\gamma-1)}-m\bigr)_{ab} \bigr\rvert^{2p}\biggr].
	\end{equation}
	To establish \eqref{eq:1G_rep_check}, it suffices to estimate $\widecheck{\Psi}_{ab}^{\gamma, p}$ uniformly, without summation in $\gamma$. 
	
	Using \eqref{eq:Gcheck_expand} for every resolvent $\widecheck{G}^{(\gamma)}$  in \eqref{eq:Psi_replace_check}, we obtain, similarly to \eqref{eq:moment_expansion} and \eqref{eq:Zreplace_def},
	\begin{equation} \label{eq:check_moment_expansion}
		\widecheck{\Psi}_{ab}^{\gamma, p} = -\sum_{r = 2}^{2p L}  \sum_{q + q' = r}  \frac{C^{p,L}_{q,q'}\mom_{q,q'}(v_{ij})}{(1 + \delta_{ij})^{r}}\sum_{(J,J^*) \in \mathfrak{J}^{p,L}_{q,q'}} \widecheck{Z}_{ab}^{\gamma}(J,J^*) \bigl( \mathcal{Y}^{(\gamma-1)}\bigr)^{p - |J|} \bigl(\overline{\mathcal{Y}^{(\gamma-1)} }\bigr)^{p - |J^*|} + \widecheck{\Omega}_{ab}^{\gamma,p}(L),
	\end{equation}
	where $\mathcal{Y}^{(\gamma-1)} := (G^{(\gamma-1)}-m)_{ab}$, the quantity $\widecheck{\Omega}_{ab}^{\gamma,p}(L)$ comprises all terms in the expansion containing $\widecheck{G}^{(\gamma)}$, and the quantities $\widecheck{Z}_{ab}^{\gamma}(J,J^*)$ are defined analogously to \eqref{eq:replaceZ}, with $\widecheck{\mathcal{Y}}^{(\gamma)}$ replaced by $\mathcal{Y}^{(\gamma-1)}$. Here we used that $\mom_{1,0}(v_{ij}) = \mom_{0,1}(v_{ij}) = 0$.
	
	Proceeding exactly as in the proof of \eqref{eq:repZ_bound} above, using \eqref{eq:1G_rep_assume} to bound $G^{(\gamma-1)}_{xy}$ and obtaining an analog of \eqref{eq:rep_partial_bound} for the derivatives of $\mathcal{Y}^{(\gamma-1)}$, we deduce that
	\begin{equation} \label{eq:rep_Z_crude}
		\bigl\lvert\widecheck{Z}_{ab}^{\gamma}(J,J^*)\bigr\rvert  \lesssim C_{p}  \bigl(\Upsilon_{ab}\bigr)^{\frac{|J| +|J^*|}{2}}.
	\end{equation}
	The error term $\widecheck{\Omega}_{ab}^{\gamma,p}(L)$ is estimated the same way as its non-checked counterpart $\Omega_{ab}^{\gamma,p}(L)$ in \eqref{eq:rep_Omega_est}, but with $\widecheck{\mathcal{Y}}^{(\gamma)}$ replaced by $\mathcal{Y}^{(\gamma-1)}$, and hence is negligible. Combining \eqref{eq:check_moment_expansion} and  \eqref{eq:rep_Z_crude}, we obtain
	\begin{equation}
		\Expv\biggl[\bigl\lvert \widecheck{\Psi}_{ab}^{\gamma, p} \bigr\rvert\biggr] \lesssim C_p \sum_{n=1}^{2p}   \bigl(\Upsilon_{ab}\bigr)^{n/2}  \bigl\lVert \mathcal{Y}^{(\gamma-1)}\bigr\rVert_{2p - n}^{2p - n} \lesssim C_p\bigl(N^{\xi}\sqrt{\Upsilon_{ab}}\bigr)^{2p}+N^{-\xi} \bigl\lVert \mathcal{Y}^{(\gamma-1)}\bigr\rVert_{2p}^{2p},
	\end{equation}
	uniformly in $\gamma \in \indset{\gamma_0}$, from which the desired \eqref{eq:1G_rep_check} immediately follows by definition of $\widecheck{\Psi}_{ab}^{\gamma, p}$ in \eqref{eq:Psi_replace_check}.
	
	Since \eqref{eq:1G_rep_check} holds for any $p\in\mathbb{N}$, \eqref{eq:rep_G-m_prec} follows straightforwardly by the relation of stochastic domination and high-moment bounds from~\eqref{eq:1G_rep_assume}.
	 This concludes the proof of Claim \ref{claim:rep_check}.
\end{proof}

\subsection{Isotropic chains of length $k\ge2$} \label{sec:k_iso_rep}
To prove Proposition \ref{prop:zag} for longer chains, we proceed by strong induction in length $k$.
As in the $k=1$ case, we first complete the proof for isotropic chains of any length $k\le K$ in the current
section. The proof for averaged chains will be performed in Section~\ref{sec:k_av_rep} relying directly
on the final result of the isotropic induction. No separate induction on $k$ for
the averaged proof is needed.  The final structure of the proof of Proposition~\ref{prop:zag} 
will be explained in Section~\ref{sec:zag_proof_for_real}. 

  By analogy with \eqref{eq:resolvent_chains} and \eqref{eq:Gk_def}, we denote, for all $u \le k' \in \indset{k}$, $\bm z \in \{z,\overline{z}\}^k$, 
and $\bm x \in \indset{N}^k$,
\begin{equation}
	G^{(\gamma)}_{[u,k']}(x_u, \dots, x_{k'-1}) := \biggl(\prod_{q = u}^{k'-1} G^{(\gamma)}_{q} S^{x_q}\biggr)G^{(\gamma)}_{k'}, \quad G^{(\gamma)}_{q} := G^{(\gamma)}(z_q), \quad q \in \indset{k},
\end{equation}
and the chains $\widecheck{G}^{(\gamma)}_{[u,k']}(\bm x')$ are defined the same way with $ G^{(\gamma)}(z_q)$ replaced by $ \widecheck{G}^{(\gamma)}(z_q)$. Recall that $\bm x = (x_1, \ldots x_{k})$.

When proving the isotropic local law for chains of length $k \in \indset{2,\maxK}$,  we assume that for all shorter chains of length $k' \in \indset{k-1}$, we already have 
\begin{equation} \label{eq:k_rep_assume}
	\bigl\lvert \bigl((G^{(\gamma)}_{[1,k']}-M_{[1,k']})(\bm x')\bigr)_{ab} \bigr\rvert \lesssim N^{\xi'} (\ell\eta)^{\alpha_{k'}} \mathfrak{s}_{k'}^\mathrm{iso}(a,\bm x', b),  \quad a,b\in\indset{N}, \quad \bm x \in \indset{N}^{k'}, \quad \bm z \in \{z,\overline{z}\}^{k'},
\end{equation}
with very high probability,  uniformly in $\gamma \in \indset{\gamma(N)}$, for some fixed tolerance exponent $\xi'$ satisfying $\xi < \xi' < \etaexp/40$. Moreover, for $\gamma =0$, we assume that the chains of length $k$ satisfy the initial condition for the replacement process, analogous to \eqref{eq:1G_GFT_init},
\begin{equation} \label{eq:k_rep_init}
	 \norm{\bigl((G^{(0)}_{[1,k]}-M_{[1,k]})(\bm x')\bigr)_{ab} }_p \lesssim C_p N^{\xi} (\ell\eta)^{\alpha_{k}} \mathfrak{s}_{k}^\mathrm{iso}(a,\bm x', b),  
\end{equation} 
uniformly in $a,b\in\indset{N}$, $\bm x \in \indset{N}^{k}$, and $\bm z \in \{z,\overline{z}\}^{k}$. 
The outcome of the strong induction in $k$ is that  \eqref{eq:k_rep_assume} and \eqref{eq:k_rep_init} imply
  \eqref{eq:k_rep_assume}  for all $k' \in \indset{k}$, in particular  \eqref{eq:k_rep_assume} holds with the same
  tolerance exponent $\xi$ even for $k'=k$. 

To prove each step of the overall induction in $k \in \indset{2,\maxK}$, we perform the one-by-one replacement of the resolvent entries via a nested induction in $\gamma \in \indset{\gamma(N)}$. We formulate
one step of the nested $\gamma$-induction for each
fixed  $k$  in Lemma~\ref{lemma:kG_replacement} below. These elementary steps  
are completely analogous  to the $k=1$ case,
formulated in Lemma \ref{lemma:1G_replacement} and proven
 in Section \ref{sec:1Grep} above.

We mention that in the sequel, all estimate on resolvent chains 
$(G^{(\gamma)}_{[\dots]})_{ab}$ and $(G^{(\gamma)}_{[\dots]}-M_{[\dots]})_{ab}$ are understood to be uniform 
in the entries $a,b\in\indset{N}$, external indices $\bm x \in \indset{N}^{k}$, and the spectral parameters $\bm z \in \{z,\overline{z}\}^{k}$ without repeating this uniformity all the times.

\begin{lemma}[$k$-Chain Replacement Lemma] \label{lemma:kG_replacement}
	Assume that \eqref{eq:k_rep_assume} and \eqref{eq:k_rep_init} hold for some $0 < \xi < \xi' < \etaexp/40$.
	There exist constants $\{\other{C}_{p,k}\}_{p\in\mathbb{N}}$, depending only on $p$, $k$, and the model parameters, such that the following holds.
	
	Let $\gamma_0 \in \indset{\gamma(N)}$, and  let $\{\mu_p\}_{p\in\mathbb{N}}$, independent of $\gamma_0$,
	be a set of positive constants, such that the resolvents $G^{(\gamma)}$ with $\gamma \le \gamma_0 - 1$ satisfy
	\begin{equation} \label{eq:kG_rep_moment_assume}
		\max_{\gamma \in \indset{\gamma_0-1}}	\norm{\bigl((G^{(\gamma)}_{[1,k]}-M_{[1,k]})(\bm x')\bigr)_{ab} }_{2p}^{2p} \le \mu_p \bigl(N^{\xi} (\ell\eta)^{\alpha_k} \mathfrak{s}_{k}^\mathrm{iso}(a,\bm x', b)\bigr)^{2p},
	\end{equation} 
	for any integer $p \in\mathbb{N}$.
	Then, the resolvent $G^{(\gamma_0)}$ satisfies
	\begin{equation} \label{eq:kG_rep_goal}
		\norm{\bigl((G^{(\gamma_0)}_{[1,k]}-M_{[1,k]})(\bm x')\bigr)_{ab} }_{2p}^{2p} \le \bigl(\other{C}_{p,k}+ 2N^{-\xi}\mu_p\bigr) \bigl(N^{\xi} (\ell\eta)^{\alpha_k} \mathfrak{s}_{k}^\mathrm{iso}(a,\bm x', b)\bigr)^{2p},
	\end{equation}
	for any integer $p \in \mathbb{N}$.
\end{lemma}
Note that the tolerance exponent $\xi$
does not deteriorate along the $k$-induction (assuming the initial condition \eqref{eq:k_rep_assume} holds with the
same $\xi$), nor does it  
not deteriorate along the $\gamma$-induction from~\eqref{eq:kG_rep_moment_assume} to
\eqref{eq:kG_rep_goal}.

The proof of Lemma \ref{lemma:kG_replacement} follows exactly the same strategy as that of Lemma \ref{lemma:1G_replacement}. We therefore outline only the key differences and leave the minor details to the reader. 
\begin{proof}[Proof of Lemma \ref{lemma:kG_replacement}]
	We fix the vector of spectral parameters $\bm z \in \{z,\overline{z}\}^k$, and suppress the dependence of all quantities of $\bm z$.
	The basic goal is to estimate the difference of high moments, defined, analogously to \eqref{eq:replace_high_moment}, as
	\begin{equation} \label{eq:kG_high_moment}
		\Psi_{k}^{\gamma, p} \equiv \Psi_{k}^{\gamma, p}(a,\bm x', b) := \Expv_{h_{ij}} \biggl[\bigl\lvert \bigl((G^{(\gamma)}_{[1,k]}-M_{[1,k]})(\bm x')\bigr)_{ab} \bigr\rvert^{2p}\biggr] - \Expv_{v_{ij}} \biggl[\bigl\lvert \bigl((G^{(\gamma-1)}_{[1,k]}-M_{[1,k]})(\bm x')\bigr)_{ab} \bigr\rvert^{2p}\biggr],
	\end{equation}
	for all $\gamma\in\indset{\gamma_0}$. The estimate on $\Psi_{k}^{\gamma, p}$ must be summable in $\gamma$, so that we can bound the telescoping summation of~\eqref{eq:kG_high_moment} from $\gamma=1$ up to $\gamma_0$,
	 analogously to the $k=1$  case  in~\eqref{eq:telescope}. 
	
	Using~\eqref{eq:resolvent_expand}, we expand every resolvent $G_q^{(\gamma)}$ and $G_q^{(\gamma-1)}$ in \eqref{eq:kG_high_moment} obtaining an expression for $\Psi_{k}^{\gamma, p}$ similar to \eqref{eq:moment_expansion},
	\begin{equation} \label{eq:Gk_moment_expansion}
		\Psi_{k}^{\gamma, p} = \sum_{r = 1}^{2p L}  \frac{1}{(1 + \delta_{ij})^{r}} \sum_{q + q' = r} \biggl(\mom_{q,q'}(h_{ij}) - \mom_{q,q'}(v_{ij})\biggr) \other{Z}_{k}^{\gamma,p,L}(q,q') + \Omega_{k}^{\gamma,p}(L),
	\end{equation}
	where $\other{Z}_{k}^{\gamma,p,L}(q,q')$ comprises all terms in the expansion that contain exactly $q$ copies of $\partial_{ij}$ and $q'$ copies of $\partial_{ji}$, while the residual term $\Omega_{k}^{\gamma,p}(L)$ comprises all terms in the expansion that contain the non-checked resolvents $G^{(\gamma)}$ or $G^{(\gamma-1)}$. Similarly to \eqref{eq:rep_Omega_est}, it is easy to show that $\Omega_{k}^{\gamma,p}(L)$ is a negligible error term. 
	
	It is straightforward to check, starting with \eqref{eq:resolvent_expand} and \eqref{eq:kG_high_moment}, that the quantities $\other{Z}_{k}^{\gamma,p,L}(q,q')$ in \eqref{eq:Gk_moment_expansion} admit the bound
	\begin{equation} \label{eq:Zk_replace_def}
		\bigl\lvert \other{Z}_{k}^{\gamma,p,L}(q,q') \bigr\rvert  \lesssim \other{C}^{p,L}_{q,q'} \sum_{(J,J^*) \in 
		\mathfrak{J}^{p,L}_{q,q'}} \bigl(S_{ij}\bigr)^{(q+q')/2} \biggl\lvert \prod_{\bm j \in J} \partial_{ij}^{\bm j} \widecheck{\mathcal{Y}}^{(\gamma)}_k \biggr\vert \biggl\lvert\prod_{\bm j \in J^*} \partial_{ij}^{\bm j} \overline{\widecheck{\mathcal{Y}}_k^{(\gamma)}} \biggr\rvert   \bigl\lvert\bigl( \widecheck{\mathcal{Y}}_k^{(\gamma)}\bigr)\bigr\rvert^{2p - |J| - |J^*|},
	\end{equation}
	where $\widecheck{\mathcal{Y}}_k^{(\gamma)} : = \bigl((\widecheck{G}^{(\gamma)}_{[1,k]}-M_{[1,k]})(\bm x')\bigr)_{ab}$
	and we recall the definition of $\mathfrak{J}^{p,L}_{q,q'}$ from \eqref{def:J}. 
	We stress the similarity between \eqref{eq:Zk_replace_def} and \eqref{eq:Zkiso_def}
	(except for the summation over $c,d$ in \eqref{eq:Zkiso_def}; whose role will ultimately be played by the final telescopic summation of $\Psi_k^{\gamma, p}$ from~\eqref{eq:kG_high_moment} over $\gamma$).
	Although the basic object $\widecheck{\mathcal{Y}}_k^{(\gamma)}$  in~\eqref{eq:Zk_replace_def} and $\mathcal{Y}_k$  in~\eqref{eq:Zkiso_def} are somewhat different, 
	they both represent long resolvent chains. And while the corresponding of differential operators $\partial_{ij}$ also slightly differ, they yield structurally similar expressions
	in terms of products of isotropic chains.  Most importantly, the spatial decay structure---encoded 
	in the off-diagonal elements of resolvent chains---is identical in both cases, and it fully connects the entry indices $a, b$
	with the derivative indices $i,j$ and the external indices $\bm x'$. In particular,  the spatial structure of
	the $\partial_{ij}^{\bm j}$-derivatives on the right-hand side of \eqref{eq:Zk_replace_def} are estimated analogously to \eqref{eq:partial_Y}. More precisely, using \eqref{eq:M_bound}, \eqref{eq:k_rep_assume} and \eqref{eq:kG_rep_moment_assume}, we obtain, for all $\bm j \in (\mathbb{Z}_{\ge 0})^2\backslash\{\bm 0\}$,
	\begin{equation} \label{eq:partial_k_check}
		\begin{split}
			\bigl\lvert\partial_{ij}^{\bm j} \widecheck{G}^{(\gamma)}_{[1,k]}(\bm x')\bigr\rvert \prec&~ 
			\sum_{1\le a_1\le\dots\le a_{|\bm j|}\le k} \sum_{\{(c_l,d_l)\}} \mathfrak{s}_{0}(a,  d_1) \mathfrak{s}_{|\bm j|}(c_{|\bm j|},  b)  \prod_{l=1}^{|\bm j|-1}\mathfrak{s}_{l}(c_l, d_{l+1}\bigr)  \\
			&\times (\ell\eta)^{|\bm j|/2 + \alpha_k} \biggl(\frac{N^{2\xi}}{\sqrt{\ell\eta}}  + \other{\delta}_1 N^{\xi} + \bigl(\other{\delta}_2  + \other{\delta}_1 \mathds{1}_{|\bm j|\ge 2} \bigr)\sqrt{\ell\eta} \biggr),
		\end{split}		
	\end{equation}   
	uniformly in $\gamma$,  where $\mathfrak{s}_l(\,\cdot\,,\,\cdot\,)$ are the decay factors corresponding to resolvent chains $\widecheck{G}^{(\gamma)}_{[a_{l},a_{l+1}]}$ for all $l \in \indset{0,|\bm j|}$, under the convention $a_0 := 1$, $a_{|\bm j| + 1}:= k $.
	Using the first bound from  \eqref{eq:SUps_comvol} and \eqref{eq:triag} repeatedly, we deduce that for $|\bm j| \ge 2$, 
	\begin{equation} \label{eq:rep_concat}
		\begin{split}
			\bigl(S_{ij}\bigr)^{|\bm j|/2-1}\prod_{l=1}^{|\bm j|-1}\mathfrak{s}_{l}(c_l, d_{l+1}\bigr)  &\lesssim \frac{1}{W^{|\bm j|/2-1}} \prod_{l=1}^{|\bm j|-2}\mathfrak{s}_{l}(c_1, c_1\bigr) \times \mathfrak{s}_{|\bm j|-1}(c_1, d_{|\bm j|}\bigr)\\
			&\lesssim \frac{1}{(W\ell\eta)^{|\bm j|/2-1}}\,\mathfrak{s}_{a_{|\bm j|} - a_1 + 1}^\mathrm{iso}\bigl(c_1, x^{a_1}, \dots, x^{a_{|\bm j|}-1} ,d_{|\bm j|}\bigr).
		\end{split}
	\end{equation}
	Plugging \eqref{eq:rep_concat} into \eqref{eq:partial_k_check} and using \eqref{eq:triag},   we obtain the corresponding 
	$\gamma$-uniform bound
	\begin{equation}  \label{eq:partial_k_check2}
			\bigl(S_{ij}\bigr)^{|\bm j|/2}\bigl\lvert\partial_{ij}^{\bm j} \widecheck{G}^{(\gamma)}_{[1,k]}(\bm x')\bigr\rvert \prec  
			  (\ell\eta)^{\alpha_k} \mathfrak{s}_{k}^\mathrm{iso} (a, \bm x', b )   \frac{1}{W^{|\bm j|/2}}\biggl(\frac{N^{2\xi}}{\sqrt{\ell\eta}}  + \other{\delta}_1 N^{\xi}   + \bigl(\other{\delta}_2  + \other{\delta}_1 \mathds{1}_{|\bm j|\ge 2} \bigr)\sqrt{\ell\eta}\biggr).
	\end{equation}
	We encourage the reader to compare equations~\eqref{eq:partial_k_check}--\eqref{eq:partial_k_check2} to the corresponding bounds~\eqref{eq:rep_partial_bound}--\eqref{eq:rep_partial_bound2} established for the case $k=1$.
	 the difference is only that $\sqrt{\Upsilon_{ab}} =
	\mathfrak{s}_{k=1}^\mathrm{iso} (a,  b )$ is replaced with $ \mathfrak{s}_{k}^\mathrm{iso} (a, \bm x', b )$
	and the factor $ (\ell\eta)^{\alpha_k} $ added which is absent for $k=1$ as $\alpha_1=0$. 
	 Even though the spatial component on the right-hand side of \eqref{eq:partial_k_check} is given by a complicated product of decay factors, the only important information about them is that these factors can be recombined into the target decay $\mathfrak{s}_k^\mathrm{iso}(a,\bm x', b)$.   
	It is straightforward to verify that the bound on $\other{Z}_{k}^{\gamma,p,L}(q,q')$ resulting from \eqref{eq:partial_k_check}--\eqref{eq:partial_k_check2} is summable in $\gamma := \phi\bigl((i,j)\bigr)$ exactly 
	as in the $k=1$ case, see~\eqref{eq:psi_rep_bound1}. 
	Hence, the remainder of the proof of \eqref{eq:kG_rep_goal} proceeds analogously to that of Lemma~\ref{lemma:1G_replacement}.  
This concludes the proof of Lemma~\ref{lemma:kG_replacement}.
\end{proof}

\subsection{Averaged chains of arbitrary length}\label{sec:k_av_rep}
To conclude the proof of Proposition \ref{prop:zag}, it remains to prove the estimates on the averaged chains. 
Throughout this part of the proof, we use the isotropic bounds established above as an input, together with the initial assumption that 
\begin{equation} \label{eq:k_av_rep_init}
	\norm{\Tr\bigl[(G^{(0)}_{[1,k]}-M_{[1,k]})(\bm x')S^{x_k}\bigr]}_p \lesssim C_p N^{\xi} (\ell\eta)^{\beta_{k}} \mathfrak{s}_{k}^\mathrm{av}(\bm x),  \quad \bm x \in \indset{N}^k, \quad \bm z \in \{z,\overline{z}\}^k,
\end{equation} 
for all $p\in\mathbb{N}$.

We proceed similarly to Section \ref{sec:k_iso_rep} above, with the a single step of the replacement process encompassed in the following lemma.
\begin{lemma}\label{lemma:av_rep}
	Fix $k \in \indset{\maxK}$. Assume that the isotropic local law
	\eqref{eq:k_rep_assume} with\footnote{In fact, for the proof of Lemma \ref{lemma:av_rep}, it suffices to assume that \eqref{eq:k_rep_assume} holds only for all $k'\in\indset{(k+1)\wedge\maxK}$.} $k' \in \indset{\maxK}$  and the initial bound 
	for the averaged law \eqref{eq:k_av_rep_init} hold for some $0 < \xi < \etaexp/40$.
	Then, there exist constants $\{\other{C}_{p,k}\}_{p\in\mathbb{N}}$, depending only on $p$, $k$,
	and the model parameters, such that the following holds.
	
	Let $\gamma_0 \in \indset{\gamma(N)}$, and  let $\{\mu_p\}_{p\in\mathbb{N}}$, independent of $\gamma_0$,
	be a set of positive constants, such that the resolvents $G^{(\gamma)}$ with $\gamma \le \gamma_0 - 1$ satisfy
	\begin{equation} \label{eq:kG_av_rep_moment_assume}
		\max_{\gamma \in \indset{\gamma_0-1}}	\norm{\Tr\bigl[(G^{(\gamma)}_{[1,k]}-M_{[1,k]})(\bm x')S^{x_k}\bigr]}_{2p}^{2p} \le \mu_p \bigl(N^{\xi} (\ell\eta)^{\beta_k} \mathfrak{s}_{k}^\mathrm{av}(\bm x)\bigr)^{2p},
	\end{equation}
	for any integer $p \in\mathbb{N}$.
	Then, the resolvent $G^{(\gamma_0)}$ satisfies
	\begin{equation} \label{eq:kG_av_rep_goal}
		\norm{\Tr\bigl[(G^{(\gamma_0)}_{[1,k]}-M_{[1,k]})(\bm x')S^{x_k}\bigr]}_{2p}^{2p} \le \bigl(\other{C}_{p,k}+ 2N^{-\xi/4}\mu_p\bigr) \bigl(N^{\xi} (\ell\eta)^{\beta_k} \mathfrak{s}_{k}^\mathrm{av}(\bm x)\bigr)^{2p},
	\end{equation}
	for any integer $p \in \mathbb{N}$.
\end{lemma}
We note that, unlike the isotropic case, comparing the averaged chains of length $k$ does not require any information about the averaged chains of shorter length  but the isotropic bounds are essential. This is because the one-by-one
replacement strategy anyway breaks up averaged chains into isotropic ones, so information about shorter averaged
chains is useless  within any given $\gamma$-replacement step.
 
The proof of Lemma \ref{lemma:av_rep} is almost entirely analogous to those of Lemmas~\ref{lemma:1G_replacement} and \ref{lemma:kG_replacement}, presented in Sections~\ref{sec:1Grep} and~\ref{sec:k_iso_rep}, respectively. Therefore, we only illustrate the key differences, leaving the routine details to the reader.
\begin{proof}[Proof of Lemma \ref{lemma:av_rep}]
	We fix the vector of spectral parameters $\bm z \in \{z,\overline{z}\}^k$, and suppress the dependence of all quantities of $\bm z$.
	For all $k \in \indset{\maxK}$ and $\gamma \in \indset{\gamma_0}$, by analogy with \eqref{eq:kG_high_moment}, we introduce the random quantities $\Psi_{k, \mathrm{av}}^{\gamma, p} \equiv \Psi_{k,\mathrm{av}}^{\gamma, p}(\bm x)$, defined as 
	\begin{equation} \label{eq:kG_av_high_moment}
		\Psi_{k,\mathrm{av}}^{\gamma, p}(\bm x) := \Expv_{h_{ij}} \biggl[\bigl\lvert \Tr\bigl[(G^{(\gamma)}_{[1,k]}-M_{[1,k]})(\bm x')S^{x_k}\bigr] \bigr\rvert^{2p}\biggr] - \Expv_{v_{ij}} \biggl[\bigl\lvert \Tr\bigl[(G^{(\gamma-1)}_{[1,k]}-M_{[1,k]})(\bm x')S^{x_k}\bigr] \bigr\rvert^{2p}\biggr].
	\end{equation}
	Recall that our goal is to estimate $\Psi_{k,\mathrm{av}}^{\gamma, p}(\bm x)$ in a way that is summable in $\gamma := \phi\bigl((i,j)\bigr)$.
	
	Using the resolvent expansion \eqref{eq:resolvent_expand}, we obtain the following expression for $\Psi_{k, \mathrm{av}}^{\gamma, p}$, similar to \eqref{eq:moment_expansion} and \eqref{eq:Gk_moment_expansion},
	\begin{equation} \label{eq:Gk_av_moment_exp}
		\Psi_{k,\mathrm{av}}^{\gamma, p} = \sum_{r = 1}^{2p L}  \frac{1}{(1 + \delta_{ij})^{r}} \sum_{q + q' = r} \biggl(\mom_{q,q'}(h_{ij}) - \mom_{q,q'}(v_{ij})\biggr) \other{Z}_{k,\mathrm{av}}^{\gamma,p,L}(q,q') + \Omega_{k,\mathrm{av}}^{\gamma,p}(L),
	\end{equation}
	where $\other{Z}_{k,\mathrm{av}}^{\gamma,p,L}(q,q')$  comprises all terms in the expansion that contain exactly $q$ copies of $\partial_{ij}$ and $q'$ copies of $\partial_{ji}$, while the residual term $\Omega_{k,\mathrm{av}}^{\gamma,p}(L)$ comprises all terms in the expansion that contain the non-checked resolvents $G^{(\gamma)}$ or $G^{(\gamma-1)}$. Similarly to \eqref{eq:rep_Omega_est}, $\Omega_{k,\mathrm{av}}^{\gamma,p}(L)$ is negligible.  
	
	Similarly to \eqref{eq:Zk_replace_def}, the quantities $\other{Z}_{k,\mathrm{av},L}^{\gamma,p}(q,q')$ 
	admit the bound
	\begin{equation} \label{eq:Zk_av_rep}
		\bigl\lvert \other{Z}_{k,\mathrm{av}}^{\gamma,p,L}(q,q') \bigr\rvert  \lesssim \other{C}^{p,L}_{q,q'} \sum_{(J,J^*) \in \mathfrak{J}^{p,L}_{q,q'}} Z_{k,\mathrm{av}}^{\gamma,p,L}(J,J^*)  \bigl\lvert\bigl( \widecheck{\mathcal{X}}_k^{(\gamma)}\bigr)\bigr\rvert^{2p - |J| - |J^*|},
	\end{equation}
	where $\widecheck{\mathcal{X}}_k^{(\gamma)} : = \Tr\bigl[(\widecheck{G}^{(\gamma)}_{[1,k]}-M_{[1,k]})(\bm x')\bigr]$ and we define 
	\begin{equation} \label{eq:Zk_av_def}
		Z_{k,\mathrm{av}}^{\gamma,p,L}(J,J^*) := \bigl(S_{ij}\bigr)^{r/2} \biggl\lvert \prod_{\bm j \in J} \partial_{ij}^{\bm j} \widecheck{\mathcal{X}}^{(\gamma)}_k \biggr\rvert \biggl\lvert\prod_{\bm j \in J^*} \partial_{ij}^{\bm j} \overline{\widecheck{\mathcal{X}}_k^{(\gamma)}} \biggr\rvert, \quad r := \sumJ.
	\end{equation}
	Similarly to \eqref{eq:repZ_bound}, our goal is to establish a very-high-probability bound on $Z_{k,\mathrm{av}}^{\gamma,p,L}(J,J^*)$ of the form 
	\begin{equation} \label{eq:av_rep_Z_goal}
		Z_{k,\mathrm{av}}^{\gamma,p,L}(J,J^*)  \lesssim C_p \psi_{k,\mathrm{av}}^{(\gamma)}(J,J^*), \quad \text{w.v.h.p.},
	\end{equation}
	for some deterministic control parameters $\psi_{k,\mathrm{av}}^{(\gamma)}(J,J^*)$,   to be defined later,  
	such that $\psi_{k,\mathrm{av}}^{(\gamma)}(J,J^*)$ are summable in $\gamma := \phi\bigl((i,j)\bigr)$. 
	Once again, the expression on the right-hand side of \eqref{eq:Zk_av_def} is reminiscent of the summands in \eqref{eq:Zkav_def}, and we estimate it using the same approach. 
	
	 Similarly to \eqref{eq:partial_X_1}, we deduce that, for $u := |\bm j| \ge 1$, we conclude that  
	\begin{equation} \label{eq:rep_partial_X1}
		\bigl\lvert \partial_{ij}^{\bm j} \widecheck{\mathcal{X}}^{(\gamma)}_k \bigr\vert  \lesssim \sum_{1\le a_1\le\dots\le a_u\le k} \sum_{\substack{(c_l,d_l)\in \{(i,j),(j,i)\}\\ l \in \indset{u}}} \biggl(\prod_{l=1}^{u-1} \bigl\lvert (\widecheck{G}_{[a_l, a_{l+1}]})_{ c_l d_{l+1}} \bigr\rvert\biggr) \bigl\lvert \bigl(\widecheck{G}_{[a_u, k]}S^{x_k}\widecheck{G}_{[1, a_{1}]}\bigr)_{ c_u d_{1}} \bigr\rvert.
	\end{equation}
	Next, similarly to \eqref{eq:partial_X2}, we estimate each chain in \eqref{eq:rep_partial_X1} using \eqref{eq:M_bound} and \eqref{eq:k_rep_assume}, obtaining  
	\begin{equation} \label{eq:rep_partial_X2}
		\bigl\lvert \partial_{ij}^{\bm j} \widecheck{\mathcal{X}}^{(\gamma)}_k \bigr\vert  \lesssim \psi_{k,\mathrm{av}}^{(\gamma)}(\bm j),\quad \text{w.v.h.p.},
	\end{equation}
	uniformly in $\gamma \in \indset{\gamma_0}$ and $\bm x \in \indset{N}^{k}$, where 
	\begin{equation} \label{eq:rep_psi_def}
		\psi_{k,\mathrm{av}}^{(\gamma)}(\bm j) :=  \bigl(1+N^{2\xi'}(\ell\eta)^{\beta_k-1/2}\bigr)(\ell\eta)^{|\bm j|/2}\sum_{ 1\le a_1\le\dots\le a_{|\bm j|}\le k } \sum_{\substack{(c_l,d_l)\in \{(i,j),(j,i)\}\\ l \in \indset{|\bm j|}}}  \prod_{q=1}^{|\bm j|}\mathfrak{s}_l(c_l, d_{l+1}).
	\end{equation}
	Here $\mathfrak{s}_l$ are the decay factors corresponding to chains $\widecheck{G}_{[a_l, a_{l+1}]}$, and we used \eqref{eq:loss_exponents}.    At this point, it is important to address the special case $k=\maxK$. Observe that 
	(see  \eqref{eq:partial_X_1}  for more details), the $\partial_{ij}^{\bm j}$ derivative of an averaged $\maxK$-chain contains isotropic resolvent chains of length\footnote{ At this point we slightly deviate from the proof structure of the global
		law in Section~\ref{sec:global}. In  Section~\ref{sec:global3} we proved the isotropic global law for chains of length $\maxK+1$ directly, before proceeding to the proof of the averaged law for $\maxK$-chains. In the present setting of Section \ref{sec:GFT}, we do not have the freedom to do this, because the zig-step (Proposition~\ref{prop:zig}), which we use as an input to initialize the replacement procedure, does not provide any isotropic local laws for $(\maxK+1)$-chains. Thus we need
		to handle the $k=K$ case separately. This situation is reminiscent of the special treatment of the $k=K$ case
		in the martingale estimates in Section~\ref{sec:mart}, where  the longest chain also needed a separate splitting 
		scheme.} 
	$\maxK+1$, which cannot be estimated directly using \eqref{eq:k_rep_assume}. 
	However, we can apply a straightforward reduction inequality, similar to \eqref{eq:1_k+1_reduction1}. That is, for all $a\in\indset{\maxK}$,
	\begin{equation} \label{eq:rep_reduction}
		\biggl\lvert \bigl(\widecheck{G}^{(\gamma)}_{[a,\maxK]}S^{x_\maxK}\widecheck{G}^{(\gamma)}_{[1,a]}\bigr)_{cd} \biggr\rvert \lesssim \sum_{q} S_{x_{b}q} \bigl\lvert \bigl(\widecheck{G}^{(\gamma)}_{[a,b]}\bigr)_{cq} \bigr\rvert \bigl\lvert \bigl(\widecheck{G}^{(\gamma)}_{[b,a]}\bigr)_{qd} \bigr\rvert \lesssim \sqrt{\ell\eta}\,\mathfrak{s}^\mathrm{iso}_{k+1}(c, \bm y_a, d) \times N^{2\xi'}(\ell\eta)^{\alpha_{\maxK/2+1}} ,
	\end{equation}
	with very high probability, where $b \equiv b(a,\maxK)\in\indset{\maxK}$ is some appropriate cutting index such that the chains\footnote{If $b\ge a$, then we interpret $\widecheck{G}^{(\gamma)}_{[b,a]}$ as $\widecheck{G}^{(\gamma)}_{[b,\maxK]}S^{x_\maxK}\widecheck{G}^{(\gamma)}_{[1,a]}$, and similar interpretation applies to $\widecheck{G}^{(\gamma)}_{[b,a]}$ is $b<a$.} $\widecheck{G}^{(\gamma)}_{[b,a]}$ and $\widecheck{G}^{(\gamma)}_{[b,a]}$ have lengths $\maxK/2$ and $\maxK/2+1$, and we used \eqref{eq:M_bound}, \eqref{eq:k_rep_assume}, \eqref{eq:Schwarz_convol} in the second step. Here $\bm y_a := (x_a, \dots, x_k, x_1, \dots, x_{a-1})$ denotes a cyclic shift of  $\bm x\in \indset{N}^{k}$ for all $a\in\indset{k}$.Recall from \eqref{eq:loss_exponents}  that $\alpha_{\maxK/2+1} = \beta_\maxK - 1/2$, hence, by \eqref{eq:M_bound}, the bound \eqref{eq:rep_partial_X2} indeed holds for $k = \maxK$.

	 Hence, \eqref{eq:av_rep_Z_goal} holds with 
	\begin{equation} \label{eq:rep_Z_psi_def}
		\psi_{k,\mathrm{av}}^{(\gamma)}(J,J^*) := \bigl(S_{ij}\bigr)^{r/2} \prod_{\bm j \in J, J^*} \psi_{k,\mathrm{av}}^{(\gamma)}(\bm j),
	\end{equation} 
	where recall $r:=\sumJ$, and $\psi_{k,\mathrm{av}}^{(\gamma)}(\bm j)$ are defined in \eqref{eq:rep_psi_def}.
	By \eqref{eq:rep_psi_def}, the decay factors on the right-hand side of \eqref{eq:rep_Z_psi_def} contain $r$ copies of both indices $i$ and $j$ among their arguments (provided $c\neq d$). Similarly to \eqref{eq:partial_X_plug_in}, when summing $\psi_{k,\mathrm{av}}^{(\gamma)}(J,J^*)$ over $i$ and $j$, we use $(r-2)$ powers of $\sqrt{S_{ij}}$ to replace all but two of $j$'s by $i$'s, with $j$-summation performed by \eqref{eq:Schwarz_convol}, obtaining
	\begin{equation} \label{eq:sum_psi_est1}
		\sum_{i,j} \psi_{k,\mathrm{av}}^{(\gamma)}(J,J^*) \lesssim \frac{(\ell\eta)^{r/2}}{W^{(r-2)/2}}\bigl(1+N^{2\xi'}(\ell\eta)^{\beta_k-1/2}\bigr)^n  \sum_i\prod_{\bm j \in J, J^*}\biggl(\sum_{ \{a_l\} } \prod_{l=1}^{|\bm j|}\mathfrak{s}_l(i, i)\biggr).
	\end{equation}
	For the $i$-summation to be effective (see the discussion around \eqref{eq:non-effective}), each summand on the right-hand side of \eqref{eq:sum_psi_est1} has to contain at least four $\sqrt{\Upsilon}$-factors (recall that each $\mathfrak{s}_l$ is a product of $\sqrt{\Upsilon}$'s as defined in \eqref{eq:sfunc_def}) that connect the index $i$ to some external index. This condition is satisfied only when $|J|+|J^*|\ge 2$,   and in this case the summation over $i$ can be performed using \eqref{eq:Schwarz_convol}, yielding
 	\begin{equation} \label{eq:sum_psi_est2}
 		\sum_{i,j} \psi_{k,\mathrm{av}}^{(\gamma)}(J,J^*) \lesssim C_r\frac{\ell\eta }{W^{(r-2)/2}\eta}\bigl(1+N^{2\xi'}(\ell\eta)^{\beta_k-1/2}\bigr)^n  \bigl( \mathfrak{s}_k^\mathrm{av}(\bm x)\bigr)^n, \quad n=|J|+|J^*|\ge 2.
 	\end{equation}
 	It is straightforward to check, by definition of $\beta_k$ in \eqref{eq:loss_exponents}	and the choice of $\xi'$, that the scalar factor $\bigl(1+N^{2\xi'}(\ell\eta)^{\beta_k-1/2}\bigr)$ in \eqref{eq:sum_psi_est2} is bounded by $C(\ell\eta)^{\beta_k}$.
 	
 	 The condition $|J|+|J^*|\ge 2$ in \eqref{eq:sum_psi_est2} is, in fact, essential. Indeed, if $|J|+|J^*|=1$, and all partial derivatives act on $\widecheck{G}_1$, then, in the resulting decay factor, the index $i$ is connected to external indices via only two $\sqrt{\Upsilon}$. This precludes the use of \eqref{eq:Schwarz_convol}, forcing us to instead apply \eqref{eq:sqrt_convol}, which incurs an additional $\sqrt{N/\ell}$-cost---a loss that is generally unaffordable (e.g., when $r=4$).   
	Hence, to estimate the sum of $\psi_{k,\mathrm{av}}^{(\gamma)}(J,J^*)$ for $|J|+|J^*| = 1$, we need to refine the bound \eqref{eq:rep_partial_X2}  by separating the last factor in~\eqref{eq:rep_partial_X1} into its deterministic approximation
	and fluctuation. The deterministic term admits a more precise bound by \eqref{eq:resum_M_bound}, 
	while the fluctuation term is smaller by a factor $1/\sqrt{\ell \eta}$ which renders the $\sqrt{N/\ell}$ loss
	from \eqref{eq:sqrt_convol} affordable. 
	
	To this end, using \eqref{eq:k_rep_assume} only for chains of length up to $k$ in \eqref{eq:rep_partial_X1}, we obtain
	\begin{equation} \label{eq:rep_bound_long}
		\begin{split}
			\bigl\lvert \partial_{ij}^{\bm j} \widecheck{\mathcal{X}}^{(\gamma)}_k \bigr\vert \lesssim &~ \psi_{k,\mathrm{short}}^{(\gamma)}(\bm j) 
			+\sum_{1\le a \le k}~ \sum_{ c\in \{i,j\} }  \bigl\lvert \bigl(M_{[a,a]}^{(x_{k})}\bigr)_{cc} \bigr\rvert\\
			&+\sum_{1\le a \le k}~ \sum_{ c,d\in \{i,j\} }  \bigl\lvert \bigl(\widecheck{G}_{[a, k]}S^{x_k}\widecheck{G}_{[1, a]}- M_{[a,a]}^{(x_{k})}\bigr)_{cd} \bigr\rvert  ,
		\end{split}
	\end{equation}
	where $M_{[a,a]}^{(x_{k})}$ is the deterministic approximation corresponding to the chain $G_{[a, k]}S^{x_k} G_{[1, a]}$, and the control parameter $\psi_{k,\mathrm{short}}^{(\gamma)}(\bm j)$ is given by
	\begin{equation} \label{eq:rep_psi_def_short}
		\psi_{k,\mathrm{short}}^{(\gamma)}(\bm j) := N^{\xi'}(\ell\eta)^{u/2}\sum_{\substack{1\le a_1\le\dots\le a_u\le k \\ a_1 \le a_u-1}} \sum_{\substack{(c_q,d_q)\in \{(i,j),(j,i)\}\\ q \in \indset{u}}}  \mathfrak{s}_q(c_q, d_{q+1}),
	\end{equation}
	with $\mathfrak{s}_q$ as below \eqref{eq:rep_psi_def}.
	
	Now we estimate the contributions of all three terms in \eqref{eq:rep_bound_long} after summation over $(i,j)$.
 Note that for any set of $\{a_q\}$, the size functions $\mathfrak{s}_l$ on the right-hand side of \eqref{eq:rep_psi_def_short} contain at least four $\sqrt{\Upsilon}$-factors connecting $i$ and $j$ to external indices.  Hence, when the control parameters $\psi_{k,\mathrm{short}}^{(\gamma)}(\bm j)$ are summed over $i,j$, together with $(S_{ij})^{r/2}$, both $i$- and $j$-summations are effective, resulting in the same bound as in \eqref{eq:sum_psi_est2}, 
	\begin{equation}
		\sum_{i,j} \bigl(S_{ij}\bigr)^{r/2} \psi_{k,\mathrm{short}}^{(\gamma)}(\bm j) \lesssim C_r \frac{\ell\eta }{W^{(r-2)/2}\eta}  \bigl( (\ell\eta)^{\beta_k}\mathfrak{s}_k^\mathrm{av}(\bm x)\bigr)^n, \quad r \ge 2.
	\end{equation}
	Moreover, it follows from \eqref{eq:resum_M_bound} that, for all $a \in \indset{N}$,
	\begin{equation} \label{eq:rep_M_resumm}
		\sum_{i,j} \bigl(S_{ij}\bigr)^{r/2} \bigl\lvert \bigl(M_{[a,a]}^{(x_{k})}\bigr)_{ii} \bigr\rvert \lesssim \frac{\ell\eta}{W^{(r-2)/2}\eta} \mathfrak{s}^\mathrm{av}_k(\bm x), \quad r \ge 2.
	\end{equation}
	Therefore, it remains to estimate the terms in the second line of \eqref{eq:rep_bound_long}.

	Using \eqref{eq:loss_exponents}, \eqref{eq:k_rep_assume} for $k\in\indset{\maxK-1}$ and \eqref{eq:rep_reduction} for $k=\maxK$, we deduce that the $\widecheck{G}S\widecheck{G}-M$ term in the second line of \eqref{eq:rep_bound_long} admits the bound
	\begin{equation} \label{eq:rep_X_long_bound}
		\bigl\lvert \bigl(\widecheck{G}_{[a, k]}S^{x_k}\widecheck{G}_{[1, a]}- M_{[a,a]}^{(x_{k})}\bigr)_{cd} \bigr\rvert 
		\lesssim 
			N^{2\xi'}(\ell\eta)^{\beta_{k}}\mathfrak{s}^\mathrm{iso}_{k+1}(c, \bm y_a, d).
	\end{equation}
	Multiplying the right-hand side of \eqref{eq:rep_X_long_bound} by $(S_{ij})^{r/2}$ and summing over $c,d \in \{i,j\}$ and over $i,j\in\indset{N}$, we obtain, by \eqref{eq:Schwarz_convol} and  \eqref{eq:sqrt_convol},
	\begin{equation} \label{eq:rep_X_long_bound_sum}
		N^{2\xi'}(\ell\eta)^{\beta_{k}}\sum_{ij} \sum_{c,d\in\{i,j\}}\bigl(S_{ij}\bigr)^{r/2}  \mathfrak{s}^\mathrm{iso}_{k+1}(c, \bm y_a, d) \lesssim N^{2\xi'}(\ell\eta)^{\beta_{k}}  \frac{\sqrt{N\eta}}{W^{(r-2)/2}\eta}\mathfrak{s}^\mathrm{av}_{k}(\bm x), \quad r \ge 2.
	\end{equation}

	Recall that $\sqrt{N\eta} \lesssim N^{\etaexp/2}\ell\eta$; hence, combining \eqref{eq:sum_psi_est2}, \eqref{eq:rep_bound_long}--\eqref{eq:rep_M_resumm}, and \eqref{eq:rep_X_long_bound}--\eqref{eq:rep_X_long_bound_sum}, we conclude that
	\begin{equation} \label{eq:sum_psi_est_final}
		\sum_{i,j} \psi_{k,\mathrm{av}}^{(\gamma)}(J,J^*) \lesssim C_r \frac{\ell\eta}{W^{(r-2)/2}\eta} \bigl( (\ell\eta)^{\beta_k}\mathfrak{s}_k^\mathrm{av}(\bm x)\bigr)^n, \quad r=\sumJ\ge 2.
	\end{equation}

	Therefore, similarly to \eqref{eq:otherZ_bound}, it follows from \eqref{eq:sum_psi_est_final}, that, for $r= q + q'\ge 2$,
	\begin{equation} \label{eq:av_otherZ_bound}
			\sum_{\gamma \le \gamma_0} \Expv \biggl[\bigl\lvert\other{Z}_{k,\mathrm{av}}^{\gamma,p,L}(q,q')\bigr\rvert\biggr]  
			\lesssim  C_p \frac{\ell\eta}{W^{(r-2)/2}\eta} \sum_{n=1}^{r\wedge (2p)}\bigl( (\ell\eta)^{\beta_k}(\mathfrak{s}_{k}^\mathrm{av}(\bm x)\bigr)^{n} \max_{\gamma'\le\gamma_0}\bigl\lVert\widecheck{\mathcal{X}}_k^{(\gamma')}  \bigr\rVert_{2p-n}^{2p-n},
	\end{equation}
	where $r:=q+q'$. Therefore, using Young's inequality and the moment-matching conditions \eqref{eq:moment_match}--\eqref{eq:m4_cond}, similarly to \eqref{eq:psi_rep_bound1} and \eqref{eq:1G_rep}, we deduce that
	\begin{equation} \label{eq:k_av_rep_final}
		\sum_{\gamma \le \gamma_0} \Expv \biggl[\bigl\lvert \Psi_{k,\mathrm{av}}^{\gamma, p}\bigr\rvert\biggr] \lesssim  \other{C}_p \biggl(1 +  \frac{\lambda\ell}{W} + \frac{\ell}{W^2}\biggr)^{2p}  \bigl(N^{\xi/4}(\ell\eta)^{\beta_k}(\mathfrak{s}_{k}^\mathrm{av}(\bm x)\bigr)^{2p} + N^{-\xi/4}\max_{\gamma'\le\gamma_0}\bigl\lVert\widecheck{\mathcal{X}}_k^{(\gamma')}  \bigr\rVert_{2p}^{2p}.
	\end{equation}
	
	Analogously to Claim \ref{claim:rep_check}, it can be shown that the high moment norm $\max\limits_{\gamma'\le\gamma_0}\bigl\lVert\widecheck{\mathcal{X}}_k^{(\gamma')}  \bigr\rVert_{2p}$ in \eqref{eq:k_av_rep_final} can be replaced by $\max\limits_{\gamma'\le\gamma_0}\bigl\lVert \Tr\bigl[(G^{(\gamma'-1)}_{[1,k]}-M_{[1,k]})(\bm x')S^{x_k} \bigr\rVert_{2p}$ up to a negligible error. Moreover, by \eqref{eq:WN} and \eqref{eq:lambda_assume}, the scalar factor $1 + \lambda W^{-1}\ell + W^{-2}\ell$ is bounded by $CN^{\xi/2}$. Therefore, the desired \eqref{eq:kG_av_rep_goal} follows from \eqref{eq:k_av_rep_init} and \eqref{eq:k_av_rep_final} by a telescopic summation argument, similar to \eqref{eq:telescope}. This concludes the proof of Lemma \ref{lemma:av_rep}.
 \end{proof}
 
\subsection{Proof of Proposition \ref{prop:zag}} \label{sec:zag_proof_for_real}
We now combine the statements established in Section \ref{sec:GFT} above to conclude Proposition \ref{prop:zag}.
\begin{proof}[Proof of Proposition \ref{prop:zag}]
	By definition of stochastic domination in \eqref{eq:stochdom}, to establish the local laws for the resolvent $G = G^{(\gamma(N))}$ of $H$ with tolerance $\xi + \nu$ in the sense of Definition \ref{def:local_laws}, it suffices to obtain the high-moment bounds
	\begin{alignat}{2}
		\norm{\bigl((G^{(\gamma(N))}_{[1,k]}-M_{[1,k]})(\bm x')\bigr)_{ab} }_p &\lesssim C_p N^{\xi} (\ell\eta)^{\alpha_{k}} \mathfrak{s}_{k}^\mathrm{iso}(a,\bm x', b),\quad  &&\text{w.v.h.p.}, \label{eq:k_iso_moment}\\
		\norm{\Tr\bigl[\bigl(G^{(\gamma(N))}_{[1,k]}-M_{[1,k]}\bigr)(\bm x')S^{x_k}\bigr]}_p &\lesssim C_p N^{\xi} (\ell\eta)^{\beta_{k}} \mathfrak{s}_{k}^\mathrm{av}(\bm x),\quad  &&\text{w.v.h.p.} \label{eq:k_av_moment}
	\end{alignat}
	for all integers $p \in \mathbb{N}$, uniformly in $a,b\in\indset{N}$, $\bm x \in \indset{N}^{k}$, and $\bm z \in \{z,\overline{z}\}^{k}$, for all $k \in \indset{\maxK}$. 
	
	For $k=1$, \eqref{eq:k_iso_moment} follows from \eqref{eq:1G_GFT_goal}, established in Section \ref{sec:1Grep}.

	Then, using \eqref{eq:1G_GFT_goal} as the base case,  the bounds \eqref{eq:k_iso_moment}  for $k \in \indset{2,\maxK}$ follow immediately from Lemma \ref{lemma:kG_replacement} by strong induction in $k$ and a nested strong induction in $\gamma_0$.
	
	Finally, once the isotropic bounds are established, \eqref{eq:k_av_moment} for all $k\in\indset{\maxK}$  follow by strong induction in $\gamma_0 \in \indset{\gamma(N)}$ using Lemma~\ref{lemma:av_rep}.
	  This concludes the proof of Proposition~\ref{prop:zag}. 
\end{proof}

\section{$M$-term analysis}\label{sec:M}
For all $1\le j\le k$, and for any diagonal matrices $\{A_i\}_{i \ge 1}$, we define the deterministic approximation terms as 
\begin{equation} \label{eq:otherM}
	M_{[j,k]} \equiv M(z_j, A_j,z_{j+1}, \dots, A_{k-1},z_k) := \other{M}(z_j, A_j,z_{j+1}, \dots, A_{k-1},z_k) \times \prod_{i=j}^{k} m(z_i), 
\end{equation}
where the quantity $\other{M}_{[j,k]} \equiv \other{M}(z_j, A_j,z_{j+1}, \dots, A_{k-1},z_k)$ satisfies the recursive formula
\begin{equation} \label{eq:M_recursion}
	\other{M}_{[j,k]} = \frac{1}{1-m(z_j)m(z_k)S}\biggl[A_j   \other{M}_{[j+1,k]} + \sum_{i=j+1}^{k-1} m(z_j)m(z_i) \mathcal{S}\bigl[\other{M}_{[j,i]}\bigr]\other{M}_{[i,k]}\biggr],
\end{equation}
where  $\other{M}_{[j,j]} := 1$ and an empty summation is zero by convention. Since all $A_j$ are diagonal matrices, so are all $M_{[j,k]}$.

Note that, for $k-j \ge 2$, the first two terms in the recursive formula \eqref{eq:M_recursion} can be combined to obtain
\begin{equation} \label{eq:M_recursion2}
	\other{M}_{[j,k]} = \frac{1}{1-m(z_j)m(z_k)S}\biggl[\other{M}_{[j,j+1]}\other{M}_{[j+1,k]} + \sum_{i=j+2}^{k-1} m(z_j)m(z_i) \mathcal{S}\bigl[\other{M}_{[j,i]}\bigr]\other{M}_{[i,k]}\biggr], \quad k \ge j+2.
\end{equation}

In the next subsections we prove certain non-trivial identities on the deterministic approximations $M_{[j,k]}$
inductively directly from their recursive definitions. We mention that alternatively these relations can also be proven
by the so-called {\it meta-argument} by tensorization  \cite{cook2018non, Cipolloni2022overlap} which exploits the fact that the analogous identities are
much easier for the original resolvent chains $G_{[j,k]}$  
and global laws are available to connect $G_{[j,k]}$  with $M_{[j,k]}$.

\subsection{Properties of the recurrence relation: Proof of Lemmas \ref{lemma:M_Ward} and \ref{lemma:dM}} \label{sec:M_rec_analysis}
We start with a technical lemma that will be needed in the sequel.
\begin{lemma}[Cyclic Property] \label{lemma:M_cyclic} 
	For any $k \ge 2$, 
	\begin{equation} \label{eq:M_cyclic}
		\Tr\bigl[\other{M}(z_1, A_1, z_2, \dots, A_{k-1}, z_k )\, A_k\bigr] = \Tr\bigl[A_{k-1} \other{M}(z_k , A_k, z_1, A_1, z_2, \dots, A_{k-2}, z_{k-1}) \bigr]
	\end{equation}
\end{lemma}
Note that it follows immediately from \eqref{eq:otherM} that \eqref{eq:M_cyclic} also holds with $\other{M}$ replaced by $M$.

\begin{proof} [Proof of Lemma \ref{lemma:M_cyclic}]
	We prove \eqref{eq:M_cyclic} by induction in $k \ge 2$. The base case $k=2$ follows immediately from \eqref{eq:M_recursion}. Indeed,
	\begin{equation}
		\begin{split}
			\Tr\bigl[\other{M}(z_1, A_1, z_2)\, A_2\bigr] &= \Tr\biggl[\bigl(1- m(z_1)m(z_2)\mathcal{S}\bigr)^{-1}[A_1]\, A_2\biggr]\\ 
			&= \Tr\biggl[A_1\,\bigl(1- m(z_1)m(z_2)\mathcal{S}\bigr)^{-1}[A_2]\biggr] = \Tr\bigl[A_1 \other{M}(z_2, A_2, z_1)\bigr].
		\end{split}
	\end{equation}
	
	Next, we prove the induction step. Assume that \eqref{eq:M_cyclic} holds for all lengths $k \le j-1$ for some $j \ge 3$. It follows from \eqref{eq:M_recursion2} that 
	\begin{equation} \label{eq:cyclic_expand1}
		\begin{split}
			\Tr\biggl[\other{M}_{[1,j]}A_j\biggr] =&~ \Tr\biggl[ \other{M}_{[j,1]} \other{M}_{[1,2]}\other{M}_{[2,j]}\biggr] + \sum_{i=3}^{j-1} m(z_1)m(z_i) \Tr\biggl[ M_{[j,1]}\mathcal{S}\bigl[\other{M}_{[1,i]}\bigr]\other{M}_{[i,j]}\biggr]\\
			=&~ \Tr\biggl[ \other{M}\bigl(z_j, \other{M}_{[j,1]} \other{M}_{[1,2]}, z_2, A_2, \dots, A_{j-2}, z_{j-1}  \bigr)A_{j-1}\biggr]\\
			&+\sum_{i=3}^{j-1} \Tr\biggl[ \other{M}\bigl(z_j, M_{[j,1]}m(z_1)m(z_i) \mathcal{S}[\other{M}_{[1,i]}], z_i, A_{i+1}, \dots,A_{j-2}, z_{j-1}\bigr)A_{j-1}\biggr]
			,
		\end{split}
	\end{equation}
	where we denote $\other{M}_{[j,1]} := \other{M}(z_j, A_j, z_1)$. Here in the second equality we used the induction hypothesis for $\{\other{M}_{[i,j]}\}_{i=2}^{j-1}$ with length at most $j-1$. Therefore, to prove \eqref{eq:M_cyclic} for $k=j$, it suffices to show that for any $z_0 \in \mathbb{C}\backslash\mathbb{R}$ and any diagonal matrix $A_0$,
	\begin{equation} \label{eq:cyclic_goal2}
		\sum_{i=2}^{j-2} \other{M}\bigl(z_0, \other{A}_{i-1}, z_i, A_{i}, \dots, z_{j-1}\bigr) + \other{M}\bigl(z_0, \other{A}_{j-2}, z_{j-1}\bigr)= \other{M}(z_0, A_0, z_1, A_1, z_2, \dots, A_{j-2}, z_{j-1}),
	\end{equation}
	where the matrices $\other{A}_b \equiv \other{A}_{b}(z_0, A_0, z_1, A_1, \dots, A_{b},z_{b+1})$ are defined as 
	\begin{equation} \label{eq:otherAs}
		\other{A}_1 := \other{M}(z_0, A_0, z_1) \other{M}_{[1,2]}, \quad \other{A}_b := \other{M}(z_0, A_0, z_1)\,m(z_1)m(z_{b+1}) \mathcal{S}[\other{M}_{[1,b+1]}], \quad b \ge 2.
	\end{equation}
	Indeed, choosing $z_0 := z_j$ and $A_0 := A_j$, we conclude from \eqref{eq:cyclic_expand1} and \eqref{eq:cyclic_goal2} that
	\begin{equation}
		\Tr\biggl[\other{M}_{[1,j]}A_j\biggr] = \Tr\biggl[ \other{M}(z_j, A_j, z_1, A_1, z_2, \dots, A_{j-2}, z_{j-1}) A_{j-1}\biggr],
	\end{equation} 
	which implies \eqref{eq:M_cyclic} for $k=j$.
	
	Therefore, it remains to prove \eqref{eq:cyclic_goal2}. Once again, we proceed by induction in $j \ge 3$.
	The base case $j=3$ holds by \eqref{eq:M_recursion2}. Indeed, for $j=3$, the left-hand side of \eqref{eq:cyclic_goal2} is given by
	\begin{equation}
		\other{M}\bigl(z_0, \other{A}_{1}, z_{2}\bigr) = \bigl(1 - m(z_0)m(z_2)\mathcal{S}\bigr)^{-1}\biggl[\other{M}(z_0, A_0, z_1) \other{M}_{[1,2]}\biggr] = \other{M}(z_0, A_0, z_1, A_1, z_2).
	\end{equation} 
	
	To prove the induction step, we analyze left-hand side of \eqref{eq:cyclic_goal2} for $j\ge 4$. Denote
	\begin{equation} \label{eq:Dj_def}
		D_j := \bigl(1-m(z_0)m(z_{j-1})\mathcal{S}\bigr)\biggl[\sum_{i=2}^{j-2} \other{M}\bigl(z_0, \other{A}_{i-1}, z_i, A_{i}, \dots, z_{j-1}\bigr) \biggr].
	\end{equation}
	Then it follows from \eqref{eq:M_recursion2} that
	\begin{equation} \label{eq:cyclic_expand2}
		\begin{split}
			D_j=&~ \sum_{i=2}^{j-2} \other{M}(z_0, \other{A}_{i-1},z_i)\other{M}_{[i,j-1]} + \sum_{i=2}^{j-2}\sum_{p=i+1}^{j-2} m(z_0)m(z_p) \mathcal{S}\bigl[\other{M}(z_0, \other{A}_{i-1}, z_i, \dots, z_p)\bigr]\other{M}_{[p,j-1]} \\
			=&~ \sum_{p=2}^{j-2} \biggl(m(z_0)m(z_p) \mathcal{S}\biggl[\sum_{i=2}^{p-1}\other{M}(z_0, \other{A}_{i-1}, z_i, \dots, z_p) + \other{M}(z_0, \other{A}_{p-1},z_p)\biggr] +\other{A}_{p-1} \biggr)\other{M}_{[p,j-1]},
		\end{split}
	\end{equation}
	where in the second step we changed the order of summation and used \eqref{eq:M_recursion} for $\other{M}(z_0, \other{A}_{p-1},z_p)$. Since $p\le j-2$, we can now use the induction hypothesis to apply the identity \eqref{eq:cyclic_goal2} to the argument of the super-operator $\mathcal{S}$ in the last line of \eqref{eq:cyclic_expand2}, thus obtaining
	\begin{equation} \label{eq:Dj_final}
		\begin{split}
			D_j=&~ \sum_{p=2}^{j-2} m(z_0)m(z_p) \mathcal{S}\bigl[\other{M}(z_0, A_0, z_1, \dots, z_p)\bigr]\other{M}_{[p,j-1]}  + \sum_{p=2}^{j-2} \other{A}_{p-1} \other{M}_{[p,j-1]}
			\\
			=&~ \bigl(1- m(z_0)m(z_{j-1})\mathcal{S}\bigr)\bigl[\other{M}_{[0,j-1]}\bigr]\\
			&+  \other{M}(z_0, A_0,z_1)\biggl(-\other{M}_{[1,j-1]} + \other{M}_{[1,2]} \other{M}_{[2,j-1]} + \sum_{p=3}^{j-2} m(z_1)m(z_{p}) \mathcal{S}[\other{M}_{[1,p]}]\other{M}_{[p,j-1]}\biggr)\\
			=&~ \bigl(1- m(z_0)m(z_{j-1})\mathcal{S}\bigr)\bigl[\other{M}_{[0,j-1]}\bigr] -  \other{M}(z_0, A_0, z_1)\,m(z_1)m(z_{j-1}) \mathcal{S}[\other{M}_{[1,j-1]}].
		\end{split}
	\end{equation}
	where in the second equality we used \eqref{eq:M_recursion2} for $\other{M}_{[0,j-1]} := \other{M}(z_0,A_0,z_1,\dots, A_{j-2},z_{j-1})$ and the definition of $\other{A}_j$ in \eqref{eq:otherAs}, while in the last step we used \eqref{eq:M_recursion2} for $\other{M}_{[1,j-1]}$.
	On the other hand, it follows immediately from the definition of $\other{A}_{j-2}$ in \eqref{eq:otherAs} that
	\begin{equation} \label{eq:Dj_residue}
		\bigl(1-m(z_0)m(z_{j-1})\mathcal{S}\bigr)\biggl[\other{M}\bigl(z_0, \other{A}_{j-2}, z_{j-1}\bigr)\biggr] = \other{A}_{j-2} = \other{M}(z_0, A_0, z_1)\,m(z_1)m(z_{j-1}) \mathcal{S}[\other{M}_{[1,j-1]}].
	\end{equation}
	Therefore, we conclude from \eqref{eq:Dj_final} and \eqref{eq:Dj_residue}, that
	\begin{equation}
		\bigl(1-m(z_0)m(z_{j-1})\mathcal{S}\bigr)^{-1}\bigl[D_j\bigr] + \other{M}\bigl(z_0, \other{A}_{j-2}, z_{j-1}\bigr) = \other{M}(z_0,A_0,z_1,\dots, A_{j-2},z_{j-1}).
	\end{equation}
	In particular, by definition of $D_j$ in \eqref{eq:Dj_def}, the identity \eqref{eq:cyclic_goal2} holds.
	This concludes the proof of Lemma \ref{lemma:M_cyclic}.	
\end{proof}

First, we establish the divided difference identity for $M$, proving Lemma \ref{lemma:M_Ward}.
\begin{proof}[Proof of Lemma \ref{lemma:M_Ward}]
	We prove \eqref{eq:M_Ward} by induction in $k \ge 2$. 
	
	For the base case $k=2$, we use \eqref{eq:otherM} and \eqref{eq:M_recursion} to obtain
	\begin{equation}
		M(z_1, I, z_2) = m(z_1)m(z_2) \bigl(1-m(z_1)m(z_2) \mathcal{S}\bigr)^{-1}\bigl[I\bigr] = \frac{m(z_1)m(z_2)}{1-m(z_1)m(z_2)} = \frac{m(z_1)-m(z_2)}{z_1-z_2},
	\end{equation}
	where in the last step we used the divided difference for $m(z) = m_{\mathrm{sc}}(z)$ that follows immediately by subtracting two copies of the Dyson equation \eqref{eq:Dyson}. 
	
	Next, we prove the induction step. Assuming $k\ge 3$, it follows from \eqref{eq:M_recursion2} that
	\begin{equation} \label{eq:M_Ward_expan}
		M_{[1,k]}\bigr\rvert_{A_j = I} = \frac{m(z_1)}{1-m(z_1)m(z_k)\mathcal{S}}\biggl[\frac{M_{[1,2]}M_{[2,k]}}{m(z_2)m(z_1)} + \sum_{i=3}^{k-1}  \mathcal{S}\bigl[M_{[1,i]}\bigr]M_{[i,k]}\biggr]\biggr\rvert_{A_j = I}~.
	\end{equation}
	
	For any $j \in \indset{2, k-1}$, the desired identity \eqref{eq:M_Ward} follows immediately applying the induction hypothesis to every $M$ in \eqref{eq:M_Ward_expan} containing the observable $A_j$, collecting the terms, and using \eqref{eq:M_recursion} to obtain the difference of two $M$'s of length $k-1$. We omit the straightforward algebraic manipulation.
	
	Therefore, it remain to consider $j=1$. Observe that the stability operator satisfies the identity
	\begin{equation} \label{eq:stab_identity}
		\frac{m(z_1)-m(z_2)}{m(z_2)}\frac{1}{1-m(z_1)m(z_k)\mathcal{S}} = \frac{m(z_1)}{m(z_2)}\frac{1-m(z_2)m(z_k)\mathcal{S}}{1-m(z_1)m(z_k)\mathcal{S}} - 1.
	\end{equation}
	Therefore, using \eqref{eq:M_Ward_expan} and the induction hypothesis for $M_{[1,i]}\rvert_{A_1 = I}$ with $i \in \indset{k-1}$, we obtain
	\begin{equation} \label{eq:Ward_j=1_case}
		\begin{split}
			M_{[1,k]}\bigr\rvert_{A_1 = I} = \frac{1}{z_1 - z_2}\biggl(&~\frac{m(z_1)-m(z_2)}{m(z_2)}\frac{1}{1-m(z_1)m(z_k)\mathcal{S}}\biggl[  M_{[2,k]}\biggr]\\
			&+ \frac{m(z_1)}{1-m(z_1)m(z_k)\mathcal{S}}\biggl[ \sum_{i=3}^{k-1}  \mathcal{S}\bigl[M_{[1,\widehat{2},i]}\bigr]M_{[i,k]} - \sum_{i=3}^{k-1}  \mathcal{S}\bigl[M_{[2,i]}\bigr]M_{[i,k]}\biggr]\biggr)\\
			= \frac{1}{z_1 - z_2}\biggl(&~\frac{m(z_1)}{1-m(z_1)m(z_k)\mathcal{S}}\biggl[  A_2M_{[3,k]} + \sum_{i=3}^{k-1}  \mathcal{S}\bigl[M_{[1,\widehat{2},i]}\bigr]M_{[i,k]}\biggr] - M_{[2,k]} \biggr), 
		\end{split}
	\end{equation}
	where we denote $M_{[1,\widehat{2},i]} := M(z_1, A_2, z_3, \dots, A_{i-1}, z_i)$, and in the second step we used the stability operator identity \eqref{eq:stab_identity}, and the recursion \eqref{eq:M_recursion} for $M_{[2,k]}$. In particular, \eqref{eq:Ward_j=1_case} together with \eqref{eq:M_recursion} imply that \eqref{eq:M_Ward} holds in the case $j=1$.
	This concludes the proof of Lemma \ref{lemma:M_Ward}.  
\end{proof}

\begin{proof} [Proof of Lemma \ref{lemma:dM}]
	We prove \eqref{eq:dM} by induction in $k$.
	For the base case $k=1$, recall that
	\begin{equation} \label{eq:dm}
		\frac{\mathrm{d}}{\mathrm{d}t}m(z_t) = \frac{1}{2}m(z_t).
	\end{equation}
	In particular, denoting $m_{i,t} := m(z_{i,t})$ for all $i\in\indset{k}$,
	\begin{equation}
		\frac{\mathrm{d}}{\mathrm{d}t}\frac{1}{1-m_{j,t}m_{k,t}S} = \frac{m_{j,t}m_{k,t}S}{\bigl(1-m_{j,t}m_{k,t}S\bigr)^2}.
	\end{equation}
	Using the definition \eqref{eq:Mt_def}, \eqref{eq:otherM} and the recursion \eqref{eq:M_recursion}, we deduce that 
	\begin{equation}
		\begin{split}
			\frac{\mathrm{d}}{\mathrm{d}t}M_{[1,k],t} =&~ \frac{1}{2}M_{[1,k],t} + \frac{m_{1,t}m_{k,t}S}{1-m_{1,t}m_{k,t}S}\bigl[M_{[1,k],t}\bigr] + \frac{m_{1,t}}{1-m_{1,t}m_{k,t}S}\biggl[A_1  \frac{\mathrm{d}}{\mathrm{d}t}M_{[2,k],t} \biggr]\\
			&+ \frac{m_{1,t}}{1-m_{1,t}m_{k,t}S}\biggl[\sum_{p=2}^{k-1} \mathcal{S}\biggl[\frac{\mathrm{d}}{\mathrm{d}t}M_{[1,p],t}\biggr]M_{[p,k],t} + \sum_{p=2}^{k-1} \mathcal{S}\bigl[M_{[1,p],t}\bigr]\frac{\mathrm{d}}{\mathrm{d}t}M_{[p,k],t}\biggr].
		\end{split}
	\end{equation}
	Using the induction hypothesis \eqref{eq:dM} for chains of length at most $k-1$ and rearranging the terms, we obtain 
	\begin{equation} \label{eq:dM_intermediate}
		\begin{split}
			\frac{\mathrm{d}}{\mathrm{d}t}M_{[1,k],t} =&~ \frac{k}{2}M_{[1,k],t}  + \sum_q\bigl(M_{[1,k],t}\bigr)_{qq} M_{[1,1],[k,k]}^{(q)}  + \Sigma_1 + \Sigma_2,
		\end{split}
	\end{equation}
	where the sums $\Sigma_1$ and $\Sigma_2$ are given by
	\begin{equation} 
		\begin{split}
			\Sigma_1 := \sum_{2 \le i < j \le k} \sum_q \bigl((M_{[i,j],t}\bigr)_{qq} \frac{m_{1,t}}{1-m_{1,t}m_{k,t}S}\biggl[&~A_1  M_{[2,i],[j,k],t}^{(q)} + \sum_{p=j}^{k-1} \mathcal{S}\bigl[ M_{[1,i],[j,p],t}^{(q)}\bigr]M_{[p,k],t}\\
			&+ \sum_{p=2}^{i}\mathcal{S}\bigl[M_{[1,p],t}\bigr]M_{[p,i],[j,k],t}^{(q)}\biggr],
		\end{split}
	\end{equation}
	\begin{equation} \label{eq:dM_Sigmas}
		\Sigma_2 := \sum_{j=2}^{k-1} \sum_q \bigl((M_{[1,j],t}\bigr)_{qq}\frac{m_{1,t}}{1-m_{1,t}m_{k,t}S}\biggl[S^qM_{[j,k],t}+ \sum_{p=j}^{k-1} \mathcal{S}\bigl[ M_{[1,1],[j,p],t}^{(q)}\bigr]M_{[p,k],t}\biggr] .
	\end{equation}
	It follows from the recursion \eqref{eq:M_recursion} for $M_{[1,j],[j,k],t}^{(q)}$ that 
	\begin{equation}
		\Sigma_1 = \sum_{2 \le i < j \le k} \sum_q \bigl((M_{[i,j],t}\bigr)_{qq}M_{[1,i],[j,k],t}^{(q)}, \quad \Sigma_2=\sum_{j=2}^{k-1} \sum_q \bigl((M_{[1,j],t}\bigr)_{qq}M_{[1,1],[j,k],t}^{(q)}
	\end{equation}
	Therefore, \eqref{eq:dM} holds by \eqref{eq:dM_intermediate} and \eqref{eq:dM_Sigmas}. This concludes the proof of Lemma \ref{lemma:dM}.
\end{proof}
We conclude this subsection by proving the re-expansion identity \eqref{eq:M_re-expansion_identity}.
\begin{proof} [Proof of \eqref{eq:M_re-expansion_identity}] Fix an index $j \in \indset{k}$, and let $\widehat{M}_{[1,k]}$ denote the right-hand side of \eqref{eq:M_re-expansion_identity}.
	 It follows from \eqref{eq:otherM}, \eqref{eq:M_recursion} and \eqref{eq:M_cyclic} that, for any diagonal matrix $A \in \mathbb{C}^{N\times N}$,
	\begin{equation} \label{eq:re-expansion_M}
		\begin{split}
			\Tr\bigl[ M_{[1,k+1]}\, A\bigr] =&~ \Tr\bigl[ M(z_{j+1}, A_{j+1}, \dots, z_{k+1}, A, z_1, \dots, z_{j})A_j\bigr]\\
			=&~ \Tr\bigl[A_j' A_{j+1}   M(z_{j+2}, A_{j+2}, \dots, z_{k+1}, A, z_1, \dots, z_{j}) \bigr]\\
			&+\sum_{i=j+2}^{k+1} \Tr\biggl[  \mathcal{S}\bigl[M_{[j+1,i]}\bigr]M(z_{i}, \dots, z_{k+1}, A, z_1, \dots, z_{j}) A_j'\biggr]\\
			&+\sum_{i=1}^{j-1} \Tr\biggl[  M(z_{j+1}, A_{j+1}, \dots, z_{k+1}, A, z_1, \dots, z_{i})\mathcal{S}\bigl[M_{[i,j]} A_j'\bigr] \biggr]\\
			=&~ \Tr \bigl[ \widehat{M}_{[1,k]} \, A  \bigr]\\
		\end{split}
	\end{equation}
	where in the second we used that $\Tr[ (1-m_{j+1}m_j\mathcal{S})^{-1}[X]Y] = \Tr[ (1-m_{j+1}m_j\mathcal{S})^{-1}[Y]X]$, the symmetry of $\mathcal{S}$, and in the last step we used \eqref{eq:M_cyclic} again. Since \eqref{eq:re-expansion_M} holds for arbitrary diagonal $A$, we conclude that $M_{[1,k+1]} = \widehat{M}_{[1,k+1]}$. This concludes the proof of \eqref{eq:re-expansion_M}. 
\end{proof}

\subsection{$M$-bounds: Proof of Lemmas  \ref{lemma:M_bounds}, \ref{lemma:extra_M}, and Claim \ref{claim:staticM}} \label{sec:M_bounds} 
We prove the estimates \eqref{eq:M_bound_av}--\eqref{eq:M_bound} following the same dynamical approach as we used for $(G-M)$-quantities in Sections \ref{sec:masters_sec}--\ref{sec:reg_sec} with the key idea being the 
observable regularization from Section \ref{sec:reg_sec}. In fact, the proof of Lemma \ref{lemma:M_bounds} can be viewed as a simplified version of Proposition \ref{prop:masters}. Recall that the deterministic $M$-quantities satisfy the evolution equation \eqref{eq:dM}, with the right-hand side of \eqref{eq:dM} for $M_{[1,k]}$ containing chains of equal or shorter length than $k$. Therefore, we prove Lemma \ref{lemma:M_bounds} by induction in $k$. Before proceeding with the proof, we fix some notation and preliminary facts about the $M$-terms.

Although the $M$-terms do not exhibit any explicit product structure, unlike the corresponding resolvent chains, we refer to the number of spectral parameters that the $M$-term depends on as its \emph{length}. For example, the length of $M_{[1,k],t}$ is $k$. This concept of length coincides with the length of the corresponding $G$-chain. Note that an $M$-term of length $k$ depends on $k-1$ external observable indices.
 
\begin{proof} [Proof of Lemma \ref{lemma:M_bounds}]
	First, observe that for $k=1$ we have $M_{[1,1],t}  = m(z_t)$,  so \eqref{eq:M_bound_av}--\eqref{eq:M_bound}
	is immediate since $|m(z)| \lesssim 1$.  For $k=2$, recall that
	 $M_{[1,2],t}$ is given by \eqref{eq:M12}
	with $A_1=S^{x_1}$. 
	  Inserting the bounds~\eqref{eq:Ups_majorates_notime} into the definition of $M_{[1,2]}$
	we immediately see that
	\eqref{eq:M_bound_av}--\eqref{eq:M_bound} hold for $k\in \{1,2\}$ by $|m(z)| \lesssim 1$.
	Therefore, we only need to prove \eqref{eq:M_bound_av}--\eqref{eq:M_bound} for $k \ge 3$. We proceed by induction in $k$.
	
	Fix an integer $k \ge 3$, and assume that for all $k' \in \indset{2,k-1}$, the bounds
	\begin{equation} \label{eq:M_bound_ind}
		\begin{split}
			\bigl\lvert  \Tr \big[ M_{[1,k'],t}(\bm z_t, \bm{x}') S^{x_{k'}} \big] \bigr\rvert &\lesssim  (\log N)^{C_{k'}+1} \times(\ell_t\eta_t)\,\mathfrak{s}_{k',t}^\mathrm{av}(\bm x),\\
			\bigl\lvert \bigl(M_{[1,k'],t}(\bm z_t, \bm{x}') \bigr)_{ab}\bigr\rvert &\lesssim \delta_{ab} \,(\log N)^{C_{k'}+1} \times \sqrt{\ell_t\eta_t}\,\mathfrak{s}_{k',t}^\mathrm{iso}(a, \bm x', b).
		\end{split}
	\end{equation}
	hold uniformly in $\bm x\in\indset{N}^{k'}$, $a,b\in \indset{N}$ and $\bm z_t\in \dom_t^{k'}$. 
	Recall that $\eta_t:= \min_j \eta_{j,t}\sim \eta_{i,t}$
	for any $i\in \indset{\maxK}$.  
	Our goal is to establish \eqref{eq:M_bound_av}--\eqref{eq:M_bound} for $M_{[1,k]}$.
 	
 	Similarly to \eqref{eq:G-M_kav} and \eqref{eq:G-M_kiso}, define the deterministic time-dependent quantities
	\begin{equation}
		X_t^k \equiv X_{\bm z_t,t}^k(\bm x) := \Tr\bigl[M_{[1,k],t}(\bm x') S^{x_k}\bigr], \quad Y_t^k \equiv Y_{\bm z_t,t}^k(\bm x', a) := \bigl(M_{[1,k],t}(\bm x') \bigr)_{aa},
	\end{equation}
	for all $\bm x \in \indset{N}^k$ and $a \in \indset{N}$. 
	 We do not follow the precise dependence of $X_t^k$ and $Y_t^k$ on $\bm z_t$, and in the sequel all bounds are understood to be uniform in $\bm z_t \in \dom_t^k$. 
	  Since $M_{[1,k],t}$’s are diagonal matrices, it suffices to analyze $(M_{[1,k],t})_{aa}$ for all $a\in\indset{N}$
	  unlike for  $\mathcal{Y}_t^k$. 
	
	Note that at the initial time $t=0$, the quantities $Y_0^k$ and $X_0^k$ are controlled by Claim \ref{claim:staticM}, which we prove in the end of Section \ref{sec:M_bounds}. Indeed, it follows from \eqref{eq:staticM} with $n := 1$ that 
	\begin{equation} \label{eq:M_bound_init}
		\bigl\lvert Y_0^k (\bm x',a)\bigr\rvert \lesssim  \sqrt{\ell_0\eta_0} \,\mathfrak{s}_{k,0}^\mathrm{iso}(a,\bm x',a), \quad \bm x \in \indset{N}^k, \quad a\in \indset{N}, 
	\end{equation}
	where we used that $\ell_0 \sim W$ and $\eta_0 \sim 1$. Hence, we the averaged quantity $X_0^k (\bm x)$ at time $t=0$ satisfies
	\begin{equation} \label{eq:av_M_bound_init}
		\bigl\lvert X_0^k (\bm x)\bigr\rvert \lesssim  \sqrt{\ell_0\eta_0} \sum_{a}S_{x_ka}\,\mathfrak{s}_{k,0}^\mathrm{iso}(a,\bm x',a) \lesssim \ell_0\eta_0 \,\mathfrak{s}_{k,0}^\mathrm{av}(\bm x), \quad \bm x \in \indset{N}^k, 
	\end{equation}
	where we used \eqref{eq:true_convol} and \eqref{eq:sfunc_def}.
	
	\vspace{5pt}
	
	\textbf{Averaged $M$-terms, proof of \eqref{eq:M_bound_av}.}\\ 
	Differentiating $X_t^k$ in time $t$, using \eqref{eq:dM}, we find
	\begin{equation} \label{eq:av_M_evol}
		\frac{\mathrm{d}}{\mathrm{d}t}X_t^k = \biggl(\frac{k}{2}I + \bigoplus_{j=1}^k \mathcal{A}_{j,t}\biggr)\bigl[X_t^k \bigr] +  F_{[1,k],t}^\mathrm{av},
	\end{equation}
	where the linear operators $\mathcal{A}_{j,t}$ are defined in \eqref{eq:lin_prop_ops}. Here the forcing terms $F_{[1,k],t}^\mathrm{av}$ are given by
	\begin{equation} \label{eq:av_M_forcing}
		F_{[1,k],t}^\mathrm{av} \equiv F_{[1,k],t}^\mathrm{av}(\bm x) : = \sum_{\substack{1 \le i < j \le k \\ 2 \le j-i \le k-2}} \sum_q \bigl(M_{[i,j],t}\bigr)_{qq} \Tr \bigl[M_{[1,i],[j,k],t}^{(q)} S^{x_k} \bigr].
	\end{equation}
	Note that the structure of \eqref{eq:av_M_evol} is similar to, yet significantly simpler than, that of \eqref{eq:k_av_evol}. This is because it lacks a martingale term, and all 
	$M$-terms on the right-hand side of \eqref{eq:av_M_forcing} have length at most $k-1$. 
	Consequently, unlike in the proof of the corresponding Lemma~\ref{lemma:forcing}, no reduction procedure or stopping time argument is necessary. 
	Hence, since \eqref{eq:M_bound_ind} hold for all lengths $2 \le k'\le k-1$ by induction assumption, we have the bound
	\begin{equation} \label{eq:av_M_forcing_bound}
		\bigl\lvert F_{[1,k],t}^\mathrm{av} \bigr\rvert \lesssim (\log N)^{\other{C}_k} \frac{1}{\eta_t} \times \ell_t\eta_t \, \mathfrak{s}_{k,t}^\mathrm{av}(\bm x), \quad t \in [0,T].
	\end{equation}
	uniformly in $\bm x \in \indset{N}^k$, where  the positive constant $\other{C}_k$ depends only on $k$ and $\etaexp$ in \eqref{eq:etabound}. Here we used \eqref{eq:Schwarz_convol} to perform the summation in $q$, as in \eqref{eq:av_insertion}.

	Next, as in Section \ref{sec:masters_sec}, we need to consider two cases: If $M_{[1,k]}$ corresponds to a non-saturated chain,  then we can apply \eqref{eq:av_prop_bound} from Lemma \ref{lemma:good_props}
	when solving \eqref{eq:av_M_evol} by Duhamel's formula.
	 On the other hand, if $M_{[1,k]}$ corresponds to a saturated chain, we need to use observable regularization.
	
	\vspace{5pt}
	
	\textbf{Case 1.} Assume that $M_{[1,k]}$ corresponds to a non-saturated chain, i.e., $\bm z_t \neq \bm z_{t,\mathrm{sat}}$, defined in \eqref{eq:sat_cahins}. Let  $\other{\mathcal{P}}^k_{s, \tau}$  denote the linear the propagator is given by \eqref{eq:k_av_propagator}, 
	\begin{equation} \label{eq:M_av_propagator}
		\other{\mathcal{P}}^k_{s, t} := \exp\biggl\{ \int_{s}^t \bigoplus_{j=1}^k \bigl(I + \mathcal{A}_{j,r}\bigr)\mathrm{d}r \biggr\}, \quad 0 \le s\le t\le T.
	\end{equation}
	Since $\bm z_t$ is not saturated, the propagator $\other{\mathcal{P}}^k_{s, t}$ satisfies the implication
	 in \eqref{eq:av_prop_bound}.
	Applying Duhamel's principle to \eqref{eq:av_M_evol}, we obtain
	\begin{equation} \label{eq:av_M_solve}
		X_t^k = \other{\mathcal{P}}^k_{0, t}\bigl[X_0^k\bigr] + \int_0^t \other{\mathcal{P}}^k_{s, t}\bigl[ F_{[1,k],s}^\mathrm{av} \bigr]\mathrm{d}s.
	\end{equation}
	Hence, plugging \eqref{eq:av_M_bound_init} and \eqref{eq:av_M_forcing_bound} into \eqref{eq:av_M_solve}, and using \eqref{eq:av_prop_bound}, we conclude that, for all $\bm x \in \indset{N}^k$,
	\begin{equation} \label{eq:av_X_bound_nonsat}
		\bigl\lvert X_t^k(\bm x)\bigr\rvert \lesssim \ell_t\eta_t \, \mathfrak{s}^\mathrm{av}_{k,t}(\bm x) \biggl( 1 + (\log N)^{\other{C}_k}\int_0^t \frac{1}{\eta_s}\mathrm{d}s\biggr) \lesssim (\log N)^{\other{C}_k+1}\,\ell_t\eta_t \, \mathfrak{s}^\mathrm{av}_{k,t}(\bm x),
	\end{equation}
	where in the last step we used \eqref{eq:int_rules}. Therefore, \eqref{eq:M_bound_av} holds for all non-saturated $M_{[1,k]}$. 
	
	\vspace{5pt}
	
	\textbf{Case 2.} Assume that $M_{[1,k]}$ corresponds to a saturated chain, hence $\mathcal{A}_{j,t} = \Theta_t$ in \eqref{eq:av_M_evol} for all $j\in\indset{k}$. In this case $k$ is even and hence $k \ge 4$. Following the same logic as in the proof of Proposition \ref{prop:masters}, we assert that (c.f., Proposition~\ref{prop:lin_term})
	\begin{equation} \label{eq:X_lin_term}
		\bigl\lvert \Theta_t^{(j)} \bigl[X_t^k\bigr] (\bm x)\bigr\rvert \lesssim \frac{1}{\eta_t}(\log N)^{\other{C}_k+1}\,\ell_t\eta_t \, \mathfrak{s}^\mathrm{av}_{k,t}(\bm x), \quad t \in [0,T],
	\end{equation}
	for all $j \in \indset{k}$, where recall that $\Theta_t^{(j)} := I^{\otimes (j-1)}\otimes\Theta_t\otimes I^{\otimes (k-j)}$.
	Equipped with \eqref{eq:X_lin_term},  $\Theta_t^{\oplus k} =\sum_{j=1}^k \Theta_t^{(j)}$, 
	and treating the linear term $\Theta_t^{\oplus k} [X_t^k]$ in \eqref{eq:av_M_evol} as an additional forcing term
	with the same bound as in~\eqref{eq:av_M_forcing_bound}  we proceed as in the non-saturated case (with all propagators replaced by harmless multiples of identity) and conclude that  \eqref{eq:av_X_bound_nonsat} also holds for saturated $M_{[1,k]}$ (see Section~\ref{sec:masters_proof} for a similar argument in the case of $G-M$).
	
	Therefore, it remains to prove~\eqref{eq:X_lin_term}.
	To this end, we adopt the strategy from Section \ref{sec:reg_sec} used for long chains, which involves regularizing two observables. A single regularization would be insufficient because the averaged $M$-term's size is given by $(\ell_t\eta_t) \times \mathfrak{s}_{k,t}^\mathrm{av}$, and only a twice-regularized propagator yields a consistent bound for this term (see \eqref{eq:two_rings}). 
	 One regularization  suffices only for a target bound of 
	  $\mathfrak{s}_{k,t}^\mathrm{av}$ (see \eqref{eq:one_ring}), which is what we primarily used when estimating $G-M$.
	  
	It might seem counterintuitive that propagating the term 
	$(\ell_t\eta_t) \mathfrak{s}_{k,t}^\mathrm{av}$,  which appears larger due to the $\ell_t\eta_t$ prefactor, 
	requires a stronger improvement (i.e., more regularization) in the propagator. However, what matters is the singularity of the bound in the 
	small $\eta$ regime. The saturated propagator $\mathcal{P}^k_{s, t}$ acts essentially as a multiplication by
	$(\eta_s/\eta_t)^k$	   on constant functions, therefore it produces a singularity of order $\eta_t^{-k}$,
	while $(\ell_t\eta_t) \mathfrak{s}_{k,t}^\mathrm{av}$ is singular only as $\eta_t^{-k+1}$.
	Recall that each regularization improves the propagator singularity by a factor $\sqrt{\eta_t/\eta_s}$
	in terms of pure $\eta$-powers. 
	Thus, to make up for the missing $\eta_t$ factor, we need at least two regularizations.  
	
	We emphasize that the rationale for using two regularizations instead of one differs between the current
	$M$-bound proof and the $G-M$ proof in Section \ref{sec:reg_sec}.  
	For the $M$-bound,  as just explained,  two regularizations are fundamental to maintain the consistency of the estimate. 
	For $G-M$, strictly speaking, one regularization would maintain the consistency, but only for the optimal (in terms of $\ell\eta$) bound;  this is the approach for shorter chains ($k\le\maxK/2$) in Section \ref{sec:reg_sec}. 
	However, for longer chains the ideal input bound is not attainable due to the reduction mechanism, which incurs a cost. In the proof, this cost is accommodated by introducing loss exponents $\alpha_k, \beta_k$ effectively resulting in weaker bounds to propagate for these long chains. Consequently, estimates with target bounds like 
    $(\ell_t\eta_t)^{\beta_k} \mathfrak{s}_{k,t}^\mathrm{av}$ for $\beta_k>0$ inherently require more than one regularization.

	Returning to the proof of   \eqref{eq:X_lin_term}, 
	 we only consider the case $j=k$, since all other cases are structurally identical.  
	It follows from \eqref{eq:av_M_evol} and \eqref{eq:dTheta} that (c.f. \eqref{eq:long_lin_tem_evol}--\eqref{eq:remainder_long})
	\begin{equation} \label{eq:lin_X_evol}
		\begin{split}
			\frac{\mathrm{d}}{\mathrm{d}t}\biggl(\Theta_t^{(k)}\bigl[ X^{k}_t \bigr]\biggr) =&~ \biggl(I+\Theta_t^{(k)}+\frac{k}{2}I + \sum_{i\neq k-2}\Theta_t^{(i)}\biggr)\circ \Theta_t^{(k)}\bigl[ X^{k}_t \bigr]\\
			&+ \dring{X}^{k}_t + R_{[1,k],t}  + \Theta_t^{(k)}\bigl[ \mathcal{F}^\mathrm{av}_{[1,k],t}\bigr],
		\end{split}
	\end{equation}
	where the mollified quantity $\dring{X}^{k}_t$ is defined as
	\begin{equation}
		\dring{X}^{k}_t := \Theta_t^{(k-2)}\circ \Theta_t^{(k)}\biggl[\Tr\bigl[ M_{[1,k],t}\bigl(S^{x_1},\dots, S^{x_{k-3}}, \reg{S}^{x_{k-2}}, S^{x_{k-1}} \bigr) \reg{S}^{x_k}\bigr]\biggr],
	\end{equation}
	(recall the definition of the regularized $\reg{S}^{x_j}$ from~\eqref{eq:circ_above}--\eqref{eq:Sring_def}) 
	and the remainder term $R_{[1,k],t}$ is given by
	\begin{equation}
		\begin{split}
			R_{[1,k],t}(\bm x) :=&~  (\Theta_t)_{x_{k-1}x_k} \mathcal{T}_{k}\circ \Theta_t^{(k-2)} \bigl[X^{k}_t\bigr](\bm x') \\
			&+ (\Theta_t)_{x_{k-2}x_{k-3}} \mathcal{T}_{k-2}\circ \Theta_t^{(k)} \bigl[X^{k}_t\bigr](\bm x''',x_{k-1},x_k)\\ 
			&- (\Theta_t)_{x_{k-1}x_k}(\Theta_t)_{x_{k-2}x_{k-3}} \mathcal{T}_{k-2, k}\bigl[X^{k}_t\bigr](\bm x''', x_{k-1}).
		\end{split}
	\end{equation}
	Here we used the definition of the regularized observables \eqref{eq:circ_above} to deduce that, analogously to \eqref{eq:long_lin_term_decomp}, 
	\begin{equation} \label{eq:lin_M_decomp}
		\Theta_t^{(k-2)}\circ \Theta_t^{(k)}\bigl[ X^{k}_t\bigr] = \dring{X}^{k}_t + R_{[1,k],t}.
	\end{equation}
	Moreover, similarly to \eqref{eq:remainder_est}, using \eqref{eq:M_Ward}, \eqref{eq:triag} and \eqref{eq:true_convol}, together with \eqref{eq:M_bound_ind} for $2 \le k' \le k-1$, that
		holds for $2 \le k' \le k-1$, we obtain, for all $\bm x \in \indset{N}^k$,
	\begin{equation} \label{eq:M_residue_bound}
		\bigl\lvert R_{[1,k],t}(\bm x) \bigr\rvert  \lesssim \frac{1}{\eta_t^2} \,\bigl((\log N)^{C_{k-1}} + (\log N)^{C_{k-2}} \bigr) \times \ell_t\eta_t\, \mathfrak{s}_{k,t}^\mathrm{av}(\bm x), \quad t\in [0,T].
	\end{equation}
	
	Completely analogously to \eqref{eq:Gcirc_solve}, the mollified quantity $\dring{X}^{k}_t$  satisfies
	\begin{equation} \label{eq:Mring_solve}
			\dring{X}^{k}_t =  \mathrm{e}^{-kt/2}\mathcal{P}^{\otimes k}_{0,t}\circ U^{(k)}_{t}  \bigl[ 
			X_{0}^{k,\mathrm{av}}
			\bigr]  + \int_{0}^{t} \mathrm{e}^{k(s-t)/2}\mathcal{P}^{\otimes k}_{s,t} \circ U^{(k)}_{t} \bigl[ F^\mathrm{av}_{[1,k],s}+ \other{F}^\mathrm{av}_{[1,k],s} \bigr]  \mathrm{d}s,
	\end{equation}
	where $U^{(k)}_{t} := \reg{\Theta}_{t}^{(k-2)} \circ  \reg{\Theta}_{t}^{(k)}$ (note that we use two 
	regularizations for all $k$ not just for $k\ge \maxK/2$, c.f.~\eqref{eq:U_operator}),
	and  the additional deterministic forcing term $\other{F}^\mathrm{av}_{[1,k],t}$ is given by
	\begin{equation}
		\begin{split}
			\other{F}^\mathrm{av}_{[1,k],t}(\bm x) :=&~  
			(\Theta_t)_{x_{k-3} x_{k-2}} \bigl(\mathcal{T}_{k-2} \bigl[X^{k}_t\bigr](\bm x^{(4)}, x_{k-2}, x_{k-1}, x_k)  + \mathcal{T}_{k-2}\bigl[X^{k}_t\bigr](\bm x''', x_{k-1},x_k)\bigr)\\
			& +(\Theta_t)_{x_{k-1}x_k} \bigl(\mathcal{T}_{k}\bigl[X^{k}_t\bigr](\bm x'', x_k)  +  \mathcal{T}_{k}\bigl[X^{k}_t\bigr](\bm x')\bigr),
		\end{split}
	\end{equation}
	exactly as~\eqref{eq:extra_forcing2}. 
	Similarly to \eqref{eq:extra_forcing}, using \eqref{eq:M_Ward}, \eqref{eq:triag} and \eqref{eq:true_convol}, together with \eqref{eq:M_bound_ind} 
	for $2 \le k' \le k-1$, we deduce that 
	\begin{equation} \label{eq:M_extra_forcing}
		\bigl\lvert \other{F}^\mathrm{av}_{[1,k],t}(\bm x) \bigr\rvert \lesssim \frac{1}{\eta_t} \,(\log N)^{C_{k-1}}\times \ell_t\eta_t \, \mathfrak{s}^\mathrm{av}_{k,t}(\bm x), \quad t\in [0,T], \quad \bm x \in \indset{N}^k.
	\end{equation}
	Plugging  \eqref{eq:av_M_bound_init}, \eqref{eq:av_M_forcing_bound}, and \eqref{eq:M_extra_forcing} into \eqref{eq:Mring_solve}, and using the bound \eqref{eq:two_rings}, we conclude that
	\begin{equation} \label{eq:Mring_bound}
		\bigl\lvert \dring{X}^{k}_t \bigr\rvert \lesssim \frac{1}{\eta_t^2} \,(\log N)^{\other{C}_{k}+1}\times \ell_t\eta_t \, \mathfrak{s}^\mathrm{av}_{k,t}(\bm x), \quad t\in [0,T], \quad \bm x \in \indset{N}^k,
	\end{equation}
	where we assume that $\other{C}_k \ge C_{k-1} \vee C_{k-2}$.
	
	Hence, plugging \eqref{eq:M_residue_bound} and \eqref{eq:Mring_bound} back into \eqref{eq:lin_X_evol}, and using Duhamel's principle, we obtain \eqref{eq:X_lin_term}, using \eqref{eq:av_prop_bound} (see \eqref{eq:long_lin_tem_solve} and the discussion below for a more detailed version of the same argument). 
	This concludes the proof of \eqref{eq:M_bound_av}.
	\vspace{5pt}
	
	\textbf{Isotropic $M$-terms, proof of \eqref{eq:M_bound}.}\\
	Next, using \eqref{eq:M_bound_av} and the induction hypothesis~\eqref{eq:M_bound_ind} as an input, we analyze the isotopic quantities $Y_t^k=\bigl(M_{[1,k],t}(\bm x') \bigr)_{aa}$ and prove \eqref{eq:M_bound}.
	Differentiating $Y_t^k$ in time $t$, using \eqref{eq:dM}, we obtain, for any $a\in\indset{N}$,
	\begin{equation} \label{eq:iso_M_evol}
		\frac{\mathrm{d}}{\mathrm{d}t}Y_t^k = \biggl(\frac{k}{2}I + \bigoplus_{j=1}^{k-1} \mathcal{A}_{j,t}\biggr)\bigl[Y_t^k \bigr] + F_{[1,k],t}^\mathrm{iso}, 
	\end{equation}
	where the linear operators $\mathcal{A}_{j,t}$ are as defined in \eqref{eq:lin_prop_ops}, and the forcing term $F_{[1,k],t}^\mathrm{iso} \equiv F_{[1,k],t}^\mathrm{iso}(\bm x', a)$ is given by 
	\begin{equation} \label{eq:MFiso}
		F_{[1,k],t}^\mathrm{iso}(\bm x', a) := \sum_{\substack{1 \le i < j \le k \\ 2 \le j-i \le k-2}} \sum_q \bigl(M_{[i,j],t}\bigr)_{qq} \bigl(M_{[1,i],[j,k],t}^{(q)}\bigr)_{aa} + \sum_q m_{1,t}m_{k,t} \bigl(I+\mathcal{A}_{k,t}\bigr)_{aq} X_{t}^k(\bm x',q).
	\end{equation}
	Again, this formula is a simplified version of the analogous equation~\eqref{eq:k_iso_evol} for the isotropic
	$G-M$ bounds, lacking a martingale term, and with a greatly simplified forcing term. 
	Similarly to \eqref{eq:av_M_forcing_bound}, it follows from \eqref{eq:M_bound_ind} 
	 for $k'\in \indset{k-1}$ and \eqref{eq:M_bound_av} for $X_t^k$, that
	\begin{equation} \label{eq:iso_M_forcing_bound}
		\bigl\lvert F_{[1,k],t}^\mathrm{iso}(\bm x', a) \bigr\rvert \lesssim \frac{1}{\eta_t} (\log N)^{C_k} \sqrt{\ell_t\eta_t}\,\mathfrak{s}_{k,t}^\mathrm{iso}(a,\bm x', a), \quad t\in[0,T], \quad \bm x' \in \indset{N}^{k-1}, \quad a \in \indset{N}.
	\end{equation}
	Hence, using Duhamel's principle for \eqref{eq:iso_M_evol}, together with \eqref{eq:M_bound_init}, \eqref{eq:iso_M_forcing_bound}, and \eqref{eq:iso_prop_bound} from Lemma \ref{lemma:good_props}, we conclude that
	\begin{equation}
		\bigl\lvert Y_{t}^k(\bm x', a) \bigr\rvert \lesssim (\log N)^{C_k+1} \sqrt{\ell_t\eta_t}\,\mathfrak{s}_{k,t}^\mathrm{iso}(a,\bm x', a), \quad t\in[0,T], \quad \bm x' \in \indset{N}^{k-1}, \quad a \in \indset{N}.
	\end{equation}
	This concludes the proof of \eqref{eq:M_bound}, and hence that of Lemma \ref{lemma:M_bounds}.
\end{proof}

Next, we prove the auxiliary $M$-bound~\eqref{eq:resum_M_bound} in 
Lemma \ref{lemma:extra_M} using the estimates from Lemma \ref{lemma:M_bounds} as an input.
\begin{proof} [Proof of Lemma \ref{lemma:extra_M}]
	Let $M_{[j,j],t}^{(x_k)}$ be an $M$-term of the form \eqref{eq:special_M}. Without loss of generality, by relabeling the observables and the spectral parameters, we can assume that $j=1$, that is, we consider
	\begin{equation} \label{eq:special_M1}
		M_{[1,k+1],t} \equiv M_{[1,1],t}^{(x_k)} = M\bigl(z_{1,t}, S^{x_1}, z_{2,t}, \dots, z_{k,t}, S^{x_{k}}, z_{1,t}\bigr).
	\end{equation}	
	Since the time $t$ is fixed along the proof, we henceforth drop it from the subscript of $z_i$, $m_i$, $\mathfrak{s}_i$, $\ell$, $\eta$, and $M_{[i,j]}$. 
	
	It follows form \eqref{eq:otherM} and \eqref{eq:M_recursion} that 
	\begin{equation} \label{eq:special_M_recur}
		\bigl(M_{[1,k+1]}\bigr)_{qq} = \sum_a \biggl(\frac{m(z_1)}{1-m(z_1)^2S}\biggr)_{qa}\biggl((S^{x_1})_{aa}  \bigl(M_{[2,k+1]}\bigr)_{aa} +   \sum_{i=2}^{k}   \bigl(\mathcal{S}\bigl[M_{[1,i]}\bigr]\bigr)_{aa}\bigl(M_{[i,k+1]}\bigr)_{aa}\biggr),
	\end{equation}
	where we identify $M_{[i,k+1]} \equiv M_{[i,1]}^{(x_k)}$, and we used the fact that the last spectral parameter in \eqref{eq:special_M1} is equal to $z_1\equiv z_{1,t}$. It follows from the second bound in \eqref{eq:Ups_majorates} 
	and~\eqref{eq:Ups_norm_bounds_notime} noting that $\eta_0\ge 1$ in the definition of $\Upsilon_0$, that
	\begin{equation} \label{eq:stabXi_sums}
		\sum_{q} \biggl\lvert \biggl(\frac{m(z_1)}{1-m(z_1)^2S}\biggr)_{qa} \biggr\rvert \lesssim \frac{1}{|m(z_1)|} \sum_{q} \bigl\lvert \bigl(I + \Xi\bigr)_{qa} \bigr\rvert \lesssim 1.
	\end{equation}
	Hence, it follows from \eqref{eq:special_M_recur} and \eqref{eq:stabXi_sums}, that 
	\begin{equation} \label{eq:special_M_sum}
		\sum_q \bigl\lvert \bigl(M_{[1,k+1]}\bigr)_{qq} \bigr\rvert \lesssim \sum_a S_{x_1a}  \bigl\lvert\bigl(M_{[2,k+1]}\bigr)_{aa}\bigr\rvert +   \sum_{i=2}^{k} \sum_{ab}  S_{ab}\bigl\lvert \bigl( M_{[1,i]}\bigr)_{bb}\bigr\rvert \bigl\lvert\bigl(M_{[i,k+1]}\bigr)_{aa}\bigr\rvert.
	\end{equation}
	Estimating every $M$-term in \eqref{eq:special_M_sum} using \eqref{eq:M_bound}, and ignoring the poly-log factors, we obtain
	\begin{equation} \label{eq:recur_bound_resum}
		\begin{split}
			\sum_q \bigl\lvert \bigl(M_{[1,k+1]}\bigr)_{qq} \bigr\rvert \lesssim&~ \sqrt{\ell\eta}\sum_a  S_{x_1a}  \,\mathfrak{s}_{k}^\mathrm{iso}(a,x_2,\dots, x_k, a) \\
			&+ \ell\eta \sum_{i=2}^{k} \sum_{ab}  S_{ab}\,\mathfrak{s}_{i}^\mathrm{iso}(b,x_1,\dots, x_{i-1}, b) \,\mathfrak{s}_{k-i+2}^\mathrm{iso}(a,x_i,\dots, x_k, a)  \\
			\lesssim&~ \ell  \,\mathfrak{s}_{k}^\mathrm{av}(\bm x),
		\end{split}
	\end{equation}
	where in the second step we used \eqref{eq:SUps_comvol} to bound $S$, followed by \eqref{eq:Schwarz_convol}, similarly to \eqref{eq:av_insertion}. This concludes the proof of Lemma~\ref{lemma:extra_M}.
\end{proof}

Finally, we prove Claim \ref{claim:staticM}.
\begin{proof}[Proof of Claim \ref{claim:staticM}]
	We conduct the proof by induction in $n$, starting with the base case $n=1$.
	\vspace{5pt}
	
	\textbf{Base case.} 
	Note that for $n=1$ the bounds \eqref{eq:V_class_bound} and \eqref{eq:V_class_av} follow immediately from \eqref{eq:V1_class}, the first bound in \eqref{eq:SUps_comvol} and \eqref{eq:sumS=1}.
	To prove \eqref{eq:staticM} for $n=1$, we proceed by a nested induction in $k$. The base case of the nested induction $k = 1$ follow immediately from \eqref{eq:Ups_majorates}. 
	
	For the induction step at a fixed $k \ge 2$, we use \eqref{eq:M_recursion} to express $M_{[1,k+1]} \equiv M_{[1,k+1]}(A_1^j, \dots, A_{k}^j)$, which has length $k+1$, in terms of shorter $M$-terms of length at most $k$,
	\begin{equation}
		M_{[1,k+1]}  = \frac{m(z_{1})}{1-m(z_1)m(z_{k+1})S}\biggl[A_1^j   M_{[2,k+1]} + \sum_{i=2}^{k}  \mathcal{S}\bigl[ M_{[1,i]}\bigr] M_{[i,k+1]}\biggr].
	\end{equation}
	We estimate the stability operator $(1-m(z_1)m(z_{k+1})S)^{-1}$ using \eqref{eq:Ups_majorates};
	note that at time $t=0$ there is no difference in the estimate for $\Theta_0$ and $\Xi_0$.
	Bounding every $M$-term by the induction hypothesis, and using \eqref{eq:true_convol} to perform the summations, similarly to \eqref{eq:recur_bound_resum}, we obtain
	\begin{equation} \label{eq:expand_resum}
		\begin{split}
			\bigl\lvert\bigl(M_{[1,k+1]}\bigr)_{aa}\bigr\rvert \lesssim&~ \sqrt{W}\sum_q \Upsilon_{aq}  S_{qx_1}\, \mathfrak{s}_{k}^\mathrm{iso}(q,x_1,x_2,\dots,x_{k}, q)\\
			&+ W \sum_{i=2}^{k} \sum_{qc} \Upsilon_{aq} S_{qc} \,\mathfrak{s}_i^\mathrm{iso}(c,x_1,\dots,x_{i-1},c) \mathfrak{s}_{k-i+2}^\mathrm{iso}(q,x_{i},\dots,x_k,q)\\
			\lesssim&~ \sqrt{W} \mathfrak{s}_{k+1}^\mathrm{iso}(a,\bm x, a),
		\end{split}
	\end{equation} 
	where we recall that $\Upsilon = \Upsilon_0$, $\ell_0 \sim W$, and $\eta_0 \sim 1$. Here, in the last step also we used \eqref{eq:triag}. Since the matrices $M$ are diagonal, $(M_{[1,k+1]})_{ab}$ is trivially zero for $b\neq a$.
	
	\vspace{5pt}
	
	\textbf{Induction step.} Fix $n \ge 2$, and assume that \eqref{eq:V_class_bound}--\eqref{eq:staticM} hold for all $n' \in \indset{n-1}$. Then the bounds \eqref{eq:V_class_bound} and \eqref{eq:V_class_av} are simple consequences of the definition \eqref{eq:Vk_class}, \eqref{eq:triag} and \eqref{eq:true_convol}. To prove \eqref{eq:staticM}, we proceed exactly as in the base case $n=1$ and argue by a nested induction in $k$. Base case $k=1$ is a trivial consequence of \eqref{eq:Ups_majorates} and \eqref{eq:true_convol}. The induction step in $k$ is proved analogously to \eqref{eq:expand_resum} by expanding $M_{[1,k+1]}$ in terms of shorter $M$-quantities and summing up the corresponding bounds, available by induction hypothesis. This concludes the proof of Claim \ref{claim:staticM}.
\end{proof}

\section{Local laws with traceless observables}\label{sec:traceless}

In this section we prove local laws with improved error bounds in case of the  traceless version of the $S^x$ observables
for spectral parameters $z$ with $\eta=|\im z|$ below the critical scale $(W/N)^2$. In this regime, the resolvent exhibits no spatial decay,  and the behavior of its fluctuations is governed by the leading tracial mode. Hence, tracelessness  of the observables results in a cancellation effect. 

We will see that each traceless observable gains an additional factor $(N/W)\sqrt{\eta}\lesssim 1$, 
corresponding to the {\it $\sqrt{\eta}$-rule} for the standard Wigner matrices, $W=N$ and $S_{ij}=1/N$, first 
found in \cite{cipolloni2021eigenstate} in a special case and proven in full 
generality  in \cite{Cipolloni2022Optimal}. We remark that the gain from tracelessness is a spectral gap effect, and the improvement factor is proportional to the square root of the ratio of the two largest eigenvalues of the stability operator $(I - |m(z)|^2 S)^{-1}$.

First, we focus on the traceless version of our standard observable of the form $S^x$, and in a later Section~\ref{sec:gen_tr} we explain how to generalize the result to any diagonal observable.  

We again prove the results for the simpler case when all spectral parameters are the same, i.e. $z_{j, t}\in \{ z_t, \bar z_t\}$
for a fixed solution $z_t$ of the characteristic equation. The extension to the general case, when $\bm z_t \in \dom_t^k$, requires minor modifications as described at the end of
Section~\ref{sec:masters_sec} and will be omitted.

Analogously  to the  (mean-field) regularization in the Wigner case, for
 any $A\in \mathbb{C}^{N\times N}$ matrix, let
$$ 
\trless{A} := A - N^{-1}\Tr[A] I
$$
denote its traceless component. In particular\footnote{We stress that the (spatial) regularization $\reg{S}$ defined in
the main proof in~\eqref{eq:circ_above} and the current mean-field or traceless regularization $\trless{S}$
 are very different objects. The latter is used  only here in Section~\ref{sec:traceless} and in the regime $\eta\le (W/N)^2$.}, 
\begin{equation} \label{eq:trless_S}
	\trless{S}^x := S^x - N^{-1} I, \quad x \in \indset{N},
\end{equation}
where we used \eqref{eq:sumS=1} to deduce that $\Tr[S^x] = 1$.

In the present section we only consider target spectral parameters $z$ satisfying $\eta=|\im z| \le (W/N)^2$ and 
$\eta \ge N^{-1+\etaexp}$. 
Let $\crit \in [0,T]$ be the critical time along the flow \eqref{eq:char_flow}, such that 
\begin{equation} \label{eq:tcrit_def}
	\eta_\crit = (W/N)^2.
\end{equation}
Recall that $\eta_t := |\im z_t|$, where $z_t$ is the trajectory of \eqref{eq:char_flow} satisfying $z_T = z$ at a fixed final time $T$. By definition \eqref{def:ell_def}, $\ell_t =\ell(\eta_t) \sim N$ for all $t \in [\crit, T]$, hence $\ell_t \eta_t \sim N\eta_t$. In particular, for all $t \in [\crit, T]$ the admissible size-functions $\Upsilon_t$ satisfy $(\Upsilon_t)_{xy} \sim (N\eta_t)^{-1}$, see~\eqref{eq:Upsilon_deloc_notime}. 
Therefore, after the critical time $\crit$, the spatial decay is completely lost, hence the size functions $\mathfrak{s}^\mathrm{av/iso}_{k,t}$, defined in \eqref{eq:sfunc_def}, can be replaced by powers of $(N\eta_t)^{-1}$, up to 
irrelevant multiplicative constants, i.e.
\begin{equation}\label{eq:simplesizefn}
	\mathfrak{s}^{\mathrm{iso}}_{k,t}(a, \bm x', b) \sim (N\eta_t)^{-k + 1/2}, \quad \mathfrak{s}^{\mathrm{av}}_{k,t}(\bm x) \sim (N\eta_t)^{-k}, \quad t\in [\crit, T],
\end{equation}
for all $a,b\in\indset{N}$, $\bm x \in \indset{N}^k$, and $k \in \indset{K}$.

Recall the notation  \eqref{eq:resolvent_chains} for resolvent chains with general observables 
and \eqref{eq:Gk_def} for special observables  $S^x$.
We now consider resolvent chains similar to  \eqref{eq:Gk_def}, where some observables $S^x$ are 
replaced by their traceless version $\trless{S}^{x}$. In the following notation, the set $\trSet'$  encodes
the location of these traceless observables. 
Let $k \in \mathbb{N}$, for any integers, $j\le l \in \indset{k}$, spectral parameters $\bm z_t \in \{z_t, \overline{z}_t\}^k$, and a subset $\trSet' \subset \indset{k-1}$, we define 
\begin{equation} \label{eq:tracelessGk}
	G_{[j,l],t}(\trSet') \equiv G_{[j,l],t}(\trSet'; x_j,\dots, x_{l-1}) := G_{[j,l],t}\bigl(\bm z_t; \other{A}_j, \dots, \other{A}_{l-1}  \bigr),
\end{equation}
where, given $\trSet'$ and $\bm x$, the observables $\{\other{A}_i\}_{i=1}^{k-1}$ are defined as
\begin{equation} \label{eq:trlessAs}
	\other{A}_i \equiv \other{A}_i(\trSet', x_i) := \begin{cases}
		\trless{S}^{x_i}, ~ &i\in\trSet',\\
		S^{x_i}, ~ &i\notin\trSet'.
	\end{cases} 
\end{equation}
We drop the explicit dependence of $G_{[j,l],t}(\trSet'; x_j,\dots, x_{l-1})$ on the external indices $x_j,\dots, x_{l-1}$ since it is essentially irrelevant for $t \ge \crit$ . We only treat $G_{[j,l],t}(\trSet')$ as a function of $\bm x$ when it is acted on by a linear operator $\mathcal{A}$ on $\mathbb{C}^{N\otimes k}$ in the sense of \eqref{eq:tensor_action}.
Recall also that $G_{[j,l],t}$ depends implicitly on the spectral parameters $z_{j,t},\dots z_{l,t}$ as indicated by the subscript ${[j,l]}$. Similarly to \eqref{eq:Mk_def}, we define the corresponding deterministic approximations
\begin{equation} \label{eq:tracelessMk}
	M_{[j,l],t}(\trSet') \equiv M_{[j,l],t}(\trSet'; x_j,\dots, x_{l-1}) := M_{[j,l],t}\bigl(\bm z_t; \other{A}_j, \dots, \other{A}_{l-1}  \bigr),
\end{equation}
where the observables $\other{A}_i$ are defined in \eqref{eq:trlessAs}. Deterministic approximations $M_{[j,l],t}(\trSet')$ satisfy the following bounds.
\begin{lemma} [Traceless $M$ bounds] \label{lemma:traceless_M_bounds}
	Let $t \in [\crit, T]$, where $\crit$ is the critical time defined by \eqref{eq:tcrit_def}, and let $k \in \mathbb{N}$ be fixed.  Then, for any (potentially empty) subset $\trSet' \subset \indset{k-1}$, the deterministic approximations $M_{[1,k],t}(\trSet')$, defined in \eqref{eq:Mt_def}, satisfy the bounds
	\begin{equation} \label{eq:trlessM_iso}
		\bigl\lvert \bigl(M_{[1,k],t}(\trSet')\bigr)_{ab}\bigr\rvert \lesssim \delta_{ab}\frac{1}{(N\eta_t)^{k-1}} \biggl(\frac{N^2\eta_t}{W^2}\biggr)^{\lceil n/2 \rceil}, \quad \text{where} \quad  n :=  |\trSet'|,\\
	\end{equation}
	uniformly in $\bm x \in \indset{N}^k$, $\bm z_t \in \{z_t, \overline{z}_t\}^k$, and $a,b \in \indset{N}$. Moreover, let  $\trSet_k \subset \{k\}$ be either an empty set or a singleton consisting of $k$, and set
	\begin{equation} \label{eq:trlessA_k}
		\other{A}_k \equiv \other{A}_k(\trSet_k, x_k) := \begin{cases}
			\trless{S}^{x_k}, \quad &k\in\trSet_k,\\
			S^{x_k}, \quad &k\notin\trSet_k.
		\end{cases} 
	\end{equation}
	Then, $M_{[1,k],t}(\trSet')$ satisfies the averaged bound
	\begin{equation} \label{eq:trlessM_av}
		\bigl\lvert \Tr\bigl[M_{[1,k],t}(\trSet')\other{A}_k\bigr] \bigr\rvert \lesssim \frac{1}{(N\eta_t)^{k-1}} \biggl(\frac{N^2\eta_t}{W^2}\biggr)^{\lceil n/2 \rceil}, \quad \text{where} \quad  n :=  |\trSet|, \quad \trSet := \trSet'\cup\trSet_k,
	\end{equation}
	uniformly in $\bm x \in \indset{N}^k$ and $\bm z_t \in \{z_t, \overline{z}_t\}^k$.
\end{lemma}
We prove Lemma~\ref{lemma:traceless_M_bounds} in Section~\ref{sec:trless_M}.
Note that the parameter $n$ counts the total number of traceless observables appearing on the left-hand side of the corresponding equation, including the test matrix $\other{A}_k$ in \eqref{eq:trlessM_av}. 
Both \eqref{eq:trlessM_iso}, \eqref{eq:trlessM_av}
exhibit an improvement over their counterparts \eqref{eq:M_bound}, \eqref{eq:M_bound_av} essentially by a factor of
$(N/W)\sqrt{\eta_t}$ for each traceless observable\footnote{More precisely, the parity of $n$ matters, even number
	$n=2m$ of traceless observables yield the same improvement as $2m-1$ traceless observables; this is expressed by
	the upper integer part $\lceil n/2 \rceil$ in the exponents.}. In the standard Wigner case ($W=N$ and $S_{ij} = N^{-1}$)
these estimates were obtained in Lemma 2.4 of \cite{Cipolloni2022Optimal}.  
The bounds \eqref{eq:trlessM_iso} and \eqref{eq:trlessM_av} are generally sharp when the corresponding resolvent chains are saturated.

Now we can formulate the local laws with traceless observables.  
\begin{theorem} [Traceless Local Laws] \label{th:traceless_laws} Fix an integer $\maxK_{\max} \in \mathbb{N}$.
	Let $t \in [\crit, T]$, where $\crit$ is the critical time defined by \eqref{eq:tcrit_def}, and let $k \in \indset{\maxK_{\max}}$. Then, for any $\bm z_t \in \{z_t, \overline{z}_t\}^k$, $\bm x \in \indset{N}^k$, and $\trSet' \subset \indset{k-1}$, the resolvent chain~$G_{[1,k],t}(\trSet')$, defined in \eqref{eq:tracelessGk}, satisfies the isotropic local law,
	\begin{equation} \label{eq:traceless_Gk_iso}
		\biggl\lvert \bigl((G-M)_{[1,k],t}(\trSet')\bigr)_{ab}\biggr\rvert \prec \frac{1}{\sqrt{N\eta_t}}\frac{1}{(N\eta_t)^{k-1}}\biggl(\frac{N^2\eta_t}{W^2}\biggr)^{n/2}, \qquad  n:= |\trSet'|,
	\end{equation}
	uniformly in $a,b \in \indset{N}$. Moreover, with $\trSet_k \subset \{k\}$ and $\other{A}_k$ from \eqref{eq:trlessA_k}, the resolvent chain~$G_{[1,k],t}(\trSet')$ satisfies the averaged local law,
	\begin{equation} \label{eq:traceless_Gk_av}
		\biggl\lvert \Tr\bigl[(G-M)_{[1,k],t}(\trSet')\other{A}_k\bigr] \biggr\rvert \prec \frac{1}{N\eta} \frac{1}{(N\eta_t)^{k-1}} \biggl(\frac{N^2\eta_t}{W^2}\biggr)^{n/2}, \quad \text{where $n:= |\trSet|$,} \quad  \trSet := \trSet'\cup\trSet_k.
	\end{equation}
\end{theorem}
Note that, compared with Theorem \ref{th:local_laws},
each traceless observable improves the error terms in \eqref{eq:traceless_Gk_iso} and \eqref{eq:traceless_Gk_av} by a small factor of $N\sqrt{\eta_t}/W$. Note that this improvement is insensitive to the parity of $n$ unlike the
bounds on the leading $M$-term in Lemma~\ref{lemma:traceless_M_bounds}.
Moreover,  the error terms in \eqref{eq:traceless_Gk_iso} and \eqref{eq:traceless_Gk_av} are smaller by a factor of $(N\eta_t)^{-1/2}$ and $(N\eta_t)^{-1}$, respectively, then the sizes of the respective deterministic terms in \eqref{eq:trlessM_iso} and \eqref{eq:trlessM_av}, at least when the number of traceless matrices $n$ is even.   This corresponds precisely to the so-called $\sqrt{\eta}$-rule in the mean-field case ($W \sim N$), see Section 1.1 of \cite{Cipolloni2022Optimal} for a more detailed discussion of the mean-field case.

We defer the proof of Theorem \ref{th:traceless_laws} to Section \ref{sec:trlaws}.

\begin{proof} [Proof of Theorem \ref{th:local_laws_traceless}]
We deduce that traceless local laws of Theorem \ref{th:local_laws_traceless} for random band matrices with general distribution from Theorem \ref{th:traceless_laws} by removing the Gaussian component introduced to $H_t$ by the OU process \eqref{eq:zigOU} and conclude Theorem~\ref{th:local_laws_traceless}.
To this end, we employ a Green function comparison argument, analogous to that in Sections~\ref{sec:k_iso_rep}--\ref{sec:k_av_rep}. 
We remark that the crucial traceless observables that carry the
$(N\sqrt{\eta}/W)$-improvement are always surrounded by resolvents on both sides and 
this structure is preserved when the surrounding resolvent chain is differentiated
in the expansion analogous to, e.g.~\eqref{eq:resolvent_expand}; 
schematically $$\partial_{ij} (GAG)  = -GE^{ij}GAG + GAGE^{ij}G.$$  Moreover, for $\eta$ below the critical $(W/N)^2$ scale, the spatial decay of the resolvent is absent, which reduces the proof of GFT for traceless observables to
a routine  power counting. We leave the technical details to the reader. This concludes the proof of Theorem \ref{th:local_laws_traceless}. 
\end{proof}

We dedicate the rest of the section to proving the local laws (Theorem~\ref{th:traceless_laws}) for resolvent chains with special traceless observables $\trless{S}^x$, defined in \eqref{eq:trless_S}, and the bounds on the corresponding deterministic approximations in Lemma~\ref{lemma:traceless_M_bounds}. 
First, we present the more involved proof of  Theorem~\ref{th:traceless_laws} in Section \ref{sec:trlaws}, then we prove the supporting technical statements in Sections~\ref{sec:tr_bootstrap}--\ref{sec:Q_bound}. We conclude this section by establishing the deterministic $M$-bounds with traceless observables in Section~\ref{sec:trless_M}, closely following the proof of Lemma~\ref{lemma:M_bounds} form Section \ref{sec:M_bounds}.

\subsection{Proof of Theorem~\ref{th:traceless_laws}} \label{sec:trlaws}
We begin by preparing some notation.
For all integers $k \ge 1$ and $0 \le n \le k-1$, consider the class of resolvent chains $\mathfrak{G}_{k,n}$, defined as
\begin{equation}
	\mathfrak{G}_{k,n} := \bigl\{ G_{[1,k],t}(\trSet'; \bm x) \, :\, |\trSet'| = n, ~ \bm z_t \in \{z_t, \overline{z}_t\}^k, ~ \bm x \in \indset{N}^{k-1} \bigr\}, 
\end{equation} 
where $\trSet' \subset \indset{k-1}$, and $ G_{[1,k],t}(\trSet'; \,\cdot\,) $ is defined in \eqref{eq:tracelessGk}. That is, the class $\mathfrak{G}_{k,n}$ consists of $k$-chains with $n$ traceless observables. Recall that the spectral parameters $\bm z_t$ and the external indices $\bm x$ are essentially irrelevant and we will drop them from the notations. 
Strictly speaking $\mathfrak{G}_{k,n}$ also depends on time but we also ignore this, and only indicate the time-dependence of the individual chains.

Let $\trl_t$ denote the improvement factor coming from each traceless observable in the error term of \eqref{eq:traceless_Gk_iso}--\eqref{eq:traceless_Gk_av}, that is,
\begin{equation} \label{eq:theta_def}
	\trl_t := \frac{N\sqrt{\eta_t}}{W}.
\end{equation}
Recall that $t\ge t_*$ implies $\theta_t\le 1$.
Since $W \le N$ and $N\eta_t \ge N^{\etaexp}$ by \eqref{eq:etabound}, there exists an integer $L \le \lceil 1/\etaexp\rceil$, such that 
\begin{equation} \label{eq:theta_lower}
	\trl_t \ge \frac{1}{(N\eta_t)^{L/2}}, \quad t \in [\crit, T].
\end{equation}

For all integers $k \ge 1$ and $0 \le n \le k-1$, we define the control quantities
\begin{equation} \label{eq:traceless_Psi_def}
	\begin{split}
		\Psi_{k,t}^{n,\mathrm{iso}} &:= \max_{\mathcal{G}_t \in \mathfrak{G}_{k,n}}  \frac{(N\eta_t)^{k-1/2}}{\trl_t^{n}}\norm{\mathcal{G}_t - M_{\mathcal{G},t}}_{\max}, \\
		 \Psi_{k,t}^{n,\mathrm{av}} &:= \max_{x\in\indset{N}} \max_{\mathcal{G}_t \in \mathfrak{G}_{k,n}}   \frac{(N\eta_t)^{k }}{\trl_t^{n}}\bigl\lvert \Tr\bigl[(\mathcal{G}_t - M_{\mathcal{G},t})S^x\bigr]\bigr\rvert,
	\end{split}
\end{equation}
where $M_{\mathcal{G},t}$ denotes the deterministic approximation corresponding to the 
chain $\mathcal{G}_t \in \mathfrak{G}_{k,n}$, defined in \eqref{eq:tracelessMk}. 
Moreover, in the special case $n=k$, we define
\begin{equation} \label{eq:traceless_Psi_av}
	\Psi_{k,t}^{k,\mathrm{av}} := \max_{x\in\indset{N}} \max_{\mathcal{G}_t \in \mathfrak{G}_{k,k-1}}   \frac{(N\eta_t)^{k }}{\trl_t^{k}}\bigl\lvert \Tr\bigl[(\mathcal{G}_t - M_{\mathcal{G},t})\trless{S}^x\bigr]\bigr\rvert.
\end{equation} 
These control quantities are the traceless analogues of $\Psi_{k,t}^{\mathrm{iso/av}}$, defined
in \eqref{eq:Psi_def}, taking into account the optimal number of improvement factors $\trl_t$. 
 Recall that since $t\ge \crit$, 
the size functions $\mathfrak{s}_{k,t}^{\mathrm{av/iso}}$ simplify to pure $(N\eta_t)$-powers according to~\eqref{eq:simplesizefn}.  
Note that no loss exponents $\alpha_k, \beta_k$ are present in \eqref{eq:traceless_Psi_av}.

Since $\trl_{\crit} = 1$ by \eqref{eq:tcrit_def} and \eqref{eq:theta_def}, Theorem \ref{th:local_laws} and Proposition \ref{prop:zig} imply that 
\begin{equation} \label{eq:trless_initial}
	\max_{n \in \indset{k-1}}\Psi_{k,\crit}^{n,\mathrm{iso}} \prec 1, \quad  \max_{n \in \indset{k}}\Psi_{k,\crit}^{n,\mathrm{av}} \prec 1, \quad \max_{t \in [\crit, T]}\Psi_{k,t}^{0,\mathrm{iso/av}} \prec 1, 
\end{equation}
for any fixed $k \in \mathbb{N}$.  The information given in~\eqref{eq:trless_initial} is assumed throughout 
Section~\ref{sec:trlaws}--\ref{sec:tr_bootstrap}. 
Our goal is to show that $\max_{n \in \indset{k-1}}\Psi_{k,T}^{n,\mathrm{iso}} \prec 1$ 
and $\max_{n \in \indset{k}}\Psi_{k,T}^{n,\mathrm{av}} \prec 1$ for any fixed $k \in \mathbb{N}$.
To this end, we proceed by induction in the number of traceless observables $n$, using the following proposition.
Note that as $n$ increases, the maximal chain length decreases;  i.e chain length is traded in for a higher $n$.
In fact, a similar mechanism was formalized  in a more subtle way earlier in Section~\ref{sec:masters_sec}: the hierarchy of master inequalities in Proposition~\ref{prop:masters} were closable
exactly because $\Psi_{k,t}^{\mathrm{iso/av}}$ encoded weaker control on longer chains
with the carefully designed  loss exponents
in the definitions~\eqref{eq:Psi_def}.

\begin{prop}  [Traceless Induction]\label{prop:n_induction} Let $n \ge 1$ be a fixed integer, 
and let $L$ be the integer such that \eqref{eq:theta_lower} holds.
	Assume that for some fixed $\maxK' \ge 2L + n+2$, the quantities $\Psi_{k,t}^{n',\mathrm{iso/av}}$, defined in \eqref{eq:traceless_Psi_def}, satisfy
	\begin{equation} \label{eq:trless_ind_assume}
		\max_{t\in[\crit, T]} \max_{n' \in \indset{0,n-1}} \max_{k \in \indset{n'+1, \maxK'}} \Psi_{k,t}^{n',\mathrm{iso/av}} \prec 1,
	\end{equation}
	i.e. we have optimal traceless law for  all chains of length at most $K'$
	with at most $n-1$ traceless observables\footnote{Note that $k\ge n'+1$ in~\eqref{eq:trless_ind_assume},
	i.e. no assumption is made for the averaged $k$-chains with $n'=k$ traceless observables (the relation $k\ge n'+1$
	is automatic for isotropic chains).}.  
	Then, the quantities $\Psi_{k,t}^{n,\mathrm{iso/av}}$, $k\ge n+1$, and  $\Psi_{n,t}^{n,\mathrm{av}}$
	from   \eqref{eq:traceless_Psi_av} satisfy
	\begin{equation} \label{eq:trless_ind_conclude}
		\max_{t\in[\crit, T]} \max_{k \in \indset{n+1, \maxK'-2L}} \Psi_{k,t}^{n,\mathrm{iso/av}} \prec 1, \quad \max_{t\in[\crit, T]} \Psi_{n,t}^{n,\mathrm{av}} \prec 1,
	\end{equation} 
	  i.e.,
	we have optimal traceless law for all chains of length at most $K'-2L$ with  $n$ traceless observables,
	including the averaged $n$-chains with maximal number of traceless observables. 
\end{prop}
We defer the proof of Proposition \ref{prop:n_induction} until the end of this section. Armed with Proposition \ref{prop:n_induction}, we are ready to prove Theorem \ref{th:traceless_laws}.
\begin{proof}[Proof of Theorem \ref{th:traceless_laws}]
	Recall that $\maxK_{\max}$ is the maximal length of the resolvent chain under consideration in~\eqref{eq:traceless_Gk_iso}--\eqref{eq:traceless_Gk_av}, and set $n_{\max} := \maxK_{\max}$ to be the maximal number of traceless observables
	in  these formulas. Define the sequence $K_j := K_{\max}+ 1 + 2(n_{\max} -j)L$, for $j \in \indset{0, n_{\max}}$, where $L$ is the integer from~\eqref{eq:theta_lower}. 
	
	Note that by the last relation in \eqref{eq:trless_initial}, the assumption \eqref{eq:trless_ind_assume} of Proposition \ref{prop:n_induction} holds 
	with $n=1$ for $K' := K_0$. 
	Using Proposition \ref{prop:n_induction} for these values, from the first formula 
	in \eqref{eq:trless_ind_conclude}, 
	it follows that the assumption \eqref{eq:trless_ind_assume}  
	also holds $n = 2$ and $K'=K_0-2L=K_1$. 
	Repeating this argument $n_{\max}$ times, 
	by induction in $j \in \indset{n_{\max}}$ using Proposition \ref{prop:n_induction},  it follows that \eqref{eq:trless_ind_conclude} holds with $n = n_{\max}$ and 
	for all chains of length $K_{n_{\max}-1}-2L = K_{\max}+2$. 
	This concludes the proof of Theorem \ref{th:traceless_laws}.
\end{proof}

Therefore, it remains to prove Proposition \ref{prop:n_induction}. To this end, we proceed by another 
induction, called {\it traceless bootstrap},  nested  inside the main induction in $n$.  In each step of this bootstrap
we gradually improve the error estimate on $\Psi_{k,t}^{n,\mathrm{iso/av}}$
by a small factor $(N\eta_t)^{-1/2}$  at the cost of decreasing the maximal chain length by two.
A new  parameter $l \in \indset{1,L}$
will track the actual gain expressed as the $l$-power of  $(N\eta_t)^{-1/2}$, we therefore also call the
traceless bootstrap as {\it $l$-induction}  to distinguish from the main {\it $n$-induction}
 step in Proposition~\ref{prop:n_induction}.
One step of the traceless bootstrap is formalized in the following proposition.  
\begin{prop} [Traceless Bootstrap] \label{prop:trless_bootstrap}
	Let $n$, $L$, and $\maxK' \ge 2L+n+2$ be the integers from Proposition \ref{prop:n_induction}, and assume that the condition \eqref{eq:trless_ind_assume} holds. Fix an $l \in \indset{L}$, and assume that
	\begin{equation} \label{eq:power_induction_assume}
		\max_{t\in[\crit, T]} \max_{k \in \indset{n+1, \maxK'-2l+2}} \min\bigl\{1, \trl_t(N\eta_t)^{(l-1)/2}\}\Psi_{k,t}^{n,\mathrm{av/iso}} \prec 1.
	\end{equation}

	Then, the quantities $\Psi_{k,t}^{n,\mathrm{av/iso}}$, defined in \eqref{eq:traceless_Psi_def}, also satisfy
	\begin{equation} \label{eq:power_induction_conclude}
		\max_{t\in[\crit, T]} \max_{k \in \indset{n+1, \maxK'-2l}} \min\bigl\{1, \trl_t(N\eta_t)^{l/2}\}\Psi_{k,t}^{n,\mathrm{av/iso}} \prec 1, 
		\quad \max_{t\in[\crit, T]} \min\bigl\{1, \trl_t(N\eta_t)^{l/2}\}\Psi_{n,t}^{n,\mathrm{av}} \prec 1.
	\end{equation}
\end{prop} 
We prove Proposition \ref{prop:trless_bootstrap} in Section \ref{sec:tr_bootstrap}.
\begin{proof}[Proof of Proposition \ref{prop:n_induction}]
	The assumption \eqref{eq:trless_ind_assume} implies that \eqref{eq:power_induction_assume} holds with the same $n$ for $l=1$, that is
	\begin{equation} \label{eq:n_trless_init}
		\max_{t\in[\crit, T]} \max_{k \in \indset{n+1, \maxK'}} \min\bigl\{1, \trl_t\}\Psi_{k,t}^{n,\mathrm{av/iso}} \prec 1.
	\end{equation}
	Indeed, by \eqref{eq:sumS=1}, any traceless observable $\trless{S}^x$ can be written as
	\begin{equation} \label{eq:trless_triv_decomp}
		\trless{S}^x= S^x -\frac{1}{N}\sum_y S^y, \quad x\in\indset{N}.
	\end{equation} 
	By \eqref{eq:trless_triv_decomp}, any chain $\mathcal{G}_{t} \in \mathfrak{G}_{k,n}$ can be expressed as a linear combination of chains in $\mathfrak{G}_{k,n-1}$ with coefficients bounded in $\ell_1$-norm, i.e. one may disregard the
	tracelessness of any observable and estimate the corresponding chain as if it had one less traceless observable. 
	Therefore, for any $k \ge n$, we have the relation
	\begin{equation}\label{eq:triv}
		\trl_t \Psi_{k,t}^{n,\mathrm{iso/av}}\lesssim \Psi_{k,t}^{n-1,\mathrm{iso/av}},
	\end{equation}
	from which  and \eqref{eq:trless_ind_assume} the conclusion \eqref{eq:n_trless_init} follows immediately.
	This trivial step, however, loses a $\trl_t$ factor that we will gain back iterating Proposition \ref{prop:trless_bootstrap}.
	
	Using \eqref{eq:n_trless_init} and Proposition \ref{prop:trless_bootstrap} for $l=1$, from the 
	first conclusion of~\eqref{eq:power_induction_conclude}
	 we deduce that \eqref{eq:power_induction_assume} also holds for $l=2$. Proceeding by induction in $l \in \indset{L}$, we deduce that the bounds \eqref{eq:power_induction_conclude} hold for $l=L$.
	It follows from \eqref{eq:theta_lower} that
	\begin{equation}
		\min\bigl\{1, \trl_t(N\eta_t)^{L/2}\} = 1.
	\end{equation}
	Hence, we obtain \eqref{eq:trless_ind_conclude} from \eqref{eq:power_induction_conclude} with $l=L$. This concludes the proof of Proposition~\ref{prop:n_induction}.
	 
\end{proof}

\subsection{Evolution equations for chains with traceless observables: Proof of Proposition~\ref{prop:trless_bootstrap}} \label{sec:tr_bootstrap}
The proof of Proposition \ref{prop:trless_bootstrap} relies on a stopping time argument and follows the same general outline as the proof of the general master inequalities in Proposition~\ref{prop:masters}. That is, we derive evolution equations for the averaged and isotropic $G-M$ quantities, apply Duhamel's principle, and bound all the emergent terms using the definition of the stopping time. In fact, since along the entire argument $t \ge \crit$, the decay of the resolvent entries becomes irrelevant, which greatly simplifies the proof.  On the other hand, the double induction makes the argument more involved,
since we have two different inputs~\eqref{eq:trless_ind_assume}  and~\eqref{eq:power_induction_assume}. 
The former is stronger as it asserts the optimal bound, but it is valid  only for chains with strictly less than $n$ traceless
observables. The latter is suboptimal due to the small factor $\min\bigl\{1, \trl_t(N\eta_t)^{(l-1)/2}\}$, 
but it is only valid for $n$ traceless observables. The goal is to improve 
this latter bound by a factor $1/\sqrt{N\eta_t}$. 

The underlying mechanism for the this improvement is based on the fact that the size of the $G_k-M_k$ fluctuations is dictated by the deterministic approximation of the $2k$-chain arising from the corresponding quadratic variation, while the fluctuations of the  $2k$-chain are suppressed by an additional $1/(N\eta)$-power. Symbolically, 
\begin{equation*}
	\bigl\lvert G_k-M_k \bigr\rvert \prec \sqrt{M_{2k}} + \sqrt{ G_{2k}-M_{2k} } \prec \frac{M_k}{\sqrt{N\eta}} \biggl(1 + \frac{1}{\sqrt{N\eta}}\biggr).
\end{equation*}
However, now we need to ensure that traceless observables are preserved along the whole procedure.

\begin{proof} [Proof of Proposition~\ref{prop:trless_bootstrap}]
	Recall the evolution equations for the averaged and  isotropic resolvent chains from \eqref{eq:k_av_evol} and \eqref{eq:k_iso_evol}, respectively. Let $G_{[1,k],t}$ be a resolvent chain with $n$ traceless observables, that is, $G_{[1,k],t} \equiv G_{[1,k],t}(\trSet'; \bm x') \in \mathfrak{G}_{k,n}$ for some $\trSet'$ with $|\trSet'|=n$, and define
	\begin{equation} \label{eq:trless_quants}
		\mathcal{X}_t^{k,n} := \Tr\bigl[(G-M)_{[1,k],t}\other{A}_k\bigr] , \quad \mathcal{Y}_t^{k,n} := \bigl((G-M)_{[1,k],t} \bigr)_{ab},
	\end{equation}
	where $\other{A}_k$ is defined in \eqref{eq:trlessA_k}.
	Similarly to the proof of Lemma \ref{lemma:av_iso_evol} in Section \ref{sec:evols}, we conclude that the quantities
	 $\mathcal{Y}_t^{k,n}, \mathcal{X}^{k,n}_t$
	satisfy the evolution equations
	\begin{equation} \label{eq:traceless_k_iso_evol}
		\begin{split}
			\mathrm{d}\mathcal{Y}^{k,n}_{t} &= \biggl(\frac{k}{2}I+\bigoplus_{j=1}^{k-1} \other{\mathcal{A}}_{j,t}\biggr)\bigl[ \mathcal{Y}^{k,n}_{t}\bigr]\mathrm{d}t
			+ \mathrm{d}\mathcal{M}^{n, \mathrm{iso} }_{[1,k],t}  + \mathcal{F}^{\, n, \mathrm{iso}}_{[1,k],t} \mathrm{d}t,\\
			\mathrm{d}\mathcal{X}^{k,n}_t &= \biggl(\frac{k}{2}I+\bigoplus_{j=1}^k \other{\mathcal{A}}_{j,t}\biggr)\bigl[ \mathcal{X}^{k,n}_t\bigr]\mathrm{d}t
			+ \mathrm{d}\mathcal{M}^{n,\mathrm{av}}_{[1,k],t}  + \mathcal{F}^{\,n,\mathrm{av}}_{[1,k],t} \mathrm{d}t.
		\end{split}
	\end{equation}
	Here the martingale
	$\mathrm{d}\mathcal{M}^{n, \mathrm{iso} }_{[1,k],t}$ and the forcing $\mathcal{F}^{\, n, \mathrm{iso}}_{[1,k],t}$ terms, as well as their averaged versions,
	 are defined completely analogously to their counterparts with no traceless observables in \eqref{eq:av_k_mart}, \eqref{eq:iso_k_mart}, \eqref{eq:F_def}, and \eqref{eq:iso_F_def}, with $G_{[i,j],t}$ being the sub-chains of $G_{[1,k],t} \in \mathfrak{G}_{k,n}$, while the linear operators  $\other{\mathcal{A}}_{j,t}$ are given by
	\begin{equation} \label{eq:traceless_lin_prop_ops}
		\other{\mathcal{A}}_{j,t} :=  \frac{m_{j,t}m_{j+1,t}S\bigl(I - \mathds{1}_{j \in \trSet'\cup \trSet_k}\Pi \bigr)}{1- m_{j,t}m_{j+1,t}S}.
	\end{equation}
	Recall $\trSet_k \subset \{k\}$ from \eqref{eq:trlessA_k}, and that $\Pi := N^{-1} \bm{1}\bm{1}^*$ denotes the orthogonal projector onto the vector of ones $\bm{1} := (1,\dots, 1) \in \mathbb{C}^{N}$. 
	Note the difference between $\other{\mathcal{A}}_{j,t}$ and $\mathcal{A}_{j,t}$  defined in \eqref{eq:lin_prop_ops};
	the subtracted projection for indices $j$ corresponding to traceless observables. 
	 We used the identity 
	\begin{equation} \label{eq:trless_lin_term}
		 \mathcal{S}[M_{[j,j+1],t}] = \sum_q \trless{S}^q \bigl(M_{[j,j+1],t}\bigr)_{qq}  = \sum_q \trless{S}^q  \bigl(\other{\mathcal{A}}_{j,t} \bigr)_{qx_j}~, \quad j \in \trSet,
	\end{equation}
	which follows from $S\Pi = \Pi S$, to recover the action of $\other{\mathcal{A}}_{j,t}$ on the quantities $\mathcal{Y}^{k,n}_{t}$ and $\mathcal{X}^{k,n}_{t}$ in the linear term of \eqref{eq:traceless_k_iso_evol} from $	G_{[1,j],t}\mathcal{S}[M_{[j,j+1],t}] G_{[j+1,k],t}$.

	Recall from \eqref{eq:iso_k_mart} and \eqref{eq:iso_F_def} that the $\mathrm{iso}$-quantities depend implicitly on the matrix entry indices $a, b$; however, we do not follow this dependence and all estimates are understood to be uniform in $a,b\in\indset{N}$. 
	
	For fixed integers $n, l$, and tolerance exponents $\xi, \nu \in (0, \etaexp/100)$, we define the stopping time $\tau\equiv\tau_{n,l,\xi, \nu}$ as\footnote{With a slight abuse of notation, 
	we use the same $\tau$ notation for the  new stopping time as in \eqref{eq:tau_def}.}
	\begin{equation}\label{eq:trless_tau_def}
		\tau :=  \inf\biggl\{ t\in[\crit,T] \, :\, \max_{k \in \indset{n+1, \maxK'-2l} }
		N^{-\xi - k\nu} \min\bigl\{1, \trl_t(N\eta_t)^{l/2}\}\Psi_{k,t}^{n,\mathrm{av/iso}} \ge 1 \biggr\}.
	\end{equation}
	Compared with~\eqref{eq:tau_def},  the control parameters are much simpler, they are just $N^{\xi +k\nu}$
	for both the isotropic and averaged quantities, hence we do not need to introduce special notation and
	  write up  sophisticated master inequalities
	for them. Nevertheless, $N^{\xi +k\nu}$ deteriorates by a factor $N^\nu$ for each additional resolvent,
	so we still need to extract a self-improving mechanism
	from the evolution equations for $\Psi_{k,t}^{n,\mathrm{av/iso}}$ since these equations also contain  chains
	longer than $k$.  The $N^\nu$ powers will play a crucial role 
	in the entire proof since the self-improvement will be manifested by showing that up to time $\tau$ the
	estimates in its definition are actually better by some power of $N^{-\nu}$, allowing to conclude that $\tau=T$ with very high probability.

	To start, it follows from \eqref{eq:trless_ind_assume} that $\tau > \crit$ with very high probability.  
	Applying Duhamel's principle to \eqref{eq:traceless_k_iso_evol}, we obtain\footnote{Here we discarded the harmless order-one $\mathrm{e}^{k(t\wedge\tau -s)/2}$ factors, see discussion below \eqref{eq:k_av_propagator}.}
	\begin{equation} \label{eq:trless_Duhamel}
		\begin{split}
			\bigl\lvert \mathcal{X}^{k,n}_{t\wedge\tau} \bigr\rvert &\lesssim\Big| \mathcal{P}^{k,n}_{\crit, t\wedge\tau}\bigl[ \mathcal{X}^{k,n}_{\crit}\bigr]  
			+\int_{\crit}^t \mathcal{P}^{k,n}_{s, t\wedge\tau}\biggl[  \mathrm{d}\mathcal{M}^{n, \mathrm{av} }_{[1,k],s} \biggr]
			+\int_{\crit}^t \mathcal{P}^{k,n}_{s, t\wedge\tau}\biggl[ \mathcal{F}^{\, n, \mathrm{av}}_{[1,k],s} \biggr] \mathrm{d}s\Big|, \\
			\bigl\lvert \mathcal{Y}^{k,n}_{t\wedge\tau} \bigr\rvert &\lesssim \Big| \mathcal{P}^{k-1,n}_{\crit, t\wedge\tau}\bigl[ \mathcal{Y}^{k,n}_{\crit}\bigr]  
			+\int_{\crit}^t \mathcal{P}^{k-1,n}_{s, t\wedge\tau}\biggl[  \mathrm{d}\mathcal{M}^{n, \mathrm{iso} }_{[1,k],s} \biggr]
			+\int_{\crit}^t \mathcal{P}^{k-1,n}_{s, t\wedge\tau}\biggl[ \mathcal{F}^{\, n, \mathrm{iso}}_{[1,k],s} \biggr] \mathrm{d}s\Big|,\\
		\end{split}	 
	\end{equation}
	where the propagators $\mathcal{P}^{k-1,n}_{s, t}$ and $\mathcal{P}^{k,n}_{s, t}$ are given by\footnote{Strictly speaking
	these propagators depend on the index set $\trSet$ and not only on its cardinality $n$, but we omit this dependence
	since in the estimates only $n$ will matter.}
	\begin{equation} \label{eq:trless_props}
		\mathcal{P}^{k',n}_{s, t} :=  \bigotimes_{j=1}^{k'} \widetilde{\mathcal{P}}^{[j]}_{s, t}, \qquad 
		\widetilde{\mathcal{P}}^{[j]}_{s, t}:=
		\exp\biggl\{ \int_{s}^t 
		 \bigl(I + \other{\mathcal{A}}_{j,r}\bigr)\mathrm{d}r \biggr\}, \quad \crit \le s\le t\le T, \quad k' \in \{k-1,k\}.
	\end{equation}
	For a function  $f(\bm x)$ of external indices $\bm x \in \indset{N}^{k}$, define the maximum norm
	 $\norm{f}_*$  
	 as 
	\begin{equation} \label{eq:star_norm}
		\norm{f}_* 
		:= \max_{\bm x \in \indset{N}^{k}}\bigl\lvert f(\bm x) \bigr\rvert. 
	\end{equation}
	Notice that since we are in the regime $t\ge t_*$, the spatial dependence is irrelevant and we can operate
	with the simple maximum norm when controlling the action of the propagator instead of the 
	more complicated weighted norms that underlie the propagator estimates  in  Lemma~\ref{lemma:good_props}. 
	
	Since $\Pi$ commutes with $S$---and hence with $\other{\mathcal{A}}_{j,t}$---it follows from \eqref{eq:P_decomp}, that for all $j\in \indset{k}$, 
	\begin{equation} \label{eq:tr_prop_decomp}
		\widetilde{\mathcal{P}}^{[j]}_{s, t} = \mathrm{e}^{t-s} I +  \bigl(\mathrm{e}^{t-s}-1\bigr) \other{\mathcal{A}}_{j,t}, \quad \crit \le s \le t \le T.
	\end{equation}
	Hence, if $\other{\mathcal{A}}_{j,t}$ is not saturated, that is, $(\im z_j) (\im z_{j+1}) > 0$, then we use the second bound in \eqref{eq:Ups_majorates} and \eqref{eq:Ups_norm_bounds} to deduce that $\lVert\widetilde{\mathcal{P}}^{[j]}_{s, t}\rVert_{\infty\to\infty} \lesssim 1$. On the other hand, if $\other{\mathcal{A}}_{j,t}$ is saturated, that is $(\im z_j) (\im z_{j+1}) < 0$, then, $\lVert\widetilde{\mathcal{P}}^{[j]}_{s, t}\rVert_{\infty\to\infty} \lesssim 1$ for $j \in \trSet$ by the second bound in \eqref{eq:supercrit_Theta} together with $(N/W)^2= 1/\eta_{t_*} \le 1/\eta_s$, and $\lVert\widetilde{\mathcal{P}}^{[j]}_{s, t}\bigr\rVert_{\infty\to\infty} \lesssim \eta_s/\eta_t$ by the first bound in \eqref{eq:supercrit_Theta}.
	  Therefore, the norm of the propagator $\mathcal{P}^{k',n}_{s, t}$ with $k'\in \{k-1,k\}$, induced by the $\norm{\cdot}_*$-norm on functions of $\bm x \in \indset{N}^{k}$, satisfies
	\begin{equation} \label{eq:tr_prop}
		\norm{\mathcal{P}^{k',n}_{s, t}}_{*\to *}  \lesssim  \biggl(\frac{\eta_s}{\eta_t}\biggr)^{k'-n}~, \quad\crit \le s\le t\le T, \quad  k' \in \{k-1,k\}.
	\end{equation}
	The estimate~\eqref{eq:tr_prop}
	expresses the key mechanism why traceless observables result in better bounds: the propagators
	$ \widetilde{\mathcal{P}}^{[j]}_{s,t}$ corresponding to indices $j$ of traceless observables are bounded, while
	they blow up as $\eta_s/\eta_t$ in general. The net effect of this improved propagator estimate on the final 
	result is analogous to the effect of the non-saturated index in~\eqref{eq:av_prop_bound} as compared to~\eqref{eq:bad_av_prop},
	as well as the effect of the regularization of the propagator in Lemma~\ref{lemma:reg_props} as compared to~\eqref{eq:noring}.
	
	Since we can work in the simple maximum norm~\eqref{eq:star_norm},  it suffices to estimate the martingale $\mathrm{d}\mathcal{M}^{n, \mathrm{av/iso} }_{[1,k],s}$ and the forcing $\mathcal{F}^{\, n, \mathrm{av/iso}}_{[1,k],s}$ terms in this norm,
	 as functions of $\bm x$ (the external index $x_k$ is redundant in the $\mathrm{iso}$ case). 
	 We start with the forcing term.
	\begin{claim} \label{claim:traceless_forcing}
		Under the assumptions of Proposition \ref{prop:trless_bootstrap}, with $\tau$ being the stopping time defined in \eqref{eq:trless_tau_def},  the quantities $\mathcal{F}^{\, n, \mathrm{av/iso}}_{[1,k],s}$  satisfy
		\begin{equation} \label{eq:tr_av_focing}
			\norm{\mathcal{F}^{\, n, \mathrm{av}}_{[1,k],s}}_* \prec N^{-\nu}\frac{1}{\eta_s} \frac{N^{\xi+k\nu}  \trl_s^{n}}{(N\eta_s)^{k}}\biggl( 1 + \frac{1}{\trl_s(N\eta_s)^{l/2}}\biggr) 
			, \quad s \in [\crit, \tau], \quad k \in \indset{n, K'-2l}
		\end{equation}
		\begin{equation} \label{eq:tr_iso_focing}
			\norm{\mathcal{F}^{\, n, \mathrm{iso}}_{[1,k],s}}_* \lesssim N^{-\nu}\frac{1}{\eta_s}\frac{N^{\xi+k\nu}  \trl_s^{n}}{(N\eta_s)^{k-1/2}}\biggl( 1 + \frac{1}{\trl_s(N\eta_s)^{l/2}}\biggr) 
			, \quad s \in [\crit, \tau], \quad k \in \indset{n+1, K'-2l}.
		\end{equation}
		In fact, the bound \eqref{eq:tr_av_focing} for $k=n$ holds for all $s \in [\crit, T]$.
	\end{claim}
	Note that this lemma controls forcing terms coming from equations for  $k$-chains with only $k\le K'-2l$,
	while the assumption \eqref{eq:power_induction_assume} involves chains up to length $K'-2l+2$. This 
	loss of chain length is necessary since a forcing term for a $k$-chain may contain $k+2$ resolvents.
	Note that the estimates  \eqref{eq:tr_av_focing}--\eqref{eq:tr_iso_focing} are natural and sufficient for
	the $l$-induction step stated in Proposition \ref{prop:trless_bootstrap}. Recall that the $1/\eta_s$ factor 
	is due to differentiation, but it will be compensated by the 
	time integration in Duhamel's  formula. The $\theta_s^n$ factor is the improvement coming from $n$ traceless observables, while
	the $N\eta_s$ powers are the natural sizes of the general local laws. The last factor with an $l$-th power
	 carries the $\sqrt{N\eta}$ improvement over the assumption \eqref{eq:power_induction_assume}.
	 Finally, notice that both bounds contain an additional  $N^{-\nu}$ compared to the target power in the definition of the stopping time; this will eventually guarantee the self-improving mechanism.

	\begin{proof}[Proof of Claim \ref{claim:traceless_forcing}]
		The proof of \eqref{eq:tr_av_focing}--\eqref{eq:tr_iso_focing} is essentially analogous to that of 
		Lemma \ref{lemma:forcing} in Section \ref{sec:forcing}, except the $M$-bounds \eqref{eq:M_bound}--\eqref{eq:M_bound_av}
		 are replaced by Lemma \ref{lemma:traceless_M_bounds}, and \eqref{eq:power_induction_assume} together with \eqref{eq:trless_tau_def} are used to bound $(G-M)$ chains. Unlike in the proof of Lemma \ref{lemma:forcing}, to bound the terms
		\begin{equation}
			\sum_{j=1}^k \Tr\bigl[\mathcal{S}[G_{j,s}-m_{j,s}] G_{[j,k],s}S^{x_k}G_{[1,j],s}\bigr]  \quad \text{and} \quad 
			\sum_{i=1}^k \bigl(G_{[1,i],s} \mathcal{S}\bigl[G_{i,s}-m_{i,s}\bigr] G_{[i, k],t}\bigr)_{ab},
		\end{equation}
		(see~\eqref{eq:1_k+1_reduction1}--\eqref{eq:1_k+1_bound3} 
		and~\eqref{eq:k+1_1_iso_term}--\eqref{eq:k+1_1_iso_bound}),
		we use \eqref{eq:power_induction_assume} for a $k+1$ chains with $n$ regular observables $G_{[j,k],s}S^{x_k}G_{[1,j],s}$ and $G_{[1,i],s} S^q G_{[i, k],s}$, respectively, instead of performing any reductions. 
		To have a suboptimal bound on longer chains available is a substantial simplification compared
		to  proof of Lemma \ref{lemma:forcing}.
		This concludes the proof of Claim \ref{claim:traceless_forcing}.
	\end{proof}
	
	Now we turn to estimating the martingale terms in~\eqref{eq:traceless_k_iso_evol}, 
	 where  simply mimicking the proof
	of the analogous Lemma~\ref{lemma:mart_est} is not sufficient. 
	Recall from the discussion in Section \ref{sec:mart}, that the key step towards estimating the contribution of the martingale terms in \eqref{eq:traceless_k_iso_evol} is bounding the $Q_{k,s}$ quantities, defined in \eqref{eq:Qkj} for the averaged and in \eqref{eq:iso_Qkj} for the isotropic case, that arise from the quadratic variation of the corresponding martingales. In Section \ref{sec:mart}, we estimated $Q_{k,s}$ with large $k$ using reduction inequalities (see \eqref{eq:reduction}). However, performing a reduction at a traceless observable would incur the loss of the tracelessness effect, which is not affordable.  	  
	To estimate  $Q_{k,s}$ without losing any traceless observables, we will derive and analyze  a separate
	evolution equation for them. To formalize this new step precisely, we
	 introduce a new class of resolvent chains;
	these  include the types of chains that arise in a quadratic variation of the chains belonging to $\mathfrak{G}_{k,n}$.
	
	For all $k \ge 1$ and $0 \le n \le k-1$, we define the class of resolvent chains $\mathfrak{Q}_{k,n}$ by
	\begin{equation} \label{eq:Q_class}
		\mathfrak{Q}_{k,n} := \bigl\{ \bigl(G_{[1,k],t}(\bm x)\bigr)^* S^x G_{[1,k],t}(\bm y) \, :\, G_{[1,k],t} \in \mathfrak{G}_{k,n}, ~ \bm x, \bm y \in \indset{N}^{k-1}, ~x \in \indset{N} \bigr\},
	\end{equation}
	where we again suppress the time-dependence of $\mathfrak{Q}_{k,n}$ in the notation.
	
	 The key observation is that, while most observables appearing in $Q_{k,s} \in \mathfrak{Q}_{k,n}$ are inherited from the original, potentially traceless, observables and must therefore be preserved,
	one observable — namely, $S^x$ — arises  from the variance $S_{ab}$ within
	  the formula for the quadratic variation in \eqref{eq:Qkj}, \eqref{eq:iso_Qkj} and it is inherently not traceless.  
	However, since $S^x$
	is placed exactly in the middle of the chain, applying reduction to it directly would not result in a self-improving bound. 
	Instead, we use a separate evolution equation for the $Q$-chains in $\mathfrak{Q}_{k,n}$.
	In turn, their quadratic variation will already contain at least three non-traceless observables distributed throughout the chain, providing sufficient freedom for an effective reduction scheme that cuts the chain specifically at these observables, which ultimately enables us to harness the optimal improvement coming from every traceless observable.  
	
	We remark that, although objects arising from the quadratic variation are explicitly symmetric (see \eqref{eq:Qkj} and \eqref{eq:iso_Qkj}), the chains in $\mathfrak{Q}_{k,n}$ are only symmetric up to the external indices $\bm x$ and $\bm y$. 
	We included these more general chains in the class $\mathfrak{Q}_{k,n}$   because 
	differentiating a symmetric $Q$-chain also gives rise to asymmetric ones (see \eqref{eq:Q_evol} below).     
	
	Let $[ \mathrm{d}\mathcal{M}^{n, \mathrm{av/iso} }_{[1,k]}]_s$ denote
	the infinitesimal quadratic variation of $\mathrm{d}\mathcal{M}^{n, \mathrm{av/iso} }_{[1,k],s}$, then, from the expressions \eqref{eq:Qkj} and \eqref{eq:iso_Qkj}, it follows that 
	\begin{equation} \label{eq:QV_Q_bounds}
		\begin{split}
			\bigl[ \mathrm{d}\mathcal{M}^{n, \mathrm{av} }_{[1,k]} \bigr]_s &\lesssim N \sum_{j=1}^k \max_{x\in \indset{N}} \bigl(G_{[1,j],s}^*S^{x_k} G_{[j,k],s}^* S^x G_{[j,k],s}S^{x_k}G_{[1,j],s}\bigr)_{xx},\\
			\bigl[ \mathrm{d}\mathcal{M}^{n, \mathrm{iso} }_{[1,k]} \bigr]_s &\lesssim  N\sum_{j=1}^k \max_{x\in \indset{N}} \bigl(G_{[1,j],s} S^x G_{[1,j],s}^* \bigr)_{aa}   \max_{y\in \indset{N}} \bigl(G_{[j,k],s}^* S^y G_{[j,k],s}\bigr)_{bb},
		\end{split}
	\end{equation} 
	where we used the identity \eqref{eq:sumS=1} to obtain the estimate
	\begin{equation} \label{eq:inset_S}
		\sum_c \bigl\lvert \bigl(G_{[1,j],s}\bigr)_{ac}\bigr\rvert^2 = \sum_{x}\bigl(G_{[1,j],s} S^x G_{[1,j],s}^* \bigr)_{aa} \le N\max_{x\in \indset{N}} \bigl(G_{[1,j],s} S^x G_{[1,j],s}^* \bigr)_{aa}.
	\end{equation} 
	Recall that $a,b$ denote the matrix-entry indices from the definition of the isotropic quantity $\mathcal{Y}^{k,n}_{t}$ in \eqref{eq:trless_quants}, that are suppressed in the notation.
	Note that all the chains on the right-hand side of \eqref{eq:QV_Q_bounds} belong to the class $\mathfrak{Q}_{k',n'}$ for some $n' \in \indset{0,n}$ and $k' \in \indset{n'+1, k+1}$, in particular, their length
	may increase by one to $k+1$. We now proceed to estimate these chains.
	
	By using the  definition of $\mathfrak{Q}_{k',n'}$ in \eqref{eq:Q_class}, and writing a matrix element of
	$Q_s\in \mathfrak{Q}_{k',n'}$ as
	$$
	 (Q_s)_{cd}:=\big[\bigl(G_{[1,k'],s}(\bm x)\bigr)^* S^x G_{[1,k'],s}(\bm y)\big]_{cd} =\sum_q S_{xq} \overline{\bigl(G_{[1,k'],s}(\bm x)\bigr)_{qc}}
	 \bigl(G_{[1,k'],s}(\bm y)\bigr)_{qd},
	 $$ 
	from Lemma~\ref{lemma:traceless_M_bounds}, and \eqref{eq:traceless_Psi_def} it follows that, for all $n' \in \mathbb{N}$ and $k' \ge n'+1$,
	\begin{equation}
		\max_{Q_{s} \in \mathfrak{Q}_{k',n'}}\norm{Q_{s}}_{\max} \lesssim \frac{\trl_s^{2n'}}{(N\eta_s)^{2k'-1}}\bigl(1+\Psi_{k',s}^{n',\mathrm{iso}}\bigr)^2,
	\end{equation}
	Hence, it follows immediately
	 from \eqref{eq:trless_ind_assume}, that  
	\begin{equation} \label{eq:lesser_Q_bounds}
		\max_{\crit\le s\le T}\max_{Q_{s} \in \mathfrak{Q}_{k',n'}} \norm{Q_{s}}_{\max} \prec \frac{\trl_s^{2n'}}{(N\eta_s)^{2k'-1}}, \quad n' \in \indset{0,n-1}, \quad k' \in \indset{n'+1, K'}.
	\end{equation}
	
	On the other hand, for $n'=n$ the assumption \eqref{eq:trless_ind_assume} cannot be used since it has less 
	than $n$ traceless
	observables, so 
	we separately estimate the chains in $\mathfrak{Q}_{k',n}$ using the following  lemma.
	\begin{lemma} [Bound on the $Q$-Chains] \label{lemma:Q_bounds} Fix $n, l\ge 1$ and $K'\ge 2L+n+2$ integers. 
		Under the assumptions of Proposition \ref{prop:trless_bootstrap}, 
		  in particular \eqref{eq:trless_ind_assume}  and \eqref{eq:power_induction_assume},   we have,
		for all $k' \in \indset{n+1, \maxK'-2l+1}$, that
		\begin{equation} \label{eq:n_Q_bounds}
			\max_{\crit\le s\le T}\max_{Q_{s} \in \mathfrak{Q}_{k',n}} \norm{Q_{s}}_{\max} \prec \frac{\trl_s^{2n}}{(N\eta_s)^{2k'-1}}\biggl(1 + \frac{1}{\trl_s^2(N\eta_s)^{l}}\biggr).
		\end{equation}
	\end{lemma}
	Note  that length of the  chain (that is, half of $Q_s$) in the conclusion~\eqref{eq:n_Q_bounds}
	 satisfies $k'\le K'-2l+1$, while the $l$-induction hypothesis
	\eqref{eq:power_induction_assume} involves slightly longer chains of length  up to $K'-2l+2$.
	The other input, the optimal traceless law~\eqref{eq:trless_ind_assume} and its main implication~\eqref{eq:lesser_Q_bounds} coming from the previous step in the $n$-induction, involve chains of length 
	up to  $K'$ but only with at most $n-1$ traceless observables. This extra room will be needed 
	since differentiating a $k$-chain naturally produces a chain of length $k+1$. 
	Note that in the proof of the Master Inequalities in Proposition~\ref{prop:masters}, the 
	maximal length of the chains was fixed, so we used a reduction step even for the forcing terms
	(see~\eqref{eq:1_k+1_reduction1}) which lost a factor $\sqrt{\ell\eta}$. 
	Since in Theorem \ref{th:local_laws} we already established the local laws for chains of arbitrarily large fixed length, we can use them as an a priori bound  and opt for the strategy of reducing the chain length along both inductions.

	We prove Lemma \ref{lemma:Q_bounds}  Section~\ref{sec:Q_bound}. 
	Combining the estimates \eqref{eq:QV_Q_bounds}, \eqref{eq:lesser_Q_bounds}, and \eqref{eq:n_Q_bounds}, we conclude that 
	\begin{equation} \label{eq:tr_mart_estimates}
		\begin{split}
			\norm{\bigl[ \mathrm{d}\mathcal{M}^{n, \mathrm{av} }_{[1,k]} \bigr]_s}_* &\prec \frac{1}{\eta_s} \frac{\trl_s^{2n}}{(N\eta_s)^{2k}}\biggl(1 + \frac{1}{\trl_s^2(N\eta_s)^{l}}\biggr), \quad k \in  \indset{n, \maxK'-2l}\\ 
			\norm{\bigl[ \mathrm{d}\mathcal{M}^{n, \mathrm{iso} }_{[1,k]} \bigr]_s}_* &\prec  \frac{1}{\eta_s} \frac{\trl_s^{2n}}{(N\eta_s)^{2k-1}}\biggl(1 + \frac{1}{\trl_s^2(N\eta_s)^{l}}\biggr), \quad k \in \indset{n+1, \maxK'-2l}.
		\end{split}
	\end{equation}
	Again, these estimates are natural and sufficient for the $l$-induction step. Time integration will 
	compensate for the  $1/\eta_s$ factor, the $(N\eta_s)$ factors express the natural sizes 
	of the corresponding $G$-chains in general, while the   $\theta_s$ factors carry the optimal improvement
	due to the traceless observables. Finally, the last additional factor in~\eqref{eq:tr_mart_estimates},
	 after taking square root, is exactly
	the (reciprocal of the) improved target factor in~\eqref{eq:power_induction_conclude}. 
	
	Indeed, plugging the first two bounds in \eqref{eq:trless_initial}, \eqref{eq:tr_av_focing}--\eqref{eq:tr_iso_focing}, and \eqref{eq:tr_mart_estimates} (together with \eqref{eq:mart_ineq}) into
	the Duhamel formula~\eqref{eq:trless_Duhamel}, we conclude, using the propagator bound~\eqref{eq:tr_prop} and 
	the time integration rules~\eqref{eq:int_rules}, that
	\begin{equation}
		\begin{split}
			\bigl\lvert \mathcal{X}^{k,n}_{t\wedge\tau} \bigr\rvert \prec   N^{-\nu}\frac{N^{\xi + k\nu}\trl_s^{n}}{(N\eta_s)^{k-1/2}}\biggl(1 + \frac{1}{\trl_s(N\eta_s)^{l/2}}\biggr), \quad k \in \indset{n, \maxK' - 2l},\\
			\bigl\lvert \mathcal{Y}^{k,n}_{t\wedge\tau} \bigr\rvert \prec  N^{-\nu}\frac{N^{\xi + k\nu}\trl_s^{n}}{(N\eta_s)^{k}}\biggl(1 + \frac{1}{\trl_s(N\eta_s)^{l/2}}\biggr), \quad k \in \indset{n+1, \maxK' - 2l}.
		\end{split}
	\end{equation}
	We direct the reader to the proof of Proposition \ref{prop:masters} in Section \ref{sec:masters_proof} for a similar argument described in full detail; in particular the forcing term bounds \eqref{eq:tr_av_focing}--\eqref{eq:tr_iso_focing}
	correspond to
	\eqref{eq:av_forcing_bound}--\eqref{eq:iso_forcing_bound}  in Lemma~\ref{lemma:forcing}
	and the martingale estimates \eqref{eq:tr_mart_estimates} correspond to
the estimates~\eqref{eq:Q_est} and~\eqref{eq:isoQ_est} within the proof of Lemma~\ref{lemma:mart_est} 
before the propagator acted on the martingales\footnote{We did not formulate the direct analogue of Lemma~\ref{lemma:mart_est} since estimating the action of the propagator here is trivial using  $*\to *$ norm.}.
	Therefore, the stopping time $\tau = T$ with very high probability, and \eqref{eq:power_induction_conclude} is established.
	This concludes the proof of Proposition~\ref{prop:trless_bootstrap}.
\end{proof}

\subsection{Equations for $Q$-chains: Proof of Lemma \ref{lemma:Q_bounds}}\label{sec:Q_bound}
\begin{proof}[Proof of Lemma \ref{lemma:Q_bounds}] 
We prove Lemma \ref{lemma:Q_bounds} by writing up a separate evolution
equation for the special $Q$-chains and analyze them by 
yet another Duhamel formula combined with a stopping time argument. We consider the integers $n, l\ge 1$ and $K'$ fixed. 
For all $k \in \indset{n+1, \maxK'-2l+1}$, define the auxiliary control quantities $\Phi_{k,t}^{n, \mathrm{av/iso}}$ as
\begin{equation} \label{eq:traceless_Phi_def}
	\Phi_{k,t}^{n,\mathrm{iso}} := \max_{Q_t \in \mathfrak{Q}_{k,n}}  \frac{(N\eta_t)^{2k-1} \norm{Q_t}_{\max}}{\trl_t^{2n }\bigl(1+\trl_t^{-2}(N\eta_t)^{-l}\bigr)}, \quad 
	\Phi_{k,t}^{n,\mathrm{av}} := \max_{Q_t \in \mathfrak{Q}_{k,n}}  \max_{y\in\indset{N}}  \frac{(N\eta_t)^{2k-1} \bigl\lvert \Tr\bigl[Q_tS^y\bigr]\bigr\rvert }{\trl_t^{2n}\bigl(1+\trl_t^{-2}(N\eta_t)^{-l}\bigr)}.
\end{equation}
Note that unlike all previous control quantities, here $\Phi$'s control the
sizes of the chains $Q_t$ and {\it not}  their fluctuations, that is, the deviation of $Q_t$  from their  deterministic approximation, given by the corresponding $M$-term (viewing $Q_t$ as a $2k$-chain). For the same reason, the 
normalization factors in definitions of $\Phi_{k,t}^{n,\mathrm{iso}}$ and $\Phi_{k,t}^{n,\mathrm{av}}$ are 
the identical: they both reflect the natural size of the corresponding $Q$-chain. 
Even though the target estimate in Lemma \ref{lemma:Q_bounds} only involves isotropic quantities, the better control (by a multiplicative factor of $N^{-\nu}$) on their averaged analogues in \eqref{eq:trless_tau2_def} below is necessary to complete the proof. 

We fix tolerance exponents $\xi, \nu \in (0, \etaexp/100)$ define the stopping time $\tau \equiv \tau_{n,l,\xi,\nu}$ 
as\footnote{This new stopping time is different from the ones in~\eqref{eq:tau_def}
and~\eqref{eq:trless_tau_def}.}
\begin{equation}\label{eq:trless_tau2_def}
	\begin{split}
		\tau := &~  \inf \biggl\{ t\in[\crit,T] \, :\, \max_{k \in \indset{n+1, \maxK'-2l+1}}
		  N^{-\xi - k\nu} \Phi_{k,t}^{n,\mathrm{av}} \ge 1 \\
		&\wedge \inf \biggl\{ t\in[\crit,T] \, :\,\max_{k \in \indset{n+1, \maxK'-2l+1}} N^{-\xi - (k+1)\nu}\Phi_{k,t}^{n,\mathrm{iso}} \ge 1  \biggr\}.
	\end{split}
\end{equation}
It follows from \eqref{eq:power_induction_assume} that $\tau > \crit$ with very high probability.  

Let $\bm x := (\bm x_1^*, x, \bm x_2, y)$ denote a vector of indices in $\indset{N}^{2k}$, where
$\bm x_1, \bm x_2\in \indset{N}^{k-1}$ and star denotes the reverse order.  
As usual, we set $\bm x' = (\bm x_1^*, x, \bm x_2)$. We often omit them from the notation. In particular, to condense the presentation, we denote 
$$
G_{[i,j],t} \equiv G_{[i,j],t}(\bm x_2), \quad \mbox{and} \quad (G_{[i,j],t})^\dagger \equiv \bigl(G_{[i,j],t}(\bm x_1)\bigr)^*,
$$
 where $(\cdot)^*$ denotes the usual Hermitian conjugation.

Let $Q_{k,t} \equiv Q_{k,t}(\bm x') := \bigl(G_{[1,k],t}(\bm x_1)\bigr)^*S^xG_{[1,k],t}(\bm x_2) \in \mathfrak{Q}_{k,n}$. Then $Q_{k,t} = G_{[1,k],t}^\dagger S^x G_{[1,k],t}$.  Define the quantities 
\begin{equation} \label{eq:Q_quantities}
	\mathcal{L}_t^{k,n} :=\Tr\bigl[Q_{k,t}S^{y}\bigr] , \quad \mathcal{K}_t^{k,n} := \bigl(Q_{k,t} \bigr)_{ab}~,
\end{equation}
where, as usual, in the notation we suppress the irrelevant dependence on $y, a, b$.
Then, $\mathcal{L}_t^{k,n}$ and $\mathcal{K}_t^{k,n}$ satisfy the evolution equations
\begin{equation}\label{eq:Q_evol}
	\begin{split}
		\mathrm{d}\mathcal{L}_t^{k,n} = \biggl(k I + \bigoplus_{j=1}^{2k} \other{\mathcal{A}}_{j,t}\biggr)\bigl[\mathcal{L}_t^{k,n}\bigr]\mathrm{d}t + \mathrm{d}\mathcal{E}_{k,t}^{n, \mathrm{av} } + \mathcal{W}_{k,t}^{n, \mathrm{av} }\mathrm{d}t,\\
		\mathrm{d}\mathcal{K}_t^{k,n} = \biggl(k I + \bigoplus_{j=1}^{2k-1} \other{\mathcal{A}}_{j,t}\biggr)\bigl[\mathcal{K}_t^{k,n}\bigr]\mathrm{d}t + \mathrm{d}\mathcal{E}_{k,t}^{n, \mathrm{iso} } + \mathcal{W}_{k,t}^{n, \mathrm{iso} }\mathrm{d}t,
	\end{split}
\end{equation}
where the linear operators $\other{\mathcal{A}}_{j,t}$ for $j \in \indset{k-1}$ are as defined in \eqref{eq:traceless_lin_prop_ops}, $\other{\mathcal{A}}_{k+j,t} = \other{\mathcal{A}}^*_{k-j,t}$, and $\other{\mathcal{A}}_{k,t} = \other{\mathcal{A}}_{2k,t} = \Theta_t$. Here, the martingale terms $\mathrm{d}\mathcal{E}_{k,t}^{n, \mathrm{av/iso} }$ are given by
\begin{equation} \label{eq:Q_mart_def}
	\mathrm{d}\mathcal{E}_{k,t}^{n, \mathrm{av} } := \sum_{ab}\sqrt{S_{ab}}\Tr\bigl[ \partial_{ab} Q_{k,t} S^{y}\bigr] \mathrm{d}\mathfrak{B}_{ab,t},  \quad \mathrm{d}\mathcal{E}_{k,t}^{n, \mathrm{iso} } := \sum_{cd}\sqrt{S_{cd}}\,\partial_{cd}\bigl( Q_{k,t}\bigr)_{ab} \mathrm{d}\mathfrak{B}_{cd,t},
\end{equation}
where we once again drop the dependence of $\mathrm{d}\mathcal{E}_{k,t}^{n, \mathrm{iso} }$ on the matrix entry indices $a,b$ as well as the dependence of $\mathrm{d}\mathcal{E}_{k,t}^{n, \mathrm{av} }$ on $y$. 
The forcing terms in \eqref{eq:Q_evol} are given by
\begin{equation} \label{eq:Q_forcing_def}
	\begin{split}
		\mathcal{W}_{k,t}^{n, \mathrm{iso} } &:= \bigl(G_{[1,k],t}^\dagger S^x \mathcal{W}_{k,t}^\mathrm{dir}  + \bigl(\mathcal{W}_{k,t}^\mathrm{dir}\bigr)^\dagger S^x G_{[1,k],t} +  \mathcal{W}_{k,t}^\mathrm{crs}\bigr)_{ab} 
		+ \delta_{ab}|m_t|^2 \sum_q \bigl(I+\Theta_t\bigr)_{aq}\Tr\bigl[Q_{k,t}S^q\bigr],\\
		\mathcal{W}_{k,t}^{n, \mathrm{av} } &:= \Tr\bigl[G_{[1,k],t}^\dagger S^x \mathcal{W}_{k,t}^\mathrm{dir}S^y\bigr] + \Tr\bigl[\bigl(\mathcal{W}_{k,t}^\mathrm{dir}\bigr)^\dagger S^x G_{[1,k],t}S^y\bigr] + \Tr\bigl[\mathcal{W}_{k,t}^\mathrm{crs}S^y\bigr],
	\end{split}
\end{equation}
where we group the forcing terms into "direct" (denoted by $\mathcal{W}^\mathrm{dir}$) and "cross" ( $\mathcal{W}^\mathrm{crs}$) terms, defined as
\begin{equation} \label{eq:direct_forcing}
	\mathcal{W}_{k,t}^\mathrm{dir} := \sum_{\substack{1\le i \le j \le k \\ j-i \le 1}}  G_{[1,i],t} \mathcal{S}\bigl[(G-M)_{[i,j],t}\bigr] G_{[j,k],t}  + \sum_{\substack{1 \le i < j \le k\\ 2 \le j-i}}  G_{[1,i],t} \mathcal{S}\bigl[G_{[i,j],t} \bigr]G_{[j+1,k],t}\\
\end{equation}
\begin{equation} \label{eq:cross_forcing}
	\begin{split}
		\mathcal{W}_{k,t}^\mathrm{crs} :=&~ G_{[1,k],t}^\dagger \mathcal{S}\bigl[G_{1,t}^* S^xG_{1,t} - \Theta^x_t\bigr]G_{[1,k],t} + \sum_q \bigl(G_{k,t}^* S^q G_{k,t} - \Theta^q_t\bigr) \bigl(Q_{k,t}\bigr)_{qq}\\
		&+ \sum_{\substack{1\le i,j \le k \\ 3 \le i+j \le 2k-1}} G_{[i,k],t}^\dagger \mathcal{S}\bigl[G_{[1,i],t}^\dagger S^x G_{[1,j],t}\bigr]G_{[j,k],t}.
	\end{split}
\end{equation}
Note that these two types of terms are structurally quite different, in particular $\mathcal{W}_{k,t}^\mathrm{dir}$
contains $k+2$ resolvents, while $\mathcal{W}_{k,t}^\mathrm{crs}$ has $2k+2$. 
Compared with~\eqref{eq:F_def} and~\eqref{eq:iso_F_def}, most $M$ terms are absent in these definitions,
since we did not subtract the corresponding $M$ from $Q$ and thus the evolution equation~\eqref{eq:dm} is not used.
However, two specific $M_{[i,j],t}$ terms remain; the  $j=i$ term comes from differentiating the spectral parameter $z_i(t)$,
while the $j=i+1$ term comes from separating the linear terms $\other{\mathcal{A}}_{j,t}$ (recall that the $M$-term
of $G^*S^xG$ is $\Theta^x$ in \eqref{eq:cross_forcing}).
Note that for the direct terms in \eqref{eq:direct_forcing}, the argument of the super-operator $\mathcal{S}$  does not contain $S^x$, meaning the two derivatives giving rise to this term acted on the same half of the chain $Q_{k,t}$. For the cross terms in \eqref{eq:cross_forcing}, on the other hand, the two derivatives acted on the different halves of $Q_{k,t}$, hence the matrix $S^x$ is inside the argument of the super-operator $\mathcal{S}$.

Similarly to \eqref{eq:trless_Duhamel} above, applying Duhamel's principle to \eqref{eq:Q_evol} generates linear propagators which can be estimated in the $\norm{\cdot}_* \to \norm{\cdot}_*$ norm, see \eqref{eq:tr_prop}. Therefore, our goal is to bound the $\norm{\cdot}_*$-norm of the forcing terms $\mathcal{W}_{k,t}^{n, \mathrm{av/iso} }$ and the infinitesimal quadratic variations of $\mathrm{d}\mathcal{E}_{k,t}^{n, \mathrm{av/iso} }$.
\begin{claim} \label{claim:Q_mart} Under the assumptions of Lemma \ref{lemma:Q_bounds}, 
in particular \eqref{eq:trless_ind_assume} 
(and its main implication~\eqref{eq:lesser_Q_bounds}), as well as~\eqref{eq:power_induction_assume}, 
the infinitesimal quadratic variation of the martingale terms in \eqref{eq:Q_evol} satisfy
	\begin{equation} \label{eq:Q_mart}
		\begin{split}
			\norm{\bigl[\mathrm{d}\mathcal{E}_{k}^{n, \mathrm{av} }\bigr]_s }_*  &\prec N^{-\nu/2}\frac{1}{\eta_s}  \biggl[\frac{N^{\xi + k\nu}\trl_s^{2n}}{(N\eta_s)^{2k-1}} \biggl( 1   +\frac{1}{\trl_s^2(N\eta_s)^{l}} \biggr)\biggr]^2, \quad s \in [\crit, \tau],\\
			\norm{\bigl[\mathrm{d}\mathcal{E}_{k}^{n, \mathrm{iso} }\bigr]_s }_*  &\prec N^{-\nu/2}\frac{1}{\eta_s}  \biggl[\frac{ N^{\xi + (k+1)\nu} \trl_s^{2n}}{(N\eta_s)^{2k-1}} \biggl( 1   +\frac{1}{\trl_s^2(N\eta_s)^{l}} \biggr)\biggr]^2, \quad s \in [\crit, \tau],
		\end{split}
	\end{equation} 
	for any $k \in \indset{n+1, \maxK'-2l+1}$, where $\tau$ is  the stopping time defined in \eqref{eq:trless_tau2_def}.
\end{claim}
\begin{claim} \label{claim:Q_forcing} Under the assumptions of Lemma \ref{lemma:Q_bounds},
in particular \eqref{eq:trless_ind_assume} 
(and its main implication~\eqref{eq:lesser_Q_bounds}), as well as~\eqref{eq:power_induction_assume},
 the forcing terms in \eqref{eq:Q_evol} satisfy
	\begin{equation} \label{eq:Q_forcing}
		\begin{split}
			\norm{\mathcal{W}_{k,s}^{n, \mathrm{av} }}_*  &\prec N^{-\nu/2} \frac{1}{\eta_s} \frac{ N^{\xi + k\nu}   \trl_s^{2n}}{(N\eta_s)^{2k-1}} \biggl( 1   +\frac{1}{\trl_s^2(N\eta_s)^{l}} \biggr), \quad s \in [\crit, \tau],\\
			\norm{\mathcal{W}_{k,s}^{n, \mathrm{iso} }}_*  &\prec N^{-\nu/2}\frac{1}{\eta_s}  \frac{ N^{\xi + (k+1)\nu}  \trl_s^{2n}}{(N\eta_s)^{2k-1}} \biggl( 1   +\frac{1}{\trl_s^2(N\eta_s)^{l}} \biggr), \quad s \in [\crit, \tau],
		\end{split}
	\end{equation} 
	for any $k \in \indset{n+1, \maxK'-2l+1}$, where $\tau$ is  the stopping time defined in \eqref{eq:trless_tau2_def}.
\end{claim}
Here we recall that $\norm{\cdot}_*$ of a $\bm x = (\bm x_1^*, x, \bm x_2, y)$ dependent function $f(\bm x)$ is given by \eqref{eq:star_norm}, and the estimate on the isotropic quantity is understood to be uniform in $a,b$.
We defer the proof of Claims~\ref{claim:Q_forcing}--\ref{claim:Q_mart} until the end of the present section. 
Note that in \eqref{eq:Q_mart}--\eqref{eq:Q_forcing}, there is no improvement in terms of the $(N\eta_s)$-power 
for the averaged quantities compared with the isotropic ones, in contrast to \eqref{eq:tr_av_focing}--\eqref{eq:tr_iso_focing} and \eqref{eq:tr_mart_estimates}. For the martingale terms, this difference arises because, when estimating $\mathrm{d}\mathcal{E}_{k}^{n, \mathrm{av/iso} }$, we perform the maximal possible number of reductions\footnote{
Without any reductions, the optimal averaged bound in \eqref{eq:Q_mart} would contain $(N\eta_s)^{4k}$, and the optimal isotropic bound would contain $(N\eta_s)^{4k-1}$,  but each reduction incurs a loss of $ N\eta_s$. 
} 
(without exceeding the natural size of the corresponding $Q_{k,s}$) to compartmentalize and preserve all traceless observables. In contrast, when estimating $\mathrm{d}\mathcal{M}_{k,s}^{n, \mathrm{av/iso} }$, no reductions were performed, allowing us to achieve the optimal bound in terms of the $(N\eta_s)$-power. For the forcing terms $\mathcal{W}_{k,s}^{n, \mathrm{av/iso} }$, the lack of improvement in the averaged case stems from the fact that the corresponding $M$-terms are not subtracted from the $Q$-chains in \eqref{eq:Q_quantities}.  Finally, notice  the additional $N^{-\nu/2}$ powers 
in \eqref{eq:Q_mart}--\eqref{eq:Q_forcing} compared to the target in the definition
of the stopping time~\eqref{eq:trless_tau2_def}, guaranteeing the self-improvement mechanism. 

Applying Duhamel's principle to \eqref{eq:Q_evol}, using the  first two
 bounds \eqref{eq:trless_initial} and \eqref{eq:Q_mart}--\eqref{eq:Q_forcing} to estimate the initial condition at $t = \crit$,  the martingale and the forcing terms, and using the analog of \eqref{eq:tr_prop} to bound the action of the propagators, we conclude that 
\begin{equation}
	\bigl\lvert \mathcal{L}_t^{k,n} \bigr\rvert + N^{-\nu}\bigl\lvert \mathcal{K}_t^{k,n}\bigr\rvert \prec N^{-\nu/4}  \frac{N^{\xi + k\nu}\trl_t^{2n}}{(N\eta_t)^{2k-1}} \biggl( 1   +\frac{1}{\trl_t^2(N\eta_t)^{l}} \biggr), \quad t \in [\crit, \tau].
\end{equation}
Here we also used the integration rules \eqref{eq:int_rules} (see the proof of  Proposition \ref{prop:masters} in Section \ref{sec:masters_proof} for more details). Therefore, the stopping time $\tau$ defined in \eqref{eq:trless_tau2_def}, satisfies $\tau = T$ with very high probability. This concludes the proof of Lemma \ref{lemma:Q_bounds}.
\end{proof}

Therefore, it remains to establish Claims \ref{claim:Q_mart}--\ref{claim:Q_forcing}. We begin by proving Claim \ref{claim:Q_mart}.
\begin{proof} [Proof of Claim \ref{claim:Q_mart}] 
	First, we prove the bound on $\mathrm{d}\mathcal{E}_{k}^{n, \mathrm{av}}$ in \eqref{eq:Q_mart}. By \eqref{eq:Q_mart_def}, we have
	\begin{equation} \label{eq:Eps_av_sum}
		 \bigl[ \mathrm{d}\mathcal{E}^{n, \mathrm{av} }_{k} \bigr]_s \lesssim \sum_{j=1}^k \sum_{ab} S_{ab} \biggl(\bigl\lvert \bigl( G_{[j,k],s} S^y G_{[1,k],s}^\dagger S^{x} G_{[1,j],s}\bigr)_{ab}\bigr\rvert^2 + \bigl\lvert \bigl( G_{[j,k],s}^\dagger S^y G_{[1,k],s} S^{x} G_{[1,j],s}^\dagger\bigr)_{ab}\bigr\rvert^2\biggr).
	\end{equation}
	We analyze the first term in the large parenthesis $(\dots)$ in detail. The other term is treated completely analogously.
	We write the first term as follows:
	\begin{equation} \label{eq:Eps_av2_id}
		\begin{split}
			\sum_{ab} S_{ab} &\bigl\lvert \bigl( G_{[j,k],s} S^y G_{[1,k],s}^\dagger S^{x} G_{[1,j],s}\bigr)_{ab}\bigr\rvert^2\\ 
			&  = \sum_a  
			\big( G_{[j,k],s} S^y G_{[1,k],s}^\dagger S^{x} G_{[1,j],s}
			S^a  G_{[1,j],s}^* S^x \big(G_{[1,k],s}^\dagger\big)^* S^y G_{[j,k],s}^* )_{aa}, 
			\end{split}
			\end{equation}
we  organize it into shorter chains with  $Q$-terms  using the reduction inequality~\eqref{eq:reduction} 
and estimate them using the two available  $Q$-estimates. Recall that
the stronger bound~\eqref{eq:lesser_Q_bounds}   is valid only if $Q\in\mathfrak{Q}_{k',n'}$ 
for $n'<n$, i.e. (half of) $Q$ contains strictly less than $n$ traceless observables, while the weaker bound, 
from~\eqref{eq:traceless_Phi_def}--\eqref{eq:trless_tau2_def}, for $s\le \tau$,
is available for $n$ traceless observables. We have the constraint that reduction inequality can only be used at some of the $S^x, S^y, S^a$
explicitly written in~\eqref{eq:Eps_av2_id}, the observables inside the $G$-chains are untouchable as they 
may be traceless. Furthermore, we need to use reduction sparingly because 
 each reduction step loses a factor $(N\eta)$.

	It turns out that the  $j = 1$ and $j=k$ cases are critical in the sense that they explicitly produce a chain of length $k+1$ that can not be further reduced without giving up the $N^{-\nu}$ improvement 
	necessary to close the stopping time argument (see discussion below \eqref{eq:Eps_av2}). However, this longer chain can be estimated using the weaker bound \eqref{eq:power_induction_assume} from the $l$-induction assumption. 
	Hence, we now distinguish the following cases:
	
	{\bf Case 1.} If $j=k$,  denoting  $Q_{k,s}:=G_{[1,k],s}^\dagger S^{x} G_{[1,k],s}$,  
	 we use \eqref{eq:reduction} splitting the chain~\eqref{eq:Eps_av2_id} at $S^x$
		    to deduce that, for $s \le \tau$,
		\begin{equation} \label{eq:Eps_av1}
			\begin{split}
				\sum_{ab} S_{ab} \bigl\lvert \bigl( G_{k,s} S^y Q_{k,s}\bigr)_{ab}\bigr\rvert^2 
				&\lesssim \sum_{ac} S_{xc} \bigl\lvert \bigl(G_{k,s} S^y G_{[1,k],s}^\dagger\bigr)_{ac}\bigr\rvert^2 \Tr \bigl[G_{[1,k],s}^* S^x  G_{[1,k],s}S^a\bigr] \\
				&\prec N\frac{\trl_s^{2n}}{(N\eta_s)^{2k}} \biggl(\frac{1}{N\eta_s}+ \frac{1}{\trl_s^2 (N\eta_s)^{l}} \biggr) \max_a  \norm{\Tr \bigl[ Q_{k,s}S^a\bigr]}_*,\\
				&\prec \frac{1}{\eta_s}\frac{ N^{\xi+k\nu} \trl_s^{4n}}{(N\eta_s)^{4k-2}} \biggl(1 + \frac{1}{\trl_s^2 (N\eta_s)^{l}} \biggr)^2,
			\end{split}
		\end{equation}
		where in the penultimate step we used the $M$-bound from Lemma \ref{lemma:traceless_M_bounds} together with the main assumption~\eqref{eq:power_induction_assume} for $\Psi^{n,\mathrm{iso}}_{k+1,s}$ to estimate the isotropic
		chain $\bigl(G_{k,s} S^y G_{[1,k],s}^\dagger\bigr)_{ac}$, 
		while in the last step we used \eqref{eq:trless_tau2_def} to bound the symmetric $2k$-chain $Q_{k,s}$ directly. 
	The case $j=1$ is estimated analogously.
	
	{\bf Case 2.} If $2 \le j \le k-1$,   then we consider two sub-cases, depending on how the $n$ traceless observables
	are distributed between the chains $G_{[1,j],s}$ and $G_{[j,k],s}$. The situation when one of the chains contains
	all traceless observables requires a separate treatment since for such chains~\eqref{eq:lesser_Q_bounds}
	is not applicable and we need to use the weaker bound coming 
	from~\eqref{eq:traceless_Phi_def}--\eqref{eq:trless_tau2_def} when $s\le \tau$. 
	 
	 {\bf Case 2.1.} If  $G_{[1,j],s} \in \mathfrak{G}_{j,n}$, 
	 i.e. all traceless observables fall into the $G_{[1,j],s}$ sub-chain,     
	 then $G_{[j,k],s} \in \mathfrak{G}_{k-j+1,0}$, and we use \eqref{eq:reduction} twice, splitting the chain~\eqref{eq:Eps_av2_id} at both $S^x$ and $S^y$, to obtain
	\begin{equation} \label{eq:Eps_av2}
		\begin{split}
			\sum_{ab} S_{ab} &\bigl\lvert \bigl( G_{[j,k],s} S^y G_{[1,k],s}^\dagger S^{x} G_{[1,j],s}\bigr)_{ab}\bigr\rvert^2\\  
			&\lesssim \sum_{a} \bigl(G_{[j,k],s} S^y G_{[j,k],s}^*\bigr)_{aa} \Tr\bigl[ G_{[1,j],s}^* S^x  G_{[1,j],s} S^a\bigr] \Tr\bigl[\bigl(G_{[1,k],s}^\dagger\bigr)^*  S^y G_{[1,k],s}^\dagger S^{x}\bigr]\\
			&\prec  \frac{1}{\eta_s}\frac{N^{2\xi + (k+j) \nu}\trl_s^{4n}}{(N\eta_s)^{4k-2}}\biggl(1 + \frac{1}{\trl_s^4(N\eta_s)^{2l} }\biggr),  
		\end{split}
	\end{equation}
	where we used the definition of $\tau$ in \eqref{eq:trless_tau2_def}  to estimate the two traces,
	 and \eqref{eq:lesser_Q_bounds}
	  to estimate the  $(G_{[j,k]}S^yG_{[j,k]}^*)_{aa}$ term.  
	The opposite case, when all traceless observables fall into the $G_{[j,k],s}$ chain, i.e. 
	$G_{[j,k],s} \in \mathfrak{G}_{k-j+1,n}$, hence
	 $G_{[1,j],s} \in \mathfrak{G}_{j,0}$,
	 is handled analogously. 
	 
	 We remark that the case separation between Case $1$ and Case $2$ is really necessary. Indeed, using the bound \eqref{eq:Eps_av2} with $j = k$ would yield no $N^{-\nu}$ improvement. On the other hand, using the strategy from Case $1$ to estimate the left-hand side of \eqref{eq:Eps_av2} with $j \ge 2$ would produce chains that are much longer than $k+1$ and can not be estimated by the $l$-induction hypothesis~\eqref{eq:power_induction_assume}.

	{\bf Case 2.2.} If $2 \le j \le k-1$ and $G_{[1,j],s} \in \mathfrak{G}_{j,n'}$ with $1 \le n' \le n-1$, then
	we can estimate the second line of~\eqref{eq:Eps_av2} differently as
	\begin{equation} \label{eq:Eps_av3} 
			\sum_{ab} S_{ab} \bigl\lvert \bigl( G_{[j,k],s} S^y G_{[1,k],s}^\dagger S^{x} G_{[1,j],s}\bigr)_{ab}\bigr\rvert^2  
			\prec  \frac{1}{\eta_s}\frac{N^{\xi + k\nu} \trl_s^{4n}}{(N\eta_s)^{4k-2}}\biggl(1 + \frac{1}{\trl_s^2(N\eta_s)^{l} }\biggr),  
	\end{equation}
	using \eqref{eq:lesser_Q_bounds} for the first two factors in the second line of~\eqref{eq:Eps_av2} and \eqref{eq:trless_tau2_def} 
	for the last trace. Note that this bound is better than~\eqref{eq:Eps_av1} from Case 1
	and~\eqref{eq:Eps_av2} from Case 2.1 since the stronger $Q$-bound~\eqref{eq:lesser_Q_bounds} could be applied
	twice.
	Combining \eqref{eq:Eps_av_sum}--\eqref{eq:Eps_av3}, we obtain the desired first bound in \eqref{eq:Q_mart}.

	Next, we prove the bound on $\mathrm{d}\mathcal{E}_{k}^{n, \mathrm{iso}}$ in \eqref{eq:Q_mart}. The proof follows the same outline; by \eqref{eq:Q_mart_def}, we obtain
	\begin{equation} \label{eq:Eps_iso_sum}
		\begin{split}
			\bigl[ \mathrm{d}\mathcal{E}^{n, \mathrm{iso} }_{k} \bigr]_s \lesssim \sum_{j=1}^k \sum_{cd} S_{cd} \biggl(&\bigl\lvert \bigl( G_{[j,k],s}\bigr)_{cb} \bigl( G_{[1,k],s}^\dagger S^{x} G_{[1,j],s}\bigr)_{ad}\bigr\rvert^2 \\
			&+ \bigl\lvert \bigl( G_{[j,k],s}^\dagger\bigr)_{ca} \bigl(G_{[1,k],s} S^{x} G_{[1,j],s}^\dagger\bigr)_{bd}\bigr\rvert^2\biggr).
		\end{split}
	\end{equation}
	With $\tau$ being the stopping time defined in \eqref{eq:trless_tau2_def}, let $R^{(j)}_{k,s}$ denote 
	\begin{equation} \label{eq:Rchain_def}
		R^{(j)}_{k,s} := \sum_{cd} S_{cd} \bigl\lvert \bigl( G_{[j,k],s}\bigr)_{cb} \bigl( G_{[1,k],s}^\dagger S^{x} G_{[1,j],s}\bigr)_{ad}\bigr\rvert^2, \quad s \in [\crit, \tau].
	\end{equation}
	
	Once again, we need to separate the critical case $j=1$ for the same reason as in the proof of the averaged bound above. However, in the isotropic case, the additional $N^{-\nu}$ factor in the control on the averaged $Q$-chains in \eqref{eq:trless_tau2_def} allows us to treat $j=k$ the same way as $2 \le j \le k-1$.
	
	For $j=1$, we  do not use any reductions and obtain for $R^{(1)}$, defined in \eqref{eq:Rchain_def}, that
	\begin{equation}
		\begin{split}
			R^{(1)}_{k,s}
			&= \sum_{d}   \bigl(G_{[1,k],s}^* S^d G_{[1,k],s}\bigr)_{bb} \bigl\lvert  \bigl( G_{[1,k],s}^\dagger S^{x} G_{1,s}\bigr)_{ad}\bigr\rvert^2\\
			&\le \max_{d} \bigl(G_{[1,k],s}^* S^{d} G_{[1,k],s}\bigr)_{bb}   \sum_{d}   \bigl\lvert  \bigl( G_{[1,k],s}^\dagger S^{x} G_{1,s}\bigr)_{ad}\bigr\rvert^2\\
			&\prec  \frac{1}{\eta_s} \frac{N^{\xi+(k+1)\nu}\trl_s^{4n}}{(N\eta_s)^{4k-2}}
			\biggl(1 + \frac{1}{\trl_s^2(N\eta_s)^{l}}\biggr)^2,
		\end{split}
	\end{equation}
	where we used \eqref{eq:trless_tau2_def} for $(G_{[1,k],s}^* S^{q} G_{[1,k],s})_{bb}$, and the main assumption~\eqref{eq:power_induction_assume} together with Lemma \ref{lemma:traceless_M_bounds} to estimate  $\lvert  ( G_{[1,k],s}^\dagger S^{x} G_{1,s})_{ad} \rvert$.
	
	Next, we consider the case $2 \le j \le k$.	 Using \eqref{eq:sumS=1} to insert $\sum_q S^q=I$ as in \eqref{eq:inset_S},  we deduce that the quantities $R^{(j)}$, defined in \eqref{eq:Rchain_def}, satisfy
	\begin{equation} \label{eq:Eps_iso2}
		\begin{split}
			R^{(j)}_{k,s}
			\le&~ \max_d\bigl(G_{[j,k],s}^* S^d G_{[j,k],s}\bigr)_{bb}\sum_{q}   \bigl( G_{[1,k],s}^\dagger S^{x} G_{[1,j],s} S^{q} G_{[1,j],s}^* S^x  (G_{[1,k],s}^\dagger)^* \bigr)_{aa}		\\
			\le&~ \max_d\bigl(G_{[j,k],s}^* S^d G_{[j,k],s}\bigr)_{bb} \bigl( G_{[1,k],s}^\dagger  S^x  (G_{[1,k],s}^\dagger)^* \bigr)_{aa} \sum_{q}   \Tr \bigl[S^{x} G_{[1,j],s} S^{q} G_{[1,j],s}^*\bigr]		\\ 
			 \prec&~ \frac{1}{\eta_s} \frac{N^{2\xi+(2k+1)\nu}\trl_s^{4n}}{(N\eta_s)^{4k-2}}\biggl(1 + \frac{1}{\trl_s^2(N\eta_s)^{l}}\biggr)^2.
		\end{split}
	\end{equation}
	Here in the second step we used the reduction inequality \eqref{eq:reduction} to split the chain at $S^x$, while in the last step we  used \eqref{eq:trless_tau2_def} to estimate $Q$-chains with $n$ traceless observables, and \eqref{eq:lesser_Q_bounds} for the $Q$-chains with
	 $n' \le n-1$ traceless observables. Note that the weaker bound~\eqref{eq:trless_tau2_def} 
	 is used at most twice since we have altogether $2n$ traceless observables. 
	 Furthermore, since the chain $G_{[j,k]}$ is strictly shorter than $G_{[1,k]}$, we gain at least a factor of $N^{-\nu}$.  Hence, combining \eqref{eq:Eps_iso_sum}--\eqref{eq:Eps_iso2}, we obtain the desired second bound in \eqref{eq:Q_mart}.
	This concludes the proof of Claim \ref{claim:Q_mart}. 
\end{proof}

We close this section by proving Claim \ref{claim:Q_forcing}.
\begin{proof}[Proof of Claim \ref{claim:Q_forcing}]
	First, we prove the estimate on $\mathcal{W}_{k,s}^{n, \mathrm{av} }$ in \eqref{eq:Q_forcing}. Recall from \eqref{eq:Q_forcing_def} that $\mathcal{W}_{k,s}^{n, \mathrm{av} }$ is comprised of two types of terms: direct and cross,
	 namely, $\mathcal{W}_{k}^{\mathrm{dir}}$   traced against $S^x G_{[1,k]} S^y$ 
	and $\mathcal{W}_{k}^{\mathrm{crs}}$,    traced against $S^y$.  
	We begin by estimating the direct terms, defined in \eqref{eq:direct_forcing}. It follows from Schwarz inequality that
	\begin{equation} \label{eq:dir_norm}
		\norm{ \Tr\bigl[G_{[1,k],s}^\dagger S^x \mathcal{W}_{k,s}^{\mathrm{dir} } S^y \bigr]}_* \le \bigl\lVert G_{[1,k],s}^\dagger\bigr\rVert_{\mathrm{dir}} \norm{\mathcal{W}_{k,s}^{\mathrm{dir} }}_{\mathrm{dir}}, \quad  \norm{X}_{\mathrm{dir}} :=  \norm{ \Tr\bigl[X^* S^x X S^y \bigr]}_*^{1/2},
	\end{equation}
	for any $\mathbb{C}^{N\times N}$-valued function $X$ of $\bm x$.  Note that
	\begin{equation} \label{eq:dir_Phi}
		  \norm{G_{[1,k],s} }_{\mathrm{dir}}^2 =
		\bigl\lVert \Tr \bigl[G_{[1,k],s}^*S^x G_{[1,k],s} S^y\bigr]\bigr\rVert_* \le \frac{\trl_s^{2n}}{(N\eta_s)^{2k-1}}\biggl(1 + \frac{1}{\trl_s^2(N\eta_s)^{l}}\biggr) \Phi_{k,s}^{n,\mathrm{av}},
	\end{equation}	
	hence it can be estimated using the stopping time $\tau$ from \eqref{eq:trless_tau2_def}.
	Therefore, it suffices to estimate $\lVert  \mathcal{W}_{k,s}^{\mathrm{dir} }\rVert_{\mathrm{dir}}$. 
	
	We proceed by estimating the terms on the right-hand side of \eqref{eq:direct_forcing} one-by-one, 
	starting with the $j=i$ term, 
	\begin{equation} \label{eq:av_dir_forcing}
		\begin{split}			
			\norm{G_{[1,i],s} \mathcal{S}\bigl[(G-M)_{i,s}\bigr] G_{[i,k],s} }_{\mathrm{dir}}^2 
			&= \sum_{ab} S_{xa}S_{yb} 	\biggl\lvert \sum_{c} \bigl(G_{i,s}-m_{i,s}\bigr)_{cc} \bigl(G_{[1,i],s} S^c G_{[i,k],s}\bigr)_{ab} \biggr\rvert^2	\\
			&\prec \frac{1}{N\eta_s}\sum_{ab} S_{xa}S_{yb} 	\biggl( \sum_{c} \bigl\lvert\bigl(G_{[1,i],s} S^c G_{[i,k],s}\bigr)_{ab}\bigr\rvert \biggr)^2	\\
			&\prec \frac{N^2}{ N\eta_s}\sum_{ab} S_{xa}S_{yb}\frac{\trl_s^{2n}}{(N\eta_s)^{2k}}\biggl(\delta_{ab} + \frac{1}{N\eta_s} + \frac{1}{\trl_s^2(N\eta_s)^{l}}\biggr)\\
			&\prec \frac{1}{ \eta_s^2}\frac{\trl_s^{2n}}{(N\eta_s)^{2k-1}}\biggl(\frac{1}{W} + \frac{1}{N\eta_s} + \frac{1}{\trl_s^2(N\eta_s)^{l}}\biggr),
		\end{split}
	\end{equation}	
	where we used \eqref{eq:trless_initial} for $(G-m)_{cc}$, and \eqref{eq:trless_ind_assume} together with Lemma \ref{lemma:traceless_M_bounds} to bound the chain $(G_{[1,i],s} S^c G_{[i,k],s})_{ab}$. 
	
	Next, we consider  the terms on the right-hand side of \eqref{eq:direct_forcing} with $j = i+1$, i.e. 
	\begin{equation} \label{eq:dir_forcing_prop}
	 	 G_{[1,i],s} \mathcal{S}\bigl[(G-M)_{[i,i+1],s}\bigr] G_{[i+1,k],s}.
	\end{equation}
	Estimating every entry of $(G-M)_{[i,i+1],s}$ by its $\norm{\cdot}_*$-norm,  we obtain
	\begin{equation} \label{eq:av_dir_forcing2}
		\begin{split}
		\bigl\lVert \, (\ref{eq:dir_forcing_prop}) \, \bigr\rVert_{\mathrm{dir}}^2 
			&\le \norm{(G-M)_{[i,i+1],s}}_*^2 \sum_{ab}  S_{xa} S_{yb} 	 \biggl(\sum_c\bigl\lvert \bigl(G_{[1,i],s} S^c G_{[i+1,k],s}\bigr)_{ab} \bigr\rvert\biggr)^2 \\
			&\le N\norm{(G-M)_{[i,i+1],s}}_*^2 \sum_{c} \sum_{ab}  S_{xa} S_{yb} 	  \bigl\lvert \bigl(G_{[1,i],s} S^c G_{[i+1,k],s}\bigr)_{ab} \bigr\rvert^2 \\
			&\le N^2\norm{(G-M)_{[i,i+1],s}}_*^2   \norm{ G_{[1,i],s} S^{x_i} G_{[i+1,k],s} }_{\mathrm{dir}}^2  \\ 
			&\prec \frac{1}{\eta_s^2} \frac{\trl_s^{2n}}{(N\eta_s)^{2k-1}} \biggl(\frac{N^{\xi+k\nu}}{N\eta_s} + \frac{1}{\trl_s^2 (N\eta_s)^{l}}\biggr).
		\end{split}
	\end{equation}
	In the second step we used Schwarz inequality, and, to bound $(G-M)_{[i,i+1],s}$ in the fourth step, we used
	 \eqref{eq:trless_ind_assume}  if $G_{[i,i+1],s} \in \mathfrak{G}_{2,n'}$ for some $n' \le n-1$ and \eqref{eq:power_induction_assume} if $G_{[i,i+1],s} \in \mathfrak{G}_{2,n}$ and $n=1$.   Here, to bound the second factor in 
	 the third line, 	   we used \eqref{eq:dir_Phi} together with \eqref{eq:trless_tau2_def},
	 noticing  that
	    if $G_{[i,i+1]} \in \mathfrak{G}_{2,0}$ then
	    $G_{[1,i],s} S^{x_i} G_{[i+1,k],s}$ has  $n$ traceless observables, while we
	    used \eqref{eq:lesser_Q_bounds} in the case $G_{[1,i],s} S^{x_i} G_{[i+1,k],s} \in \mathfrak{G}_{2,1}$, i.e. if	  $G_{[1,k],s}$
	     has only $n-1$ traceless observables.

	Finally, for the remaining terms on the right-hand side of \eqref{eq:Q_forcing_def} with $j \ge i+2$, we obtain, similarly to \eqref{eq:av_dir_forcing2} above,
	\begin{equation} \label{eq:av_dir_forcing3}
		\begin{split}
			\norm{ G_{[1,i],s} \mathcal{S}\bigl[G_{[i,j],s}\bigr] G_{[j,k],s} }_{\mathrm{dir}}^2 &\le  N^2\norm{G_{[i,j],s}}_*^2   \norm{ G_{[1,i],s} S^{x_i} G_{[j,k],s} }_{\mathrm{dir}}^2\\
			&\prec \frac{1}{\eta_s^2}\frac{N^{\xi + (k+i-j+1)\nu }\trl_s^{2n}}{(N\eta_s)^{2k-1}}\biggl(1 + \frac{1}{\trl_s^2(N\eta_s)^l}\biggr),
		\end{split}
	\end{equation}
	where we used Lemma \ref{lemma:traceless_M_bounds} together with \eqref{eq:power_induction_assume} or \eqref{eq:trless_ind_assume}  to estimate $G_{[i,j]}$, and  \eqref{eq:trless_tau2_def}  or  \eqref{eq:lesser_Q_bounds} to bound the remaining $\norm{G_{[1,i],s} S^{x_i} G_{[j,k],s}}_{\mathrm{dir}}$,  depending on whether these chains contain  exactly $n$ or less traceless observables, respectively.  Clearly, only at most one of the weaker bounds \eqref{eq:power_induction_assume}
	and  \eqref{eq:trless_tau2_def} has to be  used for any fixed $i,j$.
	   Combining \eqref{eq:av_dir_forcing}--\eqref{eq:av_dir_forcing3}, we conclude that the direct forcing terms satisfy
	\begin{equation} \label{eq:dri_forcing_est}
		\norm{ \Tr\bigl[G_{[1,k],s}^\dagger S^x \mathcal{W}_{k,s}^{\mathrm{dir} } S^y \bigr]}_* \prec \frac{1}{ \eta_s^2}\frac{N^{\xi + (k-1/2)\nu }\trl_s^{2n}}{(N\eta_s)^{2k-1}}\biggl(1 + \frac{1}{\trl_s^2(N\eta_s)^{l}}\biggr),
	\end{equation}
	where we used \eqref{eq:trless_tau2_def} and \eqref{eq:dir_norm}.
	
	Next, we estimate the cross terms, defined in \eqref{eq:cross_forcing} 
	after taking the trace against $S^y$.  Then the two terms in the first line of \eqref{eq:cross_forcing} become structurally
	 identical;  We consider the first one,
	 \begin{equation} \label{eq:cross_prop}
	 	\Tr\bigl[G_{[1,k],s}^\dagger \mathcal{S}\bigl[G_{1,s}^* S^xG_{1,s} - \Theta^x_s\bigr]G_{[1,k],s} S^y\bigr] = \sum_q \bigl(G_{1,s}^* S^xG_{1,s} - \Theta^x_s\bigr)_{qq} \Tr\bigl[G_{[1,k],s}^\dagger S^qG_{[1,k],s} S^y\bigr].
	 \end{equation}
	 Estimating the $2$-chain   $G^*S^xG-\Theta^x$ with no regular observables by 
	 the isotropic bound in \eqref{eq:trless_initial}, and using \eqref{eq:trless_tau2_def} for the remaining trace of a $Q$-chain, 	 we obtain the following bound for the first term on the right-hand side of \eqref{eq:cross_forcing},
	\begin{equation}
		\bigl\lvert \Tr\bigl[G_{[1,k],s}^\dagger \mathcal{S}\bigl[G_{1,s}^* S^xG_{1,s} - \Theta^x_s\bigr]G_{[1,k],s} S^y \bigr] \bigr\rvert \prec \frac{1}{\eta_s}\frac{N^{\xi+k\nu}}{\sqrt{N\eta_s}}\frac{\trl_s^{2n}}{(N\eta_s)^{2k-1}}\biggl(1 + \frac{1}{\trl_s^2 (N\eta_s)^{l}}\biggr).
	\end{equation}
	The second term on the right-hand side of \eqref{eq:cross_forcing} is estimated analogously.
	
	Next, we bound the terms in the second line of \eqref{eq:cross_forcing}. The super-operator $\mathcal{S}$ together with the observables $S^x$ and $S^y$ divide the two chains (separated by $\mathcal{S}$) into four segments of lengths $k-i+1, i,  j$, and $k-j +1$.
	Note that the lengths of all these segments are at most $k$. The most critical terms (see Footnote~\ref{note:2k} below for further details) are those where two of them are exactly $k$ (hence the other two are equal to one)
	 Note that the cases $i=j=k$ and $k-i+1 = k-j+1 = k$ are explicitly excluded by the condition $3\le i+j\le 2k-1$. Hence, it remains to consider the cases $i=1, j=k$ and $i=k, j=1$. 
	
	For $i=1, j=k$, using Lemma \ref{lemma:traceless_M_bounds} and \eqref{eq:power_induction_assume} to bound both $(k+1)$-chains $G_{1,s}^\dagger S^x G_{[1,k],s}$ and  $G_{k,s} S^y G_{[1,k],s}^\dagger$ (one in averaged
	sense and one in isotropic sense), we obtain
	\begin{equation}
		\bigl\lvert \Tr \bigl[G_{[1,k],s}^\dagger \mathcal{S}\bigl[G_{1,s}^\dagger S^x G_{[1,k],s}\bigr]G_{k,t} S^y\bigr] \bigr\rvert \prec \frac{1}{\eta_s} \frac{\trl_s^{2n}}{(N\eta_s)^{2k-1}}\biggl(1 + \frac{1}{\trl_s^2 (N\eta_s)^{l+1/2}}\biggr).
	\end{equation}
	The other case $i=k, j=1$ is completely analogous.

	For all remaining terms, where at most one segment has length $k$, we can divide the chain into four segments by a Schwarz inequality and estimate them individually to obtain
	\begin{equation}
		\begin{split}
			\bigl\lvert \Tr\bigl[& G_{[i,k],s}^\dagger \mathcal{S}\bigl[G_{[1,i],s}^\dagger S^x G_{[1,j],s}\bigr]G_{[j,k],s} \bigr]	\bigr\rvert \\
			=&~ \biggl\lvert \sum_{abcd} S_{ya} S_{bc}  S_{dx} \bigl(G_{[i,k],s}^\dagger\bigr)_{ab}  \bigl(G_{[1,i],s}^\dagger\bigr)_{cd}  \bigl(G_{[1,j],s}\bigr)_{dc} \bigl(G_{[j,k],s}\bigr)_{ba}  \biggr\rvert \\ 
			\lesssim &~   \biggl\lvert\sum_{abcd} S_{ya} S_{bc}  S_{dx}     \bigl\lvert \bigl(G_{[1,j],s}\bigr)_{dc} \bigr\rvert^2 \bigl\lvert \bigl(G_{[j,k],s}\bigr)_{ba}\bigr\rvert^2\biggr\rvert^{1/2}       \biggl\lvert\sum_{abcd} S_{ya} S_{bc}  S_{dx} \bigl\lvert \bigl(G_{[i,k],s}^\dagger\bigr)_{ab} \bigr\rvert^2 \bigl\lvert \bigl(G_{[1,i],s}^\dagger\bigr)_{cd} \bigr\rvert^2 \biggr\rvert^{1/2} \\ 
			\lesssim &~   N \sqrt{\bigl\lVert Q^{(1)}_{i,s}\bigr\rVert_* \, \bigl\lVert Q^{(2)}_{k-i+1,s}\bigr\rVert_* \, \bigl\lVert Q^{(3)}_{j,s}\bigr\rVert_* \, \bigl\lVert Q^{(4)}_{k-j+1,s}\bigr\rVert_* },
		\end{split}
	\end{equation}
	where each $\{Q^{(p)}_{k_p,s}\}_{p=1}^4$  belongs to $\mathfrak{Q}_{k_p, n_p}$ with $k_1 = i$, $k_2 := k-i+1$, 
	$k_3:=j$, $k_4 := k-j+1$ and $n_p \le n$ with $\sum_{p=1}^4 n_p = 2n$.
	Therefore, by \eqref{eq:trless_tau2_def}, which is used at most twice
	for two chains that have $n$ traceless observables with combined length\footnote{\label{note:2k}if the critical cases where two chains have length exactly $k$ were not explicitly excluded this combined length would be $2k$, and hence the stopping time argument would not close due to $N^{k\nu}$. } at most $2k-1$,
	and by Lemma \ref{lemma:traceless_M_bounds} together with
	 \eqref{eq:trless_ind_assume}  for the remaining chains we conclude that
	\begin{equation}
		\bigl\lvert \Tr\bigl[\dots \, S^y\bigr] \bigr\rvert \prec \frac{1}{\eta_s} \frac{N^{\xi+(k-1/2)\nu}\trl_s^{2n}}{(N\eta_s)^{2k-1}}\biggl(1 + \frac{1}{\trl_s^2 (N\eta_s)^{l}}\biggr).
	\end{equation}
	Hence, the cross forcing terms satisfy
	\begin{equation} \label{eq:cross_forcing_est}
		\norm{  \Tr\bigl[\mathcal{W}_{k,s}^\mathrm{crs} S^y\bigr]}_* \prec \frac{1}{\eta_s} \frac{N^{\xi+(k-1/2)\nu}\trl_s^{2n}}{(N\eta_s)^{2k-1}}\biggl(1 + \frac{1}{\trl_s^2 (N\eta_s)^{l}}\biggr).
	\end{equation}
	Combining \eqref{eq:dri_forcing_est} and \eqref{eq:cross_forcing_est} with \eqref{eq:Q_forcing_def}, we conclude the first bound in \eqref{eq:Q_forcing}.
	
	The proof of isotropic bound in \eqref{eq:Q_forcing} is largely analogous to that of the averaged bound proved above. 
	The last term in the definition of $\mathcal{W}_{k,t}^{n, \mathrm{iso} }$ in \eqref{eq:Q_forcing_def} admits the bound
	\begin{equation}
		\norm{|m_s|^2 \sum_q \bigl(I+\Theta_s\bigr)_{aq}\Tr\bigl[Q_{k,s}S^q\bigr]}_* \lesssim \frac{N^{\xi+k\nu}  \trl_s^{2n}}{(N\eta_s)^{2k-1}}\biggl(1 + \frac{1}{\trl_s^2 (N\eta_s)^{l}}\biggr).
	\end{equation}
	by the first bound in \eqref{eq:supercrit_Theta} and the averaged bound on $Q_{k,s}$ from \eqref{eq:trless_tau2_def}. While the critical terms from the second line of \eqref{eq:cross_forcing} are bounded by
	\begin{equation}
		 \norm{G_{[1,k],s}^\dagger \mathcal{S}\bigl[G_{1,s}^\dagger S^x G_{[1,k],s}\bigr]G_{k,s} }_*  \prec \frac{1}{\eta_s} \frac{N^{\xi/2 + (k+1)\nu/2}\trl_s^{2n}}{(N\eta_s)^{2k-1}}\biggl(1 + \frac{1}{\trl_s^2 (N\eta_s)^{l}}\biggr),
	\end{equation}
	where we used Lemma \ref{lemma:traceless_M_bounds} together with
	 \eqref{eq:power_induction_assume} for the averaged $(k+1)$-chain, to estimate the remaining $k$- and $1$- chains, we use Schwarz inequality and \eqref{eq:trless_tau2_def} together with \eqref{eq:lesser_Q_bounds}, respectively. 
	  The other terms are estimated analogously to their counterparts in the averaged case above, using the isotropic bound from \eqref{eq:trless_tau2_def}. 
	This concludes the proof of Claim \ref{claim:Q_forcing}.
	
\end{proof}
	
	\subsection{Traceless $M$-bounds: Proof of Lemma \ref{lemma:traceless_M_bounds}} \label{sec:trless_M}
	We prove  Lemma~\ref{lemma:traceless_M_bounds} using the same dynamical approach as in the proof of Lemma~\ref{lemma:M_bounds}. However, since in Lemma~\ref{lemma:traceless_M_bounds} we only consider the times $t\in[\crit, T]$, and hence  $\eta_t$ below the critical scale $(W/N)^2$, the proof simplifies greatly. Indeed, as \eqref{eq:simplesizefn} suggests, the spatial decay for such $\eta_t$ is completely lost. 
	
	Furthermore, the propagators corresponding to traceless observables are essentially bounded, see \eqref{eq:tr_prop_decomp}--\eqref{eq:tr_prop}, meaning they do not contribute the usual large $\eta_s/\eta_t$-factor to the integral in the Duhamel's formula.
	
	\begin{proof}[Proof of Lemma \ref{lemma:traceless_M_bounds}]
		It suffices to prove \eqref{eq:trlessM_iso} and \eqref{eq:trlessM_av} for chains that contain at least one traceless observable, since the other case has already been established in Lemma~\ref{lemma:M_bounds}.	Furthermore, since the $M$-objects are diagonal it suffices to consider $a=b$, since $M_{ab} = 0$ trivially for $a\neq b$.
		
		Observe that the bounds \eqref{eq:trlessM_iso}--\eqref{eq:trlessM_av} holds trivially for $k=1$, since $M_{[1,1],t} = m(z_1,t) I$, and hence vanishes when tested against a traceless observable.
		Hence, we proceed by an induction in length $k\ge 2$.
		\vspace{5pt}
		
		\textbf{Base case}. We prove the base case $k=2$ and $n=1$. 
		First, if $M_{[1,2],t}$ is not saturated, i.e., $z_1 = z_2$,  then, it follows from the second bound in \eqref{eq:Ups_majorates} and \eqref{eq:Ups_norm_bounds} that
		\begin{equation} \label{eq:trless2_iso}
			\begin{split}
				\bigl\lvert \bigl(M_{[1,k],t}(\trless{S}^x)\bigr)_{aa}\bigr\rvert &\le \bigl\lvert \bigl(M_{[1,k],t}(S^x)\bigr)_{aa}\bigr\rvert + N^{-1}\sum_j \bigl\lvert \bigl(M_{[1,k],t}(S^j)\bigr)_{aa}\bigr\rvert\\
				& \lesssim (\Upsilon_0)_{ax} + \frac{1}{N}  \sum_{j}(\Upsilon_0)_{aj} \lesssim \frac{1}{W} \lesssim \frac{N}{W^2} = \frac{1}{N\eta_t} \frac{N^2\eta_t}{W^2}.
			\end{split}
		\end{equation}
		where we used \eqref{eq:trless_S}, \eqref{eq:sumS=1}, \eqref{eq:Ups_majorates}, and \eqref{eq:Ups_norm_bounds}. Therefore, both \eqref{eq:trlessM_iso} and \eqref{eq:trlessM_av} hold for the non-saturated  $M_{[1,2],t}$  by \eqref{eq:trless2_iso}. 
		
		On the other hand, if $M_{[1,2],t}$ is saturated, i.e., $z_1 = \overline{z_2}$, then we can use the  bound in \eqref{eq:supercrit_Theta_notime}  (together with~\eqref{eq:sumgenTheta})  to deduce that
		\begin{equation}
			\bigl\lvert \bigl(M_{[1,k],t}(\trless{S}^x)\bigr)_{aa}\bigr\rvert = \biggl\lvert(\Theta_t)_{xa} - \frac{\im m(z_{1,t})}{N \im z_{1,t}}\biggr\rvert \lesssim \frac{N}{W^2} = \frac{1}{N\eta_t} \frac{N^2\eta_t}{W^2},
		\end{equation}
		where we recall that $\sum_j (\Theta_{t})_{xj} = \im m_t/(N \im z_t)$ by  \eqref{eq:sumTheta}. This concludes the proof of \eqref{eq:trlessM_iso} and \eqref{eq:trlessM_av} for $k=2$.
		
		\textbf{Induction steps}. Fix integer $k \ge 3$. Assume that \eqref{eq:trlessM_iso} and \eqref{eq:trlessM_av} holds for all $k' \le k-1$. We proceed precisely as in the proof of Lemma \ref{lemma:M_bounds} above in Section \ref{sec:M_bounds}, defining
		\begin{equation}
			X_t^{k,n} \equiv X_{\bm z_t,t}^k(\trSet'; \bm x) := \Tr\bigl[M_{[1,k],t}(\trSet'; \bm x') \other{A}_k\bigr], \quad Y_t^{k,n} \equiv Y_{\bm z_t,t}^k(\trSet'; \bm x', a) := \bigl(M_{[1,k],t}(\trSet'; \bm x') \bigr)_{aa},
		\end{equation}
		where we recall that the set $\trSet'$ contains the indices of the traceless observables among $\other{A}_i$ for $i \in \indset{k-1}$, defined in \eqref{eq:trlessAs}, and $\other{A}_k$ is defined in \eqref{eq:trlessA_k}. 
		Recall that $n$ is the total number of traceless observables, including $\other{A}_k$ in the averaged case.
		
		Since at time $t=\crit$, $\eta_\crit \sim (W/N)^2$, Lemma \ref{lemma:M_bounds} implies that
		\begin{equation} \label{eq:trl_init}
			\norm{X_{\crit}^{k,n} }_*  \lesssim \frac{1}{(N\eta_{\crit})^{k-1}} \sim \biggl(\frac{N}{W^2}\biggr)^{k-1}, \quad  \norm{Y_\crit^{k,n}}_*   \lesssim \frac{1}{(N\eta_{\crit})^{k-1}} \sim \biggl(\frac{N}{W^2}\biggr)^{k-1},
		\end{equation}
		where the $\norm{\cdot}_*$-norm is defined in \eqref{eq:star_norm}. We will use the bounds \eqref{eq:trl_init} as the initial condition for the Duhamel's principle. 
		
		Similarly to \eqref{eq:av_M_evol}, \eqref{eq:iso_M_evol}, and \eqref{eq:traceless_k_iso_evol}, we obtain the evolution equations
		\begin{equation} \label{eq:trMevol}
			\begin{split}
				\frac{\mathrm{d}}{\mathrm{d}t}X_t^{k,n} &= \biggl(\frac{k}{2}I + \bigoplus_{j=1}^k \other{\mathcal{A}}_{j,t}\biggr)\bigl[X_t^{k,n} \bigr] +  F_{[1,k],t}^{n,\mathrm{av}},\\
				\frac{\mathrm{d}}{\mathrm{d}t}Y_t^{k,n} &= \biggl(\frac{k}{2}I + \bigoplus_{j=1}^{k-1} \other{\mathcal{A}}_{j,t}\biggr)\bigl[Y_t^{k,n} \bigr] + F_{[1,k],t}^{n,\mathrm{iso}}, 
			\end{split}
		\end{equation}
		for all $\crit \le t \le T$, where the forcing terms $F_{[1,k],t}^{\mathrm{av}}$ are defined analogously to their counterparts in \eqref{eq:av_M_forcing} and \eqref{eq:MFiso}, with all observables $A_i = S^{x_i}$ replaced by $\other{A}_i$, given by \eqref{eq:trlessAs} and \eqref{eq:trlessA_k}.  	
				
		\vspace{5pt}
		\textbf{Averaged $M$-terms.}
		We start with the averaged quantity $X_t^{k,n}$. Similarly to \eqref{eq:av_M_forcing_bound}, it is straightforward to check using the induction hypothesis that the forcing term $F_{[1,k],t}^{n,\mathrm{av}}$ admits the bound 
		\begin{equation} \label{eq:trMF_bound}
			\norm{F_{[1,k],t}^{n,\mathrm{av}} }_*  \prec \frac{1}{\eta_t} \frac{1}{(N\eta_t)^{k-1}}\biggl(\frac{N^2\eta_t}{W^2}\biggr)^{\lceil n /2 \rceil},\quad \crit \le t \le T. 
		\end{equation}
		
		Similarly to \eqref{eq:tr_prop}, the propagator $ \mathcal{P}^{k,n}_{s, t}$, generated by the linear term of the averaged evolution equation in \eqref{eq:trMevol}, admits the bound 
		\begin{equation} \label{eq:trMprop}
			\norm{\mathcal{P}^{k,n}_{s, t}}_{*\to *}  \lesssim  \biggl(\frac{\eta_s}{\eta_t}\biggr)^{k-n-p}~, \quad\crit \le s\le t\le T,
		\end{equation}
		where $p \equiv p(\trSet, \bm z_t)$ denotes the number of non-saturated propagators among $\{\widetilde{\mathcal{P}}^{[j]}_{s, t}\}_{j \notin \trSet}$, defined in \eqref{eq:trless_props}, with indices $j \notin \trSet$, corresponding to non-traceless observables. We record that  $p = 0$ occurs only when $z_j = \overline{z_{j+1}}$ for all $j \notin \trSet$ (under the cyclic convention $z_{k+1} = z_1$).

		Applying Duhamel's principle to \eqref{eq:trMevol},
		and taking the $\norm{\cdot}_*$, we conclude that, for any $\crit \le t \le T$,
		\begin{equation} \label{eq:trXduhamel}
			\begin{split}
				\norm{X_t^{k,n}}_* 
				&\prec \biggl(\frac{\eta_\crit}{\eta_t}\biggr)^{k-n-p} \norm{X_\crit^{k,n}}_* 
				+ \int_{\crit}^{t} \biggl(\frac{\eta_s}{\eta_t}\biggr)^{k-n-p}\frac{1}{\eta_s} \frac{1}{(N\eta_s)^{k-1}}\biggl(\frac{N^2\eta_s}{W^2}\biggr)^{\lceil n /2 \rceil} \mathrm{d}s\\
				&\prec \frac{1}{(N\eta_t)^{k-1}}  
				\biggl(\frac{N^2\eta_t}{W^2}\biggr)^{\lceil n /2 \rceil} \times  \biggl[\biggl(\frac{\eta_t}{\eta_\crit}\biggr)^{\lfloor n/2\rfloor + p - 1} 
				+ \int_{\crit}^{t}   \biggl(\frac{\eta_t}{\eta_s}\biggr)^{\lfloor n/2\rfloor + p - 1}  \frac{\mathrm{d}s}{\eta_s}\biggr].
			\end{split}
		\end{equation}
		where we used \eqref{eq:tcrit_def},  \eqref{eq:trl_init}, and \eqref{eq:trMF_bound}. It follows from \eqref{eq:int_rules} 
		that, if $n\ge 2$ or $p\ge 1$, then the term in the square brackets $\bigl[..\bigr]$ in \eqref{eq:trXduhamel} is bounded by a constant time $\log N$, implying the desired bound \eqref{eq:trlessM_av}.

		Hence, it remains to consider the special case $n = 1$ and $p=0$, i.e., when exactly one observable is traceless and all other propagators are saturated. In this case, we once again resort to observable regularization. Recall that $k \ge 3$, and that by cyclicity of $M$ in \eqref{eq:M_cyclic}, we can assume, without loss of generality, that $\trSet = \{k\}$. Then, $(k-1) \notin \trSet$, $z_{k-1} = \overline{z_k}$, and, by decomposing  $S^{x_{k-1}}$ into its traceless part and a multiple of identity, we obtain the following decomposition for $X_t^{k,n=1}$, 
		\begin{equation} \label{eq:XWard}
			\begin{split}
				X_t^{k,n=1} &= \Tr\bigl[M_{[1,k],t}(\emptyset; \bm x') \trless{S}^{x_k}\bigr] \\
				&= \Tr\bigl[M_{[1,k],t}(\{k-1\};\bm x') \trless{S}^{x_k}\bigr] + \frac{1}{N}\Tr\bigl[M_{[1,k],t}(S^{x_1}, S^{x_2}, \dots, I)\trless{S}^{x_k}\bigr]\\
				&= \Tr\bigl[M_{[1,k],t}(\{k-1\};\bm x') \trless{S}^{x_k}\bigr] + \frac{\Tr\bigl[M_{[1,\widehat{k-1},k],t} \trless{S}^{x_k}\bigr] - \Tr\bigl[M_{[1,k-1],t} \trless{S}^{x_k}\bigr]}{N(z_{k,t}-z_{k-1,t})},
			\end{split}
		\end{equation}
		where in the second step we used \eqref{eq:trless_S}, while in the last step we used \eqref{eq:M_Ward}. Here $M_{[1,k-1],t} \equiv M_{[1,k-1],t}(\emptyset;\bm x'')$, and $M_{[1,\widehat{k-1},k],t}$ is given by
		\begin{equation}
			M_{[1,\widehat{k-1},k],t} \equiv M_{[1,\widehat{k-1},k],t}(\bm x'') := M(z_{1,t}, S^{x_1}, \dots, z_{k-2}, S^{x_{k-2}}, z_{k}),
		\end{equation}
		i.e. hat indicates that the index is omitted. 
		In particular, the first averaged quantity in the last line of \eqref{eq:XWard} contains two traceless observables, and hence admits the bound \eqref{eq:trlessM_av} already obtained by \eqref{eq:trXduhamel}, while the $M$'s in the divided difference term both have length $k-1$ and can be estimated using the induction hypothesis. The denominator $N(z_{k,t}-z_{k-1,t})$ satisfies $|N(z_{k,t}-z_{k-1,t})| \gtrsim N\eta_t$, hence if $n=1$ and $p=0$, we have
		\begin{equation}
			\norm{X_t^{k,n=1}}_* \prec \frac{ (N^2\eta_t/W^2 )^{\lceil (n+1)/2\rceil}}{(N\eta_t)^{k-1}} + \frac{1}{N\eta_t} \frac{ (N^2\eta_t/W^2 )^{\lceil n /2 \rceil}}{(N\eta_t)^{k-2}} \lesssim \frac{ (N^2\eta_t/W^2 )^{\lceil n /2 \rceil}}{(N\eta_t)^{k-1}}.
		\end{equation}
		This concludes the proof of  \eqref{eq:trlessM_av}.

		\vspace{5pt}
		\textbf{Isotropic $M$-terms.} We complete the proof by estimating the isotropic quantity  $Y_t^{k,n}$. Similarly to \eqref{eq:iso_M_forcing_bound}, it is straightforward to check that the forcing term $F_{[1,k],t}^{n,\mathrm{iso}}$ admits the bound
		\begin{equation}
			\norm{F_{[1,k],t}^{n,\mathrm{iso}} }_*  \prec \frac{1}{\eta_t} \frac{1}{(N\eta_t)^{k-1}}\biggl(\frac{N^2\eta_t}{W^2}\biggr)^{\lceil n /2 \rceil},\quad \crit \le t \le T. 
		\end{equation}
		Applying Duhamel's principle, using \eqref{eq:tr_prop} for $k' = k-1$ to estimate the linear propagators, and proceeding as in \eqref{eq:trXduhamel}, we deduce that
		\begin{equation}
			\norm{Y_t^{k,n}}_*  
			\prec \frac{ (N^2\eta_t/W^2 )^{\lceil n/2\rceil}}{(N\eta_t)^{k-1}}  
			\times  \biggl[\biggl(\frac{\eta_t}{\eta_\crit}\biggr)^{\lfloor n/2\rfloor} 
			+ \int_{\crit}^{t}   \biggl(\frac{\eta_t}{\eta_s}\biggr)^{\lfloor n/2\rfloor}  \frac{\mathrm{d}s}{\eta_s}\biggr] \prec \frac{ (N^2\eta_t/W^2 )^{\lceil n/2\rceil}}{(N\eta_t)^{k-1}}.
		\end{equation}
		Note that no special case study is needed in the isotropic case, since the isotropic evolution equation only generates a total of $(k-1)$ propagators, $n$ of which are traceless, rendering the bound on the propagator in the isotropic case effectively equivalent to \eqref{eq:trMprop} with $p=1$. Therefore, \eqref{eq:trlessM_iso} holds. This concludes the proof of Lemma \ref{lemma:traceless_M_bounds}.
	\end{proof}

\section{Local laws for general observables and test vectors} \label{sec:general}
In this section, we outline the modifications needed to adapt the proof of the local laws to include general test vectors in the isotropic laws, as well as general diagonal observables (replacing the special observables $S^x$,
i.e., the columns of the variance matrix $S$) in both isotropic and averaged laws. 
We focus primarily on the crucial zig-step in Proposition \ref{prop:zig}. As we explain later in the section, the other main components of the proof---namely, the global laws in Proposition \ref{prop:global_laws} and the zag-step in Proposition~\ref{prop:zag}---require essentially no alteration.

To illustrate the main concepts, after some preliminary estimates in Section~\ref{sec:gen_prelim}, 
we first focus on case of the general ($\bm u, \bm v$) isotropic laws with special observables $S^x$ in Section~\ref{sec:geniso}.  The extension to both isotropic and averaged laws with arbitrary diagonal observables $A_i$ will be addressed in Section~\ref{sec:genobs} using the same conceptual framework.

The main proof presented above in the case of coordinate vectors and special observables $S^x$ relies on the relations between the size functions $\mathfrak{s}_{k,t}^\mathrm{av/iso}$, which are derived from the assumption on the control function $\Upsilon_\eta$ in Definition \ref{def:adm_ups_notime}. 
We group the necessary properties of $\Upsilon_\eta$ into two categories: the basic properties \eqref{eq:Ups_majorates_notime}--\eqref{eq:triag_notime} together with \eqref{eq:true_convol_notime} (as some control on the convolution of $\Upsilon$'s is needed), and the stronger additional properties \eqref{eq:convol_notime}--\eqref{eq:suppressed_convol_notime}, which are only used exclusively in estimating the action of the linear propagators on the size functions. 

As we show in the sequel, the generalized control function $(\Upsilon_\eta)_{\bm u \bm v}$, defined in \eqref{eq:genUp} satisfies appropriate analogs of the basic properties \eqref{eq:Ups_majorates_notime}--\eqref{eq:triag_notime} and \eqref{eq:true_convol_notime}. The majority of estimates in the proof rely only on these basic properties (and their immediate consequences for the time-dependent $\Upsilon_t$): see Propositions~\ref{prop:global_laws} and \ref{prop:zag}, Lemmas~\ref{lemma:mart_est} and \ref{lemma:forcing}, which explicitly do not assume \eqref{eq:convol_notime}--\eqref{eq:suppressed_convol_notime}. Hence, it is straightforward to verify that these parts of the proof carry over to the general isotropic setup, and we leave the trivial details to the reader.

However, the additional properties \eqref{eq:convol_notime}--\eqref{eq:suppressed_convol_notime} do not hold 
verbatim
when the coordinates $x,y$ are replaced by general vectors $\bm u, \bm v\in \mathbb{C}^N $.
 In fact, the versions of these properties that are valid for general test vectors $\bm u, \bm v$ incur an unavoidable loss compared to their specialized forms for coordinate vectors. In the sequel, we amend the parts of the proof that rely on the additional properties~ \eqref{eq:convol_notime}--\eqref{eq:suppressed_convol_notime}---namely, the propagator estimates---and explain how the resulting loss can be compensated using a mechanism already established in the main proof: observable regularization.

\subsection{Preliminaries}\label{sec:gen_prelim} 
We recall the definition of the generalized control functions\footnote{
To avoid any possible confusion between the different versions of  $\Upsilon$ used here and in the sequel, we adopt the following notational conventions: uppercase letters among the indices of $\Upsilon$ (e.g., $A$, $B$) denote diagonal observables; boldface lowercase letters (e.g., $\bm{u}$, $\bm{v}$) denote test vectors in $\mathbb{C}^N$; and regular lowercase letters (e.g., $a$, $b$, $i$, $j$, $x$) denote indices in $\indset{N}$.
} 
$(\Upsilon_\eta)_{\bm u \bm v}$ and $(\Upsilon_\eta)_{AB}$ from \eqref{eq:genUp}, restating them in their time-dependent form.
For any diagonal matrices $A, B$, any test vectors $\bm u, \bm v\in \mathbb{C}^N$, and all $t\in[0,T]$, we define
\begin{equation}\label{eq:Ups_uv}
	(\Upsilon_t)_{\bm u\bm v} \equiv  (\Upsilon_t)_{\bm v\bm u} := \sum_{ij} |u_i|^2 (\Upsilon_t)_{ij} |v_j|^2,
\end{equation}
\begin{equation}\label{eq:Ups_A}
	(\Upsilon_t)_{\bm u A} \equiv (\Upsilon_t)_{A\bm u} := \sum_{ij} |u_i|^2 (\Upsilon_t)_{ij} a_j,  
	\qquad (\Upsilon_t)_{ AB} \equiv (\Upsilon_t)_{BA} := \sum_{ij} a_i (\Upsilon_t)_{ij} b_j, 
\end{equation}
where we recall that the sequences $\{ a_i\} $ and  $\{ b_i\}$ are the minimizers on which the minimum in~\eqref{eq:tri} is achieved for $A$ and $B$, respectively.

The following lemma captures the basic properties \eqref{eq:Ups_majorates_notime}--\eqref{eq:triag_notime} and \eqref{eq:true_convol_notime} for $(\Upsilon_t)_{\bm u\bm v}$, $(\Upsilon_t)_{\bm u A}$, and $(\Upsilon_t)_{ AB}$.
\begin{lemma}[Basic Properties for Generalized Control Functions] \label{lemma:gen_basic}
	The generalized control functions $\Upsilon_s$ and $\Upsilon_t$ at times $0\le s \le t \le T$ satisfy the following analogs of the basic properties:
	\begin{itemize}
		\item[(i)] Majoration of $M_{[1,2],t}(z_1, A, z_2)$.
		\begin{alignat}{2} \label{eq:gen_majorates}
			\begin{split}
				\bigl\lvert \bigl\langle \bm u, M_{[1,2],t}\, \bm v \bigr\rangle \bigr\rvert &\lesssim \sqrt{(\Upsilon_{t})_{\bm u A}(\Upsilon_{t})_{A\bm v}}, \quad \text{if } (\im z_{1,t}) (\im z_{1,t}) < 0, \\ 
				\bigl\lvert \bigl\langle \bm u, M_{[1,2],t}\, \bm v \bigr\rangle \bigr\rvert&\lesssim \sqrt{(\Upsilon_{0})_{\bm u A}(\Upsilon_{0})_{A \bm v}}, \quad \text{if } (\im z_{1,t}) (\im z_{1,t}) > 0,
			\end{split}
		\end{alignat} 
		where we recall from \eqref{eq:M12} that $M_{[1,2],t} \equiv M_{[1,2],t}(z_1, A, z_2) =  (1-m_{1,t}m_{2,t}\mathcal{S})^{-1}\bigl[ m_{1,t}m_{2,t} A\bigr]$ with $m_{j,t} := m(z_{j,t})$.
		
		\item[(ii)] Bounds. For all $\bm u, \bm v \in \mathbb{C}^N$,
		\begin{equation}\label{eq:gen_norm_bounds}
			\frac{1}{N^{2D'}}\lesssim \frac{(\Upsilon_t)_{\bm u \bm v}}{\norm{\bm u}\norm{\bm v}} \lesssim \frac{1}{\ell_t\eta_t}, \qquad \sum_q (\Upsilon_t)_{\bm u q} \lesssim \frac{1}{\eta_t}\norm{\bm u}.
		\end{equation}
		Moreover, \eqref{eq:gen_norm_bounds} holds with one or both test vectors $\bm u$ and $\bm v$ replaced by diagonal observables $A$ and $B \in \mathbb{C}^{N\times N}$, with the corresponding $\norm{\cdot}$-norms of the vectors are replaced by $ |\! | \! | \cdot  |\! | \! | $-norm of the observables, defined in \eqref{eq:tri}.
		
		\item[(iii)] Monotonicity. For all $\bm u, \bm v \in \mathbb{C}^N$,
		\begin{equation} \label{eq:gen_time_monot} 
			(\Upsilon_{s})_{\bm u \bm v}  \lesssim  (\Upsilon_{t})_{\bm u \bm v}.
		\end{equation}
		
		\item[(iv)] Triangle and Convolution inequalities. For all $\bm u, \bm v \in \mathbb{C}^N$,
		\begin{alignat}{1}
			\max_q \bigl((\Upsilon_{s})_{\bm u q} (\Upsilon_{t})_{q\bm v}\bigr) &\lesssim (\ell_s\eta_s)^{-1} (\Upsilon_{t})_{\bm u \bm v}, \label{eq:gen_triag}\\
			\sum_q (\Upsilon_{s})_{\bm u q} (\Upsilon_{t})_{q\bm v}  &\lesssim \eta_s^{-1} (\Upsilon_{t})_{\bm u \bm v}. \label{eq:gen_convol}
		\end{alignat}
		Moreover, \eqref{eq:gen_time_monot} and \eqref{eq:gen_triag}--\eqref{eq:gen_convol} hold with one or both test vectors $\bm u$ and $\bm v$ replaced by diagonal observables $A$ and $B \in \mathbb{C}^{N\times N}$.
	\end{itemize}
	Here the mixed control functions $ (\Upsilon_{t})_{\bm u q}$ and $(\Upsilon_{t})_{A q}$  with $q\in\indset{N}$  are given by
	\begin{equation}
		(\Upsilon_t)_{\bm u q} \equiv (\Upsilon_t)_{q \bm u} := \sum_{ij} |u_i|^2 (\Upsilon_t)_{iq},  
		\qquad (\Upsilon_t)_{A q} \equiv (\Upsilon_t)_{qA} := \sum_{ij} a_i (\Upsilon_t)_{iq}, 
	\end{equation} 
	for all $\bm u \in \mathbb{C}^N$ and diagonal observables $A \in \mathbb{C}^{N \times N}$.
\end{lemma}
\begin{proof} [Proof of Lemma \ref{lemma:gen_basic}]
	Since $(\Upsilon_t)_{\bm u\bm v}$, $(\Upsilon_t)_{\bm u A}$, and $(\Upsilon_t)_{ AB}$, defined in \eqref{eq:Ups_uv}--\eqref{eq:Ups_A}, are positive linear combinations of $(\Upsilon_t)_{ij}$ with indices $i,j \in \indset{N}$, the bounds \eqref{eq:gen_majorates}--\eqref{eq:gen_convol} follow immediately from \eqref{eq:Ups_majorates_notime}--\eqref{eq:triag_notime} of Definition \ref{def:adm_ups_notime}, together with \eqref{eq:true_convol_notime} (and their dependent versions \eqref{eq:Ups_majorates}--\eqref{eq:true_convol}, \eqref{eq:Ups_norm_bounds}, \eqref{eq:Ups_time_monot}, \eqref{eq:ups_lower_bound}). This concludes the proof of Lemma \ref{lemma:gen_basic}.
\end{proof}

Next, we address the additional properties \eqref{eq:convol_notime}--\eqref{eq:suppressed_convol_notime} of $\Upsilon_\eta$, which do not directly extend to the $(\Upsilon_t)_{\bm u \bm v}$ for general test vectors (or observables). Recall that these properties are the necessary key inputs for establishing the propagator bounds in Lemmas~\ref{lemma:good_props} and~\ref{lemma:reg_props},  see~\eqref{eq:P_st_bound}.

First, we consider the square-root convolution \eqref{eq:convol_notime}. It is straightforward to verify that the bound \eqref{eq:convol_notime} fails for the quantities $(\Upsilon_t)_{\bm u q}$ and $(\Upsilon_t)_{A q}$,
it holds in the following weaker form with an additional  $\sqrt{N/\ell_s}$ factor
\begin{equation} \label{eq:gen_srqt_convol}
	\sum_q \sqrt{(\Upsilon_s)_{\bm uq} (\Upsilon_t)_{q\bm v}} \lesssim \frac{1}{\eta_s} \sqrt{N\eta_s \, (\Upsilon_t)_{\bm u \bm v}},
\end{equation}
using \eqref{eq:gen_convol} and a Schwarz inequality, analogously to \eqref{eq:sqrt_convol}.
In fact, this bound 
is, in general, sharp---for example, when $\bm u = \bm v = \bm 1$.  
Similarly, one or both of the test vectors $\bm u, \bm v$ in \eqref{eq:gen_srqt_convol} can be replaced by diagonal observables $A,B$ (with saturation occurring, e.g., for $A=B=I$). 

We emphasize that the  loss of a factor $\sqrt{N/\ell_s}$ in the propagator bounds is, in general, not affordable. 
 In fact, we will never use \eqref{eq:gen_srqt_convol} to estimate propagators.  
By inspecting the proofs of Lemmas~\ref{lemma:good_props} and~\ref{lemma:reg_props}, we observe that the additional property \eqref{eq:convol_notime} (in its time-dependent form \eqref{eq:convol}) is used exclusively to establish the estimate \eqref{eq:P_st_bound} for the action of the saturated propagator $\mathcal{P}_{s,t}$ on
 $\sqrt{\Upsilon_s \Upsilon_t}$, a product of $\sqrt{\Upsilon}$'s with different times. Indeed, the corresponding bounds  \eqref{eq:P_ss_bound} and \eqref{eq:P_tt_bound} for equal times, i.e. for $\sqrt{\Upsilon_s \Upsilon_s}$ and $\sqrt{\Upsilon_t \Upsilon_t}$, as well as the estimates for the non-saturated propagators in \eqref{eq:Q_ss_bound}, rely only on \eqref{eq:Schwarz_convol}---the time-dependent version of the full-power convolution estimate \eqref{eq:true_convol_notime} combined with the Schwarz inequality---which extends naturally to $(\Upsilon_t)_{\bm u \bm v}$ and $(\Upsilon_t)_{AB}$.

Using the representation \eqref{eq:P_decomp} for the propagator $\mathcal{P}_{s,t}$, along with \eqref{eq:gen_convol} and the Schwarz inequality, similarly to the derivation of \eqref{eq:Schwarz_convol}, we  easily obtain the following generalization of \eqref{eq:P_st_bound},
\begin{equation} \label{eq:gen_P_st}
	\sum_{q} \bigl\lvert(\mathcal{P}_{s,t})_{xq}\bigr\rvert \sqrt{(\Upsilon_t)_{\bm u q}(\Upsilon_s)_{q\bm v}} \lesssim \sqrt{\frac{\ell_t}{\ell_s}}\times\sqrt{\frac{\ell_s\eta_s}{\ell_t\eta_t}} \sqrt{(\Upsilon_t)_{\bm u x}(\Upsilon_t)_{x\bm v}},  \quad 0\le s\le t \le T,
\end{equation}
for all $x \in \indset{N}$ and $\bm u, \bm v \in \mathbb{C}^N$.
Note that \eqref{eq:gen_P_st} is still weaker than \eqref{eq:P_st_bound} by a large factor of $\sqrt{\ell_t/\ell_s}$,
but this is smaller  than the factor $\sqrt{N/\ell_s}$ we would incur if  \eqref{eq:gen_srqt_convol} were used.
As we demonstrate in the sequel, this loss factor  $\sqrt{\ell_t/\ell_s}$ can be offset using observable regularization.

For completeness, we also record the generalizations of \eqref{eq:P_ss_bound}, \eqref{eq:P_tt_bound}, and \eqref{eq:Q_ss_bound}, that, as discussed above, do not rely on \eqref{eq:convol_notime} and hence incur no additional loss,
\begin{equation} \label{eq:gen_PQ_ss}
	\begin{split}
		\sum_{q} \bigl\lvert(\mathcal{P}_{s,t})_{xq}\bigr\rvert \sqrt{(\Upsilon_r)_{\bm u q}(\Upsilon_r)_{q\bm v}} &\lesssim \frac{\eta_s}{\eta_r}\sqrt{(\Upsilon_t)_{\bm u x}(\Upsilon_t)_{x\bm v}}, \quad r\in \{s,t\}, \quad 0\le s \le t \le T,\\
		\sum_{q} \bigl\lvert(\mathcal{Q}_{s,t})_{xq}\bigr\rvert \sqrt{(\Upsilon_{r_1})_{\bm u q}(\Upsilon_{r_2})_{q\bm v}} &\lesssim \sqrt{(\Upsilon_{r_1})_{\bm u q}(\Upsilon_{r_2})_{q\bm v}}, \quad r_1, r_2\in \{s,t\}, \quad 0\le s \le t \le T,
	\end{split}
\end{equation}
for all $x \in \indset{N}$ and $\bm u, \bm v \in \mathbb{C}^N$. As before, one or both test vectors $\bm u, \bm v$ in \eqref{eq:gen_P_st}--\eqref{eq:gen_PQ_ss} can be replaced by diagonal observables $A$ and $B \in \mathbb{C}^{N\times N}$.

Finally, we address the additional property \eqref{eq:suppressed_convol_notime}, which is employed—in its time-dependent form \eqref{eq:suppressed_convol}—to control the action of the regularized operator $\reg{\Theta}_t^x$ in \eqref{eq:Thetaring_action_basic}. 
Notably, this estimate is only used in conjunction with a special regularized observable $\reg{S}^{x,y}$, defined in \eqref{eq:circ_above}, which requires that two consecutive observables in the chain be $S^x$ and $S^y$. 
In the sequel, we will only regularize those observables in the chain that are themselves special (of the form $S^x$) and are surrounded by two other special $S^x$ observables. In this setting, the bound \eqref{eq:Thetaring_action_basic} is directly applicable. 

We now combine the estimates \eqref{eq:gen_P_st}--\eqref{eq:gen_PQ_ss} and \eqref{eq:Thetaring_action_basic} to bound the action of a linear propagator on the size function for general test vectors and observables, defined as the product of the corresponding $\sqrt{\Upsilon}$ decay factors and an additional scalar prefactor, analogously to \eqref{eq:sfunc_def}.
In the special case $k=1$, the size functions $\mathfrak{s}_{1,t}^\mathrm{iso/av}$ are given by
\begin{equation} \label{eq:gen_sizef1}
	\mathfrak{s}_{1,t}^\mathrm{iso}(\bm u, \bm v) := \sqrt{(\Upsilon_t)_{\bm u \bm v}}, \qquad\mathfrak{s}_{1,t}^\mathrm{av}(A) := \frac{1}{\ell_t\eta_t} |\! | \! | A  |\! | \! | , 
\end{equation}
while for all $k \in \mathbb{N}$ with $k \ge 2$, we define
\begin{equation} \label{eq:gen_sizef}
	\begin{split}
		\mathfrak{s}_{k,t}^\mathrm{iso}(\bm u, A_1, \dots, A_{k-1}, \bm v) &:= \frac{1}{(\ell_t\eta_t)^{\frac{k-1}{2}}}\sqrt{(\Upsilon_t)_{\bm u A_1}\prod_{j=2}^{k-1}(\Upsilon_t)_{A_{j-1}A_{j}} \times (\Upsilon_t)_{A_{k-1}\bm v}} , \\  
		\mathfrak{s}_{k,t}^\mathrm{av}(A_1, \dots, A_k) &:= \frac{1}{(\ell_t\eta_t)^{\frac{k}{2}}}\sqrt{\prod_{j=2}^{k}(\Upsilon_t)_{A_{j-1}A_{j}}\times (\Upsilon_t)_{A_kA_1}}, 
	\end{split}
\end{equation} 
where $t\in [0,T]$,  $\bm u, \bm v \in \mathbb{C}^N$, and $A_i \in \mathbb{C}^{N\times N}$ for $i \in \indset{k}$ are diagonal observables. Note that the size function defied in \eqref{eq:gen_sizef1}--\eqref{eq:gen_sizef} agree with the right-hand sides of the corresponding local laws \eqref{eq:avelaw} and \eqref{eq:law1}. 

For any integer $k$, let $\mathfrak{S} \subset \indset{k}$ be the subset of indices corresponding to special observables $S^x$, so that
\begin{equation} \label{eq:special_gen_obs}
	A_j = S^{x_j}, \quad j \in \mathfrak{S},
\end{equation}
for some $x_j\in\indset{N}$, 
while for $j \notin \mathfrak{S}$, the observables $A_j$ may be general diagonal matrices. 
We recall that in this special case the corresponding $\Upsilon$ functions simplify, see~\eqref{eq:UpsSS}.  
For $\{A_i\}_{i=1}^k$ satisfying \eqref{eq:special_gen_obs}, we view $\mathfrak{s}_{k,t}^\mathrm{av}(A_1,\dots, A_k)$ and $\mathfrak{s}_{k+1,t}^\mathrm{iso}(\bm u, A_1,\dots, A_k, \bm v)$ as functions of the external indices  $\{x_j\}_{j\in\mathfrak{S}}$ corresponding to the special observables, and abbreviate
\begin{equation}
	\mathfrak{s}_{k,t}^\mathrm{av}(\mathfrak{S}; \bm x) := \mathfrak{s}_{k,t}^\mathrm{av}(A_1,\dots, A_k), \quad \mathfrak{s}_{k+1,t}^\mathrm{iso}(\mathfrak{S}; \bm u, \bm x, \bm v) := \mathfrak{s}_{k+1,t}^\mathrm{iso}(\bm u, A_1,\dots, A_k, \bm v),
\end{equation}
where, with slight abuse of notation, we define $\bm x := (x_j)_{j \in \mathfrak{S}}$.  
 The other observables $A_j$'s will not play any essential role hence they are
suppressed in the notation. 

We now introduce the class of propagators that arise from the evolution equations for chains with general observables. Consider a propagator $\other{\mathcal{P}}^{k}_{s,t}(\mathfrak{S})$ acting on the external indices  $\bm x = (x_j)_{j \in \mathfrak{S}}$ corresponding to the special observables among $\{A_i\}_{i=1}^k$,
\begin{equation} \label{eq:genA_props}
	\other{\mathcal{P}}^{k}_{s,t}(\mathfrak{S}) := \bigotimes_{j \in \mathfrak{S}}\other{\mathcal{P}}^{(j)}_{s,t}, \quad \text{with}\quad  \other{\mathcal{P}}_{s,t}^{(j)} \in \bigl\{ \mathcal{P}_{s,t}, \mathcal{Q}_{s,t}, \overline{\mathcal{Q}}_{s,t},  
	I \bigr\}, \quad 0 \le s \le t\le T.
\end{equation}
We control the action of $\other{\mathcal{P}}^{k}_{s,t}(\mathfrak{S})$ on the size functions $\mathfrak{s}^{\mathrm{av/iso}}(\mathfrak{S};,\cdot,)$ using the following lemma.
\begin{lemma} [General Propagator Lemma] \label{lemma:gen_prop}
	Fix an integer $k$ and a subset $\mathfrak{S} \subset \indset{k}$ of special observable.	
	Let $s,t \in [0,T]$ be a pair of times satisfying $s\le t$, and let $\other{\mathcal{P}}^{k}_{s,t} \equiv \other{\mathcal{P}}^{k}_{s,t}(\mathfrak{S})$ be a linear operator of the form \eqref{eq:genA_props}. Let $\varphi >0 $ be a positive control parameter. 
	
	First, in the isotropic case, if $k = 1$ or $|\mathfrak{S}| \le k-1$ or there exists at least one
	index $j \in \indset{k}$, such that $\other{\mathcal{P}}_{s,t}^{(j)}$ is not saturated, then for any function $f_{\bm u \bm v}(\bm x)$ of $\bm x := (x_j)_{j \in \mathfrak{S}}$ and all test vectors $\bm u, \bm v \in \mathbb{C}^N$,
	\begin{equation} \label{eq:gen_iso_prop}
		\forall \bm x~~ \bigl\lvert f_{\bm u \bm v}(\bm x) \bigr\rvert \lesssim \varphi \, \mathfrak{s}_{k+1,s}^\mathrm{iso}(\mathfrak{S}; \bm u,\bm x, \bm v) \quad 
		\Longrightarrow \quad 
		\forall \bm x~~ \bigl\lvert \other{\mathcal{P}}^{k}_{s,t}\bigl[f_{\bm u\bm v} \bigr]  (\bm x)\bigr\rvert \lesssim \varphi \,\sqrt{\frac{\ell_t\eta_t }{\ell_s\eta_s}}\, \mathfrak{s}_{k+1,t}^\mathrm{iso}(\mathfrak{S}; \bm u,\bm x, \bm v).
	\end{equation}
	On the other hand, if $k \ge 2$ with $|\mathfrak{S}| = k$ and all propagators $\other{\mathcal{P}}_{s,t}^{(j)}$ are saturated, then
	\begin{equation} \label{eq:gen_iso_prop_bad}
		\forall \bm x~~ \bigl\lvert f_{\bm u \bm v}(\bm x) \bigr\rvert \lesssim \varphi \, \mathfrak{s}_{k+1,s}^\mathrm{iso}(\mathfrak{S}; \bm u,\bm x, \bm v) \quad 
		\Longrightarrow \quad 
		\forall \bm x~~ \bigl\lvert \other{\mathcal{P}}^{k}_{s,t}\bigl[f_{\bm u\bm v} \bigr]  (\bm x)\bigr\rvert \lesssim \varphi \, \mathfrak{s}_{k+1,t}^\mathrm{iso}(\mathfrak{S}; \bm u,\bm x, \bm v).
	\end{equation}
	However, if $k \ge 3$ with $|\mathfrak{S}| = k$ and all propagators $\other{\mathcal{P}}_{s,t}^{(j)}$ are saturated, introducing one regularized operator $\reg{\Theta}^{(2)}_t$, defined in \eqref{eq:Thetaring}, we obtain
	\begin{equation} \label{eq:iso_one_ring}
		\forall \bm x~~ \bigl\lvert f_{\bm u \bm v}(\bm x) \bigr\rvert \lesssim \varphi \, \mathfrak{s}_{k+1,s}^\mathrm{iso}(\mathfrak{S}; \bm u,\bm x, \bm v) \quad  \Longrightarrow \quad 
		\forall \bm x~~ \biggl\lvert \other{\mathcal{P}}^{k}_{s,t}\circ \reg{\Theta}_t^{(2)} \bigl[f_{\bm u \bm v} \bigr]  (\bm x)\biggr\rvert \lesssim \varphi \, \frac{1}{\eta_t} \sqrt{\frac{\ell_t\eta_t}{\ell_s\eta_s}} \, \mathfrak{s}_{k+1,t}^\mathrm{iso}(\mathfrak{S}; \bm u,\bm x, \bm v).
	\end{equation}
	Here we chose to act by a regularized operator $\reg{\Theta}^{(2)}_t$, since $x_2$ is separated from $\bm u$ and $\bm v$ by $x_1$ and $x_3$, but any $\reg{\Theta}_t^{(j)}$ with $j \in \indset{2, k-1}$ could be chosen.
	
	For the averaged case, we always assume that $|\mathfrak{S}| \le k-1$, since the case $|\mathfrak{S}| = k$ is covered by the main proof.
	If $k \le 2$ or $|\mathfrak{S}| \le k-2$ or there exists at least one
	index $j \in \indset{k}$, such that $\other{\mathcal{P}}_{s,t}^{(j)}$ is not saturated, then for any function $g(\bm x)$ of $\bm x \in \indset{N}^{k}$,
	\begin{equation} \label{eq:gen_av_prop}
		\forall \bm x~~ \bigl\lvert g(\bm x) \bigr\rvert \lesssim \varphi \, \mathfrak{s}_{k,s}^\mathrm{av}(\mathfrak{S};\bm x) \quad \Longrightarrow \quad 
		\forall \bm x~~ \bigl\lvert \other{\mathcal{P}}^{k}_{s,t}\bigl[g \bigr]  (\bm x)\bigr\rvert \lesssim \varphi \, \frac{\ell_t\eta_t}{\ell_s\eta_s}  \, \mathfrak{s}_{k,t}^\mathrm{av}(\mathfrak{S}; \bm x).
	\end{equation}
	On the other hand, if $k\ge 3$ with $|\mathfrak{S}| = k-1$ and all propagators $\other{\mathcal{P}}_{s,t}^{(j)}$ are saturated, then
	\begin{equation} \label{eq:gen_av_prop_bad}
		\forall \bm x~~ \bigl\lvert g(\bm x) \bigr\rvert \lesssim \varphi \, \mathfrak{s}_{k,s}^\mathrm{av}(\mathfrak{S};\bm x) \quad \Longrightarrow \quad 
		\forall \bm x~~ \bigl\lvert \other{\mathcal{P}}^{k}_{s,t}\bigl[g \bigr]  (\bm x)\bigr\rvert \lesssim \varphi \, \sqrt{\frac{\ell_t\eta_t}{\ell_s\eta_s}}  \, \mathfrak{s}_{k,t}^\mathrm{av}(\mathfrak{S}; \bm x).
	\end{equation}
\end{lemma}
Under their assumptions, the bounds \eqref{eq:gen_iso_prop} and \eqref{eq:gen_av_prop} match their counterparts for special observables and coordinate vectors, \eqref{eq:iso_prop_bound} and \eqref{eq:av_prop_bound}, respectively, from Lemma \ref{lemma:good_props}. However, when $k \ge 2$, all observables are special ($|\mathfrak{S}| = k$) and all propagators $\other{\mathcal{P}}_{s,t}^{(j)}$ are saturated, the bound \eqref{eq:gen_iso_prop_bad} becomes weaker than \eqref{eq:iso_prop_bound} by a large factor of $\sqrt{\ell_s\eta_s/ (\ell_t\eta_t)}$. This loss stems from the weaker bound \eqref{eq:gen_P_st} for general vectors $\bm u, \bm v \in \mathbb{C}^N$ compared to \eqref{eq:P_st_bound} for coordinate vectors. Still, even the weaker bound \eqref{eq:gen_iso_prop_bad} suffices to establish
the isotropic  local laws for shorter chains of  length $k$ with  $\alpha_{k+1} = 0$. 
For longer chains, we need to recover 
the prefactor $\sqrt{\ell_t\eta_t/ (\ell_s\eta_s)}$; this is achieved by establishing the improved propagator bound 
in \eqref{eq:iso_one_ring}. The improvement comes
 from the additional (even smaller)  factor\footnote{Note that $\ell_s\ell_t\eta_t/W^2 \lesssim  \ell_t\eta_t/ (\ell_s\eta_s)$,
 so even a full power of  $\ell_t\eta_t/ (\ell_s\eta_s)$ could have been gained in~\eqref{eq:iso_one_ring}, but we
 will not need it} $\ell_s\ell_t\eta_t/W^2$,  
 contributed by the regularized   operator  $\reg{\Theta}$   (compare \eqref{eq:Thetaring_action_basic}  to \eqref{eq:Thetaring_action_basic1}),   thus allowing to prove 
isotropic local laws for any $\alpha_{k+1} \le 1/2$.  
Note that  the same  gain is present in the main proof for averaged chains, manifested as an improvement factor  $\ell_s/\ell_t$ 
in~\eqref{eq:one_ring} compared with~\eqref{eq:bad_av_prop}.

Similarly, for $k\ge 3$, with one general observable ($|\mathfrak{S}| = k-1$) and all saturated propagators, the averaged bound \eqref{eq:gen_av_prop_bad} is again weaker than \eqref{eq:av_prop_bound} by the same factor $\sqrt{\ell_s\eta_s/ (\ell_t\eta_t)}$, due to \eqref{eq:gen_P_st}. Nevertheless, it still yields the averaged local laws for all $\beta_k \le 1/2$, i.e., for $k \le \maxK - 1$. While the special case $k = \maxK$ could, in principle, be addressed by adding a regularized operator (as in \eqref{eq:iso_one_ring}), we instead pursue a simpler approach based on the entry-wise local laws 
that have already been established in the main proof.

In the main proof, observable regularization was applied only to averaged chains, not to isotropic ones. However, when general test vectors are introduced in the isotropic law, the $\sqrt{\ell_t/\ell_s}$ deficiency arising in \eqref{eq:gen_P_st} must be offset by observable regularization in order to maintain positive loss exponents  $\alpha_k\le 1/2$. By contrast, introducing general diagonal observables into the averaged laws has the opposite effect: regularization is no longer needed to establish the averaged law for general $A_i$. This is not a coincidence. Each general observable $A_i$ effectively removes a linear propagator from the associated evolution equation, and the gain from the reduced number of propagators outweighs the $\sqrt{\ell_t/\ell_s}$ cost introduced by each general observable.

\begin{proof}[Proof of Lemma \ref{lemma:gen_prop}]

We proceed in three steps.
 First, we consider a prototypical size function with $(k+1)$ $\sqrt{\Upsilon}$-factors, given by
\begin{equation}
	\mathfrak{s}^{\bm u, \bm v}_{k+1,s}(\bm x)  := \frac{1}{(\ell_s\eta_s)^{\frac{k}{2}}}\sqrt{ (\Upsilon_{s})_{\bm u x_1} (\Upsilon_{s})_{x_1x_2} \dots (\Upsilon_{s})_{x_k \bm v}}, \quad \bm x:= (x_1,\dots, x_k),
\end{equation}
for all $0 \le s \le T$, where we consider $\bm u,\bm v\in \mathbb{C}^N$ to be fixed, and one or both of them can be replaced by diagonal observables, i.e.  $\mathfrak{s}^{A, B}_{k+1,s}(\bm x).$  For $\mathfrak{s}^{\bm u, \bm v}_{k+1,s}(\bm x)$, we    establish the bound, for all $0\le s\le t\le T$,
	\begin{equation} \label{eq:gen_basic}
		\biggl\lvert \bigotimes_{j=1}^k \other{\mathcal{P}}^{{(j)}}_{s,t}\bigl[\mathfrak{s}^{\bm u, \bm v}_{k+1,s} \bigr](\bm x)  \biggr\rvert \lesssim \sqrt{\frac{\ell_t\eta_t}{\ell_s\eta_s}} \mathfrak{s}^{\bm u, \bm v}_{k+1,t}(\bm x) \times \begin{cases}
			1,\quad & \text{if}~k=1,\\
			\sqrt{\ell_t/\ell_s},\quad &\text{if}~k\ge 2~\text{and}~n=0,\\
			\sqrt{\ell_t\eta_t/(\ell_s\eta_s)},\quad &\text{if}~k\ge 2~\text{and}~n \ge 1,\\
		\end{cases}
	\end{equation}
	for all $k \ge 1$, $\bm x \in \indset{N}^k$, and $\bm u, \bm v \in \mathbb{C}^N$, where $n$ denotes the number of non-saturated propagators in $\{\other{\mathcal{P}}_{s,t}^{(j)}\}_{j=1}^k$.  
	
	The three factors in the case separation 
	indicate the additional factors over the base estimate from \eqref{eq:iso_prop_bound} in front of the $\times$-sign. 
	The large factor $\sqrt{\ell_t/\ell_s}$ is the one lost in~\eqref{eq:gen_P_st} due to the general vectors $\bm u, \bm v$,
	however this is incurred only when all propagators are saturated.
	Once we have at least one non-saturated propagator, the estimate improves not only 
	by the usual factor $\sqrt{\ell_t\eta_t/(\ell_s\eta_s)}$ (c.f. \eqref{eq:P_st_bound} vs.
	\eqref{eq:Q_ss_bound})  but the $\sqrt{\ell_t/\ell_s}$-loss due to
	general vectors also disappears. 
	In the second step
	we will show that the bound \eqref{eq:gen_basic} can be used as a prototypical ingredient  to prove \eqref{eq:gen_iso_prop}--\eqref{eq:gen_iso_prop_bad}, \eqref{eq:gen_av_prop}--\eqref{eq:gen_av_prop_bad}.
	 Finally, in the third step, we will use \eqref{eq:Thetaring_action_basic} and \eqref{eq:gen_convol} to finish the proof of~\eqref{eq:iso_one_ring}.
	
	We now present the technical proofs. 
	
	{\bf Step 1.}  Proof of  \eqref{eq:gen_basic}. 
	The case $k=1$ follows trivially from \eqref{eq:gen_PQ_ss} with $r = r_1 = r_2 = s$. 
	
	In the case $k \ge 2$ with $n=0$ (all propagators are saturated), 
	we apply $\other{\mathcal{P}}^{(j)}_{s,t}$ one by one, similarly to the proof of \eqref{eq:iso_prop_bound}, in the decreasing order of $j$. Using the first bound in \eqref{eq:gen_PQ_ss} with $r=s$ to bound the action of $\other{\mathcal{P}}^{(k)}_{s,t}$, we obtain
	\begin{equation}
		\bigl\lvert \other{\mathcal{P}}^{(k)}_{s,t}\bigl[\mathfrak{s}^{\bm u, \bm v}_{k+1,s}\bigr](\bm x)\bigr\rvert \lesssim \frac{1}{(\ell_s\eta_s)^{\frac{k}{2}}}\sqrt{ (\Upsilon_{s})_{\bm u x_1} \dots (\Upsilon_{s})_{x_{k-2}x_{k-1}} (\Upsilon_{t})_{x_{k-1}x_k}(\Upsilon_{t})_{x_k \bm v}},
	\end{equation}
	By \eqref{eq:P_st_bound}, the next $k-2$ propagators  $\{\other{\mathcal{P}}^{(j)}_{s,t}\}_{j=2}^{k-1}$ each contribute a factor $\sqrt{\ell_s\eta_s/(\ell_t\eta_t)}$, exactly as in \eqref{eq:prop_appl_bound}, yielding
	\begin{equation}
		\bigl\lvert \other{\mathcal{P}}^{(2)}_{s,t}\circ \dots \circ \other{\mathcal{P}}^{(k)}_{s,t}\bigl[\mathfrak{s}^{\bm u, \bm v}_{k+1,s}\bigr](\bm x)\bigr\rvert \lesssim \frac{1}{\ell_s\eta_s }\frac{1}{(\ell_t\eta_t)^{\frac{k-2}{2}}}\sqrt{ (\Upsilon_{s})_{\bm u x_1} (\Upsilon_{t})_{x_1x_2} \dots (\Upsilon_{t})_{x_k \bm v}}.
	\end{equation}
	Using \eqref{eq:gen_P_st} to bound the action of the final propagator $\other{\mathcal{P}}^{(1)}_{s,t}$, we obtain
	\begin{equation}\label{eq:laststep}
		\biggl\lvert \bigotimes_{j=1}^k \other{\mathcal{P}}^{{(j)}}_{s,t}\bigl[\mathfrak{s}^{\bm u, \bm v}_{k+1,s} \bigr](\bm x)  \biggr\rvert \lesssim \sqrt{\frac{\ell_t}{\ell_s}}\times \sqrt{\frac{\ell_t\eta_t}{\ell_s\eta_s} } \times \frac{1}{(\ell_t\eta_t)^{\frac{k}{2}}}\sqrt{ (\Upsilon_{t})_{\bm u x_1} (\Upsilon_{t})_{x_1x_2} \dots (\Upsilon_{t})_{x_k \bm v}},
	\end{equation}
	thus concluding \eqref{eq:gen_basic} for $k\ge 2$ and $n=0$.
	
	In the case $k \ge 2$ and $n\ge 1$, we follow the same procedure, except we start from both
	ends of the chain, $j=1$ and $j=k$, and apply propagators from both directions 
	until we reach a fixed  site $i\in\indset{k}$ corresponding to some non-saturated $\other{\mathcal{P}}^{(i)}_{s,t}$, 
	avoiding the loss from \eqref{eq:gen_P_st} completely. This way, each of the propagators corresponding to $j\notin \{1,i,k\}$ (of which there are at most $k-2$) contributes a factor of $\sqrt{\ell_t\eta_t/(\ell_s\eta_s)}$, yielding \eqref{eq:gen_basic} for $k\ge 2$ and $n\ge 1$. This concludes the proof of \eqref{eq:gen_basic}. Note that \eqref{eq:gen_basic} also holds with one or both $\bm u, \bm v$ replaced by diagonal observables.

{\bf Step 2.}
	Now we prove \eqref{eq:gen_iso_prop}--\eqref{eq:gen_iso_prop_bad} and \eqref{eq:gen_av_prop}--\eqref{eq:gen_av_prop_bad}.  Since  $\mathfrak{s}_{k+1,s}^\mathrm{iso}(\mathfrak{S}; \bm u, \bm x, \bm v) = \mathfrak{s}_{k+1,s}^{\bm u,\bm v}(\bm x)$ for $|\mathfrak{S}| =k$, the bound \eqref{eq:gen_iso_prop_bad} follows immediately from \eqref{eq:gen_basic}. 
	
	Similarly, the bound \eqref{eq:gen_iso_prop} with $|\mathfrak{S}| = k$ also follows directly from \eqref{eq:gen_basic}, since under this assumption $k=1$ or one of the propagators is not saturated.  Therefore, it remains to consider the case $|\mathfrak{S}| \le k-1$. In this case,  
	we first split the product of $\sqrt{\Upsilon}$-factors into sub-products delineated by generic observables $A_j$ with $j\notin\mathfrak{S}$. For example, if $\mathfrak{S} = \indset{k}\backslash\{j\}$, then we write
	\begin{equation}
	\mathfrak{s}_{k+1,s}^\mathrm{iso}(\mathfrak{S};\bm u, \bm x, \bm v) = \frac{1}{\sqrt{\ell_s\eta_s}} \mathfrak{s}_{j,s}^{\bm u, A_j}(x_1,\dots, x_{j-1}) \mathfrak{s}_{k-j+1,s}^{A_j, \bm v}(x_{j+1},\dots, x_{k}). 
	\end{equation}
	Since $|\mathfrak{S}|\le k-1$, the original product is split at least once.
	The action of the corresponding propagator on each sub-product is controlled by \eqref{eq:gen_basic} yielding \eqref{eq:gen_iso_prop} by simple power counting, together with the fact that $\ell_t^2\eta_t/W^2 \lesssim 1$.
	The bounds \eqref{eq:gen_av_prop}--\eqref{eq:gen_av_prop_bad} follow by the same reasoning.
	
	{\bf Step 3.} Finally, it remains to prove \eqref{eq:iso_one_ring}. Recall that $k \ge 3$. We proceed exactly as in the proof of \eqref{eq:gen_basic} above, first applying all propagators $\{\other{\mathcal{P}}^{(j)}_{s,t}\}_{j=3}^k$ in descending order, to obtain
	\begin{equation}
		\bigl\lvert \other{\mathcal{P}}^{(3)}_{s,t}\circ \dots \circ \other{\mathcal{P}}^{(k)}_{s,t}\bigl[\mathfrak{s}_{k+1,s}^\mathrm{iso}\bigr](\bm x)\bigr\rvert \lesssim \frac{1}{(\ell_s\eta_s)^{3/2} }\frac{1}{(\ell_t\eta_t)^{\frac{k-3}{2}}}\sqrt{ (\Upsilon_{s})_{\bm u x_1} (\Upsilon_{s})_{x_1x_2} (\Upsilon_{t})_{x_2x_3} \dots (\Upsilon_{t})_{x_k \bm v}}.
	\end{equation}
	Applying $\reg{\Theta}_t^{(2)}$,  and using \eqref{eq:Thetaring_action_basic} with $r=t$, we deduce that
	\begin{equation}
		\bigl\lvert \reg{\Theta}_t^{(2)}\circ \other{\mathcal{P}}^{(3)}_{s,t}\circ \dots \circ \other{\mathcal{P}}^{(k)}_{s,t}\bigl[\mathfrak{s}_{k+1,s}^\mathrm{iso}\bigr](\bm x)\bigr\rvert \lesssim \frac{1}{\eta_s}\frac{\ell_s\ell_t\eta_t}{W^2}\times\frac{1}{ \ell_s\eta_s}\frac{1}{(\ell_t\eta_t)^{\frac{k-2}{2}}}\sqrt{ (\Upsilon_{s})_{\bm u x_1} (\Upsilon_{t})_{x_1x_2}\dots (\Upsilon_{t})_{x_k \bm v}},
	\end{equation}
	Using the first  bound in \eqref{eq:gen_PQ_ss} with $r=t$ to estimate the action of $\other{\mathcal{P}}^{(2)}_{s,t}$, and the fact that $\reg{\Theta}_t^{(2)}$ commutes with $\{\other{\mathcal{P}}^{(j)}_{s,t}\}_{j=3}^k$, we obtain
	\begin{equation}
		\bigl\lvert \other{\mathcal{P}}^{(2)}_{s,t}\circ \dots \circ \other{\mathcal{P}}^{(k)}_{s,t}\circ \reg{\Theta}_t^{(2)}\bigl[\mathfrak{s}_{k+1,s}^\mathrm{iso}\bigr](\bm x)\bigr\rvert \lesssim \frac{1}{\eta_t}\frac{\ell_s\ell_t\eta_t}{W^2}\times\frac{1}{ \ell_s\eta_s}\frac{1}{(\ell_t\eta_t)^{\frac{k-2}{2}}}\sqrt{ (\Upsilon_{s})_{\bm u x_1} (\Upsilon_{t})_{x_1x_2}\dots (\Upsilon_{t})_{x_k \bm v}},
	\end{equation}	
	For the remaining propagator $\other{\mathcal{P}}^{(1)}_{s,t}$, we proceed precisely as in the proof of \eqref{eq:gen_basic} at step~\eqref{eq:laststep}, obtaining
	\begin{equation} \label{eq:iso_pre_ring}
		\bigl\lvert \other{\mathcal{P}}^{k}_{s,t}\circ \reg{\Theta}_t^{(2)}\bigl[\mathfrak{s}_{k+1,s}^\mathrm{iso} \bigr](\bm x)  \bigr\rvert \lesssim 
		\frac{\ell_s\ell_t\eta_t}{W^2}\times \sqrt{\frac{\ell_t}{\ell_s}}\times \sqrt{\frac{\ell_t\eta_t}{\ell_s\eta_s} } \times \frac{1}{\eta_t}\frac{1}{(\ell_t\eta_t)^{\frac{k}{2}}}\sqrt{ (\Upsilon_{t})_{\bm u x_1} (\Upsilon_{t})_{x_1x_2} \dots (\Upsilon_{t})_{x_k \bm v}},
	\end{equation}
	Note that the factor 
	\begin{equation}
		\frac{\ell_s\ell_t\eta_t}{W^2}\times \sqrt{\frac{\ell_t}{\ell_s}} = \frac{\ell_t^2\eta_t}{W^2}\times \sqrt{\frac{\ell_s}{\ell_t}} \lesssim 1,
	\end{equation}
	so \eqref{eq:iso_one_ring} follows from \eqref{eq:iso_pre_ring}. This concludes the proof of Lemma \ref{lemma:gen_prop}.
\end{proof}

\bigskip  

\subsection{General test vectors} \label{sec:geniso}
First, we explain how to introduce general test vectors into the isotropic laws with special observables. 
To condense the notation, we define
\begin{equation}
	\bigl(G_{[1,k],t}\bigr)_{\bm u \bm v} := \bigl\langle \bm u, G_{[1,k],t}\,\bm v\bigr\rangle,
\end{equation}
and use analogous conventions for $M$- and $(G-M)$-quantities.

Hence, our goal is to prove that, for all $k \in \indset{K}$ and $t\in [0,T]$,
\begin{equation} \label{eq:gen_iso_goal}
	\bigl\lvert \bigl((G-M)_{[1,k],t}(\bm x')\bigr)_{\bm u \bm v} \bigr\rvert \prec (\ell_t\eta_t)^{\alpha_k} \mathfrak{s}_{k,t}^\mathrm{iso}(\bm u, \bm x', \bm v),
\end{equation}
for all $\bm x' \in \indset{N}^{k-1}$, for any $\bm u, \bm v \in \mathbb{C}^{N}$. 

To control the left-hand side of \eqref{eq:gen_iso_goal}, we modify the definition of the isotropic control quantities $\Psi_{k,t}^\mathrm{iso}$ in \eqref{eq:Psi_def} as follows,
\begin{equation} \label{eq:gen_Psi_iso}
	\Psi_{k,t}^\mathrm{iso} := 
	\max_{\bm z_t \in \dom_t^k} \max_{\bm u ,\bm v \in \mathbb{V}(\bm u_0, \bm v_0)} \max_{\bm x'\in \indset{N}^{k-1}}  \frac{\bigl\lvert \bigl( G_{[1,k],t} - M_{[1,k],t} \bigr)_{\bm u  \bm v }(\bm x') \bigr\rvert}{(\ell_t\eta_t)^{\alpha_k}\mathfrak{s}_{k,t}^\mathrm{iso}(\bm u ,\bm x',\bm v )},
\end{equation}
where, for a fixed pair of vectors $\bm u_0, \bm v_0$, the set $\mathbb{V}(\bm u_0, \bm v_0)$ is given by 
\begin{equation} \label{eq:V_set}
	\mathbb{V}(\bm u_0, \bm v_0) := \{\bm u_0, \overline{\bm u_0}, \bm v_0, \overline{\bm v_0}\}\cup \{\bm e_i\}_{i=1}^N.
\end{equation}
In the sequel, we implicitly assume that  $\bm u, \bm v \in \mathbb{V}(\bm u_0, \bm v_0)$ unless stated otherwise.
This modification of $\Psi_{k,t}^\mathrm{iso}$ is necessary because taking the supremum over all vectors $\bm u, \bm v \in \mathbb{C}^N$ is not compatible with the notion of $\prec$-bounds, unlike the entry-wise definition \eqref{eq:Psi_def}, which involves only a maximum over indices $a, b\in \indset{N}$. However, to prove \eqref{eq:gen_iso_goal} for a fixed pair of vectors $\bm u_0, \bm v_0$, we only need to evaluate resolvents on vectors belonging to the set $\mathbb{V}(\bm u_0, \bm v_0)$, which contains $O(N)$ elements. Therefore, taking a maximum over $\bm u, \bm v \in \mathbb{V}(\bm u_0, \bm v_0)$ remains within the scope of $\prec$-bounds.

Since Propositions \ref{prop:global_laws} and \ref{prop:zag} use only \eqref{eq:true_convol} and never \eqref{eq:convol} directly, it is straightforward to check that they both hold with $\Psi_{k,t}^\mathrm{iso}$, defined as in \eqref{eq:gen_Psi_iso}. Therefore, it remains to verify Proposition \ref{prop:zig}. In particular, to show that the master inequalities \eqref{eq:masters} remain valid for the modified $\Psi_{k,t}^\mathrm{iso}$. Recall that to prove \eqref{eq:masters} in Section \ref{sec:masters_proof}, we studied the evolution equation of the averaged $\mathcal{X}^k_t$  and isotropic $\mathcal{Y}^k_t$  quantities, defined in \eqref{eq:G-M_kav} and \eqref{eq:G-M_kiso}, respectively. The isotropic quantity $\mathcal{Y}_{k,t}(\bm u , \bm x', \bm v )$, corresponding to the modified $\Psi_{k,t}^\mathrm{iso}$, defined in \eqref{eq:gen_Psi_iso}, is given by
\begin{equation}
	\mathcal{Y}^{k}_{t} \equiv \mathcal{Y}^{k}_{\bm z_t,t}(\bm u ,\bm x', \bm v ) := \bigl( (G - M)_{[1,k],t}(\bm x') \bigr)_{\bm u  \bm v }, \quad \bm u ,\bm v  \in \mathbb{V}(\bm u_0, \bm v_0).
\end{equation}
It is straightforward to check that the quantity $\mathcal{Y}^{k}_{\bm z_t,t}(\bm u ,\bm x', \bm v )$ satisfies the evolution equation \eqref{eq:k_iso_evol} with all instances of $a$ and $b$ replaced by $\bm u $ and $\bm v $, respectively\footnote{
	The term in the last line of \eqref{eq:iso_Fk_def} should be interpreted as $ \sum_q m_{1,t}m_{k,t} \sum_i \overline{u_i} v_i \bigl(I+\mathcal{A}_{k,t}\bigr)_{iq}
	\mathcal{X}_{t}^k(\bm x',q)$.
	}. 
	
	Furthermore, the bounds \eqref{eq:iso_mart_bound} and \eqref{eq:iso_forcing_bound} on the martingale and forcing terms in \eqref{eq:k_iso_evol} remain valid $a$ and $b$ replaced by $\bm u $ and $\bm v $, since Lemmas~\ref{lemma:mart_est} and~\ref{lemma:forcing} rely only on \eqref{eq:true_convol_notime} but not on \eqref{eq:convol_notime}--\eqref{eq:suppressed_convol_notime}. Therefore, using the propagator bound \eqref{eq:gen_iso_prop}--\eqref{eq:gen_iso_prop_bad} from Lemma~\ref{lemma:gen_prop} instead of \eqref{eq:iso_prop_bound},  we show that the master inequalities \eqref{eq:iso_masters} hold for all $k \in \indset{\maxK/2}$ (for which $\alpha_k = 0$), precisely as in Section \ref{sec:masters_proof}.
	While~\eqref{eq:gen_iso_prop_bad} is weaker than the corresponding \eqref{eq:iso_prop_bound}, the
	extra small factor $\sqrt{\ell_t\eta_t/(\ell_s\eta_s)}$ is not necessary for $\alpha_k=0$.
	
	However, to complete the proof of $k\in \indset{\maxK/2+1, \maxK}$, for which $0 < \alpha_k \le 1/2$, 
	we need to use the stronger bound \eqref{eq:iso_one_ring} in the saturated case,  instead of \eqref{eq:gen_iso_prop_bad},  and hence we need to regularize one of the linear terms in \eqref{eq:k_iso_evol}. This is completely analogous the way we handled long ($k \ge \maxK/2$) averaged chains in Section~\ref{sec:masters_proof}, see discussion around \eqref{eq:lin_forcing} and Proposition \ref{prop:lin_term}. Similarly to \eqref{eq:Gcirc}, we consider a mollified quantity $\reg{\mathcal{Y}}_t^k$, defined as
	\begin{equation} \label{eq:Y_ring}
		\reg{\mathcal{Y}}_t^k  \equiv \reg{\mathcal{Y}}^{k}_t(\bm u , \bm x', \bm v ) :=\Theta_t^{(2)}\biggl[ \biggl\langle \bm u , (G - M)_{[1,k],t}\bigl(S^{x_1}, \reg{S}^{x_{2}}, S^{x_3}, \dots, S^{x_{k-1}}\bigr) \, \bm v  \biggr\rangle \biggr],
	\end{equation}
	where $\reg{S}^{x_{2}}$ is a regular observable, given by \eqref{eq:Sring_def}.
	Then the linear term $\Theta_t^{(2)}[\mathcal{Y}_t^k]$ in \eqref{eq:k_iso_evol} admits the decomposition  
	(c.f.~\eqref{eq:short_lin_term_decomp}--\eqref{eq:remainder_short}):
	\begin{equation} \label{eq:iso_lin_term_decomp}
		\Theta_t^{(2)} \bigl[ \mathcal{Y}^{k}_t\bigr] = \reg{\mathcal{Y}}^{k}_t + 
		\mathcal{R}^\mathrm{iso}_{[1,k],t}, \quad k\in \indset{\maxK/2+1, \maxK},
	\end{equation}
	where the remainder term $\mathcal{R}^\mathrm{iso}_{[1,k],t}\equiv \mathcal{R}^\mathrm{iso}_{[1,k],t}(\bm u , \bm x', \bm v )$ is defined as
	\begin{equation} \label{eq:iso_remainder}
		\mathcal{R}^\mathrm{iso}_{[1,k],t}(\bm u , \bm x', \bm v ) := (\Theta_t)_{x_{1}x_{2}}\mathcal{T}_{2}\bigl[\mathcal{Y}^{k}_t\bigr](x_1, x_3,\dots, x_{k-1}).
	\end{equation} 
	On the right-hand side of \eqref{eq:iso_remainder}, $\mathcal{Y}^{k}_{t} \equiv \mathcal{Y}^{k}_{t}(\bm u , \bm x', \bm v )$ is viewed as a function of $\bm x'$, and we recall that the  operator $\mathcal{T}_{2}$, defined in \eqref{eq:partial_trace},  denotes the partial trace in the variable $x_2$.

	Analogously\footnote{In \eqref{eq:Z_elov_short}, the restriction $k\in \indset{\maxK/2-1}$ stems from the definition of the quantities $\mathcal{Z}_t^k$, since for larger $k$ two observable had to be regularized instead of one. However, the same evolution equation holds for mollified quantities with one regularization for all $k\ge 2$. } to \eqref{eq:Z_elov_short} in the proof of  Lemma \ref{lemma:circ_evol}, it is straightforward to check that the quantity $\reg{\mathcal{Y}}_t^k$ satisfies the evolution equation
	\begin{equation} \label{eq:iso_ring_evol}
		\mathrm{d}\reg{\mathcal{Y}}_t^k = \biggl(\frac{k+1}{2}+\Theta_t^{(2)} + \Theta_t^{\oplus (k-1)}\biggr)\bigl[\reg{\mathcal{Y}}_t^k\bigr]\mathrm{d}t +\reg{\Theta}_t^{(2)}\biggl[
		\mathrm{d}\mathcal{M}^\mathrm{iso}_{[1,k],t}  + \mathcal{F}^\mathrm{iso}_{[1,k],t} \mathrm{d}t + \other{\mathcal{F}}^\mathrm{iso}_{[1,k],t}\mathrm{d}t\biggr],
	\end{equation} 
	where the additional forcing term $\other{\mathcal{F}}^\mathrm{iso}_{[1,k],t} \equiv \other{\mathcal{F}}^\mathrm{iso}_{[1,k],t}(\bm u , \bm x', \bm v )$ is given by (c.f., \eqref{eq:extra_forcing1}) 
	\begin{equation} \label{eq:extra_iso_forcing}
		\other{\mathcal{F}}^\mathrm{iso}_{[1,k],t}(\bm u , \bm x', \bm v ) := (\Theta_t)_{x_{1}x_{2}}  \bigl(\mathcal{T}_{2}\bigl[\mathcal{Y}^{k}_{t} \bigr](x_2, \dots, x_{k-1})  +  \mathcal{T}_{2}\bigl[\mathcal{Y}^{k}_{t} \bigr](x_1, x_3, \dots, x_{k-1})\bigr).
	\end{equation}

	Similarly to \eqref{eq:ward_local_law} and the proof of Lemma \ref{lemma:extra_forcing}, the additional forcing term $\other{\mathcal{F}}^\mathrm{iso}_{[1,k],t}$ admits the bound
	\begin{equation} \label{eq:iso_exf_bound}
		\frac{\eta_s\bigl\lvert \other{\mathcal{F}}^\mathrm{iso}_{[1,k],s}(\bm u , \bm x', \bm v ) \bigr\rvert}{(\ell_s\eta_s)^{\alpha_k}\mathfrak{s}_{k,s}^\mathrm{iso}(\bm u , \bm x', \bm v )} \lesssim \frac{\pis{k-1}}{(\ell_s\eta_s)^{\alpha_k - \alpha_{k-1}}},  \quad k \in \indset{\maxK/2 + 1, \maxK}.
	\end{equation}
	
	Equipped with \eqref{eq:iso_exf_bound} and the regularized propagator bound \eqref{eq:iso_one_ring}, we analyze \eqref{eq:iso_ring_evol} using the same procedure as in the proof of Proposition \ref{prop:circ_bound} in Section~\ref{sec:reg_evols}, concluding that 
	\begin{equation} \label{eq:Ycirc_bound}
		\max_{\tinit \le s \le  t\wedge\tau} \frac{\eta_s\bigl\lvert \reg{\mathcal{Y}}^{k}_s(\bm u , \bm x', \bm v ) \bigr\rvert}{(\ell_s\eta_s)^{\alpha_k}\mathfrak{s}^\mathrm{iso}_{k,s}(\bm u , \bm x', \bm v )} \prec  \varphi_{k,t\wedge\tau}^\mathrm{iso}.
	\end{equation}
	Estimating $\mathcal{R}^\mathrm{iso}_{[1,k],t}$, defined in \eqref{eq:iso_remainder}, analogously to $\other{\mathcal{F}}^\mathrm{iso}_{[1,k],t}$   in \eqref{eq:extra_iso_forcing}, we combine \eqref{eq:iso_lin_term_decomp} and \eqref{eq:Ycirc_bound} to obtain the following bound on the linear term (c.f., \eqref{eq:lin_term_est}),  
	\begin{equation} \label{eq:iso_lin_term}
		\max_{\tinit \le s \le  t\wedge\tau} \frac{\eta_s\bigl\lvert \Theta_s^{(2)} \bigl[ \mathcal{Y}^{k}_s\bigr](\bm x') \bigr\rvert}{(\ell_s\eta_s)^{\alpha_k}\mathfrak{s}^\mathrm{iso}_{k,s}(\bm u , \bm x', \bm v )} \prec  \varphi_{k,t\wedge\tau}^\mathrm{iso}.
	\end{equation}
	Therefore, we can proceed as in the proof of the averaged master inequalities, treating $\Theta_t^{(2)} \bigl[ \mathcal{Y}^{k}_t\bigr]$ in \eqref{eq:k_iso_evol} as an additional forcing term, and conclude that \eqref{eq:iso_masters} holds for the modified $\Psi_{k,t}^\mathrm{iso}$, defined as in \eqref{eq:gen_Psi_iso}, for all $k\in\indset{K}$.
	This completes the proof of~\eqref{eq:gen_iso_goal}.

\subsection{General observables}\label{sec:genobs}
Finally, we explain how to introduce general diagonal observables into both the averaged and the isotropic local laws. 

We begin by preparing some notation. Let $\bm A := (A_1,\dots, A_k)$ denote a vector of diagonal observables $A_i$, and let $\bm A' :=(A_1,\dots, A_{k-1})$ be the vector with the last coordinate removed. Then, we denote
\begin{equation}
	G_{[1,k],t}(\bm A' ) := G_{[1,k],t}(A_1, \dots, A_{k-1}) = G_{1,t} A_1 G_{2,t}\dots A_{k-1} G_{k-1,t},
\end{equation}
and we use analogous conventions for the corresponding $M$-terms and $(G-M)$-quantities.

Our goal is to show that, for all $k \in \indset{\maxK}$ and $t\in [0,T]$,
\begin{equation} \label{eq:genA_goal}
	\begin{split}
		\bigl\lvert \bigl((G-M)_{[1,k],t}(\bm A')\bigr)_{\bm u \bm v} \bigr\rvert &\prec (\ell_t\eta_t)^{\alpha_k} \mathfrak{s}_{k,t}^\mathrm{iso}(\bm u, \bm A', \bm v),\\
		\bigl\lvert \Tr \bigl[ (G-M)_{[1,k],t}(\bm A') A_k\bigr] \bigr\rvert &\prec (\ell_t\eta_t)^{\beta_k} \mathfrak{s}_{k,t}^\mathrm{av}(\bm A),
	\end{split}
\end{equation}
for any $\bm u, \bm v \in \mathbb{C}^N$ and diagonal observables $\bm A$.
As we discussed in Section \ref{sec:geniso} above, Propositions \ref{prop:global_laws} and \ref{prop:zag} rely only \eqref{eq:true_convol} and never use \eqref{eq:convol} directly, hence is straightforward to check that they both also for resolvent chains with general observables $\bm A$. Therefore, it remains to verify Proposition~\ref{prop:zig}. 

In the sequel, we fix a pair of deterministic vectors $\bm u_0, \bm v_0 \in \mathbb{C}^N$ and a family of diagonal observables $\{B_i\}_{i=1}^\maxK$ for which we want to prove~\eqref{eq:genA_goal}. Along the way 
we will encounter chains with observables 
 $\bm{A}$ consisting of matrices drawn from $\mathbb{M}\equiv \mathbb{M}(\{B_i\}_{i=1}^\maxK)$, defined as 
\begin{equation} \label{eq:M_set}
	\mathbb{M}\bigl(\{B_i\}_{i=1}^\maxK\bigr) := \bigl\{ B_i, B_i^*, |B_i| \bigr\}_{i=1}^\maxK \cup \{S^{x}\}_{x=1}^N,
\end{equation}
similarly as we will need to consider vectors $\bm u, \bm v\in \mathbb{V}(\bm u_0, \bm v_0)$. 
Here we include $|B_i|$ in \eqref{eq:M_set} to allow for performing reductions by splitting the chain at the general observable $B_i$.  Indeed, similarly to \eqref{eq:reduction}, we have
\begin{equation} \label{eq:genA_red}
	\begin{split}
		\bigl(X^*  B_i^* Y^* S^a Y B_i X\bigr)_{bb} &= \Tr \bigl[|B_i|^{1/2} Y^* S^a Y |B_i|^{1/2}\, U_i|B_i|^{1/2} X \lvert b \rangle \langle b \rvert X^*|B_i|^{1/2} U_i^* \bigr]\\
		&\le \Tr \bigl[|B_i|^{1/2} Y^* S^a Y |B_i|^{1/2}\bigr]~\Tr \bigl[|B_i|^{1/2} X \lvert b \rangle \langle b \rvert X^*|B_i|^{1/2} U_i^*U_i\bigr]\\
		&= \Tr \bigl[Y^* S^a Y  |B_i| \bigr]~ \bigl(X^* |B_i| X\bigr)_{bb},
	\end{split}
\end{equation}
where we used the decomposition $B_i = |B_i|^{1/2}U_i|B_i|^{1/2}$ with a diagonal unitary $U_i = (U_i^*)^{-1}$.

For any vector of diagonal observables $\bm A = (A_1,\dots, A_k) \in \mathbb{M}^k$, we define the set $\mathfrak{S}(\bm A) \subset \indset{k}$ of indices corresponding to special observables by
\begin{equation}
	\mathfrak{S}(\bm A) := \{ i \in \indset{k} \, : \, \exists\, x\in\indset{N},~  A_i = S^x \}.
\end{equation}

 We now define the modified $\Psi$-quantities. To account for the number of special observables, we introduce an additional index $n$ into the control parameters $\Psi^\mathrm{av/iso}_{k,n,t}$. This allows us to track the number of special observables explicitly and to treat all linear terms containing more special observables than the original chain as additional forcing terms. For all $k \in \indset{\maxK}$, we define  
\begin{equation} \label{eq:genA_Psi}
	\begin{split}
		\Psi_{k,n,t}^\mathrm{iso} &:=   \max_{ \bm z_{t} \in \dom_t^k}   
		\max_{\bm u, \bm v\in\mathbb{V}(\bm u_0, \bm v_0)} \max_{\bm A' \in \mathbb{M}^{k-1} \,: \, |\mathfrak{S}(\bm A')| = n}  \frac{\bigl\lvert \bigl(G_{[1,k],t} - M_{[1,k],t}\bigr)_{\bm u \bm v} (\bm A') \bigr\rvert}{(\ell_t\eta_t)^{\alpha_k}\mathfrak{s}_{k,t}^\mathrm{iso}(\bm u,\bm A',\bm v)}, \quad n \in \indset{0, k-1},\\
		\Psi_{k,n,t}^\mathrm{av} 
		&:= 
		\max_{  \bm z_{t} \in \dom_t^k} 
		\max_{\bm A\in \mathbb{M}^k\,:\, |\mathfrak{S}(\bm A)| = n} \frac{\bigl\lvert \Tr \bigl[ \bigl(G_{[1,k],t} - M_{[1,k],t}\bigr) (\bm A')  A_k\bigr]\bigr\rvert}{(\ell_t\eta_t)^{\beta_k} \mathfrak{s}_{k,t}^\mathrm{av}(\bm A)}, \quad n \in \indset{0, k},
	\end{split}
\end{equation}
where, for any $\bm A \in \mathbb{M}^k$,  $\bm A' := (A_1, \dots, A_{k-1}) \in \mathbb{M}^{k-1}$. That is, the quantities $\Psi_{k,n,t}^\mathrm{av/iso}$ control the averaged and isotropic resolvent chains, respectively, with $n$ special observables. 

To account for the new  index $n$ in the $\Psi_{k,n,t}^\mathrm{av/iso}$, we introduce the corresponding deterministic control parameter $\psi_{k,n}^\mathrm{av/iso}$ and modify the definition of the stopping time $\tau$ in \eqref{eq:tau_def} accordingly, 
\begin{equation} \label{eq:genA_tau}
	\begin{split}
		\tau :=&~  \inf\biggl\{ t\in[\tinit,T] \, :\, \max_{k\in\indset{\maxK}}\max_{n\in\indset{0,k}}\Psi_{k,n,t}^\mathrm{av} / \pav{k,n} \ge 1 \biggr\}\\
		&\wedge \inf\biggl\{ t\in[\tinit,T] \, :\, \max_{k\in\indset{2,\maxK}}\max_{n\in\indset{0,k-1}} \Psi_{k,n,t}^\mathrm{iso} / \pis{k,n} \ge 1\biggr\}.
	\end{split}
\end{equation}
We also redefine the control parameters $\pav{k}$ and $\pis{k}$, independent of $n$, as
\begin{equation}\label{eq:newpsi}
	\pav{k} := \max_{n\in\indset{0,k}}\pis{k,n}, \qquad \pis{k} := \max_{n\in\indset{0,k-1}}\pis{k,n}, \quad k\in\indset{\maxK}.
\end{equation}

Similarly to \eqref{eq:G-M_kav}, \eqref{eq:G-M_kiso}, we introduce the average $\mathcal{X}_t^k(\bm A)$, and isotropic $\mathcal{Y}_t^k(\bm u, \bm A', \bm v)$ quantities, defined as
\begin{equation}
	\mathcal{X}_t^k(\bm A) := \Tr \bigl[ (G-M)_{[1,k],t}(\bm A') A_k\bigr],\qquad \mathcal{Y}_t^k(\bm u, \bm A', \bm v) := \bigl((G-M)_{[1,k],t}(\bm A')\bigr)_{\bm u \bm v}.
\end{equation}
Note that the case $|\mathfrak{S}(\bm A)|=k$ for the averaged and $|\mathfrak{S}(\bm A')|=k-1$ for the isotropic laws have already been
covered in the main proof and Section \ref{sec:geniso} above, and we henceforth ignore them for the rest of this argument. That is, we always assume that at least one of the observables involved is a generic diagonal matrix.

When we consider the evolution equations for $\mathcal{X}_t^k(\bm A)$ and $\mathcal{Y}_t^k(\bm u, \bm A', \bm v)$, we find that they are identical to \eqref{eq:k_av_evol} and \eqref{eq:k_iso_evol}, respectively, with the corresponding observables $S^{x_i}$ replaced by $A_i$, with the only difference being the structure of the linear terms. We explain this difference using the evolution equation for the averaged quantity $\mathcal{X}_t^k(\bm A)$ as an example, with the corresponding changes for the isotropic quantity $\mathcal{Y}_t^k(\bm u, \bm A', \bm v)$ being essentially analogous. The linear term in the evolution equation for $\mathcal{X}_t^k(\bm A)$ is replaced by (ignoring the harmless constant $k/2$)
\begin{equation} \label{eq:genA_new_lin}
	\sum_{j \in \mathfrak{S}(\bm A)} \mathcal{A}_{j,t}\bigl[\mathcal{X}_t^k(\bm A)\bigr] + \sum_{j \notin \mathfrak{S}(\bm A)} \sum_q \bigl(M_{[j,j+1],t}(A_j)\bigr)_{qq} \mathcal{X}_t^k\bigl(\bm A\bigr\rvert_{A_j\mapsto S^q}\bigr),
\end{equation}
where the vector of observables $\bm A\bigr\rvert_{A_j\mapsto S^q} := (A_1, \dots, A_{j-1}, S^q, A_{j+1}, \dots, A_k)$ denotes the vector $\bm A$ with the $j$-th observable $A_j$ replaced by $S^q$.  

The terms in the first sum in \eqref{eq:genA_new_lin} are linear terms in the sense that they are given by an action of a linear operator on the original $\mathcal{X}_t^k(\bm A)$ viewed as a function of $\{x_i\}_{i \in \mathfrak{S}(\bm A)}$.  
These linear terms produce the propagators that act on the special observables, as in~\eqref{eq:genA_props}. After applying Duhamel's formula, we estimate the action of these propagators using \eqref{eq:gen_av_prop}--\eqref{eq:gen_av_prop_bad}. 
Note that the bound \eqref{eq:gen_av_prop_bad} (for $n=k-1$) is only compatible with $\beta_k \le 1/2$, hence it can not be used for saturated chains of the maximal length $k=\maxK$ with $n=\maxK-1$. We consider this special case separately at the end of the argument.

The terms in the second sum in \eqref{eq:genA_new_lin} are not linear in the sense described above. Indeed, the vector $\bm A\bigr\rvert_{A_j\mapsto S^q}$ contains one more special observable than the original $\bm A$. Hence, we can estimate the terms in the second sum in \eqref{eq:genA_new_lin} directly using the modified stopping time $\tau$, defined in \eqref{eq:genA_tau} and the bound \eqref{eq:gen_majorates} for $M_{[j,j+1],t}(A_j)$. More precisely, we obtain, for all $\tinit\le t\le \tau$
\begin{equation} \label{eq:genA_non_lin}
	\begin{split}
		\biggl\lvert \sum_q \bigl(M_{[j,j+1],t}(A_j)\bigr)_{qq} \mathcal{X}_t^k\bigl(\bm A\bigr\rvert_{A_j\mapsto S^q}\bigr) \biggr\rvert 
		&\lesssim \pav{k,n+1}(\ell_t\eta_t)^{\beta_k}\sum_q (\Upsilon_{t})_{q A_j} \mathfrak{s}_{k,t}^\mathrm{av}\bigl(\bm A\bigr\rvert_{A_j\mapsto S^q}\bigr)\\
		&\lesssim \frac{1}{\eta_t}\pav{k,n+1}(\ell_t\eta_t)^{\beta_k} \mathfrak{s}_{k,t}^\mathrm{av}(\bm A  ),
	\end{split}
\end{equation}
where in the second step we used \eqref{eq:gen_convol}. These terms contribute an additional term $\pav{k,n+1}$ to the right-hand side of the master inequality \eqref{eq:av_masters}. However, if the deterministic control parameters $\psi^\mathrm{av/iso}_{k,n}$ are chosen  (c.f., \eqref{eq:psi_choice}) as 
\begin{equation} \label{eq:new_psi_choice}
	\pis{k,n}:= N^{\xi+k\nu + (k-1-n)\nu'},~n\in\indset{0,k-1},\quad \pav{k,n}:=N^{\xi+(k+1)\nu  + (k-n)\nu'},~n\in\indset{0,k}, 
\end{equation}
for all $k\in \indset{\maxK}$, with $\nu' \in (0, \nu/(10K))$, then it is straightforward to check that the additional terms $\pav{k,n+1}$ do not impede the self-improving structure of the master inequalities. The self-improvement in terms of $k\nu$ is analogous to \eqref{eq:masters}, while the self-improvement in $(k-n)\nu'$ is guaranteed since the only new term is $\pav{k,n+1}$, coming from \eqref{eq:genA_non_lin}.

All other terms in the evolution equations are estimated the same way as in the main proof for special observables, using the new definition of the control parameters $\psi^\mathrm{av/iso}_k$ in \eqref{eq:newpsi}. Recall that Lemmas \ref{lemma:mart_est} and \ref{lemma:forcing} explicitly did not use \eqref{eq:convol_notime}, and hence their proof holds verbatim for generalized $\Upsilon$'s.   Note that the choice of $\nu' \in (0, \nu/(10K))$ guarantees that $\psi^{\mathrm{av/iso}}_{k} \ll N^{\etaexp/(10K^2)}\psi^{\mathrm{av/iso}}_{k,n}$ for all admissible $n$, which is compensated by the small factor $N^{-\etaexp/(2K^2)}$, precisely as in \eqref{eq:it}.

It remains to consider the special case of the averaged law for saturated chains of maximal length $k = \maxK$ with $n = k - 1$ special observables—that is, only one observable $A_i$ is generic. Recall that for $k=\maxK$ the loss exponent $\beta_\maxK > 1/2$, which is not compatible with the propagator bound \eqref{eq:gen_av_prop_bad}, that only allows for $\beta_k \le 1/2$.

This apparent discrepancy is easily resolved by using the optimal entry-wise local law for $k=\maxK$,
\begin{equation}
	\bigl\lvert \bigl((G-M)_{[1,\maxK],t}\bigr)_{jj} (\bm x')\bigr\rvert \prec \mathfrak{s}_{\maxK,t}^\mathrm{iso}(j,\bm x', j).
\end{equation} 
Indeed, such a local law  follow from the main proof if the maximal chain length $\maxK$ is chosen to be twice as large as $\maxK$ in the current argument.
Without loss of generality we can assume that the only generic observable is $A_\maxK$ (if not, we cyclically shift  the chain). Then, using \eqref{eq:tri} and \eqref{eq:Ups_A} to bound $|(A_\maxK)_{jj}| \lesssim (\Upsilon_0)_{jA_\maxK}$, and \eqref{eq:gen_convol}, we obtain\footnote{In the main proof this was the most critical case and it required  two regularizations, but now we can rely on previously proven isotropic laws with special observables only.}
\begin{equation}
	\begin{split}
		\bigl\lvert \Tr\bigl[(G-M)_{[1,\maxK],t}(\bm x') A_{\maxK} \bigr] \bigr\rvert &\le \sum_j \bigl\lvert \bigl((G-M)_{[1,\maxK],t}\bigr)_{jj} (\bm x') (A_{\maxK})_{jj} \bigr\rvert \\
		&\prec \sum_j  \mathfrak{s}_{\maxK,t}^\mathrm{iso}(j,\bm x', j) (\Upsilon_0)_{jA_\maxK}\\
		&\prec \sqrt{\ell_t\eta_t} \,\mathfrak{s}_{\maxK,t}^\mathrm{av}( S^{x_1},\dots, S^{x_{\maxK-1}}, A_\maxK ).
	\end{split}
\end{equation}
Hence, the averaged law for saturated chains with $k = \maxK$ and $n = k-1$ special observables (one generic observable) holds with a loss factor of $\sqrt{\ell_t \eta_t} \ll (\ell_t\eta_t)^{\beta_\maxK}$, 
without relying on the modified stopping time \eqref{eq:genA_tau}. This renders the restriction $\beta_k \le 1/2$,
needed to  use  \eqref{eq:gen_av_prop_bad},  irrelevant in this case. 
This concludes the proof of \eqref{eq:genA_goal}.

\subsection{$M$-bounds for general observables}
In this short section we explain the proof of \eqref{eq:isoM}--\eqref{eq:aveM}.

First, we observe that the isotropic bound \eqref{eq:isoM} for general test vectors $\bm u, \bm v \in \mathbb{C}^N$ follows immediately from its (time-dependent) entry-wise version (c.f., \eqref{eq:M_bound}),
\begin{equation} \label{eq:entM}
	\bigl\lvert \bigl(M_{[1,k],t}(\bm A')\bigr)_{ab} \bigr\rvert \prec \delta_{ab} \sqrt{\ell_t\eta_t}\,\mathfrak{s}_{k,t}^\mathrm{iso}(a, \bm A', a),  \quad k\ge 2.
\end{equation}
Here and in the sequel, we use $\prec$-notation to absorb the irrelevant poly-log factors. 
Indeed, assuming \eqref{eq:entM} holds, we obtain, for all $k \ge 2$,
\begin{equation}
	\begin{split}
		\bigl\lvert \bigl(M_{[1,k],t}(\bm A')\bigr)_{\bm u \bm v}  \bigr\rvert &= \biggl\lvert \sum_{ab} \overline{u_a}  \bigl(M_{[1,k],t}(\bm A')\bigr)_{ab}  v_b \biggr\rvert 
		\prec \sqrt{\ell_t\eta_t}\sum_a|u_a|\,\mathfrak{s}_{k,t}^\mathrm{iso}(a, \bm A', a) \,|v_a| \\
		&\prec (\ell_t\eta)^{1-\frac{k}{2} } \biggl(\prod_{i=1}^{k-2} (\Upsilon_t)_{A_iA_{i+1}}   \sum_a (\Upsilon_t)_{A_1a}|u_a|^2   \sum_b (\Upsilon_t)_{A_{k-1}b} |v_b|^2\biggr)^{1/2}\\
		&\prec \sqrt{\ell_t\eta_t}\,\mathfrak{s}_{k,t}^\mathrm{iso}(\bm u, \bm A', \bm v),
	\end{split}
\end{equation}
where in the penultimate we used Schwarz inequality, and the last step follows by \eqref{eq:Ups_uv} and \eqref{eq:gen_sizef}. 

Similarly, the averaged bound \eqref{eq:aveM} for $\bm A = (A_1,\dots, A_k)$ with $|\mathfrak{S}(\bm A)| = n$ special observables follows immediately from its entry-wise counterpart \eqref{eq:entM}. Indeed, if $\mathfrak{S} \neq \indset{k}$, by cyclicity of $M$ in \eqref{eq:M_cyclic}, we can assume without loss of generality that $k \notin \mathfrak{S}$, hence
\begin{equation} \label{eq:isoM_to_aveM}
	\begin{split}
		\bigl\lvert \Tr\bigl[M_{[1,k],t}(\bm A') A_k\bigr]  \bigr\rvert &= \biggl\lvert \sum_{j} \bigl(M_{[1,k],t}(\bm A')\bigr)_{jj}\,(A_k)_{jj}   \biggr\rvert 
		\prec \sqrt{\ell_t\eta_t}\sum_a\mathfrak{s}_{k,t}^\mathrm{iso}(j, \bm A', j) \bigl\lvert (A_k)_{jj} \bigr\rvert\\
		&\prec \frac{\ell_t\eta }{(\ell_t\eta)^{\frac{k}{2} }} \biggl(\prod_{i=1}^{k-2} (\Upsilon_t)_{A_iA_{i+1}}   \sum_j (\Upsilon_t)_{A_1j} (\Upsilon_0)_{jA_k}    \sum_i (\Upsilon_t)_{A_{k-1}i} (\Upsilon_0)_{iA_k}\biggr)^{1/2}\\
		&\prec \ell_t\eta_t \,\mathfrak{s}_{k,t}^\mathrm{av}( \bm A ),
	\end{split}
\end{equation}
where, in the second step, we used \eqref{eq:tri} and \eqref{eq:Ups_A} to bound $|(A_k)_{jj}| \lesssim (\Upsilon_0)_{jA_k}$, followed by Schwarz in the penultimate step, and \eqref{eq:gen_convol} in the last step.
Note that in \eqref{eq:isoM_to_aveM} the overall number of special observables in $\bm A'$ is one less than that in $\bm A$. 

Therefore, it remains to prove \eqref{eq:entM} for $\bm A'$ with $n \le k-2$ special observables (the case $n=k-1$ has already been covered in \eqref{eq:M_bound}) and $k \ge 3$ (the case $k=2$ follows from \eqref{eq:gen_majorates}).

We proceed by an overall induction in chain length $k \ge 3$, similarly to the proof of Lemma \ref{lemma:M_bounds} in Section \ref{sec:M_bounds}, and a nested induction in $k-n \in \indset{2,k}$. Note that these assumptions on $k$ and $n$ allow us to use the propagator bound \eqref{eq:gen_iso_prop} with the full $\sqrt{\ell_t\eta_t/(\ell_s\eta_s)}$ factors, compatible with the target bound \eqref{eq:entM}. In particular, we never have to perform observable regularizations. 

It follows from \eqref{eq:dM} that the quantity $Y_t^k := \bigl(M_{[1,k,t]}(\bm A')\bigr)_{aa}$ satisfies an evolution equation completely analogous to \eqref{eq:iso_M_evol} with each $S^{x_i}$ replaced by the corresponding $A_i$, with the only difference being the structure of the linear term, similarly to \eqref{eq:genA_new_lin},
\begin{equation} \label{eq:genAY_evol}
	\frac{\mathrm{d}}{\mathrm{d}t}Y_t^k = \biggl(\frac{k}{2}I + \bigoplus_{j\in \mathfrak{S}(\bm A')} \mathcal{A}_{j,t}\biggr)\bigl[Y_t^k \bigr] +\sum_{j \notin \mathfrak{S}(\bm A')} \sum_q \bigl(M_{[j,j+1],t}(A_j)\bigr)_{qq}  \bigl(M_{[1,i],[j,k],t}^{(q)}\bigr)_{aa} + F_{[1,k],t}^\mathrm{iso}.
\end{equation}
All $M$-terms in the $j \notin \mathfrak{S}(\bm A')$ summation can be estimated using the $n$-induction hypothesis since they contain strictly more special observables.
Note that all $M$-terms involved in $F_{[1,k],t}^\mathrm{iso}$ have strictly shorted length than $k$ with the exception of the last term in the definition of $F_{[1,k],t}^\mathrm{iso}$ in \eqref{eq:MFiso}, namely,
\begin{equation}
	 \sum_q m_{1,t}m_{k,t} \bigl(I+\mathcal{A}_{k,t}\bigr)_{aq} \Tr\bigl[M_{[1,k],t}(\bm A') S^q\bigr].
\end{equation}
However, using \eqref{eq:isoM_to_aveM}, we derive an bound on  $\Tr\bigl[M_{[1,k],t}(\bm A') S^q\bigr]$ from \eqref{eq:entM} with the same $k$ and $|\mathfrak{S}(\bm A')| = n+1$. More specifically, we conclude that 
\begin{equation}
	\sum_{j \notin \mathfrak{S}(\bm A')} \sum_q \bigl\lvert\bigl(M_{[j,j+1],t}(A_j)\bigr)_{qq}  \bigl(M_{[1,i],[j,k],t}^{(q)}\bigr)_{aa}\bigr\rvert + \bigl\lvert F_{[1,k],t}^\mathrm{iso}\bigr\rvert \prec \frac{1}{\eta_t} 
	\sqrt{\ell_t\eta_t} \mathfrak{s}_{k,t}^\mathrm{iso}(a, \bm A', a).
\end{equation}
Hence, applying Duhamel's principle to \eqref{eq:genAY_evol} and using \eqref{eq:gen_iso_prop}, together with $n \le k-2$, to bound the action of the propagators on the special observables, we obtain \eqref{eq:entM}. This concludes the proof of \eqref{eq:isoM}--\eqref{eq:aveM}.

\subsection{General traceless observables}  \label{sec:gen_tr}
Recall that the local laws for traceless  observables in Theorem~\ref{th:local_laws} have
already been proven in Section~\ref{sec:traceless} for special observables $S^x$
and their traceless versions $\trless{S}^x$ (Theorem~\ref{th:traceless_laws}); now we explain
how to extend them to general traceless observables and general test vectors.

Since the improvement  due to tracelessness of the observables
in the local laws  Theorem~\ref{th:traceless_laws}   
  and in corresponding $M$-bounds in Lemma \ref{lemma:traceless_M_bounds} unfold under the condition $t \ge \crit$, the spatial decay  is completely washed out, generalizing them to the case of general test vectors and diagonal (traceless) observables becomes trivial. Indeed, the difficulties we encountered in Sections \ref{sec:gen_prelim}--\ref{sec:genobs} above stemmed entirely from the loss factor $\sqrt{\ell_t/\ell_s}$ in the propagator bound \eqref{eq:gen_P_st} for general test vectors and observables compared to its specialized counterpart \eqref{eq:P_st_bound}. For $\crit \le s \le t \le T$, this  loss
  factor  is of order one,  as $\ell_t\sim\ell_s\sim N$, 
  and hence becomes irrelevant. The final ingredient for ensuring that the correct number of propagators are traceless is the identity
\begin{equation}
	 \mathcal{S}\bigl[ M_{[1,2],t} (\trless{A} ) \bigr]  =  \sum_{q} \trless{S}^q \bigl(M_{[1,2],t} (\trless{A} )\bigr)_{qq},
\end{equation}
which is an immediate consequence of the fact that $\Pi = \bm 1 \bm 1^*$ commutes with $S$, and replaces \eqref{eq:trless_lin_term}.
We leave the remaining straightforward modifications to the reader.

The corresponding $M$-bounds \eqref{eq:tr_isoM}--\eqref{eq:tr_aveM} are established following the same proof as in Section~\ref{sec:trless_M}.

\section{Real-symmetric case} \label{sec:real}
In this section, we explain how to extend the proofs from all previous sections to
incorporate the non-vanishing off-diagonal component $\mathscr{T}$ of the self-energy operator (see Section~\ref{sec:calS}), thereby including the real-symmetric case in our analysis.
 To highlight the required modifications, we focus on the crucial zig-step (Proposition \ref{prop:zig}) and detail the necessary changes to its proof. These modifications apply directly to the global laws (Proposition \ref{prop:global_laws}), while the zag-step (Proposition \ref{prop:zag}) remains valid without modification.

As established in \eqref{eq:Sdecomp} and further detailed in \eqref{eq:entries}--\eqref{eq:calT}, the self-energy operator $\mathcal{S}$
decomposes into its diagonal and off-diagonal components, $\mathscr{S}$ and $\mathscr{T}$, respectively.
Recall, in particular, the expression $\mathscr{T}[R] = S^\mathrm{od}\odot R^\mathfrak{t}$. 
Since, by assumption, $\Expv [|h_{ab}|^2] = 1$ for all $a,b \in \indset{N}$ in both in the complex-Hermitian and real-symmetric cases, the entries of the $S^\mathrm{od}$ matrix satisfy the inequality
\begin{equation} \label{eq:torS_ineq}
	\bigl\lvert \bigl(S^\mathrm{od}\bigr)_{ab} \bigr\rvert \le \delta_{a\neq b} \, S_{ab}, \quad a,b\in \indset{N}.
\end{equation}
In the remainder of this section, we specifically consider the case where  $H$ is a real symmetric random band matrix, 
$\Expv (h_{ab})^2= \Expv |h_{ab}|^2 = S_{ab}$, which implies
\begin{equation}
	  \bigl(S^\mathrm{od}\bigr)_{ab} = \delta_{a\neq b}S_{ab}, \quad a,b\in \indset{N}.
\end{equation}
While the arguments presented in this section can be easily extended to a general off-diagonal $S^\mathrm{od}$ satisfying \eqref{eq:torS_ineq}, we do not pursue this generalizations\footnote{The case of general $\Expv (h_{ab})^2$
requires studying resolvent chains that also include transposes of the resolvent $G^\mathfrak{t}$. While
transposes are usually harmless for the main singularity of our problem, e.g. $\langle GG^\mathfrak{t}\rangle$
is much less singular than $\langle GG^*\rangle$ in the small $\eta$-regime, see e.g.,
Prop. 3.3 \cite{cipolloni2023functional}, properly bookkeeping them would significantly complicate the notations.}
 for the sake of presentation.

We note that since all observables under consideration are diagonal, and $\mathscr{T}$ vanishes on any diagonal matrix, the deterministic $M$-terms (with diagonal observables) do not depend on the off-diagonal component $\mathscr{T}$. In fact, all instances of $\mathcal{S}$ in the definitions of $M$-terms in~\eqref{eq:M_recursion}, hence
in the entire Section~\ref{sec:M}) can be replaced by its diagonal components $\mathscr{S}$. Furthermore, owing to  \eqref{eq:torS_ineq}, the off-diagonal component $\mathscr{T}$ does not affect the bounds on the quadratic variation of the martingale terms in \eqref{eq:k_av_evol} and \eqref{eq:k_iso_evol}---specifically, it is easy to verify that the estimates \eqref{eq:QV_Q_bounds} remain valid (possibly with a larger implicit constant).  

However, the presence of $\mathscr{T}$ does impact other parts of the proof. Specifically, the explicit diagonal form of $\mathcal{S}$ was presumed in many estimates for forcing terms of the form $G \mathcal{S}[G]G$. These bounds remain valid only when $\mathcal{S}$ is replaced by $\mathscr{S}$.   
Therefore, we focus on estimating the new forcing terms in \eqref{eq:k_av_evol} and \eqref{eq:k_iso_evol} 
 that emerge from the off-diagonal component $\mathscr{T}$, which we  refer to these as \emph{skew-forcing} terms. 
 They are given by
\begin{equation} \label{eq:torsion_forcing}
	\Tr\bigl[\mathscr{T}[G_{[i,j],t}] \, G_{[j,k],t}A_k G_{[1,i],t}\bigr] \quad \text{or} \quad \bigl(G_{[1,i],t}\mathscr{T}[G_{[i,j],t}] \, G_{[j,k],t}\bigr)_{\bm u \bm v},
\end{equation}
in the averaged and isotropic cases, respectively, for all $1 \le i \le j \le k$.

We recall that in the main proof, only   analogous terms involving  $\mathscr{S}$ instead of $\mathscr{T}$ were present. 
There, we utilized the fact that $\mathscr{S}[G_{[i,j],t}]$ is a diagonal matrix with entries $\Tr [S^a G_{[i,j],t}]$,
directly relating the estimates of such forcing terms to averaged chains, allowing for the use of stronger estimates, as evidenced by the appearance of $\psi^{\mathrm av}_k$ quantities in, for example,
\eqref{eq:av_forcing_bound}--\eqref{eq:iso_forcing_bound}. For the skew-forcing terms we cannot achieve the same structure. However, we can still control them in a way that preserves the self-improving structure of the master inequalities, after adjusting some of the loss exponents $\beta_k$.  We now explain
these modifications in some details.

\vspace{5pt}
\textbf{Averaged chains}.
First, we consider the skew-forcing terms arising from the averaged chains.  Note that we have two different representations for them:
\begin{equation}\label{eq:2rep}
	\begin{split}
		\Tr\bigl[\mathscr{T}[G_{[i,j],t}] \, G_{[j,k],t}A_k G_{[1,i],t}\bigr] =&~ \sum_{ a,b\,:\,a\neq b} S_{ab} \bigl(G_{[i,j],t}\bigr)_{ba}  \bigl(G_{[j,k],t}A_k G_{[1,i],t}\bigr)_{ba}\\
		=&~\sum_{ ab} \bigl(G_{[i,j],t}  S^b (G_{[j,k],t}A_k G_{[1,i],t})^\mathfrak{t} \bigr)_{bb}\\
		&- \sum_{ a} \bigl(G_{[i,j],t}\bigr)_{aa} S_{aa} \bigl(G_{[j,k],t}A_k G_{[1,i],t}\bigr)_{aa}.
	\end{split}
\end{equation}
Therefore, the averaged skew-forcing terms can be treated as both $(k+2)$-chains and products of off-diagonal shorter chains.

Similarly to the proof of Lemma \ref{lemma:forcing} in Section \ref{sec:forcing}, the critical terms correspond to $i=j$. For $i=j$ the chain $G_{[i,j],t}$ contains no observables, estimating \eqref{eq:torsion_forcing} in this case requires the use of \eqref{eq:sqrt_convol}, which incurs the loss of $\sqrt{N/\ell_t}$.

Therefore, we only analyze the critical terms $i=j$ in full detail,
for the less problematic $i<j$ case we just give the final estimate
and skip the details. Without loss of generality, we assume that $i=j=1$, since all other terms are structurally identical.

\textbf{Case 1.} Assume that $k \in \indset{\maxK-2}$. Using \eqref{eq:G-M_psi_bounds} for an isotropic $(k+2)$-chain in the second representation in~\eqref{eq:2rep},  we obtain 
\begin{equation} \label{eq:skew_crit_av}
	\begin{split}
		\bigl\lvert \Tr\bigl[\mathscr{T}[G_{1,t}] \, G_{[1,k],t}A_k G_{1,t}\bigr]\bigr\rvert \prec&~ \biggl\lvert \sum_b\bigl(M(z_{1,t}, S^b, z_{1,t}, A_1, \dots, z_{k,t}, A_k, z_{1,t})\bigr)_{bb} \biggr\rvert\\
		&+ (\ell_t\eta_t)^{\alpha_{k+2}}\pis{k+2} \sum_b \mathfrak{s}_{k+2,t}^\mathrm{iso}(b,b,A_1,\dots, A_k,b)\\
		&+\frac{1}{W}\biggl(1 + \frac{\pis{1}}{\sqrt{\ell_t\eta_t}}\biggr)\biggl\lvert \sum_a\bigl(M(z_{1,t}, A_1, \dots, z_{k,t}, A_k, z_{1,t})\bigr)_{aa} \biggr\rvert\\
		&+ \frac{(\ell_t\eta_t)^{\alpha_{k+1}}\pis{k+1}}{W\sqrt{\ell_t\eta_t}}\biggl(1 + \frac{\pis{1}}{\sqrt{\ell_t\eta_t}}\biggr) \sum_a \mathfrak{s}_{k+1,t}^\mathrm{iso}(a,A_1,\dots, A_k,a)\\
		\prec&~ \biggl(1 + \frac{\sqrt{N\eta_t}}{\ell_t\eta_t}\frac{\pis{k+2}}{(\ell_t\eta_t)^{\beta_k-\alpha_{k+2}}} +  
		\frac{\sqrt{N}}{W} 
		\frac{ \pis{k+1}}{  (\ell_t\eta_t)^{\beta_k-\alpha_{k+1}} }\biggl(1 + \frac{\pis{1}}{\sqrt{\ell_t\eta_t}}\biggr)\biggr)\\
		&\times \frac{1}{\eta_t}(\ell_t\eta_t)^{\beta_k}\, \mathfrak{s}_{k,t}^\mathrm{av}(A_1,\dots, A_k),
	\end{split}
\end{equation}
where we used  \eqref{eq:true_convol} together with a Schwarz inequality (see \eqref{eq:sqrt_convol})  to estimate the sum of $\mathfrak{s}$ over $a$  and $b$.
In the last line we separated the target bound and the penultimate line contains the critical
factor that has to be added to the master inequality for $\Psi_{k}^\mathrm{av}$. 
For the deterministic $M$-term, we used the bound
\begin{equation} \label{eq:resum2_M_bound}
	\biggl\lvert \sum_b\bigl(M(z_{1,t}, S^b, z_{1,t}, A_1, \dots, z_{k,t}, A_k, z_{1,t})\bigr)_{bb} \biggr\rvert \prec \frac{1}{\eta_t}\mathfrak{s}_{k,t}^\mathrm{av}(A_1,\dots, A_k),
\end{equation}
which is proved similarly to \eqref{eq:resum_M_bound}. Note that the $M$-term in \eqref{eq:resum2_M_bound} depends on $k+2$ spectral parameters. Recall that $\beta_k = \alpha_{k+1}$ 
 for any $k \in \indset{\maxK -1}$,  hence the
most critical factor in front of  $\pis{k+2}$  is 
\begin{equation}
	 \frac{\sqrt{N\eta_t}}{\ell_t\eta_t}\frac{1}{(\ell_t\eta_t)^{\beta_k-\alpha_{k+2}}} \le \frac{\sqrt{N\eta_t}}{(\ell_t\eta_t)^{1-\sqrt{1/(2\maxK)}}} 
\end{equation}
using $\beta_k - \alpha_{k+2} \ge \alpha_{k+1}-\alpha_{k+2} \ge 1/\sqrt{2K}$. This factor 
is  much less than $N^{-\etaexp/20}\ll N^{-2\nu}$, hence  affordable in the self-improving scheme, as long as
\begin{equation} \label{eq:newK}
	\maxK \ge 8/\bandexp^2,
\end{equation}
using that $(\ell_t\eta_t)^{-1} \sim (W\sqrt{\eta_t})^{-1} + (N\eta_t)^{-1}$, $N\eta_t \ge N^{\etaexp}$, 
 $N/W^2 \le N^{-\bandexp}$ and  the relations among the small exponents   
  $\nu\le \etaexp/(10K^2) \le  \bandexp/(10K^2)$.   

\textbf{Case 2.} Assume that $k = \maxK -1$, then we use the first representation in~\eqref{eq:2rep} and we obtain
\begin{equation}
	\begin{split}
		\bigl\lvert \Tr\bigl[\mathscr{T}[G_{1,t}] \, G_{[1,\maxK -1],t}A_{\maxK -1} G_{1,t}\bigr]\bigr\rvert &\prec \frac{(\ell_t\eta_t)^{\alpha_{\maxK}}\pis{\maxK}}{\sqrt{\ell_t\eta_t}} \sum_b \mathfrak{s}_{\maxK, t}^\mathrm{iso}(b,A_1,\dots,A_{\maxK -1},b)\\
		&\prec \pis{\maxK} \sqrt{N\eta_t} \times \frac{1}{\eta_t}\mathfrak{s}_{\maxK -1,t}^\mathrm{av}(A_1,\dots, A_{\maxK -1}),
	\end{split}
\end{equation}
where in the last step we used that $\alpha_\maxK = 1/2$.  The resulting loss of order $\sqrt{N \eta_t}$
is affordable if it is smaller than the corresponding loss factor, 
$(\ell_t\eta_k)^{\beta_{\maxK-1}}$, i.e. if, in addition to assuming \eqref{eq:newK}, we redefine the exponent $\beta_{\maxK-1}$ to be 
\begin{equation} \label{eq:new_betaK-1}
	\beta_{\maxK-1} := 1 - \bandexp/3.
\end{equation}

\textbf{Case 3.} If $k=\maxK$, using again the first representation in~\eqref{eq:2rep} and proceeding
similarly to \eqref{eq:1_k+1_reduction1}, we have
\begin{equation}
	\begin{split}
		\bigl\lvert \Tr\bigl[\mathscr{T}[G_{1,t}] \, G_{[1,\maxK],t}A_\maxK G_{1,t}\bigr]\bigr\rvert &\prec \frac{(\ell_t\eta_t)^{1/2+\alpha_{\maxK/2+1}}\pis{\maxK/2}\pis{\maxK/2+1}}{\sqrt{\ell_t\eta_t}} \sum_b \mathfrak{s}_{\maxK, t}^\mathrm{iso}(b,A_1,\dots,A_\maxK,b)\\
		&\prec \pis{\maxK/2+1} \sqrt{N\eta_t}(\ell_t\eta_t)^{\alpha_{\maxK/2+1}} \times \frac{1}{\eta_t} \mathfrak{s}_{k,t}^\mathrm{av}(A_1,\dots, A_k).
	\end{split}
\end{equation}
This loss is affordable if, in addition to assuming \eqref{eq:newK} and \eqref{eq:new_betaK-1}, we modify the definition of $\beta_\maxK$ to be
\begin{equation} \label{eq:new_betaK}
	\beta_{\maxK} := 1 - \bandexp/3 + \alpha_{\maxK/2+1} = 1 - \bandexp/3 + \sqrt{\frac{1}{2\maxK}}.
\end{equation}

 Finally, for the less critical cases, 
it is straightforward to check that all other averaged skew-forcing terms in \eqref{eq:torsion_forcing}, those corresponding to $1\le i < j \le k$ with $l := j-i+1 \in [2,k]$, admit the  bound  
\begin{equation}
	\frac{\eta_t\bigl\lvert \Tr\bigl[\mathscr{T}[G_{[i,j],t}] \, G_{[j,k],t}A_k G_{[1,i],t}\bigr] \bigr\rvert}{(\ell_t\eta_t)^{\beta_k}\,\mathfrak{s}_{k,t}^\mathrm{av}(A_1,\dots, A_k) }  \lesssim   
	\begin{cases}
		1 + \frac{\pis{k+2}}{(\ell_t \eta_t)^{1/4}} + \frac{\pis{l}\pis{k-l+2}}{\sqrt{\ell_t\eta_t}}, \quad & k \in \indset{\maxK-2}\\
		\frac{\pis{l}\pis{k-l+2}}{(\ell_t\eta_t)^{\beta_k-\alpha_k}}, \quad & k \in \{\maxK-1, \maxK\},
	\end{cases}
\end{equation}
for all $0 \le t \le \tau$, where $\tau$ is the stopping time, defined in \eqref{eq:tau_def}. Here we used additionally that all deterministic approximation $M$-terms are diagonal, and hence $\mathscr{T}[M_{\dots}]\equiv 0$.
All these terms are self-improving in the sense of the master inequalities.

It is straightforward to check that the new choices \eqref{eq:new_betaK-1} and \eqref{eq:new_betaK}, together with the assumption \eqref{eq:newK}, are compatible with the rest of master inequalities \eqref{eq:masters}.
In fact we could have defined $\beta_{K-1}$ and $\beta_K$ to be  these new values 
from the very beginning in~\eqref{eq:loss_exponents}, but over there such choice would have been unmotivated.

\vspace{5pt}
\textbf{Isotropic Chains}.
Isotropic skew-forcing terms are handled similarly. In fact, their treatment is more straightforward since the target isotropic bound is weaker. Note that we again have two representations
\begin{equation}\label{eq:2repiso}
	\begin{split}
		\bigl(G_{[1,i],t}\mathscr{T}[G_{[i,j],t}] \, G_{[j,k],t}  \bigr)_{\bm u \bm v} =&~ \sum_{ a,b\,:\,a\neq b} \bigl(G_{[1,i],t}\bigr)_{\bm u a} S_{ab} \bigl(G_{[i,j],t}\bigr)_{ba}  \bigl(G_{[j,k],t} \bigr)_{b \bm v}\\
		=&~\sum_{ b} \bigl(G_{[1,i],t} S^b (G_{[i,j],t})^\mathfrak{t} \bigr)_{\bm u b} \bigl(G_{[j,k],t} \bigr)_{b \bm v}\\
		&- \sum_{ a} \bigl(G_{[1,i],t}\bigr)_{\bm u a} S_{aa} \bigl(G_{[i,j],t}\bigr)_{aa}  \bigl(G_{[j,k],t} \bigr)_{a\bm v}.
	\end{split}
\end{equation} 
This decomposition allows us to treat the isotropic skew-forcing terms as products of two shorter chains.
Similarly to the averaged case above, the critical terms correspond to $i=j$, so we analyze them first.

\vspace{5pt}
\textbf{Case 1.} Assume that $k \in \indset{\maxK-1}$, the isotropic skew-terms are handled similarly to their averaged counterparts in \eqref{eq:skew_crit_av}: based upon the second representation in~\eqref{eq:2repiso} it yields
\begin{equation}
	\frac{\eta_t\bigl\lvert \bigl(G_{[1,i],t} \mathscr{T}[G_{[i,j],t}] \, G_{[j,k],t}\bigr)_{\bm u \bm v} \bigr\rvert}{(\ell_t\eta_t)^{\alpha_k}\,\mathfrak{s}_{k,t}^\mathrm{iso}(\bm u, A_1,\dots, A_{k-1}, \bm v) } \lesssim 1 + \frac{\sqrt{N\eta_t}}{\ell_t\eta_t}\frac{\pis{j+1}\pis{k-j+1}}{(\ell_t\eta_t)^{\alpha_k-\alpha_{k+1}}} + \frac{\sqrt{N\eta_t}}{\ell_t\eta_t}\pis{j}\pis{k-j+1}\biggl(1+ \frac{\pis{1}}{\sqrt{\ell_t\eta_t}}\biggr).
\end{equation}

\vspace{5pt}
\textbf{Case 2.} For $k  = \maxK$, the same bound holds, since by symmetry we may assume that $i=j \le \maxK-1$, otherwise $k-j+1$ plays the role of $i$.

Now we consider  all other isotropic skew-forcing terms when $j-i\ge 1$.  Assuming first 
that $k \in \indset{\maxK-1}$, we have
\begin{equation}\label{eq:GTGG}
	\frac{\eta_t\bigl\lvert \bigl(G_{[1,i],t} \mathscr{T}[G_{[i,j],t}] \, G_{[j,k],t}\bigr)_{\bm u \bm v} \bigr\rvert}{(\ell_t\eta_t)^{\alpha_k}\,\mathfrak{s}_{k,t}^\mathrm{iso}(\bm u, A_1,\dots, A_{k-1}, \bm v) }  \lesssim  
				1 + \frac{\bigl(\pis{i}+\pis{j+1}\bigr)\pis{k-j+1}}{(\ell_t \eta_t)^{1/2 + \alpha_k - \alpha_{k+1}}} + \frac{\pis{i}\pis{j-i+1}\pis{k-j+1}}{\ell_t\eta_t },
\end{equation}
where we used the subadditivity of the loss exponents, $\alpha_{j+1}+ \alpha_{k-j+1} \le \alpha_{k}$. This term is affordable since $\alpha_{k+1} - \alpha_k \le \sqrt{1/(2\maxK)} \le 1/4$. An analogous bound holds for $k=\maxK$, provided $i \ge 2$ or $ j \le k-1$, since in this case all three chains in~\eqref{eq:GTGG} have length at most $k-1$, so 
all $\psi^{\mathrm{iso}}$ factors have smaller indices than $k$, ensuring the self-improving mechanism.

However, the term corresponding to $i=1$ and $j=k=\maxK$ requires a special treatment. This is the only isotropic forcing term where performing a reduction is necessary (and, as we will see, is affordable). Using Schwarz followed by reduction \eqref{eq:reduction} and another Schwarz, we obtain,
\begin{equation}
	\begin{split}
		\bigl\lvert \bigl(G_{1,t} \mathscr{T}[G_{[1,\maxK],t}] \, G_{\maxK,t}\bigr)_{\bm u \bm v} \bigr\rvert \lesssim&~
		 \sum_{b} \bigl\lvert (G_{\maxK,t})_{b\bm v} \bigr\rvert \biggl(\sum_a S_{ab} \bigl\lvert (G_{1,t})_{\bm ua} \bigr\rvert^2 \biggr)^{1/2}    \biggl(\sum_a S_{ab} \bigl\lvert( G_{[1,\maxK],t})_{ba} \bigr\rvert^2 \biggr)^{1/2}\\
		 \lesssim&~
		 \sum_{b} \bigl\lvert (G_{\maxK,t})_{b\bm v} \bigr\rvert \,\Tr\bigl[(G_{[q,\maxK],t})^*|A_q| G_{[q,\maxK],t} S^b  \bigr] ^{1/2}\\
		 &\quad\times   \biggl(\sum_a S_{ab} \bigl\lvert (G_{1,t})_{\bm ua} \bigr\rvert^2 \bigl( G_{[1,q],t}|A_q|  (G_{[1,q],t})^*\bigr)_{bb} \biggr)^{1/2}\\
		 \lesssim&~
		  \biggl(\sum_{ab} \bigl\lvert (G_{\maxK,t})_{b\bm v}\bigr\rvert^2 S_{ba} \bigl( (G_{[q,\maxK],t})^*|A_q| G_{[q,\maxK],t} \bigr)_{aa} \biggr) ^{1/2}\\
		 &\times   \biggl(\sum_{ab}  \bigl\lvert (G_{1,t})_{\bm ua} \bigr\rvert^2  S_{ab} \bigl( G_{[1,q],t}|A_q|  (G_{[1,q],t})^*\bigr)_{bb} \biggr)^{1/2}\\
		 \prec&~ \frac{1}{\eta_t} (\ell_t\eta_t)^{\alpha_\maxK}\mathfrak{s}_{\maxK,t}^\mathrm{iso}(\bm u, A_1,\dots, A_{k-1}, \bm v)   \biggl(1 +  \frac{\pav{\maxK}}{(\ell_t\eta_t)^{1-\beta_\maxK}} + \frac{\pis{2}\pis{\maxK}}{\sqrt{\ell_t\eta_t}}\biggr),
	\end{split}
\end{equation}
where $q:=\maxK/2$ (recall that $\maxK$ is even) is the cutting index. 
In the last step we used various relations on the loss exponents, and the identity,
\begin{equation}\label{eq:traceback}
	\begin{split}
	 	\sum_{ab}  \bigl\lvert (G_1)_{\bm ua} \bigr\rvert^2  S_{ab} \bigl( G_{[1,\maxK] } \bigr)_{bb} =&~
		\sum_{a c} |u_a|^2  \biggl(\frac{|m_1|^2}{1-|m_1|^2 S}\biggr)_{ac}  \Tr\bigl[ G_{[1,\maxK] }S^c \bigr]\\
		&+ \sum_{b}   \bigl(G_1S^b G_1^*-\Theta^b(z_1)\bigr)_{\bm u\bm u}\bigl( G_{[1,\maxK] } \bigr)_{bb},
	\end{split}
\end{equation}
where we wrote $G_{[1,\maxK]}= G_{[1,q],t}|A_q|  (G_{[1,q],t})^*$ for brevity.
Note that  after  using the reduction inequality followed by~\eqref{eq:traceback} 
  we recovered an averaged chain $ \Tr\bigl[ G_{[1,\maxK] }S^c \bigr]$ in the leading first  term
  in~\eqref{eq:traceback}, which was initially lost due to the action of
off-diagonal operator $\mathscr{T}$ instead of  $\mathscr{S}$, see explanation below~\eqref{eq:torsion_forcing}.

It is straightforward to check that these new terms are compatible with the master inequalities \eqref{eq:masters} given \eqref{eq:new_betaK-1} and \eqref{eq:new_betaK}, together with the assumption \eqref{eq:newK}. This concludes the analysis of the new forcing terms arising from the off-diagonal part $\mathscr{T}$ of the self-energy operator $\mathcal{S}$.

\section{Wigner-Dyson universality}\label{sec:WD}
In this section we prove Theorem~\ref{thm:WD}. 
To establish bulk universality for a random band matrix $H$, we dynamically connect
 it to the Gaussian ensemble via the Ornstein-Uhlenbeck process
\begin{equation} \label{eq:mean-field-flow}
	\mathrm{d}H_t = -\frac{1}{2}H_t + \frac{1}{\sqrt{N}}\mathrm{d}\mathfrak{B}_t, \quad H_0 = H,
\end{equation}
where $\mathfrak{B}_t$ is the standard Brownian motion in the space of $N\times N$ matrices of the same symmetry class as $H$. 
Note that, unlike \eqref{eq:zigOU}, the process \eqref{eq:mean-field-flow} introduces an order $\sqrt{1-\mathrm{e}^{-t}}$ mean-field Gaussian component to $H_t$ at time $t > 0$. In particular, at time $t=\infty$, the process \eqref{eq:mean-field-flow} equilibrates to GUE (for complex Hermitian $H$) or GOE (for real symmetric $H$).

In fact, as was shown in Theorem 2.2 of \cite{landon2019fixed} along the Ornstein-Uhlenbeck process \eqref{eq:mean-field-flow}, the correlation functions of the random matrix $H_t$ already equilibrates after a very short time $\widehat{t} := N^{-1+\varepsilon_*}$ with any fixed positive $\varepsilon_* > 0$. More precisely, for any $t \ge \widehat{t}$  the correlation functions of $H_t$ around a fixed energy $E \in [-2+\kappa, 2-\kappa]$ are universal, that is
\begin{equation}
	\int\limits_{\mathbb{R}^k} f(\bm x) \bigl(\rho_{H_{\widehat{t}}}^{(k)}-\rho_{\mathrm{GUE/GOE}}^{(k)}\bigr)
	\Big( E+ \frac{\bm x}{N\varrho_{sc}(E)}\Big) \mathrm{d}\bm x = o(1), \quad N \to \infty.
\end{equation}
for any fixed $k\in\mathbb{N}$ and all smooth compactly supported test function $f$ of $k$-variables,
 where $\rho_{H_t}^{(k)}$ denotes the $k$-th order correlation function of the eigenvalues 
 of $H_t$, and $\bm x := (x_1, \dots, x_k)$. Recall that $\varrho_{sc}$ is the semicircle density, as defined above \eqref{eq:Dyson}.

Therefore, it remains to compare $\rho_{H}^{(k)} = \rho_{H_0}^{(k)}$ to $\rho_{H_{\widehat{t}}}^{(k)}$, 
i.e. track the change when a  very small Gaussian component of order $\widehat{t}$ is added. 
This is typically done by a Green function comparison arguments with moment matching, a standard method
that is designed to bridge the deviations in the  third and higher moments between $H_0$ and $H_{\widehat{t}}$
(see e.g. Chapter 16 of \cite{erdHos2017dynamical}). More precisely, one shows that the  quantity 
\begin{equation} \label{eq:ImG_prod}
	\Expv \prod_{j=1}^{\mathfrak{n}} \bigl\langle \im G_t(z_j) \bigr\rangle, \quad 0 \le t \le \widehat{t},
\end{equation}
does not change up to leading order
along the process \eqref{eq:mean-field-flow} from $t=0$ to $t=\wh t:= N^{-1+\varepsilon_* }$, provided 
$\varepsilon_*>0$ is sufficiently small. Here, $G_t(z) := (H_t - z)^{-1}$ denotes the resolvent of $H_t$, while 
for any fixed $\mathfrak{n} \in \mathbb{N}$, the spectral parameters $\{z_{j}\}_{j=1}^\mathfrak{n}$ satisfy $|\re z_j - E| \le C N^{-1}$ and $\eta_j := \im z_j \in [N^{-1-\delta}, N^{-1+\delta_0}]$, for  $j \in \indset{\mathfrak{n}}$, with some fixed small $\delta >0$ to be chosen later. Recall that $\langle A\rangle=\frac{1}{N}\Tr A$ denotes the average trace.

It is important to note that the imaginary part of the spectral parameter $\eta_j = \im z_j$ in \eqref{eq:ImG_prod}  needs to go a bit below $1/N$.
 This places it outside the outside the regime where the local laws typically apply, as $\langle \im G_t(z_j) \bigr\rangle$ exhibits large fluctuations on
this scale; however, the expectation dampens these fluctuations.

To control  the time derivative of~\eqref{eq:ImG_prod}
within the standard proofs,  single-resolvent  local laws~\eqref{eq:law1} are routinely used in conjunction with the monotonicity  of 
the function $\eta_j\to\eta_j \langle \im G_t(z_j) \bigr\rangle$, that connects the small $\eta_j$ values
with the ones above $1/N$ where local laws are applicable.
While our main theorems state local laws for band matrices, (i.e. for $H_0$),
it is straightforward to extend the local laws in Theorems~\ref{th:local_laws} and~\ref{th:local_laws_traceless} to the resolvent $G_t$ of $H_t$, which has a small mean-field Gaussian component. This extension uses the same approach as in the proof of Proposition~\ref{prop:zig} with $\tinit = \crit$, provided $\varepsilon_* < \bandexp$ and $\im z \in (N^{-1+\etaexp}, (W/N)^2 \,)$.
Note that on this scale the spatial decay is already washed out; therefore, adding a small mean-field
component is irrelevant to the proofs.

The major difference between the current band matrix argument and its well-established counterpart in the mean-field case
 is that the OU process~\eqref{eq:mean-field-flow} does not
preserve the second moments. When tracing the change of~\eqref{eq:ImG_prod}
over a short time $\widehat{t}$, a simple power counting argument shows that discrepancies in the third and higher moments can be managed through brute-force comparison, relying solely on single-resolvent local laws.
The discrepancy in the second moment, however, necessitates  a new input. 
This input is provided by the averaged two-resolvent local law, \eqref{eq:tr_avelaw} for $k=2$, where both observables are
traceless ($n^*=2$), with the corresponding gain factor $(N\sqrt{\eta}/W)^2$ being essential for the proof.
Technically it will be used via the QUE for eigenfunctions.
This extension of standard GFT proofs to band matrices, incorporating additional information on the eigenfunctions, was first realized in \cite{xu2024bulk}, specifically in Proposition~4.17.

The difference
\begin{equation}
	\Expv \biggl[ \prod_{j=1}^{\mathfrak{n}} \bigl\langle \im G_t(z_j) \bigr\rangle \biggr] - \Expv \biggl[ \prod_{j=1}^{\mathfrak{n}} \bigl\langle \im G_{\widehat{t}}(z_j) \bigr\rangle \biggr]
\end{equation}
is controlled uniformly in $t \in [0, \widehat{t}]$ by differentiating \eqref{eq:ImG_prod} along the flow \eqref{eq:mean-field-flow} using It\^{o}'s formula together with cumulant expansion. The crucial second cumulant term in the expansion gives rise to two quantities\footnote{If $H$ is a real symmetric random band matrix, the second-order cumulants also give rise to
\begin{equation}
		\other{Z}_{2,t} := \frac{1}{N}\sum_{ab} \trless{S}^b_{aa} \bigl(G_t^2(z_u)\bigr)_{ab} \bigl(G_t(z_u)\bigr)_{ab}.
\end{equation}
The quantity $\other{Z}_{2,t}$ is estimated analogously to $Z_{1,t}$ and $Z_{2,t}$ using the fact that $G^\mathfrak{t} = G$ for real symmetric $H$.
}
\begin{subequations}
	\begin{equation} \label{eq:univ_2cum_mixed}
		Z_{1,t} := \frac{1}{N^2}\sum_{ab} \trless{S}^b_{aa} \bigl(G_t^2(z_u)\bigr)_{ab} \bigl(G_t^2(z_v)\bigr)_{ba}
	\end{equation}
	\begin{equation} \label{eq:univ_2cum_repeat}
		Z_{2,t} := \frac{1}{N}\sum_{ab} \trless{S}^b_{aa} \bigl(G_t^2(z_u)\bigr)_{aa} \bigl(G_t(z_u)\bigr)_{bb}, 
	\end{equation}
\end{subequations}
which we control using the local laws with traceless observables from Theorem \ref{th:traceless_laws}. 
Note that $G_t(z_{u/v})$ in \eqref{eq:univ_2cum_mixed}--\eqref{eq:univ_2cum_repeat} can potentially replaced by $G^*_t(z_{u/v})$, albeit this will not be relevant for our proof.  

Similarly to Lemma 4.19 and Proposition 4.17 of \cite{xu2024bulk} suffices to show that 
\begin{equation} \label{eq:univZ_goal}
	\max_{0 \le t \le \widehat{t}}~\bigl\lvert Z_{1,t} \bigr\rvert \prec N^{1  + C(\delta+\etaexp) - \zeta}, 
	\qquad \max_{0 \le t \le \widehat{t}}~\bigl\lvert Z_{2,t} \bigr\rvert \prec N^{1  + C(\delta+\etaexp) - \zeta},
\end{equation}
for some positive constant $\zeta > 0$ depending on $\bandexp$ (recall the assumption \eqref{eq:WN}). 

 For terms involving higher-order cumulants, every resolvent entry in the summation can be estimated trivially by $N^{\delta+\etaexp}$ using the single resolvent isotropic local law together
 with the standard monotonicity argument, while the entries of $G^2$ can be similarly estimated by $N^{1+2(\delta+\etaexp)}$. These trivial estimates yield the bound
\begin{equation}
	\sum_{ab} \bigl(S_{ab}\bigr)^{l/2} \partial_{ab}^l  \prod_{j=1}^{\mathfrak{n}} \bigl\langle \im G_t(z_j) \bigr\rangle	\prec C_l  N  \frac{N^{(\mathfrak{n}+l)(\delta+\etaexp)}}{W^{l/2-1}} \lesssim C_l \frac{N^{1+\mathfrak{n}(\delta+\etaexp)}  }{N^{1/24 + \bandexp/2}} , \quad l \ge 3,
\end{equation}
on the contribution for the derivative 
coming from terms involving cumulants of order $l \ge 3$.  After integrating up to time $\widehat{t}=N^{-1+\varepsilon}$,
this contribution is negligible for a sufficiently small $\delta=\delta(\mathfrak{n})$  for any fixed  $\mathfrak{n}$.

We only present the bound on the quantity $Z_{1,t}$, defined in \eqref{eq:univ_2cum_mixed}, in full detail, since the quantity $Z_{2,t}$ is estimated similarly. By spectral decomposition of the resolvent $G_t$, we obtain
\begin{equation} \label{eq:univ_Z1}
	\bigl\lvert Z_{1,t}\bigr\rvert \le \frac{1}{N^2}\sum_b \sum_{ij} \frac{ \bigl\lvert\bigl\langle \bm u_{j,t}, \trless{S}^b \bm u_{i,t} \bigr\rangle \bigr\rvert}{|\lambda_{i,t} - z_u|^2 |\lambda_{j,t} - z_v|^2 }  \bigl\lvert \bm u_{i,t}(b) \bigr\rvert\bigl\lvert \, \bm u_{j,t}(b) \bigr\rvert,
\end{equation}
where $\bm u_{j,t} := \bigl(\bm u_{j,t}(b)\bigr)_{b=1}^N$ denotes the $\ell^2$-normalized eigenvector of $H_t$ corresponding to the eigenvalue $\lambda_{j,t}$ for all $j \in \indset{N}$.

Let $\varepsilon > 0$ be a positive exponent to be chosen later, and let $\eta := N^{-1+\varepsilon}$. Note that if the indices $i,j \in \indset{N}$ are such that $|\lambda_{i,t} - z_u| \le \eta$ and $|\lambda_{j,t} - z_v| \le  \eta$, then
\begin{equation}
	\bigl\lvert\bigl\langle \bm u_{j,t}, \trless{S}^b \bm u_{i,t} \bigr\rangle \bigr\rvert^2 
	\le  \eta^2\sum_{i,j}\frac{\eta^2 \bigl\lvert\bigl\langle \bm u_{j,t}, \trless{S}^b \bm u_{i,t} \bigr\rangle \bigr\rvert^2}{|\lambda_{i,t} - z|^2 |\lambda_{j,t} - z|^2} \le \eta^2 \Tr\bigl[ \bigl(\im G_t(z)\bigr) \trless{S}^b \bigl(\im G_t(z)\bigr) \trless{S}^b \bigr],
\end{equation}
where $z := E+ \I\eta$. Recall that $|\re z_{u/v} - E| \lesssim N^{-1}$. Hence, it follows from the traceless local law \eqref{eq:tr_avelaw} and the corresponding $M$-bound \eqref{eq:tr_aveM} that (c.f. \eqref{eq:QUE})
\begin{equation} \label{eq:close_overlap}
	\bigl\lvert\bigl\langle \bm u_{j,t}, \trless{S}^b \bm u_{i,t} \bigr\rangle \bigr\rvert \prec  \eta \sqrt{\frac{1}{N\eta} \frac{N^2\eta}{W^2} \biggl(1+ \frac{1}{N\eta}\biggr)} \prec N^{-1+\varepsilon-\bandexp/2} , \quad 0\le t\le \crit,
\end{equation}
for all $b\in \indset{N}$, provided $|\lambda_{i,t} - z_u| \le N^{\zeta}\eta$ and $|\lambda_{j,t} - z_v| \le \eta$. Note that for any $u \in \indset{\mathfrak{n}}$, we have
\begin{equation} \label{eq:ImG_tr}
	N^{-1}\sum_{i}\frac{1}{|\lambda_{i,t} - z_u|^2} = \frac{1}{\eta_u}\bigl\langle \im G_t(z_u)\bigr\rangle \prec N^{1+3(\delta+\etaexp)}, \quad 0\le t\le \crit,
\end{equation}
where we used monotonicity of the map $\eta \mapsto \eta \im G(E+\I\eta)$ in $\eta \ge 0$. Similarly, we have  
\begin{equation}
	N^{-1}\sum_{i : |\lambda_{i,t} - z_u| \ge N^{\zeta}\eta }\frac{1}{|\lambda_{i,t} - z_u|^2} \lesssim  \frac{1}{N^{\zeta}\eta}\bigl\langle \im G_t(z_u + \I N^\zeta\eta)\bigr\rangle \prec N^{1-\varepsilon}, \quad 0\le t\le \crit.
\end{equation}
Finally, for all $i,j\in\indset{N}$, it follows from \eqref{eq:sumS=1} and the delocalization bound \eqref{eq:deloc} with $\bm v := \bm e_a$ that
\begin{equation} \label{eq:far_overlap}
	\bigl\lvert\bigl\langle \bm u_{j,t}, \trless{S}^b \bm u_{i,t} \bigr\rangle \bigr\rvert \prec \frac{1}{N}\sum_a \bigl\lvert \bigl(\trless{S}^b\bigr)_{aa}\bigr\rvert \prec \frac{1}{N}, \quad 0\le t\le \crit.
\end{equation}

Plugging the estimates \eqref{eq:close_overlap}--\eqref{eq:far_overlap} into \eqref{eq:univ_Z1}, and using delocalization \eqref{eq:deloc} with $\bm v := \bm e_b$ to bound~$\bm u_{i,t}(b)$ and~$\bm u_{j,t}(b)$, we obtain
\begin{equation}
	\bigl\lvert Z_{1,t}\bigr\rvert \prec  N^{1+6(\delta+\etaexp) + \varepsilon-\bandexp/2}  
	+   N^{1+3(\delta+\etaexp)-\varepsilon}, \quad 0\le t\le \crit.
\end{equation}
Therefore, choosing $\varepsilon := \bandexp/4$, we conclude that  \eqref{eq:univZ_goal} holds with $\zeta := \bandexp/4$ and $C=6$. The quantity $Z_{2,t}$ is estimated similarly. 
As in Proposition 4.17 in \cite{xu2024bulk}, we conclude that 
\begin{equation}
	\biggl\lvert \Expv \biggl[ \prod_{j=1}^{\mathfrak{n}} \bigl\langle \im G_t(z_j) \bigr\rangle \biggr] - \Expv \biggl[ \prod_{j=1}^{\mathfrak{n}} \bigl\langle \im G_{\widehat{t}}(z_j) \bigr\rangle \biggr] \biggr\rvert \lesssim N^{-c\bandexp + C\mathfrak{n}(\delta+\etaexp)},
\end{equation}
which goes to zero if $\delta=\delta(\mathfrak{n})$ 
and $\etaexp=\etaexp(\mathfrak{n})$ are chosen sufficiently small.

The remainder of proof of bulk universality follows by a standard correlation function comparison argument (see Theorem 15.3 in \cite{erdHos2017dynamical}) by choosing $\varepsilon_*$ in $\widehat{t}$ and $\delta, \delta_0$ to be much smaller than $\bandexp$.  This completes the proof of Theorem~\ref{thm:WD}.

\section{Deterministic decay profile: Proof of Proposition~\ref{prop:admS} }\label{sec:reg}

In this section we prove part (i) of Proposition~\ref{prop:admS}, i.e.
we show that translation-invariant variance profiles~\eqref{eq:trinvS} satisfy
the assumptions $(i), (vi)$ and $(vii)$ of Definition~\ref{def:adm_ups_notime}
with the given control functions in Examples~\ref{lemma:poly_Ups_notime}--\ref{lemma:exp_Ups_notime}.
 Recall that the other assumptions $(ii)-(v)$ of Definition~\ref{def:adm_ups_notime} 
 for these  control functions are elementary calculations.
 
 The proof of part (ii) of Proposition~\ref{prop:admS} will be omitted, as it  follows from part (i) after tensorization,
if one considers the profile $\sigma$  as an $L\times L$ band matrix profile with band width $W=1$.

\begin{proof}[Proof of Proposition~\ref{prop:admS} (i)] The proof consists of making the heuristic argument
via Fourier transform after Proposition~\ref{prop:admS}  rigorous. A similar argument was given 
along the proof of Proposition 2.8 of~\cite{erdHos2013delocalization} in 
Appendix A.1 of \cite{erdHos2013delocalization} with less precise estimates 
(note that the error $O(W^{-2})$ in (2.37) of  \cite{erdHos2013delocalization} is not affordable now)
and in addition to bounds on $\Theta$ now we also need regularity estimates.
To achieve these stronger bounds, we will partly rely on residue calculations.

We write 
$$
 \Xi = \frac{m_1m_2S}{1-m_1m_2S} = \frac{S}{\zeta-S}, \qquad \zeta: = \frac{1}{m_1m_2},
$$
and similarly for $\Theta$ with $\zeta=1/(m_1\overline{m_2})$, see~\eqref{def:ThetaXigen}.  The first  part of the analysis is identical for both cases, but then we need to distinguish whether $\zeta$ is away from the spectrum of $S$
at a distance order one (case of $\Xi$) or $\zeta$ can be very close to 1 (case of $\Theta$ when $\eta$ is small).
In fact, there is no other possibility. To see this, note that $|\zeta|\ge 1$, since $|m(z)|\le 1$ and 
recall that
\begin{equation}\label{SpecS}
\mbox{Spec}(S)\subset [-1+\delta, 1]
\end{equation}
 for
some small $\delta>0$, uniformly in $N, W$, 
see Lemma A.1 \cite{erdHos2012bulk}.

We recall that $\mathbb{T}_N := \mathbb{Z}/(N\mathbb{Z})$ is the discrete $N$ torus and set $Q: =(2\pi W/N)\mathbb{T}_N$,
both equipped with the periodic distance. Set $L:= N/W$.
Let $\wh S^D $ be the discrete Fourier transform of $S$, i.e.
\begin{equation}\label{def:SD}
\wh S^D(p): = \sum_{x=1}^N S_{0x}\, e^{-ipx}, \qquad  p\in P:=\frac{1}{W}Q =\frac{2\pi}{N}\mathbb{T}_N,
\end{equation}
then we have
\begin{equation}\label{def:res1}
\biggl(\frac{S}{\zeta-S}\biggr)_{0x} =  \frac{1}{N}\sum_{p\in P} \frac{\wh S^D(p)}{\zeta-\wh S^D(p)} e^{ipx},
\qquad x\in \mathbb{T}_N.
\end{equation}
Note that the definition of $\wh S^D(p)$ naturally extends to a periodic function  on the real interval $[0,2\pi]$. Thus
we may  define the  $2\pi W$-periodic function
\begin{equation}\label{def:FF}
   F(q): = \wh S^D\Bigl(\frac{q}{W}\Bigr)=\frac{1}{W}  \sum_{x=1}^N f\Bigl(\frac{ x}{W}\Bigr) e^{-ix(q/W)},  \qquad q\in [0,2\pi W],
 \end{equation}
 where we used~\eqref{eq:trinvS} and with a slight abuse of notations, we consider $f$ as $L$-periodic function
 to avoid writing $|x|_N$ instead of $x$. Sometimes we think of $F(q)$ extended periodically onto the whole $\mathbb R$.
  Clearly $\{ \wh S^D(p) \; : \; p\in P\} =\mbox{Spec}(S)$ thus using \eqref{SpecS}, essentially  the same relation holds for
the extended functions
\begin{equation}\label{SpecF}
\bigl\{ F(q)\; : \;  q\in [0, 2\pi W] \bigr\} =
\bigl\{ \wh S^D(p) \; : \; p\in  [0,2\pi] \bigr\} \subset [-1+\delta/2, 1].
\end{equation}
 
 Set
\begin{equation}\label{def:R}
  R(p) : = \frac{\wh S^D(p)}{\zeta-\wh S^D(p)} = \frac{F(pW)}{\zeta- F(pW)}, 
\end{equation}
for all $p\in [0,2\pi]$, notice that $R(p)$ is $2\pi$-periodic. Let $\wh R$ be its Fourier transform, i.e.
\begin{equation}\label{def:RF}
   \wh R(n): = \frac{1}{2\pi}\int_0^{2\pi} R(p)\, e^{-inp} \mathrm{d} p, \qquad n\in \mathbb{Z}.
\end{equation}
Then the Poisson summation formula applied to~\eqref{def:res1} gives
\begin{equation}\label{Resrep}
\biggl(\frac{S}{\zeta-S}\biggr)_{0x} = \frac{1}{N}\sum_{p\in P} R(p) e^{ipx} = \sum_{k\in \mathbb{Z}} \wh R(Nk+x), 
\qquad x\in \mathbb{T}_N.
\end{equation}
The convergence of the sum will be verified along the estimates later by proving that $\wh R$ has a strong decay
at infinity

In order to compute the Fourier transform of $R$, we first record a few basic properties of $F$
that directly follow from its definition~\eqref{def:FF}, and the basic properties of $f$, namely
that $f$ is non-negative, symmetric and together with its derivative, it has a decay~\eqref{eq:fdecay}:
\begin{equation}\label{Fprop}
    F(q)\in \mathbb{R}, \quad F(q)=F(-q), \quad |F(q)|\le 1, \quad   |\partial^k F(q)|\lesssim \langle |q| \rangle^{-2}, \;\;  
    0\le k\le 4, \quad q\in [0,2\pi W],
\end{equation}
and    
\begin{equation}\label{Fprop1}
    \qquad F(0) =1, \qquad \partial F(0)=0, \qquad -\partial^2F(0) = 
    \frac{1}{W} \sum_{x=1}^N \Bigl(\frac{ x}{W}\Bigr)^2 f\Bigl(\frac{ x}{W}\Bigr)=:D_0 \sim 1.
\end{equation}
Since $F$ is periodic, $F(-q)$ means $F(2\pi W-q)$, and $|q|$ denotes the periodic absolute value: $|q|=|q|_W=
\min\{ q, 2\pi W-q\}$. We  will omit the index $W$ from the notation. 
The proofs of~\eqref{Fprop}--\eqref{Fprop1} are elementary, here we just indicate the main mechanism.
For example, to see that $|F(q)|$ decays for large $|q|$, we employ a summation by parts in~\eqref{def:FF};
\begin{equation}\label{sumparts}
\begin{split}
    |F(q)| = & \frac{1}{W}   \Bigg| \frac{1}{1-e^{-iq/W}} \sum_{x=1}^N f\Big(\frac{ x}{W}\Big)
     \big[ e^{-ix(q/W)} - e^{-i(x+1)q/W}\big] \Bigg| \\
       = & \frac{1}{W}   \frac{1}{|1-e^{-iq/W}|} \Bigg| \sum_{x=1}^N  \Big[ f\Big(\frac{ x}{W}\Big) - f\Big(\frac{ x+1}{W}\Big)\Big]
        e^{-ix(q/W)}\Bigg| \\
        \lesssim  & \frac{1}{|q|}  \frac{1}{W} \sum_{x=1}^N  \Big| f'\Big(\frac{ x}{W}\Big)\Big| \lesssim \frac{1}{|q|},
 \end{split}
 \end{equation}
 and a second summation by parts gives $|F(q)|\le |q|^{-2}$. For small $|q|$ regime we just use $|F(q)|\lesssim 1$
 directly from \eqref{def:FF} giving the bound  $|F(q)|\le \langle |q|\rangle^{-2}$.
 The argument for higher derivatives $\partial^k F(q)$ is exactly the same; the additional $(x/W)^k$ factor
 that the $\partial^k$-derivatives in~\eqref{def:FF} bring down will be compensated by the sufficiently fast decay
 of $f$ and its derivatives.

Notice that from \eqref{Fprop1} it  follows that
\begin{equation}\label{smallq}
  F(q) = 1 - D_0|q|^2 + O(|q|^4) 
\end{equation}
for small $|q|$. Moreover, 
 for every $\epsilon>0$ there is a $c_\epsilon>0$, with $c_\epsilon\to0$ as $\epsilon\to0$,
 such that
\begin{equation}\label{largeq}
  1- \re F(q) =  \frac{1}{W} \sum_{x=1}^N  f\Big(\frac{ x}{W}\Big) \big[ 1- \cos(qx/W)\big] \ge c_\epsilon, \qquad 
  |q|\ge \epsilon,
\end{equation}
recalling that $|q|$ is understood in the $2\pi W$ periodic sense. This follows from the fact that the overwhelming mass
of $f$ cannot be supported  in the regime $|x/W|\le\epsilon$ if $\epsilon$ is small.
 We fix a sufficiently small 
 $\epsilon \ll \sqrt{D_0}$ so that the error term $O(|q|^4)$
is negligible compared with the quadratic term.

Finally, in the case when $f$ is compactly supported, i.e. $|\mbox{supp} f|$ is bounded uniformly in $L$, 
then $F(q)$ is  $2\pi W$-periodic, real analytic and can be extended to a  holomorphic function in a strip 
\begin{equation}\label{strip}
   \Omega : = \{ q\in \mathbb{C}\; : \;
|\im q|\le c_4\}  
\end{equation}
 with some constant $c_4>0$ independent of $W$,
and it remains periodic $F(q) = F(q+2\pi W)$.

\medskip

{\bf Case 1: $\zeta$ is far away from $\mbox{Spec}(S)$.} This is the case when $\zeta=1/(m_1m_2)$
with $m_i=m(z_i)$, and $\im z_i>0$, $i=1,2$. Since we are in the bulk, $\zeta$ is away from the 
entire interval $[-1+\delta/2, 1]$, so using~\eqref{SpecF},
 $|\zeta- F(q)|\ge c$ uniformly in $q\in [0, 2\pi W]$ for some small constant $c$,
independent of $W$. Now from~\eqref{def:R}--\eqref{def:RF},  we can write
\begin{equation}\label{Rhat}
  \wh R(x) = \frac{1}{2\pi W} \int_0^{2\pi W} \frac{F(q) }{\zeta - F(q)}  e^{-iqx/W} \mathrm{d} q, \qquad x\in \mathbb{Z},
\end{equation}
in particular we immediately obtain $|\wh R(x) |\lesssim 1/W$ using the decay of $F(q)$ from~\eqref{Fprop}.
 In fact $\wh R$ itself has a fast decay since $F$ has bounded
derivatives, established in \eqref{Fprop}. 
To see this fast decay, we perform an integration by parts (no boundary terms, since $F$ is periodic):
$$
    \wh R(x) = \frac{1}{2\pi W} \int_0^{2\pi W}   \biggl(\frac{W}{ix}\biggr)
    \frac{ \zeta F'(q) }{(\zeta - F(q))^2} e^{-iqx/W}\mathrm{d} q, \qquad x\in \mathbb{Z}.
$$ 
Estimating the integrand by absolute value, using $|F'(q)|\lesssim \langle |q|\rangle^{-2}$,
together with the trivial bound $|\wh R(x) |\lesssim 1$, we thus conclude
$$
   |\wh R(x)|\lesssim \frac{1}{W}\biggl\langle\frac{|x|}{W} \biggr\rangle^{-1}.
$$
Performing $D$ integration by parts, we obtain a decay of order $\langle |x|/W\rangle^{-D}$.  Plugging
this into~\eqref{Resrep}, we obtain that $|\Xi_{0x}|\lesssim  \big(\Upsilon_{\eta=1}\big)_{0x}$ with $\Upsilon_\eta$ defined in 
\eqref{eq:polyUps_notime}, recalling that $\ell(\eta=1)=W$.

If $f$ is compactly supported, hence $F(q)$ is analytic in a strip $\Omega$ of width $c_4$,
 then we can shift the integration contour in~\eqref{Rhat}
away from the real axis (either below or above depending on the sign of $x$) to pick up an exponential
factor $\exp(- c_0|x|/W)$ with $c_0= c_4/2$. The decay of $F(q)$ for large $|q|$ still holds, thus we obtain
$$
    | \wh R(x)|\lesssim \frac{1}{W} \exp(- c_0|x|/W)
$$
in this case. This proves $|\Xi_{0x}|\lesssim  \big(\Upsilon_{\eta=1}\big)_{0x}$ with $\Upsilon_\eta$ defined in 
\eqref{eq:expUps_notime}. Thus, in both cases we verified part
$(i)$ of Definition~\ref{def:adm_ups_notime} for $\Xi$.
\medskip

{\bf Case 2: $\zeta$ is close to 1.} This is the case when $\zeta=1/(m_1\overline{m}_2)$
with $m_i=m(z_i)$, and $\im z_i>0$, $i=1,2$, are both small,
 and we need to estimate $\Theta$. In the 
most important prototypical case, $z_1=z_2$, thus $\zeta -1 \sim \eta$, see the first relation in~\eqref{eq:mbound}.
For general $z_1, z_2$ we have $|\zeta|\ge 1$, thus 
\begin{equation}\label{zetaF}
\mbox{dist}(\zeta, \mbox{supp}(F)) \ge \frac{1}{2} |\zeta-1|,
\end{equation}
by elementary geometry, using~\eqref{SpecF}.
We may assume that $|\zeta-1|\le c'$ with some small $c'\le c_\epsilon/2$, otherwise 
the analysis of $\Theta$ leading to the proof of $(i)$ of Definition~\ref{def:adm_ups_notime}
 is exactly the same as that of $\Xi$ in Case 1.  
 We will explain the  proofs of $(vi)$ and $(vii)$ of Definition~\ref{def:adm_ups_notime} for $\Theta$ in the case
 $|\zeta-1|\le c'$ below. The arguments  in the case $|\zeta-1|\ge c'$ are analogous but easier, hence will be left to the reader.

The assumption  $|\zeta-1|\le c'$ implies that $|z_1-z_2|\lesssim c'$ and $\im z_1, \im z_2\lesssim c'$.
We claim that
\begin{equation}\label{zeta}
 |\im \sqrt{1-\zeta} |\gtrsim\sqrt{ \eta}, \quad |\sqrt{1-\zeta}| \lesssim  |\im \sqrt{1-\zeta} |,
 \qquad \eta:=  \min \{ \im z_1, \im z_2\}.
\end{equation}
Indeed, writing $\sqrt{1-\zeta} = a+bi$ and $m_j =|m_j|e^{i\phi_j}$ with $\phi_j\in \mathbb{R}$, we have
\begin{equation}\label{ab}
  a^2-b^2 = 1- \frac{1}{|m_1||m_2|}\cos(\phi_1-\phi_2), \qquad 2ab = -\frac{1}{|m_1||m_2|}\sin(\phi_1-\phi_2).
\end{equation}
To prove the first inequality in~\eqref{zeta},
assume to the contrary that $|b|\le c_5\sqrt{\eta}$ with some small constant $c_5$.
Then from the second equation in~\eqref{ab}
$$
 |\sin (\phi_1-\phi_2)|^2 \le 10  c_5^2\eta,
$$
and plugging this into the first equation  and using $|m_i|\le 1-c\eta_i\le 1-c\eta$ 
with some universal constant $c$ from~\eqref{eq:mbound},  
we have
$$
   a^2 \le c_5^2\eta + 1 - \frac{1}{(1-c\eta)^2}\sqrt{1- 10c_5^2\eta} \le c_5^2\eta  -(1+c\eta)^2(1-10c_5^2 \eta)
   \le (15c_5^2-2c)\eta,
  $$
which is a contradiction if $c_5$ is sufficiently small so that $15c_5^2<2c$.  This proves
the first inequality in~\eqref{zeta}. The proof of the second is analogous.

We start with proving the decay of $\wh R(x)$. Two integration by parts in \eqref{Rhat} yields
\begin{equation}\label{Rdec}
   |\wh R(x)|\lesssim \biggl( \frac{W}{|x|}\biggr)^2 \frac{1}{W} \int_0^{2\pi W} \biggl[
   \frac{ |\partial^2F(q)|}{|\zeta - F(q)|^2}  +  \frac{ |\partial F(q)|^2}{|\zeta - F(q)|^3}\biggr]  \mathrm{d} q, \qquad x\in \mathbb{Z}.
\end{equation}
The integration domain is split at $c_\epsilon$ for some sufficiently small $\epsilon$ (depending only on $D_0$) so that
the error term in \eqref{smallq} is negligible compared with the quadratic term $D_0|q|^2$.
For $|q|\ge c_\epsilon$, the denominators in~\eqref{Rdec} are harmless since
$$
 |\zeta-F(q)|\ge \re (\zeta - F(q)) = \re (\zeta-1) + 1-\re F(q)\ge  c_\epsilon/2,
$$
by \eqref{largeq} and $|\zeta-1|\le c'\le c_\epsilon/2$. From this regime we
obtain the  fast decay $W^{-1} (W/|x|)^2$ for $\wh R(x)$,
exactly as in Case 1.

 In the complementary regime, $|q|\le c_\epsilon$, we use~\eqref{smallq} and obtain the estimate
\begin{equation*}
	\begin{split}
		&\biggl(\frac{W}{|x|}\biggr)^2 \frac{1}{W} \int_{-c_\epsilon}^{c_\epsilon} \biggl[
		\frac{ |\partial^2F(q)|}{|\zeta - 1 + D_0|q|^2(1+o(1)) |^2}  +  \frac{ |\partial F(q)|^2}{|\zeta - 1 + D_0|q|^2(1+o(1))|^3}\biggr]  \mathrm{d} q \\
		&\qquad\lesssim \biggl( \frac{W}{|x|}\biggr)^2 \frac{1}{W} \frac{1}{\eta^{3/2}}
		\lesssim \frac{1}{\ell\eta}
		\biggl( \frac{\ell}{|x|}\biggr)^2.
	\end{split}
\end{equation*} 
 In estimating the integrals, we used~\eqref{zetaF},  $|\zeta-1|\gtrsim \eta$ from~\eqref{zeta},  and
 that $|\partial^2F(q)|\lesssim 1$ and $|\partial F(q)|\lesssim |q|$ in the small $|q|$ regime.
 The main contribution of order $(\ell\eta)^{-1}(\ell/|x|)^2$ to~\eqref{Rdec} comes from the regime $|q|\le c_\epsilon$,
 the estimate $W^{-1} (W/|x|)^2$ that we obtained from the large $|q|$ regime is smaller.
 Repeating the integration by parts $D$-times we eventually obtain a decay of the form
 $$
    |\wh R(x)|\lesssim \frac{1}{\ell\eta} \biggl( \frac{\ell}{|x|}\biggr)^D.
$$
The trivial bound
$$
     |\wh R(x)|\lesssim \frac{1}{\ell\eta} 
$$
may be obtained directly from \eqref{Rhat} without integration by parts just by splitting the integration regime
as above. Again, the main contribution comes from the small $|q|$ regime.
Combining these last two bounds gives $|\wh R(x)|\le (\ell\eta)^{-1} \langle |x|/\ell\rangle^{-D}$, and
 inserting it into~\eqref{Resrep} shows that $|\Theta_{0x}|\lesssim (\Upsilon_\eta)_{0x}$,
 with $\Upsilon$ given in~\eqref{eq:polyUps_notime}. 
 
 If $f$ is compactly supported, hence $F$ is analytic in a strip $\Omega$ from~\eqref{strip}, then we compute
 the behavior of~\eqref{Rhat} by residue. Using~\eqref{smallq}--\eqref{largeq} it is easy to see  that
 the equation $\zeta=F(q)$ has exactly two solutions, $q_\pm \sim \pm \sqrt{1-\zeta}$, in the strip $\Omega$
  if $c_*$ in the definition~\eqref{strip} is sufficiently small (compared with $c_\epsilon$)
 and  we can write
 $$
   \zeta -F(q) = (q-q_+)(q-q_-) h(q),
 $$
 where $h(q)$ is an analytic function that does not vanish and $|h(q)|\ge c$ with some small constant $c>0$
 in  $\Omega$. Moreover
 $|\im q_\pm|\gtrsim \sqrt{\eta}$ from~\eqref{zeta}.  
 Now we can compute the integral in \eqref{Rhat} by shifting the contour over the residue
 $q_-$ or $q_+$ depending on the sign of $x$, but still staying within $\Omega$.
  The contribution of the residue is bounded by
   $\eta^{-1/2}W^{-1}\exp(-c\sqrt{\eta} x/W)$ where we used $|q_+-q_-|\ge |\im q_+| \gtrsim \sqrt{\eta}$
   to control the contribution of the denominator at the residue.
    The integral on the shifted horizontal contour will be even smaller
   since it is exponentially small with a larger rate constant $c$.  Plugging these estimates into~\eqref{Resrep}, 
   we obtain $|\Theta_{0x}|\lesssim (\Upsilon_\eta)_{0x}$ with $\Upsilon$ from~\eqref{eq:expUps_notime}.

Thus, in both cases we verified
$(i)$ of Definition~\ref{def:adm_ups_notime} for $\Theta$.

For the proof of $(vi)$ of Definition~\ref{def:adm_ups_notime}, it is sufficient to prove that
\begin{equation}\label{1step}
  \big| (\Theta_\eta)_{0,x+1} - (\Theta_\eta)_{0x}\big| \le \frac{1}{\ell}  (\Upsilon_\eta)_{0x},  \quad x\in \mathbb{Z},
\end{equation}
then by a telescopic sum we obtain~\eqref{eq:Theta_regularity_notime}. Given~\eqref{Resrep},
we need to control $\wh R(x+1)-\wh R(x)$, and using~\eqref{Rhat} we have
\begin{equation}\label{Rxx}
     \bigl\lvert\wh R(x+1)-\wh R(x)\bigr\rvert\lesssim \frac{1}{W} 
       \biggl\lvert \int_0^{2\pi W}  (1-e^{iq/W}) \frac{F(q) }{\zeta - F(q)}  e^{-iqx/W} \mathrm{d} q  \biggr\rvert.
\end{equation}
We follow the same calculation as for the bound $|\wh R(x)|$; it is easy to see that the effect of the 
additional factor $(1-e^{iq/W})\sim |q|/W$ is the additional $1/\ell \ge \sqrt{\eta}/W$ factor. For example, 
in the more critical  $|q|\le c_\epsilon$ regime, we estimate
$$
\biggl\lvert \frac{1}{W}\int_{-c_\epsilon}^{c_\epsilon}  (1-e^{iq/W}) \frac{F(q) }{\zeta - F(q)}  e^{-iqx/W} \mathrm{d} q \biggr\rvert
\lesssim \frac{1}{W}\int_{-c_\epsilon}^{c_\epsilon}  \frac{|q|}{W}  \frac{ |F(q)|}{|\zeta - 1 + D_0|q|^2(1+o(1))|}  
\mathrm{d} q  \lesssim \frac{1}{W^2} \lesssim \frac{1}{\ell}\frac{1}{\ell\eta},
$$
since $|F(q)|\le 1$ and the denominator is bounded from below by $\eta+|q|^2$. 
To obtain the polynomial decay, we perform integration by parts as for the bound on $|\wh R(x)|$,
and we pick up the extra $1/\ell$ from the presence of $(1-e^{iq/W})$ in the same way.
This shows~\eqref{1step} in the polynomially decaying case.

If $f$ is compactly supported,
we again use the residue calculation, the additional factor $(1-e^{iq/W})$ yields  a factor of size
$$
|1-e^{iq_\pm/W}| \lesssim\frac{ |q_\pm|}{W } \lesssim \frac{ |\sqrt{1-\zeta}|}{W }
$$
at the residue. The denominator in the residue calculation is again $|q_+-q_-|\ge |\im q_+|\sim |\im \sqrt{1-\zeta}|$.
 Using the second inequality in~\eqref{zeta}, the ratio of these two factors picked up from the residue,
 $|\sqrt{1-\zeta}|/|\im \sqrt{1-\zeta}|$ is bounded. This yields an estimate of order $1/W^2$ times the 
 exponentially decaying factor for~\eqref{Rxx}. Compared with the bound $\eta^{-1/2}W^{-1}\exp(-c\sqrt{\eta} x/W)$ 
 obtained along the previous calculation of $|\wh R(x)|$ itself, we see an additional gain of order $\sqrt{\eta}/W\le 1/\ell$.
 This proves~\eqref{1step} in the exponentially decaying case as well.

Finally the proof of $(vii)$ of Definition~\ref{def:adm_ups_notime} is very similar. For any $x\in\mathbb{Z}$, we use that
\begin{equation*}
	\begin{split}
		\biggl\lvert\wh R(x) -\frac{1}{N}\sum_{y=1}^N \wh R(x+y)\biggr\rvert &\le \max_y   \bigl\lvert\wh R(x) - \wh R(x+y)\bigr\rvert\\
		&\le \max_y \frac{1}{W}  \biggl\lvert \int_0^{2\pi W} \frac{F(q) }{\zeta - F(q)} e^{-iqx/W} (e^{-iqy/W}-1) \mathrm{d} q\biggr\rvert
	\end{split}
\end{equation*}
from~\eqref{Rhat}. We again follow the same calculation as for the bound $|\wh R(x)|$ and track the 
effect of the additional factor $(e^{-iqy/W}-1)$. In the more critical  $|q|\le c_\epsilon$ regime, we have
\begin{equation*}
	\begin{split}
		\biggl\lvert  \frac{1}{W}\int_{-c_\epsilon}^{c_\epsilon}   \frac{F(q) }{\zeta - F(q)} e^{-iqx/W} (e^{-iqy/W}-1)  \mathrm{d} q \biggr\rvert
		&\lesssim \frac{1}{W}\int_{-c_\epsilon}^{c_\epsilon}  \frac{|q||y|}{W}  \frac{ |F(q)|}{|\zeta - 1 + D_0|q|^2(1+o(1))|}  
		\mathrm{d} q\\
		&\lesssim \frac{|y|}{W^2} \lesssim \frac{N}{W^2},
	\end{split}
\end{equation*}
since $|y|\le N$.  The other regime gives the same bound and is even easier.
 Using \eqref{Resrep}, this would immediately show~\eqref{eq:supercrit_Theta_notime}
if the $k\ne 0$ terms in the summation were ignored. 
In order to ensure the summability in \eqref{Resrep}, we need to recover a summable decay in $k$ as well
when $x$ is replaced with $Nk+x$. This is achieved via an integration by parts, for example, in the 
$|q|\le c_\epsilon$ regime, when the derivative hits $F/(\zeta-F)$, it gives a bound 
$$
\frac{1}{W}\int_{-c_\epsilon}^{c_\epsilon}  \frac{|q||y|}{W}  \frac{W}{|Nk+x|} \frac{ |q|}{|\zeta - 1 + D_0|q|^2(1+o(1))|^2}  
\mathrm{d} q  \lesssim  \frac{1}{|k|} \frac{1}{W\sqrt{\eta}} \le \frac{1}{|k|} \Big(\frac{N}{W^2}\Big)^{1/2},
$$
using $|F'(q)|\lesssim |q|$ from~\eqref{Fprop} and $\eta\ge 1/N$.  
When the derivative  hits $(e^{-iqy/W}-1) $ we have 
$$
\frac{1}{W}\int_{-c_\epsilon}^{c_\epsilon}  \frac{|y|}{W}  \frac{W}{|Nk+x|} \frac{1}{|\zeta - 1 + D_0|q|^2(1+o(1))|}  
\mathrm{d} q  \lesssim  \frac{1}{|k|} \frac{1}{W\sqrt{\eta}} \le \frac{1}{|k|} \Big(\frac{N}{W^2}\Big)^{1/2},
$$
yielding the same bound.
A second integration by parts achieves
the square of this bound, rendering the $k$ summation in \eqref{Resrep} finite and
resulting in  the desired $(N/W^2)$ factor. 

The estimates in the other integration regime $|q|\ge c_\epsilon$ are similar but easier.
This completes the proof of~\eqref{eq:supercrit_Theta_notime}, and thus 
the proof of Proposition~\ref{prop:admS} (i).
\end{proof}

\printbibliography
\end{document}